# Crystalline prisms: reflections and diffractions, present and past

Arthur Ogus

March 9, 2026

## Contents

---

Dedicated to my great friend and colleague, Pierre Berthelot

**Abstract**

Our goal in this manscript is to render the theory of (crystalline) prismatic cohomology explicit by showing that it can be computed (and even defined) in terms of a natural $p$-de Rham complex associated to geometrically constructed prismatic envelopes. For example, let $Y/S$ be a $p$-completely smooth morphism of "crystalline prisms," *i.e.*, $p$-torsion free $p$-adic formal schemes endowed with a Frobenius lift, and let $\overline{Y}/\overline{S}$ denote its reduction modulo $p$. We show that the category of crystals on the prismatic site of $\overline{Y}/S$ is equivalent to the category of $\mathcal{O}_Y$-modules with integrable and quasi-nilpotent $p$-connection, and that the cohomology of such a crystal is computed by the associated $p$-de Rham complex. More generally, if $X$ is a closed subscheme of $\overline{Y}$, smooth over $\overline{S}$, then the prismatic envelope $\mathbb{\Delta}_X(Y)$ of $X$ in $Y$ admits such a $p$-connection, the category of prismatic crystals on $X/S$ is equivalent to the category of $\mathcal{O}_{\mathbb{\Delta}_X(Y)}$-modules with compatible integrable and quasi-nilpotent $p$-connection, and the cohomology of such a crystal is again computed by its $p$-de Rham complex. We also give a geometric construction of the "prismatic Sen operator." Namely, we show that a lifting of $X$ (mod $p^2$) in $Y$ defines a vector field on the reduction modulo $p$ of $\mathbb{\Delta}_X(Y)$ and on a "diffracted" Higgs complex which calculates the mod $p$ prismatic and de Rham cohomologies of $X$. Surprisingly, this complex is not the reduction modulo $p$ of the afore-mentioned $p$-de Rham complex but is rather its "$\alpha$-transform." As a consequence, we get a fairly explicit description of the action of the group scheme $\mathbf{G}_m^\gamma$ on $\mathrm{R}\Gamma(X, \Omega^\bullet_{X/S})$, Drinfeld's strengthening of the Deligne-Illusie decomposition theorem. We also explain how earlier work by several authors relating Higgs fields, $p$-connections, and connections can be placed in the prismatic context.

# 0 Introduction

## 0.1 Background and motivation

This article was originally inspired by a conversation with Bhargav Bhatt in the fall of 2019 in which he explained that, if $X/k$ is a smooth scheme over a perfect field $k$ which admits a formal lifting $Y/W$ along with its Frobenius endomorphism, then the prismatic cohomology of $X$ relative to the prismatic base $(W, (p))$ can be computed as the hypercohomology of the "$p$-de Rham complex" $(\Omega^\bullet_{Y/W}, pd)$. In a further email exchange, Bhatt agreed with my speculation that, in this lifted situation, the category of prismatic crystals on $X/W$ should be equivalent to the category of $\mathcal{O}_Y$-modules with integrable and quasi-nilpotent $p$-connection. These results resonate with previous work on the Cartier transform due to several authors, for example in [34], [37], [35], and [40]. Unfortunately, the foundational manuscript on prismatic cohomology [10] neither mentions nor proves these statements. Our purpose here is to confirm them and to elucidate the relations between the prismatic theory and the afore-mentioned previous work. We discuss only crystalline prisms, omitting what is undoubtedly the most interesting aspect of the theory; it remains to be seen if any of our techniques are applicable more generally.

Although striking, the comparison between prismatic and $p$-de Rham cohomology described above is not widely applicable, since Frobenius liftings rarely exist globally. However, one can often embed $X/k$ as a closed subscheme of a smooth $Y/W$ which does admit a Frobenius lift and then hope to endow the prismatic envelope of $X$ in $Y$ with a $p$-de Rham complex which calculates the prismatic cohomology of $X/W$. This is indeed possible and leads us to reconsider the old-fashioned idea that the geometry and cohomology of a general $X$ can be understood by studying its tubular neighborhoods in a smooth $Y$. For example, Hartshorne [18] used this method to *define* the de Rham cohomology of a possibly singular variety $X$ over a field $k$ of characteristic zero by showing that the de Rham cohomology of the formal completion of $Y$ along $X$ is independent of the choice of $Y$. Similarly, if $X/k$ is a scheme over a perfect field $k$ of characteristic $p$, one can compute (and even define) the crystalline cohomology of $X/W$ by taking the de Rham complex of the (canonical connection on the) divided power envelope of $X$ in a smooth $Y/W$ [4, 20, 7]. One can similarly obtain the convergent cohomology of $X/K$ by replacing the divided power envelope by the open $p$-adic analytic tube (of radius $< 1$) [31, 5]; this works well even if $X/k$ is singular.

We shall discuss an analogous approach for prismatic cohomology. Let $S$ be $\operatorname{Spf} W$, or, more generally a "crystalline prism," *i.e.*, a $p$-torsion free $p$-adic formal scheme equipped with a lift $F_S$ of the absolute Frobenius morphism of $\overline{S}$, its reduction modulo $p$. Let $X/\overline{S}$ be a smooth scheme, embedded in a $p$-completely smooth crystalline prism $Y/S$. We shall see in Proposition 4.2.1 that the structure sheaf of the "prismatic envelope" $\mathbb{A}_X(Y)$ of $X$ in $Y$, viewed as a sheaf of $\mathcal{O}_Y$-modules, is endowed with a canonical integrable and quasi-nilpotent $p$-connection, and in Theorem 4.3.1 that the derived image of the associated $p$-de Rham complex is functorial in $X$. These results are enough to establish the existence of a functorial cohomology theory based on $p$-de Rham cohomology, without relying on the formalism of Grothendieck topologies. Since that formalism has important conceptual advantages, we go on in Theorem 6.3.5 to show that these complexes do calculate the cohomology of the structure sheaf of the prismatic topos of $X/S$. Similar results hold with coefficients: a crystal $E$ on the prismatic site of $X/S$ can be evaluated on $\mathbb{A}_X(Y)$, the resulting sheaf of $\mathcal{O}_Y$-modules has an integrable quasi-nilpotent $p$-connection which is compatible with the $p$-connection on $\mathbb{A}_X(Y)$, and the resulting $p$-de Rham complex calculates the cohomology of $E$.

After the original version of the current manuscript was submitted for consideration for publication, we learned that Drinfeld's "stacky" approach to prismatic cohomology [15] leads to a fundamentally new strengthening of the Deligne-Illusie de Rham decomposition theorem [14]. Deligne-Illusie showed that a lifting of $X$ mod $p^2$ gives rise to a direct sum decomposition (in the derived category):

$$\tau^{<p}\mathrm{R}(\Omega^\bullet_{X/k}) \sim \bigoplus_0^{p-1} \mathcal{H}^i(\Omega^\bullet_{X/k});$$

Drinfeld's strengthening gives an action of the divided power envelope $\mathbf{G}_m^\gamma$ of the multiplicative group on the derived image of the de Rham complex of $X/k$. This in turn gives an action of $\mu_p$ on all of $\mathrm{R}(\Omega^\bullet_{X/k})$, and hence a $\mathbf{Z}/p\mathbf{Z}$-grading, whose truncation agrees with that of Deligne-Illusie. Upon learning of this result, we found it impossible to resist trying to realize this $\mathbf{G}_m^\gamma$-action explicitly on the reduction modulo $p$ of prismatic envelopes and on the corresponding complexes. Preliminary investigations led to a distressing conundrum: the simplest examples showed that there could be no reasonable action on these complexes, despite the compelling geometric hope for an action on the envelopes. Finally, after several extensive and revelatory discussions with V. Vologodsky and B. Bhatt, this conundrum has been resolved. It turns out that there is a natural geometric action of $\mathbf{G}_m^\gamma$ on prismatic envelopes, and that this action extends not to the original complex, but rather to its "$\alpha$-transform," a complex which also computes prismatic (and hence also de Rham) cohomology.

Although this may all sound comfortable and reassuring to old-timers, there are new technical difficulties. First of all, an explicit description of prismatic envelopes seems difficult to come by. In [10], these envelopes are constructed using universal constructions in the category of $\delta$-rings, but, at least to this author, the geometric structure of these envelopes remains obscure. We try to address this difficulty in Theorem 2.3.3 by giving explicit formulas in local coordinates and in Theorem 2.3.8, by constructing crystalline prismatic envelopes as a limit of a sequence of "$p$-adic dilatations," which were already used in [31], [35], and [40]. A second difficulty arises from the fact that morphisms in the prismatic site must be compatible with the chosen Frobenius lifting. This makes the task of understanding coverings of the final object of the prismatic topos more complicated, an issue we explore in §6.1. The coverings constructed in the foundational paper [10], which have the advantage of working in many topologies, seem unwieldy. However, if one works in the $p$-completely flat topology, such coverings can be constructed by simply taking prismatic envelopes of embeddings in $p$-completely smooth crystalline prisms. This fact, a generalization of a technique due to Morrow and Tsuji [28], is based on a ($p$-completely flat) local lifting property explained in Theorem 1.7.11. [1] Also important is Theorem 2.3.13, which implies that the prismatic envelope of a closed immersion of crystalline prisms can be described as the PD-envelope of an associated dilatation. (A related result, [10, 2.38], plays a key role in the crystalline to prismatic comparison theorem.)

The use of extremely large "$p$-completely faithfully flat coverings" seems unavoidable. We will therefore of necessity be considering highly non-noetherian formal schemes, requiring some modifications of usual notions and techniques. We will follow the general methods suggested in [10] and [16], which we review, with some additional details and examples, in the appendix §7. Since we shall be working almost exclusively with $p$-adic formal schemes, the notion of $p$-complete flatness becomes considerably simpler than the more general notion of

[1] I learned in June 2021 at a conference in honor of Luc Illusie that this construction of coverings was independently discovered by Yichao Tian in a closely related work [39].

$I$-complete flatness, and we shall not need to consider derived completeness.

## 0.2 Overview

Our manuscript begins with a brief review of connections, $p$-connections, Higgs fields, and their associated stratifications. These can be viewed as groupoid actions and also as actions of suitable rings of differential operators. We then discuss Frobenius liftings and $\delta$-structures, reviewing the relationships among divided power structures, $\delta$-structures, and the Witt functor, especially Drinfeld's isomorphism $\mathbf{W}_F \cong \mathbf{A}^\gamma$ (Proposition 1.6.6). We then investigate deformations of subschemes of crystalline prisms, the main outtake being that the category of formal schemes endowed with a Frobenius lift is quite rigid, in that it is difficult to find morphisms which are compatible with the Frobenius lifts, even locally in the Zariski or étale topology. The most important positive result in this section is Theorem 1.7.11, which shows that this difficulty can be overcome by working in the $p$-completely flat topology.

The geometric heart of our work occurs in section §2, where we discuss various kinds of $p$-adic tubular neighborhoods of a characteristic $p$ subscheme $X$ of a $p$-adic formal scheme $Y$. In decreasing order of size, these are the formal completion $\mathbb{F}_X(Y)$ of $Y$ along $X$, the $p$-torsion free version $\mathbb{PD}_X(Y)$ of the divided power envelope of $X$ in $Y$, the dilatation $\mathbb{D}_X(Y)$ and the dilation $\mathbb{V}_X(Y)$ of $X$ in $Y$, and, if $Y$ has a suitable Frobenius lift, the prismatic envelope $\mathbb{\Delta}_X(Y)$ of $X$ in $Y$. We attempt to describe these, and especially the last two, as explicitly as possible. The dilation and dilatation agree for regular closed immersions, and we have found the formalism of dilations more convenient for our purposes. For example, $X$-valued points of $\mathbb{V}_X(Y)/Y$ correspond to liftings of $X$ in $Y$ which are flat mod $p^2$, and the restriction $\overline{\mathbb{V}}_X(Y)/X$ of $\mathbb{V}_X(Y)/Y$ to $X$ has a natural structure of a torsor under the normal bundle of $X$ in $\overline{Y}$, the reduction of $Y$ modulo $p$. We give two constructions of prismatic envelopes for regular immersions: first an explicit description in local coordinates, then a more conceptual construction as a sequence of dilations. The latter construction works more generally: $Y$ need only be endowed with an endomorphism whose restriction to $X$ is the Frobenius endomorphism. On the other hand, the former construction allows us to prove a strong universal property of prismatic envelopes which then allows a description (Proposition 2.3.7) of the functor of points of their reduction modulo $p$. Expanding on results of [10], we go on to explain some relationships between prismatic and PD envelopes. One key fact is that if $Y$ is a crystalline prism, then the scheme theoretic image $F_{\overline{Y}}(\overline{Y})$ of the Frobenius endomorphism of $\overline{Y}$ is the smallest subscheme of $Y$ whose ideal has divided powers (Proposition 2.1.4). Another is that if $X$ is regularly immersed in $\overline{Y}$, then the map $\overline{\mathbb{\Delta}}_X(Y) \to X$ from the reduction of $\mathbb{\Delta}_X(Y)$ modulo $p$ to $X$ is faithfully flat (but not of finite presentation) (see Theorem 2.3.8). Also key is the fact (Theorem 2.3.13) that, if $X \to Y$ admits a mod $p^2$ $F_Y$-invariant lift $\tilde{X} \to Y$, then $\mathbb{\Delta}_X(Y)$ identifies with the PD-envelope of the corresponding section of the dilation of $X$ in $Y$.

Section §3 explains our approach to diffraction and the Sen operator. We explain two constructions, the first based on the existence of a natural connection $\nabla$ on $\overline{\mathbb{\Delta}}_X(Y)$, viewed as a scheme over $\overline{\mathbb{D}}_X(Y)/X$, that we discovered by "fooling around," and another more geometric construction, explained to us by Bhatt, using Proposition 2.3.7 and Drinfeld's isomorphism (Proposition 1.6.6) between the divided power envelope of the origin of the affine line and the kernel of the Frobenius endomorphism of the Witt scheme. As Bhatt explained to us, the restriction of the torsorial action of $\mathbb{N}_{X/\overline{Y}}$ on $\overline{\mathbb{D}}_X(Y)$ to the divided power envelope $\mathbb{N}^{\gamma}_{X/\overline{Y}}$ of its zero section lifts to a torsorial action on $\overline{\mathbb{\Delta}}_X(Y)$. In fact, these two constructions are equivalent: the $\mathbb{N}^{\gamma}_{X/\overline{Y}}$ action coincides with the stratification corresponding to the connection $\nabla$. It appears that, at least with our framework, both methods involve some rather heavy computations. Once these are accepted, we can understand "diffraction" as follows. A mod $p^2$-lifting $\tilde{X}$ of $X$ in $Y$ trivializes the torsor $\overline{\mathbb{D}}_X(Y)/X$, thus identifying it with $\mathbb{N}_{X/\overline{Y}}$. Then the Sen operator on $\overline{\mathbb{\Delta}}_X(Y)$ is obtained by lifting the Euler vector field $\theta$ of this vector bundle via the connection $\nabla$, and the (equivalent) $\mathbf{G}^{\gamma}_m$ action is obtained by lifting the restriction of the $\mathbf{G}_m$ action corresponding to $\theta$. An amusing and useful fact is that the quasi-nilpotent part of the Sen operator and the value of the $p$-curvature of $\nabla$ at $\theta$ are mutually determinative. Examples revealed a surprising difficulty: the Sen operator is not compatible with the Higgs field on $\overline{\mathbb{\Delta}}_X(Y)$ obtained as the mod $p$-reduction of its natural $p$-connection. However, via the identification of $\Omega^1_{\mathbb{D}/X}$ with $I_{X/\overline{Y}}/I^2_{X/\overline{Y}}$ and the embedding $I_{X/\overline{Y}}/I^2_{X/\overline{Y}} \to \Omega^1_{Y/S_{|X}}$, the connection $\nabla$ on $\mathcal{O}_{\overline{\mathbb{\Delta}}_X(Y)}$ relative to $\overline{\mathbb{D}}_X(Y)$ becomes a Higgs field $\nabla'$ relative to $\overline{Y}$. The Sen operator *is* compatible with this Higgs field, and we show that it acts as multiplication by $i$ on the $i$th cohomology sheaves of the associated Higgs complex, [2] as expected. Fortunately, we are able to show in a subsequent section that this Higgs complex also calculates prismatic cohomology, unfortunately our comparison remains rather indirect. We attempt to illustrate these constructions by giving some computations which we hope shed some light on diffraction. In particular, we show in Corollary 3.4.9 that the Sen operator on $\overline{\mathbb{\Delta}}_X(Y)$ is semi-simple, equivalently, the $\mathbf{G}^{\gamma}_m$-action factors through $\mu_p$, if and only if $\tilde{X}$ is $F_Y$-invariant. Let us note that our description does not apply directly to the simplest case, when $Y$ is itself a lifting of $X$, since in that case the normal bundle is zero. (To remedy this, one can look instead at the embedding $X \to Y \times Y$, in which one obtains a diffraction along the tangent bundle.)

Section §4 is the cohomological heart of the paper. Key to crystalline cohomology is the fact that divided power envelopes carry natural connections, and we show that dilatations and prismatic envelopes carry natural $p$-connections. Then we formulate and prove Theorem 4.3.1, the prismatic Poincaré lemma, which shows that the $p$-de Rham complex of a prismatic envelope is independent of the choice of embedding, up to quasi-isomorphism. This also works for the $p$-de Rham complex of a module with integrable quasi-nilpotent $p$-connection, and we show that the category of these is, up to equivalence, also independent of the embedding. The key computation is already revealed by the case of a point $X = \mathrm{Spec}(k)$, with liftings $Y = \mathrm{Spf}(W)$ and $Z = \mathrm{Spf}(W[x]\hat{})$ (the $p$-adic completion of the polynomial algebra $W[x]$), with $\phi(x) = x^p$. Then the prismatic neighborhood of $X$ in $Z$ is the formal spectrum of the completed PD-algebra $W\langle t\rangle\hat{}$, where $x = pt$. Then $d't = dx$, the $p$-de Rham complex of $\mathbb{\Delta}_X(Z)$ identifies with the de Rham complex of the PD-algebra, and so the complex is acyclic. The prismatic Poincaré lemma follows by a reduction to what is essentially this case which requires localization in the $p$-completely flat topology. This section concludes with a description, using the language of groupoids and groupoid actions, of prismatic differential operators and stratifications, a prelude to our later treatment of the prismatic topos and prismatic crystals.

[2] It seems that an opposite sign convention has been adapted in the most recent expositions of diffraction; we have not adopted it here.

Section §5 is devoted to an explanation of what we call the "$F$," "$p$," and "$\alpha$-transforms." The F-transform has been used by many authors, beginning perhaps with Berthelot's ideas in [6], and then by others, including [17], [34], [35], [37], and [40]. We focus on Shiho's work [37] showing how a relative Frobenius lift $F\colon Y \to Y^{(p)}$ allows one to define an equivalence from the category of modules with integrable quasi-nilpotent $p$-connection on $Y^{(p)}$ to the category of modules with integrable quasi-nilpotent connection on $Y$. This can now be viewed as an equivalence from the category of prismatic crystals on $\overline{Y}^{(p)}/S$ to the category of divided power crystals on $\overline{Y}/S$. After reviewing Shiho's work, we prove Theorem 5.1.6, showing that the de Rham complex of the F-transform of a module with $p$-connection is quasi-isomorphic to the $p$-de Rham complex of the original $p$-connection. Theorem 5.1.8 shows that the F-transform takes prismatic envelopes with their $p$-connections to divided power envelopes with their connections. Then we review the $\alpha$-transform, a construction from [34], which, given a suitable Frobenius lift, takes Higgs fields to Higgs fields, and we show that the Higgs complex of a quasi-nilpotent Higgs field and that of its $\alpha$-transform become isomorphic in the derived category (although we have not found an explicit quasi-isomorphism between them). Given a Frobenius lift $F\colon Y \to Y^{(p)}$ as above, the F-transform defines a functor from the category of modules with Higgs fields on $\overline{Y}^{(p)}$ to the category of modules with integrable connection on $\overline{Y}$, and we show in Theorem 5.3.1 that the $p$-curvature of this connection is the negative of the $\alpha$-transform of the original Higgs field. (This is a generalization of a result in [34], since it does not require quasi-nilpotence of the Higgs field; our proof here is quite different.) Then we show in Theorem 5.4.3 that the $\alpha$-transform of the Higgs field obtained by reducing the canonical $p$-connection on prismatic envelopes agrees with the Higgs field $\nabla'$ constructed from the connection $\nabla$. The proof of this key result depends on a different construction of the $\alpha$-transform based on a $\delta$-structure variant of the PD-filtration familiar in crystalline cohomology. To relate all this to previous work, we explain in Theorem 5.5.2 a prismatic analog of the construction of the inverse Cartier transform [34, §2.4] [35, §1.5] and [40]: the prismatic envelope $\mathbb{\Delta}_\Gamma(Y \times_S Y^{(p)})$ of the graph of $F$ serves as a kernel for F-transform. In the prismatic context, this construction also works in the derived sense, which is not the case for the Cartier transform.

The last section is devoted to the prismatic site, its crystals, and cohomology. The key point is that, if one works in the $p$-completely flat topology, the prismatic envelope of a scheme $X$ in a completely smooth crystalline prism $Y/S$ is a covering of the final object in the prismatic topos of $X$. Using the prismatic Poincaré lemma, we show that the linearization of the $p$-de Rham complex of this envelope provides an acyclic resolution of the prismatic structure sheaf of $X$. A similar construction works for crystals. We should mention that technical difficulties seem to require the use of an auxiliary site, consisting of "small" prisms (see Definition 2.3.1). We also describe a variation of the prismatic site (based on what we call "$\phi$-prisms" or "PD-prisms"), which was inspired by work of Oyama and Xu and which gives a geometric (site-theoretic) interpretation of the F-transform. It shows that, if $S$ is a crystalline prism and $X/\overline{S}$ is smooth, then the prismatic topos of the Frobenius pullback $X^{(p)}$ of $X$ is equivalent to the PD-prismatic topos of $X/S$, which in turn is very closely related to the usual (big) crystalline topos. This gives a new "pure thought" proof of Shiho's theorem. We finish with a few applications, which are mostly illustrations of how our results can be used to give simple and explicit proofs of results that are already stated in [10]. These include the comparison between prismatic and crystalline cohomology, the fact that prismatic Frobenius is an isogeny, and the relationship between the cotangent complex and prismatic cohomology.

We end with an appendix, which the reader is invited to consult as the need arises, addressing technicalities we face when working with $p$-adic sheaves and modules without noetherian hypotheses. The point of view we adopt follows that of Drinfeld in [16], where he works with inverse systems of $p$-torsion modules. For example, the category of $p$-adically separated and complete abelian groups is not abelian, but it does form an exact subcategory of the category of such systems, which is abelian. Localization presents another problem: for example, if $A$ is a $p$-adically complete ring and $a$ is an element of $A$ then the $p$-adic completion of the localization of $A$ by $a$ may not be flat. This is overcome by the notion of "$p$-complete flatness," which for us just means $p$-torsion freeness and flatness modulo each power of $p$. There is a concomitant notion of "$p$-complete quasi-coherence," in the Zariski and $p$-completely flat topologies. We end with a some comments about regular sequences, including a quick proof of Micali's theorem on the Rees algebra of a regularly generated ideal which we found to be useful in our description of dilatations and prismatic envelopes.

## 0.3 A key example

Let us illustrate the main constructions, the main quandary, and the essentials of its resolution with an example that was a guide from the very beginning of our research. Let $Y$ be the crystalline prism $\mathrm{Spf}(W[x,y]\hat{})$, with

$$\phi(x) = x^p + py, \qquad \phi(y) = y^p,$$

and let $X := \mathrm{Spec}(k[y])$, with $X \subset Y$ defined by the ideal $(p, x)$. Then the dilatation $\mathbb{D}_X(Y)$ of $X$ in $Y$ is

$$\mathrm{Spf}(W[y,s]\hat{}\,) \to \mathrm{Spf}(W[y,x]) : x \mapsto ps,$$

and the prismatic envelope $\mathbb{\Delta}_X(Y)$ is

$$\mathrm{Spf}(W[s]\langle t\rangle\hat{}\,) \to \mathrm{Spf}(W[x,y]) : x \mapsto ps, y \mapsto s^p + pt.$$

(See Example 2.3.14 for a justification of this description.) The $p$-connections on $\mathcal{O}_Y$, $\mathcal{O}_{\mathbb{D}_X(Y)}$, and $\mathcal{O}_{\mathbb{\Delta}_X(Y)}$ are compatible, given by

$$\begin{aligned} d'x &= pdx \\ d'y &= pdy \\ d's &= dx \\ d't &= dy - s^{p-1}dx. \end{aligned}$$

Note that $dy$ and $s^{p-1}dx$ are homologous in the $p$-de Rham complex of $\mathcal{O}_{\mathbb{\Delta}_X(Y)}$.

There is also a prismatic lifting $Z := \mathrm{Spf}(W[y]\hat{}\,)$ of $X$, with $\phi(y) = y^p$, and an obvious morphism of prisms $Y \to Z$. The prismatic Poincaré lemma tells us that the corresponding map of $p$-de Rham complexes is a quasi-isomorphism, as can of course be checked directly. Reducing these complexes modulo $p$, we find Higgs complexes which calculate the cohomology of the reduction modulo $p$ of the structure sheaf of the prismatic site of $X/W$. Note that the embedding $Z \subset Y$ defined by the ideal $(x)$ does not define a morphism of prisms $Z \to Y$. However, the lifting $Z_2/W_2$ of $X$ is supposed to produce a diffraction of the prismatic cohomology of $X/k$, and since $Z_2$ is embedded in $Y$, we hoped to realize this diffraction as an action of $\mathbf{G}_m^\gamma$ on the Higgs complex $(\mathcal{O}_{\overline{\mathbb{\Delta}}_X(Y)}, d')$. Such an action is supposed to act via the negative of of the identity character on the first cohomology sheaf $\mathcal{H}^1$ of this complex , but we could find no reasonable such action compatible with the above formulas for the differential $d'$.

Despite this example, we proceeded with our attempt to construct the diffraction geometrically, first on $\overline{\mathbb{D}}_X(Y)$, then on $\overline{\mathbb{\Delta}}_X(Y)$. Recall that the lifting $Z_2 \to Y$ of $X$ in $Y$ defines an isomorphism $\overline{\mathbb{D}}_X(Y) \cong \mathbb{N}_{X/\overline{Y}}$, hence an action of $\mathbf{G}_m := \mathrm{Spec}(k[u, u^{-1}])$. In our example, the action is given by

$$r^*(s) = us, r^*(y) = y.$$

Then PD-parallel transport given by the connection $\nabla$ on $\overline{\mathbb{\Delta}}_X(Y)$, viewed over $\overline{\mathbb{D}}_X(Y)$, lifts the corresponding action $r$ of $\mathbf{G}_m^\gamma$ to $\overline{\mathbb{\Delta}}_X(Y)$. In our example, these turn out to be:

$$\begin{aligned} \nabla(s) &= ds \\ \nabla(t^{[n]}) &= -t^{[n-1]}s^{p-1}ds \\ r^*(s) &= us \\ r^*(t) &= t - y\log u \end{aligned}$$

(See Examples 3.2.10.) Thus:

$$\begin{aligned} d'(r^*(s)) &= d'(us) \\ &= udx \\ r^*(d's) &= r^*(dx) \\ &= udx, \end{aligned}$$

as desired, but:

$$\begin{aligned} d'(r^*(t)) &= d'(t - y \log u) \\ &= dy - s^{p-1}dx \\ r^*(d't) &= r^*(dy - s^{p-1}dx) \\ &= dy - u^{p-1}s^{p-1}udx \\ &= dy - s^{p-1}dx. \end{aligned}$$

Thus the presence of $dy$ destroys the compatibility with $d'$. However, if we use the Higgs field $\nabla'$ (the $\alpha$-tranform of $d'$), we find

$$\begin{aligned} \nabla'(s) &= dx \\ \nabla'(t) &= -s^{p-1}dx, \end{aligned}$$

and the difficulty disappears.

Since all constructions are natural and geometric, an explanation of the failure of the compatibility of diffraction with $d'$ is required. This naturality implies that our diffraction construction does indeed act on the prismatic Cech-Alexander complex. An explanation of why it does not act on the Higgs complex $(\mathcal{O}_{\overline{\mathbb{A}}_X(Y)}, d')$ can be seen in an even simpler example. Let $Y$ be the formal spectrum of $W[y]\hat{}$, with Frobenius lift given by $\phi(y) = y^p$, let $X$ be the reduction of $Y$ modulo $p$, and let $Y \to Y(1)$ be the diagonal embedding. Then $Y(1)$ is the formal spectrum of $W[y, \xi]\hat{}$, where we write $y$ for $y \otimes 1$ and $\xi$ for $1 \otimes y - y \otimes 1$, and $\mathbb{D}_X Y(1)$ is the formal spectrum of $W[y, \eta]\hat{}$, where $p\eta = \xi$. The lifting $Y \to Y(1)$ of $X \to Y(1)$ defines a section $Y \to \mathbb{D}_X(Y)$ defined by $\eta$, and as explained in Example 3.2.11, $\mathbb{A}_X(Y(1))$ is the PD-envelope of this section, *i.e.*: $\mathbb{A}_X(Y(1)) = \mathrm{Spf}(W[y]\langle\eta\rangle\hat{})$, It turns out that

$$\begin{aligned} r^*(\eta) &= u\eta \\ r^*(\eta^{[p]}) &= \eta^{[p]} + (1-u)y^{p-1}\eta. \end{aligned}$$

Thus the $\mathbf{G}_m^\gamma$-action preserves the ideal of $X$ in $\overline{\mathbb{A}}_X(Y(1))$. However it does not preserve the PD-filtration of this ideal and hence does not act naturally on $I_{X/\overline{Y}}/I^{[2]}_{X/\overline{Y}}$, and one needs the isomorphism

$$\Omega^1_{\overline{Y}/k} \cong I_{X/\overline{Y}}/I^{[2]}_{X/\overline{Y}}$$

to obtain a geometric interpretation of $d'$. This observation led us to introduce another filtration on $\mathcal{O}_{\overline{\mathbb{A}}_X(Y(1))}$, more closely related to the $\delta$-structure, and indeed this filtration gives another construction of the $\alpha$-transform. This is the construction that we use to prove that $\nabla'$ is the $\alpha$-transform of $d'$.

## 0.4 Acknowledgements

Enormous thanks go to Bhargav Bhatt, for the initial conversation which inspired this work, for his encouragement to carry it out, and for his frequent patient technical advice that helped clarify some of my initial difficulties. Thanks also go to Matthew Morrow, to Atsuki Shiho, to Vadim Vologodsky, to Ahmed Abbes, and especially to Luc Illusie, for additional advice and enlightenment. I must also acknowledge the contribution of two meticulous referees who pointed out numerous typographical errors, notational ambiguities, and outright mathematical mistakes in my submission of an early version of this manuscript. The current version is a complete revision of that one, attempting to incorporate their corrections and also the new results on diffraction that I obtained after its submission. Although many of the referees' corrections have been swept away by the extent of this revision, they were nevertheless invaluable. This last version was subject to an even more meticulous review, pointing out subtle punctuation and font errors, mislabeled references, unclear arguments, and several serious errors, often with suggestions on how to correct them. Finally, I wish to reiterate my indebtedness to Pierre Berthelot, whose crystal clear lectures on crystalline cohomology at Princeton in 1974 laid the foundations for my interest in $p$-adic cohomology, and whose commitment to clear exposition of fundamentals has remained an inspiration to me throughout my career.

# 1 Preliminaries

## 1.1 Notation and conventions

Throughout this paper we fix a perfect field $k$ and let $W$ be its ring of Witt vectors. If $Y/W$ is a scheme or formal scheme and $n$ is a natural number, we denote by $Y_n$ the subscheme (or formal subscheme) defined by $p^n$; we may also write $\overline{Y}$ for $Y_1$. If $X \to Y$ is a locally closed immersion of schemes or formal schemes, $I_{X/Y}$ is the ideal sheaf of $X$ in $Y$, $\mathcal{C}_{X/Y} := I_{X/Y}/I^2_{X/Y}$ is the conormal sheaf, and $\mathcal{N}_{X/Y} := \mathrm{Hom}(\mathcal{C}_{X/Y}, \mathcal{O}_X)$ is the normal sheaf. If $Y/S$ is a morphism of schemes (resp. of $p$-adic formal schemes), then $\Omega^1_{Y/S}$ is the sheaf of Kähler differentials (resp., its $p$-adic completion), which identifies with the conormal sheaf of the diagonal embedding

$$\Delta_Y : Y \to Y(1) := Y \times_S Y,$$

and

$$\mathcal{T}_{Y/S} := \mathcal{H}om_{\mathcal{O}_Y}(\Omega^1_{Y/S}, \mathcal{O}_Y)$$

is the sheaf of $S$-derivations of $\mathcal{O}_Y$. We let $\mathrm{P}_{Y/S}$ denote the first infinitesimal neighborhood of the diagonal, *i.e.*, the closed (formal) subscheme of $Y(1)$ defined by $I^2_{Y/Y(1)}$, so there is an exact sequence;

$$0 \to \Omega^1_{Y/S} \to \mathcal{O}_{\mathrm{P}_{Y/S}} \to \mathcal{O}_Y \to 0. \tag{1.1.1}$$

and, if $p_1, p_2\colon \mathrm{P}_{Y/S} \to Y$ are the two projections, then

$$d := p_2^* - p_1^*\colon \mathcal{O}_Y \to \Omega^1_{Y/S} \tag{1.1.2}$$

is the universal derivation.

If $E$ is a quasi-coherent sheaf on $Y$, then $\mathbf{V}E$ denotes the vector scheme $\operatorname{Spec}_Y S^\bullet E$, where $S^\bullet E$ is the symmetric algebra of $E$. In particular, the tangent space $\mathbb{T}_{Y/S}$ of $Y$ is $\mathbf{V}\Omega^1_{Y/S}$, the cotangent space $\mathbb{T}^*_{Y/S}$ is $\mathbf{V}\mathcal{T}_{Y/S}$, and $\mathbb{N}_{X/Y} := \mathbf{V}\mathcal{C}_{X/Y}$ is the normal bundle. If $E$ is locally free and of finite rank, then $E^\vee := \mathcal{H}om_{\mathcal{O}_Y}(E, \mathcal{O}_Y)$ is coherent, and we let $\mathbf{V}^*E := \mathbf{V}E^\vee$. In particular, if $X \to Y$ is a regular immersion, then $\mathbb{N}^*_{X/Y} := \mathbf{V}^*\mathcal{C}_{X/Y}$ is the conormal bundle.

We write $F_X$ for the absolute Frobenius endomorphism of a scheme $X$ in characteristic $p$, and if $f\colon X \to S$ is a morphism of schemes in characteristic $p$, we have a commutative diagram:

$$\begin{array}{ccccccc}
X & \xrightarrow{F_{X/S}} & X^{(p)} & \xrightarrow{\pi_{X/S}} & X & \xrightarrow{F_{X/S}} & X^{(p)} \\
 & \searrow^{f} & \downarrow{\scriptstyle f^{(p)}} & & \downarrow{\scriptstyle f} & \swarrow_{f^{(p)}} & \\
 & & S & \xrightarrow{F_S} & S, & &
\end{array} \tag{1.1.3}$$

in which the square is Cartesian and $\pi_{X/S} \circ F_{X/S} = F_X$ and $F_{X/S} \circ \pi_{X/S} = F_{X^{(p)}}$. The morphism $F_{X/S}$ is the relative Frobenius morphism. We may use a similar notation if $Y/S$ is a morphism of formal schemes endowed with compatible Frobenius lifts $F_Y$ and $F_S$.

If $C^\bullet$ is a complex in a suitable abelian category, we write $\mathrm{R}C^\bullet$ for the image of $C^\bullet$ in the associated derived category (whose defining characteristics we hope are understood from the context).

## 1.2 Connections and $p$-connections

Let us briefly review these basic notions. If $Y/S$ is a morphism of schemes or $p$-adic formal schemes, recall that the universal derivation $d\colon \mathcal{O}_Y \to \Omega^1_{Y/S}$ fits in the *de Rham complex*

$$(\Omega^\bullet_{Y/S}, d) := \mathcal{O}_Y \xrightarrow{d} \Omega^1_{Y/S} \xrightarrow{d} \Omega^2_{Y/S} \xrightarrow{d} \cdots . \tag{1.2.1}$$

Multiplying the differential by $p$, we find the *p-de Rham complex*

$$(\Omega^\bullet_{Y/S}, d') := \mathcal{O}_Y \xrightarrow{pd} \Omega^1_{Y/S} \xrightarrow{pd} \Omega^2_{Y/S} \xrightarrow{pd} \cdots . \tag{1.2.2}$$

Note that if $f$ and $g$ are sections of $\mathcal{O}_Y$, then

$$d'(fg) := pd(fg) = pfdg + pgdf = fd'g + gd'f,$$

*i.e.*, $d'$ satisfies the usual Leibniz rule.

Recall that a *connection* on a sheaf of $\mathcal{O}_Y$-modules $E$ is an additive homomorphism $\nabla\colon E \to E \otimes \Omega^1_{Y/S}$ such that $\nabla(fe) = f\nabla e + e \otimes df$ for $f \in \mathcal{O}_Y$ and $e \in E$, that such a $\nabla$ induces a collection of maps $\nabla\colon E \otimes \Omega^i_{Y/S} \to E \otimes \Omega^{i+1}_{Y/S}$ satisfying the (signed) Leibniz rule

$$\nabla(\alpha \wedge \beta) = \nabla(\alpha) \wedge \beta + (-1)^i \alpha \wedge d\beta,$$

if $\alpha \in E \otimes \Omega^i_{Y/S}$ and $\beta \in \Omega^j_{Y/S}$, and that $\nabla$ is said to be *integrable* if $\nabla^2 = 0$. In this case one can form the *de Rham complex of* $(E, \nabla)$:

$$(E \otimes \Omega^\bullet_{Y/S}, d) := E \xrightarrow{\nabla} E \otimes \Omega^1_{Y/S} \xrightarrow{\nabla} E \otimes \Omega^2_{Y/S} \xrightarrow{\nabla} \cdots .$$

A connection on $E$ gives an isomorphism of $\mathcal{O}_{\mathrm{P}_{Y/S}}$-modules:

$$p_2^*(E) \to p_1^*(E)\colon p_2^*(e) \mapsto p_1^*(e) + \nabla(e), \tag{1.2.3}$$

whose restriction to the diagonal is the identity, and such an isomorphism is equivalent to the data of a connection. If the connection is integrable, this isomorphism extends uniquely to an isomorphism between the two pullbacks of $E$ to the $p$-adic and PD-completed divided power envelope of the diagonal, and this isomorphism satisfies a cocycle condition [7, 4.8].

If $X/S$ is a smooth morphism in characteristic $p$, an integrable connection on $X/S$ has a *p-curvature* [23], which can be viewed as a linear map

$$\psi\colon E \to E \otimes F^*_{X/S}(\Omega^1_{X^{(p)}/S}). \tag{1.2.4}$$

We briefly recall the construction. For each $D \in \mathcal{T}_{X/S}$, $\nabla_D$ is an $\mathcal{O}_{X^{(p)}}$-linear endomorphism of $E$. The $p$th iterate $D^p$ of $D$ is again a derivation of $\mathcal{O}_X$, and it is easy to see that $\nabla_D^p - \nabla_{D^p}$ is in fact $\mathcal{O}_X$-linear A more difficult calculation shows that the map

$$\mathcal{T}_{X/S} \to End_{\mathcal{O}_X}(E) : D \to \nabla_D^p - \nabla_{D^p}$$

is $p$-linear and thus defines an $\mathcal{O}_X$-linear homomorphism

$$F^*_{X/S}(\mathcal{T}_{X^{(p)}/S}) \to End_{\mathcal{O}_X}(E) : 1 \otimes D \mapsto \nabla_D^p - \nabla_{D^p}. \tag{1.2.5}$$

Such a map is equivalent to a map of the form (1.2.4), and it is common to denote both it and (1.2.5) by the same letter $\psi$. One checks that the images of any two elements of $\mathcal{T}_{X^{(p)}/S}$ commute, so $\psi$ extends to a homomorphism from the symmetric algebra of $F^*_{X/S}(\mathcal{T}_{X^{(p)}/S})$. Thus, if $E$ is quasi-coherent, we can view the data $(E, \psi)$ as a quasi-coherent sheaf $E_\psi$ on $\mathbb{T}^*_{X^{(p)}} \times_{X^{(p)}} X$.

The map $\psi$ induces maps

$$E \otimes F^*_{X/S}(\Omega^i_{X^{(p)}/S}) \to E \otimes F^*_{X/S}(\Omega^{i+1}_{X^{(p)}/S})$$

which satisfy the integrability condition $\psi \circ \psi = 0$. The $p$-curvature is said to be *quasi-nilpotent* if for each local section $e$ of $E$ and each local section $D$ of $\mathcal{T}_{X'/S}$, $\psi_D^n(e) = 0$ for $n$ sufficiently large. An integrable connection is said to be quasi-nilpotent if its $p$-curvature is, and an integrable connection over a $p$-adic formal scheme is said to be quasi-nilpotent if its reduction modulo $p$ is.

Let us also recall the notion of a *p-connection.*

**Definition 1.2.1.** *If $Y/S$ is a p-completely smooth morphism of p-adic formal schemes and $E$ is a sheaf of $\mathcal{O}_Y$-modules, a p-connection on $E$ is an $\mathcal{O}_S$-linear map*

$$\nabla' \colon E \to E \otimes \Omega^1_{Y/S}$$

*such that $\nabla'(ae) = d'a \otimes e + a\nabla'(e)$ for all $a \in \mathcal{O}_Y$ and $e \in E$, where $d'a := pda$.*

A $p$-connection extends uniquely to a family of maps

$$\nabla' \colon E \otimes \Omega^i_{Y/S} \to E \otimes \Omega^{i+1}_{Y/S},$$

satisfying the signed Leibniz rule, and one has the evident notion of integrability. In the integrable case one can form the *p-de Rham complex* $(E' \otimes \Omega^\bullet_{Y/S}, d')$ of $(E, \nabla')$.

**Example 1.2.2** If $\nabla$ is a connection on an $\mathcal{O}_Y$-module $E$, then $p\nabla$ is a $p$-connection. We call $(E, p\nabla)$ the *p-transform* of $(E, \nabla)$. There is a natural morphism of complexes:

$$b \colon (E \otimes \Omega^\bullet_{Y/S}, d) \to (E \otimes \Omega^\bullet_{Y/S}, pd) \tag{1.2.6}$$

given by

$$\begin{array}{ccccccc}
E & \xrightarrow{d} & \Omega^1_{Y/S} \otimes E & \xrightarrow{d} & \Omega^2_{Y/S} \otimes E & \xrightarrow{d} & \cdots \\
\Big\downarrow{\scriptstyle \mathrm{id}} & & \Big\downarrow{\scriptstyle p} & & \Big\downarrow{\scriptstyle p^2} & & \cdots \\
E & \xrightarrow{pd} & \Omega^1_{Y/S} \otimes E & \xrightarrow{pd} & \Omega^2_{Y/S} \otimes E) & \xrightarrow{pd} & \cdots .
\end{array}$$

We note that this morphism is an isogeny if $Y$ is quasi-compact. In particular, if $Y/S$ has relative dimension $m$, there is a morphism of complexes:

$$\tilde{b} \colon (\Omega^\bullet_{Y/S} \otimes E, pd) \to (\Omega^\bullet_{Y/S} \otimes E, d), \tag{1.2.7}$$

given by:

$$\begin{array}{ccccccccc}
E & \xrightarrow{pd} & \Omega^1_{Y/S} \otimes E & \xrightarrow{pd} & \Omega^2_{Y/S} \otimes E & \xrightarrow{pd} & \cdots & \longrightarrow & \Omega^m_{Y/S} \otimes E \\
\Big\downarrow{\scriptstyle p^m} & & \Big\downarrow{\scriptstyle p^{m-1}} & & \Big\downarrow{\scriptstyle p^{m-2}} & & & & \Big\downarrow{\scriptstyle \mathrm{id}} \\
E & \xrightarrow{d} & \Omega^1_{Y/S} \otimes E & \xrightarrow{d} & \Omega^2_{Y/S} \otimes E & \xrightarrow{d} & \cdots & \longrightarrow & \Omega^m_{Y/S} \otimes E.
\end{array}$$

Evidently $b \circ \tilde{b}$ and $\tilde{b} \circ b$ are multiplication by $p^m$.

Let us also note that, if $E$ is $p$-torsion free, then $b$ factors through an isomorphism

$$(E \otimes \Omega^\bullet_{Y/S}, d) \cong L\eta(E \otimes \Omega^\bullet_{Y/S}, pd),$$

where we recall from [7, §8] and [9] that, for any $p$-torsion free complex $C^\bullet$ of abelian sheaves, $L\eta(C^\bullet)$ is the subcomplex given by

$$(L\eta(C^\bullet))^q := \{x \in p^q C^q : dx \in p^{q+1} C^{q+1}\}.$$

Indeed, since the differentials of the $p$-de Rham complex of $(E, pd)$ are divisible by $p$, $L\eta(E \otimes \Omega^\bullet_{Y/S}, pd)$ in degree $i$ is just $p^i E \otimes \Omega^i_{Y/S}$, the image of $b$ in that degree.

We shall give a geometric interpretation of this trivial-looking construction in Theorem 6.4.12.

We omit the proof of the following analog of the Gauss-Manin connection. (Follow the method of [24].)

**Proposition 1.2.3.** *Let $g \colon Y \to Z$ be a $p$-completely smooth morphism of $p$-completely smooth $p$-adic formal schemes over $S$ and $(E, \nabla')$ a sheaf of $\mathcal{O}_Y$-modules with integrable $p$-connection relative to $S$. Then the composition*

$$\nabla' \colon E \to E \otimes \Omega^1_{Y/S} \to E \otimes \Omega^1_{Y/Z}$$

*defines a $p$-connection relative to $Z$. Furthermore, the sheaves*

$$E^j := R^j g_*(E \otimes \Omega^\bullet_{Y/Z})$$

*carry a natural $p$-connection relative to $Z$, and there is a spectral sequence with*

$$E_2^{i,j} = H^i(Z, E^j \otimes \Omega^\bullet_{Z/S}) \Rightarrow H^{i+j}(Y, E \otimes \Omega^\bullet_{Y/S}). \quad \square$$

## 1.3 Higgs fields and Higgs modules

**Definition 1.3.1.** *Let $X/S$ be a smooth morphism of schemes and let $E$ be a sheaf of $\mathcal{O}_X$-modules. A Higgs field on $E$ is an $\mathcal{O}_X$-linear map*

$$\theta \colon E \to E \otimes \Omega^1_{X/S}$$

*such that the composed map*

$$E \xrightarrow{\theta} E \otimes \Omega^1_{X/S} \xrightarrow{\theta\wedge} E \otimes \Omega^2_{X/S}$$

*vanishes.*

This definition makes sense with any locally free and coherent sheaf $\Omega$ in place of $\Omega^1_{X/S}$, and it will be convenient for us to work in this more general context. We refer to such generalized data as an $\Omega$-*field* on $E$. For example, the $p$-curvature of an integrable connection in characteristic $p$ is an $F^*_{X/S}(\Omega^1_{X^{(p)}/S})$-field. Let $\mathcal{T} := \mathrm{Hom}_{\mathcal{O}_X}(\Omega, \mathcal{O}_X)$, and for $\tau \in \mathcal{T}$, let $\theta_\tau$ be the endomorphism of $E$ given by the composition of $\theta$ with the map $E \otimes \Omega \to E$ induced by pairing with $\tau$. The vanishing of the composed map in the definition is equivalent to the requirement that for every pair $(\tau, \tau')$ of local sections of $\mathcal{T}$, the endomorphisms $\theta_\tau$ and $\theta_{\tau'}$ commute. This in turn is equivalent to saying that the map $\mathcal{T} \to \mathit{End}_{\mathcal{O}_X}(E)$ induced by $\theta$ factors through an algebra homomorphism $S^\bullet\mathcal{T} \to \mathit{End}_{\mathcal{O}_X}(E)$. Thus, a $\Omega$-field on $E$ endows it with a structure of an $S^\bullet\mathcal{T}$-module, and conversely, the data of such a module is equivalent to the data of an $\mathcal{O}_X$-module endowed with a $\Omega$-field. If $E$ is quasi-coherent, the data of an $\Omega$-field on $E$ is thus equivalent to a lift of $E$ to a quasi-coherent sheaf $E_\theta$ on $\mathbf{V}\mathcal{T}$. An $\Omega$-field $\theta$ is said to be *quasi-nilpotent* if for every local section $\tau$ of $\mathcal{T}$, every local section of $E$ is annihilated by some power of $\tau$, equivalently, by $S^n\mathcal{T}$ for $n$ large enough, depending on $e$. This is the case if and only if the associated module $E_\theta$ is supported on the zero section of $\mathbf{V}T$.

The geometric interpretation of an $\Omega$-field on $E$ as a lifting of $E$ to $\mathbf{V}\mathcal{T}$ is appealing and useful, and duality provides another such interpretation. Recall from [7, A10] that the $\mathcal{O}_X$-linear dual of $S^n\mathcal{T}$ is $\Gamma^n\Omega$, and so the dual of $S^\bullet\mathcal{T}$ is

$$\hat{\Gamma}^\bullet\Omega := \prod_n \Gamma^n\Omega.$$

Thus the map $S^\bullet\mathcal{T} \otimes E \to E$ induces the map $r_\theta^*$ in the diagram

$$\begin{array}{ccc} E & \xrightarrow{r_\theta^*} & E \otimes \hat{\Gamma}^\bullet\Omega \\ & \searrow^{\theta} & \downarrow \\ & & E \otimes \Omega, \end{array} \tag{1.3.1}$$

which factors through $\Gamma^\bullet\Omega$ if the original $\Omega$-field is quasi-nilpotent. The spectrum of $\Gamma^\bullet\Omega$ (resp., of $\hat{\Gamma}^\bullet\Omega$) is the PD-envelope (resp., the completed PD-envelope) of the zero section of the vector bundle $\mathbf{V}\Omega$, which we denote by $\mathbf{V}^\gamma\Omega$ (resp. $\hat{\mathbf{V}}^\gamma\Omega$). These are group schemes, and the morphism $r_\theta^*$ (1.3.1) defines a coaction on $E$. If $(\omega_1, \ldots, \omega_r)$ is a basis for $\Omega$, with dual basis $(\tau_1, \ldots, \tau_r)$ for $\mathcal{T}$, then:

$$r_\theta^*(e) = \sum_I \theta_{\tau_1}^{I_1}\theta_{\tau_2}^{I_2} \cdots \theta_{\tau_r}^{I_r}(e) \otimes \omega_1^{[I_1]}\omega_2^{[I_2]} \cdots \omega_r^{[I_r]}. \tag{1.3.2}$$

To express this geometrically, observe that the action of of $\mathcal{T}$ on $E$ extends to an action of $\mathcal{T}$ on $S^\bullet E$ by derivations, given by the diagram

$$\begin{array}{ccc} S^\bullet E \otimes E & \xrightarrow{\mathrm{id} \otimes \theta} & S^\bullet E \otimes E \otimes \Omega \\ d \uparrow & & \downarrow \mathit{mult} \otimes \mathrm{id} \\ S^\bullet E & \longrightarrow & S^\bullet E \otimes \Omega \end{array}$$

in which $d\colon S^\bullet E \to E \otimes S^\bullet E$ is the universal derivation. Then the induced coaction of $\hat{\Gamma}^\bullet\Omega$ on $S^\bullet E$, is compatible with its ring structure, and hence corresponds to an action of $\hat{\mathbf{V}}^\gamma\Omega$ on the scheme $\mathbf{V}E$. We refer to [34, 5.19] for more details about the duality involved in these constructions.

We record this for future reference as follows, but omit details of the proof.

**Proposition 1.3.2.** *Let $\Omega$ be a locally free finite type sheaf of $\mathcal{O}_X$-modules on a scheme $X$ and let $\mathcal{T}$ be its dual.*

*(1) There is an equivalence of categories $(E,\theta) \to E_\theta$ from the category of quasi-coherent sheaves of $\mathcal{O}_X$-modules $E$ endowed with an $\Omega$-field $\theta$ to the category of quasi-coherent sheaves of $\mathcal{O}_{\mathbf{V}\mathcal{T}}$-modules. Under this equivalence, the $\Omega$-field $\theta$ is quasi-nilpotent if and only if the corresponding $\Omega$-module $E_\theta$ has support on the zero section of $\mathbf{V}\mathcal{T}$.*

*(2) If $\theta$ is a quasi-nilpotent $\Omega$-field on a quasi-coherent $\mathcal{O}_X$-module $E$, then the coaction 1.3.1 defines an action $r_\theta\colon \mathbf{V}E \times \mathbf{V}^\gamma\Omega \to \mathbf{V}E$ of the group scheme $\mathbf{V}^\gamma\Omega$ on $\mathbf{V}E$. This construction defines an equivalence from the category of quasi-nilpotent Higgs modules to the category of vector schemes with an action of $\mathbf{V}^\gamma\Omega$.*

The integrability condition of an $\Omega$-field $\theta$ on $E$ implies that we can form the *Higgs complex*:

$$(E \otimes \Omega^\bullet, \theta) := E \xrightarrow{\theta} E \otimes \Omega^1 \xrightarrow{\theta\wedge} E \otimes \Omega^2 \longrightarrow \cdots, \tag{1.3.3}$$

where $\Omega^i := \Lambda^i\Omega$. This complex has the following geometric interpretation. For the proof, we refer to [33, 2.1.2].

**Proposition 1.3.3.** *If $\theta$ is a quasi-coherent $\mathcal{O}_X$-module with $\Omega$-field $\theta$ and if $E_\theta$ is the corresponding $\mathcal{O}_{\mathbf{V}\mathcal{T}}$-module, there is a natural isomorphism in the derived category of $\mathcal{O}_X$-modules:*

$$\mathrm{R}\mathcal{H}om(i_*(\mathcal{O}_X), E_\theta) \sim \mathrm{R}(E \otimes \Omega^\bullet, \theta).$$

*where $i$ is the zero section of $\mathbf{V}\mathcal{T}$.* □

Finally, let us note that if $(E', \nabla')$ is a module with integrable $p$-connection on $Y/S$, then $\nabla'$ induces a Higgs field on the restriction of $E'$ to $\overline{Y}$, the reduction of $Y$ modulo $p$, and the restriction of the $p$-de Rham complex of $(E', \nabla')$ to $\overline{Y}$ is the Higgs complex of this Higgs field.

## 1.4 Groupoid actions, stratifications, and crystals

Here we review some constructions from [7] and its predecessor [4], recasting them in a more geometric form. To lighten our exposition, we write somewhat informally, referring to [4, II] for a more formal and rigorous presentation.

To fix ideas in a simple situation, let $G$ be a group and $X$ a set. If we view $G$ as a category with one object $y$, then a left (resp. right) action of $G$ on $X$ amounts to a functor (resp. contravariant functor) $F$ from $G$ to the category of sets such that $F(y) = X$. The analogous statement holds for monoids in places of groups. More generally, a representation of a (small) category $\mathcal{C}$ can be viewed as a functor $F$ from $\mathcal{C}$ to the category of sets. Concretely, if $Y$ is the set of objects of $\mathcal{C}$, if $A$ is the set of arrows of $\mathcal{C}$, and if $t\colon A \to Y$ and $s\colon A \to Y$ are the target and source mappings, let $X$ be the disjoint union of $\{F(y) : y \in Y\}$, and let $\pi\colon X \to Y$ be the obvious map. Then for each $a \in A$, we have a map $F(a)\colon F(s(a)) \to F(t(a))$, with $F(a \circ b) = F(a) \circ F(b)$ and $F(\mathrm{id}_y) = \mathrm{id}_{F(y)}$ for all $y \in Y$. This defines a left action of $\mathcal{C}$ on $X$; right actions correspond to contravariant functors.

We now reformulate the previous discussion by replacing the category of sets by a more general category $\mathbf{C}$ which admits fibered products. If $Y$ is an object of $\mathbf{C}$ by a "point" of $Y$ we mean a morphism $\pi_T\colon T \to Y$ in $\mathbf{C}$. We let $\mathbf{C}_Y$ denote the category whose objects are the morphisms of $\mathbf{C}$ whose target is $Y$ and whose morphisms are commutative triangles pointing to $Y$. An object (resp. morphism) of $\mathbf{C}_Y$ is often called a *$Y$-object* (resp. *$Y$-morphism*) of $\mathbf{C}$.

We recall from [4] that a "$\mathbf{C}$-category" $\mathcal{C}$ consists of the following data:

- an object $Y$ of $\mathbf{C}$, representing the objects of $\mathcal{C}$,
- an object $A$ of $\mathbf{C}$ representing the arrows of $\mathcal{C}$,
- target and source morphisms $t, s\colon A \to Y$,
- an identity morphism $\iota\colon Y \to A$, with $t \circ \iota = s \circ \iota = \mathrm{id}_Y$,
- a composition morphism: $c\colon A \times_Y A \to A$, where we view $A$ as a left object over $Y$ via $t$ and a right object over $Y$ via $s$.

Thus the $T$-valued points of $A \times_Y A$ are the pairs $(a_1, a_2) \in A(T) \times A(T)$ such that $s(a_1) = t(a_2)$, and $c(a_1, a_2)$ is defined if $s(a_1) = t(a_2)$). If $Y$ is the object of objects of $\mathcal{C}$, we say that $\mathcal{C}$ is a *category over* $Y$. The category $\mathcal{C}$ is a groupoid if all the arrows in $A$ are isomorphisms, in which case there is an automorphism $\tau$ of $A$ such that

$$\tau^2 = \mathrm{id}_A, \quad t \circ \tau = s, \quad c \circ (\mathrm{id}_A, \tau) = c \circ (\tau, \mathrm{id}_A) = \mathrm{id}_A.$$

**Example 1.4.1** If $f\colon Y \to S$ is a morphism in $\mathbf{C}$, we can form a groupoid $\mathcal{G}_{Y/S}$, with $Y$ the object of objects and $Y(1)/S$ the object of arrows, where:

$$\begin{array}{rcl}
Y(1)/S & := & Y \times_S Y \\
s\colon Y(1)/S \to Y & := & (y_1, y_2) \mapsto y_2 \\
t\colon Y(1)/S \to Y & := & (y_1, y_2) \mapsto y_1 \\
\iota\colon Y \to Y(1)/S & := & y \mapsto (y, y) \\
c\colon Y(1)/S \times_Y Y(1)/S \to Y(1)/S & := & ((y_1, y_2), (y_2, y_3)) \mapsto (y_1, y_3) \\
\tau\colon Y(1)/S \to Y(1)/S & := & (y_1, y_2) \mapsto (y_2, y_1).
\end{array}$$

This groupoid correspond to the equivalence relation given by the morphism $Y \to S$. Thus if $y_1$ and $y_2$ are points of $Y$ such that $f(y_1) = f(y_2)$, there is a unique morphism sending $y_2$ to $y_1$. If $Y \times Y$ exists in $\mathbf{C}$, then the same construction with $Y(1) := Y \times Y$ could be called the "indiscrete" groupoid structure on $Y$: all the points of $Y$ are equivalent.

**Definition 1.4.2.** *Let $\mathcal{C}$ be a* **C***-category, with $Y$ its object of objects and $A$ its object of arrows and let $\pi_V \colon V \to Y$ be a morphism in* **C**. *Then a right action of $\mathcal{C}$ on $\pi_V \colon V \to Y$ is a morphism $r \colon V \times_Y A \to V$ such that the following diagrams commute:*

$$\begin{array}{ccc} V \times_Y A & \xrightarrow{\;r\;} & V \\ {\scriptstyle pr_A}\downarrow & & \downarrow{\scriptstyle \pi_V} \\ A & \xrightarrow{\;s\;} & Y \end{array}$$

$$\begin{array}{ccc} V \times_Y A \times_Y A & \xrightarrow{\;r \times \mathrm{id}_A\;} & V \times_Y A \\ {\scriptstyle \mathrm{id}_V \times c}\downarrow & & \downarrow{\scriptstyle r} \\ V \times_Y A & \xrightarrow{\;r\;} & V \end{array} \qquad \begin{array}{ccc} V & \xrightarrow{\;(\mathrm{id}, \iota \circ \pi_V)\;} & V \times_Y A \\ & {\scriptstyle \mathrm{id}}\searrow & \downarrow{\scriptstyle r} \\ & & V. \end{array}$$

If $v$ and $a$ are $T$-valued points $v$ of $V$ and $A$ respectively, we write $va$ for $r(v, a)$, which is defined if $\pi_V(v) = t(a)$, and the first diagram tells us that $\pi_V(va) = s(a)$. Thus $a$ takes a point lying over the target of $a$ to a point lying over the source of $a$; this is because right actions correspond to contravariant functors. The remaining diagrams say that if $a_1, a_2 \in A$ and $s(a_1) = t(a_2)$ and $t(a_1) = \pi_V(v)$, then $(va_1)a_2 = v(a_1a_2)$, and also that $v\iota(y) = v$ if $y := \pi_V(v)$.

**Remark 1.4.3.** A left action of $\mathcal{C}$ on $V$ is a map $\ell \colon A \times_Y V \to A$ satisfying the analogous diagrams. In this case an arrow $a$ takes a point lying over the source of $a$ to point lying over the target of $a$. Although left actions are perhaps more natural, we have chosen to work with right actions because they seem to be more consistent with the conventions in [4]. If $\mathcal{C}$ is a groupoid, there is an equivalence between left and right actions, thanks to the isomorphism

$$\tau_V := (\tau \circ pr_A, pr_V) \colon V \times_Y A \longrightarrow A \times_Y V. \qquad (1.4.1)$$

In practice we will want to generalize the notation of categorical actions somewhat, so that $V$ belongs not to the category **C**, but rather to a category **V** fibered or cofibered over **C**. In [7], this is explained with the aid of a strict splitting of $\mathbf{V} \to \mathbf{C}$ and with resulting notation written as $\otimes$ in the cofibered case. An earlier version of this manuscript avoided such a choice, at the expense of a rather lengthy exposition that we eventually decided to forego. We will instead assume that a splitting choice has been made and, in the fibered case, just write $V \times_Y T \to V$ for the chosen Cartesian morphism covering a morphism $T \to Y$.

If $V$ is a set equipped with a right action of a monoid $Q$, there is an associated "transporter category" whose objects are the elements of $V$ and whose morphisms from an object $v$ to an object $w$ are the elements $q \in Q$ such that $wq = v$. We adopt this for category actions as follows.

**Definition 1.4.4.** *If $r\colon V \times_Y A \to V$ is a right action of a $\mathbf{C}$-category $\mathcal{C} := (Y, A, t, s, \iota, c)$ on a $Y$-object $V$ of $\mathbf{C}$, the associated transporter category is the $\mathbf{C}$-category whose object object is $V$ and whose arrow object is $V \times_Y A$, and where*

$$\begin{array}{rcl} \iota\colon V \to V \times_Y A & : & v \mapsto (v, \iota(\pi_V(v))) \\ t\colon V \times_Y A \to V & : & (v,a) \mapsto v \\ s\colon V \times_Y A \to V & : & (v,a) \mapsto va \\ c\colon V \times_Y A \times_V V \times_Y A \to V \times_Y A & : & ((v,a),(w,b)) \mapsto (v, ab). \end{array}$$

Note that $ab$ is defined because $t(b) = \pi_V(w) = \pi_V(va) = s(a)$. One checks immediately that these definitions do define a category. If the product $V \times V$ exists, there is a commutative diagam

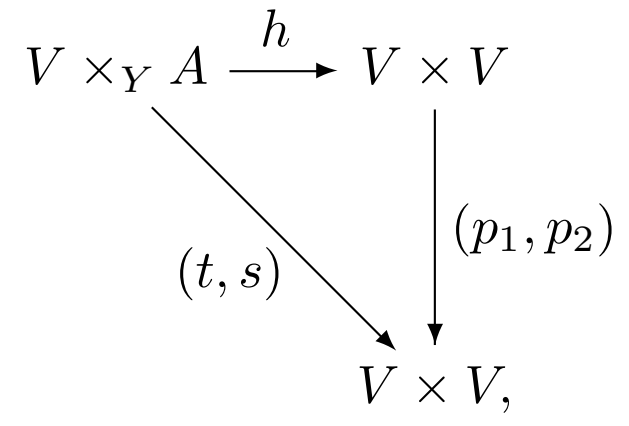

(1.4.2)

where $h := (p_1, r)$, *i.e.*, $h(v,a) = (v, va)$. In fact, $h$ is a functor from the transporter category to the indiscrete groupoid on $V$. Indeed, if $v$ is a point of $V$, then

$$h(\iota(v)) = h(v, \iota(\pi(v))) = (v, v\iota(\pi(v))) = (v, v) = \iota(v).$$

Furthermore, if $(v,a)$ and $(w,b)$ belong to $V \times_Y A$ and $(v,a)(w,b)$ is defined, then $s(v,a) = t(w,b)$, so $va = w$, and $(v,a)(w,b) = (v,ab)$. On the other hand,

$$h(v,a)h(w,b) = (v,va)(w,wb) := (v,wb) = (v,vab).$$

**Definition 1.4.5.** *An action $r$ of a $\mathbf{C}$ groupoid $(Y, A)$ on a $Y$-object $V$ is called torsorial if the map*

$$h\colon V \times_Y A \to V \times V : (v,a) \mapsto (v, va)$$

*from the transporter of $r$ to the indiscrete groupoid of $V$ is an isomorphism.*

Groupoid actions give rise to stratifications, which are generalized equivalence relations. Let us explain this first for category actions. Again we give a definition within a fixed category $\mathbf{C}$, but will apply it for fibered categories without changing notation.

**Definition 1.4.6.** *Let $\mathcal{C}$ be a $\mathbf{C}$-category and let $\pi_V : V \to Y$ be a morphism in $\mathbf{C}$. Then a $\mathcal{C}$-prestratification of $\pi_V : V \to Y$ is a morphism $\epsilon\colon A \times_Y V \to V \times_Y A$ such that the diagrams below commute. A prestratification is a stratification if $\epsilon$ is an isomorphism.*

$$
\begin{array}{ccc}
V \times_Y A & \xrightarrow{\epsilon} & A \times_Y V \\
{\scriptstyle \pi_V \times (s \circ pr_A)} \downarrow & & \downarrow {\scriptstyle (t \circ pr_A) \times \pi_V} \\
Y \times Y & \xrightarrow{\mathrm{id}} & Y \times Y
\end{array}
\qquad
\begin{array}{ccc}
V \times_Y A & \xrightarrow{\epsilon} & A \times_Y V \\
& {\scriptstyle pr_A} \searrow & \downarrow {\scriptstyle pr_A} \\
& & A
\end{array}
$$

$$
\begin{array}{ccccc}
V \times_Y A \times_Y A & \xrightarrow{\epsilon \times \mathrm{id}_A} & A \times_Y V \times_Y A & \xrightarrow{\mathrm{id}_A \times \epsilon} & A \times_Y A \times_Y V \\
{\scriptstyle \mathrm{id}_V \times c} \downarrow & & & \swarrow {\scriptstyle c \times \mathrm{id}_V} & \\
V \times_Y A & \xrightarrow{\epsilon} & A \times_Y V & &
\end{array}
$$

$$
\begin{array}{ccc}
V & \xrightarrow{(\mathrm{id}, \iota \circ \pi_V)} & V \times_Y A \\
& {\scriptstyle (\iota \circ \pi_V, \mathrm{id})} \searrow & \downarrow {\scriptstyle \epsilon} \\
& & A \times_Y V.
\end{array}
$$

If $r$ is a right action on $V \to Y$, then

$$(pr_A, r) \colon V \times_Y A \to A \times_Y V$$

is a prestratification, and if $\epsilon$ is a prestratification, then

$$pr_V \circ \epsilon \colon V \times_Y A \to V$$

is a right action. Thus the two notions are equivalent.

**Remark 1.4.7.** The prestratification associated to an action $r$ is a stratification if and only if the first diagram in Definition 1.4.2 is Cartesian, and this always holds if $\mathcal{C}$ is a groupoid. Indeed, we have a commutative diagram:

$$
\begin{array}{ccccc}
V \times_Y A & \xrightarrow{\epsilon} & A \times_Y V & \xrightarrow{pr_V} & V \\
{\scriptstyle pr_A} \downarrow & & \downarrow {\scriptstyle pr_A} & & \downarrow {\scriptstyle \pi_V} \\
A & \xrightarrow{\mathrm{id}_A} & A & \xrightarrow{s} & Y.
\end{array}
$$

The composition along the top is $r$, and the square on the right is is Cartesian. Thus the rectangle is Cartesian if and only if $\epsilon$ is an isomorphism. If $\mathcal{C}$ is a groupoid, then

$$\tau_V \circ \epsilon \circ \tau_V \circ \epsilon = \mathrm{id}_{V \times_Y A} \quad \text{and} \quad \epsilon \circ \tau_V \circ \epsilon \circ \tau_V = \mathrm{id}_{A \times_Y V}.$$

Thus $\epsilon$ is an isomorphism.

**Remark 1.4.8.** If $V$ and $W$ are $Y$-objects each of which is endowed with a right action, then the fiber product $V \times_Y W$ has a natural right action as well, defined by $r((v,w),a) := (r(v,a), r(w,a))$.

**Example 1.4.9** The projections

$$\begin{array}{ccc} Y \times_Y A \cong A & \xrightarrow{s} & Y \\ A \times_Y Y \cong A & \xrightarrow{t} & Y \end{array} \tag{1.4.3}$$

define right and left actions of $\mathcal{C}$ on $Y$, respectively. These are the "tautological" actions of $\mathcal{C}$ on the points of $Y$. If $y \in Y$ and $a \in A$, the first of these is defined if $t(a) = y$, in which case $r(y,a) = s(a)$, and the second of these is defined if $s(a) = y$, in which case $\ell(a,y) = t(a)$.

Next, take $V \to Y$ to be $s\colon A \to Y$ (resp. $t\colon A \to Y$). We write these as $A^s$ and $A^t$, respectively. The composition law defines a right (resp. left) action of $\mathcal{C}$ on $A^s$ (resp $A^t$):

$$\begin{array}{ccc} A^s \times_Y A & \xrightarrow{c} & A^s \\ A \times_Y A^t & \xrightarrow{c} & A^t \end{array} \tag{1.4.4}$$

These are the right (resp. left) "regular" representations of $\mathcal{C}$ on itself. If $\mathcal{C} = \mathcal{G}$ is a groupoid , then, as we saw in Remark 1.4.3, we also have a right action of $\mathcal{G}$ on $A^t$:

$$c^t := A^t \times_Y A \xrightarrow{\mathrm{id}\times\tau} A^t \times_Y A \xrightarrow{c} A^t \tag{1.4.5}$$

We note that these actions are compatible with the tautological actions of $\mathcal{G}$ on itself. In particular, the diagram;

$$\begin{array}{ccc} A^t \times_Y A & \xrightarrow{c^t} & A \\ {\scriptstyle t\times \mathrm{id}}\downarrow & & \downarrow{\scriptstyle t} \\ Y \times_Y A & \xrightarrow{s} & Y \end{array} \tag{1.4.6}$$

commutes.

**Example 1.4.10** A group can also act by conjugation on itself. The analog for a groupoid $\mathcal{G}$ is a bit more complicated. Assume that the category $\mathbf{C}$ has products, and let $Y(1) := Y \times Y$ and $A(1) := A \times A$. Then

$$\mathcal{G}(1) := (A(1), t \times t, s \times s, c \times c, \iota \times \iota)$$

defines a groupoid over $Y(1)$. Furthermore, the pair $(t,s)$ define a morphism $A \to Y(1)$, allowing us to view $A$ as an object over $Y(1)$. Then we have a right action of $\mathcal{G}(1)$ on $A$:

$$A \times_{Y(1)} A(1) \to A : (a, (g_1, g_2)) \mapsto g_1^{-1} a g_2, \tag{1.4.7}$$

where the fiber product is taken with the morphism $t \times t\colon A(1) \to Y(1)$. This makes sense, because for $(a, (g_1, g_2)) \in A \times_{Y(1)} A(1)$, we have $t(a) = t(g_1)$, $s(a) = t(g_2)$, so $s(g_1^{-1}) = t(a)$, and hence the product $g_1^{-1} a g_2$ is defined.

**Example 1.4.11** Assume that **C** has products, and let $\mathcal{G}_Y$ be the indiscrete groupoid on $Y$ as in Example 1.4.1. Then the tautological action (1.4.3) is given by the map

$$r\colon Y \times_Y Y(1) \longrightarrow Y : (y_1, (y_1, y_2)) \mapsto y_2. \tag{1.4.8}$$

The corresponding stratification is the map

$$\epsilon\colon Y \times_Y Y(1) \longrightarrow Y(1) \times_Y Y : (y_1, (y_1, y_2)) \mapsto ((y_1, y_2), y_2).$$

Let $Y(1)^t$ denote $Y(1)$ viewed as a $Y$-object via $t$, the left projection, and let $Y(1)^s$ denote $Y(1)$ viewed as a $Y$-object via $s$, the right projection. The morphism (1.4.1) for $V = Y(1)^t$ is given by

$$\begin{array}{rcl} \tau_{Y(1)^t}\colon (Y(1)^t)_s & \to & (Y(1)^t)_t \\ Y(1) \times_Y Y(1)^t & \to & Y(1)^t \times_Y Y(1) \\ ((y_1, y_2), (y_2, y_3)) & \mapsto & ((y_2, y_3), (y_2, y_1)) \end{array}$$

The regular representations (1.4.4) and (1.4.5) are given by:

$$\begin{array}{rcl} Y(1)^s \times_Y Y(1) \longrightarrow Y(1)^s & : & ((y_1, y_2), (y_2, y_3)) \mapsto (y_1, y_3) \\ Y(1) \times_Y Y(1)^t \longrightarrow Y(1)^t & : & ((y_1, y_2), (y_2, y_3)) \mapsto (y_1, y_3) \\ Y(1)^t \times_Y Y(1) \longrightarrow Y(1)^t & : & ((y_1, y_2), (y_1, y_3)) \mapsto (y_3, y_2). \end{array} \tag{1.4.9}$$

We shall use the following trivial-looking comparison in the proof of Shiho's theorem explained at the end of §5.6.

**Proposition 1.4.12.** *Let $Y$ be an object of a category **C** with products and fibered products, let $\mathcal{G}_Y$ (resp. $\mathcal{G}_{Y(1)}$) be the indiscrete groupoid on $Y$ (resp., on $Y(1)$), and let $\mathcal{G}_Y(1)$ be the groupoid over $Y(1)$ formed from $\mathcal{G}_Y$ as in the discussion of conjugation in Example 1.4.10. Then there is an isomorphism $\alpha\colon \mathcal{G}_{Y(1)} \to \mathcal{G}_Y(1)$ of groupoids over $Y(1)$. Furthermore, $\alpha$ takes the tautological action (1.4.8) of $\mathcal{G}_{Y(1)}$ on $Y(1)$ to the conjugation action (1.4.7) of $\mathcal{G}_Y(1)$ on $Y(1)$.*

*Proof.* We write this out using points. Define $\alpha$ in the diagram below by the formula:

$$\alpha(y_1, y_2, y_3, y_4) := (y_1, y_3, y_2, y_4).$$

$$\begin{array}{ccc} Y(1) \times Y(1) & \xrightarrow{\alpha} & Y(1) \times Y(1) \\ p_1 \Big\downarrow \Big\downarrow p_2 & & t \times t \Big\downarrow \Big\downarrow s \times s \\ Y(1) & \xrightarrow{\text{id}} & Y(1). \end{array}$$

For the sake of clarity, we check that the diagram commutes. Thus, if $y := (y_1, y_2, y_3, y_4) \in Y(1) \times Y(1)$, we have

$$(t \times t) \circ \alpha(y) = (t \times t)(y_1, y_3, y_2, y_4) = (y_1, y_2) = p_1(y)$$

$$(s \times s) \circ \alpha(y) = (s \times s)(y_1, y_3, y_2, y_4) = (y_3, y_4) = p_2(y).$$

Suppose $y$ and $z$ are elements of $Y(1) \times Y(1)$, viewed as arrows of $\mathcal{G}_{Y(1)}$. Then the composition $y \circ z$ is defined if $t(z) = s(y)$, *i.e.*, if $(y_3, y_4) = (z_1, z_2)$, in which case the composition is $(y_1, y_2, z_3, z_4)$. On the other hand, if $y'$ and $z'$ are elements of $Y(1) \times Y(1)$ viewed as arrows of $\mathcal{G}_Y(1)$, then their composition is defined if $(y'_2, y'_4) = (z'_1, z'_3)$, in which case the composition is $(y'_1, z'_2, y'_3, z'_4)$. It follows that $\alpha$ is compatibile with composition.

Finally, let us check the compatibility with the actions. The tautological action is given by

$$(y_1, y_2) \cdot (y_1, y_2, y_3, y_4) = (y_3, y_4).$$

The conjugation action is given by $h \cdot (g_1, g_2) = g_1^{-1} h g_2$. Thus if $h = (y_1, y_2)$ and $(g_1, g_2) = (y'_1, y'_2, y'_3, y'_4)$, we have

$$(y_1, y_2) \cdot (y'_1, y'_2, y'_3, y'_4) = (y'_2, y'_1)(y_1, y_2)(y'_3, y'_4) = (y'_2, y'_4).$$

Thus, if $y' = \alpha(y)$, these agree. □

Now we can describe the relations among the notions of actions, stratifications, and crystals.

**Definition 1.4.13.** *Let* $\mathbf{F} \colon \mathbf{V} \to \mathbf{C}$ *be a (split) fibered category. A crystal in* $\mathbf{F}$ *(or, somewhat abusively, in* $\mathbf{V}$*) is a Cartesian section of* $\mathbf{F}$*. Explicitly, this means a rule which assigns to each* $T \in \mathbf{C}$ *an object* $V_T$ *of* $\mathbf{V}$ *and to each morphism* $g \colon T' \to T$ *of* $\mathbf{C}$ *an isomorphism* $\epsilon_g \colon V_{T'} \to V_T \times_T T'$ *satisfying the following cocycle conditions:*

*(1)* $\epsilon_{\mathrm{id}} = \mathrm{id}$.

*(2) if* $g \colon T' \to T$ *and* $h \colon T'' \to T'$ *are morphisms in* $\mathbf{C}$*, then*

$$\epsilon_{h \circ g} = \epsilon_h \circ (\epsilon_g \times \mathrm{id}_{T''}).$$

Recall that an object $Y$ of $\mathbf{C}$ is *semifinal* if every object of $\mathbf{C}$ admits at least one morphism to $Y$. Suppose that $Y$ is such an object of $\mathbf{C}$, that $Y(1) := Y \times Y$ is representable, and that $V$ is an object of $\mathbf{V}_Y$ equipped with a right action of the indiscrete groupoid $\mathcal{G}_Y$ described in Example 1.4.11, or, equivalently, a $\mathcal{G}_Y$-stratification $\epsilon$. Then we can define a crystal in $\mathbf{F}$ as follows. For each object $T$ of $\mathbf{C}$, choose some morphism $y_T \colon T \to Y$ and let $V_T := V \times_{y_T} T$.

If $f \colon T' \to T$ is a morphism in $\mathbf{C}$, then $y_{T'}$ and $f \circ y_T$ are two morphisms $T' \to Y$, and we let $g \colon T' \to Y(1)$ be the map with $t \circ g = f \circ y_T$. and $s \circ g = y_{T'}$. Then the $Y(1)$-isomorphism $\epsilon$ induces a $T'$-isomorphism:

$$\epsilon(f) \colon V_T \times_T T' \to V_{T'}.$$

The cocycle conditions guarantees that the resulting family of isomorphisms is compatible with composition and that $\epsilon(\mathrm{id}) = \mathrm{id}$. This construction can be run backwards. We summarize these points of view as follows.

**Proposition 1.4.14.** *Let $\mathbf{F}\colon \mathbf{V} \to \mathbf{C}$ be a fibered category, where $\mathbf{C}$ is a category with fibered products. Suppose that $Y$ is a semi-final object of $\mathbf{C}$, and let $\mathcal{G}_Y$ be the indiscrete groupoid on $Y$ described in Example 1.4.1. Then the following notions are equivalent:*

*(1) An object of $\mathbf{V}_Y$ equipped with a right action of $\mathcal{G}_Y$.*

*(2) An object of $\mathbf{V}_Y$ equipped with a $\mathcal{G}_Y$-stratification.*

*(3) A crystal in $\mathbf{F}$.*

## 1.5 Differential $\mathcal{G}$-operators

We can also formulate a notion of differential operators and their linearization in the general context of groupoid actions. As motivation, we review the approach to the classical definitions explained in [7]. Suppose that $A$ is an $R$-algebra and that $E$ and $F$ are two $A$-modules. An $R$-linear homomorphism $\theta\colon E \to F$ gives rise to (and is equivalent to) an $A$-linear homomorphism $A \otimes_R E \to F$; using the $A$-module structure on $E$, we can write this as an $A$-linear homomorphism

$$(A \otimes_R A) \otimes_A E \to F.$$

Then $\theta$ is said to be a *differential operator of order at most $m$* if this homomorphism factors through a map $(A \otimes_R A)/J^{m+1} \otimes_A E$, where $J \subset A \otimes_R A$ is the ideal of the diagonal. The composition law of the indiscrete groupoid corresponding to $\mathrm{Spec}(A) \to \mathrm{Spec}(R)$ allows one to show that the composition of a differential operator of order $m$ and a differential operator of order $n$ is a differential operator of order $m + n$. In the crystalline context, a *hyper PD-differential operator* from $E$ to $F$ is by definition an $A$-linear homomorphism $D_J(A \otimes_R A) \otimes_A E \to F$, where $D_J(A \otimes_R A)$ is the divided power envelope of $J$. (In this case such an operator is not necessarily determined by its $R$-linear restriction to $E$.) We hope this discussion helps to motivate the definition which follows.

We have defined the notion of an action of a groupoid in category $\mathbf{C}$ on an object in a category $\mathbf{V}$ which is fibered over $\mathbf{C}$. Rather that working in this general context, we will restrict to the case in which $\mathbf{V}$ is the category of arrows in $\mathbf{C}$, which will allow a simplification of the notations. Before stating the definitions, let us recast the motivating discussion above, but in a geometric and more general context.

Let $\mathbf{C}$ be a category with fibered products, let $Y/S$ be a morphism in $\mathbf{C}$, and $\pi_V\colon V \to Y$ and $\pi_W\colon W \to Y$ be objects of $\mathbf{C}_Y$. Composing with the morphism $Y \to S$, we may view $W$ and $Y$ as objects of $\mathbf{C}_S$, and we suppose that we are given an $S$-morphism $f\colon V \to W$. To measure the failure of $f$ to be a $Y$-morphism, we consider the morphism

$$\tilde{f} =: V \xrightarrow{(\pi_V, f)} Y \times_S W \cong Y(1) \times_Y W,$$

which on points takes $v$ to $(\pi_V(v), \pi_W(f(v)), f(v))$. The morphism $f$ is a $Y$-morphism if and only if $\tilde{f}$ factors through the map $W \to Y(1) \times_Y W$ induced by the diagonal. If $A \to Y(1)$ is a groupoid law over $Y$, we can hope that the difference between $\pi_V$ and $\pi_W \circ f$ is mediated by $A \to Y(1)$. We can also crudely "force" $f$ to be a $Y$-morphism by forming

$$\mathcal{L}(f) := Y(1) \times_Y V \cong Y \times_S V \xrightarrow{\mathrm{id}_Y \times f} Y \times_S W \cong Y(1) \times_Y W,$$

and hope to replace $Y(1)$ by $A$. This suggests the following definition.

**Definition 1.5.1.** *Let* $\mathbf{C}$ *be a category with products and fibered products, let* $\mathcal{G} := (Y, A, t, s, \iota, c)$ *be a groupoid object in* $\mathbf{C}$. *If* $\pi_V \colon V \to Y$ *is an object of* $\mathbf{C}_Y$, *let*

$$\mathcal{L}_{\mathcal{G}}(V) := A \times_Y V \text{ and } \mathcal{R}_{\mathcal{G}}(V) := V \times_Y A,$$

*where* $A \times_Y V$ *is viewed as an object over* $Y$ *via the morphism* $t$ *and* $V \times_Y A$ *is viewed as an object over* $Y$ *via the morphism* $s$. *If* $\pi_W \colon W \to Y$ *is another object of* $\mathbf{C}_Y$, *then a differential* $\mathcal{G}$*-operator from* $V$ *to* $W$ *is a* $Y$*-morphism*

$$D \colon V \to \mathcal{L}_{\mathcal{G}}(W).$$

*If* $D$ *is such an operator, then* $\mathcal{L}_{\mathcal{G}}(D) \colon \mathcal{L}_{\mathcal{G}}(V) \to \mathcal{L}_{\mathcal{G}}(W)$ *is defined by the following diagram:*

$$\begin{array}{ccccc} \mathcal{L}_{\mathcal{G}}(V) & = & A \times_Y V & \xrightarrow{\mathrm{id}_A \times D} & A \times_Y (A \times_Y W) \\ \Big\downarrow{\scriptstyle \mathcal{L}_{\mathcal{G}}(D)} & & & & \Big\downarrow{\scriptstyle \cong} \\ {}'\mathcal{L}_{\mathcal{G}}(W) & = & A \times_Y W & \xleftarrow{c \times \mathrm{id}_W} & (A \times_Y A) \times_Y W. \end{array}$$

(Note: If $V$ and $W$ are endowed with additional structure, *e.g.*, as group or module objects, then $D$ should preserve such structure.)

Let us note that the following diagram commutes:

$$\begin{array}{ccc} \mathcal{L}_{\mathcal{G}}(V) & \xrightarrow{\mathcal{L}_{\mathcal{G}}(D)} & \mathcal{L}_{\mathcal{G}}(W) \\ \Big\downarrow & & \Big\downarrow \\ Y(1) \times_Y V & \xrightarrow{\mathcal{L}(D)} & Y(1) \times_Y W. \end{array} \tag{1.5.1}$$

Since $\mathcal{L}_{\mathcal{G}}(V)$ is viewed as an object over $Y$ using the morphism $t$, but the fiber product is formed using $s$, the left action (1.4.4) of $\mathcal{G}$ on $A_t$ induces a left action $\ell_{\mathcal{G}}$ on $\mathcal{L}_{\mathcal{G}}(V)$. We can and shall view this as a right action $r_{\mathcal{G}}$ using the twist $\tau$. Furthermore, if $D$ is a differential $\mathcal{G}$-operator from $V$ to $W$, the map $\mathcal{L}_{\mathcal{G}}(D)$ is compatible with the actions $\ell_{\mathcal{G}}$ and $r_{\mathcal{G}}$ of $\mathcal{G}$ on $\mathcal{L}_{\mathcal{G}}(V)$ and $\mathcal{L}_{\mathcal{G}}(W)$.

Suppose now that $V$ itself is endowed with a right action $r_V$ of $\mathcal{G}$, with corresponding stratification $\epsilon_V$:

$$r_V \colon \mathcal{R}_{\mathcal{G}}(V) \to V; \qquad \epsilon_V \colon \mathcal{R}_{\mathcal{G}}(V) \to \mathcal{L}_{\mathcal{G}}(V).$$

Note that the projection mapping $\pi \colon \mathcal{L}_{\mathcal{G}}(V) \to V$ is in general not compatible with the right actions of $\mathcal{G}$. However, diagram (1.4.6) shows that it *is* compatible if $V = Y$ with its canonical action. Moreover, $\mathcal{R}_{\mathcal{G}}(V)$ inherits a right action $r_{\epsilon,\mathcal{G}}$ of $\mathcal{G}$, deduced from $r_V$ and the right action (1.4.5) of $\mathcal{G}$ on $\mathcal{G}_t$, and the projection $\mathcal{R}_{\mathcal{G}}(V) \to V$ *is* compatible with this right action. Moreover, the morphism $\epsilon_V$ takes $r_{\epsilon,\mathcal{G}}$ to $r_{\mathcal{G}}$. We will explain and prove a generalization of this in Proposition 1.5.2.

Thanks to the action $r_V$, a differential $\mathcal{G}$-operator $D$ from $W$ to $W'$ induces a differential $\mathcal{G}$-operator $\epsilon(D)$ from $V \times_Y W$ to $V \times_Y W'$, given by

$$\begin{array}{ccc} V \times_Y W & \xrightarrow{\mathrm{id}_V \times D} & V \times_Y (A \times_Y W') \cong (V \times_Y A) \times_Y W' \\ \Big\downarrow \epsilon(D) & & \Big\downarrow \epsilon_V \times \mathrm{id} \\ \mathcal{L}_{\mathcal{G}}(V \times_Y W') & = & A \times_Y (V \times_Y W') \cong (A \times_Y V) \times_Y W' \end{array}$$

Also, if $W$ is any $Y$-object, we find an isomorphism:

$$\beta_\epsilon \colon V \times_Y \mathcal{L}_{\mathcal{G}}(W) \to \mathcal{L}_{\mathcal{G}}(V \times_Y W) \tag{1.5.2}$$

defined by the diagram:

$$\begin{array}{ccc} V \times_Y \mathcal{L}_{\mathcal{G}}(W) = V \times_Y (A \times_Y W) & \cong & (V \times_Y A) \times_Y W \\ \Big\downarrow \beta_\epsilon & & \Big\downarrow \epsilon_V \times \mathrm{id}_W \\ \mathcal{L}_{\mathcal{G}}(V \times_Y W) = A \times_Y (V \times_Y W) & \cong & (A \times_Y V) \times_Y W \end{array}$$

The following result is an analog (generalization) of [7, 6.15].

**Proposition 1.5.2.** *Let $\mathcal{G}$ be a $Y$-groupoid in a category* $\mathbf{C}$ *as described in Definition 1.5.1, and let $V$ be a $Y$-object endowed with a right action of $\mathcal{G}$, and let $W$ be any $Y$-object.*

(1) *The map $\beta_\epsilon \colon V \times_Y \mathcal{L}_{\mathcal{G}}(W) \to \mathcal{L}_{\mathcal{G}}(V \times_Y W)$ is an isomorphism, compatible with the right actions of $\mathcal{G}$.*

*(2) If $D$ is a differential $\mathcal{G}$-operator from $W$ to $W'$, the following diagram commutes:*

$$\begin{array}{ccc} V \times_Y \mathcal{L}_{\mathcal{G}}(W) & \xrightarrow{\beta_\epsilon} & \mathcal{L}_{\mathcal{G}}(V \times_Y W) \\ \Big\downarrow{\scriptstyle \mathrm{id} \times \mathcal{L}_{\mathcal{G}}(D)} & & \Big\downarrow{\scriptstyle \mathcal{L}_{\mathcal{G}}(\epsilon(D))} \\ V \times_Y \mathcal{L}_{\mathcal{G}}(W') & \xrightarrow{\beta_\epsilon} & \mathcal{L}_{\mathcal{G}}(V \times_Y W') \end{array}$$

*Proof.* The fact that $\beta_\epsilon$ is an isomorphism is evident from its definition. To check its compatibility with the actions of $\mathcal{G}$, we calculate with points. Suppose $(v, a, w)$ is a point of $V \times_Y \mathcal{L}_G(W)$, so $\pi_V(v) = t(a)$ and $s(a) = \pi_W(w)$. If $b \in A$ with $t(b) = \pi_V(v)$, then

$$\begin{aligned} \beta_\epsilon((v, a, w)b) &= \beta_\epsilon(vb, b^{-1}a, w) \\ &= (\epsilon(vb, b^{-1}a), w) \\ &= (b^{-1}a, va, w) \end{aligned}$$

$$\begin{aligned} (\beta_\epsilon(v, a, w))b &:= (\epsilon(a, v), w)b \\ &= (a, va, w)b \\ &= (b^{-1}a, va, w) \end{aligned}$$

This proves (1). Statement (2) is straightforward and we omit its proof. $\square$

We shall also need the following compatibility. The proof is immediate and omitted.

**Proposition 1.5.3.** *With the notations above, let $\iota_V : V \to \mathcal{L}_{\mathcal{G}}(V)$ be the morphism:*

$$V = Y \times_Y V \xrightarrow{\iota \times \mathrm{id}} A \times_Y V = \mathcal{L}_G(V)$$

*and let $\pi_V : \mathcal{L}_{\mathcal{G}}(V) \to V$ be the projection.*

*(1) $\pi_V \circ \iota_V = \mathrm{id}_V$.*

*(2) If $D$ is a differential operator from $W$ to $W'$, then*

$$\pi_{W'} \circ \mathcal{L}_{\mathcal{G}}(D) \circ \iota_W = \pi_{W'} \circ D. \quad \square$$

## 1.6 Formal $F$-schemes and formal $\delta$-schemes

W recall, and add a bit to, the following definitions from [10], all of which depend of course on a fixed prime number $p$.

**Definition 1.6.1.** *A $\phi$-ring is a pair $(R, \phi)$, where $R$ is a ring and $\phi$ is an endomorphism of $R$ such that $\phi(r) \equiv r^p \pmod{pR}$ for all $r \in R$. A $\delta$-ring is a pair $(R, \delta)$, where $\delta$ is a $\delta$-structure on $R$, i.e., an endomorphism of the set underlying $R$ satisfying the following properties:*

$$\delta(0) = \delta(1) = 0 \tag{1.6.1}$$

$$\delta(ab) = a^p\delta(b) + b^p\delta(a) + p\delta(a)\delta(b) \tag{1.6.2}$$

$$\delta(a+b) = \delta(a) + \delta(b) - C(a,b), \tag{1.6.3}$$

*where $C(A, B)$ is the homogenous polynomial of degree $p$ in $\mathbf{Z}[A, B]$ such that*

$$(A+B)^p = A^p + B^p + pC(A,B).$$

*That is,*

$$C(A,B) := \sum_{1}^{p-1} c_i A^i B^{p-i}, \text{ where } c_i := \frac{(p-1)!}{i!(p-i)!}. \tag{1.6.4}$$

*A homomorphism of $\phi$-rings (resp. $\delta$-rings) is a homomorphism of the underlying rings compatible with $\phi$ (resp. with $\delta$).*

Thus a $\phi$-ring is just a ring together with a lift of the Frobenius endomorphism of $R/pR$ to an endomorphism of $R$. If $(R, \delta)$ is a $\delta$-ring and if $\phi\colon R \to R$ is defined by $\phi(r) := r^p + p\delta(r)$, then $\phi$ is a ring homomorphism and $(R, \phi)$ is a $\phi$-ring. Conversely, if $(R, \phi)$ is a $p$-torsion free $\phi$-ring and if $\delta(r) := p^{-1}(\phi(r) - r^p)$, then $(R, \delta)$ is a $\delta$-ring. Thus one has an equivalence between the notions of $p$-torsion free $\delta$-rings and $p$-torsion free $\phi$-rings. Each of these notions has its advantages and disadvantages, and our approach will require both. Note for example that the only $\delta$-ring in which $p$ is nilpotent is the zero ring. Indeed, properties (1)–(3) imply that $\delta(0) = 0$ and that $\delta(p^n) = p^{n-1}(1 - p^{n(p-1)})$, so $p^n = 0$ implies that $p^{n-1} = 0$, hence by induction that $1 = 0$.

We defer to [10] for a discussion of the formalism of $\delta$-rings. In particular, as explained in [10, §2], the notions of $\phi$-rings and $\delta$-rings are compatible with Zariski and étale localization, provided one works with $p$-adically complete rings and $p$-complete localization. This follows from the following result, which will allow us to globalize these notions to the setting of formal schemes.

**Proposition 1.6.2.** *Let $(A, \delta)$ be a $\delta$-ring and $I$ a finitely generated ideal of $A$ containing $p$.*

*(1) The map $\delta\colon A \to A$ is $I$-adically continuous [10, 2.17].*

*(2) If $p$ belongs to the radical of $A$ and if $S$ is a multiplicative subset of $A$, then $\delta$ extends naturally to the $p$-localization of $S^{-1}A$ [10, 2.16].*

*(3) If $B$ is derived $I$-complete and $A \to B$ is an $I$-completely étale homomorphism, then the $\delta$-structure of $A$ extends uniquely to a $\delta$-structure on $B$ [10, 2.18].* □

**Definition 1.6.3.** *A formal F-scheme is a p-adic formal scheme $Y$ endowed with an endomorphism $F_Y$ lifting the absolute Frobenius endomorphism $F_{\overline{Y}}$ of $\overline{Y}$. A formal δ-scheme is a p-adic formal scheme $Y$ endowed with an endomorphism $\delta_Y$ of the sheaf of sets underlying $\mathcal{O}_Y$ satisfying properties (1)–(3) of Definition 1.6.1. A crystalline prism is a p-torsion free formal F-scheme, or, equivalently, a p-torsion free formal δ-scheme. Morphisms of these objects are defined in the obvious way.*

If $(Y, F_Y)$ is a formal F-scheme, the endomorphism $F_Y$ of the topological space underlying $Y$ is just the identity map, and the structure of the Frobenius lift just amounts to an endomorphism $\phi_Y$ of the structures sheaf $\mathcal{O}_Y$ such that $\phi_Y(a) \equiv a^p \pmod{p\mathcal{O}_Y}$ for $a \in \mathcal{O}_Y$. Thus a formal δ-scheme gives rise to a formal F-scheme, and the two notions become equivalent in the $p$-torsion free case. Thanks to the localization properties of δ-schemes explained above, we can view a formal δ-scheme with its Zariski or $p$-completely étale topology. Furthermore, if $(R, \phi)$ (resp. $(R, \delta)$ is a $\phi$-ring (resp. δ-ring) which is $p$-adically complete and has bounded $p$-torsion (and thus is $p$-derived complete) then its formal spectrum has a natural structure of a formal F-scheme (resp. formal δ-scheme). Often we will work with a fixed formal F-scheme $S$ and the category of morphisms of formal F-schemes $Y \to S$, which we shall call *S-F-schemes.*

**Example 1.6.4** The Frobenius endomorphism of $\mathbf{F}_p$ is the identity, and the identity map is its unique lifting to $\mathbf{Z}$. Thus $\mathbf{Z}$ has a unique δ-structure, with $\delta(n) = p^{-1}(n - n^p)$ for all $n$, and $(\mathbf{Z}, \delta)$ is the initial object in the category of δ-rings [10, 2.6].

Let $\mathbf{Z}[x_\bullet]$ be the polynomial ring on the variables $\{x_i : i \in \mathbf{N}\}$, and define $\phi\colon \mathbf{Z}[x_\bullet] \to \mathbf{Z}[x_\bullet]$ by $\phi(x_i) := x_i^p + px_{i+1}$. Since $\phi$ is a Frobenius lift, the map $\delta\colon x_i \mapsto x_{i+1}$ is a δ-structure on $\mathbf{Z}[x_\bullet]$. In fact, $(\mathbf{Z}[x_\bullet], \delta)$ is the free δ-ring on the variable $x_0$. Namely, if $r$ is an element of a δ-ring $(R, \delta)$, there is a unique homomorphism of δ-rings $\alpha\colon (\mathbf{Z}[x, \bullet], \delta) \to R$ sending $x_0$ to $r$, determined inductively by the formula:

$$\alpha(x_{i+1}) = \delta(\alpha(x_i)).$$

Note that the multiplicative set $S$ generated by $x_0$ is not invariant under $\phi$, and hence $\phi$ does not define a Frobenius lift of the scheme $\operatorname{Spec} \mathbf{Z}[x_\bullet]$. As Proposition 1.6.2 predicts, this problem disappears if we pass to the $p$-adic completion. For example:

$$\phi(x_0^{-1}) = x_0^{-p} - px_1x_0^{-2p} + p^2x_1^2x_0^{-3p} - \cdots .$$

More generally, if $(R, \delta)$ is any δ-ring and $X$ is any set, one can form the polynomial ring $R[X]$ on the variables in $X$ and then the polynomial ring $R[X_\bullet]$, where now $X_\bullet := \{x_j : x \in X, i \in \mathbf{N}\}$ and where we identify $x \in X$ with $x_0$. Then $R[X_\bullet]$ has a δ structure such that $\delta(x_j) = x_{j+1}$ for all $x \in X$ and all $j \in \mathbf{N}$, and $(R[X_\bullet], \delta)$ is the free $R$-δ-algebra on the variables in $X$. This object is denoted by $R\{X\}$ in [10], but we shall not use this notation, writing instead $R[X, \delta]$ for this δ-ring. Since we will be working $p$-adically, we shall work from now on with $\mathbf{Z}_p$ instead of $\mathbf{Z}$ as base ring.

**Remark 1.6.5.** As explained in [10, 2.7], the category of $\delta$-rings admits all limits and colimits, and their formation commutes with the forgetful functor to the category of rings. In particular, if $(A, \delta) \to (B, \delta)$ and $(A, \delta) \to (C, \delta)$ are homomorphisms of $\delta$-rings, then there is a unique compatible $\delta$-structure on $B \otimes_A C$.

Joyal [22] applies these facts as follows. Recall that $\mathbf{Z}_p[x, \delta]$ is the free $\delta$-ring on the variable $x$. As we saw above, the tensor product $\mathbf{Z}_p[x, \delta] \otimes_{\mathbf{Z}_p} \mathbf{Z}_p[x, \delta]$ admits a unique $\delta$-structure such that $\delta(x_i \otimes 1) = x_{i+1} \otimes 1$ and $\delta(1 \otimes x_i) = 1 \otimes x_{i+1}$ for all $i$. The scheme $\mathrm{Spec}(\mathbf{Z}_p[x])$ is a ring scheme, with coaddition and comultiplication given respectively by $x \mapsto x \otimes 1 + 1 \otimes x$ and $x \mapsto x \otimes x$. Since $(\mathbf{Z}_p[x_\bullet], \delta)$ is the free $\delta$-ring on the variable $x$, there are unique homomorphisms of $\delta$-rings:

$$
\begin{array}{rcl}
\mathbf{Z}_p[x, \delta] \to \mathbf{Z}_p[x, \delta] \otimes \mathbf{Z}_p[x, \delta] & : & x_0 \mapsto x_0 \otimes 1 + 1 \otimes x_0 \\
\mathbf{Z}_p[x, \delta] \to \mathbf{Z}_p[x, \delta] \otimes \mathbf{Z}_p[x, \delta] & : & x_0 \mapsto x_0 \otimes x_0
\end{array}
$$

These operations give $\mathbf{W} := \mathrm{Spec}(\mathbf{Z}_p[x, \delta])$ the structure of a ring-scheme. If now $A$ is any $\mathbf{Z}_p$-algebra, Joyal identifies the Witt ring $\mathrm{W}(A)$ with the set of $A$-valued points of $\mathbf{W}$, with the ring structure determined by the ring scheme structure of $\mathbf{W}$ and the $\delta$-structure taking a point $(a_0, a_1, a_2, \cdots) \in \mathrm{W}(A)$ to $(a_1, a_2, a_3, \cdots)$, in so-called "Joyal coordinates." As Joyal explains, the functor W from the category of rings to the category of $\delta$-rings is right adjoint to the forgetful functor. Thus, if $(R, \delta)$ is a $\delta$-ring and $\alpha \colon R \to A$ is a ring homomorphism, there is a $\delta$-homomorphism

$$
\alpha_\bullet \colon (R, \delta) \to (\mathrm{W}(A), \delta) : r \mapsto (\alpha(r), \alpha(\delta(r)), \alpha(\delta^2(r)), \ldots).
$$

The morphism $F \colon \mathbf{W} \to \mathbf{W}$ induced by $\phi$ is a morphism of ring-schemes, and its kernel $\mathbf{W}_F$ is the subscheme defined by the ideal $I_F$ generated by $\{x_i^p + p x_{i+1} : i \in \mathbf{N}\}$, which is invariant under $\phi$. Thus $\mathbf{W}_F$ inherits a Frobenius lifting. Furthermore, the closed immersion $\mathbf{W}_F \to \mathbf{W}$ is invariant under the action of $\mathbf{W}$ on itself, so that $\mathbf{W}_F$ is an ideal subscheme of $\mathbf{W}$ and, in particular, is a $\mathbf{W}$-module.

The relationship between $\mathbf{W}_F$ and the divided power envelope of the zero section of the affine line $\mathbf{A}$ is a key insight of Drinfeld [15, 3.2.6]. To make this explicit and precise, we shall view $\mathbf{A} = \mathrm{Spec}(\mathbf{Z}_p[x])$ as a ring scheme, representing the identity functor from the category of rings to itself. There is morphism of ring schemes $\mathbf{W} \to \mathbf{A}$ sending $(a_0, a_1, \ldots)$ to $a_0$. Note that the action of $\mathbf{W}$ on $\mathbf{W}_F$ factors through this morphism. Namely, if $x$ belongs to its kernel and $y$ belongs to $\mathbf{W}_F$, then $x = Vx'$ for some $x' \in \mathbf{W}$, and then $xy = yVx' = V(x'Fy) = 0$.

The map $\mathbf{Z}_p[x] \to \mathbf{Z}_p\langle x \rangle$ defines a morphism of schemes: $\mathbf{A}^\gamma \to \mathbf{A}$. The group scheme structure of $\mathbf{A}$ induces a unique group scheme structure on $\mathbf{A}^\gamma$, and the ring structure of $\mathbf{A}$ induces an $\mathbf{A}$-module structure on $\mathbf{A}^\gamma$. These are given by coaddition and comultiplication laws:

$$
x^{[n]} \quad \mapsto \quad \sum_{i+j=n} x^{[i]} \otimes x^{[j]}
$$

$$x^{[n]} \quad \mapsto \quad x^n \otimes x^{[n]},$$

respectively. The morphism $\mathbf{A}^\gamma \to \mathbf{A}$ is not a closed immersion, and so forms a quasi-ideal scheme (see [15, 3.1]) rather than an ideal scheme.

The following proof of Drinfeld's lemma follows a suggestion of B. Bhatt and adds a few extra details.

**Proposition 1.6.6.** *Let $I_F$ be the ideal of $\mathbf{Z}_p[x,\delta]$ generated by $\{x_i^p + px_{i+1} : i \in \mathbf{N}\}$, and let $\mathbf{A} := \operatorname{Spec} \mathbf{Z}_p[x]$ and $\mathbf{A}^\gamma := \operatorname{Spec} \mathbf{Z}_p\langle x \rangle$.*

*(1) The ring $\mathbf{Z}_p[x,\delta]/I_F$ is a free $\mathbf{Z}_p$-module, with basis the set of monomials*

$$x_\bullet^{a_\bullet} := x_0^{a_0} x_1^{a_1} \cdots x_r^{a_r} : a_i \in [0,p), r \in \mathbf{N}.$$

*(2) There is a unique isomorphism*

$$\omega \colon \mathbf{Z}_p\langle x \rangle \longrightarrow \mathbf{Z}_p[x,\delta]/I_F$$

*which fits in the commutative diagram on the left below. It is compatible with the coaddition and coaction laws and induces the isomorphism of group-schemes $\mathbf{W}_F \to \mathbf{A}^\gamma$ in the diagram on the right.*

$$\begin{array}{ccc} \mathbf{Z}_p\langle x \rangle & \xrightarrow{\cong} & \mathbf{Z}_p[x,\delta]/I_F \\ \uparrow & & \uparrow \\ \mathbf{Z}_p[x] & \longrightarrow & \mathbf{Z}_p[x,\delta] \end{array} \qquad \begin{array}{ccc} \mathbf{W}_F & \xrightarrow{\cong} & \mathbf{A}^\gamma \\ \downarrow & & \downarrow \\ \mathbf{W} & \longrightarrow & \mathbf{A}. \end{array}$$

*(3) For $j \in \mathbf{N}$, let*

$$s_j := (-p)^{1+p+\cdots+p^{j-1}} = (-p)^{\frac{p^j-1}{p-1}},$$

*and for $n \in \mathbf{N}$ with $p$-adic expansion $a_0 + a_1 p + a_2 p^2 + \cdots + a_r p^r$, let*

$$c_n := \frac{s_0^{a_0} s_1^{a_1} \cdots s_r^{a_r}}{n!}.$$

*Then $c_n$ is a $p$-adic unit, and the homomorphism $\omega$ is given by:*

$$\omega : x^{[n]} \mapsto c_n x_0^{a_0} x_1^{a_1} \cdots x_r^{a_r}.$$

*If $n$ is a power of $p$, then $c_n$ is congruent to $1$ modulo $p$.*

*Proof.* For $n \in \mathbf{N}$, let $B_n$ be the subring of $\mathbf{Z}_p[x,\delta]/I_F$ generated by $\{x_j : j \le n\}$. Thus $B_n = \mathbf{Z}_p[x_0,\ldots,x_n]/I_n$, where $I_n$ is the ideal generated by $\{x_0^p + px_1, x_1^p + px_2, \ldots, x_{n-1}^p + px_n\}$. This ring can be constructed inductively from $\mathbf{Z}_p[x_n]$ by successively adjoining variables $x_{n-j}$ in descending order and dividing by the monic polynomial $x_{n-j}^p + px_{n-j+1}$. Statement (1) can easily be checked from this construction. It follows that $\mathbf{Z}_p[x,\delta]/I_F$ is torsion free. Since the homomorphism $\mathbf{Z}_p[x] \to \mathbf{Z}_p\langle x\rangle$ is an isomorphism after tensoring with $\mathbf{Q}_p$, the homomorphism $\omega$ will be unique if it exists and will be compatible with the coaddition laws.

The morphism $\mathbf{Z}_p[x] \to \mathbf{Z}_p[x,\delta]$ sends $x$ to $x_0$. Working modulo $I_F$, we have $x_j^p + px_{j+1} = 0$, and it follows by induction that

$$(-p)^{1+p+\cdots+p^{j-1}} x_j = x_0^{p^j}.$$

Thus $x^{p^j}$ maps to $s_j x_j$ in $\mathbf{Z}_p[x,\delta]/I_F$, and if $n = a_0 + a_1 p + \cdots a_r p^r$, then $x^n$ maps to $x_0^{a_0} x_1^{a_1} \cdots x_r^{a_r} s_0^{a_0} s_1^{a_1} \cdots s_r^{a_r}$. To prove that $\omega$ exists and is an isomorphism, it will suffice to prove that each $c_n$ is a $p$-adic unit. We calculate:

$$\begin{aligned}
(p-1)\operatorname{ord}_p(s_0^{a_0} s_1^{a_1} \cdots s_r^{a_r}) &= (p-1)\left(a_0 \operatorname{ord}_p(s_0) + a_1 \operatorname{ord}_p(s_1) + \cdots + a_r \operatorname{ord}_p(s_r)\right) \\
&= a_1(p-1) + a_2(p^2-1) + \cdots a_r(p^r-1) \\
&= a_0 + a_1 p + a_2 p^2 + \cdots a_r p^r - (a_0 + a_1 + \cdots a_r) \\
&= n - (a_0 + a_1 + \cdots a_r) \\
&= (p-1)\operatorname{ord}_p(n!),
\end{aligned}$$

by statement (2) of Lemma 1.6.7 below. This shows that $c_n$ is indeed a $p$-adic unit.

For the last statement, we may as well assume that $p$ is odd. Then

$$(-1)^{1+p+\cdots p^{j-1}} = (-1)^j.$$

On the other hand, statement (1b) of Lemma 1.6.7 implies that $U_p((p^j)!) = (-1)^j$, *i.e.*, that $p^j! = p^{1+\cdots+p^{j-1}} u$ where $u \equiv (-1)^j \pmod{p}$. This proves our claim. □

**Lemma 1.6.7.** *For each natural number $n$, let $U_p(n) \in \mathbf{F}_p^*$ be the reduction modulo $p$ of $p^{-\operatorname{ord}_p(n)} n$ and let $\sigma(n)$ be the sum of the digits in the $p$-adic expansion of $n$.*

*(1) If $m$ and $n$ are natural numbers, then*

- *(a)* $\operatorname{ord}_p((p^m n)!) = \frac{n(p^m-1)}{p-1} + \operatorname{ord}_p(n!)$
- *(b)* $U_p((p^m n)!) = (-1)^{mn} U_p(n!)$.
- *(c)* $U_p((p^m)!) = (-1)^m$.

*(2) If $n \in \mathbf{N}$, then*

$$\operatorname{ord}_p(n!) = \frac{n - \sigma(n)}{p-1}.$$

*(3) If $n \in \mathbf{N}$, then* $\operatorname{ord}_p((pn)!) - \operatorname{ord}_p((n!)^p) = \sigma(n)$.

*Proof.* For a natural number $k$, let $u_k := (pk+1)(pk+2)\cdots(pk+p-1)$, and note that $u_k \equiv (p-1)! \equiv -1 \pmod{p}$. Then

$$\begin{aligned}
(pn!) &= 1\cdot 2\cdot\cdots\cdot(p-1)\cdot p\cdot(p+1)\cdot(p+2)\cdot\cdots\cdot(p+p-1)\cdot(2p)\cdot\cdots\cdot pn \\
&= u_0\cdot p\cdot u_1\cdot 2p\cdot u_2\cdot 3p\cdot u_3\cdot\cdots\cdot u_{n-1}\cdot np \\
&= p^n\cdot 1\cdot 2\cdot\cdots\cdot n\cdot u_0\cdot u_1\cdot\cdots\cdot u_{n-1} \\
&= p^n\cdot n!\cdot u_0\cdot u_1\cdot\cdots\cdot u_{n-1},
\end{aligned}$$

and $u_0\cdot u_1\cdot\cdots\cdot u_{n-1} \equiv (-1)^n \pmod{p}$. This proves statement (1b) when $m=1$.

The reader can easily prove statement (1a) by induction, using the above equality. Statement (1b) is trivial if $m=0$ or if $p=2$, so we assume that $p$ is odd, that it is true for $m$, and that $m \geq 1$. Then, using statement (1b) when $m=1$, we calculate:

$$\begin{aligned}
U_p((p^{m+1}n)!) &= U_p((p(p^m n))!) \\
&= (-1)^{p^m n}U_p((p^m n)!) \\
&= (-1)^n U_p((p^m n)!) \\
&= (-1)^n(-1)^{mn}U_p(n!) \\
&= (-1)^{n(m+1)}U_p(n!).
\end{aligned}$$

Statement (1c) is a special case of statement (1b).

Formula (2) is well-known; see for example [7, 3.3]. Let us nonetheless note that it follows easily by induction on $n$ from the case $m=1$ of formula (1a). Indeed, if $n=0$ it is obvious, and if it is true for $n-1$ and $p \nmid n$, then $\sigma(n) = \sigma(n-1)+1$ and $\mathrm{ord}_p(n!) = \mathrm{ord}_p((n-1)!)$, so the result also holds for $n$. If, on the other hand, $n = pn'$, then $\sigma(n) = \sigma(n')$, and formula (1a) and an induction assumption together imply that

$$\mathrm{ord}_p(n!) = n' + \mathrm{ord}_p(n'!) = n' + \frac{n'-\sigma(n')}{p-1} = \frac{n-\sigma(n)}{p-1}.$$

For statement (3), we just calculate:

$$\begin{aligned}
\mathrm{ord}_p(n!) &= \frac{n-\sigma(n)}{p-1} \\
\mathrm{ord}_p((pn)!) - \mathrm{ord}_p((n!)^p) &= \frac{pn-\sigma(np)-p(n-\sigma(np))}{p-1} \\
&= \frac{-\sigma(n)+p\sigma(n)}{p-1} \\
&= \sigma(n).
\end{aligned}$$

□

If $A$ is a ring and $R$ is an $A$-algebra, let $\mathrm{A}^\gamma(R)$ denote the set of $R$-valued points of $\mathbf{A}^\gamma$, *i.e.*, the set of $A$-algebra homomorphisms $A\langle x\rangle \to R$. Thanks to the **A**-module structure of $\mathrm{A}^\gamma$, $\mathrm{A}^\gamma(R)$ has a natural structure of an $A$-module. The following result shows that the functor $\mathrm{A}^\gamma$ can be viewed as the right adjoint of the functor taking an $A$-module $E$ to the $A$-algebra $\Gamma^\bullet(E)$ constructed by Roby in [36] This result is implicit in his work,

**Proposition 1.6.8.** *If $A$ is a ring, if $E$ is an $A$-module, and if $R$ is an $A$-algebra, there is a natural isomorphism*

$$\mathrm{Hom}_{A-alg}(\Gamma^\bullet(E), R) \cong \mathrm{Hom}_{A-mod}(E, \mathrm{A}^\gamma(R))$$

*taking an $A$-algebra homomorphism $\alpha\colon \Gamma^\bullet(E) \to R$ to the $A$-module homomorphism $\hat{\alpha}\colon E \to \mathrm{A}^\gamma(R)$, where, for each $e \in E$,*

$$\hat{\alpha}(e) : x^{[n]} \mapsto \alpha(e^{[n]}).$$

*Proof.* Let us check that this definition makes sense. Since $\{x^{[n]}\}$ is a basis for $A\langle x\rangle$, there is a well-defined $A$-linear homomorphism $\hat{\alpha}(e)$ sending each $x^{[n]}$ to $\alpha(e^{[n]})$. Moreover:

$$\begin{aligned}
\hat{\alpha}(e)(x^{[n]}x^{[m]}) &:= \hat{\alpha}(e)(\binom{n+m}{n}x^{[n+m]}) \\
&= \binom{n+m}{n}\hat{\alpha}(e)(x^{[n+m]}) \\
&= \binom{n+m}{n}\alpha(e^{[n+m]}), \\
\hat{\alpha}(e)(x^{[n]})\hat{\alpha}(e)(x^{[m]}) &:= \alpha(e^{[n]})\alpha(e^{[m]}) \\
&= \alpha(e^{[n]}e^{[m]}) \\
&= \alpha\left(\binom{n+m}{m}e^{[n+m]}\right) \\
&= \binom{n+m}{m}\alpha(e^{[n+m]})),
\end{aligned}$$

so $\hat{\alpha}(e)$ is in fact an algebra homomorphism.

We next check the linearity of $\hat{\alpha}$. If $e_1, e_2 \in E$:

$$\begin{aligned}
\hat{\alpha}(e_1+e_2)(x^{[n]}) &= \alpha((e_1+e_2)^{[n]}) \\
&= \alpha\left(\sum_{i+j=n} e_1^{[i]}e_2^{[j]}\right) \\
&= \sum_{i+j=n} \alpha(e_1^{[i]})\alpha(e_2^{[j]}), \\
(\hat{\alpha}(e_1)+\hat{\alpha}(e_2))(x^{[n]}) &= \mathrm{mult}_R\left(\hat{\alpha}(e_1)\otimes\hat{\alpha}(e_2)\left(\sum_{i+j=n} x^{[i]}\otimes x^{[j]}\right)\right) \\
&= \mathrm{mult}_R\left(\sum_{i+j=n}(\hat{\alpha}(e_1)(x^{[i]})\otimes\hat{\alpha}(e_2)(x^{[j]}))\right)
\end{aligned}$$

$$= \sum_{i+j=n} \alpha(e_1^{[i]})\alpha(e_2^{[j]}).$$

If $e \in E$ and $a \in A$, let $\hat{a} \in \mathbf{A}(A)$ be the homomorphism $A[X] \to A$ sending $x$ to $a$. Then:

$$\begin{aligned}
\hat{\alpha}(ae)(x^{[n]}) &= \alpha((ae)^{[n]}) \\
&= \alpha(a^n e^{[n]}) \\
&= a^n \alpha(e^{[n]}), \\
a(\hat{\alpha})(e))(x^{[n]}) &= \mathrm{mult}_R((\hat{a} \otimes \hat{\alpha}(e))(x^n \otimes x^{[n]})) \\
&= a^n \alpha(e^{[n]}).
\end{aligned}$$

Conversely, suppose that $\beta$ is an $A$-linear map $E \to \mathrm{A}^\gamma(R)$. For each $n \in \mathbf{N}$, we have a map $E \to R$ sending $e$ to $\beta(e)(x^{[n]})$. These maps are compatible with further homomorphisms $R' \to R''$, and thus define a polynomial law in the sense of Roby [36], which one checks has degree $n$. Thus we find a linear map $\Gamma^n(E) \to R$ sending $e$ to $\beta(e)(x^{[n]})$. We leave it to the reader to verify that these maps fit together to form an algebra homomorphism and that $\beta \to \hat{\beta}$ is inverse to $\alpha \mapsto \hat{\alpha}$. □

**Corollary 1.6.9.** *Let $A$ be a $\delta$-ring, let $\theta\colon A \to B$ be a ring homomorphism, and let $\theta_W\colon A \to W(B)$ be the corresponding homomorphism of $\delta$-rings. If $E$ is an $A$-module and $pE = pB = 0$, there is an isomorphism of abelian groups from the set of $A$-algebra homomorphisms: $\alpha\colon \Gamma^\bullet(E) \to B$ to the set of $A$-linear maps $E \to W(B)$. In Joyal coordinates, this isomorphism takes $\alpha$ to the map $\alpha_W\colon E \to W(B)$ given by*

$$\alpha_W\colon e \mapsto (\alpha(e), \alpha(e^{[p]}), \alpha(e^{[p^2]}), \cdots).$$

*Proof.* Since $pE = 0$, every homomorphism $E \to W(B)$ factors through the kernel of multiplication by $p$ on $W(B)$. Since $pB = 0$, the endomorphisms $F$ and $V$ of $W(B)$ commute, and since $FV = VF = p$ and $V$ is injective, the kernel of multiplication by $p$ on $W(B)$ equals the kernel of the action of $F$. Note that the group structure on this set comes from the coaddition law on $\mathbf{Z}\langle x\rangle \cong \mathbf{Z}_p[x,\delta]/I_F$, and the $A$-module structure comes from $\theta_W$ and the action of $\mathrm{W}(B)$ on $\mathrm{W}_F(B)$, and factors through the map $A \xrightarrow{\theta_W} \mathrm{W}(B) \to B$. Combining these facts with Propositions 1.6.6 and 1.6.8, we find:

$$\mathrm{Hom}_A(E, \mathrm{W}(B)) \cong \mathrm{Hom}_A(E, \mathrm{W}_F(B)) \cong \mathrm{Hom}_A(E, \mathrm{A}^\gamma(B)) \cong \mathrm{Hom}_A(\Gamma^\bullet(E), B).$$

To check the formula, let $\alpha\colon \Gamma^\bullet(E) \to B$ be an algebra homomorphism. By Proposition 1.6.8, the corresponding map $\hat{\alpha}\colon E \to \mathrm{A}^\gamma(B)$ is given by

$$\hat{\alpha}(e)(x^{[n]}) = \alpha(e^{[n]}).$$

Statement (3) of Proposition 1.6.6 implies that $\omega\colon \mathbf{Z}\langle x\rangle \to \mathbf{Z}_p[x,\delta]/I_F$ takes $x^{[p^j]}$ to $c_{p^j} x_j$ and that $c_{p^j} \equiv 1 \pmod p$. Thus $\alpha_W(e)(x_j) = \alpha(x^{[p^j]})$ as claimed. □

**Remark 1.6.10.** As is observed in [10, 2.7], the forgetful functor from $\delta$-rings to rings also has a left adjoint J, which can be described as follows. If $A$ is a ring, choose a presentation $\mathbf{Z}_p[X] \to A$ and let $K$ be its kernel. Then $\mathrm{J}(A) = \mathbf{Z}_p[X,\delta]/K_\delta$ where $K_\delta$ is the smallest $\delta$-invariant ideal of $\mathbf{Z}_p[X,\delta]$ containing $K$. This construction can be easily generalized to a relative setting, over an arbitrary $\delta$-ring $R$ in place of $\mathbf{Z}_p$, in which case we write $\mathrm{J}_R$ instead of J.

We shall want to use geometric versions of the constructions in Remarks 1.6.5 and 1.6.10, working with formal schemes in place of schemes. Let $S$ be a crystalline prism (for example, $\mathbf{Z}_p$ with its unique $\delta$-structure), Let $\Sigma_S$ denote the category of $p$-adic formal $S$-schemes and let $\Sigma_{S,\delta}$ (resp. $\Sigma_{S,F}$, resp $\mathbb{\Delta}_S$) denote the category of formal $\delta$-schemes over $S$ (resp. of formal F-schemes, of crystalline prisms, over $S$.).

**Proposition 1.6.11.** *Consider the functors:*

$$\begin{array}{ccccc} \mathbb{\Delta}_S & \longrightarrow & \Sigma_{S,\delta} & \longrightarrow & \Sigma_{S,F} \\ & {\scriptstyle \mathrm{O}}\searrow & \downarrow{\scriptstyle \mathrm{O}} & \swarrow{\scriptstyle \mathrm{O}} & \\ & & \Sigma_S, & & \end{array}$$

*where the horizontal arrows are the obvious ones and* O *is the forgetful functor.*

*(1) The functor* $\mathrm{O}\colon \Sigma_{S,\delta} \to \Sigma_S$ *admits a left adjoint*

$$\mathrm{W}\colon \Sigma_S \to \Sigma_{S,\delta}$$

*and a right adjoint*

$$\mathrm{J}_S\colon \Sigma_S \to \Sigma_{S,\delta}.$$

*(2) The inclusion functor* $\mathbb{\Delta}_S \to \Sigma_{S,\delta}$ *admits a left adjoint*

$$(Y,\delta) \mapsto (Y_{\mathrm{tf}},\delta).$$

*(3) Thus, for every* $T \in \Sigma_S$ *and* $Z \in \Sigma_{S,\delta}$, *we have*

$$\begin{array}{rcl} \mathrm{Mor}_{\Sigma_{S,\delta}}(\mathrm{W}(T),Z) & \cong & \mathrm{Mor}_{\Sigma_S}(T,\mathrm{O}(Z)) \\ \mathrm{Mor}_{\Sigma_{S,\delta}}(Z,\mathrm{J}_S(T)) & = & \mathrm{Mor}_{\Sigma_S}(\mathrm{O}(Z),T), \end{array}$$

*and if* $T$ *and* $T'$ *are objects of* $\Sigma_{\mathrm{p}}$, *we have*

$$\mathrm{Mor}_{\Sigma_S}(\mathrm{O}(W(T)),T') \cong \mathrm{Mor}_{\Sigma_{S,\delta}}(W(T),J_S(T')) \cong \mathrm{Mor}_{\Sigma_S}(T,\mathrm{O}(J_S(T')).$$

*(4) The functors* $\mathrm{O}, \mathrm{J}_S$, W, *and* $Y \mapsto Y_{\mathrm{tf}}$ *are compatible with étale localization.*

*Proof.* The construction of W and $\mathrm{J}_S$ described in Remarks 1.6.5 and 1.6.10 are compatible with $p$-adic completion and globalize to construct the functors in statement (1). The construction in (2) is straightforward. If $Y \in \Sigma_{S,\delta}$, the closure of the sheaf of $p$-torsion elements of $\mathcal{O}_Y$ is compatible with localization and in fact is $p$-completely quasi-coherent. Then the formal subscheme $Y_{\mathrm{tf}}$ of $Y$ it defines is $p$-torsion free (see Corollary 7.3.7). Furthermore, if $x \in \mathcal{O}_Y$ and $p^n x = 0$, then

$$0 = \phi(p^n x) = p^n \phi(x) = p^n x^p + p^{n+1}\delta(x) = p^{n+1}\delta(x).$$

Thus the sheaf of $p$-torsion elements of $\mathcal{O}_Y$ is invariant under $\delta$, and since $\delta$ is $p$-adically continuous, so is the closure of this ideal. Thus $Y_{\mathrm{tf}}$ inherits a $\delta$-structure, hence forms a crystalline prism. It is clear that any morphism from a crystalline prism to $Y$ factors uniquely through $Y_{\mathrm{tf}}$.

Statement (4) for W is a theorem of Van der Kallen [11, Thm.B], and the result for J is explained in [12, 11.1], although not quite in the form we need. Here is a sketch, which we write in the affine case.

Say $A \to A'$ is $p$-completely étale. Then $\mathrm{J}_S(A)\hat{\otimes}_A A'$ is $p$-completely étale over $\mathrm{J}_S(A)$, and hence by [10, 2.18], inherits a $\delta$-structure, which is evidently compatible with that of $S$. We need to check that this $\delta$-structure is universal. Suppose that $A' \to B$ is a homomorphism of $p$-adically complete and separated $A$-algebras and $B$ has a $\delta$-structure compatible with that of $S$. Then the map $A \to B$ extends uniquely to a continuous map of $\delta$-algebras $\mathrm{J}_S(A) \to B$, and hence a map of algebras $\mathrm{J}_S(A)\hat{\otimes}_A A' \to B$. We claim that this map is also compatible with the $\delta$-structures. Recall from [10, 2.4] that the $\delta$-structure on $B$ induces, and is equivalent to, a splitting of the map $W_2(B) \to B$, and similarly for $\mathrm{J}_S(A)$. Since $\mathrm{J}_S(A) \to B$ is a $\delta$-homomorphism, we have a commutative diagram:

$$\begin{array}{ccc} W_2(\mathrm{J}_S(A)) & \longrightarrow & W_2(B) \\ \uparrow & & \uparrow \\ \mathrm{J}_S(A) & \longrightarrow & B. \end{array}$$

Tensoring with $A'$ then gives the central commutative square in the diagram:

$$
\begin{array}{ccccc}
W_2(\mathrm{J}_S(A)\hat{\otimes}A') & & & \searrow & \\
\uparrow\cong & & & & \\
W_2(\mathrm{J}_S(A))\hat{\otimes}_A A' & \longrightarrow & W_2(B)\hat{\otimes}A' & \longrightarrow & W_2(B) \\
\uparrow & & \uparrow & & \uparrow \\
\mathrm{J}_S(A)\hat{\otimes}_A A' & \longrightarrow & B\hat{\otimes}_A A' & \longrightarrow & B \\
\uparrow & & \uparrow & & \\
A' & \longrightarrow & B. & &
\end{array}
$$

Here the top vertical arrow on the left is the isomorphism coming from the fact that formation of $W_2$ is compatible with étale localization. The fact that the circuit commutes then implies that our homomorphism is compatible with the $\delta$-structures. □

**Corollary 1.6.12.** *If $(S, \delta)$ is a crystalline prism and $T/S$ is a $p$-completely smooth morphism of $p$-adic formal schemes, then $\mathrm{J}_S(T)$ is $p$-torsion free, hence a crystalline prism.*

*Proof.* We may assume without loss of generality that $(S, \delta)$ is the formal spectrum of a $\delta$-ring $(R, \delta)$. Suppose that $X$ is a set of variables. Then, $\mathrm{J}_S(R[X]) = R[X, \delta]$, which is again a polynomial ring, as we saw in Remark 1.6.10. This ring is $p$-torsion free, and remains so after $p$-adic completion. Thus $\mathrm{J}_S(T)$ is $p$-torsion free whenever $T := \mathrm{Spf}(R[X]\hat{\ })$. Since any $p$-completely smooth formal scheme over $S$ is locally $p$-completely étale over such a $T$, the result follows from (4) of Proposition 1.6.11. □

**Remark 1.6.13.** The categories $\Sigma_\delta$ and $\Sigma_{\mathrm{F}}$ admit products and fiber products, and their formation commutes with the forgetful functor O. The category $\Delta$ also admits products and fiber products, but the formation of fiber products may not commute with O. Indeed, if $X \to Z$ and $Y \to Z$ are morphisms of crystalline prisms, the fiber product $X \times_Z Y$ in the category of formal $\delta$-schemes may have nontrivial $p$-torsion, and the fiber product in the category of prisms is $(X \times_Z Y)_{\mathrm{tf}}$. This can be difficult to compute and radically different from the fiber product of the underlying formal schemes.

## 1.7 Deformations of $F$-schemes and their subschemes

We begin with a review of some elementary deformation theory. Although these results are surely standard, we explain their proofs for the convenience of the reader.

For $m \in \mathbf{N}$, let $S_m := \operatorname{Spec} \mathbf{Z}/p^m\mathbf{Z}$, suppose that $Y/S_m$ is a flat morphism of schemes, and for $j \leq m$, let $Y_j$ be the reduction of $Y$ modulo $p^j$. Take $n = m-1$, and suppose we are given a closed subscheme $X$ of $Y_n$, flat over $S_n$. A *flat lifting of $X$ in $Y$* is a closed subscheme $\tilde{X} \subseteq Y$ such that $\tilde{X}_n = X$ and such that $\tilde{X}/S_m$ is flat. To describe these, let $I_{X/Y}$ be the ideal of $X$ in $Y$ and let $\mathcal{C}_{X/Y} := I_{X/Y}/I^2_{X/Y}$ be the conormal sheaf. Then $p^n \in I_{X/Y}$, and, since $Y/S_m$ is flat, multiplication by $p^n$ on $\mathcal{O}_Y$ factors through an injection $[p^n]\colon \mathcal{O}_{\overline{Y}} \to I_{X/Y}$. (Recall that $\overline{Y}$ is the reduction of $Y$ modulo $p$.) Then the top row of the following diagram is exact, and the bottom row, its restriction to $X$, is also exact.

$$\begin{array}{ccccccccc}
0 & \longrightarrow & \mathcal{O}_{\overline{Y}} & \xrightarrow{[p^n]} & I_{X/Y} & \longrightarrow & I_{X/Y_n} & \longrightarrow & 0 \\
 & & \downarrow{\scriptstyle \pi} & & \downarrow & & \downarrow & & \\
 & & \mathcal{O}_{\overline{X}} & \xrightarrow{\alpha} & \mathcal{C}_{X/Y} & \longrightarrow & \mathcal{C}_{X/Y_n} & \longrightarrow & 0
\end{array}$$

**Proposition 1.7.1.** *With the notation of the preceding paragraph, there is a bijection*

$$\tilde{X} \mapsto \rho_{\tilde{X}}$$

*from the set of flat liftings of $X$ in $Y$ to the set of $\mathcal{O}_Y$-linear homomorphisms $\rho\colon I_{X/Y} \to \mathcal{O}_{\overline{X}}$ such that $\rho(p^n) = 1$; equivalently, of splittings $\sigma_{\tilde{X}} : \mathcal{C}_{X/Y} \to \mathcal{O}_{\overline{X}}$ of $\alpha$. The inverse of this bijection takes such a homomorphism $\rho$ to the closed subscheme defined by* $\operatorname{Ker}(\rho)$. *Consequently the sheaf of flat liftings of $X$ in $Y$ has a natural structure of a pseudo-torsor under $\mathcal{N}_{\overline{X}/\overline{Y}}$, inherited from the evident pseudo-torsor on the family of splittings of $\alpha$. Explicitly: if $\eta$ is a section of $\mathcal{N}_{\overline{X}/\overline{Y}}$, and $\pi\colon I_{X/Y} \to \mathcal{C}_{\overline{X}/\overline{Y}}$ is the projection then, $\eta \circ \pi(p^n) = 0$, and if $\rho\colon I_{X/Y} \to \mathcal{O}_{\overline{X}}$ is a homomorphism with $\rho(p^n) = 1$, so is $\rho + \eta \circ \pi$.*

*Proof.* The assertion can be checked locally, so we may assume that $Y$ is affine. We will need the following elementary lemma.

**Lemma 1.7.2.** *Let $R_m := \mathbf{Z}/p^m\mathbf{Z}$ and let $M$ be an $R_m$-module. Then the following are equivalent.*

*(1) $M$ is flat.*

*(2)* $\operatorname{Tor}_1^{R_m}(M, R_i) = 0$ *for all* $i \in [1, m-1]$.

*(3)* $\operatorname{Tor}_1^{R_m}(M, R_i) = 0$ *for some* $i \in [1, m-1]$.

*(4)* $\operatorname{Tor}_1^{R_m}(M, R_1) = 0$.

*Proof.* The implications (1) $\Rightarrow$ (2) $\Rightarrow$ (3) are obvious, and (4) $\Rightarrow$ (1) follows from the local criterion for flatness. Thus it suffices to prove (3) $\Rightarrow$ (4).

The exact sequence

$$0 \longrightarrow R_{m-i} \xrightarrow{[p^i]} R_m \longrightarrow R_i \longrightarrow 0$$

yields an exact:

$$0 \to \mathrm{Tor}_1^{R_m}(M, R_i) \longrightarrow M \otimes R_{m-i} \xrightarrow{[p^i]} M \longrightarrow M \otimes R_i \to 0,$$

which shows that $\mathrm{Tor}_1^{R_m}(M, R_i) = 0$ if and only if, for all $x \in M$, $p^i x = 0$ implies that $x \in p^{m-i}M$. Assume that this holds for some $i \in [1, m-1]$; we will show that it also holds for $i-1$. Suppose that $p^{i-1}x = 0$. Then $p^i x = 0$, so $x = p^{m-i}x'$ for some $x'$. Note that $m - i - 1 \geq 0$, since $i < m$. We have

$$0 = p^{i-1}x = p^{m-1}x' = p^i p^{m-i-1}x',$$

so $p^{m-i-1}x' = p^{m-i}x''$ for some $x''$. But then

$$x = p^{m-i}x' = pp^{m-i-1}x' = pp^{m-i}x'' = p^{m-i+1}x''.$$

□

**Corollary 1.7.3.** *Let $M$ be a flat module over $R_m := \mathbf{Z}/p^m\mathbf{Z}$, let $M'$ be a submodule, and let $M'' := M/M'$. Then $M''$ is flat over $R_m$ if and only if $p^iM \cap M' = p^iM'$ for some, equivalently for all, $i \in [1, m-1]$.*

*Proof.* Indeed, the lemma shows that flatness of $M''$ is equivalent to the vanishing of $\mathrm{Tor}_1^{R_m}(M'', R_i)$, and the exact sequence:

$$0 \to \mathrm{Tor}_1^{R_m}(M'', R_i) \to M'/p^iM' \to M/p^iM \to M''/p^iM'' \to 0,$$

shows that $\mathrm{Tor}_1^{R_m}(M'', R_i) \cong (p^iM \cap M')/p^iM'$. □

To prove the proposition, assume that $n := m - 1 > 0$, that $B$ is a flat $R_m$-algebra, and that $J$ is an ideal of $B$ such that $A := B/J$ is an $R_n$-algebra. Then $p^n \in J$, and the flatness of $B$ implies that the sequence:

$$0 \longrightarrow B_1 \xrightarrow{[p^n]} J \longrightarrow JB_n \to 0,$$

is exact. We want to find the ideals $I$ of $B$ such that $IB_n = JB_n$ and such that $B/I$ is flat over $B_m$. The first condition says that $I + p^nB = J$, and the second condition is equivalent to $I \cap p^nB = p^nI$, by Corollary 1.7.3. Note that, since $p^{n+1}B = 0$, the first condition implies that $J^2 \subseteq I$.

Suppose that $I$ is such an ideal. Then

$$J/I = (I + p^nB)/I \cong p^nB/(I \cap p^nB) = p^nB/p^nI,$$

and multiplication by $p^n$ induces an isomorphism:

$$A_1 = B_1/IB_1 \xrightarrow{\cong} p^n B/p^n I \cong J/I.$$

Let $\rho\colon J \to A_1$ be the composition of the projection $J \to J/I$ with the inverse of this isomorphism. Then $I$ is the kernel of $\rho$, and we find a commutative diagram:

$$\begin{array}{ccccccccc}
0 & \longrightarrow & B_1 & \xrightarrow{[p^n]} & J & \longrightarrow & JB_n & \longrightarrow & 0 \\
 & & \downarrow \pi & \rho \swarrow & \downarrow & & \downarrow & & \\
 & & A_1 & \xrightarrow{[p^n]} & J/J^2 & \longrightarrow & JB_n/J^2B_n & \longrightarrow & 0,
\end{array}$$

Since $J^2 \subseteq I$, the map $\rho$ factors through $J/J^2$ and hence provides a splitting of the bottom sequence. In particular the bottom left map is injective.

Suppose conversely that $J$ is an ideal of $B$ containing $p^n$ and that $\rho\colon J \to A_1$ is a homomorphism such that $\rho(p^n) = 1$. Let $K$ be the kernel of $\rho$. We claim that $K + p^n B = J$. Indeed, if $x \in J$, we can find $b \in B_1$ such that $\pi(b) = \rho(x)$, and then $\rho([p^n]b) = \pi(b) = \rho(x)$, so $x - [p^n]b \in \mathrm{Ker}(\rho)$. We next claim that $K \cap p^n B = p^n K$. Indeed, if $b \in B$ and $p^n b \in K$, then

$$0 = \rho(p^n b) = \pi(b)\rho([p^n]) = \pi(b),$$

so $b \in J + pB$, say $b = x + pb'$, with $x \in J$. Write $x = x' + p^n b''$, with $x' \in K$ and $b'' \in B$. Then $p^n b = p^n x' + p^{2n} b'' = p^n x' \in p^n K$, as required. $\square$

The next proposition elucidates and exploits the naturality of the correspondence $\tilde{X} \to \rho_{\tilde{X}}$ explained in Proposition 1.7.1.

**Proposition 1.7.4.** *Let $g\colon Y \to Y'$ be a morphism of flat schemes over $S_m$, and let $n := m - 1$. Suppose we are given a commutative diagram*

$$\begin{array}{ccc}
X & \xrightarrow{i} & Y_n \\
\downarrow f_n & & \downarrow g_n \\
X' & \xrightarrow{i'} & Y'_n,
\end{array}$$

*where $i$ and $i'$ are closed immersions, and where $X$ and $X'$ are flat over $S_n$. Suppose we are also given flat lifts $\tilde{X} \subseteq Y$ and $\tilde{X}' \subseteq Y'$ of $X$ and $X'$ respectively.*

*(1) There is a unique homomorphism $\epsilon(\tilde{X}, \tilde{X}')$ fitting in the commutative*

*diagram below, in which the rightmost vertical arrow is induced by $g$.*

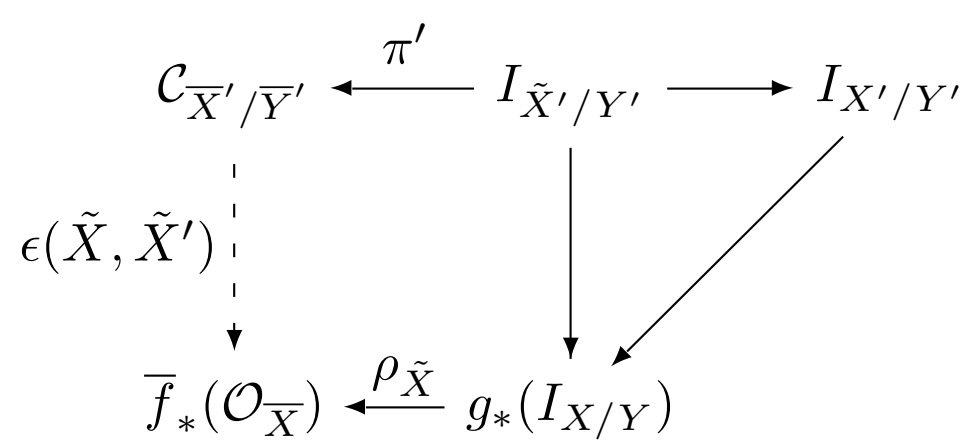

*The morphism $g$ sends $\tilde{X}$ to $\tilde{X}'$ if and only if $\epsilon(\tilde{X}, \tilde{X}')$ vanishes.*

*(2) Suppose that $\eta$ and $\eta'$ are sections of the respective normal sheaves $\mathcal{N}_{\overline{X}/\overline{Y}}$ and $\mathcal{N}_{\overline{X}'/\overline{Y}'}$, and let $\tilde{X}_\eta$ and $\tilde{X}_{\eta'}$ denote the corresponding flat lifts of $X$ and $X'$, as determined by the torsor action described in Proposition 1.7.1. Then*

$$\epsilon(\tilde{X}_\eta, \tilde{X}'_{\eta'}) = \epsilon(\tilde{X}, \tilde{X}') - \overline{f}^* \circ \eta' + \overline{f}_*(\eta) \circ g^*,$$

*where $g^* \colon \mathcal{C}_{\overline{X}'/\overline{Y}'} \to \mathcal{C}_{\overline{X}/\overline{Y}}$ is the homomorphism induced by $g$.*

*(3) There exist liftings of $X$ and $X'$ compatible with $g$ if and only if $\epsilon(\tilde{X}, \tilde{X}')$ lies in the sum of the images of*

$$\eta' \mapsto f^* \circ \eta' \colon \operatorname{Hom}(\mathcal{C}_{\overline{X}'/\overline{Y}'}, \mathcal{O}_{\overline{X}'}) \to \operatorname{Hom}(\mathcal{C}_{\overline{X}'/\overline{Y}'}, f_*(\mathcal{O}_{\overline{X}}))$$

*and*

$$\eta \mapsto \eta \circ g^* \colon \operatorname{Hom}(\mathcal{C}_{\overline{X}/\overline{Y}}, \mathcal{O}_{\overline{X}}) \to \operatorname{Hom}(\mathcal{C}_{\overline{X}'/\overline{Y}'}, f_*(\mathcal{O}_{\overline{X}}))$$

*Proof.* For (1), observe that $g^*$ takes $I_{X'/Y'}$ to $I_{X/Y}$, so we get the triangle in the diagram. The composed map $I_{\tilde{X}'/Y'} \to f_*(\mathcal{O}_{\overline{X}})$ factors through the surjection $\pi'$ because $\overline{f}_*(\mathcal{O}_{\overline{X}})$ is annihilated by $I_{\overline{X}/Y}$. Since $I_{\tilde{X}/Y}$ is the kernel of $\rho_{\tilde{X}}$, we see that $g^*$ takes $I_{\tilde{X}'/Y'}$ to $I_{\tilde{X}/Y}$ if and only if $\epsilon(\tilde{X}, \tilde{X}')$ vanishes.

dri Now suppose that $x' \in I_{\tilde{X}'_{\eta'}/Y'}$. Write $x' = x'' + p^n b'$, with $x'' \in I_{\tilde{X}'/Y'}$ and $b \in \mathcal{O}_{Y'}$. Since $\rho_{\tilde{X}'_{\eta'}}(x')$, $\rho_{\tilde{X}'}(x'')$, and $\eta' \circ \pi'(p^n b)$ all vanish, it follows that

$$\begin{aligned} 0 &= \rho_{\tilde{X}'_{\eta'}}(x') \\ &= (\rho_{\tilde{X}'} + \eta' \circ \pi')(x'' + p^n b') \\ &= \rho_{\tilde{X}'}(p^n b') + \eta' \circ \pi'(x'') \\ &= \overline{b}' + \eta'(\pi'(x'')) \end{aligned}$$

where $\overline{b}'$ is the image of $b'$ in $\mathcal{O}_{X'}$. Then

$$\begin{aligned} \epsilon(\tilde{X}_\eta, \tilde{X}'_{\eta'})(\pi'(x')) &= \rho_{\tilde{X}_\eta}(g^*(x')) \\ &= \rho_{\tilde{X}_\eta}(g^*(x'' + p^n b')) \\ &= \rho_{\tilde{X}_\eta}(g^*(x'')) + \rho_{\tilde{X}_\eta}(p^n g^*(b')) \end{aligned}$$

$$
\begin{aligned}
&= \rho_{\tilde{X}}(g^*(x'')) + \eta \circ \pi(g^*(x'')) + \overline{f}^*(\overline{b}') \\
&= \epsilon(\tilde{X}, \tilde{X}')(g^*(\pi'^*(x''))) + \eta(g^*(\pi')(x'')) - \overline{f}^*(\eta'(\pi'(x'))).
\end{aligned}
$$

Since $\pi'(x') = \pi'(x'')$, this proves the formula in statement (2). Statement (3) is an immediate consequence. □

We shall now apply these consideration to Frobenius liftings.

**Definition 1.7.5.** *Let $Y$ be a crystalline prism and let $X$ be a closed subscheme of $\overline{Y} := Y_1$ and let $n$ be a positive integer. We say that a flat lifting $\tilde{X} \subseteq Y$ is $F_n$-aligned if $F_Y$ maps $\tilde{X}_n$ to itself, and we say that $X$ is $F_n$-alignable if it admits an $F_n$-aligned lifting. We say $X$ is $F$-alignable if $X \to Y$ admits a lifting to a closed immersion of crystalline prisms $\tilde{X} \to Y$.*

In this situation, Proposition 1.7.4 simplifies somewhat to give the following criterion for $F$-alignability.

**Corollary 1.7.6.** *Let $Y$ be a crystalline prism and let $\tilde{X}$ be a closed subscheme of $Y_{n+1}$, flat over $\mathbf{Z}/p^{n+1}\mathbf{Z}$ such that $X := \tilde{X}_n$ is $F_n$-aligned. Let*

$$\epsilon(\tilde{X}, \tilde{X}) \in \operatorname{Hom}(\mathcal{C}_{\overline{X}/\overline{Y}}, F_{\overline{X}*}(\mathcal{O}_{\overline{X}}))$$

*be defined as in (1) of Proposition 1.7.4. Then the obstruction to finding an $F_{n+1}$-aligned lifting of $X_n$ is the image of $\epsilon(\tilde{X}, \tilde{X})$ in the cokernel of the natural map*

$$\operatorname{Hom}(\mathcal{C}_{\overline{X}/\overline{Y}}, \mathcal{O}_{\overline{X}}) \to \operatorname{Hom}(\mathcal{C}_{\overline{X}/\overline{Y}}, F_{\overline{X}*}(\mathcal{O}_{\overline{X}})),$$

*which, if $\overline{X}$ is regularly embedded in $\overline{Y}$, can be identified with the map*

$$\mathcal{N}_{\overline{X}/\overline{Y}} \to F_{\overline{X}*}F^*_{\overline{X}}(\mathcal{N}_{\overline{X}/\overline{Y}}).$$

*Proof.* This follows from (3) of Proposition 1.7.4 and the fact that

$$F_Y^* \colon \mathcal{C}_{\overline{X}/\overline{Y}} \to \mathcal{C}_{\overline{X}/\overline{Y}}$$

is the zero map, since $F_Y^*(I_{\overline{X}}) \subseteq I^p_{\overline{X}}$. If $\overline{X}$ is regularly embedded in $\overline{Y}$, then $\mathcal{C}_{\overline{X}/\overline{Y}}$ is locally free of finite rank, so $\mathcal{N}_{\overline{X}/\overline{Y}}$ shares this property, and

$$\mathcal{H}om(\mathcal{C}_{\overline{X}/\overline{Y}}, F_{\overline{X}*}(\mathcal{O}_{\overline{X}})) \cong F_{\overline{X}*}F^*_{\overline{X}}(\mathcal{N}_{\overline{X}/\overline{Y}}).$$

□

**Example 1.7.7** Let $Y := \operatorname{Spf}(W[x,y]\hat{\ })$, with $\phi(x) = x^p + py$ and $\phi(y) = y^p$, and let $X \subseteq Y$ be defined by $(p, x)$. Then $X$ is not $F_2$-alignable, *i.e.*, there is no lift of $X$ to $Y_2$ which is invariant under $F_Y$. Indeed, any such lift would be defined by $(p^2, x + pf)$ for some $f$ such that $\phi(x + pf) \in (p^2, x + pf)$. Say $\phi(x + pf) = (x + pf)g \pmod{p^2}$, with $g \in W[x,y]\hat{\ }$. Then $x^p + py + p\phi(f) \equiv (x + pf)g \pmod{p^2}$, and so $g \equiv x^{p-1} \pmod{p}$. Write $g = x^{p-1} + ph$, and deduce that $y + \phi(f) = hx + x^{p-1}f \pmod{p}$. But this implies that $y + f^p \equiv 0 \pmod{(x,p)}$, which is impossible. In fact, the element $\epsilon$ described in Corollary 1.7.6 is the map $\mathcal{C}_{\overline{X}/\overline{Y}} \to F^*(\mathcal{O}_X)$ taking the class of $x$ to the class of $y$, and its image in $F_{\overline{X}*}F^*_{\overline{X}}(\mathcal{N}_{\overline{X}/\overline{Y}})/\mathcal{N}_{\overline{X}/\overline{Y}}$ is nonzero.

**Example 1.7.8** Let us also note that if $Y/S$ is a flat formal scheme over $\mathrm{Spf}(W)$ endowed with two Frobenius lifts $F_Y$ and $F'_Y$, then the diagonal embedding of $\overline{Y}$ in $(Y \times_S Y, F_Y \times F'_Y)$, may not be $F_2$-alignable. For example, let $Y = \mathrm{Spf}\, W[x]\hat{}$ with $\phi(x) = x^p$ and $\phi'(x) = x^p + px$. Then the ideal of the diagonal is generated by $\xi := 1 \otimes x - x \otimes 1$, and $(\phi, \phi')$ sends $x \otimes 1$ to $x^p \otimes 1$ and $\xi$ to

$$1 \otimes x^p + 1 \otimes px - x^p \otimes 1 = \xi \left( \sum_{i=1}^{p-1} x^i \otimes x^{p-i} \right) + p \otimes x.$$

Thus $\epsilon$ is the map taking $\xi$ to $x$, which lies in $F_{\overline{Y}*} F^*_{\overline{Y}} (\mathcal{N}_{\overline{Y}/\overline{Y} \times_S \overline{Y}})$ but not in $\mathcal{N}_{\overline{Y}/\overline{Y} \times_S \overline{Y}}$.

In the next corollary we draw some conclusions about Frobenius lifts and $F$-alignment. In the first part, we "decouple" the relationship between $X_n$ and $X_n^{(p)}$, in that we don't require that $X_n^{(p)}$ be obtained by base change from $X_n$.

**Corollary 1.7.9.** *Let $Y/S$ be a morphism of crystalline prisms, where $\overline{S}$ is regular, let $X$ be a reduced closed subscheme of $\overline{Y}$, and let $X^{(p)} \subseteq \overline{Y}^{(p)}$ be its pullback by $F_{\overline{S}}$, as in the relative Frobenius diagram (1.1.3).*

*(1) For every $n$, there is at most one flat lifting $X_n^{(p)}$ of $X^{(p)}$ in $Y^{(p)}$ for which there exist some flat lifting $X_n$ of $X$ such that $F_{Y/S}$ maps $X_n$ to $X_n^{(p)}$.*

*(2) If $X \subseteq \overline{Y}$ is (Zariski) locally $F_n$-alignable, then it is $F_n$-alignable. If $X$ is locally $F$-alignable, then it is $F$-alignable.*

*Proof.* Suppose that $X_n, \tilde{X}_n \subseteq Y_n$ and $X_n^{(p)}, \tilde{X}_n^{(p)} \subseteq Y_n^{(p)}$ are flat liftings of $X$ and $X^{(p)}$ respectively and that $F_{Y/S}$ takes $X_n$ to $X_n^{(p)}$ and $\tilde{X}_n$ to $\tilde{X}_n^{(p)}$. We prove that $X_m^{(p)} = \tilde{X}_m^{(p)}$ by induction on $m$. Assuming that this is true for $m$, let $\eta$ (resp. $\eta'$) measure the difference between $X_{m+1}$ and $\tilde{X}_{m+1}$ (resp. between $X_{m+1}^{(p)}$and $\tilde{X}_{m+1}^{(p)}$). Statement (2) of Proposition 1.7.4 says that

$$\epsilon(\tilde{X}_{m+1}, \tilde{X}_{m+1}^{(p)}) = \epsilon(X_{m+1}, X_{m+1}^{(p)}) - F^*_{Y/S} \circ \eta' + F_{X/S*}(\eta) \circ F^*_{Y/S}. \qquad (1.7.1)$$

Both the $\epsilon$-terms vanish, as does $F^*_{Y/S}$, and it follows that $F^*_{X/S} \circ \eta' = 0$. Since $X$ is reduced, in fact $F^*_{X/S}$ is injective, and it follows that $\eta'$ also vanishes. Hence $\tilde{X}_{m+1}^{(p)} = X_{m+1}^{(p)}$.

Before proving (2), let us translate (1) into a uniqueness statement for $F_n$-aligned liftings. Suppose that $X_n$ and $\tilde{X}_n$ are $F_n$-aligned liftings of $X$ in $Y_n$. Then (1) implies that their Frobenius pullbacks $X_n^{(p)}$ and $\tilde{X}_n^{(p)}$ in $Y_n^{(p)}$ agree. Since $\overline{S}$ is regular, the Frobenius lifting $F_S$ of $S$ is faithfully flat, hence $Y_n^{(p)} \to Y_n$ is faithfully flat, and this implies that $X_n$ and $\tilde{X}_n$ also agree. If $X_n$ is locally $F_n$-alignable, it follows that local $F_n$-aligned liftings agree on overlaps and hence patch together to a global aligned lifting. The statement for formal schemes follows immediately. ☐

We note that in statement (1) of Corollary 1.7.9, the lifting $X_n$ of $X$ is not unique. To the contrary, if $m < n$ and if $F_{X/S}$ maps a lift $X_{m+1}$ of $X$ to a lift $X^{(p)}_{m+1}$ of $X^{(p)}$, then it also takes $\tilde{X}_{m+1}$ to $X^{(p)}_{m+1}$ for every lifting $\tilde{X}_{m+1}$ of $X_m$. Indeed, if $F_{X/S}$ maps $X_{m+1}$ to $X^{(p)}_{m+1}$, then since $F^*_{Y/S}$ and $\epsilon(X_{m+1}, X^{(p)}_{m+1})$ in equation (1.7.1) both vanish, we see that $\epsilon(\tilde{X}_{m+1}, X^{(p)}_{m+1}) = 0$ when $\eta' = 0$ for every choice of $\eta$. For an explicit example, let $Y := \operatorname{Spf} W[x]\hat{}$, and let $Y' := \operatorname{Spf} W[x']\hat{}$, with $g^\sharp(x') = x^p$, with $X$ and $X'$ defined by $(x, p)$ and $(x', p)$ respectively. Suppose $X'_2$ is defined by $(x', p^2)$ and $\tilde{X}_2$ is defined by $(\tilde{x}, p^2)$, with $\tilde{x} := x + pf$ for some $f \in W[x]\hat{}$. Then

$$g^\sharp(x') \equiv -p^p f^p \pmod{\tilde{x}}$$

and hence belongs to $(\tilde{x}, p^2)$ for every $f$. Note, however, that further lifts may not be possible. For example, if $X_2$ is defined by $(x + p, p^2)$, there is no flat lift of $X_2$ in $Y_{p+1}$ which maps to $X^{(p)}_{p+1}$. Indeed, such a lift would be defined by $(p^{p+1}, \tilde{x})$ with $\tilde{x} = x + pf$ with $f \equiv 1 \pmod{p^2}$. Then $f^p \equiv 1 \pmod{p^2}$ and $g^\sharp(x') \equiv -p^p f^p \pmod{\tilde{x}}$ and so does not belong to $(\tilde{x}, p^{p+1})$.

The following result shows that morphisms of crystalline prisms are remarkably rigid.

**Proposition 1.7.10.** *Let $T$ and $Y$ be crystalline prisms over a crystalline prism $S$, and assume that $\overline{T}$ is reduced. Then a morphism crystalline prisms $f\colon T \to Y$ over $S$ is uniquely determined by its restriction to $\overline{T}$.*

*Proof.* It is enough to prove that for each $n \geq 1$, the restriction of $f$ to $T_{n+1}$ is uniquely determined by its restriction $f_n$ to $T_n$. Suppose then that $\tilde{f}_{n+1}\colon T_{n+1} \to Y$ is a morphism compatible with the Frobenius lifts of $T_{n+1}$ and $Y$ such that $\tilde{f}_n = f_n$. Then there is a derivation $D\colon \mathcal{O}_{Y_1} \to \mathcal{O}_{T_1}$ such that

$$\tilde{f}^\sharp_{n+1} - f^\sharp_{n+1} = [p^n] \circ D \circ \pi,$$

where $\pi\colon \mathcal{O}_{Y_{n+1}} \to \mathcal{O}_{Y_1}$ is the projection and $[p^n]\colon \mathcal{O}_{T_1} \to \mathcal{O}_{T_{n+1}}$ is induced by multiplication by $p^n$. Since $\tilde{f}_{n+1}$ and $f_{n+1}$ are both compatible with the Frobenius liftings, we have:

$$\begin{aligned}
\phi_T \circ \tilde{f}^*_{n+1} &= \tilde{f}^*_{n+1} \circ \phi_Y \\
\phi_T \circ (f^*_{n+1} + [p^n] \circ D \circ \pi) &= (f^*_{n+1} + [p^n] \circ D \circ \pi) \circ \phi_Y \\
\phi_T \circ f^*_{n+1} + \phi_T \circ [p^n] \circ D \circ \pi &= f^*_{n+1} \circ \phi_Y + [p^n] \circ D \circ \pi \circ \phi^*_Y \\
\phi_T \circ [p^n] \circ D \circ \pi &= [p^n] \circ D \circ \pi \circ \phi_Y \\
F^*_{T_1} \circ D &= D \circ F^*_{Y_1}.
\end{aligned}$$

Since $D$ is a derivation, it follows that $D \circ F^*_{Y_1} = 0$, and since $T_1$ is reduced, $F^*_{T_1}$ is injective. Consequently $D = 0$ and $\tilde{f}_{n+1} = f_{n+1}$. □

The following result is inspired by [28, 3.4], which it generalizes. It will play a key role in the construction (6.1.11) of coverings of the final object in the prismatic topos and hence in the computation of prismatic cohomology.

Let us recall from [10] that a morphism of $p$-adic formal schemes is said to be *$p$-completely smooth* if it is $p$-completely flat and its reduction modulo $p$ is smooth.

**Theorem 1.7.11.** *Let $(Y, F_Y) \to (S, F_S)$ and $(T, F_T) \to (S, F_S)$ be morphisms of crystalline prisms, where the underlying morphism $Y/S$ is $p$-completely smooth.*

*(1) Suppose that $n > 0$ and that we are given an $S$-morphism $f_n \colon T_n \to Y$ which is compatible with the Frobenius liftings. Then there exists a $p$-completely faithfully flat morphism of crystalline prisms: $u \colon (\widetilde{T}, F_{\widetilde{T}}) \to (T, F_T)$ and a morphism of $S$-$F$-schemes $\tilde{f} \colon (\widetilde{T}, F_{\widetilde{T}}) \to (Y, F_Y)$ whose restriction to $\widetilde{T}_n$ agrees with the map $f_n \circ u_n$.*

*(2) In the situation of statement (1), suppose in addition that $i \colon Z \to T$ is a closed immersion of crystalline prisms and that we are given a morphism $g \colon Z \to Y$ which is compatible with the Frobenius lifts and such that $f_n \circ i_n = g_n$. Then there exists a commutative diagram of F-schemes*

$$\begin{array}{ccccc}
\tilde{Z} & \xrightarrow{\tilde{i}} & \widetilde{T} & & \\
\downarrow v & & \downarrow u & \searrow \tilde{f} & \\
Z & \xrightarrow{i} & T & & \\
& & \searrow g & & Y
\end{array}$$

*such that $\tilde{f}_n = f_n \circ u_n$, $\tilde{i}$ is a closed immersion, and the morphisms $u$ and $v$ are $p$-completely faithfully flat. (Note: the square in this diagram may not be Cartesian.)*

*Proof.* For statement (1), we first show that, after a suitable cover $\widetilde{T} \to T$, we can find a lift $\tilde{f}_{n+1}$ of $f_n$. Replacing $T$ by an open affine cover, which is certainly $p$-completely flat (although not necessarily flat!), we may and shall assume that $T$ is affine. Since the ideal $p^n \mathcal{O}_{T_{n+1}}$ of $T_n$ in $T_{n+1}$ is a square zero ideal and $Y/S$ is formally smooth, there exists an $S$-morphism $f_{n+1} \colon T_{n+1} \to Y$ extending $f_n$, not necessarily compatible with the Frobenius liftings.

Consider the following diagram, where $F_{T/S} \colon T \to T^{(p)}$ and $F_{Y/S} \colon Y \to Y^{(p)}$ are the relative lifted Frobenius maps.

$$\begin{array}{ccccc}
T_{n+1} & \xrightarrow{F_{T/S}} & T^{(p)}_{n+1} & \longrightarrow & T_{n+1} \\
\downarrow f_{n+1} & & \downarrow f'_{n+1} & & \downarrow f_{n+1} \\
Y & \xrightarrow{F_{Y/S}} & Y^{(p)} & \longrightarrow & Y
\end{array}$$

Here the right square is commutative, but the left square may not be. Since it does commute after restriction to $T_n$, the maps $f'_{n+1} \circ F_{T/S}$ and $F_{Y/S} \circ f_{n+1}$ from $T_{n+1}$ to $Y^{(p)}$ differ by an element $\epsilon$ of

$$\begin{aligned} \mathrm{Hom}(f_{n+1}^* F_{Y/S}^*(\Omega^1_{Y'/S}), p^n \mathcal{O}_{T_{n+1}}) &\cong \mathrm{Hom}(f_{n+1}^* F_{Y/S}^*(\Omega^1_{Y'/S}), \mathcal{O}_{T_1}) \\ &\cong \mathrm{Hom}(f_1^* F_{Y_1/S}^*(\Omega^1_{Y_1'/S}), \mathcal{O}_{T_1}) \\ &\cong \mathrm{Hom}(F_{T_1}^* f_1^*(\Omega^1_{Y_1/S}), \mathcal{O}_{T_1}) \\ &\cong \mathrm{Hom}(f_1^*(\Omega^1_{Y_1/S}), F_{T_1 *}(\mathcal{O}_{T_1})). \end{aligned}$$

Since we are allowed to replace $T$ by an affine cover, we may assume that $f_1^*(\Omega^1_{Y/S})$ is freely generated by elements $\omega_1, \ldots, \omega_N$. By [10, 2.12], which we restate and reprove as Proposition 2.3.10 below, we may find a $p$-completely faithfully flat map of crystalline prisms $u\colon \widetilde{T} \to T$ and sections $\tilde{b}_1, \ldots, \tilde{b}_N$ of $\mathcal{O}_{\widetilde{T}}$ such that $u^\sharp(\epsilon(\omega_i)) \equiv \phi^\sharp_{\widetilde{T}}(\tilde{b}_i) \pmod p$ for all $i$. These sections define a map $\delta\colon u^* f_1^* \Omega^1_{Y_1/S} \to \mathcal{O}_{\widetilde{T}_1}$ such that $\phi^*_{\widetilde{T}}(\delta) = u^*(\epsilon)$. Let $\tilde{f}_{n+1}\colon \widetilde{T}_{n+1} \to Y$ be the $S$-morphism corresponding to $f_{n+1} \circ u + \delta$. Then $\tilde{f}_{n+1} \circ \phi_{\widetilde{T}}$ corresponds to

$$f_{n+1} \circ u \circ \phi_{\widetilde{T}} + F^*_{\widetilde{T}_1}(\delta) = f_{n+1} \circ \phi_T \circ u + u^*(\epsilon) = \phi_Y \circ f_{n+1} \circ u,$$

while $\phi_Y \circ \tilde{f}_{n+1} \circ u$ corresponds to

$$\phi_Y \circ f_{n+1} \circ u + F^*_{Y_1}(\delta) = \phi_Y \circ f_{n+1} \circ u,$$

since $F^*_{Y_1}\colon \Omega^1_{Y_1/S} \to \Omega^1_{Y_1/S}$ is the zero map. Thus we get that $\phi_Y \circ \tilde{f}_{n+1} = \tilde{f}_{n+1} \circ \phi_T$, as desired.

Repeating this process, we find for $m \geq n$, affine and $p$-completely faithfully flat morphisms of crystalline prisms $u^{(m)}\colon T^{(m+1)} \to T^{(m)}$ and Frobenius compatible maps $f^{(m)}\colon T^{(m)}_m \to Y$, such that $f^{(m+1)}_m = f^{(m)} \circ u^{(m)}_m$. For each $m \geq n$, let $T^{(\infty)}_m := \varinjlim_{m' \geq m} T^{(m')}_m$ (in the category of schemes). Then $(T^{(\infty)}_{m'})_m = T^{(\infty)}_m$ for $m' \geq m \geq n$, and the maps $f^{(m')}_m$ induces morphisms $T^{(\infty)}_m \to Y$ which are compatible with Frobenius and each other. Then $\widetilde{T} := \varinjlim_m T^{(\infty)}_m$ is a $p$-adic formal scheme whose reduction modulo $p$ is $T^{(\infty)}_m$, by Proposition 7.1.2, and the morphism $\tilde{f} := \varinjlim_m f^{(\infty)}_m : \widetilde{T} \to Y$ is compatible with Frobenius lifts. Since each $\widetilde{T}_m \to T_m$ is faithfully flat, the morphism $\widetilde{T} \to T$ is $p$-completely faithfully flat, by Proposition 7.3.2. Furthermore, it follows from Proposition 7.3.3 that $\widetilde{T}$ is $p$-torsion free, hence is indeed a crystalline prism.

Let us sketch how to modify the proof of statement (1) to obtain statement (2). We again begin by assuming that $T$ is affine. The morphisms $g_{n+1}\colon Z_{n+1} \to Y$ and $f_n\colon T_n \to Y$ agree on $Z_{n+1} \cap T_n = Z_n$, and hence glue to define an $S$-morphism $h_n\colon Z_{n+1} \cup T_n \to Y$, again compatible with the Frobenius lifts. The ideal $I''$ of $Z_{n+1} \cup T_n$ in $T_{n+1}$ is contained in the ideal $p^n \mathcal{O}_{T_{n+1}}$ of $T_n$ in $T_{n+1}$ and, in particular, is a square zero ideal. Since $Y/S$ is formally smooth, there exists an $S$-morphism $f_{n+1}\colon T_{n+1} \to Y$ extending $h_n$, not necessarily compatible with the Frobenius liftings. Then the element $\epsilon$ measuring this incompatibility considered in the proof of statement (1) in fact belongs to

$$\mathrm{Hom}(f_{n+1}^* F_{Y/S}^*(\Omega^1_{Y'/S}), I'').$$

Since $Z$ is $p$-torsion free, in fact $I'' = p^n \mathcal{O}_{T_{n+1}} \cap I_{Z_{n+1}/T_{n+1}} \cong I_{Z_1/T_1}$, so the elements $b_i := \epsilon(\omega_i)$ belong to $F_{T_1 *}(I_{Z_1/T_1})$. Choose lifts $b'_i \in \phi_{T*}(I_{Z/T})$, and let $u\colon \tilde{T} \to T$ be the universal morphism of crystalline prisms endowed with sections $\tilde{b}_1, \ldots, \tilde{b}_N$ such that each $u^\sharp(b'_i) = \phi^\sharp(\tilde{b}_i)$ as described in Proposition 2.3.10 below. This morphism is $p$-completely faithfully flat, and, since each $i^\sharp(b'_i) = 0$, admits a section $\tilde{i}\colon Z \to \tilde{T}$. Then we can take $v := \mathrm{id}_Z$ to form the diagram in statement (2). Now let $\tilde{f}_{n+1}\colon \tilde{T}_{n+1} \to Y$ be the $S$-morphism corresponding to $f_{n+1} \circ u + \delta$. Just as before, we see that $\phi_Y \circ \tilde{f}_{n+1} \circ u = \tilde{f}_{n+1} \circ u \circ \phi_T$ as desired. Furthermore, since each $\tilde{i}^\sharp(\tilde{b}_i) = 0$, it follows that $\tilde{f}_{n+1} \circ \tilde{i} = g_{n+1} \circ v$, as claimed. The rest of the proof proceeds as before. □

Theorem 1.7.11 suggests that a morphism of F-schemes whose underlying morphism of formal schemes is formally smooth should satisfy an infinitesimal lifting property, provided one works locally in the $p$-completely flat topology. Although we shall not need this result here, let us take the time to formulate it precisely. The proof is immediate from the theorem.

**Corollary 1.7.12.** *Let $i\colon Z \to T$ be a nilpotent immersion of crystalline $S$-prisms and let $g\colon Z \to Y$ and $Y \to S$ be morphisms of crystalline $S$-prisms, where the underlying morphism $Y \to S$ is $p$-completely smooth. Then there exists a commutative diagram as in statement (2) of Theorem 1.7.11.* □

We now investigate more closely the relationship between subschemes of crystalline prisms and $\delta$-structures.

**Proposition 1.7.13.** *Let $Y$ be a crystalline prism.*

*(1) If $I \subseteq \mathcal{O}_Y$ is a sheaf of ideals, then the $\delta$-structure of $Y$ induces a $\phi$-linear map*

$$\delta_I\colon I/I^2 \to \mathcal{O}_Y/(p\mathcal{O}_Y + I^2).$$

*(2) If $X$ is a closed subscheme of $\overline{Y}$, then the $\delta$-structure of $Y$ induces a linear map:*

$$\overline{\delta}_{X/Y}\colon F_X^*(\mathcal{C}_{X/Y}) \to \mathcal{O}_X.$$

*This map sends the class of $1 \otimes p \in F_X^*(\mathcal{C}_{X/Y})$ to $1 \in \mathcal{O}_X$.*

*(3) Suppose that $X \to \overline{Y}$ is a regular immersion and $\tilde{X} \subseteq Y_2$ is a flat lifting of $X$, and let $\sigma'_{\tilde{X}}\colon \mathcal{C}_{X/\overline{Y}} \to \mathcal{C}_{X/Y}$ be the splitting associated to the splitting $\sigma_{\tilde{X}}\colon \mathcal{C}_{X/Y} \to \mathcal{O}_X$ of Proposition 1.7.1.*

*(a) The composition $\overline{\delta}_{X/Y} \circ F_X^*(\sigma'_{\tilde{X}})$ is the obstruction $\epsilon(\tilde{X}, \tilde{X})$ to the $F$-invariance of $\tilde{X}$ defined in Proposition 1.7.4.*

(b) *If $\tilde{X}$ is $F_Y$-invariant, then the $\delta$-structure of $Y$ induces a commutative diagram of $\mathcal{O}_X$-linear maps with exact rows:*

$$\begin{array}{ccccccccc}
0 & \longrightarrow & \mathcal{O}_X & \xrightarrow{[p^2]} & F_X^*(i_X^*(\mathcal{C}_{\tilde{X}/Y})) & \longrightarrow & F_X^*(\mathcal{C}_{X/\overline{Y}}) & \longrightarrow & 0 \\
 & & \downarrow F_X^* & & \downarrow \overline{\delta}_{\tilde{X}/Y} & & \downarrow \overline{\delta}_{\tilde{X}/\overline{Y}} & & \\
0 & \longrightarrow & \mathcal{O}_X & \underset{[p]}{\overset{\sigma_{\tilde{X}}}{\rightleftarrows}} & \mathcal{C}_{X/Y} & \overset{\sigma'_{\tilde{X}}}{\rightleftarrows} & \mathcal{C}_{X/\overline{Y}} & \longrightarrow & 0.
\end{array}$$

(4) *If $Y \to Z$ is a closed immersion of crystalline prisms, the $\delta$-structure of $Z$ induces linear maps:*

$$\begin{array}{rcl}
\delta_{Y/Z} \colon F_Y^*(\mathcal{C}_{Y/Z}) & \to & \mathcal{C}_{Y/Z} \\
\delta_{\overline{Y}/\overline{Z}} \colon F_{\overline{Y}}^*(\mathcal{C}_{\overline{Y}/\overline{Z}}) & \to & \mathcal{C}_{\overline{Y}/\overline{Z}}.
\end{array}$$

*Proof.* If $x, y \in I$, we have $\delta(x+y) = \delta(x) + \delta(y) - C(x,y)$ where $C(x,y) \in I^2$. It follows that the map $I \to \mathcal{O}_Y/I^2$ induced by $\delta$ is additive. Furthermore, the equation $\delta(xy) = \delta(x)y^p + x^p\delta(y) + p\delta(x)\delta(y)$ shows that $\delta(xy)$ belongs to $I^2 + p\mathcal{O}_Y$ if $x, y \in I$, and, by additivity, the same holds for any section of $I^2$. Thus we find the additive map $\delta_I$ in statement (1). Furthermore, if $x \in \mathcal{O}_Y$ and $y \in I$, then

$$\delta(xy) = \phi(x)\delta(y) + \delta(x)y^p \equiv \phi(x)\delta(y) \pmod{I^2},$$

shows that $\delta_I$ is in fact $\phi$-linear. Statement (1) follows.

Applying this to the ideal $I_{X/Y}$ we find a $\phi$-linear map

$$\mathcal{C}_{X/Y} = I_{X/Y}/I_{X/Y}^2 \xrightarrow{\delta_{I_{X/Y}}} \mathcal{O}_{\overline{Y}}/I_{X/\overline{Y}}^2 \longrightarrow \mathcal{O}_X,$$

and hence the map $\overline{\delta}_{X/Y}$ in (2). Note that this map does not factor through $\mathcal{C}_{X/\overline{Y}}$, since it sends the class of $p$ to 1.

Now suppose that $\tilde{X} \subseteq Y_2$ is a flat lifting, let $x$ be a local section of $I_{\tilde{X}/Y}$, and let $\overline{x}$ be its image in $\mathcal{C}_{X/\overline{Y}}$. Then $\sigma'_{\tilde{X}}(\overline{x})$ is the image of $x$ in $\mathcal{C}_{X/Y}$. Since $X$ is $F_Y$-invariant, we can write $\phi(x) = x' + py$ with $x' \in I_{\tilde{X}/Y}$ and $y \in \mathcal{O}_Y$, and $\epsilon(\tilde{X}, \tilde{X})(\overline{x})$ is by definition the image $\overline{y}$ of $y$ in $F_{X*}(\mathcal{O}_X)$. But in fact $\phi(x) = x^p + p\delta(x)$, so

$$\epsilon(\tilde{X}, \tilde{X})(\overline{x}) = \overline{\delta}(x) = \overline{\delta}_{X/Y}(1 \otimes \sigma'_{\tilde{X}}(\overline{x})).$$

The top row of the sequence in the diagram of (3b) is the restriction to $X$ of the exact sequence of locally free $\mathcal{O}_{\tilde{X}}$-modules:

$$0 \longrightarrow \mathcal{O}_{\tilde{X}} \xrightarrow{[p^2]} \mathcal{C}_{\tilde{X}/Y} \longrightarrow \mathcal{C}_{\tilde{X}/Y_2} \longrightarrow 0.$$

If $\tilde{X}$ is invariant under $F_Y$, then $\delta(x) \in I_{X/Y}$ whenever $x \in I_{\tilde{X}/Y}$, so we find an additive map $I_{\tilde{X}/Y} \to I_{X/Y}/I^2_{X/Y}$. If $x, y \in I_{\tilde{X}/Y}$, then

$$\delta(xy) = \delta(x)y^p + x^p\delta(y) + p\delta(x)\delta(y) \in I^2_{X/Y},$$

so this map factors through $I_{\tilde{X}/Y}/I^2_{\tilde{X}/Y}$. Furthermore, if $x \in I_{\tilde{X}/Y}$, then

$$\delta(px) = \delta(p)x^p + p\delta(x) \in I^2_{X/Y},$$

so the map also annihilates $(pI_{\tilde{X}/Y})/I^2_{\tilde{X}/Y}$. Thus it factors through $i^*_X(\mathcal{C}_{\tilde{X}/Y})$ and hence induces the arrow $\overline{\delta}_{\tilde{X}/Y}$ in the diagram. If $f \in \mathcal{O}_Y$, then

$$\delta(p^2 f) = \delta(p^2)f^p + p^2\delta(f) \equiv pf^p \pmod{p^2\mathcal{O}_Y}.$$

This implies that $[p] \circ F^*_X = \overline{\delta}_{\tilde{X}/Y} \circ [p^2]$ in the left square of the diagram, and that the right square of the diagram exists. Since $\sigma_{\tilde{X}}$ comes from the isomorphism $\mathcal{O}_X \xrightarrow{[p]} I_{X/Y}/I_{\tilde{X}/Y}$, it also shows that $\sigma_{\tilde{X}} \circ \overline{\delta}_{\tilde{X}/Y} \circ [p^2] = F^*_X$.

Finally, suppose that $Y \to Z$ is a closed immersion of crystalline prisms. We saw above that $\delta_Z$ induces an additive map $I_{Y/Z} \to \mathcal{O}_Z/I^2_{Y/Z}$, which in this case lands in $I_{Y/Z}/I^2_{Y/Z}$. Moreover, if $x, y \in I_{Y/Z}$, then $\delta(xy) \in I^2_{Y/Z}$, so our map factors through a $\phi$-linear endomorphism of $\mathcal{C}_{Y/Z}$. Since $Y \to Z$ is a closed immersion of crystalline prisms, it follows that $p\mathcal{O}_Z \cap I_{Y/Z} = pI_{Y/Z}$, hence that $\mathcal{C}_{\overline{Y}/\overline{Z}} = I_{Y/Z}/(pI_{Y/Z} + I^2_{Y/Z})$, and since $\delta(py) \equiv p\delta(y) \pmod{I^2_{Y/Z}}$ for $y \in I_{Y/Z}$, we also find a $\phi$-linear endomorphism of $\mathcal{C}_{\overline{Y}/\overline{Z}}$, hence statement (4). □

Some of the formal properties of $\delta$-structures resemble those of divided power structures[3]. We shall find it useful to adapt an analog of the PD-filtration defined by a divided power ideal $J$ into our context. Recall that $J^{[n]}$ is the ideal generated by elements of the form $z_1^{[n_1]} \dots z_r^{[n_r]}$ such that each $z_i \in J$ and $n_1 + \cdots + n_r \geq n$. Equivalently, it is the ideal generated by elements of the form $(z_1^{[p^{m_1}]})^{n_1}(z_2^{[p^{m_2}]})^{n_2} \cdots (z_r^{[p^{m_r}]})^{n_r}$ such that $\sum n_i p^{m_i} \geq n$. This motivates the following definition.

**Definition 1.7.14.** *Let $B$ be a $\delta$-ring and $I$ an ideal of $B$ which is invariant under $\delta$. Then $\mathrm{Fil}^n_\delta$ is the ideal of $B$ generated by elements of the form*

$$(\delta^{m_1}(z_1))^{n_1}(\delta^{m_2}(z_2))^{n_2} \cdots (\delta^{m_r}(z_r))^{n_r}$$

*such that $p^{m_1}n_1 + p^{m_2}n_2 + \cdots p^{m_r}n_r \geq n$ and each $z_i \in I$.*

Thus $\mathrm{Fil}^1_\delta = I$ and $\mathrm{Fil}^2_\delta$ is the ideal generated by $I^2$ and $\delta(I)$.

**Proposition 1.7.15.** *Let $I$ be a $\delta$-invariant ideal in a $\delta$-ring $B$.*

*(1)* $\mathrm{Fil}^i_\delta \mathrm{Fil}^j_\delta \subseteq \mathrm{Fil}^{i+j}_\delta$*, and* $\delta(\mathrm{Fil}^n_\delta) \subseteq \mathrm{Fil}^{pn}_\delta$ *if $n > 0$.*

[3]This would be more apparent if the definition involved division by $p!$ instead of $p$, which certainly could have been done.

*(2) If $b \in B$ and $z \in I$, then*

$$\delta(b+z) \equiv \delta(b) - b^{p-1}z \pmod{\mathrm{Fil}^2_\delta}.$$

*Proof.* It is clear from the definitions that $\mathrm{Fil}^i_\delta \mathrm{Fil}^j_\delta \subseteq \mathrm{Fil}^{i+j}_\delta$. First note that, if $a, b \in \mathrm{Fil}^n_\delta$ then $C(a,b) \in \mathrm{Fil}^{pn}_\delta$, since $C(a,b)$ a homogeneous polynomial of degree $p$ in $a, b$. Since

$$\delta(a+b) = \delta(a) + \delta(b) - C(a,b),$$

it too will belong to $\mathrm{Fil}^{pn}_\delta$ provided that $\delta(a)$ and $\delta(b)$ do. Hence, to prove the inclusion $\delta(\mathrm{Fil}^n_\delta) \subseteq \mathrm{Fil}^{pn}_\delta$, it suffices to prove that, if $z$ can be written as

$$z = (\delta^{m_1}(z_1))^{n_1} \cdots (\delta^{m_r}(z_r))^{n_r}$$

with $z_i \in I$ and $n := \sum p^{m_i} n_i$, then $\delta(z)$ belongs to $\mathrm{Fil}^{pn}_\delta$. We prove this by induction on $\sum n_i$. If this is 1, then $z = \delta^{m_1}(z_1)$ and $n = p^{m_1}$, so

$$\delta(z) = \delta^{m_1+1}(z_1) \in \mathrm{Fil}^{p^{m+1}}_\delta = \mathrm{Fil}^{pn}_\delta.$$

If $\sum n_i > 1$, then $z$ can be written as $z'z''$, where $z' \in \mathrm{Fil}^{n'}_\delta$, $z'' \in \mathrm{Fil}^{n''}_\delta$ with $n = n' + n''$ and where the induction hypothesis applies to $z'$ and $z''$. Then

$$\delta(z) = \delta(z'z'') = z'^p\delta(z'') + z''^p\delta(z') + p\delta(z')\delta(z'') \in \mathrm{Fil}^{pn'+pn''}_\delta = \mathrm{Fil}^{pn}_\delta.$$

For statement (2), we again use the equation:

$$\delta(b+z) = \delta(b) + \delta(z) - C(b,z).$$

Since $C(b,z) \equiv b^{p-1}z \pmod{z^2}$, we see that $\delta(b+z) \equiv \delta(b) - b^{p-1}z \pmod{\mathrm{Fil}^2_\delta}$. □

# 2 Enlargements, envelopes, and dilatations

Our aim in this section is to describe and compare some of the notions of $p$-adic tubular neighborhoods that have appeared in the literature and to explain their relation to the new theory of prisms. In general, $X$ will be a scheme embedded as a closed subscheme of a $p$-torsion free scheme or formal scheme $Y$, and $\mathbb{T}_X(Y)$ will denote some sort of "tubular neighborhood" of $X$ in $Y$. For example, $\mathbb{T}_X(Y)$ could be one of the following:

$\mathbb{F}_X(Y)$, the formal completion of $Y$ along $X$,

$\mathbb{PD}_X(Y)$, the $p$-torsion free divided power envelope of $X$ in $Y$, which we call the "divided power enlargement" of $X$ in $Y$,

$\mathbb{D}_X(Y)$, the $p$-adic dilatation of $X$ in $Y$ (see §2.2 for the precise defintion and some variations).

$\mathbb{\Delta}_X(Y)$, the prismatic envelope of $X$ in a crystalline prism $Y$ (see Theorem 2.3.8 below for a precise statement of a generalization).

## 2.1 Divided power enlargements

We begin by discussing divided power enlargements, which underlie the first successful attempt at integral $p$-adic cohomology. Since we are restricting our attention here to torsion free algebras, the intricacies of divided powers are considerably simplified.

Let $B$ be a torsion free ring. For each $x \in B$, let $x^{[n]} := x^n/n! \in \mathbf{Q} \otimes B$, and recall the following formulas.

$$\begin{aligned}
(x+y)^{[n]} &= \sum_{i+j=n} x^{[i]} y^{[j]} \\
(bx)^{[n]} &= b^n x^{[n]} \\
x^{[m]} x^{[n]} &= \binom{m+n}{n} x^{[m+n]} \\
(x^{[n]})^{[m]} &= \prod_{i=1}^{m-1} \binom{in+n-1}{n-1} x^{[mn]}.
\end{aligned}$$

An ideal $I$ of $B$ is said to be a *PD-ideal* if $x^{[n]} \in I$ whenever $x \in I$ and $n > 0$. If $Y$ is a $p$-torsion free scheme or formal scheme, a closed immersion $X \to Y$ is a *PD-immersion* if its defining ideal is a PD-ideal. It follows from the above formulas that the set of all $x \in B$ such that $x^{[n]}$ belongs to $B$ for all $n \in \mathbf{N}$ forms an ideal $I_{\mathbb{PD}}(B)$ of $B$ and that $x^{[m]} \in I_{\mathbb{PD}}(B)$ if $x \in I_{\mathbb{PD}}(B)$. This ideal is the largest PD-ideal of $B$. It is the unit ideal if and only if $B$ contains $\mathbf{Q}$, and if $B$ is a $\mathbf{Z}_{(p)}$-algebra, the ideal $I_{\mathbb{PD}}(B)$ contains $p$. Let us record these facts in the following definition.

**Definition 2.1.1.** *If $B$ is a torsion free ring, let*

$$I_{\mathbb{PD}}(B) := \{x \in B : x^{[n]} \in B \text{ for all } n \in \mathbf{N}\},$$

*the largest PD-ideal of $B$. If $T$ is a $p$-torsion free formal $p$-adic scheme, define $I_{\mathbb{PD}} \subseteq \mathcal{O}_T$ in the same way, and note that $p \in I_{\mathbb{PD}}$. If $I_{\mathbb{PD}}$ is quasi-coherent, let $T_{\mathbb{PD}}$ denote the closed subscheme it defines.*

If $T$ is a $p$-torsion free formal $p$-adic scheme, the ideal $I_{\mathbb{PD}}$ of $\mathcal{O}_T$ is the inverse image of its image in $\mathcal{O}_{\overline{T}}$, since it contains $p$, and quasi-coherence of $I_{\mathbb{PD}}$ amounts to quasi-coherence of this image. This does not seem to be automatic, since it is the intersection of a countable sequence of quasi-coherent sheaves, which need not be quasi-coherent. Note that, if $\overline{T}_{\text{red}}$ is the reduced subscheme of $\overline{T}$, then $I_{\mathbb{PD}} \subseteq I_{\overline{T}_{\text{red}}}$. Indeed, if $x \in I_{\mathbb{PD}}$, then $x^{[p]} \in \mathcal{O}_T$, hence $x^p = p! x^{[p]} \in p\mathcal{O}_T$, hence $x$ belongs to the ideal $I_{\overline{T}_{\text{red}}}$.

**Definition 2.1.2.** *If $X$ is a scheme (typically in characteristic $p$), a PD-enlargement of $X$ is a $p$-torsion free $p$-adic formal scheme $T$ together with a PD-immersion $i_T \colon X_T \to T$ and a morphism $z_T \colon X_T \to X$. A PD-enlargement is* small *if $z_T$ is flat. A morphism of PD-enlargements of $X$ is a pair of morphisms $(f, f_X) \colon (T, X_T) \to (T', X_{T'})$ such that $z_{T'} \circ f_X = z_T$ and $f \circ i_T = i_{T'} \circ f_X$. If $i \colon X \to Y$ is a locally closed immersion of $X$ into a $p$-torsion free formal scheme, then a PD-enlargement of $X$ in $Y$ is PD-enlargement $(T, z_T, i_T)$ of $X$ together* *with a morphism $\pi_T \colon T \to Y$ such that $i \circ z_T = \pi_T \circ i_T$. The PD-dilatation of $X$ in $Y$ is the final object in the category of PD-enlargements of $X$ in $Y$.*

Our "PD-enlargements" differ from the "PD-thickenings" usually considered in the crystalline setting in two respects: we require the ambient formal scheme $T$ to be $p$-torsion free and we allow $z_T \colon X_T \to X$ to be any morphism, not necessarily an open immersion. The "smallness" condition will allow us to partially compensate for this extra generality.

The following result establishes the existence and basic properties of PD-dilatations; these may not be the same as PD-envelopes in general.

**Proposition 2.1.3.** *Let $Y$ be a $p$-torsion free $p$-adic formal scheme and let $X$ be a locally closed subscheme of $\overline{Y}$. Then the PD-dilatation of $X$ in $Y$ is representable by a diagram:*

$$\begin{array}{ccc} X_{\mathbb{PD}} & \xrightarrow{i_{\mathbb{PD}}} & \mathbb{PD}_X(Y) \\ {\scriptstyle z_{\mathbb{PD}}}\downarrow & & \downarrow{\scriptstyle \pi_{\mathbb{PD}}} \\ X & \longrightarrow & Y. \end{array}$$

*The morphisms $\pi_{\mathbb{PD}} \colon \mathbb{PD}_X(Y) \to Y$ and $z_{\mathbb{PD}} \colon X_{\mathbb{PD}} \to X$ are affine, and $\mathbb{PD}_X(Y)$ is a closed formal subscheme of the usual PD-envelope $\mathbb{DP}_X(Y)$ of $X$ in $Y$. Furthermore, the following statements are verified.*

*(1) Formation of $\mathbb{PD}_X(Y)$ is functorial: a morphism of pairs $g \colon (Y', X') \to (Y, X)$ induces a morphism $g \colon \mathbb{PD}_{X'}(Y') \to \mathbb{PD}_X(Y)$. If $g \colon Y' \to Y$ is $p$-completely flat and $X' = g^{-1}(X)$, then the corresponding map $\mathbb{PD}_{X'}(Y') \to \mathbb{PD}_X(Y) \times_Y Y'$ is an isomorphism. If $Y$ is noetherian and $Y'$ is the formal completion of $Y$ along $X$, then the map $\mathbb{PD}_X(Y') \to \mathbb{PD}_X(Y)$ is an isomorphism.*

*(2) Suppose that $S$ is a crystalline prism and that $Y \to S$ is a $p$-completely smooth morphism of $p$-adic formal schemes. If $X \subseteq \overline{Y}$ is reduced, the map $z_{\mathbb{PD}} \colon X_{\mathbb{PD}} \to X$ is an isomorphism, and if $X \to \overline{Y}$ is a regular immersion, then the maps $z_{\mathbb{PD}} \colon X_{\mathbb{PD}} \to X$ and $\mathbb{PD}_X(Y) \to \mathbb{DP}_X(Y)$ are isomorphisms.*

*(3) If $F_Y$ is a lifting of $F_{\overline{Y}}$, then $\mathbb{PD}_X(F_Y) \colon \mathbb{PD}_X(Y) \to \mathbb{PD}_X(Y)$ is a lifting of $F_{\overline{\mathbb{PD}}_X(Y)}$.*

*Proof.* Let $\mathbb{DP}_X(Y)$ be the ($p$-adically completed) PD-envelope described in [7, 3.19], and [4, I,4.2], *i.e.*, the direct limit of the PD-envelopes of $X$ in $Y_n$, for $n > 0$. This fits into a universal commutative diagram:

$$\begin{array}{ccc} X & \longrightarrow & \mathbb{DP}_X(Y) \\ & \searrow & \downarrow \\ & & Y, \end{array}$$

where $X \to \mathbb{DP}_X(Y)$ is a PD-immersion. This diagram might not define a PD-enlargement, since $\mathbb{DP}_X(Y)$ might not be torsion free. To remedy this, let $I_{\mathrm{tor}}$ denote the sheaf of $p$-torsion sections of the structure sheaf $\mathcal{D}_X(Y)$ of $\mathbb{DP}_X(Y)$ and let $\overline{I}_{\mathrm{tor}}$ denote its closure in the $p$-adic topology. As we saw in the proof of Proposition 1.6.11, this ideal defines a closed immersion

$$\mathbb{PD}_X(Y) := \mathbb{DP}_{\mathrm{tf}} \to \mathbb{DP}_X(Y),$$

and $\mathbb{PD}_X(Y)$ is $p$-torsion free. Furthermore, the intersection of $\overline{I}_{\mathrm{tor}}$ with the ideal $I_X$ of $X$ in $\mathcal{D}_X(Y)$ is closed under divided powers, so that $\overline{I}_X/(\overline{I}_X \cap \overline{I}_{\mathrm{tor}})$ is a divided power ideal of $\mathcal{O}_{\mathbb{PD}_X(Y)}$ and hence defines a PD-immersion $X_{\mathbb{PD}} \to \mathbb{PD}_X(Y)$. It follows easily that this construction gives the universal PD-enlargement of $X$. The functoriality of the formation of $\mathbb{PD}_X(Y)$ follows from its universal property, and the affineness of $\pi_{\mathbb{PD}}$ and $z_{\mathbb{PD}}$ is clear from the construction.

If $g\colon Y' \to Y$ is $p$-completely flat, then each $g_n\colon Y'_n \to Y_n$ is flat (see Proposition 7.3.2), and if in addition $X' = g^{-1}(X)$, it follows that each map $\mathbb{DP}_{X'}(Y'_n) \to \mathbb{DP}_X(Y_n) \times_{Y_n} Y'_n$ is an isomorphism [4, 2.7.1]. Then $\mathbb{DP}_{X'}(Y') \to \mathbb{DP}_X(Y) \times_Y Y'$ is also an isomorphism, and it follows that $\mathbb{DP}_{X'}(Y') \to \mathbb{DP}_X(Y)$ is $p$-completely flat, by Proposition 7.3.3. Then $\overline{I}_{\mathrm{tor}}(X') \subset \mathcal{D}_{X'}(Y')$ is the pullback of $\overline{I}_{\mathrm{tor}}(X) \subseteq \mathcal{D}_X(Y)$, by Corollary 7.3.7. It follows that $\mathbb{PD}_{X'}(Y') \to \mathbb{PD}_X(Y) \times_Y Y'$ is also an isomorphism.

If $Y$ is noetherian and $Y'$ is the formal completion of $Y$ along $X$, then $\mathbb{DP}_X(Y') \cong \mathbb{DP}_X(Y)$ by [7, 3.34], and so $\mathbb{PD}_X(Y') \cong \mathbb{PD}_X(Y)$. This completes the proof of statement (1).

Statement (2) can be verified locally, so we may assume that $X$ and $Y$ are affine. Recall from [21, 0.§1.4] that the embedding of $X$ in $\mathrm{W}(X)$ is a PD-immersion. If $X$ is reduced, then the operator $F$ of $W\mathcal{O}_X$ is injective, and since the Verschiebung $V$ is also injective and the composition $FV$ is multiplication by $p$, the ring $W\mathcal{O}_X$ is $p$-torsion free. Since $S$ is a $\delta$-scheme, the map $X \to S$ lifts to a map $\mathrm{W}(X) \to S$. Since each $X \to \mathrm{W}_n(X)$ is nilpotent, since $Y/S$ is formally smooth, and since $\mathrm{W}(X)$ is affine, we can extend the formal embedding $X \to Y$ to a map $\mathrm{W}(X) \to Y$. Then the universal property of $\mathbb{PD}_X(Y)$ produces maps of pairs:

$$(X \subseteq \mathrm{W}(X)) \longrightarrow (X_{\mathbb{PD}} \subseteq \mathbb{PD}_X(Y)) \longrightarrow (X \subseteq Y).$$

Since the map $X_{\mathbb{PD}} \to X$ is, by construction, a closed immersion, it must in fact be an isomorphism. If $X \to \overline{Y}$ is a regular immersion, then $\mathbb{DP}_X(Y)$ is already $p$-torsion free. (To see this, one can use the compatibility of formation of PD-envelopes with flat base extension to reduce to the case in which $X \to \overline{Y}$ is the zero section of $\mathbf{A}^n_X$.) Thus $\mathbb{PD}_X(Y) \to \mathbb{DP}_X(Y)$ and $X_{\mathbb{PD}} \to X$ are isomorphisms in this case.

It suffices to check statement (3) in the affine case, and we use affine notation, with $Y = \operatorname{Spf} B$ and $X \subseteq Y$ defined by the ideal $I$. As we have seen, then $\mathbb{PD}_X(Y)$ is the $p$-adic completion $C$ of the spectrum of the $B$-algebra $D_I(B)/\overline{I}_{\mathrm{tor}}$ and hence is topologically generated by elements of the form $x^{[n]}$ for all $x \in I$.

Let us write $\phi$ for $F_Y^*$. Our claim is that $\phi(c) \equiv c^p \pmod{pC}$ for all $c \in C$. Since this is true for elements of $B$, it will suffice to check that it holds on a set of topological generators for the $B$-algebra $C$, and hence for the elements $x^{[n]}$ with $x \in I$ and $n > 0$. Since each of these belongs to a divided power ideal, $(x^{[n]})^p \in pC$. On the other hand, $\phi(x^{[n]}) = \phi(x)^{[n]} \equiv (x^p)^{[n]} \pmod{pC}$ and $(x^p)^{[n]} = (p!x^{[p]})^{[n]} = (p!)^n(x^{[p]})^{[n]} \in pC$. Thus $(x^{[n]})^p$ and $\phi(x^{[n]})$ are both congruent to zero mod $pC$. □

The following result is a reformulation of [10, 2.35].

**Proposition 2.1.4.** *Let $Y$ be a crystalline prism and let $F(\overline{Y}) \subset \overline{Y} \subseteq Y$ denote the scheme-theoretic image of $F_{\overline{Y}}$, i.e., the closed subscheme defined by $I_\phi := \{c \in \mathcal{O}_Y : c^p \in p\mathcal{O}_Y\}$. Then $F(\overline{Y}) \to Y$ is a PD-immersion. In fact, $F(\overline{Y}) = Y_{\mathbb{PD}}$, the smallest PD-subscheme of $Y$ (see Definition 2.1.1).*

*Proof.* We claim that the ideal $I_\phi$ defining the closed immersion $F_Y(\overline{Y}) \to Y$ is a PD-ideal. Assume without loss of generality that $Y = \operatorname{Spf} C$.

We check first that if $c \in I_\phi$ and $c^p = pc_1$, then also $c_1 \in I_\phi$. in fact:

$$\begin{array}{rcl} \phi(c)^p & = & p\phi(c_1) \\ (c^p + p\delta(c))^p & = & p(c_1^p + p\delta(c_1)) \\ c^{p^2} + p^2(\cdots) & = & pc_1^p + p^2\delta(c_1) \\ p^p c_1^p + p^2(\cdots) & = & pc_1^p + p^2\delta(c_1) \\ (p^{p-1} - 1)c_1^p & = & p(\delta(c_1) - \cdots). \end{array}$$

This last equation implies that $c_1^p \in pC$, as claimed.

Continuing by induction, we find a sequence $c_0, c_1, c_2 \ldots$ such that $c_0 = c$ and $c_i^p = pc_{i+1}$ for $i \geq 0$. It follows that $c^{p^i} = p^{1+p+\cdots+p^{i-1}}c_i$ for all $i$, and if $a_i \in \mathbf{N}$, that $c^{a_i p^i} \in p^{a_i(1+p+\cdots+p^{i-1})}C$. If $n$ is a natural number, let $n = \sum a_i p^i$ be its $p$-adic expansion, and let $\sigma := \sum a_i$. Then $c^n \in p^m C$, where

$$m := \sum_i a_i(1 + p + \cdots + p^{i-1}) = \sum_i a_i \frac{p^i - 1}{p - 1} = \frac{n - \sigma}{p - 1}.$$

The expression on the right is the $p$-adic ordinal of $n!$ [7, 3.3], so we conclude that $c^{[n]} := c^n/n! \in C$. Moreover, since $c \in I_\phi$, we can write $\phi(c) = pb$, and then $\phi(c^{[n]}) = \phi(c)^{[n]} = (pb)^{[n]} = b^n p^{[n]} \in pC$, so $c^{[n]}$ again belongs to $I_\phi$. This shows that $I_\phi$ is indeed a PD-ideal as claimed.

Finally, suppose that $X \subseteq \overline{Y}$ is a closed subscheme whose ideal $I \subseteq C$ is invariant under divided powers. For each $c \in I$, we have $c^p = p!c^{[p]}$, so $c^p \in pC$, hence $c \in I_\phi$. □

## 2.2 $p$-adic enlargements, dilatations, and dilations

The following special PD-enlargements were first explicitly discussed in [31], as the starting point for the theory of convergent and rigid cohomology. They were further developed in [40], in the context of a $p$-adic approach to the Cartier transform.

Recall that if $T$ is a scheme or formal scheme, we write either $T_1$ or $\overline{T}$ for the closed subscheme defined by $(p)$.

**Definition 2.2.1.** *Let $X$ be a $p$-adic formal scheme, for example an $\mathbf{F}_p$-scheme. A $p$-adic enlargement of $X$ is a $p$-torsion free $p$-adic formal scheme $T$ together with a morphism $z_T \colon \overline{T} \to X$. A $p$-adic enlargement is small if $z_T$ is flat. A morphism of $p$-adic enlargements of $X$ is a morphism $f \colon T \to T'$ such that $z_{T'} \circ \overline{f} = z_T$. If $i \colon X \to Y$ is a closed immersion of $X$ into a scheme or formal scheme $Y$, then a $p$-adic enlargement of $X$ in (or over) $Y$ is a $p$-adic enlargement $(T, z_T)$ of $X$ together with a map $\pi_T \colon T \to Y$ such that $i \circ z_T = \pi_{T|_{\overline{T}}}$:*

$$\begin{array}{ccc} \overline{T} & \longrightarrow & T \\ {\scriptstyle z_T}\downarrow & & \downarrow{\scriptstyle \pi_T} \\ X & \xrightarrow{\;i\;} & Y. \end{array}$$

As an alternate phrasing, we may say *$p$-adic $X$-enlargement* instead of "$p$-adic enlargement of $X$," and we may say a *$p$-adic $X$-enlargement over $Y$* instead of a "$p$-adic enlargement of $X$ in $Y$." Note that the categories of $p$-adic enlargements of $X$ and of $\overline{X}$ are identical. A $p$-adic enlargement of $X$ is also a PD-enlargement, since the ideal $(p)$ always has divided powers.

Suppose that $i \colon X \to Y$ is an embedding of $X$ in a scheme or formal scheme $Y$ and that $(T, z_T, \pi_T)$ is a $p$-adic enlargement of $X$ in $Y$. The map $\overline{T} \to T$ necessarily factors through $\overline{X} \times_Y T$, and the map $\overline{X} \times_Y T \to T$ necessarily factors through $\overline{T}$. It follows that these factorizations are isomorphisms, so the diagram in Definition 2.2.1 is Cartesian when $X = \overline{X}$. Thus the image of the map $\pi_T^\sharp \colon \pi_T^*(I_{\overline{X}/Y}) \to \mathcal{O}_T$ is the principal ideal $p\mathcal{O}_T$. Consequently there is a unique $\mathcal{O}_Y$-linear map:

$$\rho_T \colon I_{\overline{X}/Y} \to \pi_{T*}(\mathcal{O}_T) \tag{2.2.1}$$

such that $p\rho_T(a) = \pi_T^\sharp(a)$ for all $a \in I_{\overline{X}/Y}$. Since $\mathcal{O}_T$ is $p$-torsion free, this implies (and is equivalent to) the equality $\rho_T(p) = 1$.

**Definition 2.2.2.** *If $i \colon X \to Y$ is a locally closed immersion of a $p$-adic formal scheme $X$ into a scheme or a formal scheme $Y$, then the $p$-adic dilatation of $X$ in $Y$, denoted by $\mathbb{D}_X(Y)$, is the final object of the category of $p$-adic enlargements of $X$ in $Y$:*

$$\begin{array}{ccc} \overline{\mathbb{D}}_X(Y) & \longrightarrow & \mathbb{D}_X(Y) \\ {\scriptstyle z_{\mathbb{D}}}\downarrow & & \downarrow{\scriptstyle \pi_{\mathbb{D}}} \\ X & \xrightarrow{\;i\;} & Y. \end{array}$$

Note that since the categories of $p$-adic enlargements of $X$ and of $\overline{X}$ in $Y$ are the same, the map $\mathbb{D}_{\overline{X}}(Y) \to \mathbb{D}_X(Y)$ is an isomorphism, which we regard as an identification. Note also that, if $\hat{Y}$ is the formal completion of $Y$ along $X$ and $T \to Y$ is a $p$-adic enlargement of $X$ in $Y$, then, since the ideal $I_{X/Y}$ of $X$ in $Y$ maps to $p\mathcal{O}_T$, and $T$ is $p$-adically complete, the map $T \to Y$ factors uniquely through $\hat{Y}$. Thus $\mathbb{D}_X(\hat{Y}) \to \mathbb{D}_X(Y)$ is an isomorphism.

**Warning 2.2.3.** The $p$-adic dilatation of $X$ in $Y$ could well be empty. For example, let $V := W[x]/(x^m - p)$ with $m > 1$, let $Y := \operatorname{Spf} V$ and $X := \operatorname{Spec} k$. If $V \to B$ defines a $p$-adic enlargement of $X$ in $Y$, then $xB = pB$, hence $B = x^{m-1}B$, and since $B$ is $p$-adically complete, it follows that $B = 0$. This difficulty was addressed in [31] by considering more general enlargements of $X$, namely, $p$-adic formal schemes $T$ endowed with a map from the reduced subscheme of $\overline{T}$ to $X$. The cohomology of the corresponding topos corresponds to the so-called "convergent cohomology" of $X$ [32].

As explained in [31] and [40], the final object $\mathbb{D}_X(Y)$ is representable by a $p$-adic formal scheme. Here we review the construction and discuss some additional approaches that we shall find useful later. We begin by describing a related construction, which we describe below in the context of $p$-adic formal schemes. Of course, there are analogus constructions in the category of schemes.

**Definition 2.2.4.** *Let $Y$ be a $p$-adic formal scheme and $E$ a quasi-coherent sheaf of $\mathcal{O}_Y$-modules. If $\pi_T \colon T \to Y$ is a morphism of $p$-adic formal schemes and $e$ is a global section of $E$, let*

$$
\begin{aligned}
\hat{\mathbb{V}}E(T) &:= \{v \colon E \to \pi_{T*}(\mathcal{O}_T)\} \cong \{v \colon \pi_T^*(E) \to \mathcal{O}_T\}\\
\hat{\mathbb{V}}_{|e}E(T) &:= \{v \in \hat{\mathbb{V}}E(T) : v(e) = 0\}\\
\hat{\mathbb{V}}_{:e}E(T) &:= \{\rho \in \hat{\mathbb{V}}E(T) : \rho(e) = 1\}.
\end{aligned}
$$

*The functor $\hat{\mathbb{V}}_{:e}E$ is called the $e$-dilation of $E$. If $X \subseteq \overline{Y}$ is a closed subscheme, then $p$ belongs to $I_{X/Y}$, and $\hat{\mathbb{V}}_{:p}(I_{X/Y})$ is the dilation of $X$ in $Y$.*

The proof of the following proposition is immediate and omitted.

**Proposition 2.2.5.** *Let the notation be as in Definition 2.2.4.*

*(1) The functor $\hat{\mathbb{V}}E$ is representable by the pair $(s_E, v_E)$, where*

$$s_E \colon \operatorname{Spf}_Y(\hat{S}^\bullet E) \to Y$$

*is the formal $Y$-scheme given by the $p$-adic completion of the symmetric algebra of $E$ and*

$$v_E \colon E \to \hat{S}^\bullet E \in \hat{\mathbb{V}}E(\hat{S}^\bullet E)$$

*is the inclusion.*

*(2) The functors $\hat{\mathbb{V}}_{|e}E$ and $\hat{\mathbb{V}}_{:e}E$ are representable by the closed subschemes of $\hat{\mathbb{V}}E$ defined by the ideals of $\hat{S}^\bullet E$ generated by $v_E(e)$ and by $v_E(e) - 1$, respectively. Their formation is compatible with base change. Thus, if $f \colon Y' \to Y$ is a morphism and $E' := f^*(E), e' := f^*(e)$, then*

$$\hat{\mathbb{V}}_{|e'}E' \cong \hat{\mathbb{V}}_{|e}E \times_Y Y' \quad \text{and} \quad \hat{\mathbb{V}}_{:e'}E' \cong \hat{\mathbb{V}}_{:e}E \times_Y Y'.$$

(3) *$\hat{\mathbb{V}}_{|e}E \cong \hat{\mathbb{V}}E_{|e}$, where $E_{|e} := E/(\mathcal{O}_Y e)$, and*

$$S^\bullet E/(v_E(e)-1) \cong \varinjlim\{S^n E \xrightarrow{\cdot e} S^{n+1}E : n \in \mathbf{N}\}.$$

(4) *$\hat{\mathbb{V}}E$ has a natural structure of a (formal) vector space scheme, $\hat{\mathbb{V}}_{|e}E \subseteq \hat{\mathbb{V}}E$ is a subvector space scheme, and $\hat{\mathbb{V}}_{:e}E$ has a natural structure of a pseudo-torsor under $\hat{\mathbb{V}}_{|e}E$, given by :*

$$\hat{\mathbb{V}}_{:e}E \times_Y \hat{\mathbb{V}}_{|e}E \to \hat{\mathbb{V}}_{:e}E : (\rho, v) \mapsto v + \rho. \quad \square$$

Now suppose that $E$ is an ideal of $\mathcal{O}_Y$, or, more generally a quasi-ideal, *i.e.*, a homomorphism of $\mathcal{O}_Y$-modules $j\colon E \to \mathcal{O}_Y$ such that for all $x, y \in E$, $j(x)y = j(y)x \in E$ [15, 3.1]. Note that if $\pi_T\colon T \to Y$ is a morphism, then

$$j_T := \pi_T^*(j)\colon \pi_T^*(E) \to \mathcal{O}_T$$

is again a quasi-ideal. If $S$ is an element or a subset of $E$, we let $T_S$ be the closed formal subscheme of $T$ defined by the ideal generated by $\pi_T^*(j(S))$. We shall say that an element $e$ of $E$ is *$j$-regular* if the map $\mathcal{O}_Y \to E$ sending $a$ to $ae$ is injective, and that a sequence $(e_1, \ldots, e_r)$ is $j$-regular if for every $i$, the image of $e_i$ in $\mathcal{O}_{Y_{\{e_1,\ldots,e_{i-1}\}}} \otimes E$ is $j_{Y_{\{e_1,\ldots,e_{i-1}\}}}$-regular.

If $\rho \in \hat{\mathbb{V}}_{:e}E(T)$, then elements of $E$ become explicitly divided by $e$ in $\mathcal{O}_T$, as we see in the following proposition.

**Proposition 2.2.6.** *Suppose that $j\colon E \to \mathcal{O}_Y$ is a quasi-ideal.*

(1) *If $\pi_T\colon T \to Y$ is a morphism and $\rho \in \mathbb{V}E_{:e}(T)$, then $T_E = T_e$, and there is a commutative diagram:*

$$\begin{array}{ccc} & & \mathcal{O}_T \\ & \overset{\pi_T^*(j(e))\rho}{\nearrow} & \uparrow \pi_T^* \\ E & \xrightarrow{\ j\ } & \mathcal{O}_Y. \end{array}$$

*If $b \in \mathcal{O}_Y$ and $be = 0$, then $\pi_T^*(b) = 0$.*

(2) *If $\pi_T\colon T \to Y$ is a morphism, then the set $\hat{\mathbb{V}}_{:e}(T)$ is naturally a pseudo-torsor under $\hat{\mathbb{V}}_{|e}(T)$, and every $v \in \hat{\mathbb{V}}_{|e}(T)$ factors through the $\pi_T^*(j(e))$-torsion subsheaf of $\mathcal{O}_T$.*

(3) *There is a Cartesian diagram:*

$$\begin{array}{ccc} Y_e \times_Y \hat{\mathbb{V}}E_{:e} & \longrightarrow & \hat{\mathbb{V}}_{:e}E \\ \downarrow & & \downarrow \pi_{\mathbb{V}_{:e}} \\ Y_E & \longrightarrow & Y, \end{array}$$

*and the map $\pi^*_{V:e}$ annihilates the annihilator of $e$.*

(4) *There is a commutative diagram with exact rows*

$$\begin{array}{ccccccc}
\mathcal{O}_{Y_E} & \longrightarrow & E \otimes \mathcal{O}_{Y_e} & \longrightarrow & E_{|e} & \longrightarrow & 0 \\
\downarrow & & \downarrow & & \downarrow & & \\
\mathcal{O}_{Y_E} & \longrightarrow & E \otimes \mathcal{O}_{Y_E} & \longrightarrow & E_{|e} \otimes \mathcal{O}_{Y_E} & \longrightarrow & 0,
\end{array}$$

*where the left horizontal maps are multiplication by the class of $e$, and $E_{|e} := E/\mathcal{O}_Y e$. The first of these is injective if $e$ is $j$-regular, and they are both injective if $E$ is generated by a $j$-regular sequence $(e_1, \ldots, e_r)$ with $e_r = e$.*

*Proof.* Say $e' \in E$. Since $j$ is a quasi-ideal and $\rho$ is $\pi_T^*$-linear, we have:

$$\begin{array}{rcl}
j(e')e & = & j(e)e' \quad \in E \\
\rho(j(e')e) & = & \rho(j(e)e') \quad \in \mathcal{O}_T \\
\pi_T^*(j(e'))\rho(e) & = & \pi_T^*(j(e))\rho(e') \in \mathcal{O}_T \\
\pi_T^*(j(e')) & = & \pi_T^*(j(e))\rho(e') \in \mathcal{O}_T.
\end{array}$$

This establishes the commutative diagram in (1). It follows that $\pi_T^*(j^*(E))\mathcal{O}_T = \pi_T^*(e)\mathcal{O}_T$, and hence that $T_E = T_e$. If $b \in \mathcal{O}_Y$ and $be = 0$, then

$$\pi_T^*(b) = \pi_T^*(b)\rho(e) = \rho(be) = 0.$$

This proves statement (1), and statement (3) is the special case when $T = \hat{\mathbb{V}}E$.

The pseudo-torsor structure of $\mathbb{V}_{:e}(T)$ is described in (4) of Proposition 2.2.5. If $e' \in E$ and $v \in \mathbb{V}_{|e}(T)$, then

$$\pi_T^*(j(e))v(e') = v(j(e)e') = v(j(e')e) = \pi_T^*(j(e'))v(e) = 0.$$

This proves (2).

To prove (4), note that first that $j(e)E \subseteq \mathcal{O}_Y e \subseteq E$ since $j(e)e' = j(e')e$ if $e' \in E$. Thus we have an exact sequence

$$\mathcal{O}_Y \to E/j(e)E \to E_{|e} \to 0,$$

where the left map sends $a \in \mathcal{O}_Y$ to the class of $ae$ in $E/j(e)E$. But if $a = j(e')$ with $e' \in E$, then $ae = j(e')e = j(e)e'$ vanishes in $E/j(e)E$, so in fact this map factors through $\mathcal{O}_{Y_E}$. Now suppose $e$ is $j$-regular. Then if $a \in \mathcal{O}_Y$ and $ae = j(e)e'$ for some $e'$ in $E$, we have $ae = j(e')e$, so $a = j(e')$, which maps to zero in $\mathcal{O}_{Y_E}$. Thus the top left arrow in the diagram is injective. Finally, suppose that $(e_1, \ldots, e_r)$ is a $j$-regular sequence which generates $E$ with $e_r = e$. Let $Y' := Y_{\{e_1, \ldots, e_{r-1}\}}$ and let $j' := j \otimes \mathrm{id}_{\mathcal{O}_{Y'}} \colon E' \to \mathcal{O}_{Y'}$. Then the image $e'$ of

$e$ in $E'$ is $j'$-regular. Applying what we have just seen to the quasi-ideal $j'$, we see that the sequence

$$0 \longrightarrow \mathcal{O}_{Y'_{E'}} \longrightarrow E' \otimes \mathcal{O}_{Y'_{e'}} \longrightarrow E'_{|e'} \longrightarrow 0$$

is exact. This is the bottom row of the diagram in (3). $\square$

**Corollary 2.2.7.** *Let $Y$ be a $p$-torsion free $p$-adic formal scheme, let $X \subseteq \overline{Y}$ be a closed subscheme, and let*

$$\mathbb{V}_X(Y) := \hat{\mathbb{V}}_{:p}(I_{X/Y}).$$

(1) *The reduction modulo $p$ of the structure morphism $\pi_{\mathbb{V}} \colon \mathbb{V}_X(Y) \to Y$ factors through $X$.*

(2) *There is a natural bijection between the set of sections $X \to \mathbb{V}_X(Y)$ of $\pi_{\mathbb{V}}$ and the set of flat liftings of $X$ in $Y_2$. This bijection is compatible with the respective actions of $\mathbb{N}_{\overline{X}/\overline{Y}}$ on each of these sets (see Proposition 1.7.1).*

(3) *More generally, if $X \to \overline{Y}$ is locally closed, let*

$$\mathbb{V}_X(Y) := \mathbb{V}_{:p}(I_{X/Y'}),$$

*where $Y'$ is any open subset of $Y$ containing $X$ as a closed subset, the result being independent of the choice of $Y'$.*

*Proof.* Statement (1) follows from statement (3) of Proposition 2.2.6. For statement (2), note that a section $X \to \mathbb{V}_X(Y)$ is equivalent to a linear map $\rho \colon I_{X/Y} \to \mathcal{O}_X$ sending $p$ to 1, or to a linear map $\mathcal{C}_{X/Y} \to \mathcal{O}_X$ sending the class of $p$ to 1. Proposition 1.7.1 shows that this is equivalent to giving a flat lifting of $X$ to $Y_2$ and that this equivalence is compatible with the actions of $\mathbb{N}_{X/\overline{Y}}$. If $Y' \subseteq Y''$ are open subsets of $Y$ in which $X$ is closed, then it follows from (1) that the natural map $\mathbb{V}_X(Y') \to \mathbb{V}_X(Y'')$ is an isomorphism; this justifies statement (3). $\square$

We now discuss dilatations and relate them to dilations.

**Theorem 2.2.8.** *Let $Y$ be a scheme or formal $p$-adic scheme and let $i \colon X \to Y$ be a locally closed immersion, where $X$ is a $p$-adic formal scheme. Then $\mathbb{D}_X(Y)$ is representable, and the morphisms $\pi_{\mathbb{D}} \colon \mathbb{D}_X(Y) \to Y$ and $z_{\mathbb{D}} \colon \overline{\mathbb{D}}_X(Y) \to X$ are affine. Furthermore, the following statements are verified.*

(1) *$\mathbb{D}_X(Y) = \mathbb{V}_X(V)_{\mathrm{tf}}$, where $V \subseteq Y$ is any open set containing $X$ in which $X$ is closed. (Recall that for any $p$-adic formal scheme $T$, $T_{\mathrm{tf}} \subseteq T$ is the closed formal subscheme defined by the closure of the ideal of $p$-torsion elements of $\mathcal{O}_T$.)*

(2) *$\mathbb{D}_X(Y)$ is the open subset $D_+(p)$ of the formal blowup of $\overline{X}$ in $Y$ defined by the element $p$ of the ideal of $\overline{X}$ in $Y$.*

(3) *Formation of $\mathbb{D}_X(Y)$ is functorial: a morphism $g\colon Y' \to Y$ sending a closed subscheme $X'$ of $Y'$ to $X$ induces a morphism $\mathbb{D}_X(g)\colon \mathbb{D}_{X'}(Y') \to \mathbb{D}_X(Y)$. If $g$ is $p$-completely flat and $X' = g^{-1}(X)$, then the natural map*

$$\mathbb{D}_{X'}(Y') \to \mathbb{D}_X(Y) \times_Y Y'$$

*is an isomorphism.*

(4) *Assume that $X = \overline{X}$ and that $\bar{i}\colon X \to \overline{Y}$ is a regular immersion. Then the following hold:*

(a) *There is an exact sequence of locally free $\mathcal{O}_X$-modules:*

$$0 \to \mathcal{O}_{\overline{X}} \xrightarrow{[p]} \mathcal{C}_{\overline{X}/Y} \longrightarrow \mathcal{C}_{\overline{X}/\overline{Y}} \to 0.$$

(b) *$\mathbb{D}_X(Y) = \mathbb{V}_X(Y)$, and $\overline{\mathbb{D}}_X(Y) = \operatorname{Spec}(\varinjlim S^n(\mathcal{C}_{\overline{X}/Y}))$, where the maps*

$$S^n(\mathcal{C}_{\overline{X}/Y}) \to S^{n+1}(\mathcal{C}_{\overline{X}/Y})$$

*are given by multiplication by $[p]$. If $\pi_T\colon T \to Y$ is a morphism, then the set of liftings $T \to \mathbb{D}_X(Y)$ of $\pi_T$ is naturally a pseudo-torsor under $\operatorname{Hom}_{\mathcal{O}_Y}(i_*(\mathcal{C}_{X/\overline{Y}}), \pi_{T*}(\mathcal{O}_T))$.*

(c) *The map $z_Y\colon \overline{\mathbb{D}}_X(Y) \to X$ is faithfully flat and is naturally a torsor under the action of the normal bundle $\mathbb{N}_{X/\overline{Y}}$ of $X$ in $\overline{Y}$. There is a natural isomorphism of torsors between the set of sections of $\overline{\mathbb{D}}_X(Y) \to X$ and the set of flat liftings of $X$ in $Y_2$.*

*Proof.* Since $p$-adic enlargements of $X$ and of $\overline{X}$ are the same, we may and shall replace $X$ by $\overline{X}$. Moreover, every $p$-adic enlargement of $X$ factors through $Y_{\mathrm{tf}}$, so we may also assume that $Y$ is torsion free. If $V$ is open in $Y$ and contains $X$, then the universal property of dilatations implies that the map $\mathbb{D}_X(Y) \to Y$ factors through $V$, and it follows that the map $\mathbb{D}_X(V) \to \mathbb{D}_X(Y)$ is an isomorphism.

Statement (1) of Corollary 2.2.7 implies that the map $\overline{\mathbb{V}}_X(Y) \to Y$ factors through $\overline{X}$, and hence the same is true of the map $\overline{\mathbb{V}}_X(Y)_{\mathrm{tf}} \to Y$. Thus $\mathbb{V}_X(Y)_{\mathrm{tf}} \to Y$ defines a $p$-adic enlargement of $X$ in $Y$. Suppose that $(T, \pi_T, z_T)$ is any such. Then the map $\rho_T$ (2.2.1) sends $p$ to 1 and hence the map $\pi_T$ factors through $\mathbb{V}_X(Y)$. Since $T$ is $p$-torsion free by assumption, it also factors through $\mathbb{V}_X(Y)_{\mathrm{tf}}$. This proves (1).

For the rest of the proof, let us assume that $Y$ is affine, say $\operatorname{Spf} B$ (or $\operatorname{Spec} B$), and that $I$ is the defining ideal of $X$ in $Y$. Let $B_I$ be the Rees-algebra $B \oplus I \oplus I^2 \oplus \cdots$ and let $[p]$ be the element of degree one defined by $p$. Statement (2) asserts that $\mathbb{D}_X(Y)$ is the formal spectrum of the $p$-adic completion of the degree zero part of the localization of $B_I$ by $[p]$. This is proved in [31, 2.3] and [40, 3.5]. Alternatively, one can take the completion of the quotient of $B_I$ by the ideal generated by $[p] - 1$, or of the direct limit of the system $B \xrightarrow{[p]} I \xrightarrow{[p]} I^2 \xrightarrow{[p]} \cdots$, with the natural $B$-algebra structure.

In the situation of statement (3), let $f\colon X' \to X$ be the morphism induced by $g$. Then $(\mathbb{D}_{X'}(Y'), g \circ \pi_{\mathbb{D}'}, f \circ z_{\mathbb{D}'})$ is a $p$-adic enlargement of $X$ over $Y$, so we find the desired morphism $\mathbb{D}_{X'}(Y') \to \mathbb{D}_X(Y)$ by the universal property of $\mathbb{D}_X(Y)$. If $X' = g^{-1}(X)$, then the projection map

$$(\mathbb{D}_X(Y) \times_Y Y')_1 = \mathbb{D}_X(Y)_1 \times_{Y_1} Y'_1 \to Y'_1$$

factors through a map $z_{Y'}\colon (\mathbb{D}_X(Y) \times_Y Y')_1 \to X'$, and if $g$ is $p$-completely flat, then $\mathbb{D}_X(Y) \times_Y Y'$ is again $p$-torsion free, by statement (1) of Proposition 7.3.3. Then $(\mathbb{D}_X(Y) \times_Y Y', \pi_{Y'}, z_{Y'})$ defines a $p$-adic enlargement of $X'$ over $Y'$, and so there is a map $\mathbb{D}_X(Y) \times_Y Y' \to \mathbb{D}_{X'}(Y')$, inverse to the map $\mathbb{D}_{X'}(Y') \to \mathbb{D}_X(Y) \times_Y Y'$ coming from functoriality.

Now suppose that $X \to \overline{Y}$ is a regular immersion. Then the map $X \to Y$ is also a regular immersion. If $(p, x_1, \ldots, x_r)$ is a regular sequence generating the ideal $I_{X/Y}$, its image forms a basis for $I_{X/Y}/I^2_{X/Y}$, so (4a) is clear. Furthermore, Proposition 7.5.3 implies that the natural map from the symmetric algebra of $I$ to the Rees algebra is an isomorphism, so (4b) follows from statement (2) and Propositions 2.2.5 and 2.2.6. If $\pi_T\colon T \to X$ is a morphism, then maps $i_*(\mathcal{C}_{X/\overline{Y}}) \to \pi_{T*}(\mathcal{O}_T)$ identify with $T$-valued points of $\mathbb{N}_{X/\overline{Y}}$, and in particular $\overline{\mathbb{D}}_X(Y)$ becomes a torsor under $\mathbb{N}_{X/\overline{Y}}$, hence $\overline{\mathbb{D}}_X(Y) \to X$ is faithfully flat. The bijection between sections of $\overline{\mathbb{D}}_X(Y)$ and flat liftings follows from Corollary 2.2.7. $\square$

**Remark 2.2.9.** For yet another local construction of $\mathbb{D}_X(Y)$, suppose that $(p, x_1, \ldots, x_r)$ is a set of generators for $I$ and let $B' := B[t_1, \ldots, t_r]/(pt_1 - x_1, \ldots, pt_r - x_r)$, modulo $p$-torsion. The $p$-adic completion $B'^\wedge$ of $B'$ is again $p$-torsion free (see Lemma 7.3.4), and $IB'^\wedge$ is the principal ideal $p$. Thus $B \to B'^\wedge$ defines a $p$-adic enlargement of $X$ over $Y$, and it is clearly universal. This construction looks more down-to-earth, but computing the $p$-torsion may be difficult. See Proposition 2.2.11 below for an important case in which there is none.

**Example 2.2.10** Let $A$ be a $p$-torsion free and $p$-adically complete ring, let $B$ be the $p$-adic completion of the polynomial algebra $A[x_1, \ldots, x_r]$, and let $I$ be the ideal of $B$ generated by $(p, x_1, \ldots, x_r)$. Then the $p$-adic dilatation of $\operatorname{Spf} B/I$ in $\operatorname{Spf} B$ is given by the formal spectrum of the $p$-adic completion of

$$B' := B[t_1, \ldots, t_r]/(x_1 - pt_1, \ldots, x_r - pt_r)$$

Indeed, it is clear that this ring is already $p$-torsion free, that $IB' = pB'$, and that it is universal with this property. The following result generalizes this fact.

**Proposition 2.2.11.** *Suppose that $B$ is a $p$-torsion free $p$-adically complete ring and that $I$ is the ideal of $B$ generated by a finite regular sequence $(p, x_1, \ldots, x_r)$. Then the $p$-adic dilatation of $B/I$ in $B$ is given by the homomorphism $B \to B'$ where $B'$ is the $p$-adic completion of $B[t_1, \ldots, t_r]/(pt_1 - x_1, \ldots, pt_r - x_r)$.*

*Proof.* If $B$ is noetherian, this follows from the Example 2.2.10, because the map from the ($p$-adically completed) polynomial algebra $\mathbf{Z}_p[X_1, \ldots, X_r]\hat{}$ to $B$ sending $X_i$ to $x_i$ is flat at points of the zero set of $I$, and statement (2) of Theorem 2.2.8 tells us that formation of $p$-adic dilatations commutes with flat base change. To avoid the noetherian hypothesis, we can argue as follows. Thanks to Lemma 2.2.13 below, we can use induction on $r$ to reduce to the case in which $I$ is generated by a regular sequence $(p, x)$. Our claim is simply that $B[t]/(pt - x)\hat{}$ is $p$-torsion free. This follows from the following lemma.

**Lemma 2.2.12.** *Let $A$ be a $p$-torsion free ring, let $(p, a)$ be an $A$-regular sequence, and let $A[t]\hat{}$ be the completed polynomial algebra in $t$ over $A$. Then the ring $A[t]\hat{}/(pt - a)$ is $p$-torsion free.*

*Proof.* It is clear that $A[t]\hat{}$ is $p$-torsion free (see Lemma 7.3.4). Furthermore, its reduction modulo $p$ is the polynomial algebra $(A/pA)[t]$, which is flat over $A/pA$. Since multiplication by $a$ is injective on $A/pA$, it is also injective on $(A/pA)[t]$. Thus the sequence $(p, a)$ is $A[t]\hat{}$-regular. Now if $f$ and $g$ are elements of $A[t]\hat{}$ and $pg = f(pt - a)$, it follows that $fa = p(ft - g)$, and the $A[t]\hat{}$-regularity of $(p, a)$ implies that $f = p\tilde{f}$ for some $\tilde{f} \in A[t]\hat{}$. Then $g = \tilde{f}(pt - a)$. □

□

**Lemma 2.2.13.** *Let $Y$ be a $p$-torsion free $p$-adic formal scheme, let $X \subseteq X' \subseteq \overline{Y}$ be closed immersions, and $X'' := \pi_Y^{-1}(X) \subseteq \mathbb{D}_{X'}(Y)$. Then the natural map*

$$\mathbb{D}_{X''}(\mathbb{D}_{X'}(Y)) \to \mathbb{D}_X(Y)$$

*is an isomorphism.*

*Proof.* The map of pairs $(X'' \subseteq \mathbb{D}_{X'}(Y)) \to (X \subseteq Y)$ induces the "natural map" in the statement. On the other hand, the map of pairs $(X \subseteq Y) \to (X' \subseteq Y)$ induces a morphism $\mathbb{D}_X(Y) \to \mathbb{D}_{X'}(Y)$. Since $\overline{\mathbb{D}}_X(Y)$ is the inverse image of $X$, this map sends $\overline{\mathbb{D}}_X(Y)$ to $X''$, hence defines a $p$-adic enlargement of $X''$ over $\mathbb{D}_{X'}(Y)$, hence factors uniquely through a map $\mathbb{D}_X(Y) \to \mathbb{D}_{X''}(\mathbb{D}_Y(X))$. The two maps are inverses to each other because of the uniqueness of the factorizations. □

**Proposition 2.2.14.** *Let $j\colon Y \to Z$ be a closed immersion of $p$-torsion free $p$-adic formal schemes.*

*(1) The map $j$ factors uniquely through the $p$-adic dilatation of $Y$ in $Z$:*

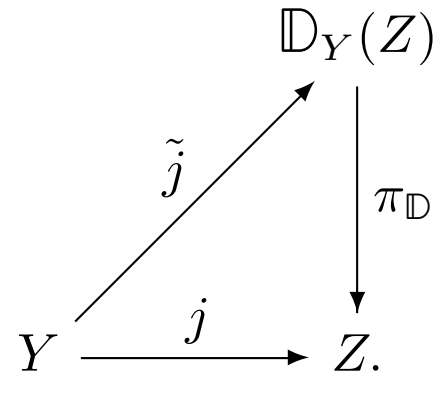

(2) *Assume that the ideal $I_{Y/Z}$ of $Y$ in $Z$ is locally finitely generated. Then the ideal $I_{Y/\mathbb{D}}$ of the closed immersion $\tilde{j}$ is generated by $\rho(I_{Y/Z})$.*

(3) *Suppose in addition that $\overline{j}\colon \overline{Y} \to \overline{Z}$ is a regular immersion. Then $j$, $\tilde{j}$, and $\overline{\tilde{j}}$ are also regular immersions, and the map $\rho$ (2.2.1) induces an isomorphism of conormal sheaves: $\mathcal{C}_{Y/Z} \to \mathcal{C}_{Y/\mathbb{D}_Y(Z)}$ whose inverse is induced by multiplication by $p$.*

*Proof.* Since $Y$ is $p$-torsion free, it is a $p$-adic enlargement of itself, so the existence and uniqueness of $\tilde{j}$ follow from the universal property of $\mathbb{D}_Y(Z)$. Let us write the rest of the proof in the affine case, with $Y = \operatorname{Spf} B$, $Z = \operatorname{Spf} C$, and $\pi_{\mathbb{D}}$ given by $\operatorname{Spf}(\theta\colon C \to C')$. If $x \in I_{Y/Z}$, then $p\rho(x) = \theta(x)$, hence

$$p\tilde{j}^\sharp(\rho(x)) = \tilde{j}^\sharp(\theta(x)) = j^\sharp(x) = 0.$$

Since $B$ is $p$-torsion free, it follows that $\tilde{j}^\sharp(\rho(x)) = 0$, hence that $\rho(x) \in I_{Y/\mathbb{D}}$. Choose a set of generators $x_1, \ldots, x_n$ for $I_{Y/Z}$. As we saw in the proof of Theorem 2.2.8, the algebra $C'$ is the $p$-adic completion of the $p$-torsion free quotient of $C[y_1, \ldots, y_n]/(py_1 - x_1, \ldots, py_n - x_n)$. In particular, $y_i = \rho(x_i)$, and these elements topologically generate $C'$. Thus every $c' \in C'$ can be written, in multi-index notation, as a $p$-adically convergent sum $c' = \sum \theta(c_I)\rho(x)^I$, with $c_I \in C$. If $c' \in I_{Y/\mathbb{D}}$, then $\tilde{j}^\sharp(c') = 0$, and hence $j^\sharp(c_0) = \tilde{j}^*(\theta(c_0)) = 0$, hence $c_0 \in I_{Y/Z}$. Then

$$c' \quad = \quad p\rho(c_0) + \sum_{I \neq 0} \theta(c_I)\rho(x)^I.$$

This completes the proof of statement (2).

Now suppose that $\overline{j}$ is a regular immersion. Then $j$ is also a regular immersion. In fact, if $(p, x_1, \ldots, x_r)$ is a regular sequence generating the ideal of $X$ in $Y$, the proof of that proposition shows that $(x_1, \ldots, x_r)$ is a regular sequence, and Proposition 2.2.11 shows that

$$C' \cong C[y_1, \ldots, y_n]\hat{\ }/(py_1 - x_1, \ldots, py_n - x_n).$$

Thus $C'_1$ identifies with the polynomial algebra $B_1[y_1, \ldots, y_n]$, in which the sequence $(y_1, \ldots, y_n)$ is regular. Then we can conclude that $(p, y_1, \ldots, y_n)$ is a regular sequence in $C'$, and in particular that $(y_1, \ldots, y_n)$ is $C'$ and $C'_1$-regular. It follows that the image of this sequence forms a basis for $\mathcal{C}_{Y/\mathbb{D}_Y(Z)}$. Since $(x_1, \ldots, x_n)$ gives a basis for $\mathcal{C}_{Y/Z}$ and $\rho(x_i) = y_i$ for all $i$, statement (3) is proved. □

Statement (2) of Theorem 2.2.8 shows that formation of $p$-adic dilatations commutes with flat base change. The flatness hypothesis is used to control possible $p$-torsion in the fiber product. The following generalization, in the case of regular immersions, will be important in applications.

**Proposition 2.2.15.** *Let $g\colon Y \to Z$ be a morphism of $p$-torsion free $p$-adic formal schemes and $i\colon X \to \overline{Z}$ a regular closed immersion.*

(1) *Suppose that the ideal of $X$ in $\overline{Z}$ is locally defined by a regular sequence which remains regular in $\mathcal{O}_{\overline{Y}}$, and let $X' := g^{-1}(X) \subseteq \overline{Y}$. Then the natural map*
$$\mathbb{D}_{X'}(Y) \to \mathbb{D}_X(Z) \times_Z Y$$
*is an isomorphism.*

(2) *Suppose that $i\colon X \to \overline{Z}$ factors as the composite of regular immersions: $X \to \overline{Y}$ and $\overline{g}\colon \overline{Y} \to \overline{Z}$. Let $\tilde{Y} \subseteq \mathbb{D}_Y(Z)$ be the image of the canonical section defined by $g$, as in (1) of Proposition 2.2.14. Then the natural map*
$$\mathbb{D}_X(Y) \to \mathbb{D}_X(Z) \times_{\mathbb{D}_Y(Z)} \tilde{Y}$$
*is an isomorphism,*

*Proof.* We can check these assertions locally, and hence we may and shall assume that $j\colon Y \to Z = \mathrm{Spf}(C \to B)$ and that $X \subseteq \overline{Z}$ is defined by a $C_1$-regular sequence $(x_1, \ldots, x_r)$. Our hypothesis in statement (1) asserts that this sequence is also $B_1$-regular. Then as we saw in Proposition 2.2.11,

$$\begin{array}{rcl} \mathbb{D}_X(Z) & = & \mathrm{Spf}\, C[t_1, \ldots, t_r]\hat{}/(pt_1 - x_1, \ldots, pt_r - x_r) \text{ and} \\ \mathbb{D}_{X'}(Y) & = & \mathrm{Spf}\, B[t'_1, \ldots, t'_r]\hat{}/(pt'_1 - x'_1, \ldots, pt'_r - x'_r), \end{array}$$

where $x'_i$ is the image of $x_i$ in $B$. Since

$$B[t'_1, \ldots, t'_r]\hat{}/(pt'_1 - x'_1, \ldots, pt'_r - x'_r) \cong B\hat{\otimes} C[t_1, \ldots, t_r]\hat{}/(pt_1 - x_1, \ldots, pt_r - x_r),$$

this proves statement (1).

In the situation of statement (2), we first find a sequence $(x_1, \ldots, x_r)$ of elements of $C$ which generates the ideal of $\overline{Y}$ in $\overline{Z}$ and is $C_1$-regular. Since $X \to \overline{Y}$ is also a regular immersion, we may then find a sequence $(y_1, \ldots, y_m)$ of elements of $C$ whose image in $B_1$ is regular and generates the ideal of $X$ in $\overline{Y}$. Then $\mathbb{D}_Y(Z) = \mathrm{Spf}\, C'$, where $C' := C[t_1, \ldots, t_r]\hat{}/(pt_1 - x_1, \ldots pt_i - x_r)$, and $\tilde{Y}$ is defined by the ideal $(t_1, \ldots, t_r)$, as we saw in (2) of Proposition 2.2.14. The ideal of $X$ in $Z$ is generated by the sequence $(p, x_1, \ldots, x_r, y_1, \ldots, y_m)$, so the ideal of $X' := X \times_Y \mathbb{D}_Y(Z)$ in $\mathbb{D}_Y(Z)$ is generated by the sequence $(p, y_1, \ldots, y_m)$ since $x_i \in pC'$. Since $C'_1$ is just a polynomial algebra over $B_1$, it follows that $(y_1, \ldots, y_r)$ is also $C'_1$-regular, so the map $X' \to \overline{\mathbb{D}_Y(Z)}$ is a regular immersion. The map $X' \cap \tilde{Y}_1 \to \tilde{Y}_1$ identifies with the map $X \to \overline{Y}$, and hence is also a regular immersion. Then statement (1), applied to the morphisms $\tilde{j}\colon \tilde{Y} \to \mathbb{D}_Y(Z)$ and $\tilde{i}\colon \tilde{X} \to \mathbb{D}_Y(Z)$, implies that the map

$$\mathbb{D}_{X' \cap \tilde{Y}}(\tilde{Y}) \to \mathbb{D}_{X'}(\mathbb{D}_Y(Z)) \times_{D_Y(Z)} \tilde{Y}$$

is an isomorphism. Lemma 2.2.13 tells us that the map $\mathbb{D}_{X'}(\mathbb{D}_Y(Z)) \to \mathbb{D}_X(Z)$ is an isomorphism, and since the map of pairs $(X' \cap \tilde{Y} \subseteq \tilde{Y}) \to (X \subseteq Y)$ is an isomorphism, so is the map $\mathbb{D}_{X' \cap \tilde{Y}}(\tilde{Y}) \to \mathbb{D}_X(Y)$. This concludes the proof of statement (2). □

**Corollary 2.2.16.** *Suppose that $X$ and $X'$ are two regularly immersed subschemes of $\overline{Y}$ which meet transversally. Then the natural map*

$$\mathbb{D}_{X\cap X'}(Y) \to \mathbb{D}_X(Y) \times_Y \mathbb{D}_{X'}(Y)$$

*is an isomorphism.*

*Proof.* We may assume that $Y = \operatorname{Spf} B$ and that $X$ (resp. $X'$) is defined by a $B/pB$-regular sequence $(x_1, \ldots, x_r)$ (resp. $(x'_1, \ldots, x'_{r'})$). Then $\mathbb{D}_{X'}(Y) = \operatorname{Spf}(B')$, where

$$B' \quad = \quad B[t'_1, \ldots, t'_{r'}]\hat{}/(pt'_1 - x'_1, \ldots, pt'_{r'} - x'_{r'}).$$

Since $X$ and $X'$ meet transversally, the sequence $(x_1, \ldots, x_r)$ remains regular in $B/(p, x'_1, \ldots, x'_{r'})$, and hence also in $B'/pB' \cong B/(p, x'_1, \ldots, x'_{r'})[t'_1, \ldots, t'_{r'}]$. Note that $X'' := X \times_Y \mathbb{D}_{X'}(Y) = \pi_Y^{-1}(X \cap X')$. Statement (1) of Proposition 2.2.15, applied to the map $\mathbb{D}_{X'}(Y) \to Y$, implies that

$$\mathbb{D}_{X''}(\mathbb{D}_{X'}(Y)) \cong \mathbb{D}_X(Y) \times_Y \mathbb{D}_{X''}(Y),$$

and Lemma 2.2.13 implies that $\mathbb{D}_{X''}(\mathbb{D}_{X'}(Y)) \cong \mathbb{D}_{X\cap X'}(Y)$. □

**Example 2.2.17** Let $i\colon X \to Y$ be the closed immersion corresponding to the map $W[x]\hat{} \to k$. Then $\mathbb{D}_X(Y) \to Y$ corresponds to the homomorphism $W[x]\hat{} \to W[t]\hat{}$ sending $x$ to $pt$, and $\mathbb{D}_X(Y) \times_Y \mathbb{D}_X(Y) \cong \operatorname{Spf} W[t_1, t_2]\hat{}/(pt_1 - pt_2)$. Thus $t_2 - t_1$ is a $p$-torsion element, and in fact the diagonal morphism $\mathbb{D}_X(Y) \to (\mathbb{D}_X(Y) \times_Y \mathbb{D}_X(Y))_{\mathrm{tf}}$ is an isomorphism.

**Remark 2.2.18.** Let $i\colon X \to Y$ be a closed immersion of a scheme into a $p$-torsion free $p$-adic formal scheme. Then we have seen that $\mathbb{D}_X(Y)$ can be empty. However, this cannot happen if $Y$ is a formal $F$-scheme. To see this, recall from [31, 2.6.1] that, for $m \gg 0$, $\mathbb{D}_{X_m}(Y)$ is not empty, where $X_m$ is the closed subscheme of $\overline{Y}$ defined by $I^m_{X/Y} + p\mathcal{O}_Y$. The map $F^m_Y$ maps $X_m$ to $X$, and hence induces a morphism $\mathbb{D}_{X_m}(Y) \to \mathbb{D}_X(Y)$. Since the first of these is not empty, neither is the second. The main idea of [31] was to consider $p$-adic enlargements of all $X_m$ as enlargements of $X$.

## 2.3 Prisms and prismatic envelopes

We now shift our attention to prisms and prismatic envelopes, as introduced in [10]. We restrict out attention to the special case of crystalline prisms, *i.e.*, the case in which the invertible ideal $I$ in Definition 3.2 of [10] is generated by $p$. We fix a a crystalline prism $(S, F_S)$ as base, and recall that, by definition. such a scheme is $p$-torsion free, so specifying the Frobenius lift is equivalent to specifying its $\delta$-structure. Since all the prisms we consider are crystalline, we may henceforth omit this qualifier.

**Definition 2.3.1.** *Let $S$ be a (crystalline) prism and let $X \to \overline{S}$ be a morphism of schemes.*

(1) *An $X/S$-prism is a morphism of (crystalline) prisms $T \to S$ endowed with an $S$-morphism $z_T \colon \overline{T} \to X$. A morphism $(T, z_T) \to (T', z_{T'})$ of $X/S$-prisms is a morphism of prisms $f \colon T \to T'$ such that $z_{T'} \circ \overline{f} = z_T$.*

(2) *An $X/S$-prism is small if $z_{\overline{T}} \colon T \to \overline{X}$ is flat.*

(3) *If $i \colon X \to Y$ is a locally closed immersion from $X$ into a crystalline prism $Y/S$, then an $X/S$-prism over $Y$ is an $X/S$-prism $(T, z_T)$ together with a morphism of prisms $\pi_T \colon T \to Y$ such that $i \circ z_T = \pi_{T|_{\overline{T}}}$.*

(4) *More generally, suppose that $i \colon X \to Y$ is a locally closed immersion, where $Y$ is a $p$-torsion free $p$-adic formal $S$-scheme endowed with an endomorphism $G_Y$ which is compatible with $F_S$ and such that $G_Y \circ i = i \circ F_X$. Then an $X/S$-prism over $Y$ is an $X$-prism $(T, z_T)$ together with a morphism of prisms $\pi_T \colon T \to Y$ such that $i \circ z_T = \pi_{T|_{\overline{T}}}$ and $G_Y \circ \pi_T = \pi_T \circ F_T$.*

$$\begin{array}{ccccc} \overline{T} & \longrightarrow & T & \xrightarrow{F_T} & T \\ \downarrow{\scriptstyle z_T} & & \downarrow{\scriptstyle \pi_T} & & \downarrow{\scriptstyle \pi_T} \\ X & \xrightarrow{i} & Y & \xrightarrow{G_Y} & Y. \end{array}$$

If $S$ is understood, we may sometimes simply say $X$-prism instead of $X/S$-prism. If $(T, z_T)$ is an $X$-prism and $g \colon T' \to T$ is a morphism of prisms, then $g$ induces a morphism $\overline{g} \colon \overline{T}' \to \overline{T}$, and $(T', \overline{g} \circ z_T)$ is also an $X$-prism. If $f \colon X' \to X$ is a morphism of $S$-schemes and $(T, z_T)$ is an $X'$-prism, then $(T, f \circ z_T)$ is an $X$-prism.

**Definition 2.3.2.** *Let $i \colon X \to Y$ be a locally closed immersion as in (4) of Definition 2.3.1. The prismatic envelope of $X$ in $Y$, denoted by $\mathbb{\Delta}_X(Y)$, is the final object of the category of $X$-prisms over $Y$.*

Slightly more generally, if $X \to S$ is a morphism, we allows ourselves to identify $\mathbb{\Delta}_{\overline{X}}(Y)$ and $\mathbb{\Delta}_X(Y)$, as we did for dilatations. Then the map $\overline{X} \to X$ induces an isomorphism on the respective categories of prisms.

The existence of prismatic envelopes in the situation (3) of Definition 2.3.1, when $X \to \overline{Y}$ is a regular immersion, is proved in [10, 3.13], using the formalism of $\delta$-rings. In attempt to make its structure more explicit, we give several constructions here, some of which work even in the situation (4). For an example of an application of this extra generality, see Proposition 2.3.10. Each of our constructions involves a sequence of $p$-adic dilatations. We note that the construction in [10] works for more general prisms (not necessarily crystalline).

We begin with an explicit description in terms of local coordinates.

**Theorem 2.3.3.** *If $Y$ is a crystalline prism and $X \to \overline{Y} \to Y$ is a regular immersion, let $\mathbb{\Delta}_X(Y)$ be the prismatic envelope of $X$ in $Y$.*

*(1)* $\Delta_X(Y)$ *can be obtained by applying the relative right adjoint* $\mathrm{J}_Y$ *(1.6.11) to the p-adic dilatation of* $X$ *in* $Y$*. Thus*

$$\Delta_X(Y) \cong \mathrm{J}_Y(\mathbb{D}_X(Y)).$$

*The map* $\overline{\Delta}_X(Y) \to X$ *is faithfully flat, and in particular,* $(\Delta_X(Y), z_\Delta)$ *is small.*

*(2) Suppose in fact that* $Y = \mathrm{Spf}(B, \phi)$ *and that* $X \subseteq Y$ *is defined by an ideal generated by a B-regular sequence* $(p, x_1, \ldots, x_r)$*.*

*(a) The prismatic envelope* $\Delta_X(Y)$ *is the formal spectrum of the p-adic completion of the ring:*

$$B^\infty := B[t_{i,j}]/(pt_{i,j+1} - \delta^j(x_i) + t_{i,j}^p : i = 1, \ldots, r, \quad j \in \mathbf{N}),$$

*where* $\{t_{i,j} : 1 \leq i \leq r, j \in \mathbf{N}\}$ *are free variables for* $j > 0$ *and* $t_{i,0} := 0$*.*

*(b) The (p-torsion free) ring* $B^\infty$ *admits a unique δ-structure compatible with that of B. It is given explicitly by the inductive formula:*

$$\delta(t_{i,j}) = t_{i,j+1} - P(t_{i,j}, t_{i,j-1}, \delta(t_{i,j-1})) \quad \text{if } j > 0,$$

*where P is a universal polynomial made explicit in Lemma 2.3.4 below.*

*(c) Consequently, if* $j > 0$*,* [4]

$$\delta(t_{i,j}) \equiv t_{i,j+1} + t_{i,j-1}^{p(p-1)}(\delta(t_{i,j-1}) - t_{i,j}) + p^{p-2}((\delta(t_{i,j-1}))^p - t_{i,j}^p) \pmod p.$$

*Proof.* We define a sequence of rings $B = B^{(0)} \subseteq B^{(1)} \subseteq B^{(2)} \subseteq \cdots$, where $B^{(j)}$ is the subring of $B^\infty$ generated by $\{t_{i,j'} : j' \leq j\}$. Then

$$B^{(j+1)} := B^{(j)}[t_{1,j+1}, \ldots, t_{r,j+1}]/(pt_{i,j+1} - \delta^j(x_i) + t_{i,j}^p : i = 1, \ldots, r),$$

and $B^\infty = \varinjlim\{B^{(j)} : j \in \mathbf{N}\}$. For each $j \in \mathbf{N}$, let $I^{(j)}$ be the ideal of $B^{(j)}$ generated by $p$ and $\{\delta^j(x_i) - t_{i,j}^p : i = 1, \ldots, r\}$. Thus $I^{(0)}$ is the ideal of $X$ in $Y$.

We shall prove the following statements by induction on $j$.

(1) For each $j > 0$, $B^{(j)}$ corresponds to the (uncompleted) $p$-adic dilatation of the ideal $I^{(j-1)}$ in $B^{(j-1)}$. In particular, each of these rings is $p$-torsion free. Furthermore, each homomorphism $B/I \to B^{(j)}/pB^{(j)}$ is faithfully flat.

[4] This congruence was pointed out by the referee, and is needed to fill a gap caused by a calculation error in a previous version.

(2) For each $j > 0$, the ring $B^{(j)}$ admits a unique endomorphism $\psi^{(j)}$ compatible with the endomorphism $\phi$ of $B$, and

$$\psi^{(j)}(t_{i,j}) \equiv \delta^j(x_i) \pmod{pB^{(j)}}.$$

(3) $\psi^{(j)}(b) \equiv b^p \pmod{I^{(j)}}$ for all $b \in B^{(j)}$.

It follows from the definition that $B^{(0)} = B$. Since $(x_1, \ldots, x_r)$ is $B/pB$-regular, Proposition 2.2.11 shows that

$$B^{(1)} := B[t_{1,1} \ldots, t_{r,1}]/(pt_{1,1} - x_1, \ldots, pt_{r,1} - x_r)$$

is $p$-torsion free and is in fact the (uncompleted) $p$-adic dilatation of $I$. Furthermore, $B^{(1)}/pB^{(1)}$ is the polynomial ring $(B/I)[t_{1,1}, \ldots, t_{1,r}]$, hence is faithfully flat over $B/I$. Thus statement (1) holds for $j = 1$. It also follows that the sequence

$$(p, t_{1,1}^p - \delta(x_1), \ldots, t_{r,1}^p - \delta(x_r))$$

is $B^{(1)}$-regular. This implies that $B^{(2)}$ is the uncompleted $p$-adic dilation of $I^{(1)}$ in $B^{(1)}$, and in particular is $p$-torsion free. Since

$$B^{(2)}/pB^{([2)} \cong \frac{(B/I)[t_{1,1}, \ldots, t_{r,1}, t_{1,2}, \ldots t_{r,2}]}{(t_{1,1}^p - \delta(x_1), \ldots, t_{r,1}^p - \delta(x_r))},$$

it too is faithfully flat over $B/I$. This same argument works for all $j$, proving (1).

To prove (2) and (3) when $j = 1$, we use the fact that $B^{(1)}$ is $p$-torsion free and that $t_{i,1} = p^{-1}x_i$. Since $\phi$ is a ring homomorphism, it induces an endomorphism of $B \otimes \mathbf{Q}$. Calculating in $B \otimes \mathbf{Q}$, we find

$$\begin{aligned} (\phi \otimes \mathrm{id})(t_{i,1}) &= p^{-1}\phi(x_i) \\ &= p^{-1}(x_i^p + p\delta(x_i)) \\ &= p^{p-1}t_{i,1}^p + \delta(x_i). \end{aligned}$$

Since $B^{(1)}$ is the subring of $B \otimes \mathbf{Q}$ generated by $\{t_{1,1}, \ldots, t_{r,1}\}$, this shows that $\phi$ induces an endomorphism $\psi^{(1)}$ of $B^{(1)}$. We also see that $\psi^{(1)}(t_{i,1}) \equiv \delta(x_i) \pmod{pB^{(1)}}$. Since $\delta(x_i) \equiv t_{i,1}^p \pmod{I^{(1)}}$, condition (2) implies that (3) holds for each of $t_{1,1}, \ldots, t_{r,1}$, and since it also holds for elements of $B^{(0)}$, it holds for all elements of $B^{(1)}$, completing the proof when $j = 1$.

To proceed by induction, suppose that (2) and (3) are satisfied for $j$. The congruence in (2) implies that

$$(\psi^{(j)}(t_{i,j}))^p \equiv (\delta^j(x_i))^p \pmod{p^2B^{(j)}}.$$

Now we compute in $B^{(j)}$:

$$\begin{aligned} \psi^{(j+1)}(pt_{i,j+1}) &= \phi(\delta^j(x_i)) - \psi^{(j)}(t_{i,j}^p) \\ &= (\delta^j(x_i))^p + p\delta^{j+1}(x_i) - \psi^{(j)}(t_{i,j}^p) \end{aligned}$$

$$\begin{aligned} &\equiv (\delta^j(x_i))^p + p\delta^{j+1}(x_i) - (\delta^j(x_i))^p \pmod{p^2B^{(j)}} \\ &\equiv p\delta^{j+1}(x_i) \pmod{p^2B^{(j)}}. \end{aligned}$$

This implies that $\psi^{(j+1)}(t_{i,j+1})$ actually lies in $B^{(j)}$ and in fact is congruent to $\delta^{j+1}(x_i) \pmod{pB^{(j)}}$. It follows that $\psi^{(j+1)}$ is well-defined on all of $B^{(j+1)}$ and satisfies (2). Since $\delta^{j+1}(x_i) - t^p_{i,j+1}$ and $\psi^{(j+1)}(t_{i,j+1}) - \delta^{j+1}(x_i)$ belong to $I^{(j+1)}$, it follows that $\psi^{(j+1)}(t_{i,j+1}) - t^p_{i,j+1} \in I^{(j+1)}$, so statement (3) follows as before.

One concludes easily that $B^\infty \subseteq B \otimes \mathbf{Q}$ is stable under $\phi \otimes \mathrm{id}$. Furthermore, since $I^{(j)}B^{(j+1)} \subseteq pB^{(j+1)}$ for all $j$, it follows that the endomorphism induced by $\phi \otimes \mathrm{id}$ is indeed a Frobenius lift. By Lemma 7.3.4, the $p$-adic completion $\hat{B}^\infty$ of $B^\infty$ is $p$-torsion free, and hence is a crystalline prism, and since $I\hat{B}^\infty = p\hat{B}^\infty$, it defines an $X$-prism over $\mathrm{Spf}\, B$.

If $T \to Y$ is a morphism of crystalline prisms and $\overline{T} \to \overline{Y}$ factors through $X$, it is easy to check, using the fact that $T$ is torsion free, that the $\delta$-morphism $T \to Y$ factors uniquely through $\mathrm{Spf}\, \hat{B}^\infty$. We leave this to the reader, proceeding with the more delicate and stronger fact that this universality holds even when $T$ has torsion. This will prove that $\Delta_X(Y) = \mathrm{J}_Y(\mathbb{D}_X(Y))$, as claimed in statement (1) of the theorem. The proof will depend on the explicit formula in statement (2b), which depends in turn on the following general fact about $\delta$-rings.

**Lemma 2.3.4.** *Define polynomials with integral coefficients as follows:*

$$\begin{aligned} C(X,Y) &:= p^{-1}\left((X+Y)^p - X^p - Y^p\right) \qquad \textit{(see also 1.6.4)} \\ C_p(X,Y) &:= p^{-1}C(X,pY) \\ P(X,Y,Z) &= C_p(Y^p,Z) - C_p(Y^p,X) + p^{p-2}(Z^p - X^p). \end{aligned}$$

*Then:*

*(1)* $C_p(X,Y) \equiv X^{p-1}Y \pmod p$

*(2)* $P(X,Y,Z) \equiv Y^{p(p-1)}(Z-X) + p^{(p-2)}(Z^p - X^p) \pmod p$

*(3) In any $\delta$-ring, we have:*

$$\delta(pa + b^p) = p\delta(a) + a^p + pP(a,b,\delta(b)).$$

*Proof.* The congruences (1) and (2) follow immediately from the definitions. Let us first prove (3) when $a = 0$. We claim that

$$\delta(b^p) = pC_p(b^p, \delta(b)) + p^{p-1}\delta(b)^p.$$

It suffices to prove this when $b$ is the element $x$ of the universal $\delta$-ring $(\mathbf{Z}[x_\bullet], \delta)$ described in Example 1.6.4. Since this ring is $p$-torsion free, we can calculate as follows:

$$\begin{aligned} p\delta(x^p) &= \phi(x^p) - (x^p)^p \\ &= \phi(x)^p - x^{p^2} \end{aligned}$$

$$
\begin{array}{rcl}
 & = & (x^p + p\delta(x))^p - x^{p^2} \\
 & = & x^{p^2} + p^p(\delta(x))^p + pC(x^p, p\delta(x)) - x^{p^2} \\
\delta(x^p) & = & C(x^p, p\delta(x)) + p^{p-1}(\delta(x))^p \\
 & = & pC_p(x^p, \delta(x)) + p^{p-1}(\delta(x))^p.
\end{array}
$$

If $b = 0$, we compute from the definiton of a $\delta$-ring:

$$
\begin{array}{rcl}
\delta(pa) & = & \delta(p)a^p + p^p\delta(a) + p\delta(p)\delta(a) \\
 & = & (1 - p^{p-1})a^p + p^p\delta(a) + (p - p^p)\delta(a) \\
 & = & a^p + p\delta(a) - p^{p-1}a^p.
\end{array}
$$

Finally:

$$
\begin{array}{rcl}
\delta(pa + b^p) & = & \delta(pa) + \delta(b^p) - C(b^p, pa) \\
 & = & a^p + p\delta(a) - p^{p-1}a^p + pC_p(b^p, \delta(b)) + p^{p-1}(\delta(b))^p - pC_p(b^p, a) \\
 & = & a^p + p\delta(a) + pP(a, b, \delta(b)).
\end{array}
$$

□

We prove statement (2b) of the theorem as follows. For $j \geq 1$, we have:

$$\delta^{j-1}(x_i) = pt_{i,j} + t_{i,j-1}^p.$$

Then the Lemma 2.3.4 implies that

$$\delta^j(x_i) = p\delta(t_{i,j}) + t_{i,j}^p + pP(t_{i,j}, t_{i,j-1}, \delta(t_{i,j-1})).$$

We also have $\delta^j(x_i) = pt_{i,j+1} + t_{i,j}^p$, so:

$$pt_{i,j+1} = p\delta(t_{i,j}) + pP(t_{i,j}, t_{i,j-1}, \delta(t_{i,j-1})).$$

Since our envelope is $p$-torsion free, we conclude that, if $j \geq 1$:

$$t_{i,j+1} = \delta(t_{i,j}) + P(t_{i,j}, t_{i,j-1}, \delta(t_{i,j-1})). \tag{2.3.1}$$

This is statement (2b) of the theorem Then (2c) follows from the congruence (2) in Lemma 2.3.4.

Statement (1) of the theorem will follow if show that for every $\delta$-ring $C$, a homomorphism $\theta\colon B^{(1)} \to C$ extends uniquely to a homomorphism of $\delta$-rings $\theta\colon (B^\infty, \delta) \to (C, \delta)$. Let $c_i := \theta(x_i)$ and let $c_{i,1} := \theta(t_{i,1})$. Then $pc_{i,1} = c_i$ for all $i$. We must show that there is then a unique family of elements $c_{i,j}$ such that $c_{i,0} = 0$ and

$$
\begin{array}{rcll}
pc_{i,j} & = & \delta^{j-1}(c_i) - c_{i,j-1}^p & (2.3.2) \\
\delta(c_{i,j}) & = & c_{i,j+1} - P(c_{i,j}, c_{i,j-1}, \delta(c_{i,j-1})). & (2.3.3)
\end{array}
$$

The second of these equations is equivalent to:

$$c_{i,j+1} = \delta(c_{i,j}) + P(c_{i,j}, c_{i,j-1}, \delta(c_{i,j-1})).$$

This is already enough to prove the uniqueness, by induction on $j$. In fact it also allows us to construct the elements $c_{i,j}$ for $j \geq 2$, by induction on $j$. Note that the first of these equations holds when $j = 1$ by definition. Now assume that we have constructed $c_{i,j}$ satisfying the first of the equations for $j' \leq j$ and the second for $j' < j$. Then define $c_{i,j+1}$ by the second equation. We know already that $pc_{i,j} = \delta^{j-1}(c_i) - c_{i,j-1}^p$, so Lemma 2.3.4 implies that

$$p\delta(c_{i,j}) = \delta^j(c_i) - c_{i,j}^p - pP(c_{i,j}, c_{i,j-1}, \delta(c_{i,j-1})).$$

Hence:

$$\begin{aligned} pc_{i,j+1} &= p\delta(c_{i,j}) + pP(c_{i,j}, c_{i,j-1}, \delta(c_{i,j-1})) \\ &= \delta^j(c_i) - c_{i,j}^p. \end{aligned}$$

Thus the first equation is also satisfied. □

**Remark 2.3.5.** The description of prismatic envelopes in Theorem 2.3.3 looks quite explicit, and perhaps appealing, but in some cases the computation of the iterates $\delta^n$ can be daunting, depending on the choices of generators for the ideal of $X$ in $Y$.

**Corollary 2.3.6.** *Suppose that $Y \to S$ is a morphism of crystalline prisms, that $X \to \overline{S}$ is flat, and that $X \to \overline{Y}$ is a regular immersion. Then $\Delta_X(Y) \to S$ is $p$-completely flat.*

*Proof.* By construction, $\Delta_X(Y)$ is $p$-torsion free, and statement (2) of Theorem 2.3.3 tells us that $z_{\Delta} \colon \overline{\Delta}_X(Y) \to X$ is flat. Then $\overline{\Delta}_X(Y) \to \overline{S}$ is also flat, and Proposition 7.3.2 tells us that $\Delta_X(Y) \to S$ is $p$-completely flat. □

Statement (1) of Theorem 2.3.3 allows us to give a useful description of the functor of points of the scheme underlying prismatic envelopes. This is a prismatic version of the functor of points description of dilations described in Definition 2.2.4, and it underlies Drinfeld's "prismatisation" functor. We shall need to apply it when the "test" scheme $T$ is not reduced, in which case $\mathrm{W}(T)$ may have torsion.

**Proposition 2.3.7.** *Let $Y$ be a crystalline prism and and let $X \to \overline{Y}$ be a regular immersion. If $t_Y \colon T \to Y$ is a morphism of formal schemes, let $t_{\mathrm{W}} \colon \mathrm{W}(T) \to Y$ be the corresponding morphism of $\delta$-schemes. Then there are natural bijections among the following:*

*(1) morphisms of formal schemes $t_{\Delta} \colon T \to \Delta_X(Y)$ such that $\pi_{\Delta} \circ t_{\Delta} = t_Y$.*

*(2) morphisms of $\delta$-schemes $w_{\Delta} \colon \mathrm{W}(T) \to \Delta_X(Y)$ such that $\pi_{\Delta} \circ w_{\Delta} = t_{\mathrm{W}}$.*

*(3) morphisms of formal schemes: $w_{\mathbb{D}} \colon \mathrm{W}(T) \to \mathbb{D}_X(Y)$ such that $\pi_{\mathbb{D}} \circ w_{\mathbb{D}} = t_{\mathrm{W}}$.*

*(4) $\mathcal{O}_Y$-linear homomorphisms $\rho_W \colon I_{X/Y} \to W(\mathcal{O}_T)$ such that $\rho_W(p) = 1$.*

*See the diagram:*

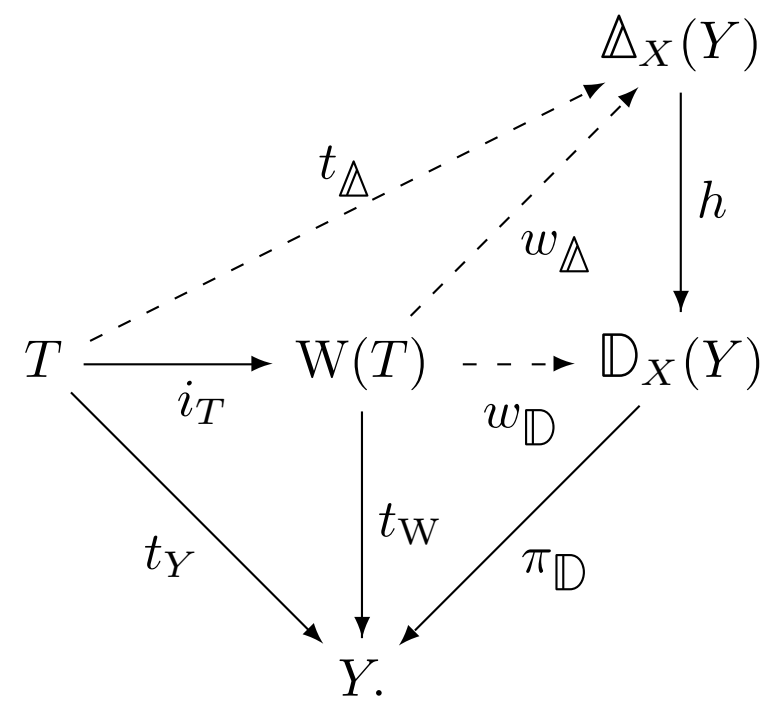

*Proof.* The equivalence of (1) and (2) is a consequence of the fact that W is left adjoint to the forgetful functor. Given a morphism $w_{\mathbb{\Delta}}$ as in (2), one obtains a morphism $w_{\mathbb{D}}$ as in (3) by composing with the map $h\colon \mathbb{\Delta}_X(Y) \to \mathbb{D}_X(Y)$. Conversely, if $w_{\mathbb{D}}$ is as in (3), there is a unique morphism of $\delta$-$Y$-schemes $W(T) \to \mathrm{J}_Y(\mathbb{D}_X(Y))$, and Theorem 2.3.3 tell us that $\mathbb{\Delta}_X(Y) \cong \mathrm{J}_Y(\mathbb{D}_X(Y))$. The equivalence of (3) and (4) comes from the isomorphism $\mathbb{D}_X(Y) \cong \mathbb{V}_X(Y)$ of statement (4) of Theorem 2.2.8 and the universal property of the latter. □

Let us now prove the existence of prismatic envelopes in the more general situation (4) of Definition 2.3.1. This extra generality is also needed because of the inductive structure of this proof.

**Theorem 2.3.8.** *Suppose that $Y$ is $p$-torsion free $p$-adic formal scheme with an endomorphism $G_Y$ and that $i\colon X \to \overline{Y} \to Y$ is a closed immersion such that $G_Y \circ i = i \circ F_X$. Then the prismatic envelope of $X$ in $Y$ is represented by an $X$-prism $(\mathbb{\Delta}_X(Y), z_{\mathbb{\Delta}})$ endowed with a map $\pi_{\mathbb{\Delta}}\colon \mathbb{\Delta}_X(Y) \to Y$. The morphisms $\pi_Y$ and $z_{\mathbb{\Delta}}$ are affine, but not necessarily of finite type. The morphism $\pi_{\mathbb{\Delta}}$ can be viewed as a completion of a limit of $p$-adic dilatations. Specifically, there exists a sequence of morphisms of closed immersions:*

$$\cdots \to (X^{(n+1)} \subseteq Y^{(n+1)}) \to (X^{(n)} \subseteq Y^{(n)}) \to \cdots \to (X \subseteq Y)$$

*such that each $Y^{(n+1)} \to Y^{(n)}$ is the $p$-adic dilatation of $X^{(n)}$ in $Y^{(n)}$ and such that the $p$-adic completion of $\varprojlim Y^{(n)}$ is the prismatic envelope $\mathbb{\Delta}_X(Y)$ of $X$ in $Y$. Furthermore, the following statements are verified.*

*(1) Formation of $\mathbb{\Delta}_X(Y)$ is functorial: a morphism $g\colon Y' \to Y$ sending $X' \subseteq Y'$ to $X \subseteq Y$ induces a morphism $\mathbb{\Delta}_{X'}(Y') \to \mathbb{\Delta}_X(Y)$. If $g$ is $p$-completely flat and $X' = g^{-1}(X)$, then the resulting morphism $\mathbb{\Delta}_{X'}(Y') \to \mathbb{\Delta}_X(Y) \times_Y Y'$ is an isomorphism.*

*(2) Suppose that $X \to \overline{Y}$ is a regular immersion. Then the same is true of each $X^{(n)} \to \overline{Y}^{(n)}$, and $\overline{\mathbb{\Delta}}_X(Y) \to X$ is faithfully flat. In particular, $\mathbb{\Delta}_X(Y)$ is small.*

*Proof.* The main computations are in the following lemma, which we formulate in the affine setting.

**Lemma 2.3.9.** *Let $I$ be an ideal in a $p$-torsion free ring $\mathbf{Z}_p$-algebra $B$. Assume that $p \in I$ and that $B$ is equipped with an endomorphism $\psi$ such that $\psi(b) \equiv b^p \pmod{I}$ for all $b \in B$; note that $\psi$ necessarily maps $I$ to $I$. Let $\theta\colon B \to B'$ be the homomorphism corresponding to the (not necessarily completed) $p$-adic dilatation of $I$ in $B$, let $\psi'$ be the endomorphism of $B'$ induced by $\psi$, and let $I' \subseteq B'$ be the ideal of $B'$ generated by $\{\psi'(b') - b'^p : b' \in B'\}$, so that $\psi'(b') \equiv b'^p \pmod{I'}$ for all $b' \in B'$. Define:*

$$\begin{array}{rcll} \rho\colon I \to B' & : & x \mapsto p^{-1}\theta(x) & \text{(see (2.2.1))} \\ \epsilon : B \to I & : & b \mapsto \psi(b) - b^p & \\ \delta : B \to B' & : & b \mapsto \rho(\epsilon(b)). & \end{array}$$

*Let $\theta'\colon B' \to B''$ be the $p$-adic dilatation of $I'$ and, analogously:*

$$\begin{array}{rcl} \rho'\colon I' \to B'' & : & x' \mapsto p^{-1}\theta'(x') \\ \epsilon' : B' \to I' & : & b' \mapsto \psi'(b') - b'^p \\ \delta' : B' \to B'' & : & b' \mapsto \rho'(\epsilon'(b)). \end{array}$$

*Then the following statements are verified.*

*(1)* $\delta(1) = 0$, *and if $a, b \in B$, then*

$$\begin{array}{rcl} \delta(a+b) & = & \delta(a) + \delta(b) - \theta(C(a,b)) \\ \delta(ab) & = & \theta(a)^p\delta(b) + \delta(a)\theta(b)^p + p\delta(a)\delta(b). \end{array}$$

*(2) If $x, y \in I$ and $b \in B$, then*

$$\begin{array}{rcl} \delta(x+y) & \equiv & \delta(x) + \delta(y) \pmod{pB'} \\ \delta(bx) & \equiv & \theta(b^p)\delta(x) \pmod{pB'}. \end{array}$$

*(3) If $x \in I$, then*

$$\begin{array}{rcl} \rho(\psi(x)) & = & \psi'(\rho(x)) \\ \rho(\psi(x)) & \equiv & \delta(x) \pmod{pB'} \\ \epsilon'(\rho(x)) & \equiv & \delta(x) - \rho(x)^p \pmod{pB'}. \end{array}$$

*(4) $\epsilon' \circ \rho$ induces a $p$-linear homomorphism*

$$\epsilon'_\rho\colon I/pB \to I'/pB'.$$

*(5) If $(p, x_1, \ldots, x_r)$ generates $I$, then $(p, x'_1, \ldots, x'_r)$ generates $I'$, where $x'_i := \epsilon'(\rho(x_i))$. If $(p, x_1, \ldots, x_r)$ is $B$-regular then $(p, x'_1, \ldots, x'_r)$ is $B'$-regular, $B'/I'$ is faithfully flat over $B/I$, and $\epsilon'_\rho$ induces an isomorphism*

$$B'/I' \otimes_B F^*(\overline{I}/\overline{I}^2) \to \overline{I}'/\overline{I}'^2,$$

*where $\overline{I}$ (resp. $\overline{I}'$) is the image of $I$ in $B/pB$ (resp., of $I'$ in $B'/pB'$).*

*Proof.* Statement (1) is straightforward:

$$\begin{aligned}
\delta(1) &= \rho(\psi(1) - 1^p) \\
&= 0. \\
\delta(a+b) &= \rho\left(\psi(a+b) - (a+b)^p\right) \\
&= \rho\left(\psi(a) + \psi(b) - (a^p + b^p + pC(a,b)\right) \\
&= \rho\left(\psi(a) - a^p + \psi(b) - b^p - pC(a,b)\right) \\
&= \rho(\psi(a) - a^p) + \rho(\psi(b) - b^p) - \theta(C(a,b))\rho(p)) \\
&= \delta(a) + \delta(b) - \theta(C(a,b)). \\
\delta(ab) &= \rho(\psi(ab) - (ab)^p) \\
&= \rho(\psi(a)\psi(b) - a^p b^p) \\
&= \rho(\psi(a)\psi(b) - a^p\psi(b) - \psi(a)b^p + a^p b^p + \psi(a)b^p - a^p b^p + a^p\psi(b) - a^p b^p) \\
&= \rho(\psi(a) - a^p)(\psi(b) - b^p) + b^p(\psi(a) - a^p) + a^p(\psi(b) - b^p)) \\
&= \theta(\psi(b) - b^p)\rho(\psi(a) - a^p) + \rho(b^p(\psi(a) - a^p) + \rho(a^p(\psi(b) - b^p)) \\
&= p\rho(\psi(b) - b^p)\rho(\psi(a) - a^p) + \theta(b^p)\rho(\psi(a) - a^p) + \theta(a^p)\rho(\psi(b) - b^p) \\
&= p\delta(b)\delta(a) + \theta(b)^p\delta(a) + \theta(a)^p\delta(b).
\end{aligned}$$

Statement (2) is an immediate consequence, since $C(x,y)$ and $x^p$ belong to $I$ if $x$ and $y$ do and $IB' = pB'$.

If $x \in I$, then $\psi(x) \in I$, and we have:

$$\begin{aligned}
\theta(\psi(x)) &= \theta(x)^p + p\delta(x) \\
p\rho(\psi(x)) &= (p\rho(x))^p + p\delta(x) \\
\rho(\psi(x)) &= p^{p-1}(\rho(x))^p + \delta(x) \\
\rho(\psi(x)) &\equiv \delta(x) \pmod{pB'}.
\end{aligned}$$

Since $\epsilon'(\rho(x)) = \psi'(\rho(x)) - \rho(x)^p = \rho(\psi(x)) - \rho(x)^p$, this proves (3). Since the maps $I \to B'/pB'$ induced by $\delta$ and $\rho^p$ are $p$-linear, so is the map induced by $\epsilon' \circ \rho$. Since $\epsilon'(\rho(p)) = \epsilon'(1) = 0$, it follows that this map annihilates $pB$ and hence factors through $I/pB$. This proves (4).

To prove (5), note first that $\epsilon'(\rho(I)) \subseteq I'$ by the definition of $I'$. Note also that $p \in I'$, since $\epsilon'(p) = p(1 - p^{p-1})$. In fact we claim that $I'$ is generated by $p$ and $\epsilon'(\rho(I))$. It will suffice to show that, for every $b' \in B'$, $\psi'(b') - b'^p$ belongs to the ideal $I''$ so generated. First suppose that $b' = \theta(b)$, with $b \in B$. Since $\psi(b) - b^p \in I$, it follows that

$$\psi'(b') - b'^p = \theta(\psi(b) - b^p) \in \theta(I) \subseteq pB'.$$

On the other hand, if $x \in I$, then $\psi'(\rho(x)) - \rho(x)^p = \epsilon'(\rho(x)) \in I''$. Since the $B$-algebra $B'$ is topologically generated by elements $\rho(x)$ with $x \in I$, the same congruence holds for all $b' \in B'$. This proves our claim. Since the map $\epsilon'_\rho$ is $p$-linear, in fact $I'$ is generated by $p$ together with the image of any set of generators for $I$.

Suppose that $I$ is generated by a $B$-regular sequence $(p, x_1, \ldots, x_r)$. Let $\overline{B} := B/pB$ and $\overline{B}' := B'/$. Then $I\overline{B}$ is generated by the $\overline{B}$-regular sequence $(\overline{x}_1, \ldots, \overline{x}_r)$, where $\overline{x}_i$ is the image of $x_i$ in $\overline{B}$, and $I'\overline{B}'$ is generated by the sequence $(\epsilon'_\rho(\overline{x}_1), \ldots, \epsilon'_\rho(\overline{x}_r))$. We claim that this sequence is $\overline{B}'$-regular. As we saw in Proposition 2.2.11, $\overline{B}'$ is isomorphic to the polynomial algebra $(B/I)[y_1, \ldots, y_r]$, where $y_i$ is the reduction modulo $pB'$ of $\rho(x_i)$. Each $\epsilon(x_i)$ belongs to $I$, and so there is a sequence $b_{i,0}, \ldots, b_{i,r}$ in $B$ such that $\epsilon(x_i) = pb_{i,0} + \sum_1^r b_{i,j}x_j$. Then

$$\begin{aligned}\delta(x_i) &= \rho(\epsilon(x_i)) \\ &= \theta(b_{i,0}) + \sum_{j=1}^{r} \theta(b_{i,j})\rho(x_j).\end{aligned}$$

In the ring $\overline{B}' \cong (B/I)[y_1, \ldots, y_r]$, we have:

$$\begin{aligned}\epsilon'_\rho(\overline{x}_i) &= \overline{\delta}(x_i) - \overline{\rho}(x_i)^p \\ &= \overline{b}_{i,0} + \sum \overline{b}_{i,j} y_j - y_i^p.\end{aligned}$$

This is a monic polynomial whose leading term is $-y_i^p$. Then one sees easily by induction that the sequence $(\epsilon'_\rho(\overline{x}_1), \ldots, \epsilon'_\rho(\overline{x}_r))$ is $\overline{B}'$-regular and that the quotient is faithfully flat over $B/IB$.

The linear map $F^*\overline{I} \to \overline{I}'$, induced by $\epsilon'_\rho$ is is surjective by construction. If $(p, x_1, \ldots, x_r)$ is a $B$-regular sequence, then $(\overline{x}_1, \ldots, \overline{x}_r)$ reduces to a basis for $\overline{I}/\overline{I}^2$, and we have just seen that $(\epsilon'_\rho(\overline{x}_1), \ldots, \epsilon'_\rho(\overline{x}_r))$ is a $\overline{B}'$-regular sequence which generates $\overline{I}'$. Hence it induces a basis for $\overline{I}'/\overline{I}'^2$, concluding the proof of statement (5). □

Let us now prove Theorem 2.3.8 assuming that $Y = \operatorname{Spf} B$ and $X = \operatorname{Spec} B/I$. Using the lemma, we inductively find a sequence:

$$(B, I, \psi) \to (B', I', \psi^{(1)}) \to (B'', I'', \psi^{(2)}) \to \cdots \to (B^{(n)}, I^{(n)}, \psi^{(n)}) \to \cdots,$$

where $\psi^{(n)}$ is an endomorphism of $B^{(n)}$ such that $\psi^{(n)}(b) - b^p \in I^{(n)}$ for all $b \in B^{(n)}$. Here each $B^{(n)} \to B^{(n+1)}$ is the (not necessarily completed) $p$-adic dilatation of $I^{(n)}$, and, for $n > 0$, $I^{(n)}$ is the ideal of $B^{(n)}$ generated by $\{\psi^{(n)}(b) - b^p : b \in B^{(n)}\}$. The endomorphisms $\psi^{(n)}$ are compatible with the transition maps $B^{(n)} \to B^{(n+1)}$. Since $I^{(n)}B^{(n+1)} = pB^{(n+1)}$, it follows that the image of $\psi^{(n)}(b) - b^p$ in $B^{(n+1)}$ is divisible by $p$. Let $(\mathcal{B}, \psi) := \varinjlim(B^{(n)}, \psi^{(n)})$. Since $\mathcal{B}$ is $p$-torsion free, its $p$-adic completion $\mathcal{B}\hat{}$ is also (see Lemma 7.3.4), and so defines a $p$-adic enlargement of $X$ in $Y$. If $b \in \mathcal{B}$ comes from an element $b_n$ of $B^{(n)}$, then $\psi^{(n)}(b_n) \equiv b_n^p \pmod{pB^{(n+1)}}$ and hence $\psi(b) \equiv b^p \pmod{p\mathcal{B}}$. It follows that this congruence also holds on all of $\mathcal{B}\hat{}$, so $\operatorname{Spf}(\mathcal{B}\hat{})$ defines an $X$-prism $Y'$ over $Y$. If $(T, \pi_T, z_T)$ is another, we claim that $\pi_T$ factors uniquely through $Y'$. We may check this locally on $T$, and so may assume that $T$ is affine, say $T = \operatorname{Spf} C$. Since $T$ is a $p$-adic enlargement of $X$, the map $B \to C$ factors through $B' = B^{(1)}$, and the endomorphisms $\phi$ of $C$ and $\psi^{(1)}$ of $B^{(1)}$ are compatible. If $c \in C$ is any element of $C$, necessarily $\phi(c) - c^p \in pC$. This applies if $c$ is the image of any element $b'$ of $B'$ and in particular if $b' = \rho(x)$ for some $x \in I$. Thus we see that $\epsilon'(x)$ becomes divisible by $p$ in $C$. Thus the homomorphism $B^{(1)} \to C$ factors through $B^{(2)}$. Continuing in this way, we find a homomorphism $\mathcal{B} \to C$ and then $\mathcal{B}\hat{} \to C$. Each step is unique, and hence so is this homomorphism. We conclude that $Y'$ is indeed a final $X$-prism over $Y$.

The proofs of functoriality of prismatic envelopes and their compatibility with base change is the same as that for $p$-adic dilatations; of course they also can be found in [10].

If $X \to \overline{Y}$ is a regular immersion, we have seen that the same is true of each of the closed immersions defining each of the successive dilatations, and it follows from Theorem 2.2.8 that each map $B^{(n)}/I^{(n)} \to B^{(n+1)}/I^{(n+1)}$ is faithfully flat. Then the same is true for the map

$$B/I \to \varinjlim B^{(n)}/I^{(n)} = \mathcal{B}/p\mathcal{B} = \mathcal{B}\hat{}/p\mathcal{B}\hat{}.$$

□

The construction in Theorem 1.7.11 used the fact that the Frobenius lift is locally surjective in the $p$-completely flat topology for crystalline prisms. This is proved by a universal construction in [10, 2.12]. It is also an easy consequence of Theorem 2.3.8, as we explain in the proof of the following proposition.

**Proposition 2.3.10.** *Let $T$ be a crystalline prism and let $b_\cdot := (b_1, \ldots, b_n)$ be a sequence of sections of $\mathcal{O}_T$. Then there exist a $p$-completely faithfully flat morphism of crystalline prisms $u\colon \widetilde{T} \to T$ and a sequence $\tilde{b}_\cdot$ of sections of $\mathcal{O}_{\widetilde{T}}$ such that $u^\sharp(b_\cdot) = \phi(\tilde{b}_\cdot)$ and which enjoys the following universal property: if $u'\colon T' \to T$ is a morphism of crystalline prisms and $b'_\cdot$ a sequence of sections of $\mathcal{O}_{T'}$ such that $u'^\sharp(b_\cdot) = \phi(b'_\cdot)$, then there is a unique morphism of prisms $v\colon T' \to \widetilde{T}$ such that $u' = u \circ v$ and $v^\sharp(\tilde{b}_\cdot) = b'_\cdot$. Formation of $\widetilde{T}$ is functorial in $(T, b_\cdot)$ and compatible with base change.*

*Proof.* Because of the universality of the construction, it will suffice to prove this when $T = \operatorname{Spf} B$. Let $C$ be the $p$-adic completion of the polynomial algebra $B[x_1, \ldots, x_n]$, extend $\phi$ to a $\phi_B$-linear endomorphism $\psi$ of $C$ by letting $\psi(x_i) := b_i$, and let $I$ be the ideal of $C$ generated by $p$ and $\{x_i^p - b_i : i = 1, \ldots, n\}$. Since $\psi(x_i) = b_i \equiv x_i^p \pmod{I}$, this endomorphism induces the Frobenius endomorphism of $C/I$. Let $\tilde{T} := \operatorname{Spf}(\tilde{B} \to C)$ be the prismatic envelope of $I$ described in Theorem 2.3.8. Since $(x_1^p - b_1, \cdots, x_n^p - b_n)$ is a $C/pC$-regular sequence, statement (2) of the theorem tells us that $C/I \to \tilde{B}/p\tilde{B}$ is faithfully flat. Since $B/pB \to C/I$ is also faithfully flat, so is the map $B/pB \to \tilde{B}/p\tilde{B}$. Then $u^\sharp\colon B \to \tilde{B}$ is $p$-completely faithfully flat by Proposition 7.3.2. If $\tilde{b}_i$ is the image of $x_i$ in $\tilde{B}$, then $\phi(\tilde{b}_i)$ is the image of $\psi(x_i)$, which is indeed $u^\sharp(b_i)$.

To check that the universality of this construction, suppose that $u'\colon T' \to T$ is a morphism of crystalline prisms and $b'_\cdot$ is a sequence of sections of $\mathcal{O}_{T'}$. Then there is a unique $T$-morphism $T' \to \operatorname{Spf} C$ sending $x_i$ to $b'_i$. If $\phi(b'_i) = u'^\sharp(b_i)$, then this homomorphism sends $x_i^p - b_i$ to $b_i'^p - u'^\sharp(b_i) = b_i'^p - \phi(b'_i) \in p\mathcal{O}_{T'}$, so $T'$ becomes a $\operatorname{Spec}(C/I)$-prism over $C$ and the map $u'$ factors uniquely through $\tilde{T}$.

The universality of $\widetilde{T}$ guarantees its functoriality. To see that its formation commutes with base change, suppose that $f\colon T' \to T$ is a morphism of crystalline prisms and that $b'_{\cdot} = f^\sharp(b_{\cdot})$. Since $\widetilde{T} \to T$ is $p$-completely flat and $T'$ is $p$-torsion free, it follows from Proposition 7.3.3 that $\widetilde{T} \times_T T'$ is again $p$-torsion free and hence a crystalline prism. We claim that the map $\widetilde{T}' \to \widetilde{T} \times_T T'$ induced by $\tilde{f}$ is an isomorphism. But

$$\phi_{\tilde{T}\times_T T''}(\pi^\sharp_{\widetilde{T}}(\tilde{b}_{\cdot})) = \pi^\sharp_{\widetilde{T}}(\phi_{\widetilde{T}}(\tilde{b}_{\cdot})) = \pi^\sharp_{\widetilde{T}}(u^\sharp(b_{\cdot})) = f^\sharp(b_{\cdot}) = b'_{\cdot},$$

and so the universal property of $\widetilde{T}'$ guarantees the existence of the desired inverse morphism. □

**Remark 2.3.11.** The construction in Proposition 2.3.10 is universal but not especially efficient. For example, if $T = \operatorname{Spf} B$ and $b_{\cdot} = 0_{\cdot}$, then $\tilde{T}$ is the prismatic envelope of $(p, x_{\cdot}^p)$ in $B[x_{\cdot}]\hat{}$, which, as we shall see in Theorem 2.3.13, is the $p$-adically completed divided power envelope $B\langle x_{\cdot}\rangle\hat{}$, with $\phi(x_{\cdot}) = 0$. Note, however, that in this case the sequence $0_{\cdot}$ in $B$ satisfies $\phi(0_{\cdot}) = 0_{\cdot}$, so no $p$-completely flat cover is necessary. The universal property of $\tilde{T}$ gives us a section $T \to \tilde{T}$ of $\tilde{T} \to T$ sending $\tilde{x}_{\cdot}$ to $0_{\cdot}$.

**Corollary 2.3.12.** *Let $Y$ be a crystalline prism and $X \to \overline{Y}$ a closed immersion. If $(p, x_1, \ldots, x_r)$ generates $I_{X/Y}$, and $m \in \mathbf{N}$, then there exists a $p$-completely faithfully flat morphism of crystalline prisms $u\colon Y' \to Y$ and a sequence $(x'_1, \ldots, x'_r)$ of elements in $\mathcal{O}_{Y'}$ generating the ideal of $X' := u^{-1}(X)$ such that each $\delta(x'_i)$ belongs to $p^m\mathcal{O}_{Y'}$. In particular, $X'$ is $F_m$-alignable in $Y'$. Furthermore, if $(p, x_1, \ldots, x_r)$ is regular, then so is $(p, x'_1, \ldots, x'_r)$.*

*Proof.* We argue by induction $m$, the case $m = 0$ being trivial. Using the fact that the composition of $p$-completely flat morphisms is again $p$-completely flat and the induction hypothesis, we may suppose that $m > 0$ and that there exist $y_i \in \mathcal{O}_Y$ such that $\delta(x_i) = p^{m-1}y_i$ for each $i$. Proposition 2.3.10 shows that there is a $p$-completely flat morphism of crystalline prisms $u\colon Y' \to Y$ and a sequence of elements $y'_1, \ldots, y'_r \in \mathcal{O}_{Y'}$ such that $\phi(y'_i) = y_i$ for all $i$. Let $x'_i := u^*(x_i) - p^m y'_i$, and note that $(p, x'_1, \ldots, x'_r)$ generates $I_{X'/Y'}$. Then

$$\begin{aligned}
\phi(x'_i) &= \phi(u^*(x_i) - p^m y'_i) \\
&= u^*(x_i^p + p\delta(x_i)) - p^m\phi(y'_i) \\
&= u^*(x_i^p + p^m y_i) - p^m\phi(y'_i) \\
&= u^*(x_i^p) \\
&= (x'_i + p^m y'_i)^p \\
&\equiv x'^p_i \pmod{p^{m+1}}.
\end{aligned}$$

Thus $\delta(x'_i) \in p^m\mathcal{O}_{Y'}$. Now suppose that $(p, x_1, \ldots, x_r)$ is $\mathcal{O}_Y$-regular. Then the sequence $(\overline{x}_1, \ldots, \overline{x}_r)$ is $\mathcal{O}_{\overline{Y}}$-regular, and since $\overline{Y}' \to \overline{Y}$ is flat, so is its pullback to $\mathcal{O}_{\overline{Y}'}$. Since $\overline{x}'_i = u^*(\overline{x}_i)$ for every $i$, we see that the sequence $(\overline{x}'_1, \ldots, \overline{x}'_r)$ is regular. □

In some cases, prismatic envelopes can be expressed as PD-dilatations. The second of the following two results appears in [10, 2.38] and is the key to its main cohomology comparison theorems. This theorem will play a similarly fundamental role here.

**Theorem 2.3.13.** *Let $Y$ be a crystalline prism and let $i\colon X \to \overline{Y}$ be a regular immersion.*

*(1) If $\tilde{i}\colon \tilde{X} \to Y_2$ is a flat lifting of $X$ in $Y_2$, let $s_{\tilde{X}}\colon X \to \mathbb{D}_X(Y)$ be the resulting section as described in statement (4c) of Theorem 2.2.8. If $\tilde{X}$ is invariant under $F_Y$,[5] then the endomorphism $\mathbb{D}_X(F_Y)$ of $\mathbb{D}_X(Y)$ preserves $s_{\tilde{X}}$, and $\mathbb{\Delta}_X(Y)$ identifies with the PD-dilatation of this section. Moreover, $s_{\tilde{X}}$ is also a regular immersion.*

*(2) Let $X^\phi \subseteq Y$ be the inverse image of $X$ under $F_Y$. The canonical map $\overline{\mathbb{PD}}_X(Y) \to \overline{Y}$ factors through a map $z\colon \overline{\mathbb{PD}}_X(Y) \to X^\phi$, and $(\mathbb{PD}_X(Y), z, \pi)$ identifies with the prismatic envelope $\mathbb{\Delta}_{X^\phi}(Y)$ of $X^\phi$ in $Y$.*

*Proof.* The statement about the regularity of the immersion $s_{\tilde{X}}$ is explained in Proposition 2.2.14. Let $\mathbb{PD}_X(\mathbb{D}_X(Y))$ be the PD-dilatation of $s_{\tilde{X}}$ defined in Proposition 2.1.3, which is $p$-torsion free by construction and which agrees with the usual PD-envelope if $s_{\tilde{X}}$ is a regular immersion. The morphism $\overline{\mathbb{PD}}_X(\mathbb{D}_X(Y)) \to Y$ factors through $\overline{\mathbb{D}}_X(Y)$ and hence through $X$, and so $\mathbb{PD}_X(\mathbb{D}_X(Y))$ can also be viewed as a $p$-adic enlargement of $X$ over $Y$. We shall show that it in fact a crystalline prism over $X$, compatible with the prism structure of $Y$.

The endomorphism $F_Y$ of $Y$ induces an endomorphism $\mathbb{D}_X(F_Y)$ of $\mathbb{D}_X(Y)$ which may not be a Frobenius lift. If $\tilde{X}$ is invariant under $F_Y$, then $\mathbb{D}_X(F_Y)$ preserves $s_{\tilde{X}}$ by functoriality, and we shall show that the resulting endomorphism of its PD-dilatation is in fact a Frobenius lift. Let us check this locally, assuming that $\mathbb{D}_X(Y) \to Y$ is given by $\theta\colon C \to C'$ and that $C' \to C''$ defines the PD-dilatation of $s_{\tilde{X}}$. Our claim is that $\phi(c'') \equiv c''^p \pmod{pC''}$ for all $c'' \in C''$. Since we know that this holds for $c \in C$, it will suffice to check this congruence for a set of topological generators for $C''$ as a $C$-algebra. The $C$-algebra $C'$ is topologically generated by $\rho(I_{\tilde{X}/Y})$, and $C''$ is topologically generated by the divided powers of the elements of $\rho(I_{\tilde{X}/Y})$. Thus it will suffice to show that $\phi(z) \equiv z^p \pmod{pC''}$ for all such $z$. These elements belong to the PD-ideal of $s_{\tilde{X}}$ in $C''$, hence their $p$th powers are automatically divisible by $p$, and it will suffice to show that the same is true of their images under $\phi$.

For the sake of further remarks, let us assume that $\tilde{X}$ is in fact flat over $\mathbf{Z}/p^r\mathbf{Z}$, with $r \geq 2$ and invariant under $F_Y$. If $x \in I_{\tilde{X}/Y}$, then

$$\phi(x) = x^p + p\delta(x) \in I_{\tilde{X}/Y}.$$

It follows that $p\delta(x) \in I_{\tilde{X}/Y}$, and since $\tilde{X}$ is flat over $\mathbf{Z}/p^r\mathbf{Z}$, that $\delta(x) \in$

[5] Recall from Corollary 1.7.9 that if $X$ is reduced, a such a lifting $\tilde{X}$ is in fact unique.

$I_{\tilde{X}/Y} + p^{r-1}C$, say $\delta(x) = x' + p^{r-1}y$, with $x' \in I_{\tilde{X}/Y}$ and $\in C$. Then:

$$\begin{array}{rcl} \phi(\rho(x)) & = & \rho(\phi(x)) \\ & = & \rho(x^p + px' + p^r y) \\ & = & p^{p-1}(\rho(x))^p + p\rho(x') + p^{r-1}y, \\ \phi(\rho(x)) & = & pz, \end{array} \tag{2.3.4}$$

where $z := p^{p-2}\rho(x)^p + \rho(x') + p^{r-2}y$. Then

$$\phi((\rho(x))^{[n]}) = z^n p^{[n]} \in pC'',$$

as required. Let us remark that, if $r \geq 3$, then $z$ belongs to the PD-ideal of $C''$ and we can write the stronger:

$$\phi((\rho(x))^{[n]} = p^n z^{[n]} \in p^n C''.$$

To complete the proof of statement (1), it remains only to show that the $X$-prism $\mathbb{PD}_X(\mathbb{D}_X(Y))$ is universal. Suppose that $T$ is another $X$-prism over $Y$. In particular it is a $p$-adic enlargement of $X$ over $Y$, and hence maps to $\mathbb{D}_X(Y)$. Equation (2.3.4) above shows that $\mathbb{D}_X(F_Y)$ maps the ideal $I_{X/\mathbb{D}}$ of $s_{\tilde{X}}$ to $p\mathcal{O}_{\mathbb{D}_X(Y)}$, and hence $\phi_T$ maps the ideal $I_{X/\mathbb{D}}\mathcal{O}_T$ to $p\mathcal{O}_T$. Thus $I_{X/\mathbb{D}}\mathcal{O}_T$ is contained in the ideal of $F_{\overline{T}}(\overline{T}) \subseteq T$. As we saw in Proposition 2.1.4, this is a PD-ideal. It follows that $T$ factors uniquely through $\mathbb{PD}_X(\mathbb{D}_X(Y))$.

To prove (2), note first that it follows from (3) of Proposition 2.1.3 that the endomorphism $\mathbb{PD}_X(F_Y)$ endows $\mathbb{PD}_X(Y)$ with the structure of a crystalline prism. Furthermore, the ideal of $X^\phi$ is generated by $p$ and the $p$th powers of elements of $I_{X/Y}$, and since $I_{X/Y}$ maps into a PD-ideal in $\mathbb{PD}_X(Y)$, all these elements become divisible by $p$ in $\mathcal{O}_{\mathbb{PD}_X(Y)}$. Thus $\mathbb{PD}_X(Y)$ defines an $X^\phi$-prism over $Y$. To see that it is universal, let us first observe that the scheme-theoretic image $F_{X^\phi}(X^\phi)$ is contained in $X$. Indeed, $I_{X^\phi/\overline{Y}}$ is generated by $\{x^p : x \in I_{X/\overline{Y}}\}$, and $I_{F(X^\phi)/\overline{Y}} = \{f \in \mathcal{O}_{\overline{Y}} : f^p \in I_{X^\phi/\overline{Y}}\}$ which certainly contains $I_{X/\overline{Y}}$. Thus if $(T, z_T, \pi_T)$ is an $X^\phi$-prism over $Y$ there is a commutative diagram:

$$\begin{array}{ccccc} \overline{T} & \longrightarrow & F_T(\overline{T}) & & \\ \downarrow{\scriptstyle z_T} & & \downarrow & \searrow{\scriptstyle z'_T} & \\ X^\phi & \longrightarrow & F_{X^\phi}(X^\phi) & \longrightarrow & X \end{array}$$

where the second vertical arrow comes by functoriality. But Proposition 2.1.4 tells us that $\phi(\overline{T}) = T_{\mathbb{PD}}$, and, since $T_{\mathbb{PD}} \to T$ is a PD-immersion, the maps $z'_T$ and $\pi_T$ then give $T$ the structure of a PD-enlargement of $X$ over $Y$ which must factor through $\mathbb{PD}_X(Y)$. □

**Example 2.3.14** Let us remark that, as observed in [10, 2.40], statement (2) of Theorem 2.3.13 implies that $\mathbb{\Delta}_{X^\phi}(Y)$ does not depend on the Frobenius lifting of $Y$, but this is not true more generally. By way of example, it is shown that if $Y = \operatorname{Spf} W[x]\hat{}$ and $X$ is defined by the ideal $(x, p)$, then the choices $\phi(x) = x^p$ and $\phi(x) = x^p + p$ give different prismatic envelopes. Statement (1) of Theorem 2.3.13 allows us to be a little more explicit about these examples, since $X$ is $F_2$-alignable in both cases. In the first case, since $\phi(x) = x^p$, one finds that $\mathbb{\Delta}_X(Y)$ is the formal spectrum of the completed PD-polynomial algebra $W\langle x/p\rangle\hat{}$. In the second case, note that $\phi(x-p) = x^p \equiv (x-p)^p \pmod{p^2}$. Thus $(x-p, p^2)$ defines an $F_2$-aligned lifting, so $\mathbb{\Delta}_X(Y)$ is the formal spectrum of $W\langle(x/p-1)\rangle\hat{}$.

For another example, consider the crystalline prism $Y$ given by $\operatorname{Spf} W[x,y]\hat{}$, with $\phi(x) = x^p + py$ and $\phi(y) = y^p$, and let $X$ be the closed subscheme defined by the ideal $(p, x)$. As we saw in Example 1.7.7, this $X$ is not $F_2$-alignable in $Y$. We claim that the prismatic envelope $\mathbb{\Delta}_X(Y) \to Y$ is given by the homomorphism:

$$W[x,y]\hat{} \to W[s]\langle t\rangle\hat{} : x \mapsto ps, y \mapsto s^p + pt.$$

Let us check this using Theorem 2.3.3, which tells us that the prismatic envelope is obtained by adjoining variables $t_j : j \in \mathbf{N}$, subject to the relations

$$pt_{j+1} = \delta^j(x) - t_j^p, \quad t_0 = 0.$$

In our case, $\delta(x) = y$ and $\delta^j(x) = 0$ for $j > 1$, so we have $t_1 = s$, $t_2 = t$, and $pt_{j+1} = -t_j^p$ for $j \geq 2$. This amounts to taking the divided power envelope of $(t)$, as claimed.

**Corollary 2.3.15.** *Let $j\colon Y \to Z$ be a closed immersion of crystalline prisms, where $\overline{Y} \to \overline{Z}$ is a regular immersion.*

*(1) There are commutative diagrams:*

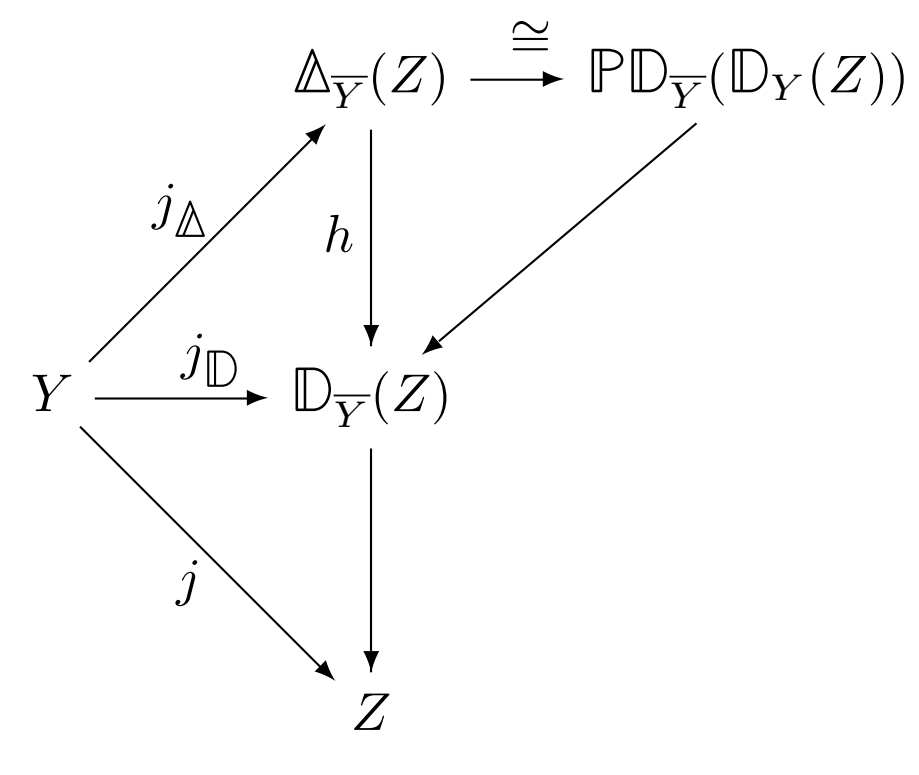

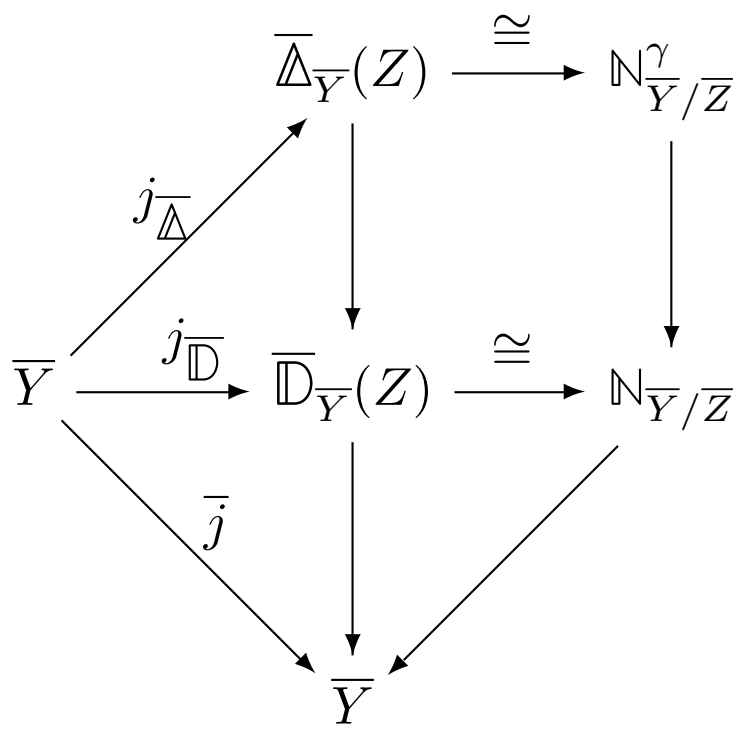

*In the second diagram, the maps $\overline{Y} \to \mathbb{N}^{\gamma}_{\overline{Y}/\overline{Z}}$ and $\overline{Y} \to \mathbb{N}_{\overline{Y}/\overline{Z}}$ are the respective zero-sections.*

(2) *In particular, if $Y/S$ is a $p$-completely smooth morphism of crystalline prisms, let $\Delta_Y(1)$ denote the prismatic envelope of the closed immersion $\overline{Y} \to Y \to Y(1)$. Then $\overline{\Delta}_Y(1)$ identifies with $\mathbb{T}^{\gamma}_{\overline{Y}/\overline{S}}$ and*

$$\mathcal{O}_{\overline{\Delta}_Y(1)} \cong \Gamma^{\bullet}(\Omega^1_{\overline{Y}/\overline{S}}).$$

*Proof.* The universal properties of $\mathbb{D}_Y(Z)$ and $\Delta_Y(Z)$ imply the existence of sections $j_{\mathbb{D}} \colon Y \to \mathbb{D}_{\overline{Y}}(Z)$ and $j_{\Delta} \colon Y \to \Delta_{\overline{Y}}(Z)$, and Theorem 2.3.13 implies that $\Delta_{\overline{Y}}(Z)$ identifies with the PD-dilatation of $j_{\mathbb{D}}$. Since the immersion $\overline{Y} \to \overline{Z}$ is regular, it follows from (4c) of Theorem 2.2.8 that $\overline{\mathbb{D}}_{\overline{Y}}(Z)$ admits the structure of an $\mathbb{N}_{\overline{Y}/\overline{Z}}$-torsor, and this torsor is trivialized by the section $j_{\overline{\mathbb{D}}}$. This proves statement (1), and statement (2) follows, with $Y(1)$ playing the role of $Z$. $\square$

In the situation of Corollary 2.3.15, the ideal $I_{Y/\Delta_Y(Z)}$ is a PD-ideal, and we shall need some information about the somewhat subtle relationship between its divided power structure and the $\delta$-structure. The following lemma contains the main computations, which we place in a more general context.

**Lemma 2.3.16.** *Let $Y$ be a crystalline prism, let $X \to \overline{Y}$ be a regular immersion, and let $x$ be a local section of $I_{X/Y}$. Suppose that $\delta(x) \in I_{X/Y}$, so that $\rho_{\Delta}(\delta(x)) \in \mathcal{O}_{\Delta_X(Y)}$, where $\rho_{\Delta} \colon I_{X/Y} \to \mathcal{O}_{\Delta_X(Y)}$ is the map for $\Delta_X(Y)$. Then $(\rho_{\Delta}(x))^{[n]} \in \mathcal{O}_{\Delta_X(Y)}$ for all $n$, and the following statements are verified.*

$$\begin{aligned}
\overline{\delta}_{\Delta}(\rho_{\Delta}(x)) &= (\overline{\rho}_{\Delta}(x))^{[p]} + \overline{\rho}_{\Delta}(\delta(x)); \\
\overline{\delta}_{\Delta}((\rho_{\Delta}(x))^{[p^m]}) &= (\overline{\rho}_{\Delta}(x))^{[p^{m+1}]} \quad \text{if } m > 0 \text{ and } p \text{ is odd;} \\
&= (\overline{\rho}_{\Delta}(x))^{[p^{m+1}]} + (\overline{\rho}_{\Delta}(\delta(x)))^{p^m} \quad \text{if } m > 0 \text{ and } p = 2\ , \\
\overline{\delta}_{\Delta}((\rho_{\Delta}(x))^{[n]}) &= 0 \quad \text{if } n \text{ is not a power of } p\ .
\end{aligned}$$

*Proof.* We start by repeating the calculations in the proof of statement (1) of Theorem 2.3.13.

$$\phi_{\Delta}(\rho_{\Delta}(x)) = p^{p-1}(\rho_{\Delta}(x))^p + \pi^*_{\Delta}(\delta(x))$$

$$= \quad p^{p-1}(\rho_{\mathbb{A}}(x))^p + p\rho(_{\mathbb{A}}\delta(x)).$$

This implies that $\phi_{\mathbb{A}}(\rho(x)) \in p\mathcal{O}_{\mathbb{A}_{Y(Z)}}$, and hence by Proposition 2.1.4, that $\rho(x)^{[n]} \in \mathcal{O}_{\mathbb{A}_{Y(Z)}}$ for all $n$. We continue as follows:

$$\begin{aligned}
p\delta_{\mathbb{A}}(\rho_{\mathbb{A}}(x)) &= p^{p-1}(\rho_{\mathbb{A}}(x))^p + p\rho_{\mathbb{A}}(\delta(x)) - (\rho_{\mathbb{A}}(x))^p \\
\delta_{\mathbb{A}}(\rho_{\mathbb{A}}(x)) &= p^{p-1}(p-1)!(\rho_{\mathbb{A}}(x))^{[p]} + \rho_{\mathbb{A}}(\delta(x)) - (p-1)!(\rho_{\mathbb{A}}(x))^{[p]} \\
\overline{\delta}_{\mathbb{A}}(\rho_{\mathbb{A}}(x)) &= \overline{\rho}_{\mathbb{A}}(\delta(x)) + (\overline{\rho}_{\mathbb{A}}(x))^{[p]}.
\end{aligned}$$

This proves this first equation. The remaining cases will require some more computation. Let $\tilde{x} := \rho_{\mathbb{A}}(x)$, $\tilde{\delta}(x) := \rho_{\mathbb{A}}(\delta(x))$, and

$$z := p^{p-2}\tilde{x}^p + \tilde{\delta}(x) = p^{p-1}(p-1)!\tilde{x}^{[p]} + \tilde{\delta}(x).$$

Then we can rewrite the result of the above calculation as:

$$\begin{aligned}
\phi_{\mathbb{A}}(\tilde{x}) &= p^{p-1}\tilde{x}^p + p\tilde{\delta}(x) \\
&= pz \\
\phi_{\mathbb{A}}(\tilde{x}^{[n]}) &= p^{[n]}z^n \\
\delta_{\mathbb{A}}(\tilde{x}^{[n]}) &= p^{-1}\left(p^{[n]}z^n - (\tilde{x}^{[n]})^p\right).
\end{aligned} \tag{2.3.5}$$

If $p$ is odd and $n > 1$, then $p^{[n]}$ is divisible by $p^2$, so

$$\delta_{\mathbb{A}}(\tilde{x}^{[n]}) \equiv -p^{-1}(\tilde{x}^{[n]})^p \pmod{p}.$$

Furthermore,

$$p^{-1}(\tilde{x}^{[n]})^p = p^{-1}(\frac{\tilde{x}^n}{n!})^p = \frac{\tilde{x}^{pn}}{p(n!)^p} = \frac{(pn)!}{p(n!)^p}\tilde{x}^{[pn]}.$$

Statement (3) of Lemma 1.6.7 shows that the $p$-adic ordinal of the coefficient of $\tilde{x}^{[pn]}$ in this expression is $\sigma(n) - 1$, where $\sigma(n)$ is the sum of the digits in the $p$-adic expansion of $n$. Thus this coefficient is divisible by $p$ unless $n$ is a power of $p$, in which case it is a unit, which we claim is congruent to $-1$ mod $p$. We use statement (1c) of Lemma 1.6.7 to compute:

$$\begin{aligned}
U\left(\frac{(pp^m)!}{p(p^m!)^p}\right) &= U(p^{m+1}!)(U(p^m!))^{-p} \\
&= (-1)^{m+1}(-1)^{-mp} \\
&= -1.
\end{aligned}$$

Thus

$$\delta_{\mathbb{A}}(\tilde{x}^{[p^m]}) \equiv -\frac{(pp^m)!}{p(p^m!)^p}\tilde{x}^{[pp^m]} \equiv \tilde{x}^{[p^{m+1}]} \pmod{p},$$

as claimed.

Now suppose $p = 2$. If $n > 1$ and is not a power of 2, then $2^{[n]}$ is divisible by 4, Lemma 1.6.7 shows that $(\tilde{x}^{[n]})^2$ is divisible by 4, and we again see that $\delta_{\Delta}(\tilde{x}^{[n]})$ is divisible by 2. Now suppose that $n = 2^m$ with $m \geq 1$. Then $\mathrm{ord}_2(2^{[n]}) = 1$, and equation (2.3.5) tells us that:

$$\begin{aligned} \delta_{\Delta}(\tilde{x}^{[n]}) &= 2^{-1}\left(2^{[n]}z^n - (\tilde{x}^{[n]})^2\right) \\ &= az^n + b\tilde{x}^{[2n]} \end{aligned}$$

where $a, b \in \mathbf{Z}_2$ are 2-adic units. Furthermore, $z \equiv \tilde{\delta}_{\Delta}(x) \pmod 2$, so

$$\delta_{\Delta}(\tilde{x}^{[2^m]}) \equiv \overline{\rho}_{\Delta}(\delta(x))^{2^m} + \tilde{x}^{[2^{m+1}]} \pmod 2.$$

This completes the proof. □

**Proposition 2.3.17.** *Let $Y \to Z$ be a closed immersion of crystalline prisms, such that $\overline{Y} \to \overline{Z}$ is a regular immersion. Since $Y \to Z$ is a morphism of prisms, the ideal $I_{Y/\Delta}$ of the canonical section $Y \to \Delta_Y(Z)$ is a PD-ideal and is invariant under $\delta$.*

*(1) The filtrations on $\mathcal{O}_{\overline{\Delta}_{Y(Z)}}$ induced by the PD-filtration and the $\delta$-filtration $\mathrm{Fil}^\bullet_\delta$ (Definition 1.7.14) on $\mathcal{O}_{\Delta_{Y(Z)}}$ satisfy the following relations.*

*(a) $I^2_{\overline{Y}/\overline{\Delta}} + I^{[p+1]}_{\overline{Y}/\overline{\Delta}} \subseteq \mathrm{Fil}^2_\delta \mathcal{O}_{\overline{\Delta}_{Y(Z)}}$.*

*(b) $(\overline{\rho}_{\Delta}(x))^{[p]} + \overline{\rho}_{\overline{\Delta}}(\delta_{\overline{Y}/\overline{Z}}(x)) \in \mathrm{Fil}^p_\delta \mathcal{O}_{\overline{\Delta}_{Y(Z)}}$ for every $x \in I_{Y/Z}$.*

*(c) $\delta_{\Delta}(I^{[p]}_{Y/\Delta}) \subseteq I^{[p^2]}_{Y/\Delta} + p\mathcal{O}_{\Delta_{Y(Z)}}$.*

*(2) There is a commutative diagram:*

$$\begin{array}{ccc} I_{Y/Z} & \xrightarrow{\rho_{\Delta}} & I_{Y/\Delta} \\ \downarrow{\scriptstyle pr_{\mathcal{C}}} & & \downarrow{\scriptstyle pr_\delta} \\ \mathcal{C}_{\overline{Y}/\overline{Z}} & \xrightarrow{\rho_\delta} & \mathrm{Gr}^1_{\mathrm{Fil}_\delta} \mathcal{O}_{\overline{\Delta}_{Y(Z)}} \end{array} \tag{2.3.6}$$

*in which $\rho_\delta$ is an isomorphism and the vertical arrows are the natural projections.*

*(3) If $x \in I_{Y/Z}$, then*

$$\rho_\delta(pr_{\mathcal{C}}(\delta(x))) = pr_\delta(\rho_{\Delta}(\delta(x))) = -pr_\delta((\rho_{\Delta}(x))^{[p]}).$$

*Proof.* Throughout the proof we write $\Delta$ for $\Delta_{Y(Z)}$. It is clear from the definition that $I^2_{Y/\Delta} \subseteq \mathrm{Fil}^2_\delta \mathcal{O}_{\Delta}$. The ideal $I_{Y/\Delta}$ is generated as a PD-ideal by elements of the form $\rho(x)$, with $x$ in the ideal $I_{Y/Z}$, which is invariant under $\delta_Z$. Lemma 2.3.16 applies to such elements and shows that $\overline{\rho}_{\mathbb{\Delta}}(x)^{[p^m]} \in \mathrm{Fil}^p_\delta \mathcal{O}_{\overline{\mathbb{\Delta}}_{Y/Z}}$ if $m \geq 2$, and hence that $I^{[p^2]}_{\overline{Y}/\mathbb{\Delta}} \subseteq \mathrm{Fil}^2_\delta \mathcal{O}_{\overline{\mathbb{\Delta}}}$. Thus

$$I^{[p+1]}_{\overline{Y}/\overline{\mathbb{\Delta}}} = I^{[p]}_{\overline{Y}/\overline{\mathbb{\Delta}}} I_{\overline{Y}/\overline{\mathbb{\Delta}}} + I^{[p^2]}_{\overline{Y}/\overline{\mathbb{\Delta}}} \subseteq \mathrm{Fil}^2_\delta \mathcal{O}_{\overline{\mathbb{\Delta}}}.$$

This proves (1a). The first equation in Lemma 2.3.16 says

$$\overline{\rho}_{\mathbb{\Delta}}(x)^{[p]} + \overline{\rho}_{\mathbb{\Delta}}(\delta(x)) = \overline{\delta}_{\mathbb{\Delta}}(\rho_{\mathbb{\Delta}}(x))$$

and this implies (1b). Lemma 2.3.16 also tells us that, if $m > 0$ and $\tilde{x} := \rho_{\mathbb{\Delta}}(x)$, then

$$\delta_{\mathbb{\Delta}}(\tilde{x}^{[p^m]}) \equiv \tilde{x}^{[p^{m+1}]} \pmod{p\mathcal{O}_{\mathbb{\Delta}}}.$$

and hence that $\delta_{\mathbb{\Delta}}(\tilde{x}^{[p^m]}) \in I^{[p^2]}_{Y/\mathbb{\Delta}} + p\mathcal{O}_{\mathbb{\Delta}}$. It also tells us that, if $n$ is not a power of $p$, then $\delta_{\mathbb{\Delta}}(\tilde{x}^{[n]}) \in p\mathcal{O}_{\mathbb{\Delta}}$. If $y$ is another element of $I_{Y/\mathbb{\Delta}}$ and $\tilde{y} := \rho_{\mathbb{\Delta}}(y)$, then

$$\delta_{\mathbb{\Delta}}(\tilde{x}\tilde{y}) = \tilde{x}^p \delta_{\mathbb{\Delta}}(\tilde{y}) + \tilde{y}^p \delta_{\mathbb{\Delta}}(\tilde{x}) + p\delta_{\mathbb{\Delta}}(\tilde{x})\delta_{\mathbb{\Delta}}(\tilde{y}) \in p\mathcal{O}_{\mathbb{\Delta}},$$

so the claim (1c) holds for every monomial in $I^{[p]}_{Y/\mathbb{\Delta}}$. Finally, if $m$ and $n$ are at least $p$, then

$$\delta_{\mathbb{\Delta}}(\tilde{x}^{[m]} + \tilde{y}^{[n]}) = \delta_{\mathbb{\Delta}}(\tilde{x}^{[m]}) + \delta_{\mathbb{\Delta}}(\tilde{y}^{[n]}) - C(\tilde{x}^{[m]}, \tilde{y}^{[n]}),$$

which belongs to $I^{[p^2]}_{Y/\mathbb{\Delta}}$ because $C(x, y)$ is homogenous of degree $p$. These statements together imply (1c).

Let us note now that statement (3) will follow from (1b) and (2), once we have proved the latter. Indeed, if $x \in I_{Y/Z}$, then $\rho_\delta(pr_{\mathcal{C}}(x)) = pr_\delta(\rho(\delta(x)))$ because of the commutativity of the diagram, and $pr_\delta(\rho(\delta(x))) = -pr_\delta((\rho_{\mathbb{\Delta}}(x))^{[p]})$ because $(\overline{\rho}_{\mathbb{\Delta}}(x))^{[p]} + \overline{\rho}_{\overline{\mathbb{\Delta}}}(\delta_{\overline{Y}/\overline{Z}}(x))$ belongs to $\mathrm{Fil}^2_{\mathrm{Fil}_\delta} \mathcal{O}_{\overline{\mathbb{\Delta}}}$.

To produce the isomorphism $\rho_\delta$, we first claim that there is a commutative diagram in which the horizontal and vertical rows are short exact and the slanted arrows are isomorphisms:

$$\begin{array}{ccccc} & & I_{\overline{Y}/\overline{\mathbb{\Delta}}}/(I_{\overline{Y}/\overline{\mathbb{D}}}\mathcal{O}_{\overline{\mathbb{\Delta}}} + I^{[p+1]}_{\overline{Y}/\overline{\mathbb{\Delta}}}) & & \\ & \nearrow^{\tilde{\gamma}_\rho} & \uparrow b & & \\ F^*_{\overline{Y}}(\mathcal{C}_{\overline{Y}/\overline{Z}}) & \xrightarrow{\tilde{\overline{\delta}}_{Y/\mathbb{\Delta}}} & I_{\overline{Y}/\overline{\mathbb{\Delta}}}/(I^2_{\overline{Y}/\overline{\mathbb{\Delta}}} + I^{[p+1]}_{\overline{Y}/\overline{\mathbb{\Delta}}}) & \xrightarrow{\tilde{\pi}_{Y/\mathbb{\Delta}}} & \mathrm{Gr}^1_{\mathrm{Fil}_\delta} \mathcal{O}_{\overline{\mathbb{\Delta}}} \\ & & \uparrow a & \nearrow_{\rho_\delta} & \\ & & \mathcal{C}_{\overline{Y}/\overline{Z}}\,. & & \end{array} \tag{2.3.7}$$

We construct this as follows. We have a diagram:

$$
\begin{array}{ccccc}
\mathcal{C}_{\overline{Y}/\overline{Z}} & \xrightarrow{\overline{\rho}_{\mathbb{D}}} & I_{\overline{Y}/\overline{\mathbb{D}}}/I^2_{\overline{Y}/\overline{\mathbb{D}}} & \longrightarrow & I_{\overline{Y}/\overline{\mathbb{\Delta}}}/(I^2_{\overline{Y}/\overline{\mathbb{\Delta}}} + I^{[p+1]}_{\overline{Y}/\overline{\mathbb{\Delta}}}) \\
& \overline{a} \searrow & & \searrow & \downarrow c \\
& & & & I_{\overline{Y}/\overline{\mathbb{\Delta}}}/I^{[2]}_{\overline{Y}/\overline{\mathbb{\Delta}}},
\end{array}
\tag{2.3.8}
$$

where $\overline{\rho}_{\mathbb{D}}$ is the mod $p$ reduction of the isomorphism described in (3) of Proposition 2.2.14, the other arrows are the obvious ones, and the composed horizontal map is the map $a$ in the diagram (2.3.7). We claim that the two slanted maps in the diagram (2.3.8) are isomorphisms, which we check locally, with the aid of a sequence $(x_1, \ldots, x_r)$ lifting a regular sequence generating $I_{\overline{Y}/\overline{Z}}$. Let $s_i := \rho_{\mathbb{D}}(x_i)$, and let $\overline{x}_i$ (resp. $\overline{s}_i$) denote their respective images in the respective conormal sheaves. Then Corollary 2.3.15 tells us that $\mathcal{O}_{\overline{\mathbb{D}}}$ (resp. $\mathcal{O}_{\overline{\mathbb{\Delta}}}$) identifies with the polynomial algebra (resp. the PD polynomial algebra) on $(\overline{s}_1, \ldots, \overline{s}_r)$, and the claim follows. We conclude that $\overline{a}^{-1} \circ c$ splits the map $a$ in diagram (2.3.7). We see also that $(\overline{s}_1, \ldots, \overline{s}_r, \overline{s}_1^{[p]}, \ldots, \overline{s}_r^{[p]})$ induces a basis for $I_{\overline{Y}/\overline{\mathbb{\Delta}}}/(I^2_{\overline{Y}/\overline{\mathbb{\Delta}}} + I^{[p+1]}_{\overline{Y}/\overline{\mathbb{\Delta}}})$.

The kernel of $c$ is

$$
\begin{aligned}
I^{[2]}_{\overline{Y}/\overline{\mathbb{\Delta}}}/(I^2_{\overline{Y}/\overline{\mathbb{\Delta}}} + I^{[p+1]}_{\overline{Y}/\overline{\mathbb{\Delta}}}) &= (I^2_{\overline{Y}/\overline{\mathbb{\Delta}}} + I^{[p]}_{\overline{Y}/\overline{\mathbb{\Delta}}})/(I^2_{\overline{Y}/\overline{\mathbb{\Delta}}} + I^{[p+1]}_{\overline{Y}/\overline{\mathbb{\Delta}}}) \\
&\cong I^{[p]}_{\overline{Y}/\overline{\mathbb{\Delta}}}/I^{[p+1]}_{\overline{Y}/\overline{\mathbb{\Delta}}} \cong I_{\overline{Y}/\overline{\mathbb{\Delta}}}/(I_{\overline{Y}/\overline{\mathbb{D}}}\mathcal{O}_{\overline{\mathbb{\Delta}}} + I^{[p+1]}_{\overline{Y}/\overline{\mathbb{\Delta}}}),
\end{aligned}
$$

since $I_{Y/\mathbb{\Delta}} = I_{Y/\mathbb{D}}\mathcal{O}_{\mathbb{\Delta}} + I^{[p]}_{Y/\mathbb{\Delta}}$. Thus $(\overline{s}_1^{[p]}, \ldots, \overline{s}_r^{[p]})$ induces a basis for this kernel, which identifies with the target of $b$. Thus the vertical column of diagram (2.3.7) is exact. Note also that the composition of $\overline{\rho}_{\mathbb{D}}$ with the $p$th divided power map induces an isomorphism

$$
\tilde{\gamma}_\rho \colon F^*_{\overline{Y}}(\mathcal{C}_{\overline{Y}/\overline{Z}}) \to I_{\overline{Y}/\overline{\mathbb{\Delta}}}/(I_{\overline{Y}/\overline{\mathbb{D}}}\mathcal{O}_{\overline{\mathbb{\Delta}}} + I^{[p+1]}_{\overline{Y}/\overline{\mathbb{\Delta}}}),
$$

sending $1 \otimes \overline{x}_i$ to $\overline{s}_i^{[p]}$.

Statement (4) of Proposition 1.7.13 implies that $\delta_{Y/\mathbb{\Delta}}$ induces a $p$-linear endomorphism of $I_{\overline{Y}/\overline{\mathbb{\Delta}}}/I^2_{\overline{Y}/\overline{\mathbb{\Delta}}}$, and we see from (1c) of the current proposition that this endomorphism takes the image of $I^{[p]}_{\overline{Y}/\overline{\mathbb{\Delta}}}$ to the image of $I^{[p+1]}_{\overline{Y}/\overline{\mathbb{\Delta}}}$. Since $I^{[2]}_{\overline{Y}/\overline{\mathbb{\Delta}}} = I^2_{\overline{Y}/\overline{\mathbb{\Delta}}} + I^{[p]}_{\overline{Y}/\overline{\mathbb{\Delta}}}$, we find a $p$-linear map

$$
\mathcal{C}_{\overline{Y}/\overline{Z}} \cong I_{\overline{Y}/\overline{\mathbb{\Delta}}}/I^{[2]}_{\overline{Y}/\overline{\mathbb{\Delta}}} \longrightarrow I_{\overline{Y}/\overline{\mathbb{\Delta}}}/(I^2_{\overline{Y}/\overline{\mathbb{\Delta}}} + I^{[p+1]}_{\overline{Y}/\overline{\mathbb{\Delta}}})
$$

hence the map $\tilde{\delta}_{Y/\overline{\mathbb{\Delta}}}$ in the diagram (2.3.7). The surjective map $\tilde{\pi}_{Y/\mathbb{\Delta}}$ exists because of statement (1a), and the horizontal row is exact in the middle because $\mathrm{Fil}^2_\delta$ is generated by $I^2_{\overline{Y}/\overline{\mathbb{\Delta}}}$ together with the image of $\delta_{\mathbb{\Delta}}$. Recall from Lemma 2.3.16 that $\overline{\delta}_{\mathbb{\Delta}}(\rho_{\mathbb{\Delta}}(x)) = (\overline{\rho}_{\mathbb{\Delta}}(x))^{[p]} + \overline{\rho}_{\overline{\mathbb{\Delta}}}(\delta_{\overline{Y}/\overline{Z}}(x))$. But $\delta_{Y/Z}(x) \in I_{Y/Z}$, so $\rho_{\mathbb{D}}(\delta_{Y/Z}(x)) \in I_{Y/\mathbb{D}}$, and thus maps to zero under $b$. This implies that the commutativity of the upper left triangle in the diagram. It follows from this that the first map in the horizontal sequence is injective and that the map $\rho_\delta$ is an isomorphism. □

## 2.4 Comparing envelopes

The various envelopes we have constructed can be thought of as having differing sizes and shapes. If $Y$ is a crystalline prism and $X$ is a closed subscheme of $\overline{Y}$, the following diagram illustrates some of their relations:

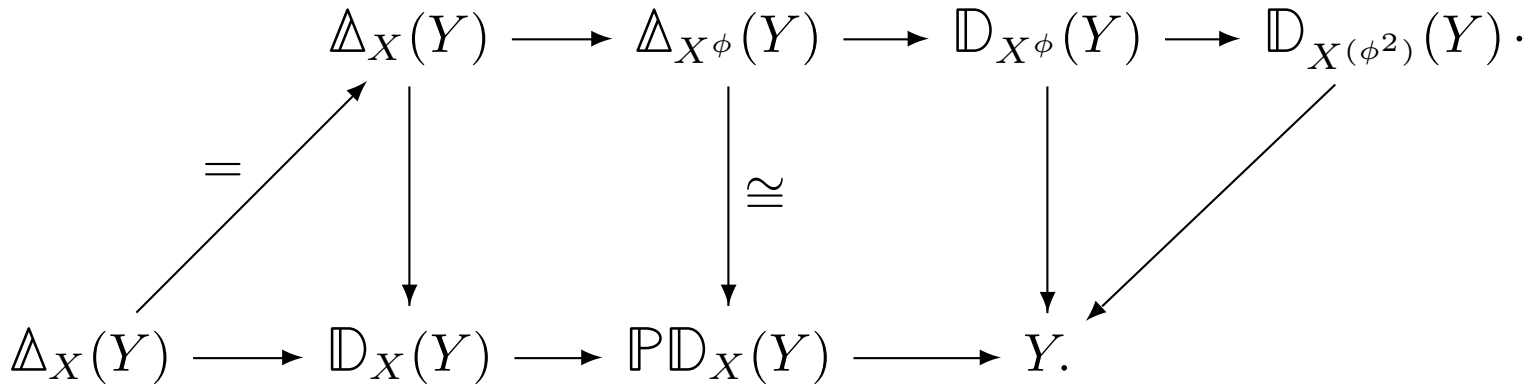

(2.4.1)

Here the central vertical isomorphism comes from statement (2) of Theorem 2.3.13.

A crude measure of the size of these envelopes is given by the radius of the analytic tubes they define. It may be interesting to compute these. For a simple example, let $Y$ be the formal spectrum of $W[x]\hat{}$, with $\phi(x) := x^p$ and let $i\colon X \to Y$ be defined by $(p, x)$. Let $\widehat{\overline{K}}$ be the completion of the algebraic closure of the fraction field of $W$. If $f = \sum a_i x^i \in W[x]\hat{}$ and $z \in \widehat{\overline{K}}$, then $\sum a_i z^i$ is guaranteed to converge if $|z| \leq 1$, so the rigid space $Y_K$ associated to $W[x]\hat{}$ corresponds to the closed unit disc $\{z : |z| \leq 1\}$. Now $\mathbb{D}_X(Y) = \operatorname{Spf} W[t]\hat{}$, with $x = pt$, and the image of $\mathbb{D}_X(Y)_{\widehat{\overline{K}}}$ in $Y_{\widehat{\overline{K}}}$ is the closed disc of radius $p^{-1}$. On the other hand, an element of the completed divided power algebra $W\langle x\rangle\hat{}$ converges at $z$ if $\operatorname{ord}_p(z^{p^n}/p^n!) \geq 0$ for all $n$. But

$$\begin{aligned}\operatorname{ord}_p(z^{p^n}/p^n!) &= p^n \operatorname{ord}_p(z) - \operatorname{ord}_p(p^n!)\\ &= p^n \operatorname{ord}_p(z) - (1 + p + p^2 + \cdots + p^{n-1})\\ &= p^n(\operatorname{ord}_p(z) - (p^{-n} + \cdots + p^{-1})).\end{aligned}$$

Thus $\operatorname{ord}_p(z^{p^n}/p^n!) \geq 0$ for all $n$ if and only if $\operatorname{ord}_p(z) \geq 1/(p-1)$, so $\mathbb{PD}_X(Y)_K$ corresponds to the closed disc of radius $p^{\frac{-1}{p-1}}$. Since $\mathbb{\Delta}_X(Y) = \operatorname{Spf} W\langle x/p\rangle\hat{}$, its image in $Y_{\widehat{\overline{K}}}$ is the closed disc of radius $p^{\frac{-p}{p-1}}$. Finally, $\mathbb{D}_{X^{(\phi^n)}}(Y) = \operatorname{Spf} W[x^{p^n}/p]\hat{}$, so the associated rigid space is the closed disc of radius $p^{-p^{-n}}$. The union of these is the open disc of radius 1, which is the neighborhood corresponding to convergent cohomology.

**Corollary 2.4.1.** *Let $Y/S$ be a morphism of crystalline prisms. (1.1.3). If $X$ is*

*a regularly immersed closed subscheme of $\overline{Y}$, there is a commutative diagram:*

$$\begin{array}{ccccc}
\mathbb{\Delta}_X(Y) & \xrightarrow{\Psi} & \mathbb{PD}_X(Y) & & \\
& \searrow^{\mathbb{\Delta}_X(F_{Y/S})} & \downarrow \Phi & \searrow^{\mathbb{PD}_X(F_{Y/S})} & \\
& & \mathbb{\Delta}_{X^{(p)}}(Y^{(p)}) & \xrightarrow{\Psi^{(p)}} & \mathbb{PD}_{X^{(p)}}(Y^{(p)}).
\end{array}$$

*If $Y/S$ is $p$-completely smooth, then $\Phi$ is $p$-completely faithfully flat.*

*Proof.* Let $X^\phi := F_Y^{-1}(X) = F_{Y/S}^{-1}(X^{(p)})$ and let $\mathrm{inc}\colon X \to X^\phi$ be the inclusion. We then have morphisms of pairs:

$$\begin{array}{ccccc}
(Y, X) & \xrightarrow{(\mathrm{id}, \mathrm{inc})} & (Y, X^\phi) & & \\
& \searrow^{(F_{Y/S}, F_{X/S})} & \downarrow & \searrow^{(F_{Y/S}, F_{X^\phi_{X/S}})} & \\
& & (Y^{(p)}, X^{(p)}) & \xrightarrow{(\mathrm{id}, \mathrm{inc})} & (Y^{(p)}, X^{\phi^{(p)}}),
\end{array}$$

where the vertical map is induced by $F_{Y/S}$. Taking the corresponding diagram of prismatic envelopes and using statement (2) of Theorem 2.3.13 to identify $\mathbb{\Delta}_{X^\phi}(Y)$ with $\mathbb{PD}_X(Y)$ and $\mathbb{\Delta}_{X^{\phi^{(p)}}}(Y^{(p)})$ with $\mathbb{PD}_{X^{(p)}}(Y^{(p)})$, we obtain the diagram in the statement of the corollary. Let us note for future reference that, since $\mathbb{D}_X(Y)$ is a PD-enlargement of $X$ over $Y$, the morphism $\Psi$ in the diagram in the corollary factors naturally as a composition:

$$\Psi : \mathbb{\Delta}_X(Y) \to \mathbb{D}_X(Y) \to \mathbb{PD}_X(Y). \tag{2.4.2}$$

If $Y/S$ is $p$-completely smooth, then $F_{Y/S}$ is $p$-completely faithfully flat, and since $X^\phi$ is the inverse image of $X^{(p)}$, the square in the diagram

$$\begin{array}{ccccc}
\Phi\colon \mathbb{PD}_X(Y) & \xrightarrow{\cong} & \mathbb{\Delta}_{X^\phi}(Y) & \longrightarrow & \mathbb{\Delta}_{X^{(p)}}(Y^{(p)}) \\
& \searrow & \downarrow & & \downarrow \\
& & Y & \xrightarrow{F_{Y/S}} & Y^{(p)}
\end{array} \tag{2.4.3}$$

is Cartesian, by statement (1) of Theorem 2.3.8. Then $\Phi$ is also $p$-completely faithfully flat. □

**Remark 2.4.2.** Suppose that $X \to Y$ admits an $F$-aligned lift $\tilde{X} \to Y$, and let $X \to \tilde{X} \to \mathbb{D}_X(Y)$ be the induced section (see Proposition 2.2.14). Then by Theorem 2.3.13, we can identify $\mathbb{\Delta}_X(Y)$ with $\mathbb{PD}_X(\mathbb{D}_X(Y))$. Furthermore, the map of pairs $(\mathbb{D}_X(Y), X) \to (Y, X)$ induces a map of divided power envelopes:

$$\mathbb{\Delta}_X(Y) \cong \mathbb{PD}_X(\mathbb{D}_X(Y)) \to \mathbb{PD}_X(Y)$$

which is none other than the morphism $\Psi$.

# 3 Diffracting Envelopes

## 3.1 Additive actions

Let $Y/S$ be a $p$-completely smooth morphism of crystalline prisms and $X \to \overline{Y}$ a regular immersion. Recall from (4c) of Theorem 2.2.8 that $\overline{\mathbb{D}}_X(Y)$ has a natural structure of a torsor under the additive group scheme $\mathbb{N}_{X/\overline{Y}}$. We will construct a compatible torsor action of the divided power envelope $\mathbb{N}^{\gamma}_{X/\overline{Y}}$ of the zero section of $\mathbb{N}_{X/\overline{Y}}$ on $\overline{\mathbb{\Delta}}_X(Y)$. This gives further insight into the geometry of $\overline{\mathbb{\Delta}}_X(Y)$ and will be key to the diffraction of prismatic envelopes.

As we shall see later, the data of an action of $\mathbb{N}^{\gamma}_{X/\overline{Y}}$ on $\overline{\mathbb{\Delta}}_X(Y)$ is equivalent to the data of an integrable and quasi-nilpotent connection

$$\nabla\colon \mathcal{O}_{\overline{\mathbb{\Delta}}_X(Y)} \to \mathcal{O}_{\overline{\mathbb{\Delta}}_X(Y)} \otimes \Omega^1_{\mathbb{N}_{X/\overline{Y}}/X}$$

which satisfies the Leibniz rule. Since the Sen operator is conveniently expressed in terms of this connection, we begin with an explicit algebraic construction of $\nabla$. We will explain a direct geometric approach to the action of $\mathbb{N}^{\gamma}_{X/\overline{Y}}$, due to Bhargav Bhatt, in Theorem 3.1.5 below.

Recall that there are affine morphisms:

$$\mathbb{\Delta}_X(Y) \xrightarrow{h} \mathbb{D}_X(Y) \xrightarrow{\pi_{\mathbb{D}}} Y.$$

Since $\overline{\mathbb{D}}_X(Y)$ is the open subset $D_+(p)$ of the exceptional divisor of the blowup of $X$ in $Y$, there is a canonical isomorphism:

$$\Omega^1_{\overline{\mathbb{D}}_X(Y)/X} \cong \pi^*_{\overline{\mathbb{D}}}(\mathcal{C}_{X/\overline{Y}})$$

fitting into a commutative diagram:

$$\begin{array}{ccc} \mathcal{O}_{\overline{\mathbb{D}}_X(Y)} & \xrightarrow{\quad d \quad} & \Omega^1_{\overline{\mathbb{D}}_X(Y)/X} \\ {\scriptstyle d'}\downarrow & & \downarrow{\scriptstyle \cong} \\ \mathcal{O}_{\overline{\mathbb{D}}_X(Y)} \otimes i^*_{X/Y}(\Omega^1_{Y/S}) & \xleftarrow{\ \mathrm{id}\otimes \overline{d}\ } & \mathcal{O}_{\overline{\mathbb{D}}_X(Y)} \otimes \mathcal{C}_{X/\overline{Y}}. \end{array} \tag{3.1.1}$$

where $\overline{d}\colon \mathcal{C}_{X/Y} \to i^*_{X/Y}(\Omega^1_{X/Y})$ is the usual map. The differential $d$ here is, of course, $\mathcal{O}_X$-linear, and since the map $\overline{\mathbb{D}}_X(Y) \to Y$ factors through $X$, the map $d'$ is $\mathcal{O}_Y$-linear, and in fact is a Higgs field.

We will construct our connection by viewing $\mathcal{O}_{\overline{\mathbb{\Delta}}_X(Y)}$ as a sheaf on $\overline{\mathbb{D}}_X(Y)$ via the affine morphism $\overline{h}$. As a matter of fact, we could just as well work with the morphism $h : \mathbb{D}_X(Y) \to Y$. It is not $p$-completely, smooth, but its reduction modulo $p$ factors through the smooth morphism $\overline{\mathbb{D}}_X(Y) \to X$, and, as is easy to see, the natural map

$$\Omega^1_{\mathbb{D}_X(Y)/Y} \to \Omega^1_{\overline{\mathbb{D}}_X(Y)/X}$$

is an isomorphism. Thus, to give a connection

$$\nabla\colon \mathcal{O}_{\mathbb{A}_X(Y)} \to \mathcal{O}_{\mathbb{A}_X(Y)} \otimes \Omega^1_{\mathbb{D}_X(Y)/Y}$$

is equivalent to giving a connection:

$$\nabla\colon \mathcal{O}_{\overline{\mathbb{A}}_X(Y)} \to \mathcal{O}_{\overline{\mathbb{A}}_X(Y)} \otimes \Omega^1_{\mathbb{D}_X(Y)/X}.$$

**Theorem 3.1.1.** *Let $Y/S$ be a $p$-completely smooth morphism of crystalline prisms and let $X \to \overline{Y}$ be a regular immersion. The structure sheaf $\mathcal{O}_{\mathbb{A}_X(Y)}$, viewed as a sheaf of $\mathcal{O}_{\mathbb{D}_X(Y)}$-modules, admits a quasi-nilpotent and integrable connection:*

$$\nabla\colon \mathcal{O}_{\mathbb{A}_X(Y)} \to \mathcal{O}_{\mathbb{A}_X(Y)} \otimes \Omega^1_{\mathbb{D}_X(Y)/Y}$$

*satisfying and uniquely determined by the following properties:*

*(1) $\nabla(h^*(a)) = 1 \otimes da$ if $a \in \mathcal{O}_{\mathbb{D}_X(Y)}$;*

*(2) $\nabla(ab) = b\nabla(a) + a\nabla(b)$ if $a, b \in \mathcal{O}_{\mathbb{A}_X(Y)}$;*

*(3) $\nabla(\delta(a)) = -a^{p-1}\nabla(a)$ if $a \in \mathcal{O}_{\mathbb{A}_X(Y)}$.*

*Proof.* Let us first remind the reader that the map $\delta\colon \mathcal{O}_{\mathbb{A}} \to \mathcal{O}_{\mathbb{A}}$ does not descend to $\mathcal{O}_{\overline{\mathbb{A}}}$, so that condition (3) in the description of $\nabla$ cannot be sensibly written directly on $\mathcal{O}_{\overline{\mathbb{A}}}$. However, it is reassuring to note that $\nabla(\delta(a))$ will in fact depend only on the reduction of $a$ modulo $p$, as the following calculation shows. If $a, b \in \mathcal{O}_{\mathbb{A}}$:

$$\begin{aligned}\nabla(\delta(a+pb)) &= \nabla\left(\delta(a) + \delta(pb) - C(a,pb)\right)\\ &= \nabla\left(\delta(a)\right) + \nabla(\delta(p)b^p + p^p\delta(b) + p\delta(a)\delta(b)) - \nabla(C(a,pb))\\ &\equiv \nabla(\delta(a)) \pmod{p\mathcal{O}_{\mathbb{A}}},\end{aligned}$$

since $C(a,pb)$ and $\nabla(b^p) = pb^{p-1}\nabla(b)$ are divisible by $p$.

We begin the proof of the theorem by showing that condition (3) is compatible with the algebra structure.

**Lemma 3.1.2.** *Suppose that $\nabla\colon \mathcal{O}_{\mathbb{A}_X(Y)} \to \mathcal{O}_{\mathbb{A}_X(Y)} \otimes \Omega^1_{\mathbb{D}_X(Y)/Y}$ is a connection satisfying the Leibniz rule and that $a, b \in \mathcal{O}_{\mathbb{A}_X(Y)}$ satisfy condition (3). Then the same is true for $a+b$ and $ab$.*

*Proof.* We write $\mathbb{D}$ for $\mathbb{D}_X(Y)$ and $\mathbb{A}$ for $\mathbb{A}_X(Y)$ in the proof. The computation for products is easy, using the facts that $p\Omega^1_{\mathbb{D}/Y} = 0$ and that, consequently, $\nabla(a^p) = \nabla(b^p) = 0$.

$$\begin{aligned}\nabla(\delta(ab)) &= \nabla\left(a^p\delta(b) + b^p\delta(a) + p\delta(a)\delta(b)\right)\\ &= a^p\nabla(\delta(b)) + b^p\nabla(\delta(a))\\ &= -a^pb^{p-1}\nabla(b) - b^pa^{p-1}\nabla(a)\end{aligned}$$

$$
\begin{aligned}
&= -a^{p-1}b^{p-1}\left(a\nabla(b)+b\nabla(a)\right)\\
&= -(ab)^{p-1}\nabla(ab).
\end{aligned}
$$

To handle the formula for the sum, first recall the following familiar computation, which holds in the polynomial ring $\mathbf{Z}[X,Y]$ and hence universally:

$$
\begin{aligned}
(X+Y)^{p-1}(d(X+Y)) &= p^{-1}d((X+Y)^p)\\
&= p^{-1}d(X^p)+p^{-1}d(Y^p)+dC(X,Y)\\
&= X^{p-1}dX+Y^{p-1}dY+dC(X,Y).
\end{aligned}
$$

We are assuming that $\nabla\colon \mathcal{O}_{\mathbb{A}} \to \mathcal{O}_{\mathbb{A}} \otimes \Omega^1_{\mathbb{D}/Y}$ satisfies the Leibniz rule, *i.e.*, is a derivation, and that it is $\mathcal{O}_Y$-linear. Thus it factors through an $\mathcal{O}_{\mathbb{A}}$-linear map $\tilde{\nabla}\colon \Omega^1_{\mathbb{A}/Y} \to \mathcal{O}_{\mathbb{A}} \otimes \Omega^1_{\mathbb{D}/Y}$. We now compute:

$$
\begin{aligned}
(a+b)^{p-1}(\nabla(a+b)) &= (a+b)^{p-1}\left(\tilde{\nabla}(da+db)\right)\\
&= \tilde{\nabla}\left((a+b)^{p-1}(da+db)\right)\\
&= \tilde{\nabla}\left(a^{p-1}da+b^{p-1}db+dC(a,b)\right)\\
&= a^{p-1}\tilde{\nabla}da+b^{p-1}\tilde{\nabla}db+\tilde{\nabla}dC(a,b)\\
&= a^{p-1}\nabla a+b^{p-1}\nabla b+\nabla C(a,b).
\end{aligned}
$$

Now if $a$ and $b$ satisfy (3), we compute;

$$
\begin{aligned}
\nabla(\delta(a+b)) &= \nabla\left(\delta(a)+\delta(b)-C(a,b)\right)\\
&= \nabla(\delta(a))+\nabla(\delta(b))-\nabla C(a,b)\\
&= -a^{p-1}\nabla(a)-b^{p-1}\nabla(b)-\nabla C(a,b)\\
&= -(a+b)^{p-1}(\nabla(a+b)),
\end{aligned}
$$

as we saw in the previous calculation. This proves our claim. □

Before attempting the construction of $\nabla$, we observe that statement (2b) of Theorem 2.3.3 implies that $\mathcal{O}_{\mathbb{A}_X(Y)}$ is generated as an $\mathcal{O}_{\mathbb{D}_X(Y)}$-algebra by the results of succesive application of the operator $\delta$. This implies that a connection $\nabla$ satisfying properties (1)–(3) is unique, and, consequently, that it will suffice to construct it locally. Thus we may assume that we are in the local situation of statement (2) of Theorem 2.3.3, and we use the notation there. As we have seen, $\mathbb{D}_X(Y)$ is the formal spectrum of the $p$-adic completion of the ring

$$
C := B[t_1,\ldots,t_r]/(pt_1-x_1,\ldots,pt_r-x_r),
$$

where $t_i=\rho(x_i)$. Thus $\overline{\mathbb{D}}_X(Y)$ is the spectrum of the ring $\overline{C}\cong A[t_1,\ldots,t_r]$, where $A=B/(p,x_1,\ldots,x_r)$ and $X=\operatorname{Spec} A$, and

$$
\Omega^1_{\mathbb{D}_X(Y)/Y}\cong\Omega^1_{\mathbb{D}_X(Y)/X}
$$

is the free $\overline{C}$-module generated by $(dt_1,\ldots,dt_r)$. Also, $\mathbb{A}_X(Y)$ is the formal spectrum of the completion of the ring $B^\infty$, obtained by adjoining variables $t_{i,j+1} : j \geq 1$ subject to the relations $pt_{i,j+1} = \delta^j(x_i) - t_{i,j}^p$. Since $d\delta^j(x_i) = 0$, this suggests defining $\nabla(t_{i,j+1})$ to be $1 \otimes dt_{i,1}$ if $j = 0$ and $-t_{i,j}^{p-1}\nabla(t_{i,j})$ if $j > 0$. This is indeed what we shall do. The construction simplifies because we are working modulo $p$.

Consider the polynomial algebra

$$C[t_{i,j+1} : 1 \leq i \leq r, j \geq 0],$$

but identify $t_{i,1}$ with $t_i$. We can define a connection $\nabla$ on this polynomial algebra by specifying its value on the free generators:

$$\nabla(t_{i,j+1}) := \begin{cases} 1 \otimes dt_{i,1} & \text{if } j = 0 \\ (-1)^j (t_{i,j} t_{i,j-1} \cdots t_{i,1})^{p-1} dt_{i,1} & \text{if } j > 0. \end{cases} \tag{3.1.2}$$

and extending to the remaining monomials by the Leibniz rule. Note that

$$\nabla(t_{i,j+1}) = -t_{i,j}^{p-1}\nabla(t_{i,j}) \qquad \text{if } j > 0. \tag{3.1.3}$$

It is evident that this connection is integrable and that its $p$-curvature is quasi-nilpotent, since it maps the subalgebra $C_n$ generated by $\{t_{i,j} : j \leq n\}$ into $C_{n-1}$ (see also Proposition 3.1.4 below). Now $\mathcal{O}_{\overline{\mathbb{\Delta}}_X(Y)}$ is the quotient $\overline{B}^\infty$ of this polynomial algebra by the relations $\delta^j(x_i) = t_{i,j}^p$, and, since $d\delta^j(x_i) = dt_{i,j}^p = 0$, the connection descends to $\overline{B}^\infty$. This connection satisfies conditions (1) and (2), and, thanks to Lemma 3.1.2, it will satisfy (3) if it does so for each $t_{i,j}$. By equation (3.1.3), (3) is equivalent to

$$\nabla(\delta(t_{i,j})) = \nabla(t_{i,j+1}) \qquad \text{if } j > 0.$$

We prove this by induction on $j$.[6] Formula (2c) of Theorem 2.3.3 tells us that, modulo $p$:

$$\delta(t_{i,j}) - t_{i,j+1} = t_{i,j-1}^{p(p-1)}(\delta(t_{i,j-1}) - t_{i,j}) + p^{p-2}((\delta(t_{i,j-1}))^p - t_{i,j}^p).$$

If $j = 1$, this says that $\delta(t_{i,1}) = t_{i,2} - p^{p-2}t_{i,1}^p$, and, since $\nabla$ kills $p$th powers, our assertion is true if $j = 1$. If $j > 1$, we use the induction hypothesis and the fact that $\nabla$ annihilates $p$th powers to conclude:

$$\begin{aligned} \nabla(\delta(t_{i,j})) - \nabla(t_{i,j+1}) &= t_{i,j-1}^{p(p-1)}(\nabla(\delta(t_{i,j-1})) - \nabla(t_{i,j})) \\ &\qquad + p^{p-2}(\nabla((\delta(t_{i,j-1}))^p) - \nabla(t_{i,j}^p)) \\ &= 0. \end{aligned}$$

□

**Remark 3.1.3.** It is straightforward to verify that the connection $\nabla$ constructed in Theorem 3.1.1 is functorial, in the following sense, Suppose that $g \colon Y' \to Y$ is a morphism of crystalline prisms and that $X' \to \overline{Y}'$ and $X \to \overline{Y}$ are regular embeddings such that $g$ maps $X'$ to $X$. Then by the functoriality of envelopes, we find a commutative diagram:

$$\begin{array}{ccc} \mathbb{A}_{X'}(Y') & \xrightarrow{h'} & \mathbb{D}_{X'}(Y') \\ \Big\downarrow{\scriptstyle \mathbb{A}(g)} & & \Big\downarrow{\scriptstyle \mathbb{D}(g)} \\ \mathbb{A}_X(Y) & \xrightarrow{h} & \mathbb{D}_X(Y). \end{array}$$

[6]Hearty thanks to the referee for this argument.

Then the compatibility of the connections on prismatic envelopes is expressed by the commutativity of the following diagram:

$$\begin{array}{ccc} \mathcal{O}_{\overline{\mathbb{A}}_{X'}(Y')} & \xrightarrow{\nabla} & \mathcal{O}_{\overline{\mathbb{A}}_{X'}(Y')} \otimes \Omega^1_{\overline{\mathbb{D}}_{X'}(Y')/X'} \\ \Big\uparrow{\scriptstyle \overline{\mathbb{A}}(g)^*} & & \Big\uparrow{\scriptstyle \overline{\mathbb{A}}(g)^* \otimes \overline{\mathbb{D}}(g)^*} \\ \mathcal{O}_{\overline{\mathbb{A}}_X(Y)} & \xrightarrow{\nabla} & \mathcal{O}_{\overline{\mathbb{A}}_X(Y)} \otimes \Omega^1_{\overline{\mathbb{D}}_X(Y)/X}. \end{array}$$

We shall find it useful to describe the $p$-curvature of the connection $\nabla$.

**Proposition 3.1.4.** *The $p$-curvature of the connection on $\mathcal{O}_{\overline{\mathbb{A}}_X(Y)}$ defined in Theorem 3.1.1 satisfies the following formula. If $f \in \mathcal{O}_{\mathbb{A}_X(Y)}$ reduces to $\overline{f} \in \mathcal{O}_{\overline{\mathbb{A}}_X(Y)}$ and if $D' \in \mathcal{T}_{X^{(p)}/S}$ is the pullback of $D \in \mathcal{T}_{X/S}$, then*

$$\psi_{D'}(\overline{\delta}(f)) = -\overline{f}^{p-1}\psi_{D'}(\overline{f}) + (\nabla_D(\overline{f}))^p.$$

*Proof.* Let us first check that the formula is reasonable. Note that, for fixed $D$, the endomorphism of the $\mathcal{O}_X$-module $\mathcal{O}_{\overline{\mathbb{A}}}$ taking an element $g$ to $(\nabla_D(g))^p$ is $p$-linear. Furthermore, the map $\mathcal{T}_{X/S} \to \mathcal{H}om_{\mathcal{O}_X}(\mathcal{O}_{\overline{\mathbb{A}}_X(Y)}, F_{X*}\mathcal{O}_{\overline{\mathbb{A}}_X(Y)})$ taking $D$ to this homomorphism is also $p$-linear, and hence can be viewed as a map

$$F^*_{X/S}(\mathcal{T}_{X^{(p)}/S}) \to \mathcal{H}om_{\mathcal{O}_X}(\mathcal{O}_{\overline{\mathbb{A}}_X(Y)}, F_{X*}\mathcal{O}_{\overline{\mathbb{A}}_X(Y)}).$$

Thus we could have used $D'$ everywhere in the formula above.

To write the proof, a change of notation will be useful. For $f \in \mathcal{O}_{\mathbb{A}}$ and $i \in [0, p)$, let $\delta_i(f)$ be the image of $(i!)^{-1}f^i$ in $\mathcal{O}_{\overline{\mathbb{A}}}$, let $\delta_p(f)$ be the image of $\delta(f)$, and let $\delta_i(f) = 0$ for $i < 0$. Then for $i \in [0, p]$, we have the formula

$$\nabla(\delta_i(f)) = \delta_{i-1}(f)\nabla(f),$$

since $(p-1)! \equiv -1 \pmod p$. Then a calculation as in the proof of [33, A.1.2] shows that

$$\nabla^p_D(\delta_i(f)) = \delta_{i-1}(f)\nabla^p_D(f) + \delta_{i-p}(f)(\nabla_D(f))^p.$$

Hence:

$$\begin{array}{rcl}\psi_{D'}(\overline{\delta}(f)) & := & \nabla_D^p(\overline{\delta}(f)) - \nabla_{D^p}(\overline{\delta}(f)) \\ & = & \nabla_D^p(\delta_p(f)) - \nabla_{D^p}(\delta_p(f)) \\ & = & \delta_{p-1}(f)\nabla_D^p(f) + (\nabla_D(f))^p - \delta_{p-1}(f)\nabla_{D^p}(f) \\ & = & -f^{p-1}\psi_{D'}(f) + (\nabla_D(f))^p. \end{array}$$

□

The connection constructed in Theorem 3.1.1 leads to a stratification, which can be interpreted as a group action, on $\overline{\mathbb{A}}_X(Y)$. Let us simplify the notation by writing $\overline{\mathbb{A}}, \overline{\mathbb{D}}$, and $\overline{\mathbb{N}}$ for $\overline{\mathbb{A}}_X(Y)$, $\overline{\mathbb{D}}_X(Y)$, and $\mathbb{N}_{X/\overline{Y}}$ respectively. We let $\overline{\mathbb{D}}(1) := \overline{\mathbb{D}} \times_X \overline{\mathbb{D}}$ and $\overline{\mathbb{D}}^\gamma(1)$ denote the divided power envelope of the diagonal, with its two projections $t := p_1$ and $s := p_2$ to $\overline{\mathbb{D}}$. The stratification associated to the connection $\nabla$ is an isomorphism $\epsilon\colon s^*(\mathcal{O}_{\overline{\mathbb{A}}}) \to t^*(\mathcal{O}_{\overline{\mathbb{A}}})$, which fits in a commutative diagram

$$\begin{array}{ccccc} s^*(\mathcal{O}_{\overline{\mathbb{A}}}) & \xrightarrow{\epsilon} & t^*(\mathcal{O}_{\overline{\mathbb{A}}}) & \xrightarrow{\cong} & \mathcal{O}_{\overline{\mathbb{A}}} \otimes \mathcal{O}_{\overline{\mathbb{D}}^\gamma(1)} \\ \uparrow s^* & & \uparrow \mathrm{id} & & \uparrow \mathrm{id} \\ \mathcal{O}_{\overline{\mathbb{A}}} & \xrightarrow{\theta} & t^*(\mathcal{O}_{\overline{\mathbb{A}}}) & \xrightarrow{\cong} & \mathcal{O}_{\overline{\mathbb{A}}} \otimes \mathcal{O}_{\overline{\mathbb{D}}^\gamma(1)}. \end{array}$$

Let $\overline{\mathbb{D}}^{(1)}$ denote the first infinitesimal neighborhood of the diagonal of $\overline{\mathbb{D}}(1)$. The ideal of the diagonal in $\overline{\mathbb{D}}^{(1)}$ identifies with $\Omega^1_{\overline{\mathbb{D}}/X}$, and if we endow this ideal with the trivial divided power structure, we find an embedding of $\overline{\mathbb{D}}^{(1)}$ in $\overline{\mathbb{D}}^\gamma(1)$. Then the connection $\nabla$ is given by the commutative diagram:

$$\begin{array}{ccc} \mathcal{O}_{\overline{\mathbb{A}}} & \xrightarrow{\nabla} & \mathcal{O}_{\overline{\mathbb{A}}} \otimes_{\mathcal{O}_{\overline{\mathbb{D}}}} \Omega^1_{\overline{\mathbb{D}}/X} \\ \downarrow \theta - t^* & & \downarrow \\ \mathcal{O}_{\overline{\mathbb{A}}} \otimes_{\mathcal{O}_{\overline{\mathbb{D}}}} \mathcal{O}_{\overline{\mathbb{D}}^\gamma(1)} & \longrightarrow & \mathcal{O}_{\overline{\mathbb{A}}} \otimes_{\mathcal{O}_{\overline{\mathbb{D}}}} \mathcal{O}_{\overline{\mathbb{D}}^{(1)}}. \end{array} \tag{3.1.4}$$

Since $\nabla$ satisfies the Leibniz rule, the map $\theta$ is compatible with multiplication, and corresponds to a morphism

$$\alpha\colon \overline{\mathbb{A}} \times_{\overline{\mathbb{D}}} \overline{\mathbb{D}}^\gamma(1) \to \overline{\mathbb{A}}.$$

The cocycle condition says that $\alpha$ is a (right) action of the groupoid $\overline{\mathbb{D}}^\gamma(1)$ on $\overline{\mathbb{A}}_X$.

On the other hand, as we have seen in (4c) of Theorem 2.2.8, $\overline{\mathbb{D}}$ is a torsor under $\mathbb{N}$, so that there is an isomorphism of groupoids:

$$(z, v) \mapsto (z, z + v) : \overline{\mathbb{D}} \times_X \overline{\mathbb{N}} \to \overline{\mathbb{D}} \times_X \overline{\mathbb{D}},$$

which induces an isomorphism:

$$\overline{\mathbb{D}} \times_X \overline{\mathbb{N}}^\gamma \to \overline{\mathbb{D}}^\gamma(1).$$

Then we have a commutative diagram:

$$\begin{array}{ccccc} & & \overline{\mathbb{\Delta}} \times_{\overline{\mathbb{D}}} \overline{\mathbb{D}}^\gamma(1) & \xrightarrow{\alpha} & \overline{\mathbb{\Delta}} \\ & \nearrow^{\cong} & \uparrow{\scriptstyle\cong} & & \uparrow{\scriptstyle \mathrm{id}} \\ \overline{\mathbb{\Delta}} \times_X \overline{\mathbb{N}}^\gamma & \xrightarrow{\cong} & \overline{\mathbb{\Delta}} \times_{\overline{\mathbb{D}}} \overline{\mathbb{D}} \times_X \overline{\mathbb{N}}^\gamma & \longrightarrow & \overline{\mathbb{\Delta}}. \end{array} \qquad (3.1.5)$$

Thus the action of the groupoid $\overline{\mathbb{D}}^\gamma$ on $\overline{\mathbb{\Delta}}$ over $X$ can be reinterpreted as an action of the group $\overline{\mathbb{N}}^\gamma$ on $\overline{\mathbb{\Delta}}$ over $\overline{\mathbb{D}}$. The following result, explained to me by Barghav Bhatt, gives a direct geometric construction of this action. and shows that it is torsorial. This last point will be important later. For another proof, which does not use the Witt vector constructions, see Remark 3.4.7.

**Theorem 3.1.5.** *Let $Y$ be a crystalline prism and let $i\colon X \to \overline{Y} \to Y$ be a regular immersion. Then $\overline{\mathbb{\Delta}}_X(Y)$ has a natural structure of an $\mathbb{N}^\gamma_{X/\overline{Y}}$-torsor, compatible with the $\mathbb{N}_{X/\overline{Y}}$-torsor structure of $\overline{\mathbb{D}}_X(Y)$. The action of $\mathbb{N}^\gamma_{X/\overline{Y}}$ on $\overline{\mathbb{\Delta}}_X(Y)$ is the stratification corresponding to the connection constructed in Theorem 3.1.1.*

*Proof.* Let $T$ be a $p$-adic formal scheme, let $t\colon T \to Y$ be a morphism, and let $t_{\mathrm{W}}\colon \mathrm{W}(T) \to Y$ be the corresponding morphism of formal $\delta$-schemes. By Proposition 2.3.7, there is a natural bijection between the set of liftings of $t$ to $\mathbb{\Delta}_X(Y)$ and the set of liftings of $t_{\mathrm{W}}$ to $\mathbb{D}_X(Y)$, equivalently, the set of homomorphisms $\rho_W\colon I_{X/Y} \to t_{\mathrm{W}*}(W(\mathcal{O}_T))$ such that $\rho_W(p) = 1$. As we saw in (4b) of Theorem 2.2.8, the set of such morphisms is naturally a pseudo-torsor under the additive group of homomorphisms $i_*(\mathcal{C}_{X/\overline{Y}}) \to t_{\mathrm{W}*}(\mathcal{O}_{\mathrm{W}(T)})$. If $p\mathcal{O}_T = 0$, Corollary 1.6.9 exhibits a natural isomorphism of groups:

$$\mathrm{Hom}_{\mathcal{O}_Y}(i_*(\mathcal{C}_{X/\overline{Y}}), t_{\mathrm{W}*}(\mathcal{O}_{\mathrm{W}(T)})) \to \mathrm{Hom}_{\mathcal{O}_X}(\Gamma^\bullet(\mathcal{C}_{X/\overline{Y}}), t_*(\mathcal{O}_T)) = \mathbb{N}^\gamma_{X/\overline{Y}}(T).$$

Furthermore, in this case, morphisms $T \to \mathbb{\Delta}_X(Y)$ are the same as morphisms $T \to \overline{\mathbb{\Delta}}_X(Y)$, so we get the desired action of $\mathbb{N}^\gamma_{X/\overline{Y}}$ on $\overline{\mathbb{\Delta}}_X(Y)$. To prove that this pseudo-torsorial action is in fact torsorial, it is enough to prove that $\overline{\mathbb{\Delta}}_X(Y) \to X$ admits a section after some fpqc base change $X' \to X$. In fact, as we saw in Theorem 2.3.8, the map $\overline{\mathbb{\Delta}}_X(Y) \to X$ is affine and faithfully flat, and thus will serve.

It remains to prove that this geometrically constructed action agrees with that provided by the connection described in Theorem 3.1.1 and that this action is torsorial. We shall do this by studying the first order deformations induced by this action.

Suppose that $T$ is scheme in characteristic $p$ and that $t_{\mathbb{A}}\colon T \to \mathbb{A}_X(Y)$ is a morphism. Then $t_Y := \pi_{\mathbb{A}} \circ t_{\mathbb{A}}\colon T \to Y$ factors through a morphism $t_X\colon T \to X$. If $E$ is a quasi-coherent sheaf of $\mathcal{O}_T$-modules, let $\mathcal{O}_T \oplus E$ be the sheaf of $\mathcal{O}_T$-algebras defined so that $J := 0 \oplus E$ is a square zero ideal, and let $T \to T'$ be the corresponding first-order thickening, with its retraction $T' \to T$. Composition with this retraction defines a morphism $t'_{\mathbb{A}}\colon T' \to \mathbb{A}_X(Y)$. We find corresponding morphisms

$$\begin{array}{rcl} t'_X\colon T' & \to & X \\ t'_Y\colon T' & \to & Y \\ w'_{\mathbb{A}}\colon \mathrm{W}(T') & \to & \mathbb{A}_X(Y) \\ w'_{\mathbb{D}}\colon \mathrm{W}(T') & \to & \mathbb{D}_X(Y) \\ \rho'_W\colon I_{X/Y} & \to & \mathcal{O}_{W(T')}, \end{array}$$

as described in Proposition 2.3.7. The map $\mathcal{O}_{T'} \to \mathcal{O}_T$ induces a surjection of sheaves of $\delta$-rings $\mathcal{O}_{W(T')} \to \mathcal{O}_{W(T)}$ whose kernel we denote by $W(J)$, a square zero ideal of $\mathcal{O}_{W(T')}$ which is stable under $\delta$.

Now suppose that $v\colon t_X^*(\mathcal{C}_{X/\overline{Y}}) \to E$ is a homomorphism of $\mathcal{O}_T$-modules. Extend it to a homomorphism of $\mathcal{O}_T$-algebras:

$$\alpha\colon t_X^*(\Gamma^\bullet(\mathcal{C}_{X/\overline{Y}})) \to \mathcal{O}_{T'} = \mathcal{O}_T \oplus E$$

with $\alpha(t_X^*(\Gamma^j(\mathcal{C}_{X/\overline{Y}})) = 0$ for $j \geq 2$ to obtain a point of $\mathbb{N}^\gamma_{X/\overline{Y}}(T')$. Then let $\tilde{t}'_{\mathbb{A}} \in \mathbb{A}_X(Y)$ be the map defined by the action of this point on $t'_{\mathbb{A}}$. Since the restriction of Spec$(\alpha)$ to $T$ vanishes, the restrictions of $\tilde{t}'_{\mathbb{A}}$ and $t'_{\mathbb{A}}$ to $T$ agree, as do the restrictions of the corresponding maps $\tilde{w}'_{\mathbb{A}}$ and $w'_{\mathbb{A}}$ and $\tilde{w}'_{\mathbb{D}}$ and $w'_{\mathbb{D}}$ to $\mathrm{W}(T)$. Thus we find derivations:

$$\begin{array}{rcl} D\colon \mathcal{O}_{\mathbb{A}} \to J \cong E & := & \tilde{t}'^*_{\mathbb{A}} - t'^*_{\mathbb{A}} \\ D_{\mathbb{A}}\colon \mathcal{O}_{\mathbb{A}} \to W(J) & := & \tilde{w}'^*_{\mathbb{A}} - w'_{\mathbb{A}}{*} \\ D_{\mathbb{D}}\colon \mathcal{O}_{\mathbb{D}} \to W(J) & := & \tilde{w}'^*_{\mathbb{D}} - w'^*_{\mathbb{D}}. \end{array}$$

It follows from the definitions that $D_{\mathbb{D}} = D_{\mathbb{A}} \circ h^*$ and that $D = i_T^* \circ D_{\mathbb{A}}$, where $h$ is the morphism $\mathbb{A}_X(Y) \to \mathbb{D}_X(Y)$ and $i_T$ is the morphism $T \to \mathrm{W}(T)$. The morphism $t_{\mathbb{A}}$ allows us to view $J \cong E$ as a sheaf of $\mathcal{O}_{\mathbb{A}}$-modules, and the morphisms $w_{\mathbb{A}}$ and $w_{\mathbb{D}}$ to view $W(J)$ as a sheaf of of $\mathcal{O}_{\mathbb{A}}$ and $\mathcal{O}_{\mathbb{D}}$-modules.

**Lemma 3.1.6.** *With the notation above, the following hold.*

*(1) $D$, $D_{\mathbb{A}}$, and $D_{\mathbb{D}}$ satisfy the Leibniz rule:*

$$\begin{array}{rcl} D(fg) & = & D_{\mathbb{A}}(g) + gD_{\mathbb{A}}(f) \\ D_{\mathbb{A}}(fg) & = & fD_{\mathbb{A}}(g) + gD_{\mathbb{A}}(f) \\ D_{\mathbb{D}}(fg) & = & fD_{\mathbb{D}}(g) + gD_{\mathbb{D}}(f) \end{array}$$

*(2) If $x \in I_{X/Y}$ and $\overline{x}$ is its image in $\mathcal{C}_{X/\overline{Y}}$, then, in Joyal coordinates:*

$$\begin{array}{rcl} D_{\mathbb{D}}(\rho_{\mathbb{D}}(x)) & = & (v(\overline{x}), 0, 0, \ldots) \qquad \textit{and} \\ D_{\mathbb{A}}(\rho_{\mathbb{A}}(x)) & = & (v(\overline{x}), 0, 0, \ldots). \end{array}$$

*(3)* $\delta(D_{\mathbb{D}}(f)) = 0$ *for all* $f \in \mathcal{O}_{\mathbb{D}}$*, and* $\delta(D_{\mathbb{A}}(f)) = 0$ *for all* $f \in \mathcal{O}_{\mathbb{A}}$.

*(4)* $D_{\mathbb{A}}(\delta(f)) = -f^{p-1}D_{\mathbb{A}}(f)$ *and* $D(\delta(f)) = -f^{p-1}D(f)$ *if* $f \in \mathcal{O}_{\mathbb{A}}$.

*Proof.* The first statement is standard consequence of the fact that $J$ and $W(J)$ are square zero ideals. The remaining statements require unwinding the definition of the torsor action. Suppose $x \in I_{X/Y}$ and $\overline{x}$ is its image in $\mathcal{C}_{X/\overline{Y}}$. The definition of $\alpha$ and Corollary 1.6.9 imply that the homomorphism $\alpha_W : \mathcal{C}_{X/\overline{Y}} \to W(\mathcal{O}_{T'})$ corresponding to $\alpha$ takes $\overline{x}$ to

$$(\alpha(\overline{x}), \alpha(\overline{x}^{[p]}), \alpha(\overline{x}^{[p^2]}), \cdots) = (v(\overline{x}), 0, 0 \cdots),$$

since $\alpha$ annihilates $\Gamma^j(\mathcal{C}_{X/\overline{Y}})$ for $j > 1$. Then

$$\begin{aligned} D_{\mathbb{D}}(\rho_{\mathbb{D}}(x)) &:= \tilde{\rho}'_W(x) - \rho'_W(x) \\ &= (v(\overline{x}), 0, 0, \cdots) \end{aligned}$$

This proves statement (2).

Statement (3) will depend on the following lemma.

**Lemma 3.1.7.** *Suppose that* $I$ *is a square zero ideal in a* $\delta$*-ring* $R$*. Then if* $x, y \in I$ *and* $a, b \in R$*, we have:*

$$\delta(ax + by) = \phi(a)\delta(x) + \phi(b)\delta(y).$$

*Proof.* Note that if $x, y \in I$, then $C(x, y) = 0$ and that $x^p = y^p = 0$, since $I^2 = 0$. Then

$$\delta(ax + by) = \delta(ax) + \delta(by) - C(ax, by) = \delta(ax) + \delta(by).$$

Furthermore:

$$\delta(ax) = x^p\delta(a) + a^p\delta(x) + p\delta(a)\delta(x) = \phi(a)\delta(x)$$

and similarly for $\delta(bx)$. □

To prove statement (3), note first that if $f = \rho(x)$, then

$$\delta(D_{\mathbb{D}}(f)) = \delta((v(\overline{x}), 0, 0, \cdots)) = 0.$$

Next, note that Lemma 3.1.7 implies that if $\delta$ annihilates $D_{\mathbb{D}}(f)$ and $D_{\mathbb{D}}(g)$, it also annihilates $D_{\mathbb{D}}(f + g)$, since $D_{\mathbb{D}}(f + g) = D_{\mathbb{D}}(f) + D_{\mathbb{D}}(g)$, and all terms belong to the square zero ideal $W(J)$. Furthermore, $D_{\mathbb{D}}(fg) = fD_{\mathbb{D}}(g) + gD_{\mathbb{D}}(f)$, and so the lemma also implies that $\delta(D_{\mathbb{D}}(fg)) = 0$. Since $\mathcal{O}_{\mathbb{D}}$ is generated as an $\mathcal{O}_Y$-algebra by the image of $\rho_{\mathbb{D}}$, this proves the first part of statment (3).

For the remaining statements of Lemma 3.1.6, suppose that $f \in \mathcal{O}_{\mathbb{A}}$ and calculate as follows:

$$D_{\mathbb{A}}(\delta(f)) \quad := \quad \tilde{w}'^*_{\mathbb{A}}(\delta(f)) - w'^*_{\mathbb{A}}\delta(f))$$

$$
\begin{aligned}
&= \delta(\tilde{w}'^*_{\mathbb{\Delta}}(f)) - \delta(w'^*_{\mathbb{\Delta}}(f)) \\
&= \delta(w'^*_{\mathbb{\Delta}}(f) + D_{\mathbb{\Delta}}(f)) - \delta(w'^*_{\mathbb{\Delta}}(f)) \\
&= \delta(D_{\mathbb{\Delta}}(f)) - C(w^*_{\mathbb{\Delta}}(f), D_{\mathbb{\Delta}}(f)) \\
&= \delta(D_{\mathbb{\Delta}}(f)) - f^{p-1} D_{\mathbb{\Delta}}(f) \qquad (3.1.6)
\end{aligned}
$$

since $D_{\mathbb{\Delta}}(f)$ belongs to the square zero ideal $W(J)$. In particular, if $\delta(D_{\mathbb{\Delta}}(f)) = 0$, we see that $D_{\mathbb{\Delta}}(\delta(f)) = -f^{p-1} D_{\mathbb{\Delta}}(f)$, and then Lemma 3.1.7 implies that

$$\delta(D_{\mathbb{\Delta}}(\delta(f))) = \delta(-f^{p-1} D_{\mathbb{\Delta}}(f)) = -\phi(f^{p-1})\delta(D_{\mathbb{\Delta}}(f)) = 0.$$

We conclude that the set of $f \in \mathcal{O}_{\mathbb{\Delta}}$ such that $\delta(D_{\mathbb{\Delta}}(f))$ vanishes is closed under addition, multiplication, and the application of $\delta$ and that it contains the image of $\mathcal{O}_{\mathbb{D}}$. Since $\mathcal{O}_{\mathbb{\Delta}}$ is generated as an $\mathcal{O}_Y$-$\delta$-algebra by this image, statement (3) holds for all $f \in \mathcal{O}_{\mathbb{\Delta}}$. Statement (4) is a consequence of this and equation (3.1.6). This concludes the proof of Lemma 3.1.6. □

To prove that the action of $\mathbb{N}_{X/\overline{Y}}$ identifies with the stratification deduced from the connection in Theorem 3.1.1, we shall apply the discussion above when $T := \overline{\mathbb{\Delta}}$ and $E := z^*_{\mathbb{\Delta}}(\mathcal{C}_{X/\overline{Y}}) \cong h^*(\Omega^1_{\overline{\mathbb{D}}/X})$. Note that the first infinitesimal neighborhood $\mathbb{N}^1_{X/\overline{Y}}$ of the zero section of $\mathbb{N}_{X/\overline{Y}}$ identifies with the relative spectrum of $\mathcal{O}_X \oplus \mathcal{C}_{X/\overline{Y}}$. Thus $T'$ in our discussion corresponds to $\overline{\mathbb{\Delta}} \times_X \overline{\mathbb{N}}^1_{X/\overline{Y}}$. We have two maps

$$T' = \overline{\mathbb{\Delta}} \times_X \overline{\mathbb{N}}^1_{X/\overline{Y}} \to \mathbb{\Delta},$$

one defined by the projection and one by the restriction of the action of $\overline{\mathbb{N}}^\gamma_{X/\overline{Y}}$ to $\overline{\mathbb{N}}^1_{X/\overline{Y}}$. As we saw, the difference between these induces a map

$$D_{\mathbb{\Delta}} \colon \mathcal{O}_{\mathbb{\Delta}} \to E \cong \mathcal{O}_{\mathbb{\Delta}} \otimes_{\mathcal{O}_{\overline{\mathbb{D}}}} \Omega^1_{\overline{\mathbb{D}}/X}.$$

Diagrams (3.1.4) and (3.1.5) show that this map identifies with the connection corresponding to the action of $\overline{\mathbb{N}}^\gamma_{X/\overline{Y}}$. Lemma 3.1.6 proves that this connection satisfies the properties which uniquely determine the connection described in Theorem 3.1.1. This completes the proof of Theorem 3.1.5. □

We can use Theorem 3.1.5 to compute the relative cohomology of the connection $\nabla$ constructed in Theorem 3.1.1. Note that since we are working relative to $X$, the differentials in the complex are $\mathcal{O}_X$-linear and that the image of $\mathcal{O}_X \to \mathcal{O}_{\overline{\mathbb{\Delta}}}$ is annihilated $\nabla$.

**Corollary 3.1.8.** *Let* $(\mathcal{O}_{\overline{\mathbb{\Delta}}_{X(Y)}} \otimes \Omega^\bullet_{\overline{\mathbb{D}}/X}, d)$ *be the de Rham complex of the connection constructed in Theorem 3.1.1, viewed as a complex of sheaves on* $X$*. Then the natural map:*

$$\mathcal{O}_X \to (\mathcal{O}_{\overline{\mathbb{\Delta}}_{X(Y)}} \otimes \Omega^\bullet_{\overline{\mathbb{D}}/X}, d)$$

*is a quasi-isomorphism.*

*Proof.* This statement is local on $X$, and in fact can be verified after base change along any faithfully flat map $X' \to X$. As we saw in Theorem 2.3.8, the map $\overline{\mathbb{\Delta}}_X(Y) \to X$ is itself faithfully flat, and after this base change, the torsor admits a section and hence becomes isomorphic to the trivial torsor corresponding to the action of $\mathbb{N}^\gamma_{X/\overline{Y}}$ on itself by translation. The map $\overline{h}\colon \overline{\mathbb{\Delta}}_X(Y) \to \overline{\mathbb{D}}_X(Y)$ also provides a compatible trivialization of the $\mathbb{N}_{X/\overline{Y}}$-torsor $\overline{\mathbb{D}}_X(Y)$. Thus, with $X' = \overline{\mathbb{\Delta}}_X(Y)$, we find a commutative diagram:

$$\begin{array}{ccc} X' \times_X \overline{\mathbb{\Delta}}_X(Y) & \xrightarrow{\cong} & X' \times_X \mathbb{N}^\gamma_{X/\overline{Y}} \\ \downarrow & & \downarrow \\ X' \times_X \overline{\mathbb{D}}_X(Y) & \xrightarrow{\cong} & X' \times_X \mathbb{N}_{X/\overline{Y}} \end{array}$$

compatible with the actions of $\mathbb{N}^\gamma_{X/\overline{Y}}$ and $\mathbb{N}_{X/\overline{Y}}$. Then the connections corresponding to these group actions correspond to the connections corresponding to the actions by translation. Thus, after the base change to $X'$, the connection on $\mathcal{O}_{\overline{\mathbb{\Delta}}_X(Y)}$ relative to $\overline{\mathbb{D}}_X(Y)/X$ identifies with the connection on $\mathcal{O}_{\mathbb{N}^\gamma_{X/\overline{Y}}}$ relative to $\mathcal{O}_{\mathbb{N}_{X/\overline{Y}}}/X$. Choosing a local trivialization of $\mathcal{C}_{X/\overline{Y}}$ we find ourselves in the situation of the standard connection on the PD-polynomial algebra $A'\langle x_1, \ldots, x_r\rangle$ relative to the $A'[x_1, \ldots, x_r]/A'$, in which case the result is classical [7]. $\square$

**Remark 3.1.9.** We have seen that, if $X \subseteq Y$ is a regular immersion, then $\overline{\mathbb{D}}_X(Y)$ and $\overline{\mathbb{\Delta}}_X(Y)$, viewed as schemes over $X$, have a natural torsor structure under an action of $\mathbb{N}_{X/\overline{Y}}$ and $\mathbb{N}^\gamma_{X/\overline{Y}}$, respectively. A similar result holds for PD envelopes. For example, Theorem 2.3.13 tells us that $\overline{\mathbb{\Delta}}_{X^\phi}(Y)$ identifies with the PD-envelope $\overline{\mathbb{PD}}_X(Y)$ of $X$ in $\overline{Y}$, and if $Y$ is regular, $F_{\overline{Y}}$ is flat, and Theorem 2.3.8 implies that $\overline{\mathbb{\Delta}}_{X^\phi}(Y)$ is the pullback of $\overline{\mathbb{\Delta}}_X(Y)$ by $F_{\overline{Y}}$. It follows that $\overline{\mathbb{PD}}_X(Y)$, viewed as a scheme over $X^\phi$, is a torsor under an action of the pullback of $\mathbb{N}^\gamma_{X/\overline{Y}}$ under the map $f_{X/Y}\colon X^\phi \to X$ induced by $F_{\overline{Y}}$. In fact, this torsor structure can be constructed directly, without reference to any lift of $\overline{Y}$ and $F_{\overline{Y}}$. Vadim Vologodsky informed me that Gleb Terentiuk has shown how to do this using Bhatt's unpublished theory of the relative de Rham stack. It can also be deduced from the standard connection on the $\mathcal{O}_Y$-module $\mathcal{O}_{\overline{\mathbb{PD}}_X(Y)}$. This connection corresponds to an action of the groupoid $\overline{\mathbb{\Delta}}_Y(1)$. The restriction of this groupoid to the diagonal of $Y(1)$ identifies with the group scheme $F^{-1}_{\overline{Y}}(\mathbb{T}^\gamma_{\overline{Y}/S})$, and the corresponding action of this group factors through an action of $f^{-1}_{X/Y}(\mathbb{N}^\gamma_{X/\overline{Y}})$. We shall not explain the details here.

## 3.2 Multiplicative actions and diffraction

Again let $Y$ be a crystalline prism and let $X \to \overline{Y} \to Y$ be a regular immersion. We explain in this section how a flat lifting $\tilde{X}$ of $X$ in $Y_2$ induces a "diffraction" of the prismatic envelope $\overline{\mathbb{\Delta}}_X(Y)$. This diffraction can be expressed in two equivalent ways: as an action of a suitable vector field (the "Sen operator"), or as an action of the divided power envelope $\mathbf{G}_m^\gamma$ of the origin of the multiplicative group $\mathbf{G}_m$. The latter is perhaps more conceptual, but the former is simpler and perhaps more applicable, so we begin with it.

By way of preparation, suppose that $E$ is a locally free sheaf of finite type on a scheme $X$, and let $\mathbf{V}E := \operatorname{Spec} S^\bullet E$, as described in the usual scheme version of Definition 2.2.4. Then $\pi\colon \mathbf{V}E \to X$ is smooth and the map $v_E\colon E \to \pi_*(\mathcal{O}_{\mathbf{V}E})$ induces an isomorphism $dv_E\colon \pi^*E \to \Omega^1_{\mathbf{V}E/X}$. Its dual is an isomorphism $\mathcal{T}_{\mathbf{V}E/X} \cong \mathcal{H}om(\pi^*E, \mathcal{O}_{\mathbf{V}E})$, so that $\Gamma(\mathbf{V}E, \mathcal{T}_{\mathbf{V}E/X}) \cong \operatorname{Hom}(E, \pi_*(\mathcal{O}_{\mathbf{V}E}))$. Thus the map $v_E$ gives rise to a canonical vector field $\theta_E$ on $\mathbf{V}E$, the so-called "Euler vector field." If $e$ is a section of $E$ and $v_E(e) \in \mathcal{O}_{\mathbf{V}E}$, it follows that $\theta_E(v_E(e)) = e$ that $\theta_E(v_E(e)^i) = ie^i$, and in fact that $\theta_E$ acts as multiplication by $i$ on $S^iE \subseteq \mathcal{O}_{\mathbf{V}E}$. Alternatively, under the identification

$$\mathcal{T}_{\mathbf{V}E/X} \cong \mathcal{H}om(\Omega^1_{\mathbf{V}E/X}, \mathcal{O}_{\mathbf{V}E}) \cong \operatorname{Hom}(E \otimes S^\bullet E, S^\bullet E),$$

provided by $dv_E$, the Euler vector field identifies with multiplication $E \otimes S^\bullet E \to S^\bullet E$. Note that if $f\colon X' \to X$ is a morphism of schemes, if $E'$ is a locally free sheaf of finite rank on $X'$, and if $h\colon E \to f_*(E')$ is a homomorphism we find a commutative diagram:

$$\begin{array}{ccc} f_*(\Omega^1_{\mathbf{V}E'/X'}) & \xrightarrow{\theta_{E'}} & f_*(\mathcal{O}_{\mathbf{V}E'}) \\ \uparrow & & \uparrow \\ \Omega^1_{\mathbf{V}E/X} & \xrightarrow{\theta_{E'}} & \mathcal{O}_{\mathbf{V}E}. \end{array} \tag{3.2.1}$$

Let us now simplify the notation by viewing $v_E$ as an inclusion. If $(e_1, \ldots, e_r)$ is a local basis for $E$, then $(de_1, \ldots, de_r)$ is a local basis for $\Omega^1_{\mathbf{V}E/X}$, and, if $(\partial/\partial e_1, \ldots, \partial/\partial e_r)$ is the dual basis,

$$\theta_E = e_1\partial/\partial e_1 + \cdots + e_r\partial/\partial e_r. \tag{3.2.2}$$

Recall that in characteristic $p$, the $p$th iterate $D^{(p)}$ of a derivation $D$ is again a derivation, and that for every $i \in \mathbf{Z}$, $i^p \equiv i \pmod p$. Thus if $X$ has characteristic $p$, the $p$th iterate of the Euler vector field $\theta_E$ is just $\theta_E$ itself.

Recall from Theorem 2.2.8 that if $X \to \overline{Y} \to Y$ is a regular immersion, the dilatation of $X$ in $Y$ is the same as the dilation, and there is a universal homomorphism $\rho_{\mathbb{D}}\colon I_{X/Y} \to \mathcal{O}_{\mathbb{D}_X(Y)}$ whose reduction modulo $p$ induces a map

$$\overline{\rho}_{\mathbb{D}}\colon \mathcal{C}_{X/Y} \to z_{\mathbb{D}*}(\mathcal{O}_{\overline{\mathbb{D}}_X(Y)}).$$

Since $z^*_{\mathbb{D}}(\mathcal{N}_{X/\overline{Y}}) \cong \mathcal{T}_{\mathbb{D}/X}$ and $\mathcal{N}_{X/\overline{Y}} \cong \mathcal{H}om(\mathcal{C}_{X/\overline{Y}}, \mathcal{O}_X)$, we find a natural isomorphism:

$$\Gamma(\overline{\mathbb{D}}_X(Y), \mathcal{T}_{\mathbb{D}/X}) \cong \operatorname{Hom}(\mathcal{C}_{X/\overline{Y}}, z_{\mathbb{D}*}(\mathcal{O}_{\overline{\mathbb{D}}_X(Y)})). \tag{3.2.3}$$

**Proposition 3.2.1.** *Let $Y$ be a $p$-torsion free $p$-adic formal scheme and $X \to \overline{Y} \to Y$ a regular immersion. Then the following sets of data are equivalent:*

*(1) A flat lifting $\tilde{X} \to Y_2$;*

*(2) A section $X \to \overline{\mathbb{D}}_X(Y)$ of the map $\overline{\mathbb{D}}_X(Y) \to X$;*

*(3) A splitting of the exact sequence*

$$0 \to \mathcal{O}_X \to \mathcal{C}_{X/Y} \to \mathcal{C}_{X/\overline{Y}} \to 0 \tag{3.2.4}$$

*(4) An isomorphism $\overline{\mathbb{D}}_X(Y) \to \mathbb{N}_{X/\overline{Y}}$.*

*These data give rise to:*

*(1) An action of the multiplicative group $\mathbf{G}_m$ on $\overline{\mathbb{D}}_X(Y)$ relative to $X$.*

*(2) A* $\mathbf{Z}$*-grading of $\mathcal{O}_{\overline{\mathbb{D}}_X(Y)}$.*

*(3) A vector field $\theta_{\tilde{X}} \in \Gamma(\overline{\mathbb{D}}_X(Y), \mathcal{T}_{\overline{\mathbb{D}}/X})$ such that $\theta_{\tilde{X}}^{(p)} = \theta_{\tilde{X}}$. In fact, $\theta_{\tilde{X}}$ is the vector field obtained by composing the splitting $\mathcal{C}_{X/\overline{Y}} \to \mathcal{C}_{X/Y}$ of the exact sequence (3.2.4) with the map $\overline{\rho}_{\mathbb{D}} \colon \mathcal{C}_{X/Y} \to z_{\mathbb{D}*}(\mathcal{O}_{\overline{\mathbb{D}}_X(Y)})$ and the isomorphism (3.2.3).*

*We call $\theta_{\tilde{X}}$ the* dilatation Sen operator *defined by the lifting $\tilde{X}$ of $X$.*

*Proof.* The equivalence among the first four sets of data follows from Theorem 2.2.8; see also Proposition 1.7.1. The isomorphism $\overline{\mathbb{D}}_X(Y) \to \mathbb{N}_{X/\overline{Y}}$ gives $\overline{\mathbb{D}}_X(Y)$ the structure of a vector bundle over $X$ and an isomorphism $\mathcal{O}_{\overline{\mathbb{D}}_X(Y)} \cong \operatorname{Spec} S^\bullet \mathcal{C}_{X/\overline{Y}}$, hence the action of $\mathbf{G}_m$, the grading of its structure sheaf, and the Euler vector field. □

**Remark 3.2.2.** If $\tilde{X}$ is a flat lifting of $X$ in $Y_2$, and if $v \in \Gamma(X, \mathcal{N}_{X/\overline{Y}})$, and if $\tilde{X} + v$ is the flat lifting of $X$ defined by the action of $\mathbb{N}_{X/\overline{Y}}$ on the set of liftings. then

$$\theta_{\tilde{X}+v} = \theta_{\tilde{X}} + v \in \Gamma(\overline{\mathbb{D}}_X(Y), \mathcal{T}_{\overline{\mathbb{D}}/X}) \cong \Gamma(\overline{\mathbb{D}}_X(Y), z^*_{\mathbb{D}}(\mathcal{N}_{X/\overline{Y}})).$$

We have seen that a flat lifting $\tilde{X}$ of $X$ in $Y_2$ induces a $\mathbf{Z}$-grading of the structure sheaf of $\overline{\mathbb{D}}_X(Y)$ and a derivation $\theta_{\tilde{X}}$ which acts as multiplication by $i$ on the $i$th graded piece. We shall use the connection $\nabla$ of Theorem 3.1.1 to lift $\theta_{\tilde{X}}$ to a derivation $\Theta_{\tilde{X}}$ of the structure sheaf of $\overline{\mathbb{\Delta}}_X(Y)$. This derivation provides the promised "diffraction" of $\overline{\mathbb{\Delta}}_X(Y)$.

**Definition 3.2.3.** *Let $Y$ be a crystalline prism and $X \to \overline{Y} \to Y$ a regular closed immersion. If $\tilde{X} \to Y_2$ is a flat lifting of $X \to Y$, the* prismatic Sen operator *is the endomorphism*

$$\Theta_{\tilde{X}} := \nabla_{\theta_{\tilde{X}}}$$

*of $\mathcal{O}_{\overline{\mathbb{\Delta}}_X(Y)}$, where $\theta_X$ is the dilatation Sen operator defined by the lifting $\tilde{X}$ and $\nabla$ is the connection constructed in Theorem 3.1.1.*

**Theorem 3.2.4.** *In the situation of Definition 3.2.3, the following results hold.*

*(1) The prismatic Sen operator $\Theta_{\tilde{X}}$ is a derivation of $\mathcal{O}_{\overline{\mathbb{\Delta}}_X(Y)}$ and is compatible with the dilatation Sen operator $\theta_{\tilde{X}}$. That is:*

$$\Theta_{\tilde{X}}(h^*(f)) = h^*(\theta_{\tilde{X}}(f))$$

*for every $f \in \mathcal{O}_{\overline{\mathbb{D}}_X(Y)}$.*

*(2) The endomorphism $\Psi_{\tilde{X}} := \Theta^p_{\tilde{X}} - \Theta_{\tilde{X}}$ of $\mathcal{O}_{\overline{\mathbb{\Delta}}_X(Y)}$ is a derivation and is quasi-nilpotent. The eigenvalues of $\Theta_{\tilde{X}}$ lie in $\mathbf{F}_p$, and $\mathcal{O}_{\overline{\mathbb{\Delta}}_X(Y)}$ admits a decomposition into generalized eigenspaces:*

$$\mathcal{O}_{\overline{\mathbb{\Delta}}_X(Y)} \cong \bigoplus_{\lambda \in \mathbf{F}_p} \mathcal{O}^{\lambda}_{\overline{\mathbb{\Delta}}_X(Y)},$$

*where the action of $\Theta_{\tilde{X}} - \lambda$ on $\mathcal{O}^{\lambda}_{\overline{\mathbb{\Delta}}_X(Y)}$ is quasi-nilpotent. This decomposition is compatible with multiplication and gives $\mathcal{O}_{\overline{\mathbb{\Delta}}_X(Y)}$ the structure of a sheaf of $\mathbf{F}_p$-graded rings. The prismatic Sen operator can be written uniquely as a sum*

$$\Theta_{\tilde{X}} = \Theta_{ss} + \Theta_{nil},$$

*where $\Theta_{ss}$ is semi-simple and $\Theta_{nil}$ is quasi-nilpotent. Each of these is a derivation, and $\Theta_{ss}$ acts as multiplication by $\lambda$ on $\mathcal{O}^{\lambda}_{\overline{\mathbb{\Delta}}_X(Y)}$.*

*Proof.* The first statement follows from properties (1) and (2) of the connection $\nabla$ as described in Theorem 3.1.1. The $p$-curvature $\psi$ of $\nabla$ is locally nilpotent, and its value at the vector field $\theta_{\tilde{X}}$ is given by:

$$\psi_{\theta_{\tilde{X}}} := \nabla^p_{\theta_{\tilde{X}}} - \nabla_{\theta^{(p)}_{\tilde{X}}} = \nabla^p_{\theta_{\tilde{X}}} - \nabla_{\theta_{\tilde{X}}} = \Theta^p_{\tilde{X}} - \Theta_{\tilde{X}}, \tag{3.2.5}$$

since $\theta^{(p)}_{\tilde{X}} = \theta_{\tilde{X}}$. Thus $\Theta^p_{\tilde{X}} - \Theta_{\tilde{X}}$ is locally nilpotent. The rest of the proposition follows formally from this. To see this, note that, in the polynomial ring $\mathbf{F}_p[t]$, we have factorizations:

$$t^p - t = \prod_{\lambda \in \mathbf{F}_p} (t - \lambda), \text{ and } (t^p - t)^n = \prod_{\lambda \in \mathbf{F}_p} (t - \lambda)^n.$$

Thus $\prod_{\lambda \in \mathbf{F}_p} (\Theta_{\tilde{X}} - \lambda)^n = (\Theta^p_{\tilde{X}} - \Theta_{\tilde{X}})^n$ annihilates the kernel $N_n$ of $\psi^n_{\theta_{\tilde{X}}}$ acting on $\mathcal{O}_{\overline{\mathbb{\Delta}}}$. Since the polynomials $(t-\lambda)^n$ for the various $\lambda$ in $\mathbf{F}_p$ are relatively prime, it follows that $N_n$ is the direct sum of submodules on which each $(\Theta_{\tilde{X}} - \lambda)^n$ acts as zero. Since $\mathcal{O}_{\overline{\mathbb{\Delta}}_X(Y)}$ is the union of the spaces $N_n$, it admits the claimed decomposition into generalized eigenspaces. Note that if $n = p^j$ for some $j$, then $(\Theta_{\tilde{X}} - \lambda)^{p^j} = \Theta^{p^j}_{\tilde{X}} - \lambda$, whose kernel $N_{p^k,\lambda}$ is $\{f : \Theta^{p^j}(f) = \lambda f\}$. Since $\Theta^{p^j}$ is again a derivation, it follows that $fg \in N_{p^k,\lambda+\lambda'}$ if $f \in N_{p^k,\lambda}$ and $g \in N_{p^k,\lambda'}$. This implies that the grading is compatible with multiplication. It also implies that the operator $\Theta_{ss}$ which is just multiplication by $\lambda$ on $\mathcal{O}^{\lambda}_{\overline{\mathbb{\Delta}}_X(Y)}$ is a derivation, and that $\Theta_{nil} := \Theta_{\tilde{X}} - \Theta_{ss}$ is also a derivation. □

**Remark 3.2.5.** The Sen operator is functorial in the following sense. Suppose that $g\colon Y' \to Y$ is a morphism of crystalline prisms, that $X' \subseteq \overline{Y}'$ and $X \subseteq \overline{Y}$ are regular immersions admitting mod $p^2$ liftings $\tilde{X}' \subseteq Y'$ and $\tilde{X} \subseteq Y$ such that $g$ maps $\tilde{X}'$ to $\tilde{X}$. Indeed, we have a commutative diagram:

$$\begin{array}{ccc} \overline{\mathbb{D}}_{X'}(Y') & \xrightarrow{\cong} & \mathbb{N}_{X'/\overline{Y}'} \\ \Big\downarrow{\scriptstyle \overline{\mathbb{D}}(g)} & & \Big\downarrow{\scriptstyle \mathbb{N}(g)} \\ \overline{\mathbb{D}}_X(Y) & \xrightarrow{\cong} & \mathbb{N}_{X/\overline{Y}}. \end{array}$$

Then it follows from the functoriality of $\nabla$ (Remark 3.1.3) and of the Euler vector field (diagram 3.2.1) that

$$\Theta_{\tilde{X}'} \circ \overline{\mathbb{\Delta}}(g)^* = \overline{\mathbb{\Delta}}(g)^* \circ \Theta_{\tilde{X}}.$$

**Remark 3.2.6.** Recall that $\Psi_{\tilde{X}}$ is the value of the $p$-curvature of $\nabla$ at the vector field $\theta_{\tilde{X}}$. In fact, $\Psi_{\tilde{X}}$ and the nilpotent part $\Theta_{nil}$ of the prismatic Sen operator are mutually determinative. We have

$$\begin{array}{rcl} \Psi_{\tilde{X}} & := & \nabla^p_{\theta_{\tilde{X}}} - \nabla_{\theta^{(p)}_{\tilde{X}}} \\ & = & \Theta^p_{\tilde{X}} - \Theta_{\tilde{X}} \\ & = & (\Theta_{ss} + \Theta_{nil})^p - (\Theta_{ss} + \Theta_{nil}) \\ & = & \Theta^p_{ss} + \Theta^p_{nil} - \Theta_{ss} - \Theta_{nil} \\ & = & \Theta^p_{nil} - \Theta_{nil}, \end{array}$$

since $\Theta^p_{ss} = \Theta_{ss}$. On the other hand:

$$\begin{array}{rcl} \Theta_{nil} & = & \Theta^p_{nil} - \Psi_{\tilde{X}} \\ & = & (\Theta^p_{nil} - \Psi_{\tilde{X}})^p - \Psi_{\tilde{X}} \\ & = & \Theta^{p^2}_{nil} - \Psi^p_{\tilde{X}} - \Psi_{\tilde{X}} \\ & = & \Theta^{p^3}_{nil} - \Psi^{p^2}_{\tilde{X}} - \Psi^p_{\tilde{X}} - \Psi_{\tilde{X}} \\ & = & -\Psi_{\tilde{X}} - \Psi^p_{\tilde{X}} - \Psi^{p^2}_{\tilde{X}} - \Psi^{p^3}_{\tilde{X}} - \cdots . \end{array}$$

As explained in [8, §3.7], the Sen operator on $\overline{\mathbb{\Delta}}_X(Y)$ can also be viewed as an action of the divided power envelope $\mathbf{G}^\gamma_m$ of the origin of $\mathbf{G}_m$. We can see this from our point of view using the the following abstract version of the argument preceding Theorem 3.1.5.

**Proposition 3.2.7.** *Let $X$ be a scheme in characteristic $p$, let $\mathbf{G}/X$ be a smooth group scheme over $X$, and let $\mathbf{G}^\gamma$ be the divided power envelope of the origin of $\mathbf{G}$. Suppose that $Y \to X$ is a morphism endowed with an action $Y \times_X \mathbf{G} \to Y$ and that $h\colon Z \to Y$ is an affine morphism. Then an integrable and quasi-nilpotent connection*

$$\nabla\colon h_*(\mathcal{O}_Z) \to h_*(\mathcal{O}_Z) \otimes \Omega^1_{Y/X}$$

*which satisfies the Leibniz rule and which extends the constant connection $d$ on $\mathcal{O}_Y$ determines an action of $\mathbf{G}^\gamma$ on $Z$ which is compatible with the action of $\mathbf{G}$ on $Y$. If the action of $\mathbf{G}$ on $Y$ is torsorial, then these two sets of data are equivalent.*

*Proof.* The data of an integrable quasi-nilpotent connection on $h_*(\mathcal{O}_Z)$ is, by standard crystalline theory [7], equivalent to the data of an HPD-stratification on $h_*(\mathcal{O}_Z)$. If the connection satisfies the Leibniz rule this stratification is compatible with the multiplication law of $h_*(\mathcal{O}_Z)$ and defines an action of the groupoid $\mathbb{DP}_{Y/X}(1)$ on $Z$. The compatibility of $\nabla$ with $d$ amounts to compatibility of this action with the action of $\mathbb{DP}_{Y/X}(1)$ on $Y$.

On the other hand, we have a commutative diagram:

$$\begin{array}{ccc} Y & \xrightarrow{\ \mathrm{id}\ } & Y \\ {\scriptstyle (\mathrm{id},\iota)}\downarrow & & \downarrow{\scriptstyle \Delta_Y} \\ Y \times_X \mathbf{G} & \longrightarrow & Y \times_X Y, \end{array}$$

where the bottom arrow is the groupoid morphism sending $(y, g)$ to $(y, yg)$ (see 1.4.5). This induces a morphism of groupoids

$$Y \times_X \mathbf{G}^\gamma \to \mathbb{DP}_Y(1). \tag{3.2.6}$$

Composing with the action of $\mathbb{DP}_Y(1)$ on $X$, we find a right action

$$Z \times_Y Y \times_X \mathbf{G}^\gamma \to Z$$

of the groupoid $Y \times_X \mathbf{G}^\gamma$ on $Z/Y$. This amounts simply to an action of the group $\mathbf{G}^\gamma$ on $Z$. If the action of $\mathbf{G}$ on $Y$ is torsorial, then the morphism (3.2.6) is an isomorphism, and we see the equivalence claimed in the proposition. □

**Corollary 3.2.8.** *Let $\mathbf{G}/X$ be a smooth commutative affine group scheme over a scheme $X$ in characteristic $p > 0$ and let $Z \to X$ be a morphism of schemes. Let*

$$\nabla\colon Lie(\mathbf{G}) \to T_{Z/X}$$

*be an $\mathcal{O}_X$-linear map such that:*

*(1) for $D, D' \in Lie(\mathbf{G})$, $[\nabla_D, \nabla_{D'}] = 0$;*

*(2) for $D \in Lie(\mathbf{G})$, $\nabla_D^p - \nabla_{D^p}$ is quasi-nilpotent.*

*Such a homomorphism is equivalent to an action $Z \times_X \mathbf{G}^\gamma \to Z$.*

*Proof.* Since we expect that the result is well-known, we provide only an informal sketch. We shall only use the case in which $\mathbf{G}$ is the multiplicative group; details in this case can be found in [8, §3.7]. We may construct the action locally on $Z$, so we may assume that $Z/X$ is affine, and we identify $\mathcal{O}_Z$ with its direct image on $X$. We apply Proposition 3.2.7 with $Y := \mathbf{G}$, the trivial torsor, and with $Z \times_X \mathbf{G}$ playing the role of $Z$.

The homomorphism $\nabla$ extends to a map

$$T_{\mathbf{G}/X} \to \mathrm{End}_{\mathcal{O}_X}(\mathcal{O}_Z \otimes_{\mathcal{O}_X} \mathcal{O}_{\mathbf{G}})$$

by

$$\nabla_D(a \otimes b) = \nabla_D(a) \otimes b + a \otimes \langle db, D \rangle.$$

This defines a connection on $\mathcal{O}_{Z \times_X \mathbf{G}}$ which satisfies the Leibniz rule, and the conditions (1) and (2) imply that it is integrable and quasi-nilpotent. The proposition then shows that these data are equivalent to an action of $\mathbf{G}^\gamma$ on $Z \times_X \mathbf{G}$ which is compatible with its action on $\mathbf{G}$. Thus, if $z \in Z, g \in \mathbf{G}$, and $h \in \mathbf{G}^\gamma$, we can write $(z, g)h = (a(z, g, h), gh)$ for some $a(z, g, h) \in Z$. The invariance of the vector fields in $Lie(\mathbf{G})$ implies that in fact this is independent of $g$. Now if $z \in Z$ and $h \in \mathbf{G}^\gamma$, let $zh := a(z, e, h)$, and one checks easily that this defines an action of $\mathbf{G}^\gamma$ on $Z$. □

For details of the proof of the following corollary, we refer again to [8].

**Corollary 3.2.9.** *Let $Z$ be a scheme in characteristic $p$ and let $\Theta$ be a vector field on $Z$ such that $\Theta^{(p)} - \Theta$ is quasi-nilpotent. Then $\Theta$ gives rise to an action of $\mathbf{G}_m^\gamma \cong \mu_p \times \mathbf{G}_a^\gamma$ on $Z$. If $\Theta^{(p)} = \Theta$, this action reduces to an action of $\mu_p$, and if $\Theta$ is quasi-nilpotent, it reduces to an action of $\mathbf{G}_a^\gamma$.*

**Example 3.2.10** Let $Y$ be the affine line over an affine scheme $X = \mathrm{Spec}(R)$ with the standard action of $\mathbf{G}_m$. We use coordinates, writing $Y = \mathrm{Spec}\, R[s]$, $\mathbf{G}_m = \mathrm{Spec}\, R[u, u^{-1}]$, and $\mathbf{G}_m^\gamma = \mathrm{Spec}\, R\langle z \rangle$, where $z := u - 1$. Then an action $r_Y$ of $\mathbf{G}_m$ on $Y$ is given by $r_Y^*(s) = us$, and the groupoid morphism $Y \times \mathbf{G}_m \to Y \times Y$ of Definition 1.4.5 is given by $s_1 \mapsto s, s_2 \mapsto us$. If we let $\sigma := s_2 - s_1$ and $z := u - 1$, we have $Y \times_X Y = \mathrm{Spec}\, R[s, \sigma]$, and the morphism is now written $\sigma \mapsto sz$.

Let $Z = \mathrm{Spec}\, R[s, t]$, with $\nabla(t) = -s^{p-1} ds$, so that

$$\begin{aligned} \Theta(t) = \nabla_{sd/ds}(t) &= -s^p \\ \Theta(s^p) = \nabla_{sd/ds}(s^p) &= 0. \end{aligned}$$

Write $\partial$ for $\nabla_{d/ds}$. Then

$$\begin{aligned} \partial(t) &= -s^{p-1} \\ \partial^i(t) &= (-1)^i (i-1)! s^{p-i} \\ \epsilon(t) &= \sum_i \sigma^{[i]} \partial^i(t) \\ &= t + \sum_{i>0} (-1)^i \sigma^{[i]} (i-1)! s^{p-i} \\ &= t + \sum_{i=1}^{p-1} (-1)^i \frac{\sigma^i s^{p-i}}{i} - \sigma^{[p]}. \end{aligned} \tag{3.2.7}$$

The induced action $r_Y$ of $\mathbf{G}_m^\gamma$ is then given by mapping $\sigma$ to $sz$, and we get

$$\begin{aligned} r^*(t) &= t + \sum_1^{p-1} (-1)^i \frac{s^i z^i s^{p-i}}{i} + \sigma^{[p]} \\ &= t + s^p \left( \sum_1^{p-1} (-1)^i \frac{z^i}{i} + z^{[p]} \right) \\ &= t - s^p \log u, \end{aligned}$$

Here we have used the fact that, in the divided power ring $R[s]\langle z \rangle$, we have

$$\log(1+z) = \sum_{i=1}^{p} (-1)^{i-1} (i-1)! z^{[i]}$$

since we are in characteristic $p$, so $(n-1)!z^{[n]} = 0$ if $n > p$.

**Example 3.2.11** Here we give the simplest possible example of the diffraction of the diagonal of a $\delta$-scheme. Let $Y$ be the formal spectrum of $W[y]\hat{}$, with Frobenius lift given by $\phi(y) = y^p$, let $X$ be the reduction of $Y$ modulo $p$, and let $Y \to Y(1)$ be the diagonal embedding. Then $Y(1)$ is the formal spectrum of $W[y, \xi]\hat{}$, where $\xi := 1 \otimes y - y \otimes 1$, and where we have identified $y$ with $y \otimes 1$. Then:

$$\begin{aligned} \phi(\xi) &= (1 \otimes y)^p - (y \otimes 1)^p \\ &= (y + \xi)^p - y^p \\ &= \xi^p + p! \sum_1^{p-1} y^{[p-i]} \xi^{[i]}. \end{aligned}$$

The dilatation $\mathbb{D}_X(Y(1))$ of $X$ in $Y(1)$ is the formal spectrum of $W[y, \eta]\hat{}$, where $p\eta = \xi$, and, by Theorem 2.3.13, $\mathbb{A}_X(Y)$ is is the formal spectrum of $W[y]\langle \eta \rangle\hat{}$. Then

$$\begin{aligned} \phi(\eta) &= p^{-1}\phi(\xi) \\ &= p^{-1}(p^p \eta^p + p! \sum_1^{p-1} y^{p-i} (p\eta)^{[i]} \\ &= p^p (p-1)! \eta^{[p]} + p(p-1)! \sum_1^{p-1} y^{[p-i]} p^i \eta^{[i]} \\ &\equiv p y^{p-1} \eta \pmod{p^2} \\ \delta(\eta) &:= p^{-1}(\phi(\eta) - \eta^p) \\ &= p^{-1}\phi(\eta) - (p-1)!\eta^{[p]}) \\ &\equiv \eta^{[p]} + y^{p-1}\eta \pmod{p}. \end{aligned}$$

We use Theorem 3.1.1 to calculate the connection $\nabla$ on $\mathcal{O}_{\overline{\mathbb{A}}_X(Y(1))}$:

$$\nabla(\eta^{[p]}) = \nabla(\delta(\eta^{[p]}) - y^{p-1}\eta)$$

$$
\begin{aligned}
&= \nabla(\delta(\eta^{[p]})) - \nabla(y^{p-1}\eta) \\
&= -\eta^{[p-1]}d\eta - y^{p-1}d\eta.
\end{aligned}
$$

Thuis, if $\partial := \nabla_{d/d\eta}$, we have

$$
\partial^i(\eta^{[p]}) = \begin{cases} -\eta^{[p-1]} - y^{p-1} & \text{if } i = 1 \\ -\eta^{[p-2]} & \text{if } 1 < i \leq p \\ 0 & \text{if } i > p. \end{cases}
$$

Now we can (partially) compute the action $r$ of $\mathbf{G}_m^\gamma$ on $\overline{\mathbb{\Delta}}_X(Y(1))$. Writing $z$ for $u-1$ as before, we get:

$$
\begin{aligned}
r^*(\eta^{[p]}) &= \sum_i (z\eta)^{[i[}\partial^i(\eta^{[p]}) \\
&= \eta^{[p]} - \sum_{i>0}(z\eta)^{[i]}\eta^{[p-i]} - z\eta y^{p-1} \\
&= \eta^{[p]} - (u-1)y^{p-1}\eta.
\end{aligned}
$$

## 3.3 Diffracted complexes

Let $Y/S$ be a $p$-completely smooth morphism of crystalline prisms and $X \to \overline{Y}$ a regular immersion, where $\overline{X} \to \overline{S}$ is smooth. Our aim in this subsection is to construct some complexes which compute (in a sense to be explained later) the prismatic cohomology of $X/\overline{S}$ and which admit an explicit diffraction expressed as either a Sen operator or an action of $\mathbf{G}_m^\gamma$.

The Cech-Alexander complex is the most obvious example of such a complex. Let $Y(n)$ be the $(n+1)$-fold product of $Y/S$ with itself and form the prismatic envelope $\mathbb{\Delta}_X(Y(n))$ of the diagonal embedding $X \to Y \to Y(n)$. (This will automatically be a regular embedding if $Y$ is noetherian.) A mod $p^2$ lifting $\tilde{X}$ of $X$ in $Y$ propagates to liftings in each $Y(n)$, so each $\overline{\mathbb{\Delta}}_X(Y(n))$ admits a Sen operator, and it follows from the functoriality of the Sen operators (Remark 3.2.5) that these are compatible with the face and degeneracy maps of the standard simplicial complex constructed from the simplicial scheme $\overline{\mathbb{\Delta}}_X(Y(\bullet))$. As we verify in Proposition 6.3.2, the resulting Cech-Alexander complex does indeed calculate prismatic cohomology (and for Drinfeld, was the very definition of prismatic cohomology). However, it is not so easy to calculate, and in fact we found it problematic to use this construction to verify even the most basic expected properties of the Sen operator.

Our goal is to find a natural differential graded algebra, constructed from a natural Higgs field on $\overline{\mathbb{\Delta}}_X(Y)$, which calculates prismatic cohomology and which is compatible with diffraction of prismatic envelopes we constructed in Definition 3.2.3. This project also turned out to be problematic, as the first Higgs field we constructed—the reduction modulo $p$ of a natural $p$-connection on $\mathbb{\Delta}_X(Y)$ (see Proposition 4.2.1) is *not* compatible with diffraction. Instead we show here how to use the connection $\nabla$ of Theorem 3.1.1 to construct a Higgs field which *is* diffraction compatible. We investigate the relationship between these two Higgs fields in Theorem 5.4.3.

We begin by investigating the compatibility between the connection $\nabla$ and the Sen operator. Note that the diffracted Sen operator $\theta_{\tilde{X}}$ on $\mathcal{O}_{\mathbb{D}_X(Y)}$ is already not strictly compatible with differential $d$ or the Higgs field $d'$ of $\mathcal{O}_{\mathbb{D}_X(Y)}$ (see diagram 3.1.1) since they decrease degrees by one. Here is the analogous statement for the connection $\nabla$ on $\overline{\Delta}_X(Y)$.

**Proposition 3.3.1.** *In the situation of Theorem 3.2.4, the prismatic Sen operator $\Theta_{\tilde{X}}$ and the connection $\nabla$ satisfy the following compatibilities.*

*(1) $\Theta_{\tilde{X}}$ fits into a commutative diagram:*

$$\begin{array}{ccc}
\mathcal{O}_{\overline{\Delta}_X(Y)} & \xrightarrow{\nabla} & \mathcal{O}_{\overline{\Delta}_X(Y)} \otimes \Omega^1_{\mathbb{D}_X(Y)/X} \\
\Big\downarrow \Theta_{\tilde{X}} & & \Big\downarrow (\Theta_{\tilde{X}} + \mathrm{id}) \otimes \mathrm{id} \\
\mathcal{O}_{\overline{\Delta}_X(Y)} & \xrightarrow{\nabla} & \mathcal{O}_{\overline{\Delta}_X(Y)} \otimes \Omega^1_{\mathbb{D}_X(Y)/X}.
\end{array}$$

*(2) The endomorphism $\Theta^{\bullet}_{\tilde{X}}$ of the de Rham complex of $(\mathcal{O}_{\overline{\Delta}_X(Y)}, \nabla)$ given by $(\Theta_{\tilde{X}} + j) \otimes \mathrm{id}$ in degree $j$ is a derivation and homotopic to zero. Namely, let*

$$r^j_{\tilde{X}} : \mathcal{O}_{\overline{\Delta}_X(Y)} \otimes \Omega^j_{\mathbb{D}/X} \to \mathcal{O}_{\overline{\Delta}_X(Y)} \otimes \Omega^{j-1}_{\mathbb{D}/X}$$

*be interior multiplication by $\theta_{\tilde{X}}$. Then*

$$\Theta^j_{\tilde{X}} = d^{j-1} \circ r^j_{\tilde{X}} + r^{j+1}_{\tilde{X}} \circ d^j.$$

*Proof.* We can check statement (1) locally, with the aid of a regular sequence $(x_1, \ldots, x_r)$ generating the ideal of $\tilde{X}$ in $Y_2$. We get a corresponding basis $(e_1, \ldots, e_r)$ for $\mathcal{C}_{X/\overline{Y}}$ and hence coordinates for $\overline{\mathbb{D}}_X(Y) \cong \mathbb{N}_{X/\overline{Y}}$ and a basis $(de_1, \ldots, de_r)$ for $\Omega^1_{\mathbb{D}_X(Y)/X}$. Let $(D_1, \ldots, D_r)$ be the dual basis. Then

$$\begin{aligned}
\theta_{\tilde{X}} &= \sum_j e_j D_j \\
D_i \theta_{\tilde{X}} &= D_i \sum_j e_j D_j \\
&= \sum_j (D_i(e_j) D_j + e_j D_i D_j) \\
&= \sum_j (\delta_{i,j} D_j + e_j D_j D_i) \\
&= D_i + \sum_j e_j D_j D_i \\
&= D_i + \theta_{\tilde{X}} D_i \\
[D_i, \theta_{\tilde{X}}] &= D_i.
\end{aligned}$$

Since the connection $\nabla$ is integrable, we conclude:

$$\begin{aligned}\nabla_{D_i}\Theta_{\tilde{X}} &= \nabla_{D_i}\nabla_{\theta_{\tilde{X}}}\\ &= \nabla_{[D_i,\theta_{\tilde{X}}]} + \nabla_{\theta_{\tilde{X}}}\nabla_{D_i}\\ &= \nabla_{D_i} + \nabla_{\theta_{\tilde{X}}}\nabla_{D_i}\\ &= \nabla_{D_i} + \Theta_{\tilde{X}}\nabla_{D_i}.\end{aligned}$$

This implies that $\nabla \circ \Theta_{\tilde{X}} = \nabla + (\Theta_{\tilde{X}} \otimes \mathrm{id}) \circ \nabla$, as asserted by the diagram.

For statement (2), first observe that both sides of the purported equation are derivations with respect to the multiplication of the complex $\mathcal{O}_{\overline{\mathbb{A}}_X(Y)} \otimes \Omega^\bullet_{\overline{\mathbb{D}}/X}$, that is, they satisfy: $\kappa(\alpha\beta) = \kappa(\alpha)\beta + \alpha\kappa(\beta)$. For example, the endomorphism $\kappa^\bullet$ given by multiplication by $j$ on each $\mathcal{O}_{\overline{\mathbb{A}}_X(Y)} \otimes \Omega^j_{\overline{\mathbb{D}}/X}$ is such a derivation, as is $\Theta_{\tilde{X}} \otimes \mathrm{id}$. Therefore $\Theta^\bullet_{\tilde{X}}$ is also a derivation. On the other hand, interior multiplication by any vector field is a derivation of degree $-1$ (so that a sign appears in the characterizing formula), and $d$ is a derivation of degree 1, and it follows that $d^\bullet r^\bullet + r^\bullet d^\bullet$ is a derivation of degree zero.

When $j = 0$, the formula asserts that $\Theta_{\tilde{X}} = r^1_{\tilde{X}} \circ \nabla$, which is just the definition of $\Theta_{\tilde{X}}$. When $j = 1$, it says that

$$\Theta_{\tilde{X}} \otimes \mathrm{id} + \mathrm{id} = \nabla \circ r^1_{\tilde{X}} + r^2 \circ d^1,$$

and so it predicts that

$$(\Theta_{\tilde{X}} \otimes \mathrm{id}) \circ \nabla + \nabla = \nabla \circ r^1_{\tilde{X}} \circ \nabla,$$

*i.e.*, that $(\Theta_{\tilde{X}} \otimes \mathrm{id}) \circ \nabla + \nabla = \nabla \circ \Theta_{\tilde{X}}$, which we have already verified. Thus our equation holds when applied to any element of $\mathcal{O}_{\overline{\mathbb{A}}}$ and also to any element in the image of $\nabla$. Since the complex is generated as an $\mathcal{O}_{\overline{\mathbb{A}}}$-algebra by such elements and since both sides are derivations, the equation holds everywhere. □

The twisted compatibility between $\Theta_{\tilde{X}}$ and $\nabla$ leads to the following decomposition of the de Rham complex $(\mathcal{O}_{\overline{\mathbb{A}}_X(Y)} \otimes \Omega^\bullet_{\overline{\mathbb{D}}/X}, d)$, *i.e.*, its "diffraction."

**Proposition 3.3.2.** *In the situation of Theorem 3.2.4, the de Rham complex of* $(\mathcal{O}_{\overline{\mathbb{A}}_X(Y)}, \nabla)$ *and the prismatic Sen operator* $\Theta_{\tilde{X}}$ *admit a canonical "diffraction" as follows:*

(1) *For each* $\lambda \in \mathbf{F}_p$, *the connection* $\nabla$ *maps* $\mathcal{O}^\lambda_{\overline{\mathbb{A}}_X(Y)}$ *to* $\mathcal{O}^{\lambda-1}_{\overline{\mathbb{A}}_X(Y)} \otimes \Omega^1_{\overline{\mathbb{D}}/X}$, *and there is a direct sum decomposition*

$$(\mathcal{O}_{\overline{\mathbb{A}}_X(Y)} \otimes \Omega^\bullet_{\overline{\mathbb{D}}/X}, d) \cong \bigoplus_\lambda (\overline{\mathbb{A}}^\bullet_\lambda/\overline{\mathbb{D}}, d),$$

*where*

$$(\overline{\mathbb{A}}^\bullet_\lambda/\overline{\mathbb{D}}, d) := \mathcal{O}^\lambda_{\overline{\mathbb{A}}_X(Y)} \xrightarrow{d} \mathcal{O}^{\lambda-1}_{\overline{\mathbb{A}}_X(Y)} \otimes \Omega^1_{\overline{\mathbb{D}}/X} \xrightarrow{d} \mathcal{O}^{\lambda-2}_{\overline{\mathbb{A}}_X(Y)} \otimes \Omega^2_{\overline{\mathbb{D}}/X} \longrightarrow \cdots.$$

(2) *The natural map*

$$\mathcal{O}_X \to (\overline{\mathbb{A}}^\bullet_0/\overline{\mathbb{D}}, d))$$

*is a quasi-isomorphism, and* $(\overline{\mathbb{A}}^\bullet_\lambda/\overline{\mathbb{D}}, d)$ *is acyclic if* $\lambda \neq 0$.

(3) *The nilpotent part $\Theta_{nil}$ of the prismatic Sen operator is a horizontal endomorphism of $(\mathcal{O}_{\overline{\mathbb{\Delta}}_X(Y)}, \nabla)$ and is compatible with the direct sum decomposition in statement (1).*

*Proof.* Let us again write $\mathbb{D}$ for $\mathbb{D}_X(Y)$ and $\mathbb{\Delta}$ for $\mathbb{\Delta}_X(Y)$. Proposition 3.3.1 implies that $\nabla \circ (\Theta_{\tilde{X}} - \lambda) = (\Theta_{\tilde{X}} - (\lambda - 1)) \circ \nabla$. Since $\mathcal{O}^{\lambda}_{\overline{\mathbb{\Delta}}} = \cup_m \operatorname{Ker}(\Theta_{\tilde{X}} - \lambda)^n$, it follows that $\nabla$ takes $\mathcal{O}^{\lambda}_{\overline{\mathbb{\Delta}}}$ into $\mathcal{O}^{\lambda-1}_{\overline{\mathbb{\Delta}}} \otimes \Omega^1_{\overline{\mathbb{D}}/X}$. This gives us the direct sum decomposition in (1). Corollary 3.1.8 says that the map

$$\mathcal{O}_X \to (\mathcal{O}_{\overline{\mathbb{\Delta}}} \otimes \Omega^{\bullet}_{\overline{\mathbb{D}}/X}, d)$$

is a quasi-isomorphism, and since this map factors through $(\overline{\mathbb{\Delta}}^{\bullet}_0/\overline{\mathbb{D}}, d)$, statement (2) follows. Since $\Theta_{nil} = \Theta_{\tilde{X}} - \lambda$ on $\mathcal{O}^{\lambda}_{\overline{\mathbb{\Delta}}}$ it follows that

$$\nabla \circ \Theta_{nil} = (\Theta_{nil} \otimes \mathrm{id}) \circ \nabla.$$

Thus $\Theta_{nil}$ is horizontal, and it is clearly compatible with the direct sum decomposition in statement (1). ☐

The reader may wonder why we are devoting so much attention to an acyclic complex. It turns out that this complex can be used as a piece of larger complex that is closely related to de Rham cohomology and prismatic cohomology. We will have to defer the explanation of that relationship till later (see Theorem 5.4.3), but it seems sensible to describe the construction here.

We assume now that $Y \to S$ is a $p$-completely smooth morphism of crystalline prisms and that $i \colon X \to Y$ is a closed subscheme such that $X \to \overline{S}$ is smooth and $\bar{i} \colon X \to \overline{Y}$ is a regular immersion. Let us also assume that $\overline{Y} \to \overline{S}$ is of finite type. We have an exact sequence of locally free sheaves on $X$:

$$0 \to \mathcal{C}_{X/\overline{Y}} \to i^*(\Omega^1_{Y/S}) \to \Omega^1_{X/\overline{S}} \to 0.$$

Using the isomorphism (3.1.1) $z^*_{\mathbb{D}}(\mathcal{C}_{X/\overline{Y}}) \cong \Omega^1_{\overline{\mathbb{D}}/X}$, we find an exact sequence:

$$0 \longrightarrow \Omega^1_{\overline{\mathbb{D}}/X} \xrightarrow{j} \mathcal{O}_{\overline{\mathbb{D}}} \otimes \Omega^1_{\overline{Y}/\overline{S}} \longrightarrow \mathcal{O}_{\overline{\mathbb{D}}} \otimes \Omega^1_{X/\overline{S}} \longrightarrow 0, \tag{3.3.1}$$

which we regard as a sequence of quasi-coherent sheaves on $X$. Composing $j$ with the differential $d \colon \mathcal{O}_{\overline{\mathbb{D}}} \to \Omega^1_{\overline{\mathbb{D}}/X}$, we find an $\mathcal{O}_{\overline{Y}}$-linear map:

$$d' \colon \mathcal{O}_{\overline{\mathbb{D}}} \to \mathcal{O}_{\overline{\mathbb{D}}} \otimes \Omega^1_{\overline{Y}/\overline{S}}, \tag{3.3.2}$$

which is in fact a Higgs field. The associated Higgs complex is:

$$\mathcal{O}_{\overline{\mathbb{D}}} \xrightarrow{d'} \mathcal{O}_{\overline{\mathbb{D}}} \otimes \Omega^1_{\overline{Y}/\overline{S}} \xrightarrow{d'} \mathcal{O}_{\overline{\mathbb{D}}} \otimes \Omega^2_{\overline{Y}/\overline{S}} \xrightarrow{d'} \cdots,$$

where $d'(f \otimes \omega) := d'f \wedge \omega$ for $f \in \mathcal{O}_{\overline{\mathbb{D}}}$ and $\omega \in \Omega^j_{\overline{Y}/\overline{S}}$. Tensoring with $\mathcal{O}_{\overline{\mathbb{\Delta}}}$ and using the connection $\nabla$ in place of $d'$, we extend $'$ to a Higgs field:

$$\nabla' \colon \mathcal{O}_{\overline{\mathbb{\Delta}}_X(Y)} \to \mathcal{O}_{\overline{\mathbb{\Delta}}_X(Y)} \otimes \Omega^1_{\overline{Y}/\overline{S}}, \tag{3.3.3}$$

fitting into a complex:

$$\mathcal{O}_{\overline{\Delta}_{X(Y)}} \xrightarrow{\nabla'} \mathcal{O}_{\overline{\Delta}_{X(Y)}} \otimes \Omega^1_{\overline{Y}/\overline{S}} \xrightarrow{\nabla'} \mathcal{O}_{\overline{\Delta}_{X(Y)}} \otimes \Omega^2_{\overline{Y}/\overline{S}} \longrightarrow \cdots$$

constructed in the same way. This leads to the following result, which will be key to our understanding of diffraction of prismatic cohomology.

**Theorem 3.3.3.** *Let $Y/S$ be a $p$-completely smooth morphism of crystalline prisms, let $X/\overline{S}$ be a smooth morphism of schemes, and suppose that $X \to \overline{Y}$ is a regular closed immersion over $S$.*

*(1) The composition $\nabla'$ of $\nabla\colon \mathcal{O}_{\overline{\Delta}} \to \mathcal{O}_{\overline{\Delta}} \otimes \Omega^1_{\mathbb{D}/X}$ with the map*

$$\mathcal{O}_{\overline{\Delta}} \otimes \Omega^1_{\mathbb{D}/X} \cong \mathcal{O}_{\overline{\Delta}} \otimes \mathcal{C}_{X/\overline{Y}} \to \mathcal{O}_{\overline{\Delta}} \otimes \Omega^1_{\overline{Y}/\overline{S}}$$

*defines a Higgs field*

$$\nabla'\colon \mathcal{O}_{\overline{\Delta}} \longrightarrow \mathcal{O}_{\overline{\Delta}} \otimes \Omega^1_{\overline{Y}/\overline{S}},$$

*fitting into a complex:*

$$(\overline{\Delta}^{\bullet}(X/Y/S), \nabla') := \mathcal{O}_{\overline{\Delta}} \xrightarrow{\nabla'} \mathcal{O}_{\overline{\Delta}} \otimes \Omega^1_{\overline{Y}/\overline{S}} \xrightarrow{\nabla'} \mathcal{O}_{\overline{\Delta}} \otimes \Omega^2_{\overline{Y}/\overline{S}} \xrightarrow{\nabla'} \cdots,$$

*where $\nabla'(f \otimes \omega) = \nabla'(f) \wedge \omega$ for $f \in \mathcal{O}_{\overline{\Delta}}$ and $\omega \in \Omega^i_{\overline{Y}/\overline{S}}$.*

*(2) For $i \in \mathbf{Z}$, there are natural isomorphisms:*

$$\mathcal{H}^i(\overline{\Delta}^{\bullet}(X/Y/S), \nabla') \cong \Omega^i_{X/\overline{S}}.$$

*(3) Let $\Theta_{\tilde{X}}$ be the prismatic Sen operator associated with a flat lifting $\tilde{X}$ of $X$ in $Y_2$, and for $\lambda \in \mathbf{F}_p$, let $\mathcal{O}^{\lambda}_{\overline{\Delta}}$ denote its generalized $\lambda$-eigenspace.*

*(a) Let*

$$\Theta^j_{\tilde{X}} := (\Theta_{\tilde{X}} + j) \otimes \mathrm{id}\colon \mathcal{O}_{\overline{\Delta}_{X(Y)}} \otimes \Omega^j_{Y/S} \to \mathcal{O}_{\overline{\Delta}_{X(Y)}} \otimes \Omega^j_{Y/S}.$$

*Then $\Theta^{\bullet}_{\tilde{X}} := \bigoplus \Theta^j_{\tilde{X}}$ defines an endomorphism of the complex $(\overline{\Delta}^{\bullet}(X/Y/S), \nabla')$ satisfying the Leibniz rule:*

$$\Theta^{\bullet}_{\tilde{X}}(\alpha \wedge \beta) = \Theta^{\bullet}_{\tilde{X}}(\alpha) \wedge \beta + \alpha \wedge \Theta^{\bullet}_{\tilde{X}}(\beta).$$

*(b) The quasi-nilpotent part $\Theta_{nil}$ of $\Theta_{\tilde{X}}$ commutes with the Higgs field $\nabla'$ and defines an endomorphism of the complex $\overline{\Delta}^{\bullet}(X/Y/S)$.*

*(c) There is a a direct sum decomposition, compatible with the action of $\Theta_{nil}$:*

$$(\overline{\Delta}^{\bullet}(X/Y/S), \nabla') = \bigoplus_{\lambda \in \mathbf{F}_p} (\overline{\Delta}^{\bullet}_{\lambda}(X/Y/S), \nabla'),$$

*where*

$$(\overline{\Delta}^{\bullet}_{\lambda}/\overline{S}, \nabla') := \mathcal{O}^{\lambda}_{\overline{\Delta}_{X(Y)}} \xrightarrow{\nabla'} \mathcal{O}^{\lambda-1}_{\overline{\Delta}_{X(Y)}} \otimes \Omega^1_{\overline{Y}/\overline{S}} \xrightarrow{\nabla'} \mathcal{O}^{\lambda-2}_{\overline{\Delta}_{X(Y)}} \otimes \Omega^2_{\overline{Y}/\overline{S}} \longrightarrow \cdots.$$

*(d) For $i \in \mathbf{Z}$, there are natural isomorphisms:*

$$\mathcal{H}^i(\overline{\mathbb{\Delta}}^\bullet_\lambda(X/Y/S), \nabla') \cong \begin{cases} \Omega^i_{X/\overline{S}} & \text{if } i \in \lambda \\ 0 & \text{otherwise.} \end{cases}$$

*Proof.* We have already explained statement (1). For the cohomology calculation, let $K^\bullet$ denote the Koszul filtration of $\mathcal{O}_{\overline{\mathbb{D}}} \otimes \Omega^\bullet_{Y/S}$ induced by the sequence (3.3.1) and let $\mathrm{Gr}_K$ be its associated graded:

$$\begin{aligned} K^i(\mathcal{O}_{\overline{\mathbb{D}}} \otimes \Omega^j_{\overline{Y}/\overline{S}}) &:= Im\left(\Omega^i_{\overline{\mathbb{D}}/X} \otimes \Omega^{j-i}_{\overline{Y}/\overline{S}} \to \mathcal{O}_{\overline{\mathbb{D}}} \otimes \Omega^j_{\overline{Y}/\overline{S}}\right) \\ \mathrm{Gr}^i_K(\mathcal{O}_{\overline{\mathbb{D}}} \otimes \Omega^j_{\overline{Y}/\overline{S}}) &\cong \Omega^i_{\overline{\mathbb{D}}/X} \otimes \Omega^{j-i}_{X/\overline{S}}. \end{aligned} \tag{3.3.4}$$

Recall that a lifting $\tilde{X}$ of $X$ in $Y_2$ produces an isomorphism $\mathcal{O}_{\overline{\mathbb{D}}_X(Y)} \cong S^\bullet \mathcal{C}_{X/\overline{Y}}$, hence a grading of $\mathcal{O}_{\overline{\mathbb{D}}}$ and of the de Rham complex:

$$(\mathcal{O}_{\overline{\mathbb{D}}} \otimes \Omega^\bullet_{\overline{\mathbb{D}}/X}, d) = \bigoplus_{m \in \mathbf{N}} (\mathcal{O}_{\overline{\mathbb{D}}} \otimes \Omega^\bullet_{\overline{\mathbb{D}}/X}, d)^m, \qquad \text{where:}$$

$$(\mathcal{O}_{\overline{\mathbb{D}}} \otimes \Omega^\bullet_{\overline{\mathbb{D}}/X}, d)^m := \mathcal{O}_{\overline{\mathbb{D}}}^{\,m} \to \mathcal{O}_{\overline{\mathbb{D}}}^{\,m-1} \otimes \Omega^1_{\overline{\mathbb{D}}/X} \to \mathcal{O}_{\overline{\mathbb{D}}}^{\,m-2} \otimes \Omega^2_{\overline{\mathbb{D}}/X} \to \cdots.$$

Note that the isomorphisms (3.3.4) shift the grading:

$$(\mathrm{Gr}^i_K(\mathcal{O}_{\overline{\mathbb{D}}} \otimes (\Omega^j_{\overline{Y}/\overline{S}})))^m \cong (\mathcal{O}_{\overline{\mathbb{D}}} \otimes \Omega^i_{\overline{\mathbb{D}}/X})^{m+i-j} \otimes \Omega^{j-i}_{X/\overline{S}}$$

for $m \in \mathbf{Z}$. In shift notation this is written:

$$\mathrm{Gr}^i_K(\mathcal{O}_{\overline{\mathbb{D}}} \otimes \Omega^j_{\overline{Y}/\overline{S}}) \cong (\Omega^i_{\overline{\mathbb{D}}/X} \otimes \Omega^{j-i}_{X/\overline{S}})(i-j). \tag{3.3.5}$$

There are analogous constructions and statements after tensoring with $\mathcal{O}_{\overline{\mathbb{\Delta}}}$. We claim that the differential of the resulting complex takes $K^i$ to $K^{i+1}$. Indeed, $\nabla'$ acts as zero on $\Omega^\bullet_{\overline{Y}/\overline{S}}$, and, if $\omega \in \mathcal{O}_{\overline{\mathbb{\Delta}}} \otimes \Omega^i_{\overline{\mathbb{D}}/X}$ and $\eta \in \Omega^{j-i}_{\overline{Y}/\overline{S}}$, then

$$\nabla'(\omega \wedge \eta) = \nabla'\omega \wedge \eta \in K^{i+1}(\mathcal{O}_{\overline{\mathbb{\Delta}}} \otimes \Omega^{j+1}_{\overline{Y}/\overline{S}}).$$

Now let $\tilde{K}$ be the dećalée of the filtration $K$ in the sense of [13]. Since $dK^i \subseteq K^{i+1}$, this and its associated graded are simply given by:

$$\begin{aligned} \tilde{K}^i(\mathcal{O}_{\overline{\mathbb{\Delta}}} \otimes \Omega^j_{\overline{Y}/\overline{S}}) &= K^{i+j}(\mathcal{O}_{\overline{\mathbb{\Delta}}} \otimes \Omega^j_{\overline{Y}/\overline{S}}) \\ \mathrm{Gr}^{-i}_{\tilde{K}}(\mathcal{O}_{\overline{\mathbb{\Delta}}} \otimes \Omega^j_{\overline{Y}/\overline{S}}) &= \mathrm{Gr}^{j-i}_K(\mathcal{O}_{\overline{\mathbb{\Delta}}} \otimes \Omega^j_{\overline{Y}/\overline{S}}) \\ &\cong \mathcal{O}_{\overline{\mathbb{\Delta}}} \otimes \Omega^{j-i}_{\overline{\mathbb{D}}/X} \otimes \Omega^i_{X/\overline{S}}. \end{aligned}$$

The differential $\nabla'$ of $\overline{\mathbb{\Delta}}^\bullet(X/Y/S)$ passes over to $\mathrm{Gr}^\bullet_{\tilde{K}}$, and we find an isomorphism of complexes (up to a sign)

$$\mathrm{Gr}^{-i}_{\tilde{K}}(\overline{\mathbb{\Delta}}^\bullet(X/Y/S), \nabla') \cong (\mathcal{O}_{\overline{\mathbb{\Delta}}} \otimes \Omega^\bullet_{\overline{\mathbb{D}}/X} \otimes \Omega^i_{X/S}, \nabla \otimes \mathrm{id})[-i].$$

Then corollary 3.1.8 shows that

$$\mathcal{H}^j(\mathrm{Gr}^{-i}_{\tilde{K}}(\overline{\mathbb{A}}^{\bullet}(X/Y/S),\nabla') \cong \begin{cases} \Omega^i_{X/\overline{S}} & \text{if } j = i \\ 0 & \text{otherwise.} \end{cases} \tag{3.3.6}$$

Thus the spectral sequence of the filtered complex $(\overline{\mathbb{A}}^{\bullet}(X/Y/S),\nabla',\tilde{K})$ degenerates. Statement (2) follows.

Now suppose that we are given a lifting $\tilde{X} \subseteq Y_2$. Parts (a) of statement (3) follows from (2) of Proposition 3.3.1, and parts (b) and (c) from Proposition 3.3.2. Thanks to the shift as we saw in (3.3.5), we find that the computation of Gr involves shifting both this grading and the placement of the complex. This yields isomorphisms:

$$(\mathrm{Gr}^{-i}_{\tilde{K}}(\overline{\mathbb{A}}^{\bullet}_{\lambda}(X/Y/S)),\nabla') \cong (\overline{\mathbb{A}}^{\bullet}_{\lambda-i}/\mathbb{D}\otimes\Omega^i_{X/\overline{S}}, d\otimes \mathrm{id})[-i].$$

Thanks to statement (2) of Proposition 3.3.2, we see that this complex is acyclic unless $i \in \lambda$, and that in that case, its cohomology agrees with that of $\mathrm{Gr}^{-i}_{\tilde{K}}(\overline{\mathbb{A}}^{\bullet}(X/Y/S,\nabla')$. This proves (d). □

**Remark 3.3.4.** The spectral sequence of the filtered complex $(\overline{\mathbb{A}}^{\bullet}(X/Y/S),K)$ used in the proof of Theorem 3.3.3 is related to the conjugate spectral sequence, and the filtration $\tilde{K}$ is related to the canonical filtration. It is convenient to use a somewhat idiosyncratic numbering of the latter when making this comparison. Namely, any complex $C^{\bullet}$ can be endowed with the so-called "trivial" filtration $T^{\bullet}$, defined by:

$$T^iC^j := \begin{cases} 0 & \text{if } i > 0 \\ C^j & \text{if } i \leq 0 \end{cases}.$$

Then if $\tilde{T}^{\bullet}$ is the décalée of of $T^{\bullet}$, we have:

$$\tilde{T}^{-i}C^{\bullet} = [\cdots \to \cdots C^{i-1} \to \mathrm{Ker}(d^i) \to 0 \to \cdots],$$

and so

$$\mathcal{H}^j(\mathrm{Gr}^{-i}_{\tilde{T}}(C^{\bullet})) = \begin{cases} \mathcal{H}^j(C^{\bullet}) & \text{if } i = j \\ 0 & \text{otherwise,} \end{cases}$$

and the associated spectral sequence is the conjugate spectral sequence.

In our case, there is an obvious map of filtered complexes:

$$(\overline{\mathbb{A}}^{\bullet}(X/Y/S),T) \to (\overline{\mathbb{A}}^{\bullet}(X/Y/S),K)$$

hence also maps

$$\begin{array}{rcl} (\overline{\mathbb{A}}^{\bullet}(X/Y/S),\tilde{T}) & \longrightarrow & (\overline{\mathbb{A}}^{\bullet}(X/Y/S),\tilde{K}) \\ \mathrm{Gr}_{\tilde{T}}(\overline{\mathbb{A}}^{\bullet}(X/Y/S)) & \longrightarrow & \mathrm{Gr}_{\tilde{K}}(\overline{\mathbb{A}}^{\bullet}(X/Y/S)) \\ \mathcal{H}^j(\mathrm{Gr}_{\tilde{T}}(\overline{\mathbb{A}}^{\bullet}(X/Y/S))) & \longrightarrow & \mathcal{H}^j(\mathrm{Gr}_{\tilde{K}}(\overline{\mathbb{A}}^{\bullet}(X/Y/S))) \end{array}$$

Thus it follows from the calculation (3.3.6) that the map

$$(\mathbb{A}^{\bullet}(X/Y/S),\tilde{T}) \to (\mathbb{A}^{\bullet}(X/Y/S),\tilde{K})$$

is a filtered quasi-isomorphism.

## 3.4 Calculations and examples

Let $Y$ be crystalline prism and $X \to \overline{Y}$ a regular closed immersion. Our goal here is to give some explicit formulas for the connection

$$\nabla\colon \mathcal{O}_{\overline{\mathbb{\Delta}}_{X(Y)}} \to \mathcal{O}_{\overline{\mathbb{\Delta}}_{X(Y)}} \otimes \Omega^1_{\overline{\mathbb{D}}_{X(Y)/X}}$$

constructed in Theorem 3.1.1, and, given a mod $p^2$ lifting $\tilde{X}$ of $X$ in $Y$, the Sen operator $\Theta_{\tilde{X}}$ acting on $\mathcal{O}_{\overline{\mathbb{\Delta}}_{X(Y)}}$ described in Definition 3.2.3.

We shall need the following notation. Recall that we have $\mathcal{O}_X$-linear maps:

$$\begin{array}{rcl}
\rho_{\mathbb{D}}\colon I_{X/Y} & \to & \pi_{\mathbb{D}*}(\mathcal{O}_{\mathbb{D}}) \\
\rho_{\mathbb{\Delta}}\colon I_{X/Y} & \to & \pi_{\mathbb{\Delta}*}(\mathcal{O}_{\mathbb{\Delta}}) \\
\rho_{\overline{\mathbb{D}}}\colon \mathcal{C}_{X/Y} & \to & \pi_{\mathbb{D}*}(\mathcal{O}_{\overline{\mathbb{D}}}) \\
\rho_{\overline{\mathbb{\Delta}}}\colon \mathcal{C}_{X/Y} & \to & \pi_{\mathbb{\Delta}*}(\mathcal{O}_{\overline{\mathbb{\Delta}}}) \\
\overline{\delta}_{X/Y}\colon F_X^*(\mathcal{C}_{X/Y}) & \to & \mathcal{O}_X,
\end{array}$$

where the maps $\rho$ are deduced from the map in equation (2.2.1) and $\overline{\delta}_{X/Y}$ is the map defined in Proposition 1.7.13. The map $\rho_{\overline{\mathbb{D}}}$ fits into a commutative diagram:

$$\begin{array}{ccccccccc}
0 & \longrightarrow & \mathcal{O}_X & \xrightarrow{[p]} & \mathcal{C}_{X/Y} & \longrightarrow & \mathcal{C}_{X/\overline{Y}} & \longrightarrow & 0 \\
 & & \Big\downarrow \text{id} & & \Big\downarrow \rho_{\overline{\mathbb{D}}} & & \Big\downarrow d\overline{\rho}_{\mathbb{D}} & & \\
0 & \longrightarrow & \mathcal{O}_X & \longrightarrow & \mathcal{O}_{\overline{\mathbb{D}}} & \xrightarrow{d} & \Omega^1_{\overline{\mathbb{D}}/X}, & &
\end{array} \tag{3.4.1}$$

where the map $d\overline{\rho}_{\mathbb{D}}$ induces the isomorphism (3.1.1):

$$\pi^*_{\overline{\mathbb{D}}}(C_{X/\overline{Y}}) \to \Omega^1_{\overline{\mathbb{D}}/X}.$$

**Proposition 3.4.1.** *Let $X \to \overline{Y} \to Y$ be a regular immersion, let $\tilde{X} \to Y_2$ be a lifting of $X$, and let $\Theta_{\tilde{X}}$ and $r_{\tilde{X}}$ be the associated Sen operator (Definition 3.2.3) and $\mathbf{G}_m^\gamma$-action (Proposition 3.2.7). If $x \in I_{\tilde{X}/Y}$ and $\overline{x}$ is its image in $\mathcal{C}_{X/Y}$, then, with the notations above, the following formulas hold:*

$$\begin{array}{rcl}
(\overline{\rho}_{\mathbb{\Delta}}(x))^p & = & \pi^*_{\overline{\mathbb{\Delta}}}(\overline{\delta}_{X/Y}(1 \otimes \overline{x})) \\
\Theta_{\tilde{X}}(\overline{\delta}(\rho_{\mathbb{\Delta}}(x))) & = & -\pi^*_{\overline{\mathbb{\Delta}}}(\overline{\delta}_{X/Y}(1 \otimes \overline{x})) \\
r^*_{\tilde{X}}(\overline{\delta}(\rho_{\mathbb{\Delta}}(x))) & = & \overline{\delta}(\rho_{\mathbb{\Delta}}(x)) - (\log u)\pi^*_{\overline{\mathbb{\Delta}}}(\overline{\delta}_{X/Y}(1 \otimes \overline{x})),
\end{array}$$

*where* $\mathbf{G}_m = \mathrm{Spec}(\mathbf{F}_p[u, u^{-1}])$.

*Proof.* The first equation follows from a calculation we have seen before:

$$\phi_Y(x) \quad = \quad x^p + p\delta_Y(x)$$

$$\begin{aligned} \phi_{\mathbb{A}}(\rho_{\mathbb{A}}(x)) &= p^{p-1}(\rho_{\mathbb{A}}(x))^p + \pi^*_{\mathbb{A}}(\delta_Y(x)) \\ &\equiv \pi^*_{\mathbb{A}}(\delta_Y(x)) \pmod{p\mathcal{O}_{\mathbb{A}}}. \end{aligned}$$

Since $\phi_{\mathbb{A}}(\rho_{\mathbb{A}}(x)) \equiv (\rho_{\mathbb{A}}(x))^p \pmod{p\mathcal{O}_{\mathbb{A}}}$, this proves the first equation.

Recall that the lifting $\tilde{X}$ gives us an isomorphism $\mathcal{O}_{\overline{\mathbb{D}}} \cong S^\bullet(\mathcal{C}_{X/\overline{Y}})$, and Definition 3.2.1 shows that the dilatation Sen operator $\theta_{\tilde{X}}$ acts as the identity on $\mathcal{C}_{X/\overline{Y}} \subseteq \mathcal{O}_{\overline{\mathbb{D}}}$. Thus $\theta_{\tilde{X}}(\rho_{\overline{\mathbb{D}}}(x)) = \rho_{\overline{\mathbb{D}}}(\overline{x})$, and hence also $\Theta_{\tilde{X}}(\rho_{\overline{\mathbb{A}}}(x)) = \rho_{\overline{\mathbb{A}}}(\overline{x})$. Now we calculate, using the properties of $\nabla$ in Definition 3.1.1.

$$\begin{aligned} \Theta_{\tilde{X}}(\overline{\delta}(\rho_{\mathbb{A}}(x))) &:= \nabla_{\theta_{\tilde{X}}}(\overline{\delta}(\rho_{\mathbb{A}}(x))) \\ &= -(\rho_{\overline{\mathbb{A}}}(x))^{p-1}\nabla_{\theta_{\tilde{X}}}(\rho_{\overline{\mathbb{A}}}(x)) \\ &= -(\rho_{\overline{\mathbb{A}}}(x))^{p-1}\theta_{\tilde{X}}(\rho_{\overline{\mathbb{A}}}(x)) \\ &= -(\rho_{\overline{\mathbb{A}}}(x))^{p} \\ &= -\pi^*_{\overline{\mathbb{A}}}(\overline{\delta}_{X/Y}(\overline{x})). \end{aligned}$$

This is the second of our claimed equations. Since $\Theta_{\tilde{X}}$ annihilates the image of $\pi^*_{\overline{\mathbb{A}}}$, we see from the calculations in Example 3.2.10 that the third equation follows. □

**Corollary 3.4.2.** *Continue with the notation of Proposition 3.4.1. If $x \in I_{\tilde{X}/Y}$, the following conditions are equivalent:*

*(1)* $\Theta^*_{\tilde{X}}(\overline{\delta}_{\mathbb{A}}(\rho_{\mathbb{A}}(x))) = 0$*;*

*(2)* $r^*_{\mathbb{A}}(\overline{\delta}_{\mathbb{A}}(\rho_{\mathbb{A}}(x))) = \overline{\delta}_{\mathbb{A}}(\rho_{\mathbb{A}}(x))$*;*

*(3)* $\pi^*_{\overline{\mathbb{A}}}(\overline{\delta}_{X/Y}(1 \otimes x)) = 0$*;*

*(4)* $\delta(x) \in I_{X/Y}$*;*

*(5)* $(\rho_{\mathbb{A}}(x))^p \in p\mathcal{O}_{\mathbb{A}}$*;*

*(6)* $(\rho_{\mathbb{A}}(x))^{[n]} \in \mathcal{O}_{\mathbb{A}}$ *for all* $n$*.*

*Proof.* The equivalence of the first three conditions follows immediately from Proposition 3.4.1. Since $\overline{\delta}_{X/Y}(1 \otimes x)$ is the image of $\delta_Y(x)$ in $\mathcal{O}_X$, it is clear that (4) implies (3), and the converse holds because $\pi_{\overline{\mathbb{A}}}$ is faithfully flat, by Theorem 2.3.3. The equivalence of (5) and (3) follows from the first equation in Proposition 3.4.1, and the equivalence of (5) and (6) follows from Proposition 2.1.4. □

Although the situation of Corollary 3.4.2 is not common, it arises often enough to warrant further study. To prepare for this, we review the standard connection on divided power algebras and explain how to twist it by a $p$-linear endomorphism.

For simplicity, we work in an affine setting. Let $R$ be a ring and $E$ a finitely generated projective $R$-module. Then $\mathbf{V}E := \operatorname{Spec} S^\bullet E$ is a vector bundle over $X := \operatorname{Spec} R$, and $\Omega^1_{\mathbf{V}E/X} \cong S^\bullet E \otimes E$. Under this identification, the exterior differential is the map $d\colon S^\bullet E \to S^\bullet E \otimes E$ sending $e^n$ to $ne^{n-1} \otimes e$. This map extends to an integrable connection on the divided power algebra $\Gamma^\bullet(E)$, sending $e^{[n]}$ to $e^{[n-1]} \otimes e$. This connection is typically denoted by $d$.

**Proposition 3.4.3.** *Let $R$ be a ring in characteristic $p$, let $E$ be a finitely generated and projective $R$-module, and let $\delta\colon E \to E$ be a $p$-linear endomorphism of $E$. Then there is a unique integrable and quasi-nilpotent connection $\nabla$ on the $S^\bullet E$-module $\Gamma^\bullet(E)$ which satisfies the Leibniz rule and such that*

$$\nabla(e^{[n]}) = e^{[n-1]} \otimes e + e^{[n-p]} \otimes \delta(e)$$

*for all $e \in E$.*

*Proof.* Since $\Gamma^\bullet(E)$ is generated as an algebra by elements of the form $e^{[n]}$ for $n \in \mathbf{N}$ and $e \in E$, the uniqueness is clear. To prove the existence, we shall use Roby's characterization [36] of $e \mapsto e^{[n]}$ as the universal polynomial law of degree $n$. If $R'$ is an $R$-algebra, let $E' := R' \otimes_R E$ and let $\delta'$ be the $p$-linear endomorphism of $E$ induced by $\delta$. By [36, III.3], the map $R' \otimes_R \Gamma^\bullet_R(E) \to \Gamma^\bullet_{R'}(E')$ is an isomorphism, and we have a set-theoretic map

$$\delta'_n\colon E' \to \Gamma^{n-p}(E') \otimes_{R'} E' : e' \mapsto {e'}^{[n-p]} \otimes \delta(e')$$

for all $e' \in E'$. This family of maps is compatible with homomorphisms $R' \to R''$. Taking $R' := R[t]$, we see that, if $e \in E$, then

$$te \mapsto (te)^{[n-p]} \otimes \delta(te) = t^{n-p} e^{[n-p]} \otimes t^p \delta(e) = t^n e^{[n-p]} \otimes \delta(e).$$

Thus this map is a "polynomial law of degree $n$" in the sense of Roby. Then by the universality of $e \mapsto e^{[n]}$, there is a unique $R$-linear map:

$$h_n\colon \Gamma^n(E) \to \Gamma^{n-p}(E) \otimes E$$

such that $h_n(e^{[n]}) = e^{[n-p]} \otimes \delta(e)$ for all $e \in E$.

**Lemma 3.4.4.** *Let*

$$h\colon \Gamma^\bullet(E) \to \Gamma^\bullet(E) \otimes_R E \cong \Gamma^\bullet(E) \otimes_{S^\bullet E} \Omega^1_{\mathbf{V}E}$$

*be the sum of the maps $h_n$. Then $h$:*

*(1) satisfies the Leibniz rule with respect to multiplication:*

$$h(fg) = h(f)g + fh(g) \quad \text{for } f, g \in \Gamma^\bullet(E);$$

*(2) is a quasi-nilpotent Higgs field on the $S^\bullet(E)$-module $\Gamma^\bullet(E)$;*

*(3) anti-commutes with the connection $d$, i.e., $hd + dh = 0$.*

*Proof.* To prove (1), let $k$ be a positive integer, let $(e_1, \ldots, e_k)$ be a sequence in $E$, and let $R'$ be the polynomial algebra $R[t_1, \ldots, t_k]$. Then if $h' := \mathrm{id}_{R'} \otimes h$, we have:

$$\begin{aligned} h'((t_1e_1 + \cdots t_ke_k)^{[n]}) &= (t_1e_1 + \cdots t_ke_k)^{[n-p]} \otimes \delta(t_1e_1 + \cdots + t_ke_k) \\ &= \sum_{I:|I|=n-p} t_1^{I_1}e_1^{[I_1]} \cdots t_k^{I_k}e_k^{[I_k]} \otimes (t_1^p\delta(e_1) + \cdots + t_k^p\delta(e_k)) \\ &= \sum_{j=1}^{k} \sum_{I:|I|=n, I_j \geq p} t_1^{I_1} \cdots t_k^{I_k} e_1^{[I_1]} \cdots e_j^{[I_j-p]} \cdots e_k^{[I_k]} \otimes \delta(e_j) \\ &= \sum_{j=1}^{k} \sum_{I:|I|=n} t_1^{I_1} \cdots t_k^{I_k} e_1^{[I_1]} \cdots\cdots e_k^{[I_k]} h(e_j^{[I_j]}). \end{aligned}$$

On the other hand:

$$\begin{aligned} (t_1e_1 + \cdots t_ke_k)^{[n]} &= \sum_{I:|I|=n} t_1^{I_1}e_1^{[I_1]} \cdots t_k^{I_k}e_k^{[I_k]} \\ h'((t_1e_1 + \cdots t_ke_k)^{[n]}) &= \sum_{I:|I|=n} h'(t_1^{I_1}e_1^{[I_1]} \cdots t_k^{I_k}e_k^{[I_k]}) \\ &= \sum_{I:|I|=n} t_1^{I_1} \cdots t_k^{I_k} h(e_1^{[I_1]} \cdots e_k^{[I_k]}) \end{aligned}$$

Comparing coefficients of $t_1^{I_1} \cdots t_k^{I_k}$, we get

$$h(e_1^{[I_1]} \cdots e_k^{[I_k]}) = \sum_{j=1}^{k} e_1^{[I_1]} \cdots\cdots e_k^{[I_k]} h(e_j^{[I_j]}).$$

Now suppose that $1 \leq i \leq k$ and let $f := e_1^{[I_1]} \cdots e_i^{[I_i]}$ and $g := e_{i+1}^{[I_{i+1}]} \cdots e_k^{[I_k]}$. Then the equation above can be rewritten::

$$\begin{aligned} h(fg) &= \sum_{j=1}^{i} e_1^{[I_1]} \cdots e_k^{[I_k]} h(e_j^{[I_j]}) + \sum_{j=i+1}^{k} e_1^{[I_1]} \cdots e_k^{[I_k]} h(e_j^{[I_j]}) \\ &= h(f)g + fh(g). \end{aligned}$$

This shows that the Leibniz rule holds for the product of any two divided power monomials. Since every element of $\Gamma^\bullet(E)$ is a sum of such, it holds in general. This proves (1).

It is clear from the definition that $h$ annihilates the image of $E \to \Gamma^\bullet(E)$, and it follows from the Leibniz rule that it also annihilates the image of $S^\bullet(E)$. Thus it is in fact $S^\bullet(E)$-linear. To show that it is a Higgs field, we check that the composition $h \wedge h \colon \Gamma^\bullet(E) \to \Gamma^\bullet(E) \otimes E \to \Gamma^\bullet(E) \otimes \Lambda^2 E$ vanishes. First note that if $e \in E$, then

$$h(e^{[n]}) = e^{[n-p]} \otimes \delta(e)$$

$$\begin{aligned}(h \wedge h)(e^{[n]}) &= h(e^{[n-p]} \wedge \delta(e)) \\ &= e^{[n-2p]} \otimes (\delta(e) \wedge \delta(e)) \\ &= 0.\end{aligned}$$

Next, if $f, g \in \Gamma^\bullet(E)$, then

$$\begin{aligned}h(fg) &= h(f)g + fh(g) \\ &= gh(f) + fh(g) \\ (h \wedge h)(fg) &= h(gh(f)) + h(fh(g)) \\ &= h(g) \wedge h(f) + g(h \wedge h)(f) + h(f) \wedge h(g) + f(h \wedge h)(g) \\ &= g(h \wedge h)(f) + f(h \wedge h)(g).\end{aligned}$$

Thus $h \wedge h$ is again a derivation, and since it annihilates the generators of the algebra $\Gamma^\bullet(E)$, it vanishes.

Next we check $hd + dh$ vanishes. We compute:

$$\begin{aligned}(hd + dh)(e^{[n]}) &= h(e^{[n-1]} \otimes e) + d(e^{[n-p]} \otimes \delta(e)) \\ &= e^{[n-p-1]} \otimes \delta(e) \wedge e + e^{[n-p-1]} \otimes e \wedge \delta(e) \\ &= 0.\end{aligned}$$

Since $hd + dh$ is also a derivation, it therefore annihilates all of $\Gamma^\bullet(E)$. □

The proposition now follows easily. Namely, if $\nabla\colon \Gamma^\bullet(E) \to E \otimes \Gamma^\bullet(E)$ is $d + h$, then $\nabla$ satisfies the Leibniz rule because $d$ and $h$ do, $\nabla$ is a connection because $d$ is a connection and $h$ is linear, and $\nabla$ is integrable because $\nabla^2 = (d+h)^2 = d^2 + dh + hd + h^2 = 0$. Finally it is apparent from the formulas that $\nabla$ is quasi-nilpotent, since its action decreases the PD-filtration level. □

As a simple example, suppose that $E$ is free with basis $\{e\}$. Then the lemma implies that the map $h\colon e^{[n]} \mapsto e^{[n-p]} \otimes e$ satisfies the Leibniz rule. In particular, if $m$ and $n$ are natural numbers, then:

$$\begin{aligned}h(e^{[m]}e^{[n]}) &= h(e^{[m]})e^{[n]} + e^{[m]}h(e^{[n]}) \\ \binom{m+n}{m} h(e^{[m+n]}) &= e^{[m-p]}e^{[n]} \otimes e + e^{[m]}e^{[n-p]} \otimes e \\ \binom{m+n}{m} e^{[m+n-p]} \otimes e &= \binom{m+n-p}{n} e^{[m-p+n]} \otimes e + \binom{m+n-p}{m} e^{[m-p+n]} \otimes e.\end{aligned}$$

If $m + n \geq p$, this leads to the following congruence, which (of course) can be checked directly (but for us, not painlessly), and which we will need later:

$$\binom{m+n}{m} \equiv \binom{m+n-p}{m} + \binom{m+n-p}{n} \pmod p \quad \text{if } m+n \geq p. \tag{3.4.2}$$

We shall use Proposition 3.4.3 to describe the connection $\nabla$ on $\overline{\mathbb{A}}_X(Y)$ in some special cases.

**Proposition 3.4.5.** *Let $Y$ be a crystalline prism, let $X \to \overline{Y}$ be a regular closed immersion, and let $x$ be a local section of $I_{X/Y}$ such that $\delta(x) \in I_{X/Y}$, so that $(\rho_{\mathbb{A}}(x))^{[n]} \in \mathcal{O}_{\overline{\mathbb{A}}_{X/Y}}$ for all $n$, as we saw in Lemma 2.3.16. Then for $n \in \mathbf{N}$:*

$$\nabla(\overline{\rho}_{\mathbb{A}}(x)^{[n]}) = \overline{\rho}_{\mathbb{A}}(x)^{[n-1]} d\overline{\rho}_{\mathbb{D}}(x) - \overline{\rho}_{\mathbb{A}}(x)^{[n-p]} d\rho_{\mathbb{D}}(\overline{\delta}(x)).$$

*In particular, if $\delta(x) \in p^2\mathcal{O}_Y$, then*

$$\nabla(\overline{\rho}_{\mathbb{A}}(x)^{[n]}) = \overline{\rho}_{\mathbb{A}}(x)^{[n-1]} d\overline{\rho}_{\mathbb{D}}(x).$$

*for all $n$.*

*Proof.* For simplicity of notation, let

$$\tilde{x} := \overline{\rho}_{\mathbb{A}}(x), \quad d\tilde{x} := d\overline{\rho}_{\mathbb{D}}(x), \quad \tilde{\delta}(x) := d\overline{\rho}_{\mathbb{D}}(\overline{\delta}(x)).$$

Thus our claim is that

$$\nabla(\tilde{x}^{[n]}) = \tilde{x}^{[n-1]} d\tilde{x} - \tilde{x}^{[n-p]}\tilde{\delta}(x).$$

(compare with 3.4.3). We first check the formula when $n = p^m$ for some $m \in \mathbf{N}$, by induction on $m$. If $m = 0$, then $n = 1$, and the formula is true from the defining property (1) of the definition of $\nabla$ in Theorem 3.1.1. If $m = 1$, we have, from Lemma 2.3.16 and property (3) of the definition of $\nabla$:

$$\begin{aligned}
\tilde{x}^{[p]} &= \overline{\delta}_{\mathbb{A}}(\rho_{\mathbb{A}}(x)) - \overline{\rho}_{\mathbb{A}}(\delta(x)) \\
\nabla(\tilde{x}^{[p]}) &= \nabla(\overline{\delta}_{\mathbb{A}}(\rho_{\mathbb{A}}(x))) - \nabla(\overline{\rho}_{\mathbb{A}}(\delta(x))) \\
&= -\tilde{x}^{p-1} d\overline{\rho}_{\mathbb{D}}(x) - d\overline{\rho}_{\mathbb{D}}(\delta(x)) \\
&= \tilde{x}^{[p-1]} d\tilde{x} - \tilde{\delta}(x)
\end{aligned}$$

as claimed.

Assume the formula for $m$, with $m \geq 1$. Lemma 2.3.16 tells us that

$$\tilde{x}^{[p^{m+1}]} = \overline{\delta}_{\mathbb{A}}(\tilde{x}^{[p^m]}) + e^p.$$

where the term $e^p$ can only occur if $p = 2$. Then using the defining property of $\nabla$, and the induction assumption, we compute

$$\begin{aligned}
\nabla(\tilde{x}^{[p^{m+1}]}) &= \nabla(\overline{\delta}_{\mathbb{A}}(\tilde{x})^{[p^m]}) + \nabla(e^p) \\
&= -(\tilde{x}^{[p^m]})^{p-1}\nabla(\tilde{x}^{[p^m]}) + pe^{p-1}\nabla(e) \\
&= -(\tilde{x}^{[p^m]})^{p-1}\left(\tilde{x}^{[p^m-1]} d\tilde{x} - \tilde{x}^{[p^m-p]}\tilde{\delta}(x)\right) \\
&= -\left(\frac{\tilde{x}^{p^m}}{p^m!}\right)^{p-1}\left(\frac{\tilde{x}^{p^m-1}}{(p^m-1)!} d\tilde{x} - \frac{\tilde{x}^{p^m-p}}{(p^m-p)!}\tilde{\delta}(x)\right) \\
&= -\frac{(p^{m+1}-1)!\tilde{x}^{[p^{m+1}-1]}}{((p^m)!)^{p-1}(p^m-1)!} d\tilde{x} + \frac{(p^{m+1}-p)!\tilde{x}^{[p^{m+1}-p]}}{((p^m)!)^{p-1}(p^m-p)!}\tilde{\delta}(x)
\end{aligned}$$

Strictly speaking, these computations don't make sense in characteristic $p$, and should be regarded as "formal." We shall use Lemma 1.6.7 to justify the calculation and to compute these terms. First of all,

$$\begin{aligned}
\mathrm{ord}_p((p^m-1)!) &= \frac{p^m-1}{p-1} - m \\
\mathrm{ord}_p((p^m!)^{p-1}) &= p^m - 1 \\
\mathrm{ord}_p((p^m-1)!) + \mathrm{ord}_p((p^m!)^{p-1}) &= \frac{p^{m+1}-1}{p-1} - (m+1) \\
&= \mathrm{ord}_p((p^{m+1}-1)!),
\end{aligned}$$

so the coefficient of the first term in the previous equation is a $p$-adic unit. A similar calculation shows that the same is true for the second term. Furthermore, since $p^m! = p^m(p^m-1)!$, we have $U_p((p^m-1)!) = U_p((p^m)!)$, and, by Lemma 1.6.7, $U_p((p^m)!) = (-1)^m$. It follows that

$$\frac{U_p((p^{m+1}-1)!)}{U_p((p^m)!)^{p-1}U_p((p^m-1)!} = \frac{(-1)^{m+1}}{(-1)^{m(p-1)}(-1)^m} = -1.$$

Similarly, $U_p((p^{m+1}-p)!) = U_p((p^m-1)!) = U_p((p^m)!) = (-1)^m$, so

$$\frac{U_p((p^{m+1}-p)!)}{U_p((p^m)!)^{p-1}U_p((p^m-p)!} = \frac{(-1)^m}{(-1)^{m(p-1)}(-1)^{m-1}} = -1.$$

Thus, the expression above reduces to

$$\nabla(\tilde{x}^{[p^{m+1}]}) = \tilde{x}^{[p^{m+1}]}d\tilde{x} - \tilde{x}^{[p^{m+1}-p]})\tilde{\delta}(x)$$

as required.

For the general case, let:

$$\tilde{\nabla}(\tilde{x}^{[n]}) := \tilde{x}^{[n-1]}d\tilde{x} + \tilde{x}^{[n-p]}\tilde{\delta}(x).$$

Then

$$\begin{aligned}
\tilde{\nabla}(\tilde{x}^{[m]}\tilde{x}^{[n]}) &= \tilde{\nabla}\left(\binom{m+n}{n}\tilde{x}^{[m+n]}\right) \\
&= \binom{m+n}{n}\left(\tilde{x}^{[m+n-1]}d\tilde{x} + \tilde{x}^{[m+n-p]}\tilde{\delta}(x)\right).
\end{aligned}$$

On the other hand:

$$\begin{aligned}
\tilde{\nabla}(\tilde{x}^{[m]})\tilde{x}^{[n]} + \tilde{x}^{[m]}\tilde{\nabla}(\tilde{x}^{[n]}) &= (\tilde{x}^{[m-1]}d\tilde{x} + \tilde{x}^{[m-p]}\tilde{\delta}(x))\tilde{x}^{[n]} + \tilde{x}^{[m]}(\tilde{x}^{[n-1]}d\tilde{x} + \tilde{x}^{[n-p]}\tilde{\delta}(x)) \\
&= (\tilde{x}^{[m-1]}\tilde{x}^{[n]} + x^{[m]}\tilde{x}^{[n-1]})d\tilde{x} + (\tilde{x}^{[m-p]}\tilde{x}^{[n]} + \tilde{x}^{[m]}\tilde{x}^{[n-p]})\tilde{\delta}(x) \\
&= \left(\binom{m+n-1}{n} + \binom{m+n-1}{m}\right)\tilde{x}^{[m+n-1]}d\tilde{x} \quad + \\
&\qquad \left(\binom{m+n-p}{n} + \binom{m+n-p}{m}\right)\tilde{x}^{[m+n-p]}\tilde{\delta}(x)
\end{aligned}$$

$$= \binom{m+n}{n}\tilde{x}^{[m+n-1]}d\tilde{x} + \binom{m+n}{m}\tilde{x}^{[m+n-p]}\tilde{\delta}(x),$$

by Pascal's formula and equation (3.4.2).

Thus $\nabla$ and $\tilde{\nabla}$ are both derivations, and they agree on elements of the form $x^{[n]}$ whenever $n$ is a power of $p$. Since the divided power algebra is generated by such elements, this implies that $\tilde{\nabla}$ and $\nabla$ agree.

Now suppose that $\delta(x) \in p^2\mathcal{O}_Y$. Since $\rho_{\mathbb{D}}(\delta(x)) \in p\mathcal{O}_{\mathbb{D}}$, it follows that $d\overline{\rho}_{\mathbb{D}}(x) = 0$, and the formula simplifies as claimed. □

Proposition 3.4.5 immediately implies the following.

**Corollary 3.4.6.** *Suppose that $X \subseteq \overline{Y}$ is as in Proposition 3.4.5 and admits an $F_Y$-invariant lift $\tilde{X} \subseteq Y_2$. Make the identification:*

$$\overline{\mathbb{A}}_X(Y) \cong \mathbb{N}^{\gamma}_{X/\overline{Y}} = \mathrm{Spec}_X(\Gamma^{\bullet}(\mathcal{C}_{X/\overline{Y}})),$$

*as explained in Theorem 2.3.13 and Proposition 3.2.1. Then the connection $\nabla$ on $\mathcal{O}_{\overline{\mathbb{A}}_X(Y)}$ is given as in Proposition 3.4.3, with $E = \mathcal{C}_{X/\overline{Y}}$ and $\delta = \overline{\delta}_{\tilde{X}/\overline{Y}}$ (see (3) of Proposition 1.7.13). If in fact $\tilde{X}$ admits an $F_Y$-invariant lift in $Y_3$, then $\overline{\delta}_{\tilde{X}/\overline{Y}}$ vanishes, and $\nabla$ is the standard connection on $\Gamma^{\bullet}(\mathcal{C}_{X/\overline{Y}})$.* □

**Remark 3.4.7.** We can use these calculations to give another proof that the action of $\mathbb{N}^{\gamma}_{X/\overline{Y}}$ on $\overline{\mathbb{A}}_X(Y)$ defined by $\nabla$ is torsorial, without the use of Theorem 3.1.5. Working fpqc locally on $X$, we may assume that there is a regular sequence $(p, x_1, \ldots, x_r)$ generating $I_{X/Y}$. Corollary 2.3.12 tells us that, after a $p$-completely flat covering, we can assume that $\delta(x_i) \in p^2\mathcal{O}_Y$ for all $i$. Then the lifting of $X$ in $Y$ defined by $(x_1, \ldots, x_r)$ is $F_3$-aligned, and Theorem 2.3.13 identifies $\overline{\mathbb{D}}_X(Y)$ with $\mathbb{N}_{X/\overline{Y}}$ and $\overline{\mathbb{A}}_X(Y)$ with $\mathbb{N}^{\gamma}_{X/\overline{Y}}$. Furthermore, Corollary 3.4.6 tells us that the connection on $\mathcal{O}_{\overline{\mathbb{A}}_X(Y)}$ is just the standard connection on $\mathcal{O}_{\mathbb{N}^{\gamma}_{X/\overline{Y}}}$, and hence that the corresponding action of $\mathbb{N}^{\gamma}_{X/\overline{Y}}$ is the standard action by translation.

The following result, a follow-up to Proposition 3.4.1, describes the action of the Sen operator when our lifting $\tilde{X}$ is invariant under $F_Y$.

**Proposition 3.4.8.** *In the situation of Proposition 3.4.5, suppose that $\tilde{X}$ is a flat lifting of $X$ in $Y_2$, that $x \in I_{\tilde{X}/Y}$, and that $\delta(x) \in I_{X/Y}$, so that $\rho_{\mathbb{A}}(x)$ and $\rho_{\mathbb{A}}(\delta(x))$ belong to $\mathcal{O}_{\mathbb{A}}$. Let $\tilde{\delta}_2(x)$ be the image of $\delta(x)$ under the map*

$$I_{X/Y} \longrightarrow \mathcal{C}_{X/\overline{Y}} \xrightarrow{s_{\tilde{X}}} \mathcal{C}_{X/Y} \xrightarrow{\overline{\rho}_{\overline{\mathbb{A}}}} \mathcal{O}_{\overline{\mathbb{A}}},$$

*where $s_{\tilde{X}}$ is the splitting defined by $\tilde{X}$. Then the following statements are verified.*

*(1) $\overline{\rho}_{\mathbb{A}}(\delta(x)) - \tilde{\delta}_2(x)$ belongs to the image of the map $\pi^*_{\overline{\mathbb{A}}}\colon \mathcal{O}_X \to \mathcal{O}_{\overline{\mathbb{A}}_X(Y)}$.*

*(2) The diffraction operators act as follows:*

$$\begin{aligned}
\Theta_{\tilde{X}}(\tilde{\delta}_2(x)) &= \tilde{\delta}_2(x) \\
\Theta_{\tilde{X}}(\overline{\delta}_{\mathbb{\Delta}}(x)) &= \Theta_{\tilde{X}}(\tilde{\delta}(x)) \\
&= \tilde{\delta}_2(x) \\
\Theta_{\tilde{X}}((\overline{\rho}_{\mathbb{\Delta}}(x))^{[p]}) &= -\tilde{\delta}_2(x) \\
r^*_{\tilde{X}}((\overline{\rho}_{\mathbb{\Delta}}(x))^{[p]}) &= (\overline{\rho}_{\mathbb{\Delta}}(x))^{[p]} - (u-1)\tilde{\delta}_2(x) \\
\Theta_{\tilde{X}}((\overline{\rho}_{\mathbb{\Delta}}(x))^{[n]}) &= n(\overline{\rho}_{\mathbb{\Delta}}(x))^{[n]} - (\overline{\rho}_{\mathbb{\Delta}}(x))^{[n-p]}\tilde{\delta}_2(x)).
\end{aligned}$$

*(3)* $\Theta^{nil}_{\tilde{X}}$ *annihilates* $(\overline{\rho}_{\mathbb{\Delta}}(x))^{[n]}$ *for all* $n$.

*Proof.* Since we are assuming that $\delta(x)$ belongs to $I_{X/Y}$ we can write $\delta(x) = pf + g$ with $f \in \mathcal{O}_Y$ and $g \in I_{\tilde{X}/Y}$. Then $\delta(x)$ maps to the image of $(f, g)$ in $\mathcal{C}_{X/Y} \cong \mathcal{O}_X \oplus \mathcal{C}_{X/\overline{Y}}$, so $\tilde{\delta}_2(x) = \overline{\rho}_{\mathbb{\Delta}}(g)$. Furthermore, $\rho_{\mathbb{D}}(g)$ belongs to the ideal of the embedding of $X$ in $\mathbb{D}_X(Y)$ determined by the lifting $\tilde{X}$. Since

$$\rho_{\mathbb{\Delta}}(\delta(x)) = \rho_{\mathbb{\Delta}}(pf) + \rho_{\mathbb{\Delta}}(g) = \pi^*_{\mathbb{\Delta}}(f) + \rho_{\mathbb{\Delta}}(g),$$

we see that $\overline{\rho}_{\mathbb{\Delta}}(\delta(x)) - \tilde{\delta}_2(x)$ is the image of $f$ in $\mathcal{O}_{\overline{\mathbb{\Delta}}_X(Y)}$, proving statement (1).

Using the isomorphism $\mathcal{O}_{\overline{\mathbb{D}}} \cong S^\bullet(\mathcal{C}_{X/\overline{Y}})$ provided by $\tilde{X}$, we find that $\pi^*_{\overline{\mathbb{D}}}(f)$ has degree zero and $\rho_{\overline{\mathbb{D}}}(g)$ has degree one, hence the former is annihilated by the Sen operator and the latter is fixed. We conclude

$$\Theta_{\tilde{X}}(\overline{\rho}_{\mathbb{\Delta}}(\delta(x))) = \Theta_{\tilde{X}}(\overline{\rho}_{\mathbb{\Delta}}(g)) = \tilde{\delta}_2(x). \tag{3.4.3}$$

Since $\delta(x) \in I_{X/Y}$, it follows that $\overline{\delta}_{X/Y}(1 \otimes \overline{x}) = 0$, and then, by Proposition 3.4.1, that $\overline{\delta}(\rho_{\mathbb{\Delta}}(x))$ is killed by $\Theta_{\tilde{X}}$. Then we can use Proposition 2.3.16 to compute:

$$\begin{aligned}
(\overline{\rho}_{\mathbb{\Delta}}(x))^{[p]} &= \overline{\delta}_{\mathbb{\Delta}}(\rho_{\mathbb{\Delta}}(x)) - \overline{\rho}_{\mathbb{\Delta}}(\delta(x)) \\
\Theta_{\tilde{X}}((\overline{\rho}_{\mathbb{\Delta}}(x))^{[p]}) &= \Theta_{\tilde{X}}(\overline{\delta}_{\mathbb{\Delta}}(\rho_{\mathbb{\Delta}}(x))) - \Theta_{\tilde{X}}(\overline{\rho}_{\mathbb{\Delta}}(\delta(x))) \\
&= -\tilde{\delta}_2(x),
\end{aligned}$$

and the formula for $r^*_{\tilde{X}}((\overline{\rho}_{\mathbb{\Delta}}(x))^{[p]})$ follows because

$$r^*_{\tilde{X}}(\overline{\delta}_{\mathbb{\Delta}}(\rho_{\mathbb{\Delta}}(x))) = \overline{\delta}_{\mathbb{\Delta}}(\rho_{\mathbb{\Delta}}(x)) \text{ and } r^*_{\tilde{X}}(\overline{\rho}_{\mathbb{\Delta}}(\delta(x))) = u\overline{\rho}_{\mathbb{\Delta}}(\delta(x)).$$

Proposition 3.4.5 tells us that

$$\begin{aligned}
\nabla_{\theta_{\tilde{X}}}((\overline{\rho}_{\mathbb{\Delta}}(x))^{[n]}) &= (\overline{\rho}_{\mathbb{\Delta}}(x))^{[n-1]}\Theta_{\tilde{X}}(\overline{\rho}_{\mathbb{\Delta}}(x)) - (\overline{\rho}_{\mathbb{\Delta}}(x))^{[n-p]}\Theta_{\tilde{X}}(\rho_{\mathbb{\Delta}}(\overline{\delta}(x))) \\
&= (\overline{\rho}_{\mathbb{\Delta}}(x))^{[n-1]}\overline{\rho}_{\mathbb{\Delta}}(x) - (\overline{\rho}_{\mathbb{\Delta}}(x))^{[n-p]}\tilde{\delta}_2(x) \\
&= n(\overline{\rho}_{\mathbb{\Delta}}(x))^{[n]} - (\overline{\rho}_{\mathbb{\Delta}}(x))^{[n-p]}\tilde{\delta}_2(x).
\end{aligned}$$

This completes the proof of statement (2). Note also that when $n = p$, it reaffirms our previous computation.

To prove statement (3), let $\Psi$ be the value of the $p$-curvature of $\nabla$ at the Euler vector field $\theta_{\tilde{X}}$, and recall from Equation (3.2.5) that $\Psi = \Theta^p_{\tilde{X}} - \Theta_{\tilde{X}}$ and from Remark 3.2.6 that $\Theta^{nil}_{\tilde{X}} = -\Psi - \Psi^p - \Psi^{p^2} - \Psi^{p^3} - \cdots$. Thus it will suffice to prove that $\Psi$ annihilates $\overline{\rho}(x)^{[n]}$ for all $n$. Since $\Psi$ is a derivation, it is enough to treat the case in which $n$ is a power of $p$.

To save ink, we sometimes abuse notation by writing $\Psi$ and $\Theta$ for their respective compositions with reduction modulo $p$. Since $\nabla$ is just the standard connection on $\mathcal{O}_{\overline{\mathbb{D}}}$, it follows that the $p$-curvature vanishes on $\mathcal{O}_{\overline{\mathbb{D}}}$, and in particular that $\Psi(\tilde{x}) = 0$. Now using statement (2), we calculate:

$$\begin{aligned} \Psi(\tilde{x}^{[p]}) &= \Theta^p(\tilde{x}^{[p]}) - \Theta(\tilde{x}^{[p]}) \\ &= -\Theta^{p-1}(\tilde{\delta}_2(x)) + \tilde{\delta}_2(x) \\ &= -\delta_2(x) + \tilde{\delta}_2(x) \\ &= 0. \end{aligned}$$

We proceed by induction. Recall that Lemma 2.3.16 tells us that, if $m > 0$,

$$\tilde{x}^{[p^{m+1}]} = \overline{\delta}_{\mathbb{\Delta}}(\tilde{x}^{[p^m]}) + e^p,$$

where the term $e^p$ can only occur if $p = 2$. Then, using Proposition 3.1.4, we compute:

$$\begin{aligned} \Psi(\tilde{x}^{[p^{m+1}]}) &= \Psi(\overline{\delta}_{\mathbb{\Delta}}(\tilde{x}^{[p^m]})) + \Psi(e^p) \\ &= -(\tilde{x}^{[p^m]})^{p-1}\Psi(\tilde{x}^{[p^m]}) + (\Theta(\tilde{x}^{[p^m]}))^p \end{aligned}$$

The first term vanishes by the induction assumption, and the second because $\Theta(\tilde{x}^{[p^m]})$ belongs to the PD-ideal of $\mathcal{O}_{\overline{\mathbb{\Delta}}}$, as follows from the formulas in (2). □

**Corollary 3.4.9.** *In the situation of Proposition 3.4.8, the following are equivalent.*

*(1) $\tilde{X}$ is invariant under $F_Y$.*

*(2) The action $r_{\tilde{X}}$ of $\mathbf{G}_m^\gamma$ on $\overline{\mathbb{\Delta}}_X(Y)$ factors through $\mu_p$.*

*(3) The Sen operator $\Theta_{\tilde{X}}$ on $\mathcal{O}_{\overline{\mathbb{\Delta}}_X(Y)}$ is semi-simple.*

*Proof.* The equivalence of (2) and (3) follows from the generalities in Corollary 3.2.9. Since $\tilde{X}$ is flat modulo $p^2$ it follows that it is invariant under $F_Y$ if and only if $\delta_Y$ maps $I_{\tilde{X}/Y}$ to $I_{X/Y}$. If this fails, there is some $x \in I_{\tilde{X}/Y}$ such that the image of $\delta_Y(x)$ in $\mathcal{O}_X$ is not zero. Since $\pi_{\overline{\mathbb{\Delta}}}\colon \overline{\mathbb{\Delta}} \to X$ is faithfully flat, it follows that $\pi^*_{\overline{\mathbb{\Delta}}}(\delta_{X/Y}(\overline{x}))$ is not zero, and then Proposition 3.4.1 shows that

$$r^*_{\tilde{X}}(\overline{\rho}_{\mathbb{\Delta}}(x)) = r^*_{\tilde{X}}(\overline{\rho}_{\mathbb{\Delta}}(x)) + (\log u)\pi^*_{\overline{\mathbb{\Delta}}}(\delta_{X/Y}(\overline{x})).$$

Thus $r_{\tilde{X}}$ does not factor through $\mu_p$ in this case.

Suppose on the other hand that $\tilde{X}$ is $F_Y$-invariant. Then Theorem 2.3.13 tells us that $\overline{\mathbb{\Delta}}_X(Y)$ identifies with the PD-envelope of the zero section of $\mathbb{N}_{X/\overline{Y}}$. Thus it is enough to show that $\Theta_{nil}$ annihilates $(\overline{\rho}(x))^{[n]}$ for every $x \in I_{\tilde{X}/Y}$ and every $n > 0$. This follows from Proposition 3.4.8. □

# 4 Prismatic connections and cohomology

Our goal here is to show that, if $S$ is a crystalline prism, if $X/\overline{S}$ is smooth, and if $X \to Y$ is an embedding of $X$ is a $p$-completely smooth crystalline prism over $S$, then the prismatic envelope of $X$ in $Y$ carries a natural $p$-connection and that the resulting $p$-de Rham complex, viewed as an object in the derived category, is independent of the choice of the embedding. A similar result will hold with a suitable notion of coefficients. This will allow us to give a relatively concrete construction of prismatic cohomology, which we will reinterpret using the prismatic topos in §6.

## 4.1 Dilatations and $p$-connections

The $p$-differential $d'$ of the $p$-de Rham complex (1.2.2) has a geometric interpretation, analogous to the interpretation (1.1.2) of the usual differential $d$. Suppose that $Y/S$ is a smooth morphism of $p$-torsion free schemes (or formal schemes), let $\mathbb{D}_Y(1) := \mathbb{D}_{\overline{Y}}(Y(1))$, and consider the commutative diagram

$$\begin{array}{ccccc} & & \mathbb{D}_Y(1) & & \\ & \nearrow^{\tilde{\Delta}_{Y/S}} & \downarrow{\scriptstyle \pi_{\mathbb{D}}} & \searrow^{\tilde{p}_i} & \\ Y & \xrightarrow{\Delta_{Y/S}} & Y(1) & \xrightarrow{p_i} & Y, \end{array}$$

where $\tilde{\Delta}_{Y/S}$ is the map induced from the universal property of the dilatation as discussed in Proposition 2.2.14 and, for $i = 1, 2$, $p_i$ is the natural projection. Let $J$ (resp. $\tilde{J}$) be the ideal of the locally closed immersion $\Delta_{Y/S}$ (resp. $\tilde{\Delta}_{Y/S}$). Proposition 2.2.14 implies that $\rho$ induces an isomorphism $\tilde{\rho}\colon J/J^2 \to \tilde{J}/\tilde{J}^2$, fitting into a diagram:

$$\begin{array}{ccccccc} \mathcal{O}_Y & \xrightarrow{p_2^* - p_1^*} & J & \xrightarrow{[\,]} & J/J^2 & \xrightarrow{\cong} & \Omega^1_{Y/S} \\ & & \downarrow{\scriptstyle \rho} & & \downarrow{\scriptstyle \tilde{\rho}} & \nearrow_{\cong} & \\ \mathcal{O}_Y & \xrightarrow{\tilde{p}_2^* - \tilde{p}_1^*} & \tilde{J} & \xrightarrow{[\,]} & \tilde{J}/\tilde{J}^2. & & \end{array}$$

The map $d\colon \mathcal{O}_Y \to \Omega^1_{Y/S}$ is the composition of the top arrows. Letting $d'$ be the composition along the bottom, we find, for $f \in \mathcal{O}_Y$:

$$\begin{aligned} [(\tilde{p}_2^*(f) - \tilde{p}_1^*(f)] &= [\pi_{\mathbb{D}}^*(p_2^*(f)) - \pi_{\mathbb{D}}^*(p_1^*(f))] \\ &= [\pi_{\mathbb{D}}^*(p_2^*(f) - p_1^*(f))] \\ &= [p(\rho(p_2^*(f) - p_1^*(f)))] \\ &= p\tilde{\rho}[p_2^*(f) - p_1^*(f)], \qquad \text{hence} \\ d'f &= pdf. \end{aligned}$$

If we continue to use the same identification $J/J^2 \cong \tilde{J}/\tilde{J}^2 \cong \Omega^1_{Y/S}$ and its exterior powers, we find a commutative diagram:

$$\begin{array}{ccccccc} \mathcal{O}_Y & \xrightarrow{d'} & \tilde{J}/\tilde{J}^2 & \xrightarrow{d'} & \Lambda^2 \tilde{J}/\tilde{J}^2 & \cdots \\ \uparrow{\scriptstyle \mathrm{id}} & & \uparrow{\scriptstyle \cong} & & \uparrow{\scriptstyle \cong} & \\ \mathcal{O}_Y & \xrightarrow{pd} & \Omega^1_{Y/S} & \xrightarrow{pd} & \Omega^2_{Y/S} & \cdots . \end{array} \tag{4.1.1}$$

This point of view also gives a geometric interpretation of the data of a $p$-connection on an $\mathcal{O}_Y$-module $E$: it becomes equivalent to the data of an isomorphism between the two pullbacks of $E$ to the first infinitesimal neighborhood of $Y$ in $\mathbb{D}_Y(1)$ whose restriction to $Y$ is the identity.

These constructions can be used to amplify Mazur's interpretation [26, B] of the Cartier isomorphism as follows.

**Proposition 4.1.1.** *Let $Y \to S$ be a $p$-completely smooth morphism of crystalline prisms, and let $\mathrm{P}^1_{Y/S}$ (resp. $\mathrm{D}^1_{Y/S}$ ) be the first infinitesimal neighborhood of $Y$ in $Y(1)$ (resp., of $\tilde{\Delta}_{Y/S}$ in $\mathbb{D}_{Y/S}(1)$). Then the composition*

$$f_{Y/S} \colon \mathrm{P}^1_{Y/S} \subseteq Y(1) \xrightarrow{F_{Y(1)}} Y(1)$$

*defines a $p$-adic enlargement of $Y$ over $Y(1)$ and factors uniquely through a map $\tilde{F} \colon \mathrm{P}^1_{Y/S} \to \mathrm{D}^1_{Y/S}$, and there is a commutative diagram*

$$\begin{array}{ccccc} & & \tilde{J}/\tilde{J}^2 & \xrightarrow{\tilde{F}^*} & J/J^2 \\ & {\scriptstyle \tilde{\rho}}\nearrow & \uparrow & & \uparrow{\scriptstyle \cong} \\ J/J^2 & \xrightarrow{\cong} & \Omega^1_{Y/S} & \xrightarrow{\zeta} & \Omega^1_{Y/S}, \end{array} \tag{4.1.2}$$

*where the top horizontal arrow is induced by $\tilde{F}$, and $p\zeta = dF_Y^*$.*

*Proof.* We claim that the morphism $F_{Y(1)} \colon Y(1) \to Y(1)$ fits into a diagram:

$$\begin{array}{ccccc} \overline{\mathrm{P}}^1_{Y/S} & \longrightarrow & \mathrm{P}^1_{Y/S} & \longrightarrow & Y(1) \\ {\scriptstyle z}\,\downarrow & & \downarrow{\scriptstyle f_Y} & \swarrow{\scriptstyle F_{Y(1)}} & \\ \overline{Y} & \xrightarrow{\overline{\Delta}} & Y(1). & & \end{array}$$

The ideal of $\overline{\Delta}$ is generated by $p$ and the set of sections of $\mathcal{O}_{Y(1)}$ of the form $1 \otimes f - f \otimes 1$ for $f \in \mathcal{O}_Y$. Then

$$F_Y(1)^*(1 \otimes f - f \otimes 1) \equiv 1 \otimes f^p - f^p \otimes 1 \equiv (1 \otimes f - f \otimes 1)^p \pmod{p}.$$

Thus $F_Y(1)^*(1 \otimes f - f \otimes 1)$ belongs to the ideal of $\overline{\mathrm{P}}^1_{Y/S}$ in $Y(1)$, proving the existence of the dotted arrow $z$. This map gives $\mathrm{P}^1_{Y/S}$ the structure of a $p$-adic enlargement of $Y$ in $Y(1)$, and hence the map $f_Y$ factors

$$\mathrm{P}^1_{Y/S} \to \mathbb{D}_Y(1) \to Y(1).$$

The first map takes the diagonal of $\mathrm{P}^1_{Y/S}$ to the diagonal section of $\mathbb{D}_Y(1)$ and hence induces a map $\tilde{F}\colon \mathrm{P}^1_{Y/S} \to \mathrm{D}^1_{Y/S}$ as claimed. It is then immediate to check the commutativity of the diagram (4.1.2). $\square$

## 4.2 Connections and $p$-connections on envelopes

If $X$ is a closed subscheme of a smooth scheme $Y$ over $\mathbf{C}$, a fundamental ingredient in the theory of de Rham cohomology is the fact that the differential of the de Rham complex of $Y/\mathbf{C}$ extends naturally to the formal completion $\mathbb{F}_X(Y)$ of $Y$ along $X$, or, equivalently, that its structure sheaf $\mathcal{O}_{\mathbb{F}_X(Y)}$, viewed as a sheaf of $\mathcal{O}_Y$-algebras, admits an integrable connection compatible with its structure as an $\mathcal{O}_Y$-algebra. The crystalline incarnation of this fact says that if $Y/W$ is a smooth scheme and $X \subseteq Y$ is a closed subscheme, then the structure sheaf of the divided power envelope of $X$ in $Y$, viewed as a sheaf of $\mathcal{O}_Y$-algebras, admits an integrable and quasi-nilpotent connection. In this section we explain that analogs hold for $p$-adic dilatations and prismatic envelopes, but one must replace connections by $p$-connections. We will give a more conceptual geometric proof of these results later, in section §4.4.

**Proposition 4.2.1.** *Let $Y/S$ be a $p$-completely smooth morphism of $p$-torsion free $p$-adic formal schemes and let $i\colon X \to \overline{Y}$ be a closed immersion.*

*(1) The $p$-connection $d' := pd$ of $\mathcal{O}_Y$ extends uniquely to a multiplicative integrable and quasi-nilpotent $p$-connection $d'$ on the $\mathcal{O}_Y$-algebra $\mathcal{O}_{\mathbb{D}_X(Y)}$.*

*(2) If $Y/S$ is a morphism of crystalline prisms and $X \to \overline{Y}$ is a regular immersion, then $d'$ extends uniquely to a multiplicative integrable and quasi-nilpotent $p$-connection on the $\mathcal{O}_Y$-algebra $\mathcal{O}_{{\mathbb{\Delta}}_X(Y)}$. This $p$-connection is compatible with $F_Y$, in the sense that the following diagram commutes:*

$$\begin{array}{ccc}
\mathcal{O}_{{\mathbb{\Delta}}_X(Y)} & \xrightarrow{\quad d' \quad} & \mathcal{O}_{{\mathbb{\Delta}}_X(Y)} \otimes \Omega^1_{Y/S} \\
\Big\downarrow{\scriptstyle F_Y^*} & & \Big\downarrow{\scriptstyle F_Y^* \otimes dF_Y} \\
F_{Y*}(\mathcal{O}_{{\mathbb{\Delta}}_X(Y)}) & & F_{Y*}(\mathcal{O}_{{\mathbb{\Delta}}_X(Y)}) \otimes_{\mathcal{O}_Y} F_{Y*}(\Omega^1_{Y/S}) \\
& \searrow^{F_{Y*}(d')} & \Big\downarrow \\
& & F_{Y*}(\mathcal{O}_{{\mathbb{\Delta}}_X(Y)} \otimes_{\mathcal{O}_Y} \Omega^1_{Y/S}).
\end{array}$$

*Proof.* For (1), we suppose without loss of generality that $Y$ is affine, with $Y/S = \mathrm{Spf}(B/R)$. Let $I$ be the ideal defining $X$ and and let $B'$ be the uncompleted version of the $p$-adic dilatation of $I$. Then $B' \subseteq B_{\mathbf{Q}} := \mathbf{Q} \otimes B$ and is generated as a $B$-algebra by elements of the form $\rho(x)$ with $x \in I$. The $p$-connection $d'_B := pd_B$ extends uniquely to an integrable $p$-connection $d'_{\mathbf{Q}} \colon B_{\mathbf{Q}} \to B_{\mathbf{Q}} \otimes \Omega^1_{B/R}$. We claim that $d'_{\mathbf{Q}}$ maps $B'$ to $B' \otimes \Omega^1_{B/R}$. In fact, for $x \in I$:

$$d'_{\mathbf{Q}}\rho(x) = pd_{\mathbf{Q}}\rho(x) = d_{\mathbf{Q}}(p\rho(x)) = dx \in \Omega^1_{B/R}. \tag{4.2.1}$$

Since $d'$ satisfies the Leibniz rule and $B'$ is generated as an algebra over $B$ by $\rho(I)$, it follows that $d'_{\mathbf{Q}}b' \in \Omega^1_B \otimes B'$ for all $b' \in B'$. We write $d' \colon B' \to \Omega^1_{B/R} \otimes B'$ for the induced map, which is evidently an integrable $p$-connection on $B'$ and is the unique such extending the canonical $p$-connection on $B$ and which is multiplicative. This $p$-connection extends uniquely to the $p$-adic completion of $B'$ by continuity. To see the quasi-nilpotence, choose an element $\tau$ of $\mathcal{T}_{B/R}$, and let $d'_\tau$ be the resulting derivation. Observe from equation (4.2.1) that $d'_\tau(\rho(x)) \in B$ for all $x \in I$. Then $(d'_\tau)^2(\rho(x)) \in pB$, and, since $B'$ is generated as a $B$-algebra by $\rho(I)$, the multiplicativity of $d'_\tau$ implies that it is quasi-nilpotent modulo $p$.

Now suppose that $Y/S$ is a morphism of crystalline prisms. We are assuming that $X$ is regularly immersed in $\overline{Y}$, which will allow us to apply the explicit description of $\mathbb{A}_X(Y)$ from Theorem 2.3.3. With the notation there, we calculate:

$$\begin{aligned} pd't_{i,j+1} &= d'\delta^j(x_i) - d'(t^p_{i,j}) \\ &= pd\delta^j(x_i) - pt^{p-1}_{i,j}d't_{i,j} \\ d't_{i,j+1} &= d\delta^j(x_i) - t^{p-1}_{i,j}d't_{i,j}. \end{aligned} \tag{4.2.2}$$

This shows that the subrings $B^{(j)}$ and $B^\infty$ of $\mathbf{Q} \otimes B$ are stable under the $p$-connection $d'_{\mathbf{Q}}$. We shall show by induction on $j$ that the resulting $p$-connection $d'$ on each $B^{(j)}$ (and hence also on $B^\infty$) is quasi-nilpotent. This is clear for $j = 0$, since $d'$ is divisible by $p$ on $B = \mathcal{O}_Y$. It is enough to show that, if $\tau$ is a tangent vector and $d'_\tau$ is the resulting derivation, then a sufficiently high power of $d'_\tau$ maps each generator of $B^{(j+1)}$ to $B^{(j)}$ modulo $p$. In fact, it follows from (4.2.2) that

$$d'_\tau(t_{i,j+1}) = d_\tau(\delta^j(x_i)) - t^{p-1}_{i,j}d'_\tau(t_{i,j}) \in B^{(j)}.$$

We conclude that the $p$-connection on $B^\infty$ is quasi-nilpotent. Furthermore, the diagram in statement (2) with $B_{\mathbf{Q}}$ in place of $B^\infty$ commutes by functoriality, and it follows that the same holds with subring $B^\infty$ in place of $B_{\mathbf{Q}}$. These results extend to the $p$-adic completion by continuity. □

**Remark 4.2.2.** We should point out that the more general prismatic envelopes constructed in Theorem 2.3.8 may not inherit a $p$-connection. For example, let $Y := \mathrm{Spf}(\mathbf{Z}_p[x]\hat{\ })$, with $G_Y(x) := \psi(x) = x$ and let $X$ be the closed subscheme defined by $(p, x)$. Since $a \equiv a^p \pmod{p}$ for every $a \in \mathbf{Z}_p$, it follows that $f(x) \equiv (f(x))^p \pmod{I_{X/Y}}$ for every $f \in \mathbf{Z}_p[x]\hat{}$, so the Theorem applies. Let us follow the procedure there to compute $\mathbb{A}_X(Y)$. We find:

$$\begin{aligned} ps_1 &= x \\ p\psi(s_1) &= \psi(x) \\ &= ps_1 \\ \psi(s_1) &= s_1. \end{aligned}$$

Hence

$$\begin{aligned} ps_2 &= \psi(s_1) - s_1^p \\ &= s_1 - s_1^p \\ \psi(ps_2) &= \psi(s_1) - (\psi(s_1))^p \\ &= s_1 - s_1^p \\ &= ps_2 \\ \psi(s_2) &= s_2. \end{aligned}$$

Continuing in this way, we find that $\mathbb{A}_X(Y)$ is the formal spectrum of the $p$-adic completion of

$$\mathbf{Z}_p[s_1, s_2, \ldots]/(ps_2 = s_1 - s_1^p,\ ps_3 = s_2 - s_2^p,\ \ldots),$$

with $\psi(s_i) = s_i$ for all $i$, and with $\mathbb{A}_X(Y) \to Y$ defined by $x \mapsto ps_1$.

Suppose that the $p$-connection $d'$ of $\mathcal{O}_Y$ did extend to to $\mathcal{O}_{\mathbb{A}_X(Y)}$. Then for $i > 0$ we would have:

$$\begin{aligned} s_i &= s_i^p + ps_{i+1} \\ d's_i &= ps_i^{p-1} d's_i + pd's_{i+1}. \end{aligned}$$

This equation implies that $d's_i \in p\mathcal{O}_{\mathbb{A}} \otimes \Omega^1_{Y/S}$ for all $i > 0$. Using this fact and the same equation again, we see that in fact $d's_i \in p^2\mathcal{O}_{\mathbb{A}} \otimes \Omega^1_{Y/S}$, and, by induction, that $d's_i \in p^n\mathcal{O}_{\mathbb{A}} \otimes \Omega^1_{Y/S}$ for all $n$. Since $d's_1 = dx$, we see that $dx$ vanishes in the $p$-adic completion of $\mathcal{O}_{\mathbb{A}} \otimes \Omega^1_{Y/S}$, which is absurd.

**Example 4.2.3** Suppose that $(x_1, \ldots, x_n)$ is a local system of coordinates for $Y/S$, that $X \subseteq Y$ is defined by $(x_1, \ldots, x_r, p)$, and that $\phi(x_i) \in (x_1, \ldots, x_r, p^2)$ for $1 \leq i \leq r$. Then by (1) of Theorem 2.3.13, $\mathbb{A}_X(Y)$ is the $p$-adic completion of $\mathcal{O}_Y\langle t_1, \ldots, t_r\rangle$, where $pt_i = x_i$. Then $d't_i = pdt_i = dx_i$, and, more generally, $d't_i^{[n]} = t_i^{[n-1]}dx_i$.

**Remark 4.2.4.** Our next goal is to define a suitable notion of module with prismatic connection. Suppose again that $Y/S$ is a $p$-completely smooth morphism of crystalline prisms and that $i \colon X \to \overline{Y} \to Y$ is a regular immersion. Let $(\mathbb{A}_X(Y), z, \pi)$ denote the prismatic envelope of $X$ in $Y$. We shall sometimes find it convenient to work with sheaves on $Y$, sometimes with sheaves on $\mathbb{A}_X(Y)$, and sometimes with sheaves on $X$. Let us explain how we can shift among these points of view. As we have seen, the projection $\pi\colon \mathbb{\Delta}_X(Y) \to Y$ is an affine morphism of crystalline prisms. Since formation of prismatic envelopes commutes with $p$-completely flat base change $Y' \to Y$, the sheaf $\pi_*(\mathcal{O}_{\mathbb{\Delta}_X(Y)})$ is $p$-completely quasi-coherent, in the sense of Definition 7.4.1. Moreover, thanks to a (relative version of) Proposition 7.4.2, there is a natural equivalence from the category of sheaves of $p$-completely quasi-coherent $\mathcal{O}_Y$-modules on $Y$ endowed with a (compatible) $\mathcal{O}_{\mathbb{\Delta}_X(Y)}$-module structure, and the category of $p$-completely quasi-coherent $\mathcal{O}_{\mathbb{\Delta}_X(Y)}$-modules on $\mathbb{\Delta}_X(Y)$. We shall allow ourselves to identify these categories when we feel that no confusion should result. Note also that, if $E$ is a $p$-completely quasi-coherent sheaf of $\mathcal{O}_{\mathbb{\Delta}_X(Y)}$-modules, Proposition 7.4.2 tells us that $R^q\pi_*(E)$ vanishes for $q > 0$, so that cohomology computed with these two points of view coincide.

Next, observe that the underlying topological spaces of $\mathbb{\Delta}_X(Y)$ and of $\overline{\mathbb{\Delta}}_X(Y)$ coincide, that $F_{\mathbb{\Delta}}$ is the identity on this space, and that $z\colon \overline{\mathbb{\Delta}}_X(Y) \to X$ is an affine morphism of schemes. Thus, if $E$ is a $p$-completely quasi-coherent sheaf of $\mathcal{O}_{\mathbb{\Delta}_X(Y)}$-modules, Proposition 7.4.2 also implies that $R^q z_*(E)$ vanishes for $q > 0$. Furthermore the support of any sheaf of $\pi_*(\mathcal{O}_{\mathbb{\Delta}_X(Y)})$-modules is contained in the subset $X$; and if $F$ is any sheaf on $Y$ with this property, then the natural map $F \to i_*i^{-1}F$ is an isomorphism, so we may indifferently view such a sheaf as living on $Y$ or on $X$.

**Definition 4.2.5.** *Let $Y/S$ be a $p$-completely smooth morphism of crystalline prisms, let $X \to \overline{Y}$ be a regular closed immersion, and let $E$ be a $p$-completely quasi-coherent sheaf of $\mathcal{O}_{\mathbb{\Delta}_X(Y)}$-modules (see the previous paragraph).*

*(1) An integrable $p$-connection $\nabla'\colon E \to E \otimes \Omega^1_{Y/S}$ on $E$ is compatible if*

$$\nabla'(ae) = a\nabla'(e) + e \otimes d'a.$$

*for every pair of local sections $a \in \mathcal{O}_{\mathbb{\Delta}_X(Y)}, e \in E$.*

*(2) A module with integrable prismatic $p$-connection on $X/Y/S$ is a $p$-completely quasi-coherent sheaf of $\mathcal{O}_{\mathbb{\Delta}_X(Y)}$-modules endowed with an integrable, quasi-nilpotent, and compatible $p$-connection.*

*We denote the category of modules with prismatic $p$-connection on $X/Y/S$ by* $\mathrm{MICP}(X/Y/S)$.

The analogous category $\mathrm{MIC}(X/Y/S)$ of $p$-completely quasi-coherent $\mathcal{O}_{\mathbb{PD}_X(Y)}$-modules with compatible quasi-nilpotent connection is discussed, for example, in [7, §6]. In the next section, we shall show that the category $\mathrm{MICP}(X/Y/S)$ is, in a suitable sense, independent of $Y$.

## 4.3 The prismatic Poincaré Lemma

Let $S$ be a crystalline prism and $X/\overline{S}$ a smooth morphism. Although we have not yet discussed the prismatic site, we can use "bare hands" methods to prove the existence of a functorial cohomology theory $H^*_{\mathbb{\Delta}}(X/S)$, including a theory of coefficients, that will agree with the cohomology of the prismatic topos $(X/S)_{\mathbb{\Delta}}$. The key is to show that, if $X$ is embedded in $p$-completely smooth $Y/S$, the category of modules with prismatic $p$-connection on $X/Y/S$ is independent of the choice of $Y$, and that the same independence holds for the cohomology of any such module. We accomplish this using our detailed study of the nature of prismatic envelopes in §2.3. (When the base $S$ is perfect, another approach is possible, by reduction to the usual crystalline theory using the F-transform, but we shall not explain this here.) In the following theorem, we show that the category $\mathrm{MICP}(X/Y/S)$ and the cohomology of an object therein are independent of the choice of $Y$. Here we view an object $(E, \nabla')$ of this category as living on the topological space $X$, as explained in Remark 4.2.4.

**Theorem 4.3.1** (The prismatic Poincaré Lemma)**.** *Let $S$ be a crystalline prism, let $X/\overline{S}$ be a smooth morphism, let $i\colon X \to \overline{Y} \to Y$ be a regular closed immersion, and let $g\colon Y \to Z$ be a morphism of $p$-completely smooth crystalline prisms over $S$ such that $j := g \circ i$ is again a regular closed immersion.*

*(1) The morphism $g$ induces an equivalence of categories:*

$$\mathbb{\Delta}(g)^*\colon \mathrm{MICP}(X/Z/S) \to \mathrm{MICP}(X/Y/S).$$

*(2) If $(E, \nabla')$ is an object of $\mathrm{MICP}(X/Z/S)$, then the natural map*

$$(E \otimes \Omega^\bullet_{Z/S}, d') \to (\mathbb{\Delta}(g)^*(E) \otimes \Omega^\bullet_{Y/S}, d')$$

*is a strict quasi-isomorphism (see Definition 7.2.2).*

*Proof.* The functor induced by the morphism $g$ is obtained as follows. Since $g$ takes the embedding $i\colon X \to Y$ to the embedding $j\colon X \to Z$, it induces a morphism of envelopes: $\mathbb{\Delta}(g)\colon \mathbb{\Delta}_X(Y) \to \mathbb{\Delta}_X(Z)$. This morphism is compatible with the corresponding $p$-connections: if $b \in \mathcal{O}_{\mathbb{\Delta}_X(Z)}$, then

$$(\mathbb{\Delta}(g)^* \otimes g^*)(d'b) = d'\mathbb{\Delta}(g)^*(b) \in \mathcal{O}_{\mathbb{\Delta}_X(Y)} \otimes \Omega^1_{Y/S}.$$

If $(E, \nabla')$ is an object of $\mathrm{MICP}(X/Z, S)$, there is a unique a $p$-connection:

$$\nabla'\colon \mathcal{O}_{\mathbb{\Delta}_X(Y)} \otimes_{\mathbb{\Delta}_X(Z)} E \to \mathcal{O}_{\mathbb{\Delta}_X(Y)} \otimes_{\mathbb{\Delta}_X(Z)} E \otimes_{\mathcal{O}_Y} \Omega^1_{Y/S}$$

such that

$$a \otimes e \mapsto pe \otimes d'a + a \otimes (\mathrm{id}_E \otimes \mathbb{\Delta}(g)^*)(\nabla'(e)).$$

(One easily checks that this is well-defined, compatible with the $p$-connection on $\mathcal{O}_{\mathbb{\Delta}_X(Y)}$, and quasi-nilpotent.)

The proof will proceed by reducing to calculations in the case of prismatic neighborhoods of a point. Thus the key case occurs when $Y = S$, when $g$ is a closed immersion, which we now denote by $s\colon S \to Z$, and when $X = \overline{S}$. Under these conditions, the category $\mathrm{MICP}(X/Y/S) = \mathrm{MICP}(\overline{S}/S/S)$ is just the category of $p$-completely quasi-coherent $\mathcal{O}_S$-modules, and the theorem reduces to the following lemma.

**Lemma 4.3.2.** *Let $h\colon Z \to S$ be a $p$-completely smooth morphism of crystalline prisms. Assume that $h$ admits a section $s\colon S \to Z$, compatible with the Frobenius lifts.*

*(1) The functor*

$$\mathbb{\Delta}(s)^*\colon \mathrm{MICP}(\overline{S}/Z/S) \to \mathrm{MICP}(\overline{S}/S/S) : (E, \nabla') \mapsto s^*(E)$$

*is an equivalence, with quasi-inverse*

$$\mathbb{\Delta}(h)^*\colon \mathrm{MICP}(\overline{S}/S/S) \to \mathrm{MICP}(\overline{S}/Z/S) : E' \mapsto (\mathcal{O}_{\mathbb{\Delta}_S(Z)}\hat{\otimes} E', d'\hat{\otimes}\mathrm{id}).$$

*In particular, if $E \in \mathrm{MICP}(\overline{S}/Z/S)$, the natural maps*

$$(\mathcal{O}_{\mathbb{\Delta}_S(Z)}\hat{\otimes}_{\mathcal{O}_S} E^{\nabla'}, d'\hat{\otimes}\mathrm{id}) \to (E, \nabla') \quad \text{and} \quad E^{\nabla'} \to s^*(E)$$

*are isomorphisms, and if $E' \in \mathrm{MICP}(\overline{S}/S/S)$ and $(E, \nabla') := \mathbb{\Delta}(h)^*(E')$, then the maps*

$$E' \to E^{\nabla'} \quad \text{and} \quad E' \to s^*(E)$$

*are isomorphisms.*

*(2) If $(E, \nabla')$ is an object of $\mathrm{MICP}(\overline{S}/Z/S)$, then*

*(a) The family $(E_n^{\nabla'} : n \in \mathbf{N})$ forms a $p$-adic module, and $E^{\nabla'} \cong \varprojlim E_n^{\nabla'}$.*

*(b) The map $E^{\nabla'} \to (E \otimes \Omega^\bullet_{Z/S}, d')$ is a strict quasi-isomorphism.*

*Proof.* Since $(S, \mathrm{id}_{\overline{S}}, \mathrm{id}_S)$ is an $\overline{S}$-prism over $S$, the map $\mathbb{\Delta}_{\overline{S}}(Z) \to S$ admits a section $\tilde{s}'\colon S \to \mathbb{\Delta}_{\overline{S}}(Z)$, and the functor $\mathbb{\Delta}(s)^*$ in statement (1) takes an object $(E, \nabla')$ of $\mathrm{MICP}(\overline{S}/Z/S)$ to its pullback via $\tilde{s}'$. As we saw in statement (1) of Theorem 2.3.13, the prismatic envelope $\mathbb{\Delta}_{\overline{S}}(Z)$ of $S$ in $Z$ can be identified with the divided power envelope $\mathbb{PD}_{\overline{S}}(\mathbb{D}_{\overline{S}}(Z))$ of the section of the $p$-adic dilatation of $\overline{S}$ in $Z$ defined by $s$. Working locally, we may assume that $Z \to S$ is $\mathrm{Spf}(B \to C)$ and that there is a set of local coordinates $(x_1, \ldots, x_r)$ for $Z/S$ which generates the ideal of $S$ in $Z$. Then the formal completion of $S$ along $Z$ is $\mathrm{Spf}\, B[[x_1, \ldots, x_r]]$, and $\mathbb{D}_{\overline{S}}(Z)$ is $\mathrm{Spf}\, B[t_1, \ldots, t_r]\hat{}$, with $x_i$ mapping to $pt_i$. Thus the map $\mathbb{D}_{\overline{S}}(Z) \to S$ is also $p$-completely smooth.

Continuing to work locally with the aid of coordinates, we write $\mathbb{PD}_{\overline{S}}(\mathbb{D}_{\overline{S}}(Z))$ as the formal spectrum of the $p$-adic completion of the PD-algebra $B\langle t_1, \ldots, t_r\rangle$. Since $pt_i = x_i$ and $d'x_i = pdx_i$, it follows that $d't_i = dx_i$, and since $\mathcal{O}_{\mathbb{\Delta}_S(Z)}$ is $p$-torsion free, that $d't_i^{[n]} = t_i^{[n-1]}dx_i$ for all $n$. (See Example 4.2.3.) Because the $p$-connection $\nabla'$ on $E$ is compatible with the $p$-connection on $\mathcal{O}_{\mathbb{\Delta}_S(Z)}$, it behaves like a usual connection on $E$, viewed now as a $B\langle t_1, \ldots, t_r\rangle\hat{}$-module; furthermore, it is quasi-nilpotent as a connection. In this light the lemma becomes a standard fact from the usual crystalline theory. Let us nevertheless write some details.

Let $\partial_i := \nabla'_{\partial/\partial x_i}$, and let

$$r := \sum_I (-t)^{[I]} \partial^I \in \mathrm{End}_{\mathcal{O}_S}(E). \tag{4.3.1}$$

This formula converges modulo every power of $p$, and hence $p$-adically, because of the quasi-nilpotence of $\nabla'$: if $e \in E$ and $n \in \mathbf{N}$, then $\partial^I e \in p^n E$ for $|I| \gg 0$. Then a standard formal verification shows that

$$\begin{aligned} \nabla' \circ r &= 0 \\ \sum_I t^{[I]} r \circ \partial^I &= \mathrm{id}. \end{aligned} \tag{4.3.2}$$

These equationsimply that every element $e$ of $E$ can be written as a $p$-adically convergent sum

$$e = \sum t^{[I]} e_I$$

where $e_I := r(\partial_I(e)) \in E^{\nabla'}$. We claim that any such expression is unique. To see this suppose that $e = \sum_I t^{[I]} e'_I$ is another such expression and calculate:

$$\begin{aligned} r(e) &= r\left(\sum t^{[I]} e'_I\right) \\ &= \sum_J (-t)^{[J]} \partial^J \left(\sum_I t^{[I]} e'_I\right) \\ &= \sum_{I,J} (-t)^{[J]} t^{[I-J]} e'_I \\ &= \sum_{I,J} \binom{I}{J} (-1)^J t^{[I]} e'_I \\ &= \sum_I \sum_{0 \le J \le I} \binom{I}{J} (-1)^{|J|} t^{[I]} e'_I \\ &= \sum_I (1-1)^{|I|} t^{[I]} e'_I \\ &= e'_0. \end{aligned}$$

This holds modulo every power of $p$ and in the limit, and it proves that the map $E^{\nabla'} \to s^*(E)$ is an isomorphism modulo every power of $p$ and in the limit. It follows that $(E_n^{\nabla'} : n \in \mathbf{N})$ forms a $p$-adic module. The equality $E^{\nabla'} = \varprojlim E_n^{\nabla'}$ follows from the left exactness of inverse limits. We have now proved statement (2a).

Applying the above calculations to $\partial^I e$, we see that $r(\partial^I e) = e'_I$ for all $I$. The proves the claimed uniqueness, and shows that the map

$$\mathbb{A}(h)^*(E^{\nabla'}) = E^{\nabla'} \hat{\otimes}_{\mathcal{O}_S} \mathcal{O}_{\mathbb{A}_{\overline{S}}(Z)} \longrightarrow E \tag{4.3.3}$$

is an isomorphism modulo every power of $p$ and in the limit. Furthermore, if $E'$ is an object of $\mathrm{MICP}(\overline{S}/S/S)$ and $E := \mathbb{A}(h)^*(E')$, it is clear that $E^{\nabla'} \cong E'$. Statement (1) follows.

We continue in our localized situation to address statement (2b). Statement (1) and the identification $\mathcal{O}_{\mathbb{A}_{\overline{S}}(Z)} \cong B\langle t_1, \ldots, t_r\rangle\hat{}$ together imply that we have a natural isomorphism:

$$(E \otimes \Omega^\bullet_{Z/S}, d') \cong (E^{\nabla'} \hat{\otimes}_B B\langle t_1, \ldots, t_r\rangle\hat{}, d)$$

where $d$ is the usual differential of the de Rham complex of the divided power polynomial algebra. Thus the standard divided power Poincaré lemma [7, 6.12] implies that the map is a strict quasi-isomorphism. □

The remainder of the proof of the theorem will be series of reductions to the case handled in Lemma 4.3.2.

**Lemma 4.3.3.** *Theorem 4.3.1 holds if $Z \to S$ is the identity map and $g\colon Y \to S$ admits a section $s$ such that $i = s \circ j$. In particular,*

$$\mathbb{A}^*(i)\colon \mathrm{MICP}(X/Y/S) \to \mathrm{MICP}(X/S/S)$$

*is quasi-inverse to $\mathbb{A}^*(g)$. If $(E, \nabla') \in \mathrm{MICP}(X/Y/S)$, then $E^{\nabla'} \to i^*(E)$ is an isomorphism in the limit and modulo every power of $p$. Thus $(E_n^{\nabla'} : n \in \mathbf{N})$ forms a $p$-adic module and $E^{\nabla'} = \varprojlim E_n^{\nabla'}$.*

*Proof.* Since $Y/S$ was by assumption $p$-completely smooth, and in this case $g$ is the structural morphism of $Y/S$, it follows that $g$ is $p$-completely smooth. Since $X \to \overline{S}$ is smooth and factors through the closed immersion $\overline{j}\colon X \to \overline{S}$, necessarily $X/\overline{S}$ is étale, and in fact $j$ is an open immersion. Since our statements are local along $X$, we may as well assume that $X = \overline{S}$. Then our claim reduces to Lemma 4.3.2, with $g$ here playing the role of $h$ there. □

The next step requires localization in the $p$-completely flat topology.

**Lemma 4.3.4.** *Theorem 4.3.1 holds if $Z \to S$ is the identity map. Moreover,*

$$\mathbb{A}^*(i)\colon \mathrm{MICP}(X/Y/S) \to \mathrm{MICP}(X/S/S)$$

*is quasi-inverse to $\mathbb{A}^*(g)$. If $(E, \nabla') \in \mathrm{MICP}(X/Y/S)$, then $E^{\nabla'} \to i^*(E)$ is an isomorphism in the limit and modulo every power of $p$. Thus $(E_n^{\nabla'} : n \in \mathbf{N})$ forms a $p$-adic module and $E^{\nabla'} = \varprojlim E_n^{\nabla'}$.*

*Proof.* Thanks to the descent results in section §7, this can be checked $p$-completely flat locally on $S$, so we may and shall assume that $S$ is affine. As we saw in the proof of the previous lemma, our assumption that $Z = S$ implies that $X \to \overline{S}$ is an open immersion, and we may as well assume that it is an isomorphism. As each map $S_n \to S_{n+1}$ is a nilpotent immersion of affine schemes and $g$ is formally smooth, there exists a section of $g$ which extends the embedding $S_1 = X \to Y$. This map is not necessarily compatible with the Frobenius lifts, but, by Theorem 1.7.11, there exist a $p$-completely flat covering $\tilde{S} \to S$ of crystalline prisms and a section $\tilde{s}\colon \tilde{S} \to Y \times_S \tilde{S}$ which is compatible with the embeddings as well at the prism structures. Then Lemma 4.3.3 applies. □

**Lemma 4.3.5.** *Theorem 4.3.1 holds if $g$ is $p$-completely smooth.*

*Proof.* Let $(E, \nabla')$ be an object of $\mathrm{MICP}(X/Y/S)$ and let $(E_n, \nabla')$ be its reduction modulo $p$. The relative cohomology sheaves

$$E_n^q := R^q g_*(E_n \otimes \Omega^{\bullet}_{Y/Z}, d'),$$

are just the cohomology sheaves of the complex $(E_n \otimes \Omega^{\bullet}_{Y/Z}, d')$. Lemma 4.3.4 applies to the morphism $g$, with $S$ replaced by $Z$ and tells us that these sheaves vanish for $q > 0$. It also tells us that the family $(E_n^0 : n \in \mathbf{N})$ forms a $p$-adic module, that $E^0 = \varprojlim E_n^0$, and that each $g^*(E_n^0) \to E_n$ is an isomorphism.

As we saw in Proposition 1.2.3, the sheaves $E_n^0$ admit a Gauss-Manin $p$-connection, and there is a spectral sequence with

$$E_2^{p,q} \cong \mathcal{H}^p(E_n^q \otimes \Omega^{\bullet}_{Z/S}, d'_{GM}) \Rightarrow \mathcal{H}^{p+q}(E_n \otimes \Omega^{\bullet}_{Y/S}, d').$$

Since $E_n^q = 0$ for $q > 0$, in fact we have isomorphisms

$$\mathcal{H}^n(E_n \otimes \Omega^{\bullet}_{Y/S}, d') \cong \mathcal{H}^n(E_n^0 \otimes \Omega^{\bullet}_{Z/S}, d'_{GM}). \tag{4.3.4}$$

Furthermore it follows from the naturality of the Gauss-Manin construction that the isomorphism $\mathbb{\Delta}(g)^*(E^0) \to E$ is compatible with the $p$-connections, so in fact we have an isomorphism:

$$\mathbb{\Delta}(g)^*(E^0, \nabla'_{GM}) \to (E, \nabla').$$

This implies that the functor $\mathbb{\Delta}(g)^*$ is essentially surjective.

Now suppose that $(E', \nabla')$ is an object of $\mathrm{MICP}(X/Z/S)$, and let $(E, \nabla') := \mathbb{\Delta}(g)^*(E', \nabla')$. It follows from Lemma 4.3.4 that $E' \cong E^0$, and hence that the functor $\mathbb{\Delta}(g)^*$ is fully faithful. The fact that $E' \to (\mathbb{\Delta}(g)^*(E'))^0 = E^0$ is an isomorphism, compatibly with the $p$-connections $\nabla'$ on the source and $\nabla'_{GM}$ on the target, implies that

$$(E' \otimes \Omega^{\bullet}_{Z/S}, d') \to (E \otimes \Omega^{\bullet}_{Y/S}, d')$$

is a strict quasi-isomorphism, as claimed in the theorem. □

We can now prove the general case of Theorem 4.3.1. First note that if $f$ and $g$ are composable morphisms and the result holds for any two of $f, g$, and $g \circ f$, then it also holds for all three, as is easily seen. Now in the situation of the theorem, consider the graph $\Gamma_g \colon Y \to Y \times_S Z$ of $g$. The projections $\pi_Y \colon Y \times_S Z \to Y$ and $\pi_Z \colon Y \times_S Z \to Y$ are $p$-completely smooth, and Lemma 4.3.5 implies that the theorem is true for each of these. Since $\pi_Y \circ \Gamma_g = \mathrm{id}_Y$, the result for $\pi_Y$ and for $\mathrm{id}_Y$ imply the result for $\Gamma_g$. Since $g = \pi_Z \circ \Gamma_g$, the result for $\pi_Z$ and for $\Gamma_g$ imply the result for $g$. □

**Corollary 4.3.6.** *Let $Y/S$ be a $p$-completely smooth morphism of crystalline prisms and let $X \to \overline{Y}$ be a closed immersion, with $X/\overline{S}$ smooth. Then if $(E, \nabla')$ is an object of $\mathrm{MICP}(X/Y/S)$ and $n > 0$, the cohomology sheaves $\mathcal{H}^q(E_n \otimes \Omega^{\bullet}_{Y/S}, d')$ vanish if $q > \dim(X/S)$.*

*Proof.* We can check this statement locally on $X$, and so we may and shall assume that there is a $p$-completely smooth crystalline prism $\tilde{X}/S$ lifting $X/\overline{S}$. The corollary is certainly true with $\tilde{X}$ in place of $Y/S$. Now we have morphisms of embeddings:

$$(\tilde{X}, X) \leftarrow (\tilde{X} \times Y, X) \rightarrow (Y, X)$$

and so Theorem 4.3.1 allows us to carry the result for $(\tilde{X}, X)$ over to $(Y, X)$. □

## 4.4 Prismatic stratifications and differential operators

In this section we adapt some standard constructions relating groupoids, stratifications, and differential operators to prismatic crystals, adapting the general framework described in §1.5 to our current context. To facilitate the comparison to other theories, we let $\mathbb{T}$ stand for any one of $\mathbb{PD}$, $\mathbb{D}$, or $\mathbb{\Delta}$, and sometimes also for $\mathbb{F}$. We work over a $p$-torsion free $p$-adic formal scheme $S$, endowed with a $\delta$-scheme structure in the prismatic case. If $X$ is a smooth $\overline{S}$-scheme, we let $\mathbb{T}(X/S)$ be the category of PD-enlargements of $X/S$, of $p$-adic enlargements of $X/S$, or of $X/S$-prisms, respectively. If $Y/S$ is a $p$-completely smooth morphism of $p$-adic formal schemes (or of crystalline prisms if $\mathbb{T} = \mathbb{\Delta}$) and if $X \subseteq Y$ is a closed embedding, we write $\mathbb{T}_X(Y)$ for the appropriate envelope of $X$ in $Y$.

Let $f \colon Y \to S$ be a $p$-completely smooth morphism of $p$-adic formal schemes, or, in the prismatic context, of crystalline prisms. If $n \in \mathbf{N}$, we write $Y(n)$ for the $(n+1)$-fold product $Y \times_S Y \times_S \cdots \times_S Y$, noting that this product is again $p$-torsion free and that it inherits a Frobenius lifting in the prismatic case. If $X \to Y$ is a closed immersion, we let $\mathbb{T}_{X/Y}(n)$ denote the $\mathbb{T}$-envelope of $X$ in $Y(n)$, embedded via the diagonal. In the important case in which $X = \overline{Y}$, we abbreviate this to $\mathbb{T}_Y(n)$.

**Proposition 4.4.1.** *Let $f \colon Y \to S$ and $g \colon Z \to S$ be $p$-completely smooth morphisms of $p$-torsion free $p$-adic formal schemes (resp. of formal $\delta$-schemes if $\mathbb{T} = \mathbb{\Delta}$), and let $i \colon X \to \overline{Y}$ and $j \colon X \to \overline{Z}$ be closed $\overline{S}$-immersions.*

*(1) The envelope $\mathbb{T}_X(Y \times_S Z)$ represents the product $\mathbb{T}_X(Y) \times \mathbb{T}_X(Z)$ in the category $\mathbb{T}(X/S)$. Similarly, for each $n$, $\mathbb{T}_{X/Y}(n)$ represents the $(n+1)$-fold product of $\mathbb{T}_X(Y)$ with itself in the category $\mathbb{T}(X/S)$.*

*(2) For every $n \in \mathbf{N}$ and for $0 \leq k \leq n$, the diagram*

$$\begin{array}{ccc} \mathbb{T}_{X/Y}(n) & \longrightarrow & \mathbb{T}_Y(n) \\ {\scriptstyle \mathbb{T}_X(p_k)} \downarrow & & \downarrow {\scriptstyle p_k} \\ \mathbb{T}_X(Y) & \xrightarrow{\ \pi\ } & Y \end{array}$$

*is Cartesian.*

*Proof.* The proof of statement (1) is purely formal and omitted.

We focus on the prismatic context in the proof of statement (2), but continue to use the generalized notation. First note that each map

$$\mathbb{T}_Y(n) \longrightarrow Y(n) \xrightarrow{pr_k} Y$$

is $p$-completely flat, by statement (3) of Proposition 6.1.2. It follows that the product $\mathbb{T}_X(Y) \times_Y \mathbb{T}_Y(n)$ is $p$-torsion free, and hence that the maps

$$z_Y \times z_{Y(n)} \colon \overline{\mathbb{T}}_X(Y) \times_Y \overline{\mathbb{T}}_Y(n) \to X \times_Y \overline{Y} = X$$

and

$$\pi_Y \times \pi_{Y(n)} \colon \mathbb{T}_X(Y) \times_Y \mathbb{T}_Y(n) \to Y \times_Y Y(n) = Y(n)$$

endow $\mathbb{T}_X(Y) \times_Y \mathbb{T}_Y(n)$ with the structure of an $X$-prism over $Y(n)$. Hence there is a unique map

$$\mathbb{T}_X(Y) \times_Y \mathbb{T}_Y(n) \to \mathbb{T}_{X/Y}(n) \tag{4.4.1}$$

of $X$-prisms over $Y(n)$. On the other hand, the map of pairs

$$(\mathrm{id}, i) \colon (Y(n), X) \to (Y(n), \overline{Y})$$

defines a $Y$-morphism $\mathbb{T}_X(Y(n)) \to \mathbb{T}_Y(n)$, and the map

$$(p_k, i) \colon (Y(n), X) \to (Y, X)$$

defines a $Y$-morphism $\mathbb{T}_{X/Y}(n) \to \mathbb{T}_X(Y)$. These assemble to define a morphism

$$\mathbb{T}_{X/Y}(n) \to \mathbb{T}_X(Y) \times_Y \mathbb{T}_Y(n)$$

which is inverse to the mapping (4.4.1). □

We find morphisms

$$\begin{aligned} t, s : \mathbb{T}_{X/Y}(1) &\longrightarrow \mathbb{T}_X(Y) \\ \iota : Y &\longrightarrow \mathbb{T}_Y(1) \\ c : \mathbb{T}_{X/Y}(1) \times_{\mathbb{T}_X(Y)} \mathbb{T}_{X/Y}(1) &\longrightarrow \mathbb{T}_{X/Y}(1) \\ \tau \colon \mathbb{T}_{X/Y} &\longrightarrow \mathbb{T}_{X/Y} \end{aligned}$$

covering the corresponding structural morphisms of the "indiscrete groupoid" $\mathcal{G}_Y$ defined by $Y$, as described in Example 1.4.1. (The morphism $c$ is obtained as the composition of the isomorphism

$$\mathbb{T}_{X/Y}(1) \times_{\mathbb{T}_X(Y)} \mathbb{T}_{X/Y}(1) \cong \mathbb{T}_{X/Y}(2)$$

of statement (1) of Proposition 4.4.1 with the map $\mathbb{T}_{X/Y}(2) \to \mathbb{T}_{X/Y}(1)$ induced by functoriality from the composition law of $\mathcal{G}_Y$.) These assemble to define a groupoid $\mathcal{G}_{\mathbb{T}_{X/Y}}$ over $Y$, with $\mathbb{T}_{X/Y}(1)$ representing the class of arrows. In fact, statement (1) of Proposition 4.4.1 tells us that $\mathbb{T}_{X/Y}(1)$ is the fiber product of $\mathbb{T}_X(Y)$ with itself over $S$ in the category $\mathbb{T}(X/S)$, so $\mathcal{G}_{\mathbb{T}_{X/Y}}$ can be viewed as the indiscrete groupoid on the object $\mathbb{T}_X(Y)$ when viewed in this category. In the important case when $X = \overline{Y}$, we just write $\mathcal{G}_{\mathbb{T}_Y}$ for this indiscrete groupoid on the object $Y$ of the category $\mathbb{T}(\overline{Y}/S)$.

If $X$ is a closed subscheme of $\overline{Y}$, there is a unique morphism $\epsilon$ making the following diagram commute:

$$\begin{array}{ccc} \mathbb{T}_X(Y) \times_Y \mathbb{T}_Y(1) & \xrightarrow{\epsilon} & \mathbb{T}_Y(1) \times_Y \mathbb{T}_X(Y) \\ \downarrow \cong & \swarrow \cong & \\ \mathbb{T}_X(Y(1)), & & \end{array}$$

where the isomorphisms come from statement (2) of Proposition 4.4.1. On points, the morphism $\epsilon$ takes $(y_1, y_1, y_2)$ to $(y_1, y_2, y_2)$; its nontrivial content is that if $y_1$ is in $\mathbb{T}_X(Y)$ and $(y_1, y_2)$ in in $\mathbb{T}_Y(1)$, then $y_2$ is also in $\mathbb{T}_X(Y)$. It is easy to see from this that $\epsilon$ defines a $\mathcal{G}_{\mathbb{T}_Y}$-stratification on $\mathbb{T}_X(Y)$, or, equivalently, a right action:

$$r\colon \mathbb{T}_X(Y) \times_Y \mathbb{T}_Y(1) \to \mathbb{T}_X(Y). \tag{4.4.2}$$

In every case we are considering, the morphisms $s$ and $t : \mathbb{T}_{X/Y}(1) \to Y$ are affine. Let $\mathcal{A}_{\mathbb{T}_{X/Y}} := t_*(\mathcal{O}_{\mathbb{T}_{X/Y}(1)})$. Then the groupoid $\mathcal{G}_{\mathbb{T}_{X/Y}}$ comes from a Hopf algebra structure:

$$\begin{array}{rcl} c^\sharp\colon \mathcal{A}_{\mathbb{T}_{X/Y}} & \to & \mathcal{A}_{\mathbb{T}_{X/Y}} \hat{\otimes}_{\mathcal{O}_{\mathbb{T}_X(Y)}} \mathcal{A}_{\mathbb{T}_{X/Y}} \\ \iota^\sharp\colon \mathcal{A}_{\mathbb{T}_{X/Y}} & \to & \mathcal{O}_{\mathbb{T}_X(Y)} \\ \tau^\sharp\colon \mathcal{A}_{\mathbb{T}_{X/Y}} & \to & \mathcal{A}_{\mathbb{T}_{X/Y}}. \end{array}$$

If $X$ is closed in $\overline{Y}$, the right action of $\mathcal{G}_{\mathbb{T}_Y}$ on $\mathbb{T}_X(Y)$ corresponds to a coaction of $\mathcal{A}_{\mathbb{T}_Y}$ on $\mathcal{O}_{\mathbb{T}_X(Y)}$:

$$r^\sharp\colon \mathcal{O}_{\mathbb{T}_X(Y)} \to \mathcal{O}_{\mathbb{T}_X(Y)} \hat{\otimes}_{\mathcal{O}_Y} \mathcal{A}_{\mathbb{T}_Y}.$$

The ring of *(hyper-)$\mathbb{T}$-differential operators* is defined by taking the $\mathcal{O}_{\mathbb{T}_X(Y)}$-linear (topological) dual:

$$\mathcal{HD}_{\mathbb{T}_{X/Y}} := \mathcal{H}om(\mathcal{A}_{\mathbb{T}_{X/Y}}, \mathcal{O}_{\mathbb{T}_X(Y)}) := \varprojlim \mathcal{H}om(\mathcal{A}_{\mathbb{T}_{X/Y,n}}, \mathcal{O}_{\mathbb{T}_X(Y)_n}), \tag{4.4.3}$$

with ring structure induced by the dual of $c^\sharp$. If $X = \overline{Y}$, we just write $\mathcal{HD}_{\mathbb{T}_Y}$ for $\mathcal{HD}_{\mathbb{T}_{X/Y}}$. The right action of $\mathcal{G}_{\mathbb{T}_Y}$ on $\mathbb{T}_X(Y)$ corresponds to a left action of $\mathcal{HD}_{\mathbb{T}_Y}$, given by

$$D(e) := (\mathrm{id} \hat{\otimes} D)(r^\sharp(e)),$$

for $e \in \mathcal{O}_{\mathbb{T}_X(Y)}$ and $D \in \mathcal{HD}_{\mathbb{T}_Y}$. If $\mathbb{T} = \mathbb{PD}$, the ring $\mathcal{HD}_{\mathbb{T}_Y}$ is topologically generated by derivations, and if $\mathbb{T} = \mathbb{A}$, by $p$-derivations; see [37], which we shall review below. This is not the case if $\mathbb{T} = \mathbb{F}$ or $\mathbb{D}$.

Statement (2) of Proposition 4.4.1 tells us that $\mathbb{T}_{X/Y}(1) \cong \mathbb{T}_X(Y) \times_Y \mathbb{T}_Y(1)$, hence that

$$\mathcal{A}_{\mathbb{T}_{X/Y}} \cong \mathcal{O}_{\mathbb{T}_X(Y)} \hat{\otimes}_{\mathcal{O}_Y} \mathcal{A}_{\mathbb{T}_Y}$$

and that

$$\mathcal{HD}_{\mathbb{T}_{X/Y}} \cong \mathcal{O}_{\mathbb{T}_X(Y)} \hat{\otimes}_{\mathcal{O}_Y} \mathcal{HD}_{\mathbb{T}_Y} \tag{4.4.4}$$

where the completion here means the inverse limit over the duals of the formal neighborhoods of the diagonal ideal in $\mathcal{A}_{\mathbb{T}_Y}$.

Let us make these constructions explicit in a local situation with the aid of a system of étale local coordinates. We suppose that $X/\overline{S}$ is smooth and $Y/S$ is $p$-completely smooth, that $S = \operatorname{Spf} R$, and that we have a system of local coordinates $(x_1, \ldots, x_n)$ for $Y$ such that the ideal of $X$ in $Y$ is generated by $(p, x_1, \ldots, x_r)$. The following straightforward computations appear in the literature, for example in [7], [37], [35], [40].

First of all, there are $p$-completely étale maps:

$$\begin{array}{rcl} Y & \to & \operatorname{Spf} R[x_1, \ldots, x_n]\hat{} \\ \mathbb{F}_X(Y) & \to & \operatorname{Spf} R[[x_1, \ldots, x_r]][x_{r+1}, \ldots, x_n]\hat{} \\ \mathbb{PD}_X(Y) & \to & \operatorname{Spf} R\langle x_1, \ldots, x_r\rangle[x_{r+1}, \ldots, x_n]\hat{} \\ \mathbb{D}_X(Y) & \to & \operatorname{Spf} R[t_1 \ldots, t_r][x_{r+1}, \ldots, x_n]\hat{} \\ , & & \text{where } x_i = pt_i \text{ for } 1 \leq i \leq r. \end{array} \tag{4.4.5}$$

The prismatic case is not so familiar or explicit; we can only say that, as a consequence of Theorem 2.3.3, there is a $p$-completely étale map:

$$\mathbb{\Delta}_X(Y) \to \operatorname{Spf} B^\infty[x_{r+1}, \ldots, x_n]\hat{}, \quad \text{where}$$

$$B^\infty := R[t_{i,j}]\hat{}/(t_{i,0}, pt_{i,j+1} - \delta^j(x_i) + t_{i,j}^p : i = 1, \ldots, r, j \in \mathbf{N}).$$

Let $\xi_i := 1 \otimes x_i - x_i \otimes 1$ in $\mathcal{O}_{Y(1)}$. Viewing $\mathcal{O}_{Y(1)}$ as an $\mathcal{O}_Y$-module via the first projection $t$, we identify $x_i \in \mathcal{O}_Y$ with $x_i \otimes 1 \in \mathcal{O}_{Y(1)}$. Then we have coordinates $x_1, \ldots, x_n, \xi_1, \ldots, \xi_n$ for $Y(1)$, and $(p, \xi_1, \ldots, \xi_n)$ is a sequence of generators for the ideal of $\overline{Y}$ in $Y(1)$. Thus:

$$\begin{array}{rcl} \mathbb{F}_Y(1) & = & \operatorname{Spf} \mathcal{O}_Y[[\xi_1, \ldots, \xi_n]] \\ \mathbb{PD}_Y(1) & = & \operatorname{Spf} \mathcal{O}_Y\langle \xi_1, \ldots, \xi_n\rangle\hat{} \\ \mathbb{D}_Y(1) & = & \operatorname{Spf} \mathcal{O}_Y[\eta_1, \ldots, \eta_n]\hat{}, \quad \text{where } \xi_i = p\eta_i \\ \mathbb{\Delta}_Y(1) & = & \operatorname{Spf} \mathcal{O}_Y\langle \eta_1 \ldots, \eta_n\rangle\hat{}, \quad \text{where } \xi_i = p\eta_i. \end{array} \tag{4.4.6}$$

Recall that the composition law $c\colon Y(1) \times_Y Y(1) \to Y(1)$ corresponds to the map $p_{13}\colon Y \times_S Y \times_S Y \to Y \times_S Y$. One then gets the following formulas for the Hopf algebra structures:

$$\begin{array}{rcl} c^\sharp\colon \mathcal{O}_{\mathbb{F}_Y(1)} \to \mathcal{O}_{\mathbb{F}_Y(2)} & : & \xi_i \mapsto \xi_i \otimes 1 + 1 \otimes \xi_i \\ c^\sharp\colon \mathcal{O}_{\mathbb{PD}_Y(1)} \to \mathcal{O}_{\mathbb{PD}_Y(2)} & : & \xi_i \mapsto \xi_i \otimes 1 + 1 \otimes \xi_i \\ c^\sharp\colon \mathcal{O}_{\mathbb{D}_Y(1)} \to \mathcal{O}_{\mathbb{D}_Y(2)} & : & \eta_i \mapsto \eta_i \otimes 1 + 1 \otimes \eta_i \\ c^\sharp\colon \mathcal{O}_{\mathbb{\Delta}_Y(1)} \to \mathcal{O}_{\mathbb{\Delta}_Y(2)} & : & \eta_i \mapsto \eta_i \otimes 1 + 1 \otimes \eta_i. \end{array} \tag{4.4.7}$$

For each multi-index $I := (I_1, \dots, I_n)$, define differential operators:

$$\begin{aligned} \partial_I \in \mathcal{HD}_{\mathbb{PD}_Y} \colon t_* \mathcal{O}_{\mathbb{PD}_Y(1)} \to \mathcal{O}_Y &\quad : \quad \xi^{[J]} \mapsto \delta_{I,J} \\ \partial_I \in \mathcal{HD}_{\mathbb{D}_Y} \colon t_* \mathcal{O}_{\mathbb{D}_Y(1)} \to \mathcal{O}_Y &\quad : \quad \eta^{J} \mapsto \delta_{I,J} \\ \partial_I \in \mathcal{HD}_{\mathbb{A}_Y} \colon t_* \mathcal{O}_{\mathbb{A}_Y(1)} \to \mathcal{O}_Y &\quad : \quad \eta^{[J]} \mapsto \delta_{I,J}. \end{aligned} \tag{4.4.8}$$

The set of these operators forms a formal basis for each $\mathcal{HD}_{\mathbb{T}_Y}$, in that every operator can be written uniquely as a formal sum $\sum_I a_I \partial_I$ with $a_I \in \mathcal{O}_Y$ (with no convergence conditions). Operators in $\mathcal{HD}_{\mathbb{T}_{X/Y}}$ can be a sum of the same form, now with coefficients in $\mathcal{O}_{\mathbb{T}_X(Y)}$. Furthermore we have the following composition laws:

$$\begin{aligned} \partial_I \partial_J &= \partial_{I+J} && \in \mathcal{HD}_{\mathbb{PD}_Y} \\ \partial_I \partial_J &= \frac{(I+J)!}{I!J!} \partial_{I+J} && \in \mathcal{HD}_{\mathbb{D}_Y} \\ \partial_I \partial_J &= \partial_{I+J} && \in \mathcal{HD}_{\mathbb{A}_Y} \end{aligned} \tag{4.4.9}$$

In particular, the rings $\mathcal{HD}_{\mathbb{PD}_Y}$ and $\mathcal{HD}_{\mathbb{A}_Y}$ are formally generated by operators of degree at most one, but this is not the case for $\mathcal{HD}_{\mathbb{D}_Y}$.

The right actions of $\mathcal{G}_{\mathbb{T}_Y}$ on $\mathbb{T}_X(Y)$ cover the action of $\mathcal{G}_Y$ on $Y$. These actions are given by the following formulas:

$$\begin{aligned} r^\sharp \colon \mathcal{O}_Y \to \mathcal{O}_{\mathbb{F}_Y(1)} &\quad : \quad x_i \mapsto x_i + \xi_i \\ r^\sharp \colon \mathcal{O}_{\mathbb{PD}_X(Y)} \to \mathcal{O}_{\mathbb{PD}_X(Y)} \hat{\otimes} \mathcal{O}_{\mathbb{PD}_Y(1)} &\quad : \quad x_i \mapsto x_i + \xi_i \\ r^\sharp \colon \mathcal{O}_{\mathbb{D}_X(Y)} \to \mathcal{O}_{\mathbb{D}_X(Y)} \hat{\otimes} \mathcal{O}_{\mathbb{D}_Y(1)} &\quad : \quad x_i \mapsto x_i + p\eta_i, \quad t_i \mapsto t_i + \eta_i \\ r^\sharp \colon \mathcal{O}_{\mathbb{A}_X(Y)} \to \mathcal{O}_{\mathbb{A}_X(Y)} \hat{\otimes} \mathcal{O}_{\mathbb{A}_Y(1)} &\quad : \quad x_i \mapsto x_i + p\eta_i, \quad t_i \mapsto t_i + \eta_i. \end{aligned}$$

These actions give rise to corresponding stratifications and consequently to actions of rings of differential operators. These in turn correspond to connections on $\mathcal{O}_Y$ and on $\mathcal{O}_{\mathbb{PD}_X(Y)}$ and to $p$-connections on $\mathcal{O}_{\mathbb{D}_X(Y)}$ and $\mathcal{O}_{\mathbb{A}_X(Y)}$, obtained by restricting to the action of first-order differential operators, or, equivalently, projecting from the ideal of the ideal $I$ of the diagonal to $I/I^2$, or to $I/I^{[2]}$ in the case of $\mathbb{PD}$ and $\mathbb{A}$. These connections are given by by the following formulas:

$$\begin{aligned} d \colon \mathcal{O}_Y \to \Omega^1_Y &\quad : \quad x_i^n \mapsto n x^{n-1} dx_i \\ d \colon \mathcal{O}_{\mathbb{PD}_X(Y)} \to \mathcal{O}_{\mathbb{PD}_X(Y)} \otimes \Omega^1_{Y/S} &\quad : \quad x_i^{[n]} \mapsto x_i^{[n-1]} dx_i \\ d' \colon \mathcal{O}_{\mathbb{D}_X(Y)} \to \mathcal{O}_{\mathbb{D}_X(Y)} \otimes \Omega^1_{Y/S} &\quad : \quad t_i^n \mapsto n t_i^{n-1} dx_i \\ d' \colon \mathcal{O}_{\mathbb{A}_X(Y)} \to \mathcal{O}_{\mathbb{A}_X(Y)} \otimes \Omega^1_{Y/S} &\quad : \quad t_{i,j} \mapsto d\delta^{j-1}(x_i) - t_{i,j-1}^{p-1} d' t_{i,j-1}. \end{aligned}$$

This last formula follows from the formula (4.2.2) for $d'$ in the proof of Proposition 4.2.1.

The *right* regular representations of $\mathcal{G}_Y$ on $Y(1)_s$ and $Y(1)_t$ described in Example 1.4.11 give rise to stratifications and connections which we claim are given by the following formulas:

$$
\begin{array}{rcll}
\nabla\colon \mathcal{O}_{\mathbb{F}_Y(1)_s} \to \mathcal{O}_{\mathbb{F}_Y(1)_s} \otimes \Omega^1_{Y/S} & : & x_i \otimes 1 \mapsto 0 & \\
 & & \xi_i \mapsto 1 \otimes 1 \otimes dx_i & (4.4.10) \\
\nabla\colon \mathcal{O}_{\mathbb{F}_Y(1)_t} \to \mathcal{O}_{\mathbb{F}_Y(1)_t} \otimes \Omega^1_{Y/S} & : & 1 \otimes x_i \mapsto 0 & \\
 & & \xi_i \mapsto -1 \otimes 1 \otimes dx_i & (4.4.11) \\
\nabla\colon \mathcal{O}_{\mathbb{A}_Y(1)_s} \to \mathcal{O}_{\mathbb{A}_Y(1)_s} \otimes \Omega^1_{Y/S} & : & x_i \otimes 1 \mapsto 0 & \\
 & & \eta_i \mapsto 1 \otimes 1 \otimes dx_i & (4.4.12) \\
\nabla\colon \mathcal{O}_{\mathbb{A}_Y(1)_t} \to \mathcal{O}_{\mathbb{A}_Y(1)_t} \otimes \Omega^1_{Y/S} & : & 1 \otimes x_i \mapsto 0 & \\
 & & \eta_i \mapsto -1 \otimes 1 \otimes dx_i. & (4.4.13)
\end{array}
$$

It suffices to check these formulas when $\mathbb{T} = \mathbb{F}$. We begin with the following formulas for the right actions of $\mathcal{G}_Y$ on $\mathbb{F}_Y(1)_s$ and $\mathbb{F}_Y(1)_t$, which are consequences of the descriptions in equation (1.4.9).

$$
\begin{array}{rcl}
r_s^*\colon \mathcal{O}_{\mathbb{F}_Y(1)_s} \to \mathcal{O}_{\mathbb{F}_Y(1)_s} \otimes \mathcal{O}_{\mathbb{F}_Y(1)} & : & a \otimes b \mapsto a \otimes 1 \otimes 1 \otimes b \\
r_t^*\colon \mathcal{O}_{\mathbb{F}_Y(1)_t} \to \mathcal{O}_{\mathbb{F}_Y(1)_t} \otimes \mathcal{O}_{\mathbb{F}_Y(1)} & : & a \otimes b \mapsto 1 \otimes b \otimes 1 \otimes a.
\end{array}
$$

Recalling that $1 \otimes x_i = x_i \otimes 1 + \xi_i$, we find

$$
\begin{array}{rcl}
r_s^*\colon x_i \otimes 1 & \mapsto & x_i \otimes 1 \otimes 1 \otimes 1 \\
1 \otimes x_i & \mapsto & 1 \otimes 1 \otimes 1 \otimes x_i \\
\xi_i & \mapsto & 1 \otimes 1 \otimes 1 \otimes x_i - x_i \otimes 1 \otimes 1 \otimes 1 \\
 & = & 1 \otimes 1 \otimes x_i \otimes 1 + 1 \otimes 1 \otimes \xi_i - x_i \otimes 1 \otimes 1 \otimes 1 \\
 & = & 1 \otimes x_i \otimes 1 \otimes 1 + 1 \otimes 1 \otimes \xi_i - x_i \otimes 1 \otimes 1 \otimes 1 \\
 & = & \xi_i \otimes 1 \otimes 1 + 1 \otimes 1 \otimes \xi_i.
\end{array}
$$

When calculating $r_t$, we must remember to use the structure on $\mathbb{F}_Y(1)$ coming from the left factor. Thus:

$$
\begin{array}{rcl}
r_t\colon 1 \otimes x_i & \mapsto & 1 \otimes x_i \otimes 1 \otimes 1 \\
x_i \otimes 1 & \mapsto & 1 \otimes 1 \otimes 1 \otimes x_i \\
 & = & 1 \otimes 1 \otimes x_i \otimes 1 + 1 \otimes 1 \otimes \xi_i \\
 & = & x_i \otimes 1 \otimes 1 \otimes 1 + 1 \otimes 1 \otimes \xi_i \\
\xi_i & \mapsto & 1 \otimes x_i \otimes 1 \otimes 1 - x_i \otimes 1 \otimes 1 \otimes 1 - 1 \otimes 1 \otimes \xi_i \\
 & = & \xi_i \otimes 1 \otimes 1 - 1 \otimes 1 \otimes \xi_i.
\end{array}
$$

We deduce that:

$$
\begin{array}{rcl}
\nabla_s(x_i \otimes 1) & = & 0 \\
\nabla_s(\xi_i) & = & 1 \otimes 1 \otimes dx_i \\
\nabla_t(1 \otimes x_i) & = & 0 \\
\nabla_t(\xi_i) & = & -1 \otimes 1 \otimes dx_i,
\end{array}
$$

as claimed.

The following proposition, a consequence of the above discussion, summarizes what we shall need.

**Proposition 4.4.2.** *Let $Y/S$ be a $p$-completely smooth morphism of $p$-torsion free formal $\delta$-schemes.*

(1) *The reduction of* $\mathcal{HD}_{\mathbb{A}_Y}$ *modulo* $p$ *is canonically isomorphic to the completion of* $S^\bullet \mathcal{T}_{\overline{Y}/S}$ *along the ideal* $S^+\mathcal{T}_{\overline{Y}/S}$ *of the zero section. Thus there are ring homomorphisms:*

$$\mathcal{HD}_{\mathbb{A}_Y} \to \hat{S}^\bullet \mathcal{T}_{\overline{Y}/S} \to \mathcal{O}_{\overline{Y}}.$$

(2) *If* $X$ *is a closed subscheme of* $\overline{Y}$ *which is smooth over* $\overline{S}$*, there is a natural isomorphism*

$$\mathcal{O}_{\mathbb{A}_{X(Y)}} \hat{\otimes}_{\mathcal{O}_Y} \mathcal{HD}_{\mathbb{A}_Y} \to \mathcal{HD}_{\mathbb{A}_{X/Y}},$$

*where the completion is taken with respect to the kernel of the homomorphism* $\mathcal{HD}_{\mathbb{A}_Y} \to \mathcal{O}_Y$.

*Proof.* Definition 4.4.3 implies that

$$\mathcal{HD}_{\mathbb{A}_Y}/p\mathcal{HD}_{\mathbb{A}_Y} \cong \mathcal{H}om(t_*(\mathcal{O}_{\overline{\mathbb{A}}_{Y(1)}}), \mathcal{O}_{\overline{Y}}).$$

By Theorem 2.3.13, $\overline{\mathbb{A}}_Y(1)$ is isomorphic to the PD-envelope of the canonical section of the dilatation of the diagonal of $\overline{Y}(1)$, which is in turn isomorphic (via the map induced by multiplication by $p$), to the PD-envelope of the diagonal of $\overline{Y}(1)$. The dual of this PD-envelope is the completion of the symmetric algebra of the tangent sheaf, thus proving statement (1). Then statement (2) follows from equation (4.4.4). □

To prepare for a discussion of prismatic crystals, we let $\mathbf{F}\colon \mathbf{V} \to \mathbb{T}(X/S)$ be the fibered category whose objects are pairs $(T, E_T)$, where $T$ is an object of $\mathbb{T}(X/S)$, where $E_T$ is a $p$-completely quasi-coherent sheaf of $\mathcal{O}_T$-modules on $T$ (Definition 7.4.1), and where $\mathbf{F}(T, E_T) := T$. Since we may also want to work geometrically, we let $\mathbf{V}E_T = \mathrm{Spec}_T\, S^\bullet E_T$.

The following result, whose (omitted) proof follows from the above discussions and the methods described, for example in [7] and [37], summarizes various ways of understanding stratifications. For clarity, we state it just for prismatic crystals.

**Theorem 4.4.3.** *Let* $Y/S$ *be a* $p$*-completely smooth morphism of crystalline prisms and let* $X \to Y$ *be a closed immersion, where* $X/\overline{S}$ *is smooth. Suppose that* $E \in \mathbf{V}_{\mathbb{A}_{X(Y)}}$ *is a* $p$*-completely quasi-coherent sheaf of* $\mathcal{O}_{\mathbb{A}_{X(Y)}}$*-modules on* $Y$*. Then the following sets of data are equivalent:*

(1) *A right action of* $\mathcal{G}_{\mathbb{A}_{X/Y}}$ *on* $\mathbf{V}E$.

(2) *A* $\mathcal{G}_{\mathbb{A}_{X/Y}}$*-stratification on* $\mathbf{V}E$.

(3) *A left action of the ring* $\mathcal{HD}_{\mathbb{A}_{X/Y}}$ *on* $E$,

(4) *A left action of the ring* $\mathcal{HD}_{\mathbb{A}_{X/Y}}$ *on* $E$ *compatible with its action on* $\mathcal{O}_{\mathbb{A}_{X(Y)}}$.

(5) *An integrable quasi-nilpotent $p$-connection $\nabla'$ on $E$, compatible with the $p$-connection on $\mathcal{O}_{\mathbb{\Delta}_X(Y)}$.*

□

Let us now discuss more general prismatic differential operators and their linearization, following the framework described in §1.5. Note that there we used geometric notation, so the arrows in our algebraic constructions here are the reverse of those there.

**Definition 4.4.4.** *Let $Y/S$ be a $p$-completely smooth morphism of crystalline prisms, and let $t, s\colon \mathbb{\Delta}_Y(1) \to Y$ be the two projections.*

(1) *If $\Omega$ is a sheaf of $\mathcal{O}_Y$-modules, let*

$$\mathcal{L}_{\mathbb{\Delta}}(\Omega) := t_* s^*(\Omega) = \mathcal{O}_{\mathbb{\Delta}_Y(1)} \hat{\otimes} \Omega.$$

(2) *If $\Omega$ and $\Omega'$ are sheaves of $\mathcal{O}_Y$-modules, a (hyper) prismatic differential operator from $\Omega$ to $\Omega'$ is an $\mathcal{O}_Y$-linear map $t_*s^*(\Omega) \to \Omega'$, i.e., a map*

$$D\colon \mathcal{L}_{\mathbb{\Delta}}(\Omega) = \mathcal{O}_{\mathbb{\Delta}_Y(1)} \hat{\otimes}_{\mathcal{O}_Y} \Omega \to \Omega'.$$

(3) *If $D$ is a prismatic differential operator from $\Omega$ to $\Omega'$, then $\mathcal{L}_{\mathbb{\Delta}}(D)$ is the $\mathcal{O}_Y$-linear map:*

$$\mathcal{L}_{\mathbb{\Delta}}(D) := \mathcal{O}_{\mathbb{\Delta}_Y(1)} \hat{\otimes} \Omega \xrightarrow{\delta \otimes \mathrm{id}} \mathcal{O}_{\mathbb{\Delta}_Y(1)} \hat{\otimes} \mathcal{O}_{\mathbb{\Delta}_Y(1)} \hat{\otimes} \Omega \xrightarrow{\mathrm{id} \otimes D} \mathcal{O}_{\mathbb{\Delta}_Y(1)} \hat{\otimes} \Omega'.$$

(4) *If $E$ is an $\mathcal{O}_Y$-module with prismatic stratification $\epsilon$ and $D$ is a prismatic differential operator from $\Omega$ to $\Omega'$, then $\epsilon(D)$ is the prismatic differential operator from $E \hat{\otimes} \Omega$ to $E \hat{\otimes} \Omega'$ defined by the following diagram:*

$$\begin{array}{ccc} \mathcal{L}_{\mathbb{\Delta}}(E \hat{\otimes} \Omega) & \xrightarrow{\epsilon(D)} & E \hat{\otimes} \Omega' \\ \cong \downarrow & & \uparrow \mathrm{id} \hat{\otimes} D \\ \mathcal{O}_{\mathbb{\Delta}_Y(1)} \hat{\otimes} E \otimes \Omega & \xrightarrow{\epsilon \hat{\otimes} \mathrm{id}} & E \hat{\otimes} \mathcal{O}_{\mathbb{\Delta}_Y(1)} \hat{\otimes} \Omega. \end{array}$$

Since $\mathcal{L}_{\mathbb{\Delta}}(\Omega)$ is computed using $s$ but is endowed with the $\mathcal{O}_Y$-module structure coming from $t$, the $p$-connection on $t_*\mathcal{O}_{\mathbb{\Delta}_Y(1)}$ described in Proposition 4.2.1 carries over to $\mathcal{L}_{\mathbb{\Delta}}(\Omega)$ and $\mathcal{L}_{\mathbb{\Delta}}(\Omega')$, and if $D$ is a prismatic differential operator from $\Omega$ to $\Omega'$, the corresponding $\mathcal{O}_Y$-linear map $\mathcal{L}_{\mathbb{\Delta}}(D)$ is horizontal.

There is a commutative diagram:

$$\begin{array}{ccc} \Omega & & \\ s^* \downarrow & \searrow^{\mathrm{id}} & \\ \mathcal{L}_{\mathbb{\Delta}}(\Omega) & \xrightarrow{\iota^*} & \Omega, \end{array} \tag{4.4.14}$$

where $\iota^*$ is induced by the diagonal mapping and $s^*(\omega) := 1 \otimes \omega$. When $\Omega = \mathcal{O}_Y$, the map $s$ is compatible with the prismatic connections. (See diagram (1.4.6)).

The following prismatic variant of the crystalline construction [7, 6.15] then follows from Proposition 1.5.2.

**Proposition 4.4.5.** *If $E$ and $\Omega$ are $p$-completely quasi-coherent $\mathcal{O}_Y$-modules and $\epsilon$ is a prismatic stratification on $E$, define $\beta$ by the following commutative diagram:*

$$\begin{array}{ccc} E \hat{\otimes} \mathcal{L}_{\Delta}(\Omega) & \xrightarrow{\beta} & \mathcal{L}_{\Delta}(E \hat{\otimes} \Omega) \\ \downarrow \cong & & \downarrow \cong \\ E \hat{\otimes} \mathcal{O}_{\overline{\Delta}_Y(1)} \hat{\otimes} \Omega & \xrightarrow{\epsilon \otimes \mathrm{id}} & \mathcal{O}_{\overline{\Delta}_Y(1)} \hat{\otimes} E \hat{\otimes} \Omega \end{array}$$

*(1) The map $\beta$ is a horizontal isomorphism.*

*(2) If $D$ is a prismatic differential operator from $\Omega$ to $\Omega'$, the following diagram commutes:*

$$\begin{array}{ccc} E \hat{\otimes} \mathcal{L}_{\Delta}(\Omega) & \xrightarrow{\beta} & \mathcal{L}_{\Delta}(E \hat{\otimes} \Omega) \\ {\scriptstyle \mathrm{id} \hat{\otimes} \mathcal{L}_{\Delta}(D)} \downarrow & & \downarrow {\scriptstyle \mathcal{L}_{\Delta}(\epsilon(D))} \\ E \hat{\otimes} \mathcal{L}_{\Delta}(\Omega') & \xrightarrow{\beta} & \mathcal{L}_{\Delta}(E \hat{\otimes} \Omega'). \end{array}$$

The differentials in the $p$-de Rham complex of $Y/S$ are prismatic differential operators; see [4, Theorem 3.1.2] for the proof of the analogous statement for PD differential operators. As in the crystalline model explained in [7, 6.11], we have the following description of their linearization.

**Proposition 4.4.6.** *Let $Y/S$ be a $p$-completely smooth morphism of crystalline prisms. View $\Delta_Y(1)$ as a formal scheme over $Y$ via the morphism $t$. Then there is a natural isomorphism of complexes:*

$$(\Omega^\bullet_{\Delta_Y(1)/Y}, d') \to \mathcal{L}_{\Delta}(\Omega^\bullet_{Y/S}, d').$$

*In particular, if $Y/S$ admits a system of local coordinates $(x_1, \ldots, x_n)$ and $(\eta_1, \ldots, \eta_n) \in \mathcal{O}_{\Delta_Y(1)}$ are as in Equation (4.4.6), then the homomorphism*

$$\mathcal{L}_{\Delta}(d') \colon \mathcal{L}_{\Delta}(\Omega^i_{Y/S}) \to \mathcal{L}_{\Delta}(\Omega^{i+1}_{Y/S})$$

*is given by*

$$\eta_1^{[m_1]} \cdots \eta_n^{[m_n]} \otimes \omega \mapsto \sum_k \eta_1^{[m_1]} \cdots \eta_k^{[m_k - 1]} \eta_n^{[m_n]} \otimes dx_k \wedge \omega + \eta_1^{[m_1]} \cdots \eta_n^{[m_n]} \otimes d\omega.$$

*Proof.* As noted in diagram (1.5.1) in the discussion of the the general construction, there is a commutative diagram:

$$
\begin{array}{ccccccc}
\mathcal{O}_{\mathbb{\Delta}_{Y(1)}} & \xrightarrow{\mathcal{L}(d')} & \mathcal{O}_{\mathbb{\Delta}_Y} \otimes \Omega^1_{Y(1)/Y} & \xrightarrow{\mathcal{L}(d')} & \mathcal{O}_{\mathbb{\Delta}_Y} \otimes \Omega^2_{Y(1)/Y} & \xrightarrow{d'} & \cdots \\
\uparrow & & \uparrow & & \uparrow & & \\
\mathcal{O}_{Y(1)} & \xrightarrow{d'} & \Omega^1_{Y(1)/Y} & \xrightarrow{d'} & \Omega^2_{Y(1)/Y} & \xrightarrow{d'} & \cdots \\
{\scriptstyle t^\sharp}\uparrow & & {\scriptstyle t^*}\uparrow & & \uparrow{\scriptstyle t^*} & & \\
\mathcal{O}_Y & \xrightarrow{d'} & \Omega^1_{Y/S} & \xrightarrow{d'} & \Omega^2_{Y/S} & \xrightarrow{d'} & \cdots .
\end{array}
$$

As we have argued before, $\mathcal{O}_{\mathbb{\Delta}_{Y(1)}}$ is locally topologically generated by elements $a$ such that $p^n a$ lies in $\mathcal{O}_Y$ for some $n$, and since the terms of the top complex are $p$-torsion free and $p$-adically separated, its differentials $\mathcal{L}(d')$ are uniquely determined by the commutativity of the diagram. Thus they must agree with the differentials of the complex $(\Omega^{\bullet}_{\mathbb{\Delta}_{Y(1)/Y}}, d')$. In the presence of local coordinates, we have $p\eta_i = s^\sharp(x_i) - t^\sharp(x_i)$, so $d'\eta_i = pd\eta_i = dx_i$. The formula in the proposition follows.

□

# 5 The $F$- and $\alpha$-transforms

In this section we discuss two operations relating $p$-connections, connections, and Higgs fields. The $F$-transform, relates $p$-connections and connections, and is somehow inspired by Simpson's groundbreaking work on nonabelian Hodge theory [38]. The $\alpha$-transform first appeared in [33] and [34] in the context of $p$-curvature calculations, but it plays an additional and important role here, as we shall see in Theorem 5.4.3.

## 5.1 The $F$-transform

We recall the following construction, which has been used by many authors, for example, in [34], [37], [17], and [40]. The starting point is easily checked.

**Proposition 5.1.1.** *Let $Y/S$ be a $p$-completely smooth morphism of crystalline prisms, let $F_{Y/S}\colon Y \to Y^{(p)}$ be the relative Frobenius morphism,* [7] *let $(E', \nabla')$ be a $p$-completely quasi-coherent module with integrable $p$-connection on $Y^{(p)}/S$, and let*

$$\zeta \quad := \quad p^{-1}F^*_{Y/S} \quad : \quad \Omega^1_{Y^{(p)}/S} \quad \to \quad F_{Y/S*}(\Omega^1_{Y/S})$$

[7] In fact the construction requires less than what is provided by prisms: one only needs a lifting of the relative Frobenius morphism.

$$\eta \quad := \quad e' \mapsto 1 \otimes e' \quad : \quad E' \quad \to \quad F^*_{Y/S}(E') \qquad (5.1.1)$$

*be the natural maps.*

*(1) There is a unique integrable connection*

$$\nabla \colon F^*_{Y/S}(E') \to F^*_{Y/S}(E') \otimes \Omega^1_{Y/S}$$

*with the property that*

$$\nabla(\eta(e')) = (\zeta \otimes \eta)(\nabla'(e'))$$

*for every local section $e'$ of $E'$. We call $(F^*_{Y/S}(E'), \nabla)$ the $F$-transform of $(E', \nabla')$.*

*(2) There is a morphism of complexes*

$$\zeta^{\bullet}_E \colon (E' \otimes \Omega^{\bullet}_{Y^{(p)}/S}, d') \to F_{Y/S*}(F^*_{Y/S}(E') \otimes \Omega^{\bullet}_{Y/S}, d)$$

*defined by*

$$\begin{array}{ccccccc}
E' & \xrightarrow{\nabla'} & E' \otimes \Omega^1_{Y^{(p)}/S} & \xrightarrow{d'} & E' \otimes \Omega^2_{Y^{(p)}/S} & \longrightarrow & \cdots \\
\Big\downarrow \eta & & \Big\downarrow \eta \otimes \zeta & & \Big\downarrow \eta \otimes \wedge^2 \zeta & & \\
F^*_{Y/S}(E') & \xrightarrow{\nabla} & F^*_{Y/S}(E') \otimes \Omega^1_{Y/S} & \xrightarrow{d} & F^*_{Y/S}(E') \otimes \Omega^2_{Y/S} & \longrightarrow & \cdots
\end{array}$$

□

Connections, $p$-connections, and the $F$-transform can also be interpreted geometrically. Recall from Proposition 4.1.1 that we are writing $P^1_{Y/S}$ (resp. $D^1_{Y/S}$) for the first infinitesimal neighborhood of $Y \to Y(1)$ (resp. of $Y \to \mathbb{D}_Y(1)$. As is well-known, to give a connection on an $\mathcal{O}_Y$-module $E$ is the same as to give an isomorphism of $\mathcal{O}_{P^1_{Y/S}}$-modules $\epsilon \colon p_2^*(E) \to p_1^*(E)$ whose restriction to the diagonal is the identity. The isomorphism $\epsilon$ and the connection $\nabla$ are related by the formula

$$\epsilon(p_2^*(e)) = p_1^*(e) + \nabla(e).$$

Similarly, giving a $p$-connection on an $\mathcal{O}_Y$-module $E$ is equivalent to giving an isomorphism of $\mathcal{O}_{D^1_{Y/S}}$-modules $p_2^*(E) \to p_1^*(E)$ whose restriction to the diagonal is the identity. Using these interpretations of connections and $p$-connections, we find the following geometric interpretation of the F-transform. The proof is immediate from the definitions and diagram (4.1.2) of Proposition 4.1.1.

**Proposition 5.1.2.** *With the notation of Proposition 5.1.1, let*

$$\epsilon' \colon {p'_2}^*(E') \to {p'_1}^*(E')$$

*be the isomorphism of $\mathcal{O}_{D^1_{Y^{(p)}/S}}$-modules corresponding to a p-connection $\nabla'$ on $E'$. Let $\tilde{F}\colon P^1_{Y/S} \to D^1_{Y^{(p)}/S}$ be the morphism*[8] *defined in Proposition 4.1.1. Then the F-transform of $(E', \nabla')$ corresponds to the isomorphism of $\mathcal{O}_{P^1_{Y/S}}$-modules:*

$$\tilde{F}^*(\epsilon')\colon \tilde{F}^*({p'_2}^*(E')) \to \tilde{F}^*({p'_1}^*(E')),$$

*after the identifications $\tilde{F}^*{p'_i}^* \cong p_i^* F^*$.* □

**Remark 5.1.3.** The $F$-transform of $E'$ depends on the choice of Frobenius lifting. Two Frobenius liftings $F_1$ and $F_2$ agree modulo $p$, and since the ideal $(p)$ has a unique PD-structure, an integrable connection on $E'$ would define an isomorphism $F_2^*E' \cong F_1^*E'$, but a $p$-connection is not enough to do this on its own. This issue was studied in [34], [35], and [40]. See Proposition 5.6.1 for a review.

The following important theorem is due to A. Shiho [37, 3.1]. As he explains, it is a "descent to level minus one" version of Berthelot's Frobenius descent ideas [6]. We shall review his proof, using the point of view developed in our current context, in subsection 5.6.

**Theorem 5.1.4** (Shiho)**.** *Let $F\colon Y \to Y^{(p)}$ be as in Proposition 5.1.1, and let $(E, \nabla)$ be the F-transform of an integrable p-connection $(E', \nabla')$ on $Y^{(p)}/S$. Then $\nabla$ is quasi-nilpotent if $\nabla'$ is. The F-transform $(E', \nabla') \mapsto (E, \nabla)$ defines an equivalence from the tensor category* $\mathrm{MICP}(Y^{(p)}/S)$ *of quasi-nilpotent integrable p-connections on $Y^{(p)}$ to the tensor category* $\mathrm{MIC}(Y/S)$ *of quasi-nilpotent integrable connections on $Y$.* □

**Remark 5.1.5.** Shiho shows in [37] that the F-transform of an integrable quasi-nilpotent $p$-connection is a quasi-nilpotent connection. We shall give a more precise version of his argument by computing the $p$-curvature of the F-transform explicitly; see Theorem 5.3.1 below. We should also point out that the converse result is not true: the F-transform of an integrable non-quasi-nilpotent $p$-connection can be quasi-nilpotent. In particular, the F-transform functor is not fully faithful on the category of all integrable $p$-connections.

The following result is a strengthening of a special case of a result of Shiho [37, 4.4] relating the cohomology of a module with integrable quasi-nilpotent $p$-connection to the cohomology of its F-transform.

**Theorem 5.1.6.** *With the notation of Theorem 5.1.4, let $(E', \nabla')$ be a p-completely quasi-coherent sheaf of $\mathcal{O}_{Y^{(p)}}$-modules with quasi-nilpotent integrable p-connection and let $(E, \nabla)$ be its F-transform. Then the natural morphism of complexes*

$$\zeta_E^\bullet\colon (E' \otimes \Omega^\bullet_{Y^{(p)}/S}, d') \to (E \otimes \Omega^\bullet_{Y/S}, d)$$

*is a strict quasi-isomorphism (see Definition 7.4.1).*

[8]The treatment there used the absolute instead of relative Frobenius, but the construction here is essentially the same.

*Proof.* First suppose that $(E', \nabla') = (\mathcal{O}_{\overline{Y}}, 0)$. Then $(E, \nabla) = (\mathcal{O}_{\overline{Y}}, d)$, and our assertion is that

$$\zeta_E^\bullet \colon (\Omega^\bullet_{\overline{Y}^{(p)}/S}, 0) \to F_{Y/S*}(\Omega^\bullet_{\overline{Y}/S}, d)$$

is a quasi-isomorphism. This follows from the fact that in each degree $i$, $\zeta_E$ reduces to the inverse Cartier isomorphism:

$$C^{-1}_{\overline{Y}/S} \colon \Omega^i_{\overline{Y}^{(p)}/S} \to F_{Y/S*}\mathcal{H}^i(\Omega^\bullet_{\overline{Y}/S}).$$

Now let $E'$ be any quasi-coherent sheaf of $\mathcal{O}_{\overline{Y}^{(p)}}$-modules. Since all the terms of the complex $F_{Y/S*}(\Omega^\bullet_{\overline{Y}/S})$, as well as its cohomology sheaves, are locally free over $\mathcal{O}_{\overline{Y}^{(p)}}$, it follows that the map

$$\zeta_E^\bullet \colon (E' \otimes \Omega^\bullet_{\overline{Y}^{(p)}/S}, 0) \to F_{Y/S*}(E' \otimes_{\mathcal{O}_{Y^{(p)}}} \Omega^\bullet_{\overline{Y}/S}, d \otimes \mathrm{id})$$

is still a quasi-isomorphism. The right side is the de Rham complex of the F-transform of $(E', 0)$. Thus we have proved the theorem whenever $E'$ is annihilated by $p$ and by $\nabla'$.

Note next that, since $F$ is flat, the $F$-transform preserves exact sequences. Hence if $0 \to (E'_1, \nabla') \to (E'_2, \nabla') \to (E'_3, \nabla') \to 0$ is an exact sequence of quasi-coherent $p^n$-torsion modules with integrable and quasi-nilpotent $p$-connection and the theorem is true for any two of the terms, then it is true for the third. Consequently the result also holds if $E'$ is $p^n$-torsion and has a finite filtration whose associated graded satisfies the claimed result.

Let us now check that the theorem is true whenever $E'$ is annihilated by $p$. We may work locally, and in particular we may assume that $Y$ is affine. Let $E'_1 := \mathrm{Ker}(\nabla')$, let $E'_2$ be the inverse image in $E'$ of the kernel of $\nabla'$ acting on $E'/E'_1$, etc. We obtain an increasing filtration $E'_1 \subseteq E'_2 \subseteq \cdots$ on whose associated graded $\nabla'$ is trivial. It follows that the result holds for each $E'_n$. If $E'$ is noetherian, this sequence terminates, and $E'_n = E'$ for $n \gg 0$ because $\nabla'$ was assumed to be locally nilpotent. In any case, $E' = \varinjlim E'_n$, and since formation of F-transforms and of cohomology commute with direct limits, the result also holds for $E'$.

It now follows by induction that the theorem is true for any $E'$ annihilated by a finite power of $p$. To check the induction step, assume that $p^{n+1}E' = 0$, so that there is an exact sequence:

$$0 \to p^n E' \to E' \to E'/p^n E' \to 0.$$

The first term in this sequence is annihilated by $p$ and is a quotient of $E'/pE'$, hence its $p$-connection is quasi-nilpotent, and hence the theorem holds for this term. Since the $p$-connection of the last term is also quasi-nilpotent, the theorem also holds for it, and hence also for $E'$.

To complete the proof of the theorem, recall that, by definition, a morphism of complexes $p$-completely quasi-coherent sheaves is a strict quasi-isomorphism if and only if its reduction modulo each power $p$ is a quasi-isomorphism. Since the $F$-transform commutes with reduction modulo each power of $p$, and we have just checked that the result is true for every such reduction, the theorem is proved. □

The following result relates the Frobenius pull-back of a module with $p$-connection $(E', \nabla')$ on $Y^{(p)}/S$ to the $p$-transform (explained in Example 1.2.2) of the F-transform of $(E', \nabla')$. It will allow us to see that the prismatic Frobenius is an isogeny on cohomology.

**Proposition 5.1.7.** *Let $Y/S$ be a $p$-completely smooth morphism of crystalline prisms, and let $(E', \nabla')$ be an object of* $\mathrm{MICP}(Y^{(p)}/S)$.

*(1) The $F_{Y/S}$-pullback $(E'', \nabla'')$ of $(E', \nabla')$ as a module with $p$-connection is the $p$-transform (see Example 1.2.2) of the F-transform $(E, \nabla)$ of $(E', \nabla')$.*

*(2) There exists a commutative diagram of complexes:*

$$\begin{array}{ccc} (E' \otimes \Omega^{\bullet}_{Y^{(p)}/S}, d') & \xrightarrow{\zeta^{\bullet}_E} & F_{Y/S*}(E \otimes \Omega^{\bullet}_{Y/S}, d) \\ & \searrow^{c} & \downarrow F_{Y/S*}(b) \\ & & F_{Y/S*}(E'' \otimes \Omega^{\bullet}_{Y/S}, d'), \end{array}$$

*where $\zeta^{\bullet}_E$ is the quasi-isomorphism of Theorem 5.1.6, where $b$ is the isogeny defined in Example 1.2.2, and where $c$ is the canonical map induced by the morphism $F_{Y/S}$.*

*Proof.* Let $\eta\colon E' \to F_{Y/S*}(E) = E''$ be the adjunction map, as in (5.1.1). Then we have a commutative diagram:

$$\begin{array}{ccc} E' & \xrightarrow{\nabla'} & E' \otimes \Omega^1_{Y^{(p)}/S} \\ \downarrow \eta & & \downarrow \eta \otimes \zeta \\ F_{Y/S*}(E) & \xrightarrow{F_{Y/S*}(\nabla)} & F_{Y/S*}(E \otimes \Omega^1_{Y/S}) \\ \downarrow \mathrm{id} & & \downarrow \mathrm{id} \otimes p \\ F_{Y/S*}(E'') & \xrightarrow{F_{Y/S*}(\nabla'')} & F_{Y/S*}(E'' \otimes \Omega^1_{Y/S}). \end{array}$$

The composition of the vertical arrows on the left is $\eta$, and on the right is $\eta \otimes F^*_{Y/S}$, since $F^*_{Y/S} = p\zeta$. Then the rectangle commutes, by the definition of the pullback $p$-connection $\nabla''$. The top square commutes by the definition of the F-transform connection $\nabla$, and it follows that the bottom square commutes as well. This implies that $\nabla'' = p\nabla$, so $\nabla''$ is the $p$-transform of $\nabla$.

The diagram above extends to a diagram of complexes:

$$\begin{array}{ccccccc}
E' & \xrightarrow{d'} & E'\otimes\Omega^1_{Y^{(p)}/S} & \xrightarrow{d'} & E\otimes\Omega^2_{Y^{(p)}/S} & \xrightarrow{d'} & \cdots \\
\Big\downarrow{\scriptstyle \eta} & & \Big\downarrow{\scriptstyle \eta\otimes p^{-1}F^*_{Y/S}} & & \Big\downarrow{\scriptstyle \eta\otimes p^{-2}F^*_{Y/S}} & & \\
F_{Y/S*}(E) & \xrightarrow{d} & F_{Y/S*}(E\otimes\Omega^1_{Y/S}) & \xrightarrow{d} & F_{Y/S*}(E\otimes\Omega^2_{Y/S}) & \xrightarrow{d} & \cdots \\
\Big\downarrow{\scriptstyle \mathrm{id}} & & \Big\downarrow{\scriptstyle p} & & \Big\downarrow{\scriptstyle p^2} & & \\
F_{Y/S*}(E) & \xrightarrow{pd} & F_{Y/S*}(E\otimes\Omega^1_{Y/S}) & \xrightarrow{pd} & F_{Y/S*}(E\otimes\Omega^2_{Y/S}) & \xrightarrow{pd} & \cdots,
\end{array}$$

which is the content of statement (2). □

We now suppose that $X$ is a closed subscheme of $\overline{Y}$ and is smooth over $\overline{S}$. We have seen in Proposition 4.2.1 that the structure sheaf of the prismatic neighborhood $\mathbb{\Delta}_{X^{(p)}}(Y^{(p)})$ carries a $p$-connection, while that of the divided power neighborhood $\mathbb{PD}_X(Y)$ carries a connection. We shall see that the F-transform takes the former to the latter.

Recall from Theorem 2.3.13 that if $X^\phi := F^{-1}_{Y/S}(X^{(p)})$, then $\mathbb{\Delta}_{X^\phi}(Y) \cong \mathbb{PD}_X(Y)$, and, from Corollary 2.4.1 and diagram (2.4.3) in its proof, that there is a commutative diagram:

$$\begin{array}{ccccc}
\mathbb{\Delta}_{X^\phi}(Y)\cong\mathbb{PD}_X(Y) & \xrightarrow{\Phi_{Y/S}} & \mathbb{\Delta}_{X^{(p)}}(Y^{(p)}) & \longrightarrow & \mathbb{\Delta}_X(Y) \\
\Big\downarrow & & \Big\downarrow & & \Big\downarrow \\
Y & \xrightarrow{F_{Y/S}} & Y^{(p)} & \xrightarrow{\pi_{Y/S}} & Y,
\end{array} \tag{5.1.2}$$

whose left square is Cartesian. (If $S$ is perfect or regular, the right square is also Cartesian, since in those cases $F_S$ is flat.)

**Theorem 5.1.8.** *With the notation and hypotheses described in the previous paragraphs, the canonical connection*

$$d\colon \mathcal{O}_{\mathbb{PD}_X(Y)} \to \mathcal{O}_{\mathbb{PD}_X(Y)} \otimes \Omega^1_{Y/S}$$

*on $\mathcal{O}_{\mathbb{PD}_X(Y)}$ is the F-transform (5.1.1) of the canonical $p$-connection*

$$d'\colon \mathcal{O}_{\mathbb{\Delta}_{X^{(p)}}(Y^{(p)})} \to \mathcal{O}_{\mathbb{\Delta}_{X^{(p)}}(Y^{(p)})} \otimes \Omega^1_{Y^{(p)}/S}$$

*on $\mathcal{O}_{\mathbb{\Delta}_{X^{(p)}}(Y^{(p)})}$. Furthermore, the F-transform defines an equivalence of categories:*

$$\mathrm{MICP}(X^{(p)}/Y^{(p)}/S) \to \mathrm{MIC}(X/Y/S).$$

*Proof.* Since the left square in the diagram is Cartesian, the map

$$\tau\colon F^*_{Y/S}(\mathcal{O}_{\mathbb{\Delta}_{X^{(p)}}(Y^{(p)})}) \to \mathcal{O}_{\mathbb{\Delta}_{X^\phi(Y)}} \cong \mathcal{O}_{\mathbb{PD}_X(Y)}$$

is an isomorphism. We claim that this isomorphism identifies the F-transform connection with the standard connection $d$. We must show that if $e'$ is a local section of $\mathcal{O}_{\mathbb{\Delta}_{X^{(p)}}(Y^{(p)})}$, then

$$(\mathrm{id} \otimes F^*_{Y/S})(d'(e')) = pd(\tau(e')).$$

The construction of $\mathbb{\Delta}_{X^{(p)}}(Y^{(p)})$ shows that $\mathcal{O}_{\mathbb{\Delta}_{X^{(p)}}(Y^{(p)})}$ is topologically generated by elements $e'$ such that $p^m e'$ belongs to $\mathcal{O}_{Y^{(p)}}$ for some $m > 0$, and so it suffices to check this equality for such $e'$. Since the target is $p$-torsion free, in fact it suffices to check the result for elements of $\mathcal{O}_{Y^{(p)}}$. Since the $p$-connection $d'$ on $\mathcal{O}_{\mathbb{\Delta}_{X^{(p)}}(Y^{(p)})}$ is compatible with the $p$-connection $pd$ on $\mathcal{O}_{Y^{(p)}}$ and the connection on $\mathcal{O}_{\mathbb{PD}_X(Y)}$ is compatible with the connection $d$ on $\mathcal{O}_Y$, the result is obvious for such elements.

To prove the last statement of the theorem, observe first that Shiho's theorem already tells us that the F-transform is an equivalence from the category $\mathrm{MICP}(Y^{(p)}/S)$ to the category $\mathrm{MIC}(Y/S)$. It remains only to check the compatibility condition. Suppose that $E'$ is an object of $\mathrm{MICP}(X^{(p)}/Y^{(p)}/S)$ and that $E$ is its $F$-transform. The $\mathcal{O}_{\mathbb{\Delta}_{X^{(p)}/S}(Y^{(p)})}$-module structure of $E'$ is a horizontal map:

$$\mathcal{O}_{\mathbb{\Delta}_{X^{(p)}}(Y^{(p)})} \otimes_{\mathcal{O}_{Y^{(p)}}} E' \to E'$$

which induces a horizontal map

$$F^*_{Y/S}(\mathcal{O}_{\mathbb{\Delta}_{X^{(p)}}(Y^{(p)})} \otimes_{\mathcal{O}_{Y^{(p)}}} E') \to F^*_{Y/S}(E').$$

Since the $F$-transform is compatible with tensor products, this gives us a horizontal map

$$\mathcal{O}_{\mathbb{PD}_X(Y)} \otimes_{\mathcal{O}_Y} F^*_{Y/S}(E') \to F^*_{Y/S}(E'),$$

endowing $F^*_{Y/S}(E')$ with the structure of an object of $\mathrm{MIC}(X/Y/S)$. This construction works in the opposite direction, proving the theorem. $\square$

**Remark 5.1.9.** As we saw in Theorem 4.3.1, the category $\mathrm{MICP}(X^{(p)}/Y^{(p)}/S)$ is, in a suitable sense, independent of the choice of $Y$, and the same is true for the category $\mathrm{MIC}(X/Y/S)$. One can check that the F-transform is also independent of this choice. Namely, if $g\colon Y \to Z$ is a morphism of $p$-completely smooth crystalline prisms such that $g \circ i$ is again a closed immersion, then one finds a 2-commutative diagram

$$\begin{array}{ccc} \mathrm{MICP}(X^{(p)}/Z'/S) & \longrightarrow & \mathrm{MIC}(X/Z/S) \\ \downarrow & & \downarrow \\ \mathrm{MICP}(X'/Y^{(p)}/S) & \longrightarrow & \mathrm{MIC}(X/Y/S). \end{array}$$

In section 6.4 we shall give a topos-theoretic version of the F-tranform that does not refer to any embedding.

## 5.2 The $\alpha$-transform

Let us revisit a construction from [34, §2.2]. Since everything takes place in characteristic $p$, we change the notation, letting $X \to S$ be a smooth morphism of schemes in characteristic $p > 0$.

We begin by recalling how to turn a $p$-linear homomorphism into an Artin-Schreier covering. Suppose for simplicity that $\Omega$ is a locally free coherent sheaf on $X$ and that

$$\beta \colon F_X^*(\Omega) \to \Omega$$

is an $\mathcal{O}_X$-linear map. Let $\mathcal{T}$ be the $\mathcal{O}_X$-linear dual of $\Omega$. We deduce $\mathcal{O}_X$-linear maps:

$$\begin{array}{rcl} \hat{\beta} & : & \mathcal{T} \longrightarrow F_X^*(\mathcal{T}) \\ z^* := S^\bullet(\hat{\beta}) & : & S^\bullet\mathcal{T} \longrightarrow F_X^*(S^\bullet\mathcal{T}) \\ F^*_{\mathbf{V}\mathcal{T}/X} & : & F_X^*(S^\bullet\mathcal{T}) \longrightarrow S^\bullet\mathcal{T} \\ h^* := F^*_{\mathbf{V}^*\Omega/X} \circ \hat{\beta}^* & : & S^\bullet\mathcal{T} \longrightarrow S^\bullet\mathcal{T} \end{array} \tag{5.2.1}$$

Here $\hat{\beta}$ is the dual of $\beta$ and induces the algebra homomorphism $z^*$, and $F^*_{\mathcal{T}/X}$ is the relative Frobenius homomorphism for $\mathbf{V}\mathcal{T}/X$, sending $1 \otimes D$ to $D^p$ for each $D \in \mathcal{T}$. Thus we have a commutative diagram of $X$-schemes:

$$\begin{array}{ccc} \mathbf{V}\mathcal{T}^{(p)} & \xrightarrow{\;z\;} & \mathbf{V}\mathcal{T} \\ {\scriptstyle F_{\mathbf{V}\mathcal{T}/X}} \uparrow & \nearrow {\scriptstyle h} & \\ \mathbf{V}\mathcal{T}, & & \end{array} \tag{5.2.2}$$

and we define:

$$\alpha \colon \mathbf{V}\mathcal{T} \to \mathbf{V}\mathcal{T} := \mathrm{id} - h. \tag{5.2.3}$$

Note that $h$ and $\alpha$ are linear over $X$ and compatible with the group structure of $\mathbf{V}\mathcal{T}$, although not with its vector bundle structure.

In particular, suppose that we are given a splitting of the inverse Cartier isomorphism, as in the following diagram.

$$
\begin{array}{ccc}
 & & F_{X/S*}\Omega^1_{X/S} \\
 & \overset{\zeta}{\nearrow} & \uparrow \text{inc} \\
\Omega^1_{X^{(p)}/S} & \xrightarrow{\zeta'} & F_{X/S*}\mathcal{Z}^1_{X/S} \\
 & \underset{C^{-1}_{X/S}}{\searrow} & \downarrow \text{proj} \\
 & & F_{X/S*}\mathcal{H}^1(\Omega^\bullet_{X/S}),
\end{array}
\tag{5.2.4}
$$

where $\mathcal{Z}^1_{X/S}$ is the sheaf of closed one-forms on $X/S$. (Recall, e.g. from Proposition 4.1.1, that a lifting of the Frobenius morphism to a smooth morphism of $p$-torsion free schemes provides such a splitting.) We deduce a map

$$\tilde{\zeta}\colon F_X^*(\Omega^1_{X/S}) \to \Omega^1_{X/S}, \tag{5.2.5}$$

the composite of the adjoint to $\zeta$ with the isomorphism $F_X^*(\Omega^1_{X/S}) \cong F_{X/S}^*(\Omega^1_{X^{(p)}/S})$.

The Artin-Schreier construction (5.2.2) then becomes:

$$
\begin{array}{ccc}
\mathbb{T}^{*(p)}_{X/S} & \xrightarrow{z} & \mathbb{T}^*_{X/S} \\
F_{\mathbb{T}^*/X} \uparrow & \overset{h}{\nearrow} & \\
\mathbb{T}^*_{X/S}, & &
\end{array}
\tag{5.2.6}
$$

and again:

$$\alpha := \mathrm{id} - h\colon \mathbb{T}^*_{X/S} \to \mathbb{T}^*_{X/S} \tag{5.2.7}$$

**Definition 5.2.1.** *If $(E, \theta)$ is a Higgs module on $X/S$ and $E_\theta$ is the corresponding sheaf on the cotangent bundle, then the $\alpha$-transform of $(E, \theta)$ is the Higgs module corresponding to $\alpha_*(E_\theta)$.*

Since $\alpha$ is $\mathcal{O}_X$-linear, the $\alpha$-transform does not change the underlying $\mathcal{O}_X$-module, so we could write $(E, \alpha(\theta))$ for the $\alpha$-transform of $(E, \theta)$. The following proposition gives an explicit local formula for $\alpha(\theta)$ which we shall need later.

**Proposition 5.2.2.** *Let $(E, \theta)$ be a Higgs module on $X/S$. Suppose that $(\omega_1, \ldots, \omega_n)$ is a basis for $\Omega^1_{X/S}$ with dual basis $(\tau_1, \ldots, \tau_n)$ for $\mathcal{T}_{X/S}$, and let $\theta_i := \theta_{\tau_i}$. Then*

$$
\begin{array}{rcl}
\theta\colon E & \to & E \otimes \Omega^1_{X/S} \\
e & \mapsto & \sum_i \theta_i(e) \otimes \omega_i,
\end{array}
$$

$$\begin{array}{rcl}\alpha(\theta)\colon E & \to & E \otimes \Omega^1_{X/S} \\ e & \mapsto & \sum_i \theta_i(e) \otimes \omega_i - \sum_i \theta_i^p(e) \otimes \zeta(1 \otimes \omega_i),\end{array}$$

*where $\theta_i^p$ is the pth iterate of $\theta_i$.*

*Proof.* The formula for $\theta$ is immediate from the definitions, and the formula for $\alpha(\theta)$ is straightforward. If $\tau \in \mathcal{T}_{X/S}$, let $\theta_\tau$ be the endomorphism of $E$ given by $\theta$ and let $\tilde{\theta}_\tau$ be the endomorphism given by the right side of the formula. Say $\tau = \sum a_i \tau_i$ with $a_i \in \mathcal{O}_X$ and $\hat{\zeta}(\tau) = \sum b_j \otimes \tau_j$, with $b_j \in \mathcal{O}_X$, where $\hat{\zeta}$ is the dual of $\zeta$. Then $F^*_{\mathbb{T}^*/X}(\hat{\zeta}(\tau)) = \sum b_j \tau_j^p$, and

$$\begin{array}{rcl}\tilde{\theta}_\tau(e) & = & \sum_i \theta_i(e)\langle \omega_i, \tau\rangle - \sum_i \langle \zeta(1 \otimes \omega_i), \tau\rangle \theta_i^p(e) \\ & = & \sum_i a_i \theta_i(e) - \sum_i \langle 1 \otimes \omega_i, \hat{\zeta}(\tau)\rangle \theta_i^p(e) \\ & = & \theta_\tau(e) - \sum_{i,j} \langle 1 \otimes \omega_i, b_j \otimes \tau_j\rangle \theta_i^p(e) \\ & = & \theta_\tau(e) - \sum_j b_j \theta_j^p(e) \\ & = & \theta_\tau(e) - \theta_{F_{\mathbb{T}^*/X}(\hat{\zeta}(\tau))}(e) \\ & = & \theta_{\alpha^*(\tau)}(e).\end{array}$$

□

Note that $\alpha$ induces an isomorphism when restricted to the formal completion of $\mathbb{T}^*_{Y/S}$, and so the $\alpha$-transform becomes an equivalence when restricted to the category of quasi-nilpotent Higgs modules. Furthermore, one has the following result.

**Proposition 5.2.3.** *Let $(E, \theta)$ be a quasi-nilpotent Higgs module on $X/S$ and let $(\tilde{E}, \tilde{\theta})$ be its $\alpha$-transform. Then there is a natural isomorphism in the derived category of $\mathcal{O}_X$-modules between their respective Higgs complexes.*

*Proof.* Let $E_\theta$ be the sheaf on $\mathbb{T}^*_{X/S}$ corresponding to $(E, \theta)$, so that $\alpha_*(E_\theta)$ is the sheaf corresponding to its $\alpha$-transform. Both of these are quasi-nilpotent, and so can be viewed as living on $\hat{\mathbb{T}}^*_{X/S}$. The map $\hat{\alpha}\colon \hat{\mathbb{T}}^*_{X/S} \to \hat{\mathbb{T}}^*_{X/S}$ is an isomorphism, and the inverse image of the zero section $i$ under $\alpha$ is again the zero section. Thus $\hat{\alpha}^*(i_*(\mathcal{O}_X)) \cong i_*(\mathcal{O}_X)$. Now we have a commutative diagram in the derived category of $\mathcal{O}_X$-modules:

$$
\begin{array}{ccc}
(E \otimes \Omega^{\bullet}_{X/S}, \alpha(\theta)) & & (E \otimes \Omega^{\bullet}_{X/S}, \theta) \\
\downarrow \sim & & \downarrow \sim \\
\mathrm{RHom}_{\mathcal{O}_{\mathbb{T}^*_{X/S}}}(i_*(\mathcal{O}_X), \alpha_*(E_\theta)) & & \mathrm{RHom}_{\mathcal{O}_{\mathbb{T}^*_{X/S}}}(i_*(\mathcal{O}_X), E_\theta) \\
\downarrow \sim & & \downarrow \sim \\
\mathrm{RHom}_{\mathcal{O}_{\hat{\mathbb{T}}^*_{X/S}}}(i_*(\mathcal{O}_X), \alpha_*(E_\theta)) & & \mathrm{RHom}_{\mathcal{O}_{\hat{\mathbb{T}}^*_{X/S}}}(i_*(\mathcal{O}_X), E_\theta) \\
\downarrow \sim & \nearrow \sim & \\
\mathrm{RHom}_{\mathcal{O}_{\hat{\mathbb{T}}^*_{X/S}}}(\alpha^* i_*(\mathcal{O}_X), E_\theta). & &
\end{array}
$$

Here the top two vertical isomorphisms come from Proposition 1.3.3, the next two are consequences of the flatness of $\hat{\mathbb{T}}^*_{X/S} \to \mathbb{T}^*_{X/S}$, and the third vertical isomorphism on the left is by adjunction. The slanted arrow comes from the isomorphism $\alpha^*(i_*(\mathcal{O}_X)) \cong i_*(\mathcal{O}_X)$. Note that this argument does not produce an explicit morphism between the two Higgs complexes. However, since all the arrows in the diagram above are natural so is the resulting isomorphism in the derived category. □

**Remark 5.2.4.** It is also possible, and sometimes convenient, to compute the $\alpha$-transform using a dual point of view. Recall from [7, A10] that $\Gamma^{\bullet}\Omega$ is the Matlis dual of $S^{\bullet}\mathcal{T}$. Thus, the dual of the relative Frobenius morphism

$$F_{\mathbf{V}\mathcal{T}/X} \colon \mathbf{V}\mathcal{T} \to \mathbf{V}\mathcal{T}^{(p)}$$

is the Verschiebung morphism

$$V_{\Omega/X} \colon \mathbf{V}^{\gamma}\Omega^{(p)} \to \mathbf{V}^{\gamma}\Omega,$$

corresponding to the algebra homomorphism $\Gamma^{\bullet}\Omega \to F_X^*(\Gamma^{\bullet}\Omega)$ which is zero on $\Gamma^n\Omega$ if $p \not| n$ and which takes $\omega^{[mp]}$ to $\omega^{[m]}$. Then given $\beta \colon F_X^*\Omega \to \Omega$ as above, we can form:

$$
\begin{array}{rcl}
\hat{z}^* := F_X^*(\Gamma^{\bullet}\beta) & : & F_X^*(\Gamma^{\bullet}\Omega) \longrightarrow \Gamma^{\bullet}\Omega \\
V^*_{\Omega/X} & : & \Gamma^{\bullet}\Omega \longrightarrow F_X^*(\Gamma^{\bullet}\Omega) \\
\hat{h}^* := \hat{z}^* \circ V^*_{\Omega/X} & : & \Gamma^{\bullet}\Omega \longrightarrow \Gamma^{\bullet}\Omega.
\end{array}
\tag{5.2.8}
$$

This leads to the following diagram:

$$
\begin{array}{ccc}
\mathbf{V}^\gamma\Omega^{(p)} & \xleftarrow{\hat{z}} & \mathbf{V}^\gamma\Omega \\
\Big\downarrow V_{\Omega/X} & \swarrow \hat{h} & \\
\mathbf{V}^\gamma\Omega & &
\end{array}
$$

and we can set

$$\hat{\alpha} := \mathrm{id} - \hat{h} \colon \mathbf{V}^\gamma\Omega \to \mathbf{V}^\gamma\Omega.$$

One could also work with the completed PD-algebra $\hat{\Gamma}^{\bullet}\Omega$, but since we will be concentrating on quasi-nilpotent fields, we shall not need this.

As we saw in Proposition 1.3.2, a quasi-nilpotent $\Omega$-field $\theta$ on $E$ can be viewed as an action $r_\theta \colon \mathbf{V}E \times \mathbf{V}^\gamma\Omega \to \mathbf{V}E$ of $\mathbf{V}^\gamma\Omega$ on $\mathbf{V}E$. and it follows easily that the $\alpha$ transform of $\theta$ corresponds to the action given by composition with $\hat{\alpha}$. That is:

$$r_{\alpha(\theta)} = r_\theta \circ (\mathrm{id} \times \hat{\alpha}).$$

Proposition 5.4.1 below will give yet another construction of the $\alpha$-transform.

## 5.3 F-transform and $p$-curvature

Our goal here is to give an explicit formula for the $p$-curvature of the reduction modulo $p$ of the F-transform of a module with integrable $p$-connection on a crystalline prism. Since everything takes place in characteristic $p$, we change the notation, letting $X \to S$ be a smooth morphism of schemes in characteristic $p > 0$.

We can generalize the F-transform construction slightly as follows. Suppose that we are given a splitting $\zeta$ of the inverse Cartier isomorphism, as in the diagram (5.2.4). (Recall that a lifting of $F_{X/S} \colon X/S \to X^{(p)}/S$ to a morphism of smooth schemes $\tilde{X}/\tilde{S} \to \tilde{X}'/\tilde{S}$, where $\tilde{S}$ is flat over $\mathbf{Z}/p^2\mathbf{Z}$, provides such a splitting.) If $\theta \colon E' \to E' \otimes \Omega^1_{X^{(p)}/S}$ is a Higgs field on a sheaf of $\mathcal{O}_{X^{(p)}}$-modules $E'$, there is a unique connection $\nabla$ on $E := F^*_{X/S}(E')$ such that $\nabla(1 \otimes e') = (\zeta \otimes \mathrm{id})\theta(e')$ for every local section $e'$ of $E'$. We call this $(E, \nabla)$ the *$\zeta$-transform* of $(E', \theta)$. This construction generalizes the F-transform (5.1.1) in the obvious way. Our goal is to compute the $p$-curvature of $(E, \nabla)$. The result is compatible with the formulas in [34, 2.5, 2.11, and 2.13] and generalizes the formula of Katz [23, 7.1.2], which is the case in which $E = \mathcal{O}_X$.

Since $F^*_{X/S}(\Omega^1_{X^{(p)}/S}) \cong F^*_X(\Omega^1_{X/S})$, we can and shall view the $p$-curvature of $(E, \nabla)$ as an $F^*_X(\Omega^1_{X/S})$-field on $E$, *i.e.*, as a sheaf $E_\psi$ on $\mathbb{T}^{*(p)}_{X/S}$. The relative

Frobenius diagram for $\mathbb{T}^*_{X/S}$ is

$$
\begin{array}{ccccc}
\mathbb{T}^*_{X/S} & \xrightarrow{F_{\mathbb{T}^*/S}} & \mathbb{T}^{*(p)}_{X/S} & \xrightarrow{\pi_{\mathbb{T}^*/X}} & \mathbb{T}^*_{X/S} \\
 & \searrow & \downarrow & & \downarrow \\
 & & X & \xrightarrow{F_{X/S}} & X.
\end{array}
\tag{5.3.1}
$$

The square in this diagram factors into a diagram with two Cartesian squares:

$$
\begin{array}{ccccc}
\mathbb{T}^{*(p)}_{X/S} & \xrightarrow{\pi'} & \mathbb{T}^*_{X^{(p)}/S} & \xrightarrow{\pi} & \mathbb{T}^*_{X/S} \\
\downarrow & & \downarrow & & \downarrow \\
X & \xrightarrow{F_{X/S}} & X^{(p)} & \xrightarrow{\pi_{X/S}} & X,
\end{array}
\tag{5.3.2}
$$

where $\pi$ is defined by the map $\mathcal{T}_{X/S} \to \mathcal{T}_{X^{(p)}/S} = \pi^*_{X/S}(\mathcal{T}_{X/S})$ sending $D$ to $1 \otimes D$. The endomorphisms $z$ and $\alpha$ of $\mathbb{T}^*_{X/S}$ constructed in (5.2.6) pull back to endomorphisms $h^{(p)}$ and $\alpha^{(p)}$ of $\mathbb{T}^*_{X^{(p)}/S}$, where $h^{(p)} \circ \pi = \pi \circ h$ and $\alpha^{(p)} \circ \pi = \pi \circ \alpha$.

**Theorem 5.3.1.** *Let $X/S$ be a smooth morphism of schemes in characteristic $p$ endowed with a splitting $\zeta$ of the inverse Cartier isomorphism as in diagram (5.2.4), let $\alpha$ be the resulting endomorphism of the cotangent space $\mathbb{T}^*_{X/S}$ as in in Definition 5.2.1, and let $\alpha^{(p)}$ be its pullback to $\mathbb{T}^*_{X^{(p)}/S}$.*

*(1) If $(E', \theta)$ is Higgs module on $X^{(p)}/S$, there is a unique connection $\nabla$ on $E := F^*_{X/S}(E')$ such that the following diagram*

$$
\begin{array}{ccc}
E' & \xrightarrow{\theta} & E' \otimes \Omega^1_{X^{(p)}/S} \\
\downarrow \eta & & \downarrow \eta \otimes \zeta \\
F^*_{Y/S}(E') & \xrightarrow{\nabla} & F^*_{Y/S}(E' \otimes \Omega^1_{X/S})
\end{array}
$$

*commutes. The resulting module with connection $(E, \nabla)$ is called the $\zeta$-transform of $(E', \theta)$.*

*(2) The $p$-curvature $\psi$ of $(E, \nabla)$ is the negative of the pullback of the $\alpha^{(p)}$-transform of $(E', \theta)$ along $F_{X/S}$. This can be made explicit in the following equivalent ways.*

*(a) For every local section $D$ of $\mathcal{T}_{X/S}$, we have, if $D' := \pi^*_{X/S}(D)$,*

$$\psi_{D'} = -F^*_{X/S}(\theta^*_{\alpha^{(p)}(D')}).$$

(b) *If $E'_\theta$ is the sheaf on $\mathbb{T}^*_{X^{(p)}/S}$ corresponding to $(E', \theta)$ and $E_\psi$ is the sheaf on $\mathbb{T}^{*(p)}_{X/S}$ corresponding to $(E, \psi)$, then*

$$E_\psi = \pi'^*(E'_{-\alpha^{(p)}(\theta)}),$$

*where $\pi'$ is the morphism in diagram (5.3.2).*

(c) *If $\theta' := \alpha^{(p)}(\theta)$, then there is a commutative diagram:*

$$\begin{array}{ccc} F^*_{X/S}(E') & \xrightarrow{-\theta'} & F^*_{X/S}(E' \otimes \Omega^1_{X^{(p)}/S}) \\ \downarrow \cong & & \downarrow \cong \\ E & \xrightarrow{\psi} & E \otimes F^*_{X/S}(\Omega^1_{X^{(p)}/S}). \end{array}$$

*Proof.* The first statement has already been explained, but, for the proof of remaining statements, we shall need to make the formula for the $\zeta$-transform $\nabla$ more explicit. The Higgs field $\theta$ gives us an algebra homomorphism:

$$\theta^\bullet \colon S^\bullet \mathcal{T}_{X^{(p)}/S} \longrightarrow End_{\mathcal{O}_{X^{(p)}}}(E').$$

Since $E = F^*_{X/S}(E')$, we can define a homomorphism $\Theta^\bullet$ by the following diagram.

$$\begin{array}{ccccc} F^*_X(S^\bullet \mathcal{T}_{X/S}) & \xrightarrow{\cong} & F^*_{X/S}(S^\bullet \mathcal{T}_{X^{(p)}/S}) & \xrightarrow{F^*_{X/S}(\theta^\bullet)} & End_{\mathcal{O}_X}(E) \\ & \nwarrow^{z^*} & & \nearrow_{\Theta^\bullet} & \\ & & S^\bullet \mathcal{T}_{X/S}, & & \end{array}$$

where $z^*$ is the map in diagram (5.2.6). This defines a Higgs field $\Theta$ on the $\mathcal{O}_X$-module $E$, and it follows from the definitions that the $\zeta$-transform connection $\nabla$ corresponds to the map

$$\nabla \colon \mathcal{T}_{X/S} \to \mathrm{End}_{\mathcal{O}_{X^{(p)}}}(E) \colon D \mapsto D \otimes \mathrm{id} + \Theta_D,$$

where $D \otimes \mathrm{id}$ is the endomorphism of $E$ corresponding to the Frobenius pull-back connection on $F^*_{X/S}(E')$.

**Lemma 5.3.2.** *If $D \in \mathcal{T}_{X/S}$, then*

$$\begin{aligned} \Theta^p_D &= F^*_{X/S}(\theta^\bullet(\pi^*(h^*(D)))) \\ &= F^*_{X/S}(\theta^\bullet(h^{(p)*}(\pi^*(D)))). \end{aligned}$$

*Proof.* Consider the following diagram:

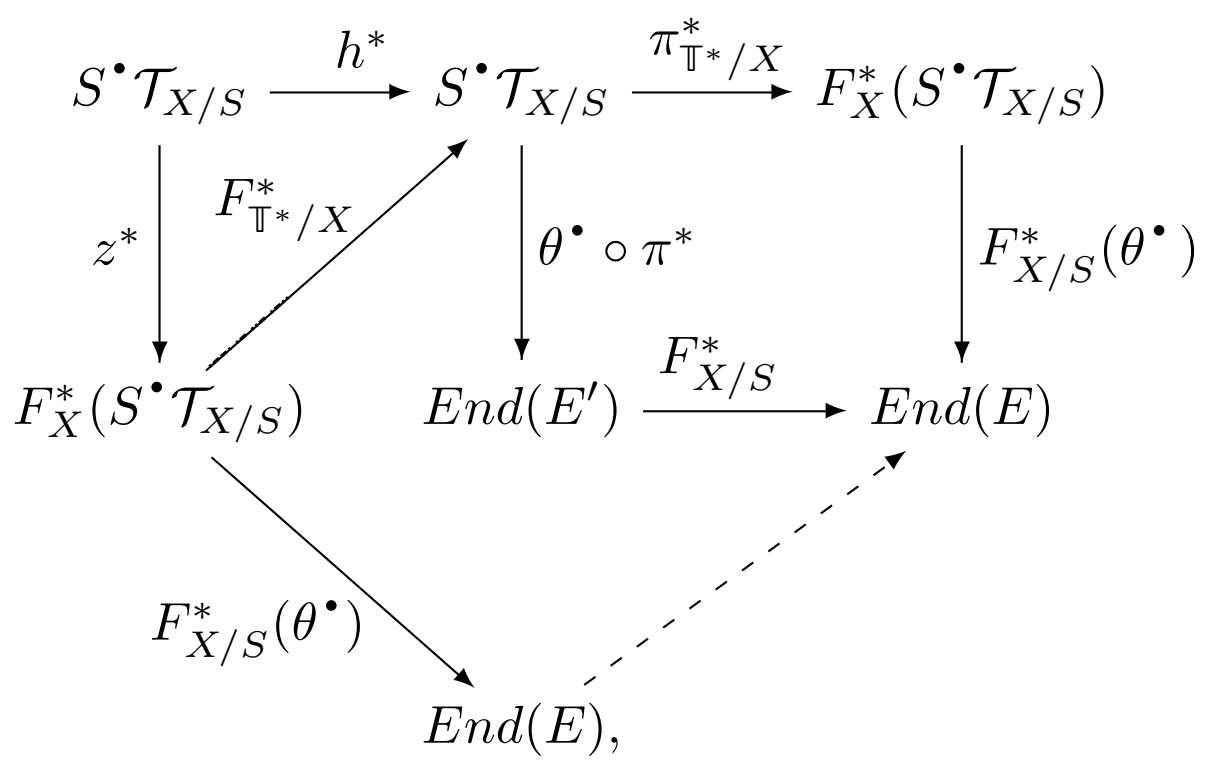

where the dotter arrow takes an element to its $p$th power. The triangle commutes by the definition of $h$, and the square commutes by construction and the $F_X$-linearity of $\pi_{\mathbb{T}^*/X}$. Furthermore, $F_{\mathbb{T}^*/X} \circ \pi_{\mathbb{T}^*/X}$ is the absolute Frobenius endomorphism of $\mathbb{T}^{*(p)}_{X/S}$, taking an element $D'$ to its $p$th power. Since $F^*_{X/S}(\theta^\bullet)$ is an algebra homomorphism, the entire diagram commutes. Since $\Theta = F^*_{X/S}(\theta^*) \circ z^*$, this proves the first equation of the lemma. The second follows from the fact that $\pi^* \circ h^* = h^{(p)*} \circ \pi^*$. □

The next lemma is where we use the fact that $\zeta$ splits the inverse Cartier isomorphism.

**Lemma 5.3.3.** *If $D$ is any section $\mathcal{T}_{X/S}$ and $D' := \pi^*(D)$, then as endomorphisms of $F^*_{X/S}(E')$, we have*

$$\mathrm{ad}^{p-1}_{D\otimes \mathrm{id}}(\Theta_D) - \Theta_{D^p} = -F^*_{X/S}(\theta_{D'}).$$

*Furthermore, $\Theta_D$ commutes with $\mathrm{ad}^j_{D\otimes\mathrm{id}}(\Theta_D)$ for every $j$.*

*Proof.* We work locally, with the aid of a system of coordinates $(t_1, \ldots, t_n)$. Let $(D_1, \ldots, D_n)$ be the basis for $\mathcal{T}_{X/S}$ dual to $(dt_1, \ldots, dt_n)$, let $\omega_i := \zeta(dt'_i)$, let $D'_i = \pi^*(D_i) \in \mathcal{T}_{X^{(p)}/S}$, and let $\theta_i := \theta_{D'_i}$. The integrability condition implies that $\theta_1, \ldots, \theta_n$ is a family of commuting endomorphisms of $E'$. Then for each $e' \in E'$, we have

$$\begin{aligned} \theta(e') &= \sum_i dt'_i \otimes \theta_i(e') \\ \nabla(1 \otimes e') &= \sum_i \omega_i \otimes \theta_i(e'). \end{aligned}$$

For $D \in \mathcal{T}_{X/S}$, we have:

$$\Theta_D = F^*_{X/S}(\theta)_{\hat{\zeta}(D)} \in \mathrm{End}_{\mathcal{O}_X}(E)$$

$$
\begin{aligned}
&= \sum_i \langle 1\otimes dt_i', \hat\zeta(D)\rangle \otimes \theta_i \\
&= \sum_i \langle \zeta(1\otimes dt_i'), D\rangle \otimes \theta_i \\
&= \sum_i \langle \omega_i, D\rangle \otimes \theta_i.
\end{aligned}
$$

This formula implies that $\Theta_D$ commutes with each $\theta_i$. Furthermore, if $j\in\mathbf{N}$, then

$$
\mathrm{ad}^j_{D\otimes \mathrm{id}}(\Theta_D) = \sum_i D^j(\langle \omega_i, D\rangle)\otimes\theta_i,
$$

so $\Theta_D$ also commutes with $\mathrm{ad}^j_{D\otimes\mathrm{id}}(\Theta_D)$.

Recall from [23, 7.1.2.6] that the Cartier operator $C\colon F_{Y/S*}(\mathcal{Z}^1_{X/S}) \to \Omega^1_{X^{(p)}/S}$ satisfies the following formula

$$
F^*_{X/S}(\langle C(\omega), \pi^*_{X/S}D\rangle) = \langle \omega, D^p\rangle - D^{p-1}(\langle\omega, D\rangle),
$$

for $\omega\in F_{X/S*}\mathcal{Z}^1_{X/S}$ and $D\in\mathcal{T}_{X/S}$. Using the fact that $\Theta_D$ commutes with each $\theta_i$, we compute:

$$
\begin{aligned}
\mathrm{ad}^{p-1}_{D\otimes\mathrm{id}}(\Theta_D) - \Theta_{D^p} &= \mathrm{ad}^{p-1}_{D\otimes\mathrm{id}}\left(\sum_i \langle\omega_i, D\rangle\otimes\theta_i\right) - \sum_i \langle\omega_i, D^p\rangle\otimes\theta_i \\
&= \sum_i \left(D^{p-1}(\langle\omega_i, D\rangle)\otimes\theta_i - \langle\omega_i, D^p\rangle\otimes\theta_i\right) \\
&= -\sum_i F^*_{X/S}(\langle C(\omega_i), D'\rangle)\otimes\theta_i.
\end{aligned}
$$

Since $\omega_i = \zeta(dt_i')$ and $\zeta$ is a splitting of the inverse Cartier isomorphism, in fact $C(\omega_i) = dt_i'$, and we find that

$$
\begin{aligned}
\mathrm{ad}^{p-1}_{D\otimes\mathrm{id}}(\Theta_D) - \Theta_{D^p} &= -\sum F^*_{X/S}\langle dt_i', D'\rangle\otimes\theta_i \\
&= -F^*_{X/S}(\theta_{D'})
\end{aligned}
$$

as claimed. □

To prove the the theorem, recall that the $p$-curvature is the $F_{X/S}$-linear map:

$$
\mathcal{T}_{X^{(p)}/S} \to End_{\mathcal{O}_X}(E) : D' \mapsto \nabla^p_D - \nabla_{D^p},
$$

where $D' := \pi^*(D)$. Since $\Theta_D$ commutes with each $\mathrm{ad}^j_{D\otimes\mathrm{id}}(\Theta_D)$, Jacobson's formula for the $p$th power of a sum simplifies:

$$
\begin{aligned}
\psi_{D'} &:= \nabla^p_D - \nabla_{D^p} \\
&= (D\otimes\mathrm{id} + \Theta_D)^p - (D^p\otimes\mathrm{id} + \Theta_{D^p}) \\
&= (D\otimes\mathrm{id})^p + \mathrm{ad}^{p-1}_{D\otimes\mathrm{id}}(\Theta_D) + \Theta^p_D - (D^p\otimes\mathrm{id} + \Theta_{D^p}) \\
&= \Theta^p_D + \mathrm{ad}^{p-1}_{D\otimes\mathrm{id}}(\Theta_D) - \Theta_{D^p}
\end{aligned}
$$

Using Lemma 5.3.2 to compute $\theta^p_D$ and Lemma 5.3.3 for the remaining terms, we get:

$$
\begin{aligned}
\psi_{D'} &= F^*_{X/S}(\theta^\bullet(h^{(p)*}(D')) - F^*_{X/S}(\theta_{D'}) \\
&= F^*_{X/S}\left(\theta^\bullet(h^{(p)*}(D')) - \theta_{D'}\right) \\
&= -F^*_{X/S}\left(\theta^\bullet(\alpha^{(p)}(D'))\right).
\end{aligned}
$$

This completes the proof. □

**Remark 5.3.4.** Theorem 2.11 of [34] also defines a functor, called the "local inverse Cartier transform," from Higgs modules to modules with connection. As explained in Theorem 2.13 and Corollary 2.14 there, this functor is obtained by changing the sign of the Higgs field, then applying $\alpha^*$ (the inverse of the $\alpha$-transform), and then what we call here the F-transform.

## 5.4 The $\alpha$-transform, $d'$, and $\nabla'$

We again suppose that $Y/S$ is a $p$-completely smooth of crystalline prisms and that $X \to \overline{Y}$ is a (regular) closed immersion, with $X/\overline{S}$ smooth. The morphism $\mathbb{A}_X(Y) \to Y$ is affine, and we write $\mathcal{O}_{\mathbb{A}}$ for the structure sheaf of $\mathbb{A}_X(Y)$, now viewed as a sheaf of $\mathcal{O}_Y$-modules on $Y$. We saw in Proposition 4.2.1 that this sheaf is naturally endowed with a $p$-connection $d'$ whose reduction modulo $p$

$$\overline{d}' : \mathcal{O}_{\overline{\mathbb{A}}} \to \mathcal{O}_{\overline{\mathbb{A}}} \otimes \Omega^1_{Y/S} \tag{5.4.1}$$

is a Higgs field, and as we shall see in Theorem 6.3.5, the corresponding Higgs complex calculates prismatic cohomology mod $p$. On the other hand, we found another Higgs field $\nabla'$ (3.3.3), which, unlike $\overline{d}'$, is compatible, in a suitable sense, with the Sen operator (Definition 3.2.3). The main result of this section, Theorem 5.4.3, shows that $\nabla'$ is the $\alpha$-transform of $\overline{d}'$. The proof of this key result will depend on another construction of the $\alpha$-transform, which we explain next.

Recall that, since $Y$ is a crystalline prism, the diagonal embedding $Y \to Y(1)$ lifts to morphism $Y \to \mathbb{A}_Y(1) := \mathbb{A}_{\overline{Y}}(Y(1))$, and that the ideal $I_{Y/\mathbb{A}}$ is invariant under the $\delta$-structure. We denote by $\mathrm{Fil}_\delta$ the corresponding filtration of $\mathcal{O}_{\mathbb{A}_Y(1)}$, as in Definition 1.7.14. There are natural maps:

$$
\begin{aligned}
\pi_\omega \colon I_{\overline{Y}/\overline{Y}(1)} \to I_{\overline{Y}/\overline{Y}(1)}/I^2_{\overline{Y}/\overline{Y}(1)} &\xrightarrow{\cong} \Omega^1_{\overline{Y}/\overline{S}} \\
\pi_\gamma \colon I_{\overline{Y}/\overline{\mathbb{A}}_Y(1)} \to I_{\overline{Y}/\overline{\mathbb{A}}_Y(1)}/I^{[2]}_{\overline{Y}/\overline{\mathbb{A}}_Y(1)} &\xrightarrow{\cong} \Omega^1_{\overline{Y}/\overline{S}} \\
\pi_\delta \colon I_{\overline{Y}/\overline{\mathbb{A}}_Y(1)} \to \mathrm{Gr}^1_{\mathrm{Fil}_\delta} \mathcal{O}_{\overline{\mathbb{A}}_Y(1)} &\xrightarrow{\cong} \Omega^1_{\overline{Y}/\overline{S}},
\end{aligned}
$$

where $\pi_\omega$ is the usual map, where $\pi_\gamma$ is obtained from $\pi_\omega$ and the isomorphism

$$I_{\overline{Y}/\overline{Y}(1)}/I^2_{\overline{Y}/\overline{Y}(1)} \cong I_{\overline{Y}/\overline{\mathbb{D}}_Y(1)}/I^2_{\overline{Y}/\overline{\mathbb{D}}_Y(1)} \cong I_{\overline{Y}/\overline{\mathbb{A}}_Y(1)}/I^{[2]}_{\overline{Y}/\overline{\mathbb{A}}_Y(1)},$$

and the second map in the definition of $\pi_\delta$ is the inverse of the isomorphism $\rho_\delta$ defined in Proposition 2.3.17.

**Proposition 5.4.1.** *Let $Y/S$ be a $p$-completely smooth morphism of crystalline prisms. Then the $\alpha$-transform $(E, \alpha(\theta))$ of a quasi-nilpotent Higgs module $(E, \theta)$ on $\overline{Y}/S$ can be computed as follows. Let*

$$r_\theta^* \colon E \longrightarrow E \otimes \Gamma^\bullet(\Omega^1_{\overline{Y}/\overline{S}}) \cong E \otimes \mathcal{O}_{\overline{\mathbb{\Delta}}_{Y(1)}}$$

*be the coaction (1.3.2) associated to $\theta$ , and let*

$$p_1^* \colon E \longrightarrow E \otimes \Gamma^\bullet(\Omega^1_{\overline{Y}/\overline{S}}) \cong E \otimes \mathcal{O}_{\overline{\mathbb{\Delta}}_{Y(1)}}$$

*be the map $e \mapsto e \otimes 1$. Then $r_\theta^* - p_1^*$ factors through $E \otimes I_{\overline{Y}/\overline{\mathbb{\Delta}}_Y(1)}$, and*

$$\begin{array}{rcl} \theta & = & (\mathrm{id} \otimes \pi_\gamma) \circ (r_\theta - p_1^*) \\ \alpha(\theta) & = & (\mathrm{id} \otimes \pi_\delta) \circ (r_\theta - p_1^*). \end{array}$$

*Proof.* Statement (1) for $\pi_\gamma$ has already been explained in Proposition 1.3.2. The statement for $\pi_\delta$ will depend on the comparison between the PD-filtration and the $\delta$-filtration of the ideal $I_{Y/\mathbb{\Delta}_Y(1)}$ described in Proposition 2.3.17.

**Lemma 5.4.2.** *Let $Y/S$ be a $p$-completely smooth morphism of crystalline prisms and let $Z = Y(1) := Y \times_S Y$, with its inherited $\delta$-structure $\delta_Z$.*

*(1) For $f \in \mathcal{O}_Y$, let $\xi_f := 1 \otimes f - f \otimes 1 \in I_{Y/Z}$ and let $\eta_f := \rho_{\mathbb{\Delta}}(\xi_f) \in I_{Y/\mathbb{\Delta}_Y(1)}$. Then:*

*(a) $\pi_\omega(\overline{\delta}_Z(\eta_f)) = \tilde{\zeta}(1 \otimes d\overline{f})$.*

*(b) $\pi_\delta(\overline{\eta}_f) = df$.*

*(c) $\pi_\delta(\overline{\eta}_f^{[p]}) = -\tilde{\zeta}(1 \otimes d\overline{f})$.*

*(2) The standard identification of the conormal sheaf $\mathcal{C}_{\overline{Y}/\overline{Z}}$ with $\Omega^1_{\overline{Y}/\overline{S}}$ fits into a commutative diagram*

$$\begin{array}{ccc} F^*_{\overline{Y}}(\mathcal{C}_{\overline{Y}/\overline{Z}}) & \xrightarrow{\delta_{\overline{Y}/\overline{Z}}} & C_{\overline{Y}/\overline{Z}} \\ \cong \downarrow & & \downarrow \cong \\ F^*_{\overline{Y}}(\Omega^1_{\overline{Y}/\overline{S}}) & \xrightarrow{\tilde{\zeta}} & \Omega^1_{\overline{Y}/\overline{S}} \end{array}$$

*where $\delta_{\overline{Y}/\overline{Z}}$ is the map induced by $\delta_Z$ (see (4) of Proposition 1.7.13) and $\tilde{\zeta}$ is the map (5.2.5).*

*Proof.* We must calculate $\delta_Z$. Let $p_1$ and $p_2$ be the projections $Y(1) \to Y$, so that $\xi_f := p_2^*(f) - p_1^*(f)$. Then:

$$\phi_Z(\xi_f) \quad = \quad \phi_Z(p_2^*(f)) - \phi_Z(p_1^*(f))$$

$$
\begin{aligned}
&= p_2^*(\phi_Z(f)) - p_1^*(\phi_Z(f)) \\
&= p_2^*(f^p + p\delta_Y(f)) - p_1^*(f^p + p\delta_Y(f)) \\
&= p_2^*(f^p) + pp_2^*(\delta_Y(f)) - p_1^*(f^p) - pp_1^*(\delta_Y(f)) \\
&= (p_1^*(f) + \xi_f)^p + pp_2^*(\delta_Y(f)) - (p_1^*(f))^p - pp_1^*(\delta_Y(f)) \\
&= (p_1^*(f))^p + pC(p_1^*(f), \xi_f) + \xi_f^p - (p_1^*(f))^p + p\xi_{\delta_Y(f)} \\
&= \xi_f^p + pC(p_1^*(f), \xi_f) + p\xi_{\delta_Y(f)}.
\end{aligned}
$$

Thus we can conclude that

$$
\begin{aligned}
\delta_Z(\xi_f) &= C(p_1^*(f), \xi_f) + \xi_{\delta_Y(f)} \\
&\equiv p_1^*(f^{p-1})\xi_f + \xi_{\delta_Y(f)} \pmod{I_{Y/Z}^2}.
\end{aligned}
$$

The image of $\xi_f$ under the identification $\mathcal{C}_{Y/Z} \cong \Omega^1_{Y/S}$ takes $\xi_f$ to $df$, so we see that

$$
\delta(\xi_f) \mapsto f^{p-1}df + d\delta_Y(f) \in I_{Y/Z}/I_{Y/Z}^2 \cong \Omega^1_{Y/S}. \tag{5.4.2}
$$

On the other hand,

$$
\begin{aligned}
\phi_Y(f) &= f^p + p\delta_Y(f) \\
\phi_Y^*(df) &= pf^{p-1}df + pd\delta_Y(f) \\
p^{-1}\phi_Y^*(df) &= f^{p-1}df + d\delta_Y(f).
\end{aligned}
$$

Reducing these equations modulo $p$ proves the first equation in statement (1).

For the second equation, note that, by definition, $\pi_\delta = \rho_\delta^{-1} \circ pr_\delta$, and from the diagram in statement (2) of Proposition 2.3.17 that $\pi_\omega = \rho_\delta^{-1} \circ pr_\delta \circ \rho_{\mathbb{A}}$. Then

$$
\begin{aligned}
\pi_\delta(\eta_f) &= \rho_\delta^{-1}(pr_\delta(\eta_f)) \\
&= \rho_\delta^{-1}(pr_\delta(\rho_{\mathbb{A}}(\xi_f))) \\
&= \pi_\omega(\xi_f) \\
&= df.
\end{aligned}
$$

Using the same method, statement (1a), and (3) of Proposition 2.3.17, we get:

$$
\begin{aligned}
\pi_\delta(\eta_f^{[p]}) &= \rho_\delta^{-1}(pr_\delta(\eta_f^{[p]})) \\
&= \rho_\delta^{-1}(pr_\delta((\rho_{\mathbb{A}}(\xi_f))^{[p]})) \\
&= \rho_\delta^{-1}(\rho_\delta(\pi_\omega(-\overline{\delta}_Z(\xi_f)))) \\
&= -(\pi_\omega(\overline{\delta}_Z(\xi_f))) \\
&= -\zeta(1 \otimes df).
\end{aligned}
$$

This proves statement (1). Statement (2) is just a reformulation of statement (1a). □

To prove the proposition, consider the following diagram. The commutativity of its upper part is clear, the commutativity of the lower square follows from that of diagram (2.3.6), and the commutativity of the lower triangle is the definition of $\tilde{\theta}$, which we claim is the $\alpha$-transform of $\theta$.

$$
\begin{array}{ccccc}
 & & E \otimes \Omega^1_{\overline{Y}/\overline{S}} & \longleftarrow & E \otimes I_{\overline{Y}/\overline{\mathbb{\Delta}}_Y(1)}/I^{[2]}_{\overline{Y}/\overline{\mathbb{\Delta}}_Y(1)} \\
 & \overset{\theta}{\nearrow} & \Big\uparrow {\scriptstyle \mathrm{id}\otimes\pi_\gamma} & & \Big\uparrow \\
E & \xrightarrow{r_\theta - p_1^*} & E \otimes I_{\overline{Y}/\overline{\mathbb{\Delta}}_Y(1)} & \longrightarrow & E \otimes I_{\overline{Y}/\overline{\mathbb{\Delta}}}/(I^2_{\overline{Y}/\overline{\mathbb{\Delta}}_Y(1)} + I^{[p+1]}_{\overline{Y}/\overline{\mathbb{\Delta}}_Y(1)}) \\
 & \underset{\tilde{\theta}}{\searrow} & \Big\downarrow {\scriptstyle \mathrm{id}\otimes\pi_\delta} & & \Big\downarrow \\
 & & E \otimes \Omega^1_{\overline{Y}/\overline{S}} & \xleftarrow{\mathrm{id}\otimes\rho_\delta^{-1}} & E \otimes \mathrm{Gr}^1_{\mathrm{Fil}_\delta}\, \mathcal{O}_{\overline{\mathbb{\Delta}}_Y(1)}.
\end{array}
$$

We can check this locally, with the aid of a local coordinate system $(x_1, \ldots, x_n)$ for $Y/S$. We get a corresponding set of generators $(\xi_1, \ldots, \xi_n)$ for the ideal $I_{Y/Y(1)}$ and $(\eta_1, \ldots, \eta_n)$ for the ideal $I_{Y/\mathbb{D}_Y(1)}$, where $\eta_i := \rho_{\mathbb{D}}(\xi_i)$. Let $\overline{\xi}_i$ (resp. $\overline{\eta}_i$) be the image of $\xi_i$ in $I_{\overline{Y}/\overline{Y}(1)}$ (resp., of $\eta_i$ in $I_{\overline{Y}/\overline{\mathbb{\Delta}}}$). These elements $\overline{\eta}_i$, together with their $p$th divided powers, induce a basis for $I_{\overline{\mathbb{\Delta}}}/(I^2_{\overline{\mathbb{\Delta}}} + I^{[p+1]}_{\overline{\mathbb{\Delta}}})$. Statement (1) of Lemma 5.4.2 implies that

$$
\begin{aligned}
\pi_\delta(\overline{\eta}_i) &= dx_i \\
\pi_\delta(\overline{\eta}_i^{[p]}) &= -\zeta(dx_i).
\end{aligned}
$$

Let $D'_i := pd/dx_i$, and for a multi-index $I$, let $\theta^I_{D'} := \theta^{I_1}_{D'_1} \cdots \theta^{I_n}_{D'_n}$. Then for $e \in E$,

$$
\begin{aligned}
r^*_\theta(e) - e \otimes 1 &= \sum_{I \neq 0} \theta^I_{D'}(e) \otimes \overline{\eta}^{[I]} \\
&\equiv \sum_i \theta_{D'_i}(e) \otimes \overline{\eta}_i + \sum_i \theta^p_{D'_i}(e) \otimes \overline{\eta}_i^{[p]} \pmod{E \otimes (I^2_{\overline{\mathbb{\Delta}}} + I^{[p+1]}_{\overline{\mathbb{\Delta}}})} \\
\pi_\delta(r^*_\theta(e) - e \otimes 1) &= \sum_i \theta_{D'_i}(e) \otimes \pi_\delta(\overline{\eta}_i) + \sum_i \theta^p_{D'_i}(e) \otimes \pi_\delta(\overline{\eta}_i^{[p]}) \\
&= \sum_i \theta_{D'_i}(e) \otimes dx_i - \sum_i \theta^p_{D'_i}(e) \otimes \zeta(dx_i).
\end{aligned}
$$

According to Proposition 5.2.2, this is indeed the $\alpha$ transform of $\theta$. □

Now we can prove the main result of this section.

**Theorem 5.4.3.** *Let $Y/S$ be a $p$-completely smooth morphism of crystalline prisms and let $X \to \overline{Y}$ be a regular closed immersion. Then the $\alpha$-transform of the Higgs field (5.4.1)*

$$
\overline{d}' : \mathcal{O}_{\overline{\mathbb{\Delta}}_X(Y)} \to \mathcal{O}_{\overline{\mathbb{\Delta}}_X(Y)} \otimes \Omega^1_{\overline{Y}/\overline{S}}
$$

*is the Higgs field (3.3.3):*

$$\nabla' \colon \mathcal{O}_{\overline{\mathbb{A}}_X(Y)} \to \mathcal{O}_{\overline{\mathbb{A}}_X(Y)} \otimes \Omega^1_{\overline{Y}/\overline{S}}.$$

*Consequently there is a canonical isomorphism*

$$\mathrm{R}(\mathcal{O}_{\overline{\mathbb{A}}_X(Y)}, \overline{d}') \sim \mathrm{R}(\mathcal{O}_{\overline{\mathbb{A}}_X(Y)}, \nabla') \tag{5.4.3}$$

*in the derived category of sheaves of $\mathcal{O}_X$-modules.*

*Proof.* Recall that $\nabla \colon \mathcal{O}_{\overline{\mathbb{A}}_X(Y)} \to \mathcal{O}_{\overline{\mathbb{A}}_X(Y)} \otimes \Omega^1_{\overline{\mathbb{D}}/X}$ is the unique connection such that

(1) $\nabla(f) = df$ if $f \in \mathcal{O}_{\overline{\mathbb{D}}}$;

(2) $\nabla$ satisfies the Leibniz rule;

(3) $\nabla(\overline{\delta}(f)) = -\overline{f}^{p-1}\nabla(\overline{f})$ if $f \in \mathcal{O}_{\mathbb{A}}$.

Let $\tilde{d}'$ denote the $\alpha$-transform of $\overline{d}'$ and let $\nabla'$ be the Higgs field associated to $\nabla$ as explained in Theorem 3.3.3. This satisfies, and is determined by, the the conditions (1)–(3) characterizing the connection $\nabla$, with $\Omega^1_{\overline{Y}/\overline{S}}$ in place of $\Omega^1_{\overline{\mathbb{D}}/X}$. To prove that $\tilde{d}' = \nabla'$, we shall show that $\tilde{d}'$ shares these properties.

First note that if $f \in \mathcal{O}_{\overline{Y}}$, then $\overline{d}' f = 0$, and if $f \in \mathcal{O}_{\overline{\mathbb{D}}}$, then $\overline{d}' f$ identifies with the usual exterior differential of $f$. This implies that, if we work in local coordinates with associated derivations $D_1, \ldots, D_n$, then $D_i^p(f) = 0$ for all $i$ and all $f \in \mathcal{O}_{\overline{\mathbb{D}}}$, and hence that $\alpha(\overline{d}') = \overline{d}'$ on $\mathcal{O}_{\overline{\mathbb{D}}}$. This implies condition (1). For (2) and (3), we use Proposition 5.4.1 to compute the $\alpha$-transform $\tilde{d}'$. Since $\overline{d}'$ satisfies the Leibniz rule, the action $r_{\overline{d}'}$ is a ring homomorphism, and since $I^2_{\overline{Y}/\overline{\mathbb{A}}} \subseteq \mathrm{Fil}^2_\delta$ it follows that $\alpha(\overline{d}') = (\mathrm{id} \otimes \pi_\delta)(p_1^* - r_{\overline{d}'})$ satisfies the Leibniz rule (2).

To prove (3), we compute as follows. If $f \in \mathcal{O}_{\mathbb{A}}$, let $g := p_2^*(f) - p_1^*(f) \in I_{Y/\mathbb{A}_Y(1)}$, whose image $\overline{g}$ in $I_{\overline{Y}/\overline{\mathbb{A}}_Y(1)}$ is equal to $r_{\overline{d}'}(\overline{f}) - p_1^*(\overline{f})$ modulo $\mathrm{Fil}^2_\delta$, where $\overline{f}$ is the image of $f$ in $\mathcal{O}_{\overline{\mathbb{A}}}$. Thus Proposition 5.4.1 tells us that $\tilde{d}'(\overline{f}) = \pi_\delta(\overline{g})$. Calculating in $\mathcal{O}_{\mathbb{A}_Y(1)}$, we have $p_2^*(f) = p_1^*(f) + g$, and by Proposition 1.7.15:

$$\delta(p_2^*(f)) \equiv \delta(p_1^*(f)) - p_1^*(f)^{p-1}g \pmod{\mathrm{Fil}^2_\delta}.$$

Reducing modulo $p$,

$$\begin{aligned} \overline{\delta}(p_2^*(f)) - \overline{\delta}(p_1^*(f)) &\equiv -p_1^*(\overline{f})^{p-1}\overline{g} \pmod{\mathrm{Fil}^2_\delta} \\ \pi_\delta(\overline{\delta}(p_2^*(f)) - \overline{\delta}(p_1^*(f))) &= -\overline{f}^{p-1}\pi_\delta(\overline{g}) \\ \tilde{d}'(\overline{\delta}(f)) &= -\overline{f}^{p-1}\tilde{d}'(\overline{f}). \end{aligned}$$

This concludes the proof that $\nabla'$ is the $\alpha$-transform of $\overline{d}'$. The existence of the derived isomorphism of the corresponding Higgs complexes then follows from Proposition 5.2.3. $\square$

## 5.5 A kernel for the F-transform

Let $Y/S$ be a $p$-completely smooth crystalline prism. Shiho's Theorem 5.1.4, gives an equivalence of categories

$$\mathrm{MICP}(Y^{(p)}/S) \to \mathrm{MIC}(Y/S).$$

The description of the functor in this direction is explicit and simple, but that of the inverse is less so. Here we shall describe a more symmetric construction of this equivalence, as a geometric transform. Namely, we will see that there is a crystalline prism $\mathbb{\Delta}$ with affine morphisms $\pi_Y \colon \mathbb{\Delta} \to Y$ and $\pi_{Y^{(p)}} \colon \mathbb{\Delta} \to Y^{(p)}$, and that the F-transform and its inverse are given by (suitably defined) functors $\pi_{Y*}\pi^*_{Y^{(p)}}$ and $\pi_{Y^{(p)}*}\pi^*_Y$. This prism $\mathbb{\Delta}$ can be viewed as a prismatic analog of the "torsor of Frobenius liftings" used in [34, §2.4] and its $p$-adic versions in [35] and [40].

If $Y/S$ is a $p$-completely smooth morphism of crystalline prisms, we let $Z := Y \times_S Y^{(p)}$, noting that there are isomorphisms:

$$\begin{array}{rcl} \Omega^i_{Z/Y^{(p)}} & \cong & pr_Y^*(\Omega^i_{Y/S}) \\ \Omega^i_{Z/Y} & \cong & pr^*_{Y^{(p)}}(\Omega^i_{Y^{(p)}/S}), \end{array}$$

which we may use without further mention in what follows. Let $\mathbb{\Delta}_\Gamma(Z)$ be the prismatic envelope of the graph $\Gamma \subset Z$ of the relative Frobenius map $F_{Y/S} \colon Y \to Y^{(p)}$. There are morphisms

$$\begin{array}{rcl} \pi_Z \colon \mathbb{\Delta}_\Gamma(Z) & \to & Z \\ \pi_Y \colon \mathbb{\Delta}_\Gamma(Z) & \to & Z \to Y \\ \pi_{Y^{(p)}} \colon \mathbb{\Delta}_\Gamma(Z) & \to & Z \to Y^{(p)}. \end{array}$$

Note that these morphisms are affine, because the maps $\overline{\mathbb{\Delta}}_\Gamma(Z) \to \overline{\Gamma}$, $\overline{\Gamma} \to \overline{Y}$, and $\overline{\Gamma} \to \overline{Y}^{(p)}$ are all affine.

**Proposition 5.5.1.** *With the notation above, let $\mathcal{B}_{Y/S} := \mathcal{O}_{\mathbb{\Delta}_\Gamma(Z)}$, which we view as a sheaf of $\mathcal{O}_Z$-algebras, a sheaf of left $\mathcal{O}_Y$-modules, and a sheaf of right $\mathcal{O}_{Y^{(p)}}$-modules.*

*(1) The maps $\pi_Y$ and $\pi_{Y(p)}$ are affine and $p$-completely flat. Furthermore, there is a natural isomorphism:*

$$F^*_{Y/S}(\mathcal{O}_{\mathbb{\Delta}_{Y^{(p)}}(1)}) \cong \mathcal{B}_{Y/S},$$

*where $\mathbb{\Delta}_{Y^{(p)}}(1)$ is the prismatic envelope of the diagonal of $Y^{(p)} \times_S Y^{(p)}$.*

*(2) The $p$-differential $d' \colon \mathcal{O}_Z \to \Omega^1_{Z/Y} \cong \mathcal{O}_Z \otimes \Omega^1_{Y^{(p)}/S}$ extends uniquely to an integrable quasi-nilpotent $p$-connection on $\mathcal{B}_{Y/S}$ viewed as a right $\mathcal{O}_{Y^{(p)}}$-module:*

$$d' \colon \mathcal{B}_{Y/S} \to \mathcal{B}_{Y/S} \otimes_{\mathcal{O}_{Y^{(p)}}} \Omega^1_{Y^{(p)}/S}.$$

(3) *The differential* $d\colon \mathcal{O}_Z \to \Omega^1_{Z/Y^{(p)}} \cong \Omega^1_{Y/S} \otimes \mathcal{O}_Z$ *extends uniquely to an integrable quasi-nilpotent connection on* $\mathcal{B}_{Y/S}$ *viewed as a left* $\mathcal{O}_Y$*-module:*

$$d\colon \mathcal{B}_{Y/S} \to \Omega^1_{Y/S} \otimes_{\mathcal{O}_Y} \mathcal{B}_{Y/S}.$$ [9]

(4) *The map* $d'$ *in (2) is* $\mathcal{O}_Y$*-linear, the map* $d$ *in (3) is* $\mathcal{O}_{Y^{(p)}}$*-linear, and the following diagram commutes:*

$$\begin{array}{ccc}
\mathcal{B}_{Y/S} & \xrightarrow{\quad d \quad} & \Omega^1_{Y/S} \otimes \mathcal{B}_{Y/S} \\
\Big\downarrow{\scriptstyle d'} & & \Big\downarrow{\scriptstyle \mathrm{id} \otimes d'} \\
\mathcal{B}_{Y/S} \otimes \Omega^1_{Y^{(p)}/S} & \xrightarrow{d \otimes \mathrm{id}} & \Omega^1_{Y/S} \otimes \mathcal{B}_{Y/S} \otimes \Omega^1_{Y^{(p)}/S}.
\end{array}$$

*Proof.* We have already seen that $\pi_Y$ and $\pi_{Y^{(p)}}$ are affine. For the flatness, observe that the map $F_{Y/S} \times \mathrm{id}\colon Z \to Y^{(p)} \times_S Y^{(p)}$ is $p$-completely flat and that the diagram below is Cartesian

$$\begin{array}{ccc}
Y & \xrightarrow{\quad \Gamma \quad} & Z \\
\Big\downarrow{\scriptstyle F_{Y/S}} & & \Big\downarrow{\scriptstyle F_{Y/S} \times \mathrm{id}} \\
Y^{(p)} & \xrightarrow{\ \Delta\ } & Y^{(p)} \times_S Y^{(p)},
\end{array}$$

where the bottom arrow is the diagonal embedding of $Y^{(p)}$ in $Y^{(p)} \times_S Y^{(p)}$. Since $Y^{(p)}/S$ is $p$-completely smooth, $\Delta$ is regular, hence so is $\Gamma$. Let $\mathbb{A}_{Y^{(p)}}(1)$ denote the prismatic envelope of this embedding. It follows from Theorem 2.3.8 that the induced map

$$\mathbb{A}_\Gamma(Z) \to Z \times_{Y^{(p)}(1)} \mathbb{A}_{Y^{(p)}}(1) \cong Y \times_{Y^{(p)}} \mathbb{A}_{Y^{(p)}}(1)$$

is an isomorphism. Thus we have found the isomorphism claimed in statement (1). Since $\Gamma \to Z$ is also a regular immersion, Theorem 2.3.8 implies that $\overline{\mathbb{A}}_\Gamma(Z) \to \overline{\Gamma} = \overline{Y}$ is flat, and hence $\overline{\mathbb{A}}_\Gamma(Z) \to \overline{Y}^{(p)}$ is also flat. Thus $\pi_Y$ and $\pi_{Y^{(p)}}$ are $p$-completely flat.

The existence of the $p$-connection in statement (2) can be deduced from Proposition 4.2.1. However, the existence of the connection in statement (3) needs another justification, and, for the sake of symmetry, we explain how to construct $d$ and $d'$ simultaneously. We begin by observing that, since $\Gamma \subset Z$ is an immersion of $p$-completely smooth crystalline prisms and $\overline{\Gamma} \to \overline{Z}$ is a regular immersion, Theorem 2.3.13 implies that $\mathbb{A}_\Gamma(Z)$ is isomorphic to the divided power envelope $\mathbb{PD}_{\tilde{\Gamma}}(\mathbb{D}_\Gamma(Z))$ of the immersion $\tilde{\Gamma} \subset \mathbb{D}_\Gamma(Z)$ induced by $\Gamma \subset Z$ (see Proposition 2.2.14). If $f$ is a local section of $\mathcal{O}_{Y^{(p)}}$, then

$$\epsilon_f := 1 \otimes f - F^*_{Y/S}(f) \otimes 1$$

belongs to the ideal $I_{\Gamma/Z}$, and in fact this ideal is generated by the family of such elements. Then the algebra $\mathcal{O}_{\mathbb{D}_\Gamma(Z)}$ is topologically generated by elements of the form $\rho(\epsilon_f)$ for all such $f$, and the algebra $\mathcal{O}_{\mathbb{A}_\Gamma(Z)}$ is topologically generated by the divided powers of such elements. Then, in $\mathbf{Q} \otimes \mathcal{B}_{Y/S} \otimes_{\mathcal{O}_Z} \Omega^1_{Z/S}$, we have

$$\begin{aligned} d\rho(\epsilon_f) &= p^{-1} d(1 \otimes f - F^*_{Y/S}(f) \otimes 1) \\ &= p^{-1} \otimes df - p^{-1} dF^*_{Y/S}(f) \otimes 1 \\ &= p^{-1} \otimes df - 1 \otimes \zeta(df). \end{aligned}$$

[9]Here we must break our convention of writing differentials on the right.

Applying the projections

$$\begin{aligned} \mathcal{B}_{Y/S} \otimes_{\mathcal{O}_Z} \Omega^1_{Z/S} &\to \mathcal{B}_{Y/S} \otimes_{\mathcal{O}_Z} \Omega^1_{Z/Y^{(p)}} &\cong \mathcal{B}_{Y/S} \otimes_{\mathcal{O}_Y} \Omega^1_{Y/S} \\ \mathcal{B}_{Y/S} \otimes_{\mathcal{O}_Z} \Omega^1_{Z/S} &\to \mathcal{B}_{Y/S} \otimes_{\mathcal{O}_Z} \Omega^1_{Z/Y} &\cong \mathcal{B}_{Y/S} \otimes_{\mathcal{O}_{Y^{(p)}}} \Omega^1_{Y^p/S}, \end{aligned}$$

we find:

$$\begin{aligned} d\rho(\epsilon_f) &= -\zeta(df) \otimes 1 \in \Omega^1_{Y/S} \otimes \mathcal{B}_{Y/S} \\ d'\rho(\epsilon_f) &= 1 \otimes df \in \mathcal{B}_{Y/S} \otimes \Omega^1_{Y^{(p)}/S}. \end{aligned} \tag{5.5.1}$$

We see that $d\rho(\epsilon_f) \in \Omega^1_{Y/S} \otimes \mathcal{B}_{Y/S}$ and that $d'\rho(\epsilon_f) \in \mathcal{B}_{Y/S} \otimes \Omega^1_{Y^{(p)}/S}$, and hence the same is true for all the divided powers of $\rho(\epsilon_f)$. This implies statements (2) and (3).

The linearity assertions in (4) are immediate from the definitions, and the diagram commutes when restricted to $\mathcal{O}_Z$. It follows that it commutes on all of $\mathcal{B}_{Y/S}$. □

We can now describe the F-transform in the following way. Recall that $\mathcal{B}_{Y/S}$ is regarded as a left $\mathcal{O}_Y$-module and a right $\mathcal{O}_{Y^{(p)}}$-module, and, if $E$ (resp. $E'$) is an $\mathcal{O}_Y$-module (resp., $\mathcal{O}_{Y^{(p)}}$-module), we make identifications:

$$E \otimes_{\mathcal{O}_Y} \mathcal{B}_{Y/S} \cong \pi^*_Y(E), \quad \mathcal{B}_{Y/S} \otimes_{\mathcal{O}_{Y^{(p)}}} E' \cong \pi^*_{Y^{(p)}}(E').$$

Suppose $(E', \nabla')$ is an object of MICP$(Y^{(p)}/S)$. Endow $\mathcal{B}_{Y/S} \otimes E'$ with the tensor product $p$-connection coming from the $p$-connections on $E'$ and on $\mathcal{B}_{Y/S}$:

$$\nabla' \colon \mathcal{B}_{Y/S} \otimes E' \to \mathcal{B}_{Y/S} \otimes E' \otimes \Omega^1_{Y^{(p)}/S}.$$

Since the connection $d$ on $\mathcal{B}_{Y/S}$ is $\mathcal{O}_{Y^{(p)}}$-linear, it induces a connection

$$\nabla := d \otimes \mathrm{id} \colon \mathcal{B}_{Y/S} \otimes E' \to \Omega^1_{Y/S} \otimes \mathcal{B}_{Y/S} \otimes E'.$$

It follows from the commutative diagram in statement (4) that $\nabla$ and $\nabla'$ commute and that $\nabla$ induces a connection on the cohomology sheaves of the complex $(\mathcal{B}_{Y/S} \otimes E' \otimes \Omega^{\bullet}_{Y^{(p)}/S}, d')$ and in particular on $(\mathcal{B}_{Y/S} \otimes E')^{\nabla'}$.

Similarly, if $(E, \nabla)$ is an object of $\mathrm{MIC}(Y/S)$, we endow $E \otimes \mathcal{B}_{Y/S}$ with the tensor product connection

$$\nabla\colon E \otimes \mathcal{B}_{Y/S} \to \Omega^1_{Y/S} \otimes E \otimes \mathcal{B}_{Y/S}$$

coming from the connections on $\mathcal{B}_{Y/S}$ and on $E$. Since the $p$-connection on $\mathcal{B}_{Y/S}$ is $\mathcal{O}_Y$-linear, we also find a $p$-connection:

$$\nabla' := \mathrm{id} \otimes d'\colon E \otimes \mathcal{B}_{Y/S} \to E \otimes \mathcal{B}_{Y/S} \otimes \Omega^1_{Y^{(p)}/S},$$

which commutes with $\nabla$. Then $\nabla'$ induces a $p$-connection on the cohomology sheaves of the complex $(E \otimes \mathcal{B}_{Y/S} \otimes \Omega^{\bullet}_{Y/S}, d)$ and in particular on $(E \otimes \mathcal{B}_{Y/S} \otimes E)^{\nabla}$.

**Theorem 5.5.2.** *With the constructions described in the previous paragraph, the following results hold.*

*(1) If $(E', \nabla') \in \mathrm{MICP}(Y^{(p)}/S)$, let*

$$(E, \nabla) := ((\mathcal{B}_{Y/S} \otimes E')^{\nabla'}, d \otimes \mathrm{id}).$$

*Then the map*

$$E \to (\mathcal{B}_{Y/S} \otimes E' \otimes \Omega^{\bullet}_{Y^{(p)}/S}, d')$$

*is a strict quasi-isomorphism, and the natural $\mathcal{O}_Y$-linear map*

$$E \otimes \mathcal{B}_{Y/S} \to \mathcal{B}_{Y/S} \otimes E'$$

*is an isomorphism. Furthermore, $(E, \nabla)$ is isomorphic to the F-transform of $(E', \nabla')$.*

*(2) If $(E, \nabla) \in \mathrm{MIC}(Y/S)$, let*

$$(E', \nabla') := ((E \otimes \mathcal{B}_{Y/S})^{\nabla}, \mathrm{id} \otimes d').$$

*Then the map*

$$E' \to (\Omega^{\bullet}_{Y/S} \otimes E \otimes \mathcal{B}_{Y/S}, d)$$

*is a strict quasi-isomorphism, and the natural $\mathcal{O}_{Y^{(p)}}$-linear map*

$$\mathcal{B}_{Y/S} \otimes E' \to E \otimes \mathcal{B}_{Y/S}$$

*is an isomorphism.*

*(3) The functor $\mathrm{MICP}(Y^{(p)}/S) \to \mathrm{MIC}(Y/S)$*

$$(E', \nabla') \mapsto ((\mathcal{B}_{Y/S} \otimes E')^{\nabla'}, \nabla)$$

*is an equivalence of categories, with quasi-inverse*

$$(E, \nabla) \mapsto ((E \otimes \mathcal{B}_{Y/S})^{\nabla}, \nabla').$$

*Proof.* The map $\pi_Y\colon Z \to Y$ is $p$-completely smooth, with a section $\Gamma\colon Y \to Z$. Thus we are in the situation of Lemma 4.3.2, with a change of notation: $S \mapsto Y$, $s \mapsto \Gamma$, and $E := \pi^*_{Y^{(p)}}(E')$. The first sentence of statement (1) follows.

To see that $((\mathcal{B}_{Y/S} \otimes E')^{\nabla'}, \nabla)$ is the F-transform of $(E', \nabla')$, let $\Gamma'\colon Y \to \mathbb{A}_\Gamma(Z)$ be the canonical section of $\pi_Y$, arising from the fact that $Y$ is itself a $Y$-prism. Lemma 4.3.2 tells us that the map

$$E := (\mathcal{B}_{Y/S} \otimes E')^{\nabla'} \to \Gamma'^*(\mathcal{B}_{Y/S} \otimes E')$$

is an isomorphism. Since $\Gamma'^*(\mathcal{B}_{Y/S} \otimes E') \cong F^*_{Y/S}(E')$, we find an isomorphism $E \to F^*_{Y/S}(E')$. It remains to show that the connection on $(\mathcal{B}_{Y/S} \otimes E')^{\nabla'}$ agrees with the F-transform connection.

We have a split exact sequence

$$0 \to I_{Y/\mathbb{A}} \to \mathcal{B}_{Y/S} \to \mathcal{O}_Y \to 0,$$

hence also:

$$0 \to I_{Y/\mathbb{A}} \hat{\otimes} E' \to \mathcal{B}_{Y/S} \hat{\otimes} E' \to F^*_{Y/S}(E') \to 0.$$

There are also maps:

$$\begin{array}{rcl} \Omega^1_{Y/S} \otimes \mathcal{B}_{Y/S} \hat{\otimes} E' & \xrightarrow{\mathrm{id}\otimes\Gamma'^*} & \Omega^1_{Y/S} \otimes F^*_{Y/S}(E') \\ \mathcal{B}_{Y/S} \hat{\otimes} E' \otimes \Omega^1_{Y^{(p)}/S} & \xrightarrow{\Gamma'^*\otimes\mathrm{id}} & F^*_{Y/S}(E') \otimes \Omega^1_{Y^{(p)}/S} \\ F^*_{Y/S}(E') \otimes \Omega^1_{Y^{(p)}/S} & \xrightarrow{switch} & \Omega^1_{Y^{(p)}/S} \otimes F^*_{Y/S}(E'), \end{array}$$

the last of which is an isomorphism.

**Lemma 5.5.3.** *The following diagram commutes:*

$$\begin{array}{ccccc} \mathcal{B}_{Y/S} \otimes E' & \xrightarrow{\nabla'} & \mathcal{B}_{Y/S} \otimes E' \otimes \Omega^1_{Y^{(p)}/S} & \xrightarrow{\Gamma'^* \otimes \mathrm{id}} & F^*_{Y/S}(E') \otimes \Omega^1_{Y^{(p)}/S} \\ \downarrow \nabla & & & & \downarrow switch \\ \Omega^1_{Y/S} \otimes \mathcal{B}_{Y/S} \otimes E' & \xrightarrow{\mathrm{id} \otimes \Gamma'^*} & \Omega^1_{Y/S} \otimes F^*_{Y/S}(E') & \xleftarrow{-\zeta \otimes \mathrm{id}} & \Omega^1_{Y^{(p)}/S} \otimes F^*_{Y/S}(E'). \end{array}$$

*Proof.* We have seen in the proof of Theorem 2.3.13 that the ideal $I_{Y/\mathbb{A}_\Gamma(Z)}$ is topologically generated by the divided powers of elements of the form $g := \rho(\epsilon_f)$, where $f \in \mathcal{O}_{Y^{(p)}}$ and $\epsilon_f := 1 \otimes f - F^*_{Y/S}(f) \otimes 1$. Hence it will suffice to check the formula for elements of the form $g^{[n]}e'$, for all $f \in \mathcal{O}_{Y^{(p)}}$ and $e' \in E'$.

Note first that

$$\nabla(g^{[n]}e') = g^{[n-1]}dg \otimes e' + g^{[n]}\nabla(e'),$$

which again belongs to $I_{Y/\mathbb{A}_\Gamma(Z)}E'$ if $n > 1$, and the analogous statement holds for $\nabla'$. Thus $\nabla(g^{[n]}e')$ and $\nabla'(g^{[n]}e')$ both vanish when pulled back via $\Gamma'$, and the claim becomes trivially true on such elements. Thus it suffices to treat the case $n = 1$. By equation (5.5.1),

$$\begin{aligned}(\mathrm{id}\otimes\Gamma'^{*})(\nabla(ge')) &= (\mathrm{id}\otimes\Gamma'^{*})(dg\otimes e')\\ &= (\mathrm{id}\otimes\Gamma'^{*})(d\rho(\epsilon_f)\otimes e')\\ &= -\zeta(df)\otimes e'.\end{aligned}$$

Similarly:

$$\begin{aligned}(\Gamma'^{*}\otimes\mathrm{id})(\nabla'(ge')) &= (\Gamma'^{*}\otimes\mathrm{id})(e'\otimes d'g + g\otimes\nabla'(e))\\ &= (\Gamma'^{*}\otimes\mathrm{id})(e'\otimes d'(\rho(\epsilon_f))) + \Gamma'^{*}(g)\nabla'(e')\\ &= e'\otimes df.\end{aligned}$$

This proves the lemma. □

Now suppose that $e' \in E'$. By Lemma 4.3.2, we can uniquely write $1\otimes e' = e + e''$, with $\nabla'(e) = 0$ and $e'' \in I_{Y/\mathbb{A}_\Gamma(Z)}E'$. Then, using Lemma 5.5.3, we have:

$$\begin{aligned}\nabla(e) &= \nabla(1\otimes e') - \nabla(e'')\\ &= -\nabla(e'')\\ &= (\zeta\otimes\mathrm{id})\nabla'(e'')\\ &= (\zeta\otimes\mathrm{id})(\nabla'(1\otimes e') - \nabla'(e))\\ &= (\zeta\otimes\mathrm{id})\nabla'(e').\end{aligned}$$

This concludes the proof of statement (1).

We begin the proof of statement (2) with the following special case, which is worth stating separately.

**Lemma 5.5.4.** *The de Rham complex $(\Omega^{\bullet}_{Y/S}\otimes\mathcal{B}_{Y/S}, d)$ of the $\mathcal{O}_Y$-module with connection $(\mathcal{B}_{Y/S}, d)$, viewed as a complex of $\mathcal{O}_{Y^{(p)}}$-modules, is a resolution of $\mathcal{O}_{Y^{(p)}}$, each term of which is $p$-completely flat.*

*Proof.* This statement can be verified étale locally on $Y$. Thus we may assume that $Y/S$ admits a sequence of local coordinates $(x_1, \ldots, x_n)$, and we find a corresponding sequence $(x'_1, \ldots, x'_n)$ for $Y^{(p)}/S$. Let $\xi_i := 1\otimes x'_i - x'_i\otimes 1 \in \mathcal{O}_{Y^{(p)}(1)}$ and $\tilde{\xi}_i := \rho(\xi_i) \in \mathcal{O}_{\mathbb{D}_{Y^{(p)}(1)}}$. Corollary 2.3.15 implies that the structure sheaf $\mathcal{B}'_{Y/S}$ of $\mathbb{A}_{Y^{(p)}(1)}$ identifies with the completed PD-polynomial algebra $\mathcal{O}_{Y^{(p)}}\langle\tilde{\xi}_1, \ldots, \tilde{\xi}_n\rangle$. Similarly, the structure sheaf $\mathcal{B}_{Y/S}$ of $\mathbb{A}_\Gamma(Z)$ identifies with the completed PD-polynomial algebra $\mathcal{O}_Y\langle\tilde{\eta}_1, \ldots, \tilde{\eta}_n\rangle$, where

$$\tilde{\eta}_i = (F_{Y/S}\times\mathrm{id})^*(\tilde{\xi}_i) \quad\text{and}\quad p\tilde{\eta}_i = 1\otimes x'_i - F^*_{Y/S}(x'_i)\otimes 1.$$

The prismatic envelope $\mathcal{B}'_{Y/S}$ has a natural $p$-connection relative to $Y^{(p)}(1)/S$, and if we work relative to the second factor, we find a $p$-connection

$$\nabla'\colon \mathcal{B}'_{Y/S} \to \Omega^1_{Y^{(p)}(1)/Y}\otimes\mathcal{B}'_{Y/S} : \tilde{\xi}\mapsto -dx'_i\otimes 1.$$

Statement (1) of Proposition 5.5.1 provides us with a canonical isomorphism

$$F_{Y/S}^*(\mathcal{B}'_{Y/S}) \to \mathcal{B}_{Y/S} : F_{Y/S}^*(\tilde{\xi}_i) \mapsto \tilde{\eta}_i,$$

and equation (5.5.1) shows that this isomorphism identifies $(\mathcal{B}_{Y/S}, d)$ as the $F$-transform of $(\mathcal{B}'_{Y/S}, \nabla')$. Then Theorem 5.1.6 implies that the map

$$(\mathcal{B}'_{Y/S} \otimes \Omega^{\bullet}_{Y^{(p)}/S}, d') \to F_{Y/S*}(\mathcal{B}_{Y/S} \otimes \Omega^{\bullet}_{Y/S}, d)$$

is a strict quasi-isomorphism. The prismatic Poincaré lemma (Theorem 4.3.1) tells us that the complex on the left is a strict resolution of $\mathcal{O}_{Y^{(p)}}$, and hence the same is true of the complex on the right. It is clear that each of its terms is $p$-completely flat over $\mathcal{O}_Y$ and hence also over $\mathcal{O}_{Y^{(p)}}$. □

To prove the general case of statement (2), we use Shiho's Theorem 5.1.4, which tells us that every $(E, \nabla) \in \mathrm{MIC}(Y/S)$ is the F-transform of some $(E', \nabla') \in \mathrm{MICP}(Y^{(p)}/S)$. Then, by statement (1), we have isomorphisms:

$$\begin{array}{rcl} (E, \nabla) & \cong & ((\mathcal{B}_{Y/S} \otimes E')^{\nabla'}, d \otimes \mathrm{id}) \\ E \otimes \mathcal{B}_{Y/S} & \cong & \mathcal{B}_{Y/S} \otimes E' \\ (\Omega^{\bullet}_{Y/S} \otimes E \otimes \mathcal{B}_{Y/S}, d) & \cong & (\Omega^{\bullet}_{Y/S} \otimes \mathcal{B}_{Y/S} \otimes E', d \otimes \mathrm{id}_{E'}). \end{array}$$

Lemma 5.5.4 tells us that $(\Omega^{\bullet}_{Y/S} \otimes \mathcal{B}_{Y/S}, d)$ is a strict resolution of $\mathcal{O}_{Y^{(p)}}$ all of whose terms are $p$-completely flat, and hence $(\Omega^{\bullet}_{Y/S} \otimes \mathcal{B}_{Y/S} \otimes E', d \otimes \mathrm{id}_{E'})$ is a strict resolution of $E'$. Statement (2) follows, and Theorem 5.1.4 now also implies statement (3). □

## 5.6 Comparing tubular groupoids

Let us explain the relationships among the various groupoids we have now constructed, as well as how they relate to the relative Frobenius morphism. As we shall see, these constructions allow for a geometric framework for a proof of Shiho's Theorem 5.1.4.

Let $Y/S$ be a $p$-completely smooth morphism of crystalline prisms and $X \to Y$ a closed subscheme, with $X/\overline{S}$ smooth. The morphisms of tubes, some of which were illustrated in diagram (2.4.1), prolong to morphisms of groupoids:

$$\mathcal{G}_{\mathbb{A}X/Y} \to \mathcal{G}_{\mathbb{D}X/Y} \to \mathcal{G}_{\mathbb{P}\mathbb{D}X/Y} \to \mathcal{G}_{\mathbb{F}X/Y}. \tag{5.6.1}$$

In particular, an action of $\mathcal{G}_{\mathbb{P}\mathbb{D}X/Y}$ restricts to an action of $\mathcal{G}_{\mathbb{D}X/Y}$, which in turn restricts to an action of $\mathcal{G}_{\mathbb{A}X/Y}$. When $X = \overline{Y}$ and we are working with actions on $\mathcal{O}_Y$-modules, the latter restriction amounts to restricting an action of the ring $\mathcal{HD}_{\mathbb{D}_Y}$ to the actions of the $p$-derivations, *i.e.*, to the underlying $p$-connection. The following result, inspired by an argument in [25], gives a criterion for extending actions of $\mathcal{G}_{\mathbb{A}_Y}$ to the ring of dilatation differential operatros $\mathcal{G}_{\mathbb{D}_Y}$ defined in section § 4.4.

**Proposition 5.6.1.** *Suppose that $E$ is a $p$-torsion free sheaf of $\mathcal{O}_Y$-modules endowed with an integrable $p$-connection $\nabla'$ whose action on $E/pE$ has level less than $p-1$ (resp. $p-2$). Then $\nabla'$ extends uniquely to an action of the ring $\mathcal{D}_{\mathbb{D}_Y}$ (resp.to $\mathcal{HD}_{\mathbb{D}_Y}$).*

*Proof.* The uniqueness follows from the assumption that $E$ is $p$-torsion free, and it implies that we can work locally. Choose local coordinates $(t_1, \ldots, t_n)$ on $Y/S$, let $(\partial_1, \ldots, \partial_n)$ be the basis for $\mathcal{T}_{Y/S}$ dual to $(dt_1, \ldots, dt_n)$, and write $\nabla_i$ for the corresponding endomorphism of $E'$ induced by $\nabla'$. To say that $(E/pE, \nabla')$ has level less than $\ell$ means that the ideal $\{\oplus S^i T_{Y/S} : i \geq \ell\}$ annihilates $E/pE$, *i.e.*, that $\nabla_1^{I_1} \cdots \nabla_n^{I_n} E \subseteq pE$ whenever $I_1 + \cdots + I_n \geq \ell$. In particular, if the level is less than $p-1$, then $\nabla_i^{p-1} E \subseteq pE$. It follows that for each $m > 0$

$$\nabla_i^{p^m} E = (\nabla_i^{p-1})^{p^{m-1}} (\nabla_i^{p-1})^{p^{m-2}} \cdots (\nabla_i^{p-1}) E \subseteq p^{p^{m-1}} p^{p^{m-2}} \cdots pE.$$

Recalling that $\mathrm{ord}_p(p^m!) = p^{m-1} + \cdots + 1$, we see that $\nabla_i^{p^m} E \subseteq p^m! E$ and then that $\nabla_i^n \subseteq n! E$ for all $n$. It then follows from the formulas in (4.4.9) that the ring $\mathcal{D}_{\mathbb{D}}$ operates on $E$. If the level is less than $p-2$, we have $\nabla_i^{p-2} E \subseteq pE$, and then

$$\begin{aligned}
\nabla_i^{p^m} E &= (\nabla_i^{p-2}\nabla_i)^{p^{m-1}} (\nabla_i^{p-2}\nabla_i)^{p^{m-2}} \cdots (\nabla_i^{p-2})\nabla_i E \\
&\subseteq p^{p^{m-1}} p^{p^{m-2}} \cdots p \nabla_i^{p^{m-1}} E \\
&\subseteq p^m! p^{m-2} E.
\end{aligned}$$

This implies that the operation of $\mathcal{D}_{\mathbb{D}_Y}$ is quasi-nilpotent, and hence that it extends to $\mathcal{HD}_{\mathbb{D}_Y}$. □

In the prismatic context, the Frobenius lifting induces additional morphisms of groupoids. This gives another approach to the F-transform and Shiho's theorem. The main point is to see how a quasi-nilpotent connection gives rise to descent data for the relative Frobenius morphism, as explained by the morphism of groupoids $u$ in the following proposition.

**Proposition 5.6.2.** *Let $Y/S$ be a $p$-completely smooth morphism of crystalline prisms, let $F_{Y/S} : Y \to Y^{(p)}$ be the associated relative Frobenius morphism, and let $\mathcal{G}_{\Phi_Y}$ be the groupoid over $Y$ corresponding to $F_{Y/S}$ as explained in Example 1.4.1 (and justified by Lemma 5.6.4 below.) There are corresponding*

*morphism of groupoids:*

$$\begin{array}{ccccc} & & \mathcal{G}_{\Phi Y} & & \\ & & \downarrow{\scriptstyle u} & & \\ \mathcal{G}_{\Delta Y} & \xrightarrow{\Psi} & \mathcal{G}_{\mathbb{PD}Y} & & \\ & {\scriptstyle \Delta(F_{Y/S})}\searrow & \downarrow{\scriptstyle \Phi} & \searrow{\scriptstyle \mathbb{PD}(F_{Y/S})} & \\ & & \mathcal{G}_{\Delta Y^{(p)}} & \xrightarrow{\Psi^{(p)}} & \mathcal{G}_{\mathbb{PD}Y^{(p)}}. \end{array}$$

*The morphisms defining $\Phi$ are $p$-completely faithfully flat. Moreover, the composite $\Phi \circ u$ factors through the identity section.*

*Proof.* The existence of the lower rhombus is a consequence of Corollary 2.4.1, applied with $Y$ replaced by $Y(1)$ and $X$ by the diagonal embedding of $\overline{Y}$. That corollary also explains why $\Phi$ is $p$-completely faithfully flat. Using the same construction with $Y(2)$ in place of $Y(1)$, one can check that $\Phi$ and $\Psi$ are compatible with the composition laws and identity sections.

To construct the morphism $u$, we use the following lemmas.

**Lemma 5.6.3.** *Let $X/\overline{S}$ be a morphism of schemes in characteristic $p$, and let $F_{X/\overline{S}}\colon X \to X^{(p)}$ be the relative Frobenius morphism, and $\pi\colon X^{(p)} \to X$ the base change map. Then the scheme theoretic image of the Frobenius endomorphism of $X \times_{X^{(p)}} X$ is contained in the diagonal.*

*Proof.* We have a commutative diagram:

$$\begin{array}{ccccc} & & X^{(p)} & \xrightarrow{\pi} & X \\ & \nearrow & \downarrow{\scriptstyle \Delta_{X^{(p)}}} & & \downarrow{\scriptstyle \Delta_X} \\ X \times_{X^{(p)}} X & \xrightarrow{F_{X/S} \times F_{X/S}} & X^{(p)} \times_{X^{(p)}} X^{(p)} & \xrightarrow{\pi \times \pi} & X \times_{X^{(p)}} X. \end{array}$$

The composition along the bottom is the absolute Frobenius endomorphism of $X \times_{X^{(p)}} X$. □

**Lemma 5.6.4.** *Let $Y/S$ be a $p$-completely smooth morphism of crystalline prisms, with relative Frobenius morphism $F_{Y/S}\colon Y \to Y^{(p)}$. Then $Y \times_{Y^{(p)}} Y$ inherits the structure of a crystalline prism, and the diagonal embedding $Y \to Y \times_{Y^{(p)}} Y$ is defined by a PD-ideal.*

*Proof.* Let $(Z, F_Z) := (Y \times_{Y^{(p)}} Y, F_Y \times_{F_{Y^{(p)}}} F_Y)$. Since $F_{Y/S}$ is $p$-completely flat, this fiber product is again $p$-torsion free, by Proposition 7.3.3, and $F_Z$ is a Frobenius lift. By Proposition 2.1.4, the scheme theoretic image of $\overline{Z}$ by $F_Z$ is $Z_{\mathbb{PD}}$, the smallest PD-subscheme of $Z$. Applying the previous lemma with $X = \overline{Y}$, we see that this image is contained in the diagonal embedding of $Y$ in $Z$. Thus the ideal of $Y$ is contained in the maximal PD-ideal $I_{\mathbb{PD}}$ of $\mathcal{O}_Z$. Since $Y$ is $p$-torsion free, it follows that its ideal is also a PD-ideal. □

The map $Y \times_{Y^{(p)}} Y \to Y \times_S Y$ sends the diagonal of $Y \times_{Y^{(p)}} Y$ to the diagonal of $Y \times_S Y$, and since the former is defined by a PD-ideal, this map factors uniquely through a map $u$ as claimed. Repeating this argument with $Y \times_{Y^{(p)}} Y \times_{Y^{(p)}} Y$ in place of $Y \times_{Y^{(p)}} Y$ shows that $u$ is compatible with the composition laws.

The composite $\Phi \circ u$ is given by the following diagram,

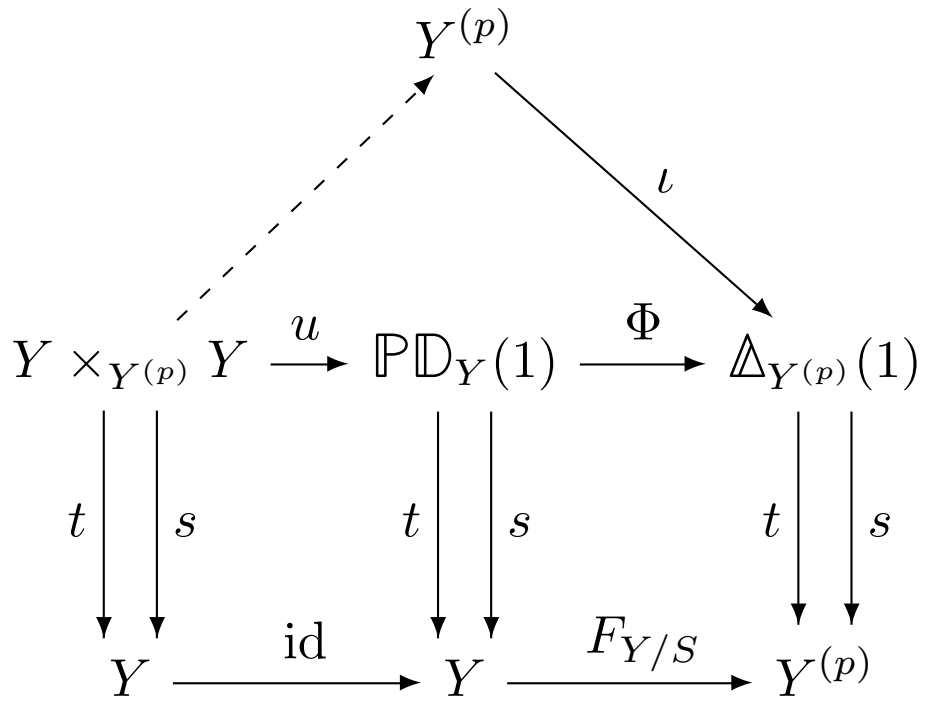

which shows that $\Phi \circ u$ factors through the identity section. □

**Remark 5.6.5.** It may be enlightening to write formulas for the morphisms in Proposition 5.6.2 in terms of the local coordinate description of these groupoids given in equations (4.4.6). The formals for $\Psi$ and $\Phi$ are straightforward:

$$\Psi^\sharp\colon \mathcal{O}_Y\langle \xi_1, \ldots, \xi_n\rangle^\wedge \to \mathcal{O}_Y\langle \eta_1, \ldots, \eta_n\rangle^\wedge \quad : \quad \xi_i \mapsto p\eta_i \tag{5.6.2}$$

$$\Phi^\sharp\colon \mathcal{O}_{Y^{(p)}}\langle \eta'_1, \ldots \eta'_n\rangle^\wedge \to \mathcal{O}_Y\langle \xi_1, \ldots, \xi_n\rangle^\wedge \quad : \quad \eta'_i \mapsto \phi_{Y/S}(\xi_i) \tag{5.6.3}$$

It is apparent from the formulas that $\Psi^\sharp$ is well-defined, but this is less apparent for $\Phi^\sharp$. The point is that $\phi_{Y/S}(\xi'_i) \equiv \xi_i^p$ modulo $pI_Y$, not just modulo $p$. This will imply that $\phi_{Y/S}(\eta'_i) \equiv (p-1)!\xi_i^{[p]}$ modulo $I_Y$, which belongs to the PD ideal $I_Y$. To check the claimed congruence, we compute as follows:

$$\begin{aligned}
\phi_{Y/S}(\xi'_i) &= \phi_{Y/S}(1 \otimes x'_i - x'_i \otimes 1)\\
&= 1 \otimes \phi_{Y/S}(x'_i) - \phi_{Y/S}(x'_i) \otimes 1\\
&= 1 \otimes x_i^p - x_i^p \otimes 1 + p \otimes \delta(x_i) - \delta(x_i) \otimes p\\
&= (x_i \otimes 1 + \xi_i)^p - x_i^p \otimes 1 + p \otimes \delta(x_i) - \delta(x_i) \otimes p\\
&= \sum_{j=0}^{p} \binom{p}{j} (x_i^j \otimes 1)\xi_i^{p-j} - x_i^p \otimes 1 + p \otimes \delta(x_i) - \delta(x_i) \otimes p\\
&= \xi_i^p + \sum_{j=1}^{p-1} \binom{p}{j} (x_i^j \otimes 1)\xi_i^{p-j} + p(1 \otimes \delta(x_i) - \delta(x_i) \otimes 1).
\end{aligned}$$

Thus $\phi_{Y/S}(\xi_i') \equiv \xi_i^p \pmod{pI_Y}$, as claimed.

The morphism $u$ is more difficult to write explicitly, as we explain in the remark below.

**Remark 5.6.6.** The morphism $u$ does not factor through the PD-completion of the diagonal of $Y$ in $\mathbb{PD}_Y(1)$; this explains why the inverse F-transform is only defined for quasi-nilpotent connections. For example, let $Y := \operatorname{Spf} W[t]\hat{}$ and let $\phi(t) = t^p$. Then $Y \times_S Y = \operatorname{Spf} W[t_1, t_2]\hat{}$, where $t_1 := t \otimes 1, t_2 := 1 \otimes t$. Let $\xi := t_2 - t_1$ and view $Y(1)$ as a $Y$-scheme via $p_1$. Thus $t$ identifies with $t_1$, and the ideal of $\mathcal{G}_{\Phi_Y} := Y \times_{Y^{(p)}} Y$ in $Y \times_S Y$ is generated by

$$t_2^p - t_1^p = (t+\xi)^p - t^p = \xi^p + \sum_{i=1}^{p-1} \binom{p}{i} \xi^i t^{p-i} = \xi^p - p\xi y,$$

where

$$y := -(p-1)! \sum_{i=1}^{p-1} \frac{\xi^{i-1} t^{p-i}}{i!(p-i)!} \equiv -t^{p-1} \pmod{\xi}.$$

Then $Y \times_{Y^{(p)}} Y = \operatorname{Spf} W[t, \xi]\hat{}/(\xi^p - p\xi y)$, and in this quotient, we have

$$\begin{aligned} \xi^p &= p\xi y \\ \xi^{pn} &= p^n \xi^n y^n. \end{aligned}$$

Note that the reduction of $\mathcal{G}_{\Phi_Y}$ modulo $p$ is the spectrum of $k[\overline{t}, \overline{\xi}]/(\overline{\xi}^p)$, so that $\overline{\xi}$ and $\overline{y}^n$ for $n \in \mathbf{N}$, are not zero. We have inductive formulas for $u^*$:

$$\begin{aligned} u^*(\xi^{[p]}) &= \frac{\xi y}{(p-1)!} \\ u^*(\xi^{[pn]}) &= \frac{p^n n!}{(pn)!} u^*(\xi^{[n]}) y^n. \end{aligned}$$

Statement (1) of Lemma 1.6.7 tells us that $\operatorname{ord}_p((pn)!) = n + \operatorname{ord}_p(n!)$, so the coefficient in the last equation is a unit. Since $\xi^{[n]}$ is not divisible by $p$ if $n < p$ and $y^n$ is not divisible by $p$ for all $n$, we see by induction on $n$ that $u^*(\xi^{[pn]})$ is not divisible by $p$.

Before explaining our proof of Shiho's theorem, we should make explicit the relationship between the F-transform and the morphism $\Phi$.

**Lemma 5.6.7.** *Let $Y/S$ be $p$-completely smooth and $(E', \nabla')$ an object of $\mathrm{MICP}(Y^{(p)}/S)$, with corresponding action $r'$ of $\mathcal{G}_{\mathbb{A}_{Y^{(p)}}}$. Then the F-transform of $(E', \nabla')$ corresponds to the pullback of $(E', r')$ via the morphism $\Phi$ of Proposition 5.6.2. Furthermore, the action of $\mathcal{G}_{\Phi_Y}$ on $F_{Y/S}^*(E')$ corresponding to its descent data is the pullback of $\Phi^*(r')$ by $u$.*

*Proof.* Let $(E, \nabla)$ be the F-transform of $(E', \nabla')$. By definition, $E = F_{Y/S}^*(E')$, and it remains only to prove that the connection $\nabla$ corresponds to the action $\Phi^*(r')$. This can be checked after restricting to the first infinitesimal neighborhood of the diagonal, where the requisite compatibility is proved in Proposition 5.1.2. The last statement follows from the fact that $\Phi \circ u$ factors through the identity section. $\square$

We can now explain how Shiho's proof of theorem 5.1.4 fits into our current context. His proof is quite detailed, so here we give only a sketch.

*Proof of Theorem 5.1.4.* Let $(E, \nabla)$ be a $p$-completely quasi-coherent sheaf of $\mathcal{O}_Y$-modules with integrable and quasi-nilpotent connection, and let

$$r\colon \mathbf{V}E \times_Y \mathcal{G}_{\mathbb{PD}_Y} \to \mathbf{V}E$$

be the corresponding action of $\mathcal{G}_{\mathbb{PD}_Y}$ on $\mathbf{V}E$. (Here we abusively use the same notation for the groupoid $\mathcal{G}_{\mathbb{PD}_Y}$ and its object of arrows $\mathbb{PD}_Y(1)$.) We claim that there is an object $(E', \nabla') \in \mathrm{MICP}(Y^{(p)}/S)$ whose F-tranform is $(E, \nabla)$. Thanks to Lemma 5.6.7, it will suffice to show that there is a sheaf of $\mathcal{O}_{Y^{(p)}}$-modules $E'$ with an action $r'$ of $\mathcal{G}_{\mathbb{\Delta}_{Y^{(p)}}}$ such that $\Phi^*(E', r') \cong (E, r)$.

The pullback of $r$ by the morphism $u\colon \mathcal{G}_{\mathbb{\Phi}_Y} \to \mathcal{G}_{\mathbb{PD}_Y}$ of Proposition 5.6.2 induces an action $u^*(r)$ of $\mathcal{G}_{\mathbb{\Phi}_Y}$ on $\mathbf{V}E$. This action is in fact descent data for the $p$-completely faithfully flat morphism $F_{Y/S}$. Thanks to Proposition 7.3.8, these data give rise to a $p$-completely quasi-coherent sheaf of $\mathcal{O}_{Y^{(p)}}$-modules $E'$ and an isomorphism $F^*_{Y/S}(E') \cong E$, which we view as a linear scheme $\mathbf{V}E'$ over $Y^{(p)}$ together with an isomorphism $Y \times_{Y^{(p)}} \mathbf{V}E' \cong \mathbf{V}E$. We claim that $r$ induces an action $r'$ of $\mathcal{G}_{\mathbb{\Delta}_{Y^{(p)}}}$ on $\mathbf{V}E'$, *i.e.*, that there is a commutative diagram:

$$\begin{array}{ccc} \mathbf{V}E \times_Y \mathcal{G}_{\mathbb{PD}_Y} & \xrightarrow{\ r\ } & \mathbf{V}E \\ {\scriptstyle \pi \times \Phi}\downarrow & & \downarrow{\scriptstyle \pi} \\ \mathbf{V}E' \times_{Y^{(p)}} \mathcal{G}_{\mathbb{\Delta}_{Y^{(p)}}} & \overset{r'}{\dashrightarrow} & \mathbf{V}E', \end{array}$$

Such an $r'$ will automatically define an action since $r$ does.

Applying Theorem 5.1.8 to the diagonal embedding of $\overline{Y}$ in $Y(1)$, we see that $\mathcal{G}_{\mathbb{PD}_Y} := \mathbb{PD}_Y(1)$ is the F-transform of $\mathcal{G}_{\mathbb{\Delta}_{Y^{(p)}}} := \mathbb{\Delta}_{Y^{(p)}}(1)$, viewed over the crystalline prism $Y^{(p)}(1)$, via the morphisms $t$ and $s$. Thus, Lemma 5.6.7 shows that the morphism $\Phi$ presents $\mathcal{G}_{\mathbb{\Delta}_{Y^{(p)}}}$ as the quotient of $\mathcal{G}_{\mathbb{PD}_Y}$ by an action of the groupoid $\mathcal{G}_{\mathbb{\Phi}_Y}(1)$. By Lemma 5.6.7, this action is the restriction by $u$ of the action of $\mathcal{G}_{\mathbb{PD}_{Y(1)}}$ corresponding to the canonical connection on $\mathbb{PD}_{Y(1)}$. By construction, this action corresponds to the tautological action of $\mathcal{G}_{\mathbb{PD}_{Y(1)}}$ on $\mathbb{PD}_Y(1)$. As explained in Proposition 1.4.12, there is an isomorphism of groupoids $\mathcal{G}_{\mathbb{PD}_{Y(1)}} \cong \mathcal{G}_{\mathbb{PD}_Y}(1)$, and the tautological action of the former corresponds to the conjugation action by the latter. Thus the action of $\mathcal{G}_{\mathbb{\Phi}_Y}(1)$ on $\mathcal{G}_{\mathbb{PD}_Y}$ is the restriction by $u$ of the action of $\mathcal{G}_{\mathbb{PD}_Y}(1)$ given by $h(g_1, g_2) = g_1^{-1} h g_2$. Similarly, by its very definition, the action of $\mathcal{G}_{\mathbb{\Phi}_Y}$ on $\mathbf{V}E$ is given by the restriction of the action of $\mathcal{G}_{\mathbb{PD}_Y}$ via the morphism $u$, *i.e.*, $vg = r(v, u(g))$ for $g \in \mathcal{G}_{\mathbb{PD}_Y}$.

Since $\pi$ is the quotient morphism by an action of $\mathcal{G}_{\mathbb{\Phi}_Y}$ and $\Phi$ is the quotient morphism by an action of $\mathcal{G}_{\mathbb{\Phi}_Y}(1)$, we find that $\pi \times \Phi$ is a quotient morphism by an action of $\mathcal{G}_{\mathbb{\Phi}_Y}(1)$, where $(v, h)(g_1, g_2) := (vg_1, h(g_1, g_2))$. Thus, to show that $r$ descends, it will suffice to check that $\pi(r(v,h)) = \pi(r((v,h)(g_1,g_2)))$ for all $(g_1,g_2) \in \mathcal{G}_{\Phi_Y}(1)$. In fact, since the action of $\mathcal{G}_{\Phi_Y}$ on $\mathbf{V}E$ is the restriction of the action $r$ of $\mathcal{G}_{\mathbb{P}\mathbb{D}_Y}$, we compute:

$$\begin{aligned} r((v,h)(g_1,g_2)) &= r(vg_1, h(g_1,g_2)) \\ &= r(vg_1, g_1^{-1}hg_2) \\ &= r(v, hg_2) \\ &= r(v,h)g_2 \end{aligned}$$

Thus, if $g_2 \in \mathcal{G}_{\Phi_Y}$, we see that $\pi((r(v,h)(g_1,g_2))) = \pi(r(v,h))$, as required. □

# 6 The prismatic topos and its cohomology

We finally turn to study of the prismatic topos per se. We shall need to introduce an additional topos, involving only small prisms, which is technically easier to work with, but which does not seem to be adequately functorial. This topos can be viewed as an analog of the small version of the crystalline topos. It is perhaps worth recalling that functoriality of the small crystalline topos is a delicate matter.

In order to focus ideas, we repeat the main definitions of [10], in the special case we are studying here. Throughout this section, we let $S$ be a fixed crystalline prism. If $S = \operatorname{Spf} A$, then $(A,(p))$ is a bounded prism in the sense of [10].

## 6.1 Prismatic sites

Recall from Definition 2.3.1 that an $X/S$-prism (resp. a small $X$-prism) is a pair $(T, z_T)$, where $T$ is a crystalline prism over $S$ and where $z_T \colon \overline{T} \to X$ is an $S$-morphism (resp. an $S$-morphism such that the induced map $\overline{T} \to \overline{X}$ is flat).

**Definition 6.1.1.** *Let $S$ be a crystalline prism and $X/S$ a morphism of $p$-adic formal schemes. Then we define $\Delta(X/S)$ to be the category of $X/S$-prisms, and $\Delta_s(X/S)$ to be the full subcategory of small $X/S$-prisms. We endow these categories with the topology in which the coverings are the $\delta$-morphisms $f \colon (\widetilde{T}, z_{\widetilde{T}}) \to (T, z_T)$ where $f$ is quasi-compact and $p$-completely faithfully flat (7.3.1). We may also sometimes consider $\Delta(X/S)$ and $\Delta_s(X/S)$ with the Zariski or $p$-completely étale topology.*

Typically $X$ will be a scheme over $\overline{S}$. In any case, the category of $X/S$-prisms is the same as the category of $\overline{X}/S$-prisms, but we will find the additional notational flexibility convenient.

To establish the workability of Definition 6.1.1, we shall need to discuss finite inverse limits of prisms; as we shall see, these are better behaved if we restrict our attention to small prisms.

**Proposition 6.1.2.** *Let $S$ be a crystalline prism and $X/\overline{S}$ a morphism of finite type.*

(1) *In the category of $X/S$-prisms, inverse limits over nonempty finite index sets are representable.*

(2) *If $X \to Y$ is a closed $S$-immersion from $X$ into a crystalline prism $Y$, then in the category of $X/S$-prisms over $Y$, inverse limits over finite index sets are representable.*

(3) *Suppose that $X/\overline{S}$ is smooth and that $T'$ and $T''$ are small $X$-prisms. Then their product $T' \times_S T''$ (in the category of $X/S$-prisms) is also small, and its projection mappings to $T'$ and $T''$ are $p$-completely flat.*

*Proof.* To show that a category admits finite nonempty inverse limits, it is enough to show that the product of every pair of objects and the fibered product of every pair of morphisms are representable. With this strategy, statement (1) is straightforward. If $T'$ and $T''$ are $X/S$-prisms, let $T''' := (T' \times_S T'')_{\mathrm{tf}}$ be their product in the category of crystalline prisms over $S$ as described in Remark 1.6.13, and let $X'''$ be the inverse image of the diagonal under the natural map $\overline{T}''' \to \overline{T}' \times_S \overline{T}'' \to X \times_S X$. Then it is straightforward to verify that $\mathbb{\Delta}_{X'''}(T''')$ represents the product of $T'$ and $T''$ in the category of $X$-prisms over $S$. Note, however, that computation of $T'''$ may not be easy, because of the possibility of torsion in the fiber product. The analogous construction works for fibered products. Since the category of $X/S$-prisms over $Y$ is equivalent to the category of morphisms of $X$-prisms $T \to \mathbb{\Delta}_X(Y)$, statement (1) implies statement (2).

Now suppose that $T'$ and $T''$ are small. Since $X \to \overline{S}$, $\overline{T}' \to X$, and $\overline{T}'' \to X$ are flat, the maps $\overline{T}' \to \overline{S}$ and $\overline{T}'' \to \overline{S}$ are flat. Thus $T' \to S$ and $T'' \to S$ are $p$-completely flat, by Proposition 7.3.2, and it follows that $T''' := T' \times_S T''$ is $p$-torsion free (Proposition 7.3.3). Furthermore, the map $\overline{T}''' = \overline{T}' \times_S \overline{T}'' \to X \times_S X$ is flat, and since $X/\overline{S}$ is smooth, the diagonal $X \to X \times_S X$ is a regular immersion. Then $X''' \to \overline{T}'''$ is also a regular immersion, and it follows from Theorem 2.3.8 that $\overline{\mathbb{\Delta}}_{X'''}(T''') \to X'''$ is flat. Since $X''' \to X$ is also flat, we conclude that $\overline{\mathbb{\Delta}}_{X'''}(T''') \to X$ is flat, so $\mathbb{\Delta}_{X'''}(T''')$ forms a small $X$-prism. Since $\overline{\mathbb{\Delta}}_{X'''}(T''') \to X'''$ and $X''' \to \overline{T}'$ are flat, so is the map $\overline{\mathbb{\Delta}}_{X'''}(T''') \to \overline{T}'$, and it follows that $\mathbb{\Delta}_{X'''}(T''') \to T'$ is $p$-completely flat. The case of the projection to $T''$ follows by symmetry. □

**Remark 6.1.3.** If $T' \to T$ and $T'' \to T$ are morphisms of small $X$-prisms, and if $T' \times_T T''$ is $p$-torsion free, then it is also small. For example, this holds if $T' \to T$ or $T'' \to T$ is $p$-completely flat. However I do not know if this is true in general.

**Proposition 6.1.4.** *Let $g \colon Y \to Z$ be a morphism of crystalline prisms and $i \colon X \to \overline{Z}$ a closed subscheme..*

(1) *Suppose that $X \to \overline{Z}$ and $X' := g^{-1}(X) \to \overline{Y}$ are regular immersions. Then the natural map*

$$\mathbb{\Delta}_{X'}(Y) \to \mathbb{\Delta}_X(Z) \times_Z Y$$

*is an isomorphism.*

(2) *Suppose that $i$ factors as a composite of regular immersions $X \to \overline{Y}$ and $\overline{g}\colon \overline{Y} \to \overline{Z}$. Then the natural map*

$$\mathbb{\Delta}_X(Y) \to \mathbb{\Delta}_X(Z) \times_{\mathbb{\Delta}_Y(Z)} Y$$

*is an isomorphism.*

*Proof.* We shall deduce statement (1) from its analog in Proposition 2.2.15 and the explicit construction of prismatic envelopes in Theorem 2.3.8. Recall that $\mathbb{\Delta}_X(Z)$ is a (completed) limit of $p$-completely flat $W$-schemes $Z^{(n)}$, where each $Z^{(n+1)} \to Z^{(n)}$ is obtained as the dilatation of a subscheme $X^{(n)}$. The envelope $\mathbb{\Delta}_{X'}(Y)$ is constructed in the same way, and in fact there is a commutative diagram

$$\begin{array}{ccc} Y^{(n+1)} & \longrightarrow & Z^{(n+1)} \\ \downarrow & & \downarrow \\ Y^{(n)} & \longrightarrow & Z^{(n)}. \end{array}$$

When $n = 0$, the vertical map on the right (resp. left) is the dilatation of $X$ in $Z$ (resp. of $X'$ in $Y$), and the diagram is Cartesian by Proposition 2.2.15. It follows from Statement (2) of Theorem 2.3.8 that each $X^{(n)} \to Z^{(n)}$ and each ${X'}^{(n)} \to Y^{(n)}$ is again a regular immersion, and from the construction of the ideal defining $X^{(n)}$ described in Lemma 2.3.9, that ${X'}^{(n)} = X^{(n)} \times_{Z^{(n)}} Y^{(n)}$. Then Proposition 2.2.15 implies that that the diagram is Cartesian for every $n$. We deduce that $Y^{(n)} \cong Y \times_Z Z^{(n)}$ for every $n$, and statement (1) follows.

For statement (2) we shall need the following analog of Lemma 2.2.13, which is proved in the same way as that result.

**Lemma 6.1.5.** *Let $Y$ be a crystalline prism, let $X \subseteq X' \subseteq \overline{Y}$ be closed immersions, and let $X'' := \pi_Y^{-1}(X) \subseteq \mathbb{\Delta}_{X'}(Y)$. Then the natural map*

$$\mathbb{\Delta}_{X''}(\mathbb{\Delta}_{X'}(Y)) \to \mathbb{\Delta}_X(Y)$$

*is an isomorphism.* □

We prove statement (2) following the method of Proposition 2.2.15, working locally in an affine setting with the same notation, We first find a sequence $(x_1, \ldots, x_r)$ in $C$ which generates the ideal $I_{\overline{Y}/\overline{Z}}$ and is $C/pC$-regular. Since $X \to \overline{Y}$ is also a regular immersion, we may then find a $B/pB$-regular sequence $(y_1, \ldots, y_m)$ in $C$ whose image generates the ideal $I_{X/\overline{Y}}$ of $X$ in $\overline{Y}$. Then by Theorem 2.3.13, the prismatic envelope $\mathbb{\Delta}_Y(Z)$ of $Y$ in $Z$ corresponds to the PD-envelope $\mathbb{PD}_{\tilde{Y}}(\mathbb{D}_Y(Z))$ of the section $\tilde{Y}$ of $\mathbb{D}_Y(Z)$ defined by $Y \to Z$. The ideal of this section is generated by the regular sequence $(t_1, \ldots, t_r)$, where $t_i := \rho(x_i)$, by Proposition 2.2.14, and hence this PD-envelope is given by the formal spectrum of the completion $\tilde{C}$ of the PD-polynomial algebra $B\langle t_1, \ldots, t_r\rangle$. The ideal of the inverse image of $X$ in $\mathbb{D}_Y(Z)$ is generated by the sequence $(p, x_1, \ldots, x_r, y_1, \ldots, y_m)$, or in fact just by $(p, y_1, \ldots, y_m)$ since $x_i \in p\tilde{C}$. Since $\tilde{C}/p\tilde{C}$ is a PD-polynomial algebra over $B/pB$, it follows that $(y_1, \ldots, y_r)$ is also $\tilde{C}/p\tilde{C}$-regular. Thus the maps $X' := \pi_Z^{-1}(X) \to \overline{\mathbb{\Delta}}_Y(Z)$ and $X' \cap Y \to Y$ are regular immersions. Then statement (1), applied to the morphisms $j'\colon Y \to \mathbb{\Delta}_Y(Z)$ and $\tilde{i}\colon X' \to \mathbb{\Delta}_Y(Z)$, implies that the map

$$\mathbb{\Delta}_X(Y) = \mathbb{\Delta}_{X'\cap\tilde{Y}}(\tilde{Y}) \to \mathbb{\Delta}_{X'}(\mathbb{\Delta}_Y(Z)) \times_{\mathbb{\Delta}_Y(Z)} j'(Y).$$

is an isomorphism. Lemma 6.1.5 tells us that $\mathbb{\Delta}_{X'}(\mathbb{\Delta}_Y(Z)) \cong \mathbb{\Delta}_X(Z)$, proving statement (2). □

**Corollary 6.1.6.** *If $Y$ is a crystalline prism and $X$ and $X'$ are regularly immersed in $\overline{Y}$ and meet transversally, then the natural map*

$$\mathbb{\Delta}_{X\cap X'}(Y) \to \mathbb{\Delta}_X(Y) \times_Y \mathbb{\Delta}_{X'}(Y)$$

*is an isomorphism.*

*Proof.* This is proved in the same way as Corollary 2.2.16 □

The following result explains the relation between the two prismatic sites and establishes functoriality of the prismatic topos. Note that we have not established functoriality of the small prismatic topos.

**Proposition 6.1.7.** *Let $X/\overline{S}$ be a smooth morphism of schemes. Then $\mathbb{\Delta}(X/S)$ and $\mathbb{\Delta}_s(X/S)$, endowed with any of the topologies above, form sites. Furthermore, the inclusion functor $u\colon \mathbb{\Delta}_s(X/S) \to \mathbb{\Delta}(X/S)$ is continuous and cocontinuous and hence induces a morphism of topoi:*

$$u\colon (X/S)_{\mathbb{\Delta}_s} \to (X/S)_{\mathbb{\Delta}}.$$

*A morphism $f\colon \overline{X} \to \overline{X}'$ of smooth schemes over $\overline{S}$ induces a continuous and cocontinuous map $\mathbb{\Delta}(f)\colon \mathbb{\Delta}(X/S) \to \mathbb{\Delta}(X'/S)$ and hence a morphism of topoi*

$$f_{\mathbb{\Delta}}\colon (X/S)_{\mathbb{\Delta}} \to (X'/S)_{\mathbb{\Delta}}.$$

*Proof.* It is verified in [10, 3.12] that $\mathbb{\Delta}(X/S)$ does in fact form a site with the $p$-completely flat topology. The Zariski topology and étale topologies are even easier. Let us check that the same is true for $\mathbb{\Delta}_s(X/S)$. It is clear that isomorphisms are coverings and that the composition of coverings is a covering (see Proposition 7.3.3). Furthermore, if $T' \to T$ is a morphism and $\widetilde{T} \to T$ is a covering in $\mathbb{\Delta}_s(X/S)$, then Proposition 7.3.3 implies that the fiber product $T' \times_T \widetilde{T}$ is $p$-torsion free. Furthermore, the map $\widetilde{T}_1 \to T_1$ is flat, hence so is the map $T'_1 \times_T \widetilde{T}_1 \to T'_1$, and since $T'_1 \to X$ is flat, so is the map $(T' \times_T \widetilde{T})_1 \to X$. Thus $T' \times_T \widetilde{T}$ is again small.

To prove that $u$ is cocontinuous, let $\widetilde{T} \to T$ be a cover in $\mathbb{\Delta}_s(X/S)$ and let $\widetilde{T}' \to u(\widetilde{T})$ be a cover in $\mathbb{\Delta}(X/S)$. Then $\widetilde{T}' \to u(\widetilde{T})$ is $p$-completely flat, hence $\widetilde{T}'_1 \to T_1$ is flat, hence $\widetilde{T}'_1 \to X$ is flat, so $T'$ is again small, hence in the image of $u$. To check that $u$ is continuous, observe that if $\widetilde{T} \to T$ is a covering in $\mathbb{\Delta}_s(X/S)$, then $u(\widetilde{T}) \to u(T)$ is a covering in $\mathbb{\Delta}(X/S)$, and if $T' \to T$ is a morphism in $\mathbb{\Delta}_s(X/S)$, then the fiber product $\widetilde{T} \times_T T'$ is again small, hence the map $u(\widetilde{T} \times_T T') \to u(\widetilde{T}) \times_{u(T)} u(T')$ is an isomorphism.

It is also checked in [10, 4.3] that a morphism $f \colon X \to X'$ of $S$-schemes induces a morphism of topoi. The argument there is rather abstract; here is another. If $(T, z_T)$ is an $X$-prism, then

$$\mathbb{\Delta}(f)(T, z_T) := (T, f \circ z_T)$$

is an $X'$-prism. Thus we find a functor

$$\mathbb{\Delta}(f) \colon \mathbb{\Delta}(X/S) \to \mathbb{\Delta}(X'/S) : (T, z_T) \mapsto (T, f \circ z_T).$$

To check that this functor is cocontinuous, let $g \colon (\widetilde{T}, z_{\widetilde{T}}) \to (T, f \circ z_T)$ be a covering of $X'$-prisms. The diagram

$$\begin{array}{ccc} T_1 & \xleftarrow{\;g_1\;} & \widetilde{T}_1 \\ {\scriptstyle z_T}\downarrow & & \downarrow{\scriptstyle z_{\widetilde{T}}} \\ X & \xrightarrow{\;f\;} & X' \end{array}$$

commutes. Then $(\widetilde{T}, z_T \circ g_1)$ is an $X$-prism, and the map $g$ defines a covering of $X$-prisms $(\widetilde{T}, z_T \circ g_1) \to (T, z_T)$ with $\mathbb{\Delta}(f)(\widetilde{T}, z_{\widetilde{T}} \circ g_1) = (\widetilde{T}, z_{\widetilde{T}})$. To check that $\mathbb{\Delta}(f)$ is continuous, suppose instead that $(\widetilde{T}, z_{\widetilde{T}}) \to (T, z_T)$ is a covering of $X$-prisms. Then $(\widetilde{T}, f \circ z_{\widetilde{T}}) \to (T, f \circ z_T)$ is a covering of $X'$-prisms. Furthermore, if $T' \to T$ is a morphism in $\mathbb{\Delta}(X/S)$, then from the description of fiber products we saw in Proposition 6.1.2, it is clear that the map

$$\mathbb{\Delta}(f)\left((\widetilde{T}, z_{\widetilde{T}}) \times_{(T, z_T)} (T', z_{T'})\right) \to (\widetilde{T}, f \circ z_{\widetilde{T}}) \times_{(T, f \circ z_T)} (T', f \circ z_{T'})$$

is an isomorphism. □

The following result illustrates an important advantage of working with small prisms. Before stating it, we introduce a useful, if somewhat abusive, abbreviated notation. Suppose that $f \colon X \to X'$ is a closed immersion of $S$-schemes and $(T', z_{T'})$ is an $X'$-prism. Then, $z_{T'}^{-1}(X) \to T'$ is also a closed immersion, and we can form the prismatic envelope $\mathbb{\Delta}_{z_{T'}^{-1}(X)}(T')$. The map $\overline{\mathbb{\Delta}}_{z_{T'}^{-1}(X)}(T') \to z_{T'}^{-1}(X) \to X$ endows $\mathbb{\Delta}_X(T')$ with the structure of an $X$-prism, and if no confusion is likely, we write

$$\mathbb{\Delta}_X(T') := f_{\mathbb{\Delta}}^{-1}(T') := \mathbb{\Delta}_{z_{T'}^{-1}(X)}(T'). \tag{6.1.1}$$

**Proposition 6.1.8.** *Suppose that $Y/S$ is a $p$-completely smooth morphism of crystalline prisms and $i\colon X \to \overline{Y}$ is a regular closed immersion.*

*(1) If $T$ is a small $Y$-prism, the prismatic envelope $\mathbb{\Delta}_X(T)$ of $z_T^{-1}(X)$ in $T$ is a small $X$-prism. Furthermore, if $T' \to T$ is a morphism in $\mathbb{\Delta}_s(Y/S)$, the natural map*

$$\mathbb{\Delta}_X(T') \to T' \times_T \mathbb{\Delta}_X(T)$$

*is an isomorphism.*

*(2) The prismatic envelope $\mathbb{\Delta}_X(Y(1))$ of the composition $X \to Y \to Y(1)$ of $i$ with the diagonal embedding in $Y(1)$ is small. Furthermore, if $T \to Y$ is any morphism in $\mathbb{\Delta}_s(Y/S)$, the natural map*

$$\mathbb{\Delta}_X(T \times_S Y) \to T \times_Y \mathbb{\Delta}_X(Y(1))$$

*is an isomorphism.*

*Proof.* If $T$ is a small $Y$-prism, the map $z_T\colon \overline{T} \to \overline{Y}$ is flat, and since $X \to \overline{Y}$ is a regular immersion, the same is true of the map $z_T^{-1}(X) \to \overline{T}$. Then it follows from statement (2) of Theorem 2.3.8 that the map $\overline{\mathbb{\Delta}}_X(T) \to X$ is flat, so $\mathbb{\Delta}_X(T)$ is a small $X$-prism. Suppose $T' \to T$ is a morphism in $\mathbb{\Delta}_s(Y/S)$. Working locally, choose an $\mathcal{O}_{\overline{Y}}$-regular sequence which generates the ideal of $X$ in $\overline{Y}$. Since $z_T$ and $z_{T'}$ are flat, this sequence remains regular in $\overline{T}$ and $\overline{T}'$. Then statement (1) of Proposition 6.1.4 implies that the map $\mathbb{\Delta}_X(T') \to T' \times_T \mathbb{\Delta}_X(T)$ is an isomorphism.

Since $\overline{Y}/\overline{S}$ is smooth, the diagonal $\overline{Y} \to \overline{Y} \times_S \overline{Y}$ is a regular immersion, and hence so is the the map $X \to \overline{Y}(1)$. Then statement (1) implies that $\mathbb{\Delta}_X(Y(1))$ is small. Furthermore, if $T \to Y$ is a morphism of small $Y/S$-prisms, the map $\overline{T} \to \overline{Y}$ is flat, hence $\overline{T} \times_S \overline{Y} \to \overline{Y} \times_S \overline{Y}$ is also flat, and hence $T \times_S Y \to Y(1)$ is $p$-completely flat. Thus it follows from Theorem 2.3.8 that the natural map

$$\mathbb{\Delta}_X(T \times_S Y) \to (T \times_S Y) \times_{Y(1)} \mathbb{\Delta}_X(Y(1)) \cong T \times_Y \mathbb{\Delta}_X(Y(1))$$

is an isomorphism. □

Note that Proposition 6.1.8 is not true if we work with full prismatic sites. For example, let $S = \operatorname{Spf} W$, let $X = \operatorname{Spec} k$, and let $Y = \operatorname{Spf}(W[x]\hat{})$, with $\phi(x) := x^p$ and the inclusion $X \to Y$ defined by $(p, x)$. Then the prismatic envelope $\mathbb{\Delta}_X(Y) = \operatorname{Spf}(W\langle x/p\rangle\hat{})$ of $X$ in $Y$ can also be viewed as a $Y$-prism $T'$ via the (not flat) map $\overline{\mathbb{\Delta}}_X(Y) \to X \to \overline{Y}$, and the natural map $f\colon T' \to Y$ is a morphism in $\mathbb{\Delta}(Y/S)$. However, the map

$$\mathbb{\Delta}_X(T') \to T' \times_Y \mathbb{\Delta}_X(Y) \tag{6.1.2}$$

is not an isomorphism. Indeed, $\mathbb{\Delta}_X(T') \cong \mathbb{\Delta}_X(Y) = \operatorname{Spf}(W\langle x/p\rangle\hat{})$, but

$$T' \times_Y \mathbb{\Delta}_X(Y) \cong \operatorname{Spf}\left(W\langle x/p\rangle\hat{}\,\hat{\otimes}_{W[x]\hat{}}\, W\langle x/p\rangle\hat{}\right),$$

which is not isomorphic to $\mathrm{Spf}(W\langle x/p\rangle\hat{})$. (For example, $1 \otimes x/p - x/p \otimes 1$ is a nonzero $p$-torsion element in the tensor product.) It could be argued that (6.1.2) is an isomorphism if the fiber product appearing is taken in the category of crystalline prisms but in general this would involve killing an incomputable amount of $p$-torsion.

**Proposition 6.1.9.** *Let $i\colon X \to X'$ be a closed immersion of smooth schemes over $\overline{S}$. If $T' \in \Delta(X'/S)$, the sheaf $i^{-1}(T')$ in $(X/S)_{\Delta}$ is represented by $\Delta_X(T')$, which is small if $T'$ is small.*

*Proof.* For any morphism $f\colon X \to X'$, the functor $f^{-1}(T')$ is by definition the sheaf on $\Delta(X/S)$ sending an object $T \in \Delta(X/S)$ to the set of $\phi$-morphisms $v\colon T \to T'$ such that $z_{T'} \circ \overline{v} = f \circ z_T$, *i.e.*, such that $(z_T, \overline{v})$ factors through $X \times_{X'} \overline{T}'$. If $f = i$ is a closed immersion, then this functor is represented by $\Delta_X(T')$. Furthermore, since $X$ and $X'$ are smooth, $i$ is a regular immersion, and if $\overline{T}' \to X'$ is flat, the immersion $X \times_{X'} \overline{T}' \to \overline{T}'$ is also regular. Then Theorem 2.3.8 tells us that $\Delta_X(T')$ is also small. ☐

Let $S$ be a crystalline prism and $X/\overline{S}$ a smooth morphism. We shall describe two constructions of coverings of the final object of the topos $(X/S)_{\Delta}$. The first of these, based on [10, 4.16, 4.17] and a conversation with A. Shiho, is in fact a covering in the Zariski topology. However, it seems to be quite unwieldy in practice.

**Proposition 6.1.10.** *Let $S$ be a crystalline prism, let $Y \to S$ be a $p$-completely smooth morphism of $p$-adic formal schemes, and let $r\colon \mathrm{J}_{S,\mathrm{tf}}(Y) \to Y$ be the universal morphism from a crystalline $S$-prism to $Y$ described in Proposition 1.6.11. If $X \to \overline{Y}$ is a closed immersion, let $\tilde{X} := r^{-1}(X)$, and let $(\Delta_{\tilde{X}}(\mathrm{J}_{S,\mathrm{tf}}(Y), z_{\tilde{X}})$ be the prismatic envelope of $\tilde{X}$ in $\mathrm{J}_{S,\mathrm{tf}}(Y)$. Then*

$$(\Delta_X(\mathrm{J}_{S,\mathrm{tf}}(Y)), r_{|\tilde{X}} \circ z_{\tilde{X}})$$

*is a small $X$-prism and is a covering of the final objects of $(X/S)_{\Delta}$ and $(X/S)_{\Delta_s}$ respectively, if these are endowed with the Zariski topology.*

*Proof.* We first claim that $r\colon \mathrm{J}_{S,\mathrm{tf}}(Y) \to Y$ is $p$-completely flat. This can be checked étale locally on $Y$, so we may assume that $Y$ and $S$ are affine and that $Y$ is $p$-completely étale over an affine space over $S$. It follows from the explicit construction of the functor $\mathrm{J}_{S,\mathrm{tf}}$ in Example 1.6.4, and the fact that it commutes with $p$-completely étale localization that the space $\mathrm{J}_{S,\mathrm{tf}}(Y)$ is $p$-completely étale over a (possibly infinite dimensional) affine space over $Y$, hence $p$-completely flat over $Y$. Since $X \to \overline{Y}$ is a closed immersion of smooth $\overline{S}$-schemes, it is a regular immersion, and since $\overline{\mathrm{J}}_{S,\mathrm{tf}}(Y) \to \overline{Y}$ is flat, $\tilde{X} \to \overline{\mathrm{J}}_{S,\mathrm{tf}}(Y)$ is also a regular immersion. Then it follows from Theorem 2.3.8 that $\overline{\Delta}_X(\mathrm{J}_{S,\mathrm{tf}}(Y)) \to \tilde{X}$ is flat. Since $r_{|\tilde{X}}$ is also flat, we can conclude that $\Delta_X(\mathrm{J}_{S,\mathrm{tf}}(Y))$ is small.

We claim that every affine $X$-prism $T$ admits a morphism to $\Delta_X(\mathrm{J}_{S,\mathrm{tf}}(Y))$. Since $T$ is $p$-adically complete and affine and $Y/S$ is formally smooth, the map $\overline{T} \to X \to Y$ extends to an $S$-morphism $T \to Y$, which in turn lifts to a unique morphism of crystalline prisms $T \to \mathrm{J}_{S,\mathrm{tf}}(Y)$. This morphism necessarily maps $\overline{T}$ to $r^{-1}(X)$ and hence factors through $\mathbb{\Delta}_X(\mathrm{J}_{S,\mathrm{tf}}(Y))$. □

We should explain our characterization of $\mathbb{\Delta}_X(\mathrm{J}_{S,\mathrm{tf}}(Y))$ as "unwieldy." This is partly because $\mathrm{J}_{S,\mathrm{tf}}(Y)$ is itself somewhat unwieldy, but there is more to the story. For example, suppose that $X = \operatorname{Spec} k[x]$ and $Y = \operatorname{Spf}(W[x]^\wedge)$. Then $\mathrm{J}_{S,\mathrm{tf}}(Y) = \operatorname{Spf} W[x_0, x_1, \ldots]^\wedge$ with $\phi(x_i) = x_i^p + px_{i+1}$, as we saw in Example 1.6.4. Here $\tilde{X} := r^{-1}(X) = \overline{\mathrm{J}}_{S,\mathrm{tf}}(Y)$, which does not seem so unwieldy. However to do cohomology calculations, one needs to understand $\mathbb{\Delta}_X(\mathrm{J}_{S,\mathrm{tf}}(Y)(1))$, the prismatic envelope of $\tilde{X} \times_X \tilde{X}$ in $\mathrm{J}_{S,\mathrm{tf}}(Y) \times \mathrm{J}_{S,\mathrm{tf}}(Y)$, which seems very difficult to describe explicitly.

Fortunately, it turns out that, if one is willing to use the $p$-completely flat topology, the construction of coverings of the final object becomes much simpler. The following proposition was inspired by a result [28, 3.4] of Morrow and Tsuji, which it generalizes.

**Proposition 6.1.11.** *Let $Y/S$ be a $p$-completely smooth morphism of crystalline prisms and let $X \to \overline{Y}$ be a closed immersion, where $X/\overline{S}$ is smooth. Then the prismatic envelope $\mathbb{\Delta}_X(Y)$ of $X$ in $Y$ is small and is a covering of the final objects of $(X/S)_{\mathbb{\Delta}}$ and $(X/S)_{\mathbb{\Delta}_s}$ respectively, if these are endowed with the $p$-completely flat topology.*

*Proof.* Since $X \to \overline{Y}$ is a regular immersion, statement (2) of Theorem 2.3.8 shows that $\mathbb{\Delta}_X(Y)$ is small. To show that it covers the final object of $(X/S)_{\mathbb{\Delta}}$, we shall show that if $T$ is an affine $X$-prism, then there exists a $p$-completely faithfully flat morphism of crystalline prisms $\widetilde{T} \to T$ and a morphism of $X$-prisms $\widetilde{T} \to \mathbb{\Delta}_X(Y)$.

As in the previous proposition, the formal smoothness of the underlying morphism $Y \to S$ implies that the map $f_1 \colon T_1 \to X \to Y$ lifts to a morphism $T \to Y$. This map may not be compatible with the Frobenius liftings, but Proposition 1.7.11 shows that we may find a $p$-completely faithfully flat morphism of crystalline prisms $u \colon \widetilde{T} \to T$ and a morphism of crystalline prisms $\tilde{f} \colon \widetilde{T} \to Y$ such that $\tilde{f}_1 = u_1 \circ f_1$. Then $\widetilde{T}$ becomes an $X$-prism over $Y$ and hence $\tilde{f}$ factors through $\mathbb{\Delta}_X(Y)$. □

## 6.2 Prismatic crystals

Continuing to follow [10], we endow $(X/S)_{\mathbb{\Delta}}$ with the structure of a ringed topos, with structure sheaf $\mathcal{O}_{X/S}$ given by $T \mapsto \mathcal{O}_T$. This sheaf is $p$-completely quasi-coherent, in the sense of Definition 7.4.1.

**Definition 6.2.1.** *If $S$ is a crystalline prism and $X/S$ is a scheme, a crystal of $\mathcal{O}_{X/S}$-modules on $\mathbb{\Delta}(X/S)$ (resp. $\mathbb{\Delta}_s(X/S)$) is a $p$-completely quasi-coherent sheaf $E$ of $\mathcal{O}_{X/S}$-modules such that for each morphism $f \colon T' \to T$ in $\mathbb{\Delta}(X/S)$ (resp. $\mathbb{\Delta}_s(X/S)$), the corresponding transition map $f^*(E_T) \to E_{T'}$ is an isomorphism, where $f^*(E_T)$ means the $p$-completed pullback.*

A crystal $E$ of $\mathcal{O}_{X/S}$-modules can also be viewed geometrically. If $T$ is an object of $\Delta(X/S)$, then $E_T$ is a $p$-completely quasi-coherent sheaf of $\mathcal{O}_T$-modules, and we let $\mathbf{V}E_T$ denote the corresponding formal scheme, as described in Remark 7.4.3. Then $T \mapsto \mathbf{V}E_T$ is a crystal of formal schemes. More precisely, let $\mathbf{V}\Delta(X/S)$ be the category of morphisms of $p$-adic formal schemes $V \to T$, where $T$ is an object of $\Delta(X/S)$, and let $\mathbf{F}\colon \mathbf{V}\Delta(X/S) \to \Delta(X/S)$ be the functor taking a morphism to its target. Then $\mathbf{V}\Delta(X/S)$ is fibered over $\Delta(X/S)$, and the assignment $T \mapsto (\mathbf{V}E_T \to T)$ defines a crystal in $\mathbf{V}\Delta(X/S)$ in the sense of Definition 1.4.13.

**Proposition 6.2.2.** *Let $Y/S$ be a $p$-completely smooth morphism of crystalline prisms and $X \to \overline{Y}$ a closed immersion, where $X/\overline{S}$ is smooth.*

*(1) The morphism $u\colon \Delta_s(X/S) \to \Delta(X/S)$ induces an equivalence between the corresponding categories of crystals of $\mathcal{O}_{X/S}$-modules.*

*(2) The functor*

$$u^* \circ i_{\Delta *}\colon (X/S)_{\Delta} \to (Y/S)_{\Delta} \to (Y/S)_{\Delta_s}$$

*takes crystals of $\mathcal{O}_{X/S}$-modules on $\Delta(X/S)$ to crystals of $\mathcal{O}_{Y/S}$-modules on $\Delta_s(Y/S)$. In particular, $\mathcal{A}_{X/Y/S} := u^* i_{\Delta *}(\mathcal{O}_{\Delta})$ is a crystal of $\mathcal{O}_{Y/S}$-algebras on $\Delta_s(Y/S)$. If $E$ is a crystal of $\mathcal{O}_{X/S}$-modules on $\Delta(X/S)$, then $u^* i_{\Delta *}(E)$ inherits the structure of a $p$-completely quasi-coherent $\mathcal{A}_{X/Y/S}$-module.*

*(3) The functor $E \mapsto u^* i_{\Delta *}(E)$ is an equivalence from the category of crystals of $\mathcal{O}_{X/S}$-modules on $\Delta(X/S)$ to the category of crystals of $p$-completely quasi-coherent crystals of $\mathcal{A}_{X/Y/S}$-modules on $\Delta_s(Y/S)$.*

*Proof.* The general formalism of crystals, reviewed in Proposition 1.4.14, tells us that if $T$ is a covering of the final object, the category of crystals of $\mathcal{O}_{X/S}$-modules is equivalent to the category of $p$-completely quasi-coherent $\mathcal{O}_T$-modules endowed with a stratification; this holds both on $\Delta(X/S)$ and on $\Delta_s(X/S)$. Proposition 6.1.11 tells us that $\Delta_X(Y)$ is such a covering, and Proposition 6.1.8 tells us that $\Delta_X(Y)$ and $\Delta_X(Y(1))$ are small. Statement (1) follows.

Recall from Proposition 6.1.9 that if $T$ is a $Y$-prism, the sheaf $i_{\Delta}^{-1}(T)$ is represented by the prismatic envelope $(\Delta_X(T), z_{\Delta}, \pi_T)$ of $z_T^{-1}(X)$ in $T$. If $E$ is a sheaf in $(X/S)_{\Delta}$, it follows that

$$(i_{\Delta *}(E))_T = \pi_{T*}(E_{\Delta_X(T)}). \tag{6.2.1}$$

The morphism $\pi_T$ is affine, and we can conclude that, if $E$ is a sheaf of $p$-completely quasi-coherent $\mathcal{O}_{X/S}$-modules, then $i_{\Delta *}(E)$ becomes a $p$-completely quasi-coherent sheaf of $\mathcal{O}_{Y/S}$-modules. Moreover, $(\mathcal{A}_{X/Y/S})_T = \pi_{T*}(\mathcal{O}_{\Delta_X(T)})$, and thus $(i_{\Delta *}(E))_T$ inherits the structure of a $p$-completely quasi-coherent $\mathcal{A}_{X/Y/S}$-module.

Suppose that $E$ is a crystal of $\mathcal{O}_{X/S}$-modules in $(X/S)_{\Delta}$ and $f\colon T' \to T$ is a morphism in $\Delta(Y/S)$. This map induces a morphism $\Delta(f)\colon \Delta_X(T') \to \Delta_X(T)$, fitting into a commutative diagram:

$$\begin{array}{ccc} \Delta_X(T') & \xrightarrow{\Delta(f)} & \Delta_X(T) \\ \downarrow{\scriptstyle \pi_{T'}} & & \downarrow{\scriptstyle \pi_T} \\ T' & \xrightarrow{f} & T. \end{array}$$

Proposition 6.1.8 tells us that this diagram is Cartesian if $T'$ and $T$ are small. Since $\pi_T$ and $\pi_{T'}$ are affine, this implies that the map

$$f^* \pi_{T*}(E_{\Delta_X(T)}) \to \pi_{T'*} \Delta(f)^*(E_{\Delta_X(T)})$$

is an isomorphism. Since $E$ is a crystal on $\Delta(X/S)$, the map

$$\Delta(f)^*(E_{\Delta_X(T)}) \to E_{\Delta_X(T')}$$

is also an isomorphism. As we have seen, $i_{\Delta *}(E)_T = \pi_{T*}(E_{\Delta_X(T)})$, and similarly for $T'$. Thus $u^* i_{\Delta *}(E)$ forms a crystal of $\mathcal{O}_{Y/S}$-modules on $\Delta_s(Y/S)$. This proves statement (2).

To prove statement (3), suppose that $E$ is a crystal of $p$-completely quasi-coherent $\mathcal{A}_{X/Y/S}$-modules on $\Delta_s(Y/S)$. Then $E_Y$ and $(\mathcal{A}_{X/Y/S})_Y$ are sheaves of $\mathcal{O}_Y$-modules endowed with prismatic stratifications; in fact $E_Y$ has a structure of a $p$-completely quasi-coherent $(\mathcal{A}_{X/Y/S})_Y$-module, compatible with the given stratifications. But $(\mathcal{A}_{X/Y/S})_Y = \pi_{T*}(\mathcal{O}_{\Delta_X(Y)})$, and so $E_Y$ can be viewed as a sheaf of $\mathcal{O}_{\Delta_X(Y)}$-modules on $\Delta_X(Y)$, and is endowed with prismatic stratification. As we recalled in the proof of statement (1), this gives rise to a crystal of $\mathcal{O}_{X/S}$-modules on $\Delta(X/S)$. It is clear that this construction gives a quasi-inverse to the functor $u^* i_{\Delta *}$. □

Note that Proposition 6.2.2 is not true if we work with full prismatic sites, as the discussion after Proposition 6.1.8 shows.

The construction of coverings in Proposition 6.1.11 enables us to establish the speculated relationship between prismatic crystals and $p$-connections sketched in the introduction.

**Theorem 6.2.3.** *With the hypothesis of Proposition 6.2.2, there is a natural equivalence, made explicit below, from the category* MICP$(X/Y/S)$ *given in Definition 4.2.5, to the category of crystals of $\mathcal{O}_{X/S}$-modules on $\Delta(X/S)$.*

*Proof.* Recall that an object $(E, \nabla')$ of MICP$(X/Y/S)$ is a $p$-completely quasi-coherent sheaf of $\pi_{Y*}(\mathcal{O}_{\Delta_X(Y)})$-modules $E$ together with an integrable quasi-nilpotent $p$-connection $\nabla'\colon E' \to E' \otimes \Omega^1_{Y/S}$ which is compatible with the canonical $p$-connection on $\pi_{Y*}(\mathcal{O}_{\Delta_X(Y)})$. As we saw in Theorem 4.4.3, the connection $\nabla'$ induces a prismatic stratification on $E'$, which is compatible with the prismatic stratification on $\pi_{Y*}(\mathcal{O}_{\mathbb{A}_X(Y)})$. These data lead to a crystal of $p$-completely quasi-coherent $\mathcal{A}_{X/Y/S}$-modules on $Y/S$, and so the theorem follows from statement (3) of Proposition 6.2.2. □

The following special case is worth stating.

**Corollary 6.2.4.** *If $Y/S$ is a $p$-completely smooth morphism of crystalline prisms, the category of $\mathcal{O}_{Y/S}$-modules on $\mathbb{A}(Y/S)$ is equivalent to the category of $p$-completely quasi-coherent $\mathcal{O}_Y$-modules with integrable quasi-nilpotent $p$-connection.* □

## 6.3 Cohomology of the prismatic topos

We are now prepared to show that the $p$-de Rham and Higgs complexes we constructed in Proposition 4.2.1 and Theorem 3.3.3. do indeed calculate prismatic cohomology in the sense of [10]. Let $S$ be a crystalline prism and $X/\overline{S}$ a smooth morphism. A $p$-completely quasi-coherent sheaf of $\mathcal{O}_{X/S}$-modules on $\mathbb{A}(X/S)$ assigns to each $X$-prism $T$ a $p$-completely quasi-coherent sheaf of $\mathcal{O}_T$-modules $E_T$. In particular, each of these is a "$p$-adic sheaf" in the sense of Definition 7.1.2, and we shall identify it with the inverse system $\{E_{T,n} : n \in \mathbf{N}\}$ of its reductions modulo powers of $p$. As explained in §7.2, the category of $p$-adic sheaves is not abelian, but is an exact subcategory of the abelian category of inverse systems of $p$-torsion sheaves. In particular, in forming the derived category $D^+_{\mathbb{A}}(X/S)$ of sheaves of $\mathcal{O}_{X/S}$-modules on $\mathbb{A}(X/S)$, we localize by morphisms of complexes which are strict quasi-isomorphisms, *i.e.*, are quasi-isomorphisms modulo every power of $p$; (see Definition 7.2.2). As explained in [10, 4.4], we have a morphism of topoi

$$v_{X/S} \colon (X/S)_{\mathbb{A}} \to X_{\acute{e}t},$$

and hence derived functors

$$Rv_{X/S*} \colon D^+_{\mathbb{A}}(X/S) \to D^+(X_{\acute{e}t}),$$

where the latter is the derived category of inverse systems of abelian sheaves on $X_{\acute{e}t}$.

If $T$ is an $X/S$-prism, let $\mathbb{A}(X/S)_{|T}$ denote the localized site whose objects are morphisms $T' \to T$ in $\mathbb{A}(X/S)$. There is a morphism $s \colon \mathbb{A}(X/S)_{|T} \to \mathbb{A}(X/S)$ taking an object $T' \to T$ of $\mathbb{A}(X/S)_{|T}$ to the object $T'$ of $\mathbb{A}(X/S)$, and as explained in [7, 5.23], a corresponding morphism of topoi

$$j_T \colon ((X/S)_{\mathbb{A}})_{|T} \to (X/S)_{\mathbb{A}},$$

where $j_T^{-1}$ takes the (sheaf represented by) $T'$ to (the sheaf represented by)

$$\mathbb{A}_X(T' \times_S T) := \mathbb{A}_{X'}(T' \times_S T),$$

where $X' := (z_{T'} \times z_T)^{-1}(X)$. Note that if $T$ is small, then $T/S$ is $p$-completely flat, hence $T' \times_S T$ is $p$-torsion free, and the product is computed in the category of formal schemes.

We shall need the following analog of [7, 5.26]. We state it for $\mathbb{\Delta}(X/S)$, but the same arguments show it also holds for $\mathbb{\Delta}_s(X/S)$.

**Proposition 6.3.1.** *If $T$ is an object of $\mathbb{\Delta}(X/S)$, let $T_{\mathrm{pcf}}$ denote the topos of sheaves on the site defined by the $p$-completely flat coverings of $T$. Then there is a 2-commutative diagram of morphisms of topoi:*

$$\begin{CD} ((X/S)_{\mathbb{\Delta}})_{|T} @>{\psi}>> T_{\mathrm{pcf}} \\ @V{j_T}VV @VV{\lambda}V \\ (X/S)_{\mathbb{\Delta}} @>{v_{X/S}}>> X_{\acute{e}t}. \end{CD}$$

*Furthermore, the following statements are verified.*

*(1) The functor $\psi_*$ is exact. If $E_T$ is a sheaf on $T_{\mathrm{pcf}}$, then*

$$v_{X/S*} j_{T*} \psi^*(E_T) \cong \lambda_*(E_T).$$

*(2) If $E_T$ is a $p$-completely quasi-coherent sheaf of $\mathcal{O}_T$-modules on $T_{\mathrm{pcf}}$, then $\psi^*(E_T)$ is a crystal of $\mathcal{O}_{X/S}$-modules on $((X/S)_{\mathbb{\Delta}})_{|T}$.*

*(3) If $z_T \colon \overline{T} \to X$ is affine and $E$ is a crystal of $\mathcal{O}_{X/S}$-modules on $((X/S)_{\mathbb{\Delta}})_{|T}$, then $R^q j_{T*}E$ and $R^q v_{X/S*}(j_{T*}E)$ vanish for $q > 0$.*

*Proof.* If $E$ is a sheaf in $((Y/S)_{\mathbb{\Delta}})_{|T}$, then $\psi_*(E) := E_{(T,\mathrm{id})}$, and if $E_T$ is a sheaf on $T$ and $(T', h') \in \mathbb{\Delta}(X/S)_{|T}$, then $\psi^{-1}(E)_{(T',h)} = h'^{-1}(E_T)$. It is clear that $\psi^{-1}$ is left adjoint to $\psi_*$, that $\psi^{-1}$ and $\psi_*$ are exact, and that $\psi^*(E_T)$ is a crystal if $E$ is $p$-completely quasi-coherent. The morphism $\lambda$ is defined as the composition

$$\lambda \colon T_{\mathrm{pcf}} \longrightarrow T_{\acute{e}t} \xrightarrow{z_T} X_{\acute{e}t}.$$

The 2-commutativity of the diagram is easy to check from the definitions. Then if $E_T$ is a sheaf on $T_{\mathrm{pcf}}$,

$$v_{X/S*} j_{T*} \psi^*(E_T) \cong \lambda_* \psi_* \psi^*(E_T) = \lambda_*(\psi^*(E)_{(T,\mathrm{id})}) = \lambda_*(E_T).$$

This completes the proof of statements (1) and (2).

To prove (3), let $T'$ be an object of $\mathbb{\Delta}(X/S)$, let $T'' := \mathbb{\Delta}_X(T' \times_S T)$, with projections $h \colon T'' \to T$, and $h' \colon T'' \to T'$. Then $(T'', h) \in \mathbb{\Delta}(X/S)_{|T}$, and if $E$ is a sheaf in $\mathbb{\Delta}(X/S)_{|T}$, then $(j_{T*}E)_{T'} = h'_*(E_{(T'',h)})$. Let $X'$ be the fiber product in the following Cartesian diagram:

$$\begin{CD} X' @>>> \overline{T} \\ @VVV @VV{z_T}V \\ \overline{T}' @>{z_{T'}}>> X, \end{CD}$$

Since $z_T : \overline{T} \to X$ is by assumption affine, the map $X' \to \overline{T}'$ is also affine. Since the map $\overline{T}'' \to X'$ is affine by the construction in Theorem 2.3.8, we conclude that the composition $\overline{h}' : \overline{T}'' \to X' \to \overline{T}$ is also affine. If $E$ is a $p$-completely quasi-coherent sheaf in $((X/S)_{\mathbb{\Delta}})_{|T}$, then $E_{(T'',h)}$ is $p$-completely quasi-coherent on $T''$, and it follows from Proposition 7.4.2 that $R^q h'_*(E_{(T'',h)})$ vanishes for $q > 0$. Since this holds for all $T' \in \mathbb{\Delta}(X/S)$, we can conclude that $R^q j_{T*}E$ vanishes for $q > 0$. Since $\psi_*(E)$ is also $p$-completely quasi-coherent and $z_T$ is affine, it follows also that $R^q\lambda_*\psi_*E$ also vanishes for $q > 0$. Since $\psi_*$ is exact, we conclude that $R^q(\lambda \circ \psi)_*E = 0$ for $q > 0$. Then the vanishing of $R^q j_{T*}E$ implies that $R^q v_{X/S*}(j_{T*}E)$ also vanishes for $q > 0$. □

Proposition 6.3.1 will allow us to use Čech-Alexander complexes to compute prismatic cohomology. Let $T$ be a small $X/S$-prism, which we can view as an object of $\mathbb{\Delta}(X/S)$ or of $\mathbb{\Delta}_s(X/S)$. Since $T$ is small, it is $p$-completely flat over $S$, so by Proposition 7.3.3, the $(n+1)$-fold product $T(n) := T \times_S T \times_S \cdots \times_S T$, computed in the category of $p$-adic formal schemes is $p$-torsion free. We endow it with its natural structure of a crystalline prism. The map $z_{T(n)} : \overline{T}(n) \to X(n)$ is again flat, and since $X/\overline{S}$ is smooth, the diagonal embedding $X \to X(n)$ is a regular immersion. Then $z^{-1}_{T(n)}(X) \to \overline{T}(n)$ is also a regular immersion, so the prismatic envelope $\mathbb{\Delta}_X(T(n))$ of $z^{-1}_{T(n)}(X)$ in $T(n)$ is again small, by Theorem 2.3.8. Now if $E$ is a sheaf of $\mathcal{O}_{X/S}$-modules on $\mathbb{\Delta}(X/S)$, let

$$C^n_T(E) := j_{T(n)*}j^*_{T(n)}(E),$$

and let $C^\bullet_T(E)$ be the Čech-Alexander complex:

$$C^\bullet_T(E) := C^0_T(E) \to C^1_T(E) \to \cdots \to C^n_T(E) \to \cdots$$

with the usual boundary maps.

**Proposition 6.3.2.** *With the notations above, suppose that $E$ is a $p$-completely quasi-coherent sheaf of $\mathcal{O}_{X/S}$-modules on $\mathbb{\Delta}(X/S)$ or $\mathbb{\Delta}_s(X/S)$.*

*(1) For each $n$, the sheaf $C^n_T(E)$ is acyclic for the functors*

$$v_{X/S} : D^+(X/S)_{\mathbb{\Delta}} \to X_{\acute{e}t} \text{ and } D^+(X/S)_{\mathbb{\Delta}_s} \to X_{\acute{e}t}.$$

*(2) If $T$ is a covering of the final object of $(X/S)_{\mathbb{\Delta}}$ or $(X/S)_{\mathbb{\Delta}_s}$, then the natural maps $E \to C^\bullet_T(E)$ and $Rv_{X/S*}E \to v_{X/S*}C^\bullet_T(E)$ are strict quasi-isomorphisms.*

*Proof.* Statement (1) of this proposition follows from statement (3) of Proposition 6.3.1, applied to the prism $T(n)$. If $T$ is a covering of the final object of $(X/S)_{\mathbb{\Delta}}$, then, as explained for example in [7, 5.29], the natural map $E \to C^\bullet_T(E)$ is a quasi-isomorphism, and statement (2) follows. □

**Corollary 6.3.3.** *Let $E$ be a $p$-completely quasi-coherent sheaf of $\mathcal{O}_{X/S}$-modules on $\mathbb{\Delta}(X/S)$ and let $u^*(E)$ be its restriction to $\mathbb{\Delta}_s(X/S)$. Then the natural map*

$$Rv_{X/S*}(E) \to Rv_{X/S*}(u^*(E))$$

*is an isomorphism.*

*Proof.* Without loss of generality we may assume that $X$ is affine and choose a lifting $Y$ of $X$ along with its Frobenius. Then $Y$ is a covering of the final object of $(X/S)_{\mathbb{\Delta}}$ and of $(X/S)_{\mathbb{\Delta}_s}$, and $Rv_{X/S*}(E)$ and $Rv_{X/S*}(u^*(E))$ are both computed by the same Čech-Alexander complex. □

Recall from Definition 4.4.4 the functor $\mathcal{L}_{\mathbb{\Delta}}$ from the category of $\mathcal{O}_Y$-modules to the category of $\mathcal{O}_Y$-modules with prismatic stratification. If $E$ is a $p$-completely quasi-coherent $\mathcal{O}_Y$-module, we write $L(E)$ for the crystal of $\mathcal{O}_{Y/S}$-modules on $\mathbb{\Delta}_s(Y/S)$ corresponding to $\mathcal{L}_{\mathbb{\Delta}}(E)$. A prismatic differential operator $D\colon \mathcal{L}_{\mathbb{\Delta}}(E) \to E'$ induces a map $\mathcal{L}_{\mathbb{\Delta}}(D)\colon \mathcal{L}_{\mathbb{\Delta}}(E) \to \mathcal{L}_{\mathbb{\Delta}}(E')$ (see Definition 1.5.1)) and hence a map $L(D)\colon L(E) \to L(E')$. The following result is the prismatic analog of [7, 6.10] in the crystalline context. Here it seems to be important to work in the small site.

**Proposition 6.3.4.** *Let $j_Y\colon ((Y/S)_{\mathbb{\Delta}_s})_{|Y} \to (Y/S)_{\mathbb{\Delta}_s}$ and $\psi\colon ((Y/S)_{\mathbb{\Delta}_s})_{|Y} \to Y_{\mathrm{pcf}}$ be the morphisms as in Proposition 6.3.1, and let $E$ be a $p$-completely quasi-coherent sheaf of $\mathcal{O}_Y$-modules on $Y_{\mathrm{pcf}}$.*

*(1) There is a natural isomorphism*

$$L(E) \cong j_{Y*}(\psi^*(E))$$

*of crystals of $\mathcal{O}_{Y/S}$-modules on $\mathbb{\Delta}_s(Y/S)$.*

*(2) If $\lambda\colon Y_{\mathrm{pcf}} \to Y_{\acute{e}t}$ is the natural map, there is a natural isomorphism:*

$$v_{Y/S*}(L(E)) \cong \lambda_*(E).$$

*If $D\colon \mathcal{L}_{\mathbb{\Delta}}(E) \to E'$ is a prismatic differential operator, let $s^*$ be the natural map $E \to t_*s^*E = \mathcal{L}_{\mathbb{\Delta}}(E)$ (see diagram 4.4.14). Then the diagram*

$$\begin{array}{ccc} v_{Y/S*}(L(E)) & \xrightarrow{v_{Y/S*}(L(D))} & v_{Y/S*}(L(E')) \\ \cong\downarrow & & \downarrow\cong \\ \lambda_*(E) & \xrightarrow{D \circ s^*} & \lambda_*(E') \end{array}$$

*commutes.*

*(3) $R^q v_{Y/S*} L(E) = 0$ for $q > 0$.*

*Proof.* We first define a morphism $j_Y^* L(E) \to \psi^*(E)$ as follows. For $(T', h) \in \mathbb{\Delta}_s(X/S)_{|Y}$, the definitions tell us that

$$(j_Y^* L(E))_{(T',h)} := L(E)_{T'} = h^*(\mathcal{L}_{\mathbb{\Delta}}(E)),$$

where as usual $h^*$ means the completed pullback. Composing with the pullback of the map (4.4.14) $\mathcal{L}_{\mathbb{\Delta}}(E) \to E$, we find the desired:

$$(j_Y^* L(E))_{(T',h)} = h^*(\mathcal{L}_{\mathbb{\Delta}}(E)) \to h^*(E) = (\psi^*(E))_{(T',h)}.$$

To prove the proposition, it will be enough to check that the adjoint to this construction is an isomorphism. Our claim is that for each $T \in \mathbb{\Delta}_s(Y/S)$, the map $L(E)_T \to (j_{Y*}(\psi^* E))_T$ is an isomorphism. We can check this $p$-completely flat locally on $T$, and so by Proposition 6.1.11, we may assume that there is a morphism of prisms $h\colon T \to Y$. We find a commutative diagram:

$$\begin{array}{ccccc}
\mathbb{\Delta}_{\overline{Y}}(T \times_S Y) & \xrightarrow{f} & \mathbb{\Delta}_Y(1) & \xrightarrow{p_2} & Y \\
{\scriptstyle p_T}\downarrow & & \downarrow{\scriptstyle p_1} & & \\
T & \xrightarrow{h} & Y, & &
\end{array}$$

in which the square is Cartesian by Proposition 6.1.8. Then, using the definitions and the fact that $p_1\colon \mathbb{\Delta}_Y(1) \to Y$ is affine, we get that:

$$\begin{aligned}
L(E)_T &= h^*(\mathcal{L}_{\mathbb{\Delta}}(E)) \\
&= h^*(\mathcal{O}_{\mathbb{\Delta}_Y(1)} \hat{\otimes} E) \\
&= h^*(p_{1*}(p_2^*(E))) \\
&= p_{T*}(f^*(p_2^*(E))).
\end{aligned}$$

On the other hand, since $\mathbb{\Delta}_Y(T \times_S Y) \to Y$ represents $j_Y^{-1}(T)$ we have

$$j_{Y*}(\psi^*(E))_T = p_{T*}((\psi^*(E))_{\mathbb{\Delta}_Y(T \times_S Y)}) = p_{T*}(f^*(p_2^*(E))),$$

proving statement (1). The remaining statements then follow from Propositions 6.3.1 and 1.5.3. $\square$

We are now ready to prove that our complexes do indeed calculate prismatic cohomology.

**Theorem 6.3.5.** *Let $Y/S$ be a $p$-completely smooth morphism of crystalline prisms, and let $i\colon X \to \overline{Y}$ be a closed immersion, where $X/\overline{S}$ is smooth. If $E$ is a crystal of $\mathcal{O}_{X/S}$-modules on $\mathbb{\Delta}(X/S)$, let $(E_{\mathbb{\Delta}_X(Y)}, \nabla')$ be the corresponding object of* MICP$(X/Y/S)$*, as described in Corollary 6.2.4. Then there is a canonical and functorial strict quasi-isomorphism:*

$$Rv_{X/S*}E \cong (E_{\mathbb{\Delta}_X(Y)} \otimes \Omega^{\bullet}_{Y/S}, d').$$

*Proof.* In fact we shall deduce the theorem from the special case in which $X = \overline{Y}$, stated as Lemma 6.3.7 below. This deduction uses the following lemma.

**Lemma 6.3.6.** *If $E$ is an abelian sheaf on $\mathbb{A}(X/S)$, then $R^q i_{\mathbb{A}*}(E) = 0$ for $q > 0$.*

*Proof.* If $T$ is an object of $\mathbb{A}(Y/S)$ and $E$ is a sheaf on $\mathbb{A}(X/S)$, then the value of $i_{\mathbb{A}*}(E)$ on $T$ is $E_{\mathbb{A}_X(T)}$ (using the notation in Equation 6.1.1). It follows that the functor $i_{\mathbb{A}*}$ is exact, and hence that $R^q i_{\mathbb{A}*}$ vanishes if $q > 0$. □

Now let $E$ be a crystal of $\mathcal{O}_{X/S}$-modules. Since $v_{X/S} = v_{Y/S} \circ i_{\mathbb{A}}$, the lemma implies that $Rv_{X/S*}(E) \cong Rv_{Y/S*}(i_{\mathbb{A}*}(E))$. Unfortunately $i_{\mathbb{A}*}(E)$ is not a crystal of $\mathcal{O}_{Y/S}$-modules, but, by statement (2) of Proposition 6.2.2, its restriction $\tilde{E} := u^* i_{\mathbb{A}*} E$ to $\mathbb{A}_s(Y/S)$ is. Lemma 6.3.7 below, together with Corollary 6.3.3, will tell us that $Rv_{Y/S*}(\tilde{E})$ is represented by the $p$-de Rham complex of $(\tilde{E}_Y, \nabla')$. Since $\tilde{E}_Y = E_{\mathbb{A}_X(Y)}$, this will prove the theorem. □

**Lemma 6.3.7.** *Let $Y/S$ be a $p$-completely smooth morphism of crystalline prisms. let $E$ be a crystal of $\mathcal{O}_{Y/S}$-modules on $\mathbb{A}(Y/S)$ (resp., on $\mathbb{A}_s(Y/S)$), and let $(E_Y, \nabla')$ be the corresponding $\mathcal{O}_Y$-module with integrable $p$-connection (6.2.3). Then there is a canonical strict quasi-isomorphism*

$$Rv_{Y/S*}(E) \cong (E_Y \otimes \Omega^{\bullet}_{Y/S}, d').$$

*Proof.* Our proof will follow the method of proof of its crystalline analog as carried out in [7]. We first explain the case when $E = \mathcal{O}_Y$, with its canonical $p$-connection. The map $s^*$ in diagram (4.4.14) defines a horizontal morphism $\mathcal{O}_Y \to \mathcal{L}_{\mathbb{A}}(\mathcal{O}_Y)$. It follows from diagram (1.5.1) that this extends to a morphism of complexes of modules with $p$-connection:

$$(\mathcal{O}_Y, d') \longrightarrow \mathcal{L}_{\mathbb{A}}(\Omega^{\bullet}_{Y/S}, d'). \tag{6.3.1}$$

Tensoring with $E_Y$ and composing with the isomorphism $\beta$ of Proposition 4.4.5, we get:

$$(E_Y, \nabla') \longrightarrow E_Y \hat{\otimes} \mathcal{L}_{\mathbb{A}}(\Omega^{\bullet}_{Y/S}, d') \xrightarrow[\cong]{\beta} \mathcal{L}_{\mathbb{A}}(E_Y \otimes \Omega^{\bullet}_{Y/S}, d'), \tag{6.3.2}$$

and hence a corresponding morphism of complexes of crystals on $\mathbb{A}(Y/S)$:

$$E \to L(E_Y \otimes \Omega^{\bullet}_{Y/S}, d'). \tag{6.3.3}$$

Lemma 6.3.8 below will tell us that this map is a strict quasi-isomorphism, after restricting to $\mathbb{A}_s(Y/T)$, and Corollary 6.3.3 tells that this restriction is anodyne. By statement (3) of Proposition 6.3.4, each term of the complex $L(E \otimes \Omega^{\bullet}_{Y/S})$ is acyclic for $v_{Y/S*}$, and by statement (2) of that proposition, $v_{Y/S*} L(E \otimes \Omega^{\bullet}_{Y/S}, d')$ identifies with $(E_Y \otimes \Omega^{\bullet}_{Y/S}, d')$. Thus the following lemma will complete the proof of the theorem. □

**Lemma 6.3.8.** *Let $(E_Y, \nabla')$ be an $\mathcal{O}_Y$-module with integrable and quasi-nilpotent $p$-connection, and let $E$ be the corresponding crystal on $\mathbb{\Delta}_s(Y/T)$. Then the map (6.3.3) is a strict quasi-isomorphism.*

*Proof.* The isomorphism $\beta$ induces an isomorphism of complexes of crystals of $\mathcal{O}_{Y/S}$-modules on $\mathbb{\Delta}_s(Y/S)$:

$$E\hat{\otimes}L(\Omega^\bullet_{Y/S}, d') \cong L(E_Y \otimes \Omega^\bullet_{Y/S}, d').$$

Thus it will be enough to prove that the morphism $E \to E\hat{\otimes}L(\Omega^\bullet_{Y/S}, d')$ induced by the first arrow in equation (6.3.2) is a strict quasi-isomorphism. Since the terms of the complex $L(\Omega^\bullet_{Y/S}, d')$ are $p$-completely flat, it is enough to prove this statement when $E = \mathcal{O}_{Y/S}$. Indeed, this case will show that the augmented complex $\mathcal{O}_T \to L(\Omega^\bullet_{Y/S})_T$ is strictly exact, and this will remain the case after the completed tensor product with $E_T$.

The claim is that for every $T \in \mathbb{\Delta}_s(Y/S)$, the map $\mathcal{O}_T \to L(\Omega^\bullet_{Y/S})_T$ is a strict quasi-isomorphism. Thanks to Proposition 7.3.5, this can be verified after replacing $T$ by a $p$-completely flat cover, so by Proposition 6.1.11 we may without loss of generality that there is a morphism of $Y$-prisms $T \to Y$. Then, as we saw in Propositions 6.3.4 and 4.4.5

$$L(\Omega^\bullet_{Y/S})_T \cong \mathcal{O}_T\hat{\otimes}_{\mathcal{O}_Y}\mathcal{L}_{\mathbb{\Delta}}(\Omega^\bullet_{Y/S}).$$

Since $T$ is small, the map $T \to Y$ is $p$-completely flat, so again by Proposition 7.3.5, we are reduced to checking our claim when $T = Y$. Moroever, it follows from Proposition 4.4.6 that the complex $\mathcal{L}_{\mathbb{\Delta}}(\Omega^\bullet_{Y/S})$ identifies with the complex $(\Omega^\bullet_{\mathbb{\Delta}_Y(1)/Y}, d')$. If we write $Z$ for $Y(1)$, then the first projection $Z \to Y$ is a $p$-completely smooth morphism of crystalline prisms, with a section defined by the diagonal, and $\mathbb{\Delta}_Y(1)$ is the prismatic envelope of this section. Then the simplest form of prismatic Poincaré lemma, statement (2) of Lemma 4.3.2, implies that the map $\mathcal{O}_Y \to \Omega^\bullet_{Z/Y}$ is a strict quasi-isomorphism. This concludes the proof. ☐

**Corollary 6.3.9.** *In the situation of Theorem 6.3.5, there is a canonical and functorial quasi-isomorphism:*

$$Rv_{X/S*}(\mathcal{O}_{X/\overline{S}}) \cong (\mathcal{O}_{\overline{\mathbb{\Delta}}} \otimes \Omega^\bullet_{\overline{Y}/\overline{S}}, \nabla'),$$

*where $(\mathcal{O}_{\overline{\mathbb{\Delta}}}\otimes\Omega^\bullet_{\overline{Y}/\overline{S}}, \nabla')$ is the Higgs complex of the Higgs field $\nabla'$ on $\mathcal{O}_{\overline{\mathbb{\Delta}}}$ defined in Theorem 3.3.3.*

*Proof.* This follows from Theorem 6.3.5 and Proposition 5.2.3. ☐

## 6.4 PD-prisms and the prismatic F-transform

In this section we introduce a variant of the prismatic site which lies between the prismatic and crystalline theories. Inspired by work of Oyama [35] and Xu [40], it will allow us to give a topos-theoretic interpretation of the F-transform. It also provides a canonical factorization of the prismatic Frobenius endomorphism which clarifies why it is an isogeny. Although this construction is not formalized explicitly in [10], some of its key aspects are used in some of the comparison theorems there.

Recall that if $X$ is a scheme in characteristic $p$, the scheme theoretic image $F_X(X)$ of its absolute Frobenius endomorphism is the closed subscheme defined by the ideal of sections of $\mathcal{O}_X$ whose $p$th power is zero. There is a canonical factorization

$$X \xrightarrow{F_X} X \quad = \quad X \xrightarrow{f_X} F_X(X) \xrightarrow{j_X} X, \tag{6.4.1}$$

where $j_X$ is a closed immersion. If $S$ is a crystalline prism, we may write $F_S(\overline{S})$ or $F(\overline{S})$ instead of $F_{\overline{S}}(\overline{S})$. If $X/S$ is an $\overline{S}$-scheme, there is a commutative diagram:

$$\begin{array}{ccccc} X & \xrightarrow{f_X} & F_X(X) & \xrightarrow{j_X} & X \\ \downarrow & & \downarrow & \swarrow & \\ S & \xrightarrow{F_S} & S, & & \end{array} \tag{6.4.2}$$

so that $j_X$ is a morphism of $S$-schemes. Note that, since $F_X$ is a homeomorphism and $f_X^\sharp$ is injective, the morphism $f_X$ is an epimorphism in the category of schemes.

If $X$ is embedded as a closed subscheme of a crystalline prism $Y$, we let $X^\phi := F_{\overline{Y}}^{-1}(X)$. The following commutative diagram illustrates two factorizations of $F_{X^\phi}$:

$$\begin{array}{ccccc} X^\phi & \xrightarrow{f_{X^\phi}} & F_{X^\phi}(X^\phi) & & \\ {\scriptstyle f_{X/Y/S}}\downarrow & & {\scriptstyle k_{X/Y}}\downarrow & \searrow^{j_{X^\phi}} & \\ X^{(p)} & \xrightarrow{\pi_{X/S}} & X & \xrightarrow{i_{X^\phi}} & X^\phi. \end{array} \tag{6.4.3}$$

**Definition 6.4.1.** *Let $S$ be a crystalline prism and $X/\overline{S}$ a morphism of schemes. Then an $X/S$-$\phi$-prism is a pair $(T, y_T)$, where $T$ is a crystalline $S$-prism and $y_T \colon F_T(\overline{T}) \to X$ is an $\overline{S}$-morphism. (We sometimes just say $X$-$\phi$-prism instead of $X/S$-$\phi$ prism.) We denote by $\mathbb{A}_\phi(X/S)$ the category of $X/S$-$\phi$-prisms and endow it with the $p$-completely flat topology. If $X \to Y$ is a closed immersion of $X$ in a crystalline prism, then the $\phi$-prismatic envelope of $X$ in $Y$ is the universal morphism $\pi_{\mathbb{A}} \colon (\mathbb{A}_X^\phi(Y), y_{\mathbb{A}}) \to Y$ from an $X/S$-$\phi$-prism to $Y$.*

We omit the verifications that $\mathbb{A}_\phi(X/S)$ forms a site and that its formation is functorial. Crystals on $\mathbb{A}_\phi(X/S)$ are defined in the obvious way.

**Remark 6.4.2.** If $T$ is a crystalline prism, Proposition 2.1.4 tells us that $F_{\overline{T}}(\overline{T})$, is equal to $T_{\mathbb{PD}}$, the smallest PD-subscheme of $T$. Thus we could just as well have defined $\mathbb{A}_\phi(X/S)$ to be the site whose objects are "PD-prisms," *i.e.*, pairs $(T, y_T)$, where $T$ is a crystalline $S$-prism and $y_T$ is an $S$-morphism $T_{\mathbb{PD}} \to X$, and denote this site by ${\mathbb{\Delta}}_{\mathbb{PD}}(X/S)$.

Theorem 6.4.6 below implies that, when $S$ is perfect, the sites ${\mathbb{\Delta}}_\phi(X/S)$ and ${\mathbb{\Delta}}(X/S)$ give rise to equivalent topoi. The key geometric construction we shall need is the following description of $\phi$-prismatic envelopes. Recall that if $X \subseteq \overline{Y}$, then $X^\phi := F_{\overline{Y}}^{-1}(X)$, and, from Theorem 2.3.13, that $\mathbb{PD}_X(Y)$ identifies with ${\mathbb{\Delta}}_{X^\phi}(Y)$. This will allow us to relate $\phi$-prismatic envelopes to divided power envelopes.

**Proposition 6.4.3.** *If $Y$ is a crystalline prism and $i\colon X \to Y$ is a closed immersion, let $X^\phi$ be the inverse image of $X$ under $F_{\overline{Y}}$ and let $({\mathbb{\Delta}}_{X^\phi}(Y), z_{{\mathbb{\Delta}}^\phi}, \pi_{{\mathbb{\Delta}}^\phi})$ be the prismatic envelope of $X^\phi$ in $Y$. Then there is a unique morphism $y_{\mathbb{\Delta}}\colon F(\overline{{\mathbb{\Delta}}}_{X^\phi}(Y)) \to X$ such that the left square below commutes, and the resulting outer diagram represents $({\mathbb{\Delta}}_X^\phi(Y), y_{\mathbb{\Delta}}, \pi_{\mathbb{\Delta}})$.*

$$\begin{array}{ccccc}
F(\overline{{\mathbb{\Delta}}}_{X^\phi}(Y)) & \xrightarrow{j_{\overline{\mathbb{\Delta}}}} & \overline{{\mathbb{\Delta}}}_{X^\phi}(Y) & \longrightarrow & {\mathbb{\Delta}}_{X^\phi}(Y) \\
\Big\downarrow{\scriptstyle y_{\mathbb{\Delta}}} & & \Big\downarrow{\scriptstyle z_{{\mathbb{\Delta}}^\phi}} & \swarrow{\scriptstyle \pi_{{\mathbb{\Delta}}^\phi}} & \\
X & \xrightarrow{\quad i \quad} & Y. & &
\end{array}$$

*If $Y/S$ is $p$-completely smooth and $X/\overline{S}$ is smooth, then $({\mathbb{\Delta}}_X^\phi(Y), y_{\mathbb{\Delta}}$ is small.*

*Proof.* The factorization (6.4.1) of the Frobenius endomorphism of $\overline{{\mathbb{\Delta}}}_{X^\phi}(Y)$ gives the top row of the following commutative diagram:

$$\begin{array}{ccccccc}
\overline{{\mathbb{\Delta}}}_{X^\phi}(Y) & \xrightarrow{f_{\overline{\mathbb{\Delta}}}} & F({\mathbb{\Delta}}_{X^\phi}(Y)) & \xrightarrow{j_{\overline{\mathbb{\Delta}}}} & \overline{{\mathbb{\Delta}}}_{X^\phi}(Y) & & \\
\Big\downarrow{\scriptstyle z_{{\mathbb{\Delta}}^\phi}} & & \Big\downarrow & \searrow{\scriptstyle y_{\mathbb{\Delta}}} & & \searrow{\scriptstyle z_{{\mathbb{\Delta}}^\phi}} & \\
X^\phi & \xrightarrow{f_{X^\phi}} & F_{X^\phi}(X^\phi) & \xrightarrow{k_{X/Y}} & X & \xrightarrow{i_{X^\phi}} & X^\phi.
\end{array} \tag{6.4.4}$$

The bottom row comes from the factorization 6.4.3 of the Frobenius endomorphism of $X^\phi$, and the remaining solid arrows exist because of the commutativity of the diagram (6.4.3). The dotted arrow $y_{\mathbb{\Delta}}$ is defined to make the triangle (and hence the diagram) commute, and endows ${\mathbb{\Delta}}_{X^\phi}(Y)$ with the structure of an $X/S$-$\phi$-prism.

If $(T, y_T, \pi_T)$ is an $X/S$-$\phi$-prism over $Y$, the commutative diagram

$$\begin{array}{ccccc}
\overline{T} & \xrightarrow{f_T} & F(\overline{T}) & \xrightarrow{y_T} & X \\
\Big\downarrow{\scriptstyle \overline{\pi}_T} & & \Big\downarrow{\scriptstyle \pi_T \circ j_{\overline{T}}} & \swarrow{\scriptstyle \overline{i}} & \\
\overline{Y} & \xrightarrow{F_{\overline{Y}}} & \overline{Y}, & &
\end{array}$$

shows that $\overline{\pi}_T$ factors through $X^\phi$, hence that $\pi_T$ factors through $\mathbb{A}_{X^\phi}(Y)$. This proves that $(\mathbb{A}_{X^\phi}(Y), y_{\mathbb{A}}, \pi_{\mathbb{A}^\phi})$ does enjoy the requisite universal property. If $Y/S$ is $p$-completely smooth and $X/\overline{S}$ is smooth, then $X^{(p)}$ is regularly immersed in $\overline{Y}^{(p)}$, and since $F_{\overline{Y}/S}$ is flat and $X^\phi = F^{-1}_{\overline{Y}/S}(X^{(p)})$, it is also regularly immersed in $\overline{Y}$. Then it follows from Theorem 2.3.8 that its prismatic neighborhood is small. □

Combining this result with Theorem 2.3.13, we see in the following corollary that $\phi$-prismatic envelopes are essentially the same as PD-envelopes, at least in the smooth case.

**Corollary 6.4.4.** *In the situation of Proposition 6.4.3, suppose that $Y/S$ is $p$-completely smooth and that $X/\overline{S}$ is smooth. Then there is a natural $Y$-isomorphism $\mathbb{PD}_X(Y) \to \mathbb{A}^\phi_X(Y)$ fitting into a commutative diagram*

$$
\begin{array}{ccccc}
\mathbb{PD}_X(Y) & \xrightarrow{\cong} & \mathbb{A}^\phi_X(Y) & \xrightarrow{\cong} & \mathbb{A}_{X^\phi}(Y) \\
\uparrow & & \uparrow & & \uparrow \\
\mathbb{A}^\phi_X(Y)_{\mathbb{PD}} & \xrightarrow{j_{\overline{\mathbb{A}}}} & \overline{\mathbb{A}}^\phi_X(Y) & \xrightarrow{\cong} & \overline{\mathbb{A}}_{X^\phi}(Y) \\
 & \searrow^{y_{\mathbb{A}}} & & & \downarrow{\scriptstyle z_{\mathbb{A}^\phi}} \\
 & & X & \xrightarrow{i_{X^\phi}} & X^\phi.
\end{array}
$$
□

We are now ready to define a pair of functors $A$ and $B$ which allow us to factor the prismatic Frobenius morphism. We shall see later that $B$ is a geometric incarnation of the F-transform.

**Definition 6.4.5.** *Suppose that $S$ is a crystalline prism and that $X/\overline{S}$ is a smooth morphism.*

*(1) $A\colon \mathbb{A}(X/S) \to \mathbb{A}_\phi(X/S)$ is the functor*

$$A\colon (T, z_T) \mapsto (T, z_T \circ j_{\overline{T}}),$$

*where $j_{\overline{T}}\colon F_{\overline{T}}(\overline{T}) \to \overline{T}$ is the inclusion (6.4.1).*

*(2) $B\colon \mathbb{A}_\phi(X/S) \to \mathbb{A}(X^{(p)}/S)$ is the functor*

$$B\colon (T, y_T) \mapsto (T, \tilde{z}_T),$$

*where $\tilde{z}_T\colon \overline{T} \to X^{(p)}$ is the unique $S$-morphism making the following*

*diagram commute:*

$$\begin{CD} \overline{T} @>f_{\overline{T}}>> F(\overline{T}) \\ @V{\tilde{z}_T}VV @VV{y_T}V \\ X^{(p)} @>{\pi_{X/S}}>> X. \end{CD}$$

*(Note: This morphism exists because $y_T \circ f_{\overline{T}}$ is an $F_{\overline{S}}$-morphism.)*

**Theorem 6.4.6.** *Suppose that $S$ is a crystalline prism and that $X/\overline{S}$ is a smooth morphism.*

*(1) There is a 2-commutative diagram of continuous and cocontinuous functors:*

$$\begin{array}{ccccc} \mathbb{\Delta}(X/S) & \xrightarrow{A} & \mathbb{\Delta}_\phi(X/S) & & \\ & \searrow^{\mathbb{\Delta}(F_{X/S})} & \downarrow{\scriptstyle B} & \searrow^{\mathbb{\Delta}_\phi(F_{X/S})} & \\ & & \mathbb{\Delta}(X^{(p)}/S) & \xrightarrow{A^{(p)}} & \mathbb{\Delta}_\phi(X^{(p)}/S). \end{array}$$

*(2) If $\mathcal{F}$ is a sheaf on $\mathbb{\Delta}_\phi(X/S)$ (resp. on $\mathbb{\Delta}(X^{(p)}/S)$), then the presheaf*

$$(T, z_T) \mapsto \mathcal{F}(A(T, z_T)) \quad \text{resp.} \quad (T, y_T) \mapsto \mathcal{F}(B(T, y_T))$$

*is the sheaf $A^{-1}(\mathcal{F})$ (resp. $B^{-1}(\mathcal{F})$). Moreover, there are isomorphisms:*

$$\begin{array}{rcl} A^{-1}(\mathcal{O}_{X/S}) & \xrightarrow{\cong} & \mathcal{O}_{X/S} \\ B^{-1}(\mathcal{O}_{X^{(p)}/S}) & \xrightarrow{\cong} & \mathcal{O}_{X/S}, \end{array}$$

*and hence morphisms of ringed topoi:*

$$\begin{array}{rcl} A_{\mathbb{\Delta}} \colon ((X/S)_{\mathbb{\Delta}}, \mathcal{O}_{X/S}) & \to & ((X/S)_{\mathbb{\Delta}_\phi}, \mathcal{O}_{X/S}) \\ B_{\mathbb{\Delta}} \colon ((X/S)_{\mathbb{\Delta}_\phi}, \mathcal{O}_{X/S}) & \to & ((X^{(p)}/S)_{\mathbb{\Delta}}, \mathcal{O}_{X^{(p)}/S}). \end{array}$$

*(3) The morphism $B_{\mathbb{\Delta}}$ is an equivalence of ringed topoi and induces an equivalence from the category of crystals of $\mathcal{O}_{X^{(p)}/S}$-modules on $\mathbb{\Delta}(X^{(p)}/S)$ to the category of crystals of $\mathcal{O}_{X/S}$-modules on $\mathbb{\Delta}_\phi(X/S)$.*

*Proof.* To check the commutativity of the diagram in (1), let $(T, z_T)$ be an object of $\mathbb{\Delta}(X/S)$. Then $A(T, Z_T) = (T, y_T)$, where $y_T := z_T \circ j_T$ and $B(A(T, z_T)) = (T, \tilde{z}_T)$, where $\tilde{z}_T$ is the unique $S$-morphism such that $\pi \circ \tilde{z}_T = y_T \circ f_{\overline{T}}$. Since $\mathbb{\Delta}(F_{X/S})(T, z_T) = (T, F_{X/S} \circ z_T)$, we must check that the two maps $F_{X/S} \circ z_T$ and $\tilde{z}_T\colon \overline{T} \to X^{(p)}$ from $\overline{T}$ to $X^{(p)}$ agree. Since $F_{X/S} \circ z_T$ is an $S$-morphism, it will suffice to check that $\pi_{X/S} \circ F_{X/S} \circ z_T = y_T \circ f_{\overline{T}}$. This follows from the following commutative diagram.

$$\begin{array}{ccccc}
\overline{T} & \xrightarrow{f_{\overline{T}}} & F(\overline{T}) & \xrightarrow{j_{\overline{T}}} & \overline{T} \\
{\scriptstyle z_T}\downarrow & & & \searrow^{y_T} & \downarrow{\scriptstyle z_T} \\
X & \xrightarrow{F_{X/S}} & X^{(p)} & \xrightarrow{\pi_{X/S}} & X.
\end{array}$$

This shows $B \circ A = \mathbb{\Delta}(F_{X/S})$. We leave the second triangle for the reader.

We omit the proof that $A$ is continuous and cocontinuous. To see that $B$ is continuous, observe first that it takes coverings to coverings. Furthermore, if $(\widetilde{T}, y_{\widetilde{T}}) \to (T, y_T)$ and $(T', y'_T) \to (T, y_T)$ are morphisms in $\mathbb{\Delta}_\phi(X/S)$ and $(\widetilde{T}, y_{\widetilde{T}}) \to (T, y_T)$ is $p$-completely flat, then it follows from the construction in Proposition 6.1.2 that the fiber product $(\widetilde{T}, y_{\widetilde{T}}) \times_{(T,y_T)} (T', y_{T'})$ in the category $\mathbb{\Delta}_\phi(X/S)$ is given by the usual fiber product in the category of $p$-adic formal schemes. Since the same holds in the category $\mathbb{\Delta}(X^{(p)}/S)$, we see that $B$ preserves fiber products (at least) in this case. To see that $B$ is cocontinuous, let $(T, y_T)$ be an object of $\mathbb{\Delta}_\phi(X/S)$ and let $u\colon (\widetilde{T}, z_{\widetilde{T}}) \to (T, \tilde{z}_T) := B(T, y_T)$ be a covering. We shall see that this covering is induced by a covering of $(T, y_T)$. In fact, $u\colon \widetilde{T} \to T$ is a $p$-completely faithfully flat morphism and $z_{\widetilde{T}}\colon \widetilde{T}_1 \to X^{(p)}$ is a morphism such that $z_{\widetilde{T}} = \tilde{z}_T \circ u_1$. Let $y_{\widetilde{T}} := y_T \circ u_\phi\colon F(\widetilde{T}_1) \to X$, where $u_\phi\colon F(\widetilde{T}_1) \to F(T_1)$ is the morphism induced by $u$. Then $(\widetilde{T}, y_{\widetilde{T}})$ is an object of $\mathbb{\Delta}_\phi(X/S)$, and $u$ defines a $p$-completely flat covering $\tilde{u}\colon (\widetilde{T}, y_{\widetilde{T}}) \to (T, y_T)$. Let $(\widetilde{T}, \tilde{z}_{\widetilde{T}}) := B(\widetilde{T}, y_{\widetilde{T}})$, so that $\tilde{z}_{\widetilde{T}}\colon \widetilde{T}_1 \to X^{(p)}$ is the unique $S$-morphism such that $\pi_{X/S} \circ \tilde{z}_{\widetilde{T}} = y_{\widetilde{T}} \circ f_{\widetilde{T}_1}$. But then

$$\begin{aligned}
\pi_{X/S} \circ \tilde{z}_{\widetilde{T}} &= y_{\widetilde{T}} \circ f_{\widetilde{T}_1} \\
&= y_T \circ u_\phi \circ f_{\widetilde{T}_1} \\
&= y_T \circ f_{\overline{T}} \circ u_1 \\
&= \pi_{X/S} \circ \tilde{z}_T \circ u_1 \\
&= \pi_{X/S} \circ z_{\widetilde{T}}.
\end{aligned}$$

This implies that $\tilde{z}_{\widetilde{T}} = z_{\widetilde{T}}$, so in fact $B(\tilde{u})$ coincides with the original covering $u\colon (\widetilde{T}, \tilde{z}_{\widetilde{T}}) \to (T, \tilde{z}_T)$.

If $\mathcal{F}$ is a sheaf on $\mathbb{\Delta}_\phi(X/S)$, recall that $A^{-1}(\mathcal{F})$ is the sheaf associated to the presheaf which takes an object $(T, z_T)$ of $\mathbb{\Delta}(X/S)$ to $\mathcal{F}(A(T, z_T))$. Since $A$ is continuous, this presheaf is in fact a sheaf. Note that $A(T, z_T) = (T, y_T)$, and so the sheaf $A^{-1}(\mathcal{F})_{(T,z_T)}$ on $T$ is the same as the sheaf $\mathcal{F}_{(T,y_T)}$. Applying this to $\mathcal{O}_{X/S}$, we see that $A^{-1}(\mathcal{O}_{X/S}) = \mathcal{O}_{X/S}$. It follows that if $\mathcal{F}$ is a sheaf of $\mathcal{O}_{X/S}$-modules on $\mathbb{\Delta}_\phi(X/S)$, then $A^{-1}(\mathcal{F}) \cong A^*(\mathcal{F})$, and that $A^*(\mathcal{F})$ is $p$-completely quasi-coherent if and only if $\mathcal{F}$ is. The same argument works for $B$.

**Lemma 6.4.7.** *In the situation of Theorem 6.4.6,*

*(1) The functor $B$ is fully faithful.*

*(2) Every object of $\mathbb{\Delta}(X^{(p)}/S)$ admits a cover by an object in the image of $B$.*

*Proof.* The proof is a modification of Oyama's argument, as explained by Xu [40, 9.8]. To see that $B$ is fully faithful, suppose that $(T, y_T)$ and $(\widetilde{T}, y_{\widetilde{T}})$ are objects of $\mathbb{\Delta}_\phi(X/S)$. Then a morphism $(\widetilde{T}, y_{\widetilde{T}}) \to (T, y_T)$ is a morphism $g\colon \widetilde{T} \to T$ such that

$$y_T \circ g_\phi = y_{\widetilde{T}}, \tag{6.4.5}$$

and a morphism $B(\widetilde{T}, y_{\widetilde{T}}) \to B(T, y_T)$ is a morphism $g\colon \widetilde{T} \to T$ such that

$$\tilde{z}_T \circ \overline{g} = \tilde{z}_{\widetilde{T}}. \tag{6.4.6}$$

Thus it is obvious that $B$ is faithful. To see that it is full, suppose that $g\colon \widetilde{T} \to T$ defines a morphism $B(\widetilde{T}, y_{\widetilde{T}}) \to B(T, y_T)$. Then

$$\begin{aligned} \tilde{z}_T \circ \overline{g} &= \tilde{z}_{\widetilde{T}} \\ \pi_{X/S} \circ \tilde{z}_T \circ \overline{g} &= \pi_{X/S} \circ \tilde{z}_{\widetilde{T}} \\ y_T \circ f_{\overline{T}} \circ \overline{g} &= y_{\widetilde{T}} \circ f_{\widetilde{T}_1} \\ y_T \circ g_\phi \circ f_{\widetilde{T}_1} &= y_{\widetilde{T}} \circ f_{\widetilde{T}_1}. \end{aligned}$$

Since $f_{\widetilde{T}_1}$ is a scheme-theoretic epimorphism, it follows that $y_T \circ g_\phi = y_{\widetilde{T}}$, so $g$ defines a morphism $(\widetilde{T}, y_{\widetilde{T}}) \to (T, y_T)$. This proves statement (1).

To prove (2), suppose that $(T', z_{T'})$ is an object of $\mathbb{\Delta}(X^{(p)}/S)$. Without loss of generality, we assume that $X^{(p)}$ and $T'$ are affine. Choose a $p$-completely smooth crystalline $S$-prism $Y/S$ lifting $X/\overline{S}$. Then $Y^{(p)}/S$ is again $p$-completely smooth, and by Theorem 1.7.11, we may, after replacing $T'$ by a $p$-completely flat covering, assume that there is a morphism $(T', z_{T'}) \to (Y^{(p)}, \mathrm{id}_{X^{(p)}})$. The map $F_{Y/S}\colon Y \to Y^{(p)}$ is $p$-completely flat, and hence so is the map $v\colon T := T' \times_{Y^{(p)}} Y \to T'$. Note that $\overline{T} \cong \overline{T}' \times_{X^{(p)}} X$, and that we have a commutative diagram:

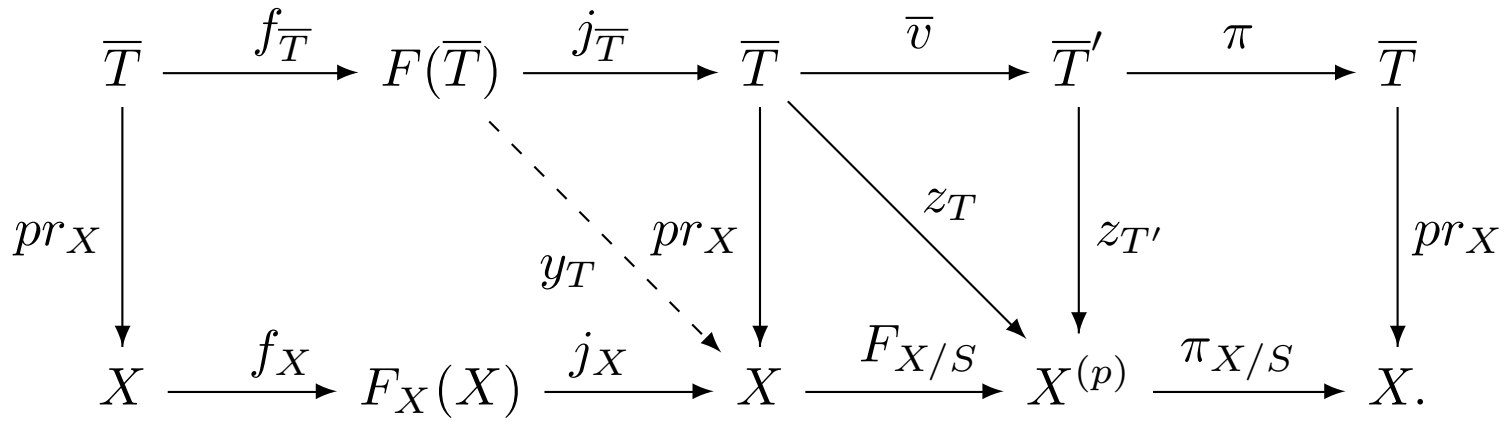

The map $z_T$ defined by the diagram endows $T$ with the structure of an $X^{(p)}/S$-prism, the map $y_T$ gives $T$ the structure of an $X/S$-$\phi$ prism, and $v$ defines a $p$-completely flat cover $(T, z_T) \to (T', z_{T'})$. Note that $y_T \circ f_{\overline{T}} = F_X \circ pr_X$. We claim that $B(T, y_T)$ is equal to $(T, z_T)$, which will prove (2). We have:

$$\pi_{X/S} \circ z_T \quad = \quad \pi_{X/S} \circ z_{T'} \circ \overline{v}$$

$$\begin{aligned} &= F_X \circ pr_X \\ &= y_T \circ f_{\overline{T}}. \end{aligned}$$

By the definition of $B(T, y_T)$ this proves the claim. □

Since the $p$-completely flat topology comes from a pre-topology, Lemma 6.4.7 and a general theorem of Oyama [35, 4.2.1] imply that $B$ induces an equivalence of topoi. As we have seen, $B^{-1}(\mathcal{O}_{X^{(p)}/S}) \cong \mathcal{O}_{X/S}$, so in fact $B_{\mathbb{\Delta}}$ is an equivalence of ringed topoi. If $E'$ is a crystal of $\mathcal{O}_{X^{(p)}/S}$-modules on $\mathbb{\Delta}(X^{(p)}/S)$, it is easy to check that $B^*(E')$ is a crystal of $\mathcal{O}_{X/S}$-modules on $\mathbb{\Delta}_\phi(X/S)$. To prove the converse, one can argue locally, using Lemma 6.4.7. We leave details to the reader, following the method of Oyama. □

The following result expresses the naturality of the functors $B$ and $A$. Its proof is immediate from the definitions.

**Proposition 6.4.8.** *Let $S$ be a crystalline prism and $f\colon X \to X'$ a morphism of smooth $\overline{S}$-schemes. Then there are commutative diagrams:*

$$\begin{array}{ccc} \mathbb{\Delta}_\phi(X/S) & \xrightarrow{B_X} & \mathbb{\Delta}(X^{(p)}/S) \\ {\scriptstyle \mathbb{\Delta}_\phi(f)}\big\downarrow & & \big\downarrow{\scriptstyle \mathbb{\Delta}(f^{(p)})} \\ \mathbb{\Delta}_\phi(X'/S) & \xrightarrow{B_{X'}} & \mathbb{\Delta}(X'^{(p)}/S) \end{array}$$

$$\begin{array}{ccc} \mathbb{\Delta}(X/S) & \xrightarrow{A_X} & \mathbb{\Delta}_\phi(X/S) \\ {\scriptstyle \mathbb{\Delta}(f)}\big\downarrow & & \big\downarrow{\scriptstyle \mathbb{\Delta}_\phi(f)} \\ \mathbb{\Delta}(X'/S) & \xrightarrow{A_{X'}} & \mathbb{\Delta}_\phi(X'/S). \end{array}$$

□

**Remark 6.4.9.** If $(T, y_T)$ is an $X/S$-$\phi$-prism, then $T_{\mathbb{PD}} \to T$ is a PD-thickening, and thus by forgetting the $\delta$-structure of $T$ we can view $(T, y_T)$ as an object of the site $\mathbb{PD}(X/S)$ consisting of the PD-enlargements of $X/S$. This defines a functor

$$\mathbb{\Delta}_\phi(X/S) \to \mathbb{PD}(X/S). \tag{6.4.7}$$

There are similarly defined functors

$$\mathbb{\Delta}_\phi(X/S) \to \mathbb{D}_\phi(X/S), \quad \mathbb{\Delta}(X/S) \to \mathbb{D}(X/S). \tag{6.4.8}$$

Here $\mathbb{D}_\phi(X/S)$ and $\mathbb{D}(X/S)$ are the sites considered by Oyama and Xu. The objects of $\mathbb{D}_\phi(X/S)$ are pairs $(T, y_T)$, where $T$ is a $p$-torsion free $p$-adic formal scheme and $y_T\colon F_{\overline{T}}(T) \to X$ is a morphism from the scheme theoretic image of $F_{\overline{T}}$ to $X$. Then there is a 2-commutative diagram

$$\begin{array}{ccc} \mathbb{A}_\phi(X/S) & \xrightarrow{B} & \mathbb{A}(X^{(p)}/S) \\ \downarrow & & \downarrow \\ \mathbb{D}_\phi(X/S) & \xrightarrow{C} & \mathbb{D}(X^{(p)}/S), \end{array}$$

where $C$ is the equivalence defined by Oyama and Xu.

**Proposition 6.4.10.** *Suppose that $Y/S$ is a $p$-completely smooth morphism of crystalline prisms.*

*(1) Let $j_Y\colon F(\overline{Y}) \to \overline{Y}$ be the inclusion. Then $(Y, j_Y) \in \mathbb{A}_\phi(\overline{Y}/S)$ covers the final object of the topos $(\overline{Y}/S)_{\mathbb{A}_\phi}$. More generally, if $X$ is closed in $Y$, then $(\mathbb{A}^\phi_X(Y), y_{\mathbb{A}})$ covers the final object of $(X/S)_{\mathbb{A}_\phi}$.*

*(2) If $E$ is a crystal of $\mathcal{O}_{Y/S}$-modules on $\mathbb{A}_\phi(\overline{Y}/S)$, then its value on $(Y, j_Y)$ is canonically endowed with an integrable and quasi-nilpotent nilpotent connection, and this correspondence induces an equivalence from the category of such crystals to the category of $p$-completely quasi-coherent $\mathcal{O}_Y$-modules with integrable quasi-nilpotent connection.*

*Proof.* Suppose that $(T, y_T)$ is an affine object of $\mathbb{A}_\phi(Y/S)$. Then $y_T\colon F(\overline{T}) \to \overline{Y}$ is an $\overline{S}$-morphism, and since $\overline{Y}/\overline{S}$ is smooth and $F(\overline{T}) \to \overline{T}$ is a nil-immersion, $y_T$ can be extended to an $S$-morphism $\overline{T} \to \overline{Y}$. This morphism necessarily takes $F(\overline{T})$ to $F(\overline{Y})$, and hence $y_T$ factors through $F(\overline{Y})$. Now using the formal smoothness of $Y/S$, one can find a lifting $T \to Y$, which, after a $p$-completely faithfully flat cover, can be chosen to be compatible with the Frobenius liftings, thanks to Theorem 1.7.11. The generalization to the case of $X \subseteq Y$ is straightforward, thanks to Proposition 6.4.3.

It follows from statement (1) that the category of crystals of $\mathcal{O}_{Y/S}$-modules on $\mathbb{A}_\phi(Y/S)$ is equivalent to the category of $p$-completely quasi-coherent $\mathcal{O}_Y$-modules endowed with a right action of the groupoid $t, s\colon \mathbb{A}^\phi_Y(1) \to Y$. Corollary 6.4.4 shows that there is an isomorphism of groupoids:

$$\begin{array}{ccc} \mathbb{A}^\phi_Y(1) & \xrightarrow{\cong} & \mathbb{PD}_Y(1) \\ t \downarrow\downarrow s & & t \downarrow\downarrow s \\ Y & \xrightarrow{\text{id}} & Y, \end{array} \tag{6.4.9}$$

and hence an equivalence between the categories of crystals of $\mathcal{O}_{Y/S}$-modules on $\mathbb{PD}(Y/S)$ and on $\mathbb{A}_\phi(Y/S)$. Since crystals on $\mathbb{PD}(Y/S)$ are given by modules with integrable quasi-nilpotent connection, statement (2) follows. ☐

**Remark 6.4.11.** The functor $\mathbb{\Delta}_\phi(X/S) \to \mathbb{PD}(X/S)$ defined in Remark 6.4.9 induces an equivalence on the category of crystals of modules, compatibly with cohomology. Here is a sketch of an argument. If $i\colon X \to Y$ is a closed immersion from $X$ into a $p$-completely smooth $S$-prism $Y/S$, then, as we saw in Proposition 6.4.10, $\mathbb{\Delta}^\phi_X(Y)$ covers the final object of $(X/S)_{\mathbb{\Delta}_\phi}$. Thus crystals on $\mathbb{\Delta}_\phi(X/S)$ correspond to sheaves of $\mathcal{O}_{\mathbb{\Delta}^\phi_X(Y)}$-modules endowed with a compatible action of the groupoid $\mathbb{\Delta}^\phi_Y(1)$. Since $\mathbb{\Delta}^\phi_X(Y) \cong \mathbb{PD}_X(Y)$ and $\mathbb{\Delta}^\phi_Y(1) \cong \mathbb{PD}_Y(1)$, these correspond to crystals on $\mathbb{PD}(X/S)$. A comparison of Cech-Alexander complexes gives the claimed result about cohomology. We leave details to the reader.

We can now explain why $B^*$ corresponds to the F-transform and why $A^*$ corresponds to the $p$-transform. Since $B_{\mathbb{\Delta}}$ is an equivalence, this gives another proof of Shiho's theorem 5.1.4.

**Theorem 6.4.12.** *Suppose that $i\colon X \to Y$ is closed immersion of $X$ into a $p$-completely smooth crystalline $S$-prism $Y/S$.*

1. *If $E$ is a crystal of $\mathcal{O}_{X/S}$-modules on $\mathbb{\Delta}_\phi(X/S)$, its value on $\mathbb{\Delta}^\phi_X(Y)$ is endowed with a canonical integrable and quasi-nilpotent connection.*
2. *If $E'$ is a crystal of $\mathcal{O}_{X^{(p)}/S}$-modules on $\mathbb{\Delta}(X^{(p)}/S)$, then $B^*(E')$ is a crystal of $\mathcal{O}_{X/S}$-modules on $\mathbb{\Delta}_\phi(X/S)$ whose value on $\mathbb{\Delta}^\phi_X(Y) \cong \mathbb{PD}_X(Y)$ is the $F$-transform of the value of $E'$ on $\mathbb{\Delta}_{X^{(p)}}(Y^{(p)})$.*
3. *If $E$ is a crystal of $\mathcal{O}_{Y/S}$-modules on $\mathbb{\Delta}_\phi(Y/S)$, then $A^*(E)$ is a crystal of $\mathcal{O}_{Y/S}$-modules on $\mathbb{\Delta}(Y/S)$, whose value on $Y$ is the $p$-transform (Example 1.2.2) of its value on $Y$.*

*Proof.* The geometry behind this result comes from the following relationships between the functors $A$ and $B$ and the formation of envelopes.

**Lemma 6.4.13.** *In the situation of Theorem 6.4.12, the following statements hold.*

1. *Let $\Phi\colon \mathbb{\Delta}^\phi_X(Y) = \mathbb{\Delta}_{X^\phi}(Y) \longrightarrow \mathbb{\Delta}_{X^{(p)}}(Y^{(p)})$ be the morphism induced by*
$$(F_{Y/S}, f_{X/Y/S})\colon (Y, X^\phi) \to (Y^{(p)}, X^{(p)}),$$
*where $f_{X/YS}$ is the map defined in (6.4.3). Then $\Phi$ fits into a morphism of $X^{(p)}/S$-prisms*
$$\Phi\colon B(\mathbb{\Delta}^\phi_X(Y), y_{\mathbb{\Delta}}) \to (\mathbb{\Delta}_{X^{(p)}}(Y^{(p)}), z_{\mathbb{\Delta}^{(p)}}).$$
2. *The morphism $\Psi\colon \mathbb{\Delta}_X(Y) \to \mathbb{\Delta}^\phi_X(Y)$ induced by the map*
$$(\mathrm{id}, i_{X^\phi})\colon (Y, X) \to (Y, X^\phi)$$
*fits into a morphism of $X/S$-$\phi$-prisms:*
$$\Psi\colon A(\mathbb{\Delta}_X(Y), z_{\mathbb{\Delta}}) \to (\mathbb{\Delta}^\phi_X(Y), y_{\mathbb{\Delta}}).$$

*Proof.* Let $(\mathbb{\Delta}^{\phi}_X(Y), \tilde{z}_{\mathbb{\Delta}}) := B(\mathbb{\Delta}^{\phi}_X(Y), y_{\mathbb{\Delta}})$. To prove that $\Phi$ fits into a morphism of prisms as shown, we must show that $z_{\mathbb{\Delta}^{(p)}} \circ \overline{\Phi} = \tilde{z}_{\mathbb{\Delta}}$, *i.e.*, that the triangle in the diagram below commutes.

$$\begin{array}{ccccccc}
\overline{\mathbb{\Delta}}_{X^{\phi}}(Y) & \xrightarrow{\overline{\Phi}} & \overline{\mathbb{\Delta}}_{X^{(p)}}(Y^{(p)}) & \xrightarrow{\overline{\mathbb{\Delta}}(\pi_{X/S})} & \overline{\mathbb{\Delta}}_X(Y) & \xrightarrow{\overline{\mathbb{\Delta}}(i_{X^{\phi}})} & \overline{\mathbb{\Delta}}_{X^{\phi}}(Y) \\
& \searrow^{\tilde{z}_{\mathbb{\Delta}}} & \downarrow{\scriptstyle z_{\mathbb{\Delta}^{(p)}}} & & \downarrow{\scriptstyle z_{\mathbb{\Delta}}} & & \downarrow{\scriptstyle z_{\mathbb{\Delta}^{\phi}}} \\
& & X^{(p)} & \xrightarrow{\pi_{X/S}} & X & \xrightarrow{i_{X^{\phi}}} & X^{\phi}.
\end{array}$$

Recall from Proposition 6.4.3 that $i \circ y_{\mathbb{\Delta}} = z_{\mathbb{\Delta}^{\phi}} \circ j_{\overline{\mathbb{\Delta}}}$. Since the composed morphism along the top of the diagram is $F_{\overline{\mathbb{\Delta}}}$ and its squares commute, we have:

$$\begin{aligned}
i \circ y_{\mathbb{\Delta}} \circ f_{\overline{\mathbb{\Delta}}} &= z_{\mathbb{\Delta}^{\phi}} \circ j_{\overline{\mathbb{\Delta}}} \circ f_{\overline{\mathbb{\Delta}}} \\
&= z_{\mathbb{\Delta}^{\phi}} \circ F_{\overline{\mathbb{\Delta}}} \\
&= z_{\mathbb{\Delta}^{\phi}} \circ \overline{\mathbb{\Delta}}(i) \circ \overline{\mathbb{\Delta}}(\pi_{X/S}) \circ \overline{\Phi} \\
&= i \circ \pi_{X/S} \circ z_{\mathbb{\Delta}^{(p)}} \circ \overline{\Phi}.
\end{aligned}$$

Since $i$ is a monomorphism, it follows that $y_{\mathbb{\Delta}} \circ f_{\overline{\mathbb{\Delta}}} = \pi_{X/S} \circ z_{\mathbb{\Delta}^{(p)}} \circ \overline{\Phi}$. Since $\tilde{z}_{\mathbb{\Delta}}$ is the unique $S$-morphism such that $\pi_{X/S} \circ \tilde{z}_{\mathbb{\Delta}} = y_{\mathbb{\Delta}} \circ f_{\overline{\mathbb{\Delta}}}$, we conclude that $\tilde{z}_{\mathbb{\Delta}} = z_{\mathbb{\Delta}^{(p)}} \circ \overline{\Phi}$ as claimed.

For (2), we recall that $\mathbb{\Delta}^{\phi}_X(Y) = \mathbb{\Delta}_{X^{\phi}}(Y)$ and note the diagrams:

$$\begin{array}{ccc}
\overline{\mathbb{\Delta}}_X(Y) & \xrightarrow{\overline{\Psi}} & \overline{\mathbb{\Delta}}_{X^{\phi}}(Y) \\
\downarrow{\scriptstyle z_{\mathbb{\Delta}}} & & \downarrow{\scriptstyle y_{\mathbb{\Delta}}} \\
X & \xrightarrow{i_{X^{\phi}}} & X^{\phi}
\end{array}$$

$$\begin{array}{ccccc}
F(\overline{\mathbb{\Delta}}_X(Y)) & \xrightarrow{F(\overline{\Psi})} & F(\overline{\mathbb{\Delta}}_{X^{\phi}}(Y)) & \xrightarrow{\mathrm{id}} & F(\overline{\mathbb{\Delta}}_{X^{\phi}}(Y)) \\
\downarrow{\scriptstyle F(z_{\mathbb{\Delta}})} & & \downarrow{\scriptstyle F(z^{\phi}_{\mathbb{\Delta}})} & & \downarrow{\scriptstyle y_{\mathbb{\Delta}}} \\
F(X) & \longrightarrow & F(X^{\phi}) & \longrightarrow & X.
\end{array}$$

Thus we find a commutative diagram

$$\begin{array}{ccc}
F(\overline{\mathbb{\Delta}}_X(Y)) & \xrightarrow{F(\overline{\Psi})} & F(\overline{\mathbb{\Delta}}_{X^{\phi}}(Y)) \\
\downarrow{\scriptstyle z_{\mathbb{\Delta}} \circ j_{\mathbb{\Delta}}} & & \downarrow{\scriptstyle y_{\mathbb{\Delta}}} \\
X & \xrightarrow{\mathrm{id}} & X,
\end{array}$$

which fills in to define the morphism described in statement (2). $\square$

Suppose that $E'$ is a crystal of $\mathcal{O}_{X^{(p)}/S}$-modules on $\mathbb{\Delta}_{X^{(p)}}(Y^{(p)})$. Then if $(T, z_T)$ is an object of $\mathbb{\Delta}_\phi(X/S)$, the definitions tell us that $B^*(E')_T = E'_{B(T)}$, and the lemma implies that

$$B^*(E')_{\mathbb{\Delta}^\phi_X(Y)} \cong \Phi^*(E'_{\mathbb{\Delta}_{X^{(p)}}(Y^{(p)})}). \tag{6.4.10}$$

Similarly, if $E$ is a crystal of $\mathcal{O}_{X/S}$ crystals on $\mathbb{\Delta}_\phi(X/S)$, we find that

$$A^*(E)_{\mathbb{\Delta}_X(Y)} \cong \Psi^*(E_{\mathbb{\Delta}^\phi_X(Y)}), \tag{6.4.11}$$

where in both cases we use completed pullbacks.

Turning to the proof of the theorem, we first consider the case when $X = \overline{Y}$. Then a crystal of $\mathcal{O}_{X/S}$-modules on $\mathbb{\Delta}_\phi(Y/S)$ is given by a $p$-completely quasi-coherent $\mathcal{O}_Y$-module $E_Y$ endowed with a $\mathbb{\Delta}^\phi_Y(1)$-stratification which, thanks to the isomorphism (6.4.9), is the same as a $\mathbb{PD}_Y(1)$-stratification, which in turn is equivalent to the data of an integrable quasi-nilpotent connection on $E_Y$. This explains (1).

For (2), we obtain from Proposition 6.4.3 and functoriality a commutative diagram of groupoids:

$$\begin{array}{ccccc} \mathbb{\Delta}^\phi_Y(1) & \xrightarrow{\cong} & \mathbb{\Delta}_{Y^\phi}(Y(1)) & \longrightarrow & \mathbb{\Delta}_{Y^{(p)}}(1) \\ & & t \downarrow\downarrow s & & t \downarrow\downarrow s \\ & & Y & \xrightarrow{F_{Y/S}} & Y^{(p)}, \end{array}$$

which in fact we already saw in Proposition 5.6.2. Equation (6.4.10) shows that

$$B^*(E')_Y \cong F^*_{Y/S}(E') \text{ and that } B^*(E')_{\mathbb{\Delta}^\phi_Y(1)} \cong \Phi^*(E')_{\mathbb{\Delta}_{Y^{(p)}}(1)}.$$

Since these isomorphisms are compatible with (completed, as always) pullback by $t$ and $s$, we conclude that the $\mathbb{\Delta}_\phi$-stratification of $B^*(E'_Y)$ is the $\Phi$-pullback of the $\mathbb{\Delta}$-stratification of $E'$. Thanks to Lemma 5.6.7, this shows that the corresponding connection on $F^*_{Y/S}(E)$ is the F-transform of the $p$-connection on $E'$, proving (2) in this case.

The proof of (3) is similar. We have a morphism of groupoids (see again Proposition 5.6.2):

$$\begin{array}{ccc} \mathbb{\Delta}_Y(1) & \xrightarrow{\Psi} & \mathbb{\Delta}^\phi_Y(1) \\ t \downarrow\downarrow s & & t \downarrow\downarrow s \\ Y & \xrightarrow{\text{id}} & Y, \end{array}$$

Equation 6.4.11 implies that

$$A^*(E)_Y \cong E_Y \text{ and that } A^*(E)_{\mathbb{\Delta}_Y(1)} \cong \Psi^*(E_{\mathbb{\Delta}^\phi_Y(1))}).$$

It follows that the $\mathbb{\Delta}$-stratification $\epsilon'$ on $A^*(E)_Y$ is the $\Psi$-pullback of the $\mathbb{\Delta}_\phi$-stratification $\epsilon$ on $E_Y$. Using equation (5.6.2), it is easy to see that the $p$-connection corresponding to $\epsilon'$ is the $p$-transform of connection corresponding to $\epsilon$.

To deal with the case when $i\colon X \to \overline{Y}$ is a general closed immersion, we will use the following compatibility between the functor $B$ and the functors $i_{\mathbb{\Delta}_\phi}\colon \mathbb{\Delta}_\phi(X) \to \mathbb{\Delta}_\phi(\overline{Y})$ and $i^{(p)}_{\mathbb{\Delta}}\colon \mathbb{\Delta}(X^{(p)}) \to \mathbb{\Delta}(\overline{Y}^{(p)})$.

**Lemma 6.4.14.** *In the situation of Theorem 6.4.12, the map (coming from the commutative diagram of Proposition 6.4.8)*

$$B^*_Y(i^{(p)}_{\mathbb{\Delta}*}(E')) \to i_{\mathbb{\Delta}_\phi *}(B^*_X(E'))$$

*is an isomorphism.*

*Proof.* This follows from the fact that $B_X$ and $B_Y$ induce equivalences of topoi, but can also be verified directly by checking the value of both sides on each object of $\mathbb{\Delta}_\phi(X/S)$. $\square$

The general case of the theorem now follows easily. If $E$ is a crystal of $\mathcal{O}_{X/S}$-modules on $\mathbb{\Delta}_\phi(X/S)$, then $i_{\mathbb{\Delta}_\phi *}(E)$ forms a crystal of $\mathcal{O}_{Y/S}$-modules on the small version of $\mathbb{\Delta}_\phi(Y/S)$ whose value on $Y$ is $E_{\mathbb{\Delta}^\phi_X(Y)}$, and the special case we have already discussed shows that this sheaf carries an integrable quasi-nilpotent connection. If $E = B^*_X(E')$ for some crystal $E'$ of $\mathcal{O}_{X^{(p)}/S}$-modules on $\mathbb{\Delta}(X^{(p)}/S)$, then

$$E_{\mathbb{\Delta}^\phi_X(Y)} = i_{\mathbb{\Delta}_\phi *}(E)_Y \cong \Big(B^*_Y(i^{(p)}_{\mathbb{\Delta}*}(E'))\Big)_Y,$$

which we have shown is the F-transform of $i^{(p)}_{\mathbb{\Delta}_\phi *}(E')_{Y^{(p)}} = E'_{\mathbb{\Delta}^\phi_{X^{(p)}}(Y^{(p)})}$. $\square$

**Proposition 6.4.15.** *Let $S$ be a crystalline prism and $X/\overline{S}$ a smooth morphism. There is a 2-commutative diagram of topoi:*

$$\begin{array}{ccc}
(X/S)_{\mathbb{\Delta}_\phi} & \xrightarrow{\;B\;} & (X^{(p)}/S)_{\mathbb{\Delta}} \\
\Big\downarrow{\scriptstyle u_{X/S}} & & \Big\downarrow{\scriptstyle v_{X^{(p)}/S}} \\
X_{\acute{e}t} & \xrightarrow{\;F_{X/S}\;} & X^{(p)}_{\acute{e}t},
\end{array}$$

*where $u_{X/S}\colon (X/S)_{\mathbb{\Delta}_\phi} \to X_{\acute{e}t}$ is defined in the same way as $v_{X/S}$. If $E'$ is a crystal of $\mathcal{O}_{X^{(p)}/S}$-modules on $\mathbb{\Delta}(X^{(p)}/S)$, there is a natural isomorphism*

$$Rv_{X^{(p)}/S*}(E') \to F_{X/S*}Ru_{X/S*}(B^*(E')).$$

*Proof.* The existence of the diagram follows immediately from the definitions. Since $F_{X/S*}$ is exact, we find an isomorphism:

$$Rv_{X^{(p)}/S*} \circ RB_*(B^*(E')) \cong F_{X/S*} \circ Ru_{X/S*}(B^*(E')).$$

Since $B_{\mathbb{\Delta}}$ is an equivalence of topoi, the functor $B_*$ is exact, so $RB_*B^*(E') \cong B_*B^*(E') \cong E'$, and the result is proved. □

Our next result shows that the topos-theoretic formulation of the F-transform is compatible with derived direct images.

**Corollary 6.4.16.** *Let $S$ be a crystalline prism, let $f\colon X \to X'$ be a morphism of smooth $\overline{S}$-schemes, and let $\mathbb{\Delta}_{s\phi}(X'/S)$ denote the small version of $\mathbb{\Delta}_{\phi}(X'/S)$.*

*(1) If $E$ is a crystal on $\mathbb{\Delta}_{\phi}(X/S)$, and $i \geq 0$, the natural map*

$$R^i f^{(p)}_{\mathbb{\Delta}*}(B_{X*}(E)) \to B_{X'*}(R^i f_{\mathbb{\Delta}_{\phi}*}(E))$$

*is an isomorphism, and if $E'$ is a crystal of $\mathcal{O}_{X^{(p)}/S}$-modules on $\mathbb{\Delta}(X^{(p)}/S)$, the natural map*

$$B^*_{X'}(R^i f^{(p)}_{\mathbb{\Delta}*}(E')) \to R^i f_{\mathbb{\Delta}_{\phi}*}(B^*_X(E'))$$

*is an isomorphism.*

*(2) Suppose that the restriction of $R^i f^{(p)}_{\mathbb{\Delta}*}(E')$ to $\mathbb{\Delta}_s(X'^{(p)}/S)$ forms a crystal of $\mathcal{O}_{X'^{(p)}/S}$-modules. Then $R^i f_{\mathbb{\Delta}_{\phi}*}(B^*_X(E'))$ forms a crystal of $\mathcal{O}_{X'/S}$ modules on $\mathbb{\Delta}_{s\phi}(X'/S)$, and if $Y'/S$ is a crystalline prism lifting $X'/\overline{S}$, its value on $Y'$, endowed with its corresponding connection, is the F-transform of $R^i f^{(p)}_{\mathbb{\Delta}*}(E')_{Y'^{(p)}}$.*

*Proof.* Theorem 6.4.6 tells us that the morphisms $B_X$ and $B_{X'}$ induce equivalences of topoi. It follows that $B_{X*}$ and $B_{X'*}$ are exact. Then if $E$ is a crystal on $\mathbb{\Delta}_{\phi}(X/S)$:

$$R^i f^{(p)}_{\mathbb{\Delta}*}(B_{X*}(E)) \cong R^i(f^{(p)}_{\mathbb{\Delta}} \circ B_X)_*(E') \cong R^i(B_{X'*} \circ f_{\mathbb{\Delta}_{\phi}*})(E) \cong B_{X'*}(R^i f_{\mathbb{\Delta}_{\phi}*}(E)).$$

This proves the first part of statement (1). Now let $E = B^*_X(E')$. Then $B_{X*}(E) = E'$, and applying $B_{X'*}$ to both sides of the previous equation yields the second part. If $R^i f^{(p)}_{\mathbb{\Delta}*}(E')$ forms a crystal on $\mathbb{\Delta}_s(X'^{(p)}/S)$, then $B^*_{X'}(R^i f^{(p)}_{\mathbb{\Delta}*}(E'))$ forms a crystal on $\mathbb{\Delta}_{s\phi}(X'/S)$, and, by statement (1), the same is true of $R^i f_{\mathbb{\Delta}_{\phi}*}(B^*_X(E'))$. Statement (2) follows, thanks to Theorem 6.4.12 □

## 6.5 Applications

In this section we explain a few applications of our comparison theorems. These are not really new, but our results allow us to present formulations that are more explicit than has been possible previously. Throughout this section we assume that $S$ is a crystalline prism and that $\overline{f}\colon X \to \overline{S}$ is a smooth morphism of schemes. We denote by $f\colon X \to S$ the composition of $\overline{f}$ with the inclusion $\overline{S} \to S$, and we refer to diagram (1.1.3) for the notation surrounding the relative Frobenius morphism $F_{X/S}\colon X \to X^{(p)}$. Recall that we have canonical morphisms:

$$v_{X/S}\colon (X/S)_{\Delta} \to X_{\acute{e}t} \text{ and } u_{X/S}\colon (X/S)_{\mathbb{PD}} \to X_{\acute{e}t}.$$

Compatibility of prismatic cohomology with base change, based on its computation by Čech-Alexander complexes, is discussed in [10, 4.18]. The use of $p$-de Rham complexes provides stronger results. This has been carried out by Y. Tian [39] who also deals with more general prisms. Here is our (slightly different) version.

**Theorem 6.5.1.** *If $E$ is a crystal of $\mathcal{O}_{X/S}$-modules on $\Delta(X/S)$, the following statements hold.*

*(1) $R^q v_{X/S*}(E) = 0$ if $q > \dim(X/\overline{S})$.*

*(2) If $E$ is locally free of finite rank, then formation of $\mathrm{R}v_{X/S*}(E)$ commutes with base change $S' \to S$.*

*(3) If $X/\overline{S}$ is proper of relative dimension $n$, if $E$ is locally free of finite rank, and if $S = \mathrm{Spf}\, A$, then $\mathrm{R}f_{\Delta*}E$ is a perfect complex of $p$-adic $A$-modules, of amplitude in $[0, 2n]$.*

*Proof.* The first statement may be verified Zariski locally on $X$, so we may without loss of generality assume that there is a crystalline prism $Y/S$ lifting $X/\overline{S}$. Then $\mathrm{R}v_{X/S*}(E)$ is calculated by the $p$-de Rham complex of $(E_Y, \nabla')$, which has no terms in degrees greater than $n$. Statement (2) is proved in a similar manner, using the fact that formation of the $p$-de Rham complex is compatible with base change, and that its terms are locally free. The last statement also follows; see [39, 2.10] for details. □

A key result of Bhatt-Scholze, [10, 5.3], compares prismatic and crystalline cohomology. It asserts the existence of a canonical quasi-isomorphism:

$$\mathrm{L}F_S^*\mathrm{R}f_{\Delta*}(\mathcal{O}_{X/S}) \to \mathrm{R}f_{\mathbb{PD}*}(\mathcal{O}_{X/S}), \tag{6.5.1}$$

where, as usual, all pullbacks are completed. Strictly speaking, this result is only stated and proved in [10] when $X$ and $S$ are affine; presumably the general case follows by standard simplicial methods. The following result provides a generalization to the case of crystals.

**Theorem 6.5.2.** *Let $E'$ be a crystal of $\mathcal{O}_{X/S}$-modules on $\Delta(X/S)$ and let $E$ be the crystal of $\mathcal{O}_{X/S}$-modules on $\mathbb{PD}(X/S)$ corresponding to the F-transform (5.1.9) of $\pi^*_{X/S}(E')$.*

*(1) There is a canonical strict quasi-isomorphism*

$$\mathrm{R}v_{X^{(p)}/S*}(\pi^*_{X/S}(E')) \to F_{X/S*}Ru_{X/S*}(E).$$

*(2) If $E'$ is locally free or if $F_S$ is p-completely flat, there is a canonical strict quasi-isomorphism*

$$F_S^*(Rf_{\Delta *}E') \to \mathrm{R}f_{\mathbb{PD}*}E.$$

*(3) In particular, there is a canonical strict quasi-isomorphism:*

$$F_S^*(\mathrm{R}f_{\Delta *}\mathcal{O}_{X/S}) \to \mathrm{R}f_{\mathbb{PD}*}\mathcal{O}_{X/S}.$$

*Proof.* Statement (1) here is just a restatement of (2) of Proposition 6.4.15. The remaining statements then follow from the base change results of Theorem 6.5.1. □

Next we show that Frobenius induces an isogeny on prismatic cohomology. We can again include a version with coefficients in a prismatic crystal.

**Proposition 6.5.3.** *Let $E'$ be a crystal of $\mathcal{O}_{X^{(p)}/S}$-modules on $\Delta(X^{(p)}/S)$ and let $F_{X/S}^*(E')$ be its pullback to a crystal of $\mathcal{O}_{X/S}$-modules on $\Delta(X/S)$. If $X$ is quasi-compact, the natural map*

$$Rv_{X^{(p)}/S*}(E') \to F_{X/S*}Rv_{X/S*}(F_{X/S}^*(E'))$$

*is an isogeny. In particular, the Frobenius morphism*

$$F_{X/S_{\Delta}}^*\colon R\Gamma((X^{(p)}/S)_{\Delta}, \mathcal{O}_{X^{(p)}/S}) \to R\Gamma((X/S)_{\Delta}, \mathcal{O}_{X/S})$$

*is an isogeny.*

*Proof.* We explain this first assuming that $X$ admits an embedding in a $p$-completely smooth crystalline $S$-prism $Y/S$. Let $E'_{Y^{(p)}}$ denote the value of $E'$ on the prismatic envelope of $X^{(p)}$ in $Y^{(p)}$. Then the morphism in the statement is represented by the morphism of complexes:

$$(E'_{Y^{(p)}} \otimes \Omega^{\bullet}_{Y^{(p)}/S}, d') \to F_{Y/S*}(F_{Y/S}^*(E'_{Y^{(p)}}) \otimes \Omega^{\bullet}_{Y/S}, d').$$

This is the morphism $c$ in Proposition 5.1.7, which is the composite of the quasi-isomorphism $\zeta_E$ and the morphism $b$, which is an isogeny when $X$ is quasi-compact. It follows that the map in question becomes an isogeny in the derived category. The general case will follow by the usual simplicial methods. □

The following result is stated in [10, 4.15] and [1, 3.2.1]. The treatments there only explicitly discuss the affine case, and the task of globalization (presumably using simplicial methods) is left to the reader. For quasi-projective schemes, our methods lead to a very explicit construction.

**Theorem 6.5.4.** *Let $S$ be a crystalline prism, let $X/\overline{S}$ be a smooth morphism, together with a closed immersion $X \to Y$, where $Y/S$ is a p-completely smooth crystalline $S$-prism, and let $\mathbb{L}_{X/S}$ denote the (completed) cotangent complex of $X$ relative to $S$.*

*(1) There are natural quasi-isomorphisms:*

$$\mathbb{L}_{X/S} \to \tau^{\leq 0}\left(\mathcal{O}_{\overline{\mathbb{A}}_X(Y)} \otimes \Omega^\bullet_{Y/S}[1], \overline{d}'\right)$$

$$\mathbb{L}_{X/S} \to \tau^{\leq 0}\left(\mathcal{O}_{\overline{\mathbb{A}}_X(Y)} \otimes \Omega^\bullet_{Y/S}[1], \overline{\nabla}'\right),$$

*compatible with the isomorphism (5.4.3) of the derived images of the corresponding Higgs complexes.*

*(2) A mod $p^2$ flat lifting $\tilde{X} \to Y$ of $X \to Y$ induces a splitting of $\mathbb{L}_{X/S}$ compatible with the Sen operator on $\left(\mathcal{O}_{\overline{\mathbb{A}}_X(Y)} \otimes \Omega^\bullet_{Y/S}[1], \overline{\nabla}'\right)$.*

*Proof.* Strictly speaking, our formulation does not make sense without referring to an explicit complex representing $\mathbb{L}_{X/S}$, and we shall indeed realize the morphisms in the theorem as explicit morphisms of complexes. In fact, the closed immersion $X \to \overline{Y} \to Y$, which is necessarily a regular immersion, provides us with such a representation: by [3, 2.1], the cotangent complex of $X/S$ can be identified with the complex:

$$\mathbb{L}_{X/Y/S} := \mathcal{C}_{X/Y} \xrightarrow{\overline{d}} \Omega^1_{Y/S} \otimes \mathcal{O}_X,$$

placed in degrees $-1$ and 0. Here $\overline{d}$ is the map induced by the Kähler differential $d\colon I_{X/Y} \to \Omega^1_{Y/S}$. It follows from the theory of cotangent complexes that this construction is functorial in the derived sense. (In fact, it is elementary to check directly that if $g\colon Y \to Z$ is a morphism of smooth $S$-schemes such that $g \circ i$ is also a closed immersion, then the induced map of complexes $\mathbb{L}_{X/Z/S} \to \mathbb{L}_{X/Y/S}$ is a quasi-isomorphism.)

Recall from (2.2.1) that we have an $\mathcal{O}_Y$-linear map $\rho\colon I_{X/Y} \to \mathcal{O}_{\mathbb{A}_X(Y)}$, with $p\rho(x) = \pi^*_{\mathbb{A}}(x)$ for all $x \in I_{X/Y}$. Thus $\rho$ induces a morphism

$$\overline{\rho}\colon \mathcal{C}_{X/Y} \to \mathcal{O}_{\overline{\mathbb{A}}_X(Y)}.$$

Furthermore,

$$pd'\rho(x) = d'\pi^*_{\mathbb{A}}(x) = pd\pi^*_{\mathbb{A}}(x) = p\pi^*_{\mathbb{A}}(dx),$$

so $d'\rho(x) = \pi^*_{\mathbb{A}}(dx)$. We see that there is a commutative diagram of $\mathcal{O}_X$-linear maps:

$$\begin{array}{ccccc}
\mathcal{C}_{X/Y} & \xrightarrow{\overline{d}} & \Omega^1_{Y/S} \otimes \mathcal{O}_X & \xrightarrow{\overline{d}} & \Omega^2_{Y/S} \otimes \mathcal{O}_X \\
\downarrow{\scriptstyle \overline{\rho}} & & \downarrow{\scriptstyle \pi^*_{\mathbb{A}}} & & \downarrow{\scriptstyle p\pi_{\mathbb{A}^*}} \\
\mathcal{O}_{\overline{\mathbb{A}}_X(Y)} & \xrightarrow{\overline{d}'} & \mathcal{O}_{\overline{\mathbb{A}}_X(Y)} \otimes \Omega^1_{Y/S} & \xrightarrow{\overline{d}'} & \mathcal{O}_{\overline{\mathbb{A}}_X(Y)} \otimes \Omega^2_{Y/S}
\end{array}$$

Since the rightmost vertical arrow is in fact zero, the diagram defines the first morphism of complexes as stated in the theorem. In fact, since $\overline{d}'$ and $\overline{\nabla}'$ agree on the image of $\mathcal{O}_{\overline{\mathbb{D}}_X(Y)}$ in $\mathcal{O}_{\overline{\mathbb{A}}_X(Y)}$, we also get the second morphism of complexes. It follows that the map $\alpha$ (5.2.3) acts as the identity on this image and hence from the construction of the derived isomorphism (5.4.3) in Proposition 5.2.3 that it is compatible with the second complex as well.

Suppose that $X \to Y'$ is another embedding of $X$ into a $p$-completely smooth crystalline $S$-prism. Let $Z := Y \times Y'$ and let $X \to Z$ be the map induced by these two embeddings. We find a commutative diagram of complexes:

$$\begin{array}{ccc}
\mathbb{L}_{X/Y/S} & \longrightarrow & \tau^{\leq 0}(\mathcal{O}_{\overline{\mathbb{A}}_X(Y)} \otimes \Omega^\bullet_{Y/S}[1], \overline{d}') \\
\downarrow & & \downarrow \\
\mathbb{L}_{X/Z/S} & \longrightarrow & \tau^{\leq 0}(\mathcal{O}_{\overline{\mathbb{A}}_X(Z)} \otimes \Omega^\bullet_{Z/S}[1], \overline{d}') \\
\uparrow & & \uparrow \\
\mathbb{L}_{X/Y'/S} & \longrightarrow & \tau^{\leq 0}(\mathcal{O}_{\overline{\mathbb{A}}_X(Y')} \otimes \Omega^\bullet_{Y'/S}[1], \nabla').
\end{array}$$

The vertical maps on the right are quasi-isomorphisms by the prismatic Poincaré lemma (Theorem 4.3.1), and the maps on the left are quasi-isomorphisms as we have already explained. We can conclude that our construction gives a well-defined morphism in the derived category, as stated in the theorem. To show that it is a quasi-isomorphism, we may work locally, and in particular we may assume that $Y$ is a lifting of $X$. In that case, the complex $\mathbb{L}_{X/Y/S}$ is $(p) \otimes \mathcal{O}_X$ in degree $-1$ and $\Omega^1_{X/S}$ in degree zero, with $\overline{d} = 0$, while the complex $\tau^{\leq 0}(\mathcal{O}_{\overline{\mathbb{A}}_{X(Y)}} \otimes \Omega^\bullet_{Y/S}[1])$ is $\mathcal{O}_X$ in degree $-1$ and $\Omega^1_{X/S}$ in degree zero, with $\overline{d}' = 0$. Furthermore, $\rho$ sends the class of $p$ to 1, and we see that the map of complexes is the identity map in this case.

As we have seen in Proposition 1.7.1, a lifting $\tilde{X} \to Y$ of $X \to Y$ provides a splitting of the exact sequence

$$0 \to \mathcal{O}_X \to \mathcal{C}_{X/Y} \to \mathcal{C}_{X/\overline{Y}} \to 0.$$

By construction, the Sen operators $\theta_{\tilde{X}}$ and $\Theta_{\tilde{X}}$ are compatible with this splitting: they act as 0 on $\mathcal{O}_X$ and as the identity on the image of $\mathcal{C}_{X/\overline{Y}}$. Now the

diagram

$$
\begin{array}{ccccc}
\mathbb{L}_{\overline{Y}/Y/S_{|X}} & = & \mathcal{O}_X & & \\
\downarrow & & \downarrow & & \\
\mathbb{L}_{X/Y/S} & = & \mathcal{C}_{X/Y} & \longrightarrow & \Omega^1_{Y/S_{|X}} \\
\downarrow & & \downarrow & & \downarrow \scriptstyle{\mathrm{id}} \\
\mathbb{L}_{X/\overline{Y}/S} & = & \mathcal{C}_{X/\overline{Y}} & \longrightarrow & \Omega^1_{Y/S_{|X}}
\end{array}
$$

shows that a splitting of the vertical sequence on the left is equivalent splitting the exact sequence of complexes provided by the diagram. Thus, the splitting of $\mathbb{L}_{X/Y/S}$ is compatible with the action of the Sen operator on $(\mathcal{O}_{\overline{\mathbb{\Delta}}_{X}(Y')} \otimes \Omega^{\bullet}_{Y'/S}[1], \nabla')$. $\square$

# 7 Appendix: $p$-adic sheaves and modules

Let $A$ be a $p$-adically separated and complete ring, not necessarily noetherian. We will need to work with $A$-modules which are also $p$-adically complete but not necessarily finitely generated. There are several possible strategies: "derived complete modules" as suggested in [10], topological modules, filtered modules, or $p$-adic inverse systems as in [16]. All these categories are additive, but only the first is abelian, which makes it a tempting choice, and is the strategy followed in [1]. However, this reference only studies derived $p$-complete modules which have bounded $p$-torsion, and as it turns out, these are necessarily $p$-adically separated and complete. Furthermore, we have been unable to prove the key exactness properties in Proposition 7.3.3 for derived $p$-complete modules in general. Consequently we will work instead with the category of $p$-adically separated and complete modules. This category can be profitably viewed as living either in the category of all modules, or in the category of inverse systems of torsion modules, each of which has its advantages and disadvantages. The following exposition, which is perhaps overly detailed and which includes many well-known results, summarizes our approach.

## 7.1 $p$-adic modules

**Definition 7.1.1.** *Let $A$ be a $p$-adically separated and complete ring. We write $\mathbf{M}_A$ for the category of all $A$-modules, and $\mathbf{M}_{\bullet}$ for the category of inverse systems $\{M_n \to M_{n-1} : n \in \mathbf{N}\}$ of $A$-modules such that $p^n M_n = 0$ for all $n$. A $p$-adic $A$-module is an object of $\mathbf{M}_{\bullet}$ such that each map $M_n \otimes_A A_{n-1} \to M_{n-1}$ is an isomorphism.*

We have functors

$$\begin{array}{rcll} \varprojlim : \mathbf{M}_\bullet & \to & \mathbf{M}_A : & M_\bullet \mapsto \hat{M} \\ S_\bullet : \mathbf{M}_A & \to & \mathbf{M}_\bullet : & M \mapsto \{M/p^n M \to M/p^{n-1}M\}. \end{array} \qquad (7.1.1)$$

For the following result, see also Temkin's [16, 2.3.8].

**Proposition 7.1.2.** *With the definitions above, the following statements hold.*

*(1) For every $M \in \mathbf{M}_A$, the object $S_\bullet(M) \in \mathbf{M}_\bullet$ is a $p$-adic $A$-module. For every $M_\bullet$ in $\mathbf{M}_\bullet$, the object $\varprojlim M_\bullet$ of $\mathbf{M}_A$ is $p$-adically separated and complete.*

*(2) If $M_\bullet$ is a $p$-adic $A$-module and $\hat{M} := \varprojlim M_\bullet$, then for each $n$, the natural map $\hat{M} \otimes_A A_n \to M_n$ is an isomorphism.*

*(3) The functor $\varprojlim$ induces an equivalence from the category of $p$-adic $A$-modules to the category of $p$-adically separated and complete $A$-modules.*

*Proof.* The first part of statement (1) is obvious. For the second part, let $\hat{M} := \varprojlim M_\bullet$, and for each $n$, let $F^n\hat{M}$ be the kernel of the map $\hat{M} \to M_n$. It follows from the construction that $\hat{M}$ is separated and complete for the $F$-adic topology. Since $p^n M_n = 0$ for every $n$, this topology is weaker than the $p$-adic topology, and the following lemma shows that $\hat{M}$ is also $p$-adically separated and complete. This is proved in [30, A1] for any finitely generated ideal in place of $(p)$; since the proof in this case is considerably simpler, we explain it again here.

**Lemma 7.1.3.** *Let $M$ be an $A$-module, separated and complete with respect to an $A$-linear topology $F^\bullet$ such that $p^n M \subseteq F^n M$ for all $n$. Then $M$ is also $p$-adically separated and complete.*

*Proof.* The separation is automatic, and for the completeness we check that every $p$-adically Cauchy sequence $(x_\bullet)$ in $M$ has a $p$-adic limit. Since $(x_\bullet)$ is $p$-adically Cauchy, it suffices to find a subsequence which converges $p$-adically, so we may assume that $x_{n+1} - x_n \in p^n M$ for all $n$. Thus it suffices to prove that every series of the form $\sum p^n y_n$ converges $p$-adically. Since this series is also $F$-adically Cauchy, it has an $F$-adic limit $y$, and we claim that this is also the $p$-adic limit. Fix $N$, and consider the series $\sum_{n \geq N} p^{n-N} y_n$. This series is also $F$-adically Cauchy and therefore has an $F$-adic limit $z_N$, and it is clear that

$$p^N z_N = \sum_{n \geq N} p^n y^n = y - \sum_{n<N} p^n y^n.$$

Thus, the sequence of partial sums converges $p$-adically to $y$. □

To prove (2), we repeat the argument of [16, 2.3.8] for the convenience of the reader. See also [2, 05GG] for a more general statement. Since the transition maps $M_n \to M_{n-1}$ are all surjective, so is each projection map $\hat{M} \to M_n$;

let $F^n M$ denote its kernel. By construction, the module $\hat{M}$ is complete with respect to the $F$-adic topology. We know that $p^n\hat{M} \subseteq F^n\hat{M}$, and we must show that the reverse inclusion also holds. We have

$$\hat{M}/F^n\hat{M} \cong M_n \cong M_{n+1}/p^n M_{n+1} \cong \hat{M}/(F^{n+1}\hat{M} + p^n\hat{M}),$$

and it follows that that $F^n\hat{M} = p^n\hat{M} + F^{n+1}\hat{M}$. An element $x$ of $F^n\hat{M}$ can then be written as $p^n x'_0 + x_1$ with $x_0 \in \hat{M}$ and $x_1 \in F^{n+1}\hat{M}$, and then we can write $x_1 = p^{n+1}x'_1 + x_2$ with $x_2 \in F^{n+2}\hat{M}$. Continuing in this way, we find sequences $(x'_0, x'_1, x'_2, \ldots)$ and $(x_1, x_2, \ldots)$, with $x'_i \in \hat{M}$, $x_i \in F^{n+i}\hat{M}$ and $x_i = p^{n+i}x'_i + x_{i+1}$ for all $i$. Let $x'$ be the limit of the series $\sum p^i x'_i$, which converges $F$-adically. Then $x = p^n x'$, since $\hat{M}$ is $F$-adically separated.

If $M$ is $p$-adically separated and complete, then the natural map $M \to \varprojlim S_\bullet(M)$ is an isomorphism, by definition. On the other hand, if $M_\bullet$ is a $p$-adic $A$-module, statements (1) and (2) tell us that $\varprojlim M_\bullet$ is $p$-adically separated and complete and that the natural map of inverse systems $S_\bullet(\varprojlim M_\bullet) \to M_\bullet$ is an isomorphism. Statement (3) follows. □

## 7.2 $p$-adic exactness

**Example 7.2.1** The category of $p$-adically separated and complete $A$-modules is not abelian: a proper inclusion with dense image is a monomorphism and an epimorphism but not an isomorphism. For an example, let $A$ be $\mathbf{Z}_p$, let $M$ be the $\mathbf{Z}_p$-module consisting of the set of sequences $(a_0, a_1, \ldots)$ in $\mathbf{Z}_p$ which tend to zero $p$-adically, and let $N$ be the set of all sequences in $\mathbf{Z}_p$. These module are $p$-adically separated and complete, as is easy to verify. Define a map $\alpha\colon M \to M$ by $\alpha(a_0, a_1, a_2, \ldots) := (a_0, pa_1, p^2a_2, \ldots)$ and a map $\beta\colon N \to M$ by the same formula. We have a commutative diagram with exact rows:

$$\begin{array}{ccccccccc}
0 & \longrightarrow & M & \xrightarrow{\alpha} & M & \longrightarrow & M'' & \longrightarrow & 0 \\
 & & \downarrow{\scriptstyle\gamma} & & \downarrow{\scriptstyle\mathrm{id}} & & \downarrow{\scriptstyle\delta} & & \\
0 & \longrightarrow & N & \xrightarrow{\beta} & M & \longrightarrow & N'' & \longrightarrow & 0,
\end{array}$$

where $\gamma$ is the inclusion. The image of $\beta$ is the set of sequences $c_\bullet$ such that each $c_i$ is divisible by $p^i$, and is evidently closed in $M$. If $b_\bullet \in N$ and $n > 0$, then

$$\begin{aligned}
\beta(b_\bullet) &= (b_0, pb_1, p^2b_2, \ldots) \\
&= (b_0, pb_1, \ldots, p^{n-1}b_{n-1}, 0, \ldots) + p^n(0, 0, \ldots, b_n, pb_{n+1}, \ldots)
\end{aligned}$$

which is of the form $\alpha(a_\bullet) + p^n a'_\bullet$, with $a_\bullet$ and $a'_\bullet \in M$. Thus $\alpha(M)$ is dense in $\beta(M)$, and in fact $\beta(M)$ is the closure of $\alpha(M)$ in $M$.

Note that $M''$ is the standard example of a module which is derived $p$-complete but not $p$-adically separated [2, 0G3F]. We also note that $M''$ can be written as an extension of two $p$-adically separated and complete modules, showing that the class of such modules does not form a thick subcategory of the category of all modules. To see this, we check first that the image of $\gamma$ is closed in $N$. If $b_\bullet \in N$ is in the closure of $\gamma(M)$, then for every $n > 0$, there exist $a_\bullet \in M$ and $b'_\bullet \in N$ such that $b_\bullet = \gamma(a_\bullet) + p^n b'_\bullet$. Then there exists $m_n$ such that $a_i \in p^n \mathbf{Z}_p$ for all $i \geq m_n$, which implies that the same holds for $b_i$. Thus $b_\bullet$ also belongs to the image of $\gamma$. The snake lemma implies that $\mathrm{Ker}(\delta) \cong \mathrm{Cok}(\gamma)$. Since the image of $\gamma$ is closed, these modules are $p$-adically separated, and since $\beta(N)$ is closed in $N$, the quotient $N''$ is also $p$-adically separated. Then

$$0 \to \mathrm{Ker}(\delta) \to M'' \to N'' \to 0$$

is the desired extension of separated and complete modules which is derived complete but not separated. We note further that $\mathrm{Ker}(\delta) = \cap_n p^n M''$ and that $N''$ is the $p$-adic completion of $M''$. Thus, by the exact sequence in [2, 0BKG], $\mathrm{Ker}(\delta) = R^1 \lim\{M''[p^n] : n \in \mathbf{N}\}$. Note also that this module is $p$-torsion free.

Since the category of $p$-adic $A$-modules is not abelian, it is not safe to speak about exact sequences and the cohomology of complexes. We therefore introduce the following terminology.

**Definition 7.2.2.** *Let $A$ be a $p$-adically separated and complete ring.*

*(1) A morphism $M'_\bullet \to M_\bullet$ of $p$-adic $A$-modules is a strict monomorphism if each $M'_n \to M_n$ is a monomorphism.*

*(2) A sequence of $p$-adic $A$-modules $M'_\bullet \to M_\bullet \to M''_\bullet$ is strictly short exact if each sequence $M'_n \to M_n \to M''_n$ is short exact.*

*(3) A complex $M^\bullet_\bullet$ of $p$-adic $A$-modules is strictly acyclic if each $M^\bullet_n$ is acyclic.*

*(4) A morphism $u_\bullet$ of complexes of $p$-adic $A$-modules is a strict quasi-isomorphism if each $u_n$ is a quasi-isomorphism.*

**Remark 7.2.3.** The category $\mathbf{M}_\bullet$ is an abelian category, and a sequence in this category is exact if and only if the corresponding sequences for each $n$ are all exact. The category of $p$-adic $A$-modules is a full subcategory of $\mathbf{M}_\bullet$, and we claim that it is closed under extension. Indeed, if $0 \to M'_\bullet \to M_\bullet \to M''_\bullet \to 0$ is exact in $\mathbf{M}_\bullet$ and $M'_\bullet$ and $M''_\bullet$ are $p$-adic systems, then for each $n$ we have a commutative diagram

$$\begin{array}{ccccccccc}
 & & M'_n \otimes A_{n-1} & \longrightarrow & M_n \otimes A_{n-1} & \longrightarrow & M''_n \otimes A_{n-1} & \longrightarrow & 0 \\
 & & \downarrow & & \downarrow & & \downarrow & & \\
0 & \longrightarrow & M'_{n-1} & \longrightarrow & M_{n-1} & \longrightarrow & M''_{n-1} & \longrightarrow & 0
\end{array}$$

whose rows are exact. The vertical arrows on the left and right are isomorphisms, and it follows that the central one is also. Thus our definition of "strict exact sequences" is consistent with Quillen's notion of an "exact category." We have two functors from the category of $p$-adic $A$-modules to an abelian category: the inclusion into the category $\mathbf{M}_\bullet$ and the inverse limit functor into the category $\mathbf{M}_A$ (or into the category of derived $p$-complete modules). By definition, a sequence $0 \to M'_\bullet \to M_\bullet \to M''_\bullet \to 0$ of $p$-adic $A$-modules is (strictly) exact if and only if it is so when viewed in $\mathbf{M}_\bullet$. The inverse limit of such a sequence is again exact, but converse does not hold.

**Lemma 7.2.4.** *Let $M$ be a $p$-adically separated and complete $A$-module and let $K$ be a submodule.*

*(1) If $K$ is closed in the $p$-adic topology of $M$, then it is $p$-adically separated and complete.*

*(2) Let $K^-$ be the closure of $K$ in $M$. Then $M/K^-$ is the $p$-adic completion of $M/K$.*

*Proof.* The topology of $K$ induced by the $p$-adic topology of $M$ could be weaker than the $p$-adic topology of $K$ but it follows nonetheless that $K$ is $p$-adically separated and complete, by Lemma 7.1.3. We have an inverse system of short exact sequences:

$$0 \to K/(p^n M \cap K) \to M/p^n M \to M/(p^n M + K) \to 0.$$

Since the transition maps on the left are surjective, we find an exact sequence

$$0 \to \varprojlim K/(p^n M \cap K) \to \varprojlim M/p^n M \to \varprojlim M/(p^n M + K) \to 0.$$

Since $K^- = \varprojlim K/(p^n M \cap K)$, the lemma is proved. □

**Proposition 7.2.5.** *Let $A$ be a $p$-adically separated and complete ring.*

*(1) The category of $p$-adically separated and complete $A$-modules admits kernels, whose formation is compatible with the inclusion in the category $\mathbf{M}_A$.*

*(2) Suppose that the kernel $K$ of a homomorphism $u\colon M \to N$ of $p$-adically separated and complete $A$-modules agrees with the kernel of the corresponding homomorphism of $u_\bullet\colon M_\bullet \to N_\bullet$ in $\mathbf{M}_\bullet$. Then*

*(a) The inclusion $K \to M$ is a strict monomorphism.*

*(b) The image $M''$ of $M$ in $N$ is $p$-adically separated and complete and is closed in $N$, and the inclusion $M'' \to N$ is a strict monomorphism.*

*A homomorphism of $p$-adically separated and complete $A$-modules is said to be* strict *if it satisfies these properties.*

*Proof.* If $u\colon M \to N$ is a homomorphism of $p$-adically separated and complete $A$-modules, let $K := \{x \in M : u(x) = 0\}$. This module is evidently closed in the $p$-adic topology of $M$, and so it is $p$-adically separated and complete by Lemma 7.2.4. As an example, if $u$ is multiplication by $p$ on $A := \mathbf{Z}_p$, we see that $K = \{0\}$, although $u_n$ is not injective for $n > 0$.

Now suppose that $K_n = \mathrm{Ker}(u_n)$ for all $n$. Then necessarily each $K_n \to M_n$ is injective, so $K \to M$ is a strict monomorphism. For each $n$, we have a commutative diagram

$$\begin{array}{ccccccccc} & & K \otimes_A A_n & \longrightarrow & M \otimes_A A_n & \longrightarrow & M'' \otimes_A A_n & \longrightarrow & 0 \\ & & \downarrow & & \downarrow & & \downarrow & & \\ 0 & \longrightarrow & K_n & \longrightarrow & M_n & \longrightarrow & \mathrm{Im}(u_n) & \longrightarrow & 0 \end{array}$$

with exact rows. Since the first two vertical arrows are isomorphisms, so is the third, and the left horizontal top arrow is injective. Since the transition maps for the system $K_\bullet$ are surjective, we can deduce that the map $M = \varprojlim M_\bullet \to \varprojlim M''_\bullet$ is surjective, and since $M''_n = \mathrm{Im}(u_n) \subseteq N_n$, that the map $\varprojlim M''_\bullet \to \varprojlim N_\bullet = N$ is injective. Thus $M'' = \hat{M}''$ is $p$-adically separated and complete. Moreover, since $M''_n \subseteq N_n$ for all $n$, the $p$-adic topology of $M''$ agrees with the induced topology from that of $N$, and the map $M'' \to N$ is a strict monomorphism. It follows also that $M''$ is closed in the $p$-adic topology of $N$. □

**Remark 7.2.6.** We should remark that the category of $p$-adic $A$-modules also admits cokernels, whose formation is compatible with the inclusion into the category $\mathbf{M}_\bullet$. but not into the category $\mathbf{M}_A$. Namely, if $u_\bullet \colon M_\bullet \to N_\bullet$ is a homomorphism of $p$-adic $A$-modules, then the inverse system $\mathrm{Cok}(u_\bullet)$ is again a $p$-adic $A$-module, but its inverse limit need not be the cokernel of the map $\hat{u} := \varprojlim u_\bullet$. Indeed, if $M''_n := \mathrm{Im}(u_n)$, then the inverse system $M''_\bullet$ has surjective transition maps, although it might not form a $p$-adic $A$-module. We thus get a surjective map $\varprojlim N_\bullet \to \varprojlim \mathrm{Cok}(u_\bullet)$, although the map $\varprojlim M_\bullet \to \varprojlim \mathrm{Im}(u_\bullet)$ might not be surjective. In fact it is easy to verify that $\varprojlim \mathrm{Im}(u_\bullet)$ is the closure of $\mathrm{Im}(\hat{u})$ in $\varprojlim N_\bullet$. Thus $\varprojlim \mathrm{Cok}(u_\bullet)$ is the quotient of $\varprojlim N_\bullet$ by the closure of $\mathrm{Im}(\hat{u})$.

**Proposition 7.2.7.** *Let $A$ be a $p$-adically separated and complete ring.*

(1) *The inverse limit of a a strictly exact sequence $0 \to M'_\bullet \to M_\bullet \to M''_\bullet$ (resp. $0 \to M'_\bullet \to M_\bullet \to M''_\bullet \to 0$) of $p$-adic $A$-modules is an exact sequence of $A$-modules.*

(2) *The inverse limit of a strictly acyclic complex of $p$-adic $A$-modules is an acyclic complex of $A$-modules.*

(3) *The inverse limit of a strict quasi-isomorphism of $p$-adic $A$-modules is a quasi-isomorphism of complexes of $A$-modules.*

*Proof.* If each sequence $0 \to M'_n \to M_n \to M''_n$ is exact, then so so is the limit sequence, since the inverse limit functor is left exact. Surjectivity is also preserved on the right, since the transition maps in the kernel sequence are all surjective. This proves statement (1).

For (2), suppose that $K^{\bullet}_{\bullet}$ is a complex of $p$-adic $A$-modules and that each complex $K^{\bullet}_n$ is acyclic For each $q \in \mathbf{Z}$, let $\overline{K}^q_n$ be the image of $d^{q-1}_n$, so that we have exact sequences:

$$0 \to \overline{K}^q_n \to K^q_n \to \overline{K}^{q+1}_n \to 0$$

for each $n$. Since we have a surjection $K^{q-1}_{\bullet} \to \overline{K}^q_{\bullet}$ and the transition maps of $K^{q-1\bullet}$ are surjective, the same is true of the transition maps of $\overline{K}^q_{\bullet}$, and hence the sequence

$$0 \to \varprojlim \overline{K}^q_{\bullet} \to \varprojlim K^q_{\bullet} \to \varprojlim \overline{K}^{q+1}_{\bullet} \to 0$$

is also exact. This holds for all $q$, and so the map $\varprojlim K^{q-1}_{\bullet} \to \varprojlim \overline{K}^q_{\bullet}$ is also surjective. Thus $\varprojlim \overline{K}^q_{\bullet}$ is the image of the map $\varprojlim d^{q-1}_{\bullet}$. Since we also know that it is the kernel of the map $\varprojlim d^q_{\bullet}$, we can conclude that $\varprojlim K^{\bullet}_{\bullet}$ is exact. Since the cone of a map of complexes of $p$-adic $A$-modules is again a complex of $p$-adic $A$-modules, statement (3) follows. □

We should also mention the following simple version of the derived Nakayama Lemma (see [2, Tag 0G1U]).

**Proposition 7.2.8.** *Let $K^{\bullet}_{\bullet}$ be a complex of $p$-adic $A$-modules. Suppose that each $K^q_n$ is flat over $\mathbf{Z}/p^n\mathbf{Z}$. If $H^m(K^{\bullet}_1)$ vanishes for some $m$, then $H^m(K^{\bullet}_n)$ vanishes for every $n$.*

*Proof.* Since $K^q_1 \cong K^q_n/pK^q_n$ and $K^q_n$ is flat, we have an exact sequence of complexes:

$$0 \to K^{\bullet}_{n-1} \to K^{\bullet}_n \to K^{\bullet}_1 \to 0$$

Then the result follows by induction on $n$. □

## 7.3 $p$-complete flatness

Since the Artin-Rees lemma is not available in our context, the map from a module to its formal completion might not be flat, and as far as we know, even the map from a ring $A$ to the $p$-adic completion of a localization of $A$ might not be flat. The notion of "complete flatness" helps to overcome this difficulty. We add some details to the discussion in [10, Corollary 3.14].

**Definition 7.3.1.** *If $I$ is an ideal in a ring $A$, then an $A$-module $M$ is said to be $I$-completely flat if $\operatorname{Tor}^A_i(M, N)$ vanishes whenever $i > 0$ and $N$ is an $A$-module with $IN = 0$. An $A$-module is $I$-completely faithfully flat if it is $I$-completely flat and in addition $M/IM$ is faithfully flat as an $A/I$-module.*

In our case, the ideal $I$ will be principally generated by $p$, and the notion simplifies. In particular, if $A$ is $p$-torsion free and $p$-adically separated and complete, then a $p$-adically separated and complete $A$-module $M$ is $p$-completely flat if and only if the corresponding $p$-adic $A$-module $M_\bullet$ is flat as an inverse system of $A_\bullet$-modules.

**Proposition 7.3.2.** *Let $A$ be a $p$-torsion free ring and $M$ an $A$-module; write $M_n$ for $M/p^nM$ if $n \in \mathbf{N}$. Then conditions (1)–(3) below are equivalent and imply condition (4). If $M$ is $p$-adically separated, then all four conditions are equivalent.*

*(1) $M$ is $p$-completely flat (resp. faithfully flat).*

*(2) $M$ is $p^n$-completely flat (resp. faithfully flat) for all $n > 0$.*

*(3) $M$ is $p$-torsion free and $M/pM$ is flat (resp. faithfully flat) over $A/pA$.*

*(4) $M_n$ is flat (resp. faithfully flat) over $A_n$ for all $n > 0$.*

*Proof.* Since a flat $A_n$-module is faithfully flat if and only if its reduction modulo $p$ is faithfully flat as an $A_1$-module, it suffices to prove the equivalence for flatness.

We prove that (1) implies (2) by induction on $n$. Suppose that $M$ is $p^n$-completely flat and that $p^{n+1}N = 0$. Let $N' := \{x \in N : px = 0\}$ and let $N'' := N/N'$. Then $p^nN'' = p^nN' = 0$, and the exact sequence

$$\mathrm{Tor}_i^A(M, N') \to \mathrm{Tor}_i^A(M, N) \to \mathrm{Tor}_i^A(M, N'')$$

shows that $\mathrm{Tor}_i^A(M, N) = 0$ if $i > 0$, so $M$ is $p^{n+1}$-completely flat.

If (2) holds, let $F^\bullet \to M$ be a flat resolution of $M$. Condition (2) implies that $\mathrm{Tor}_i^A(M, A_n) = 0$ for $i > 0$, so $F_n^\bullet \to M_n$ is a flat resolution of the $A_n$-module $M_n$. If $N$ is an $A_n$-module, the isomorphism

$$\mathrm{Tor}_i^{A_n}(M_n, N) \cong H_i(F_n^\bullet \otimes_{A_n} N) \cong H_i(F^\bullet \otimes_A N) \cong \mathrm{Tor}_i^A(M, N) \qquad (7.3.1)$$

shows that $\mathrm{Tor}_i^{A_n}(M_n, N) \cong \mathrm{Tor}_i^A(M, N)$. Condition (2) implies that this vanishes if $i > 0$, and since $N$ is an arbitrary $A_n$-module, it follows that $M_n$ is $A_n$-flat. Thus (2) implies (4).

If (2) holds, then $\mathrm{Tor}_1^A(M, A/pA)$ vanishes. Then the $p$-torsion freeness of $A$ implies that $M$ is also $p$-torsion free, and since (4) also holds, we see that (2) implies (3).

To see that (3) implies (2), we again choose a flat resolution $F^\bullet$ of $M$. Condition (3) implies that $M$ is $p$-torsion free and hence that $\mathrm{Tor}_i^A(M, A/pA) = 0$ for all $i > 0$, so $F^\bullet \otimes_A A/pA$ is a flat resolution of $M/pM$. In particular, equation (7.3.1) holds with $n = 1$, *i.e.*, $\mathrm{Tor}_i^A(M, N) \cong \mathrm{Tor}_i^{A/pA}(M/pM, N)$ for every $A/pA$-module $N$. Since $M/pM$ is assumed to be flat over $A/pA$, this module vanishes for every $i > 0$, so $M$ is $p$-completely flat. Thus conditions (1)–(3) are equivalent.

It remains only to show that (4) implies (3) if $M$ is $p$-adically separated. The flatness of $M_n$ and $p$-torsion freeness of $A$ imply that the sequence

$$0 \longrightarrow M/pM \xrightarrow{p^{n-1}} M_n \xrightarrow{p} M_n$$

is exact and hence that $\varprojlim M_n$ is $p$-torsion free. (See the proof of Lemma 7.3.4 below.) Thus any $p$-torsion element of $M$ maps to zero in this limit and necessarily vanishes if $M$ is $p$-adically separated. □

**Proposition 7.3.3.** *Suppose that $A$ is $p$-torsion free and $p$-adically separated and complete and that $M$ is a $p$-completely flat $A$-module.*

*(1) If $N$ is a $p$-torsion free $A$-module, then $N \otimes_A M$ and $N \hat{\otimes}_A M$ are also $p$-torsion free.*

*(2) If $N$ is a $p$-completely flat $A$-module, then $N \hat{\otimes}_A M$ is also $p$-completely flat.*

*(3) If $A \to A'$ is a homomorphism of $p$-torsion free $p$-adically separated and complete rings, then $A' \hat{\otimes}_A M$ is $p$-completely flat over $A'$, and is faithfully flat if $M$ is $p$-completely faithfully flat over $A$.*

*(4) The composition of $p$-completely flat (resp. faithfully flat) homomorphisms is again $p$-completely flat (resp. faithfully flat). If $A \to A'$ is $p$-completely faithfully flat, then a $p$-adically separated and complete $A$-module $N$ is $p$-completely flat (resp. faithfully flat) over $A$ if and only if $A' \hat{\otimes}_A N$ is so over $A'$.*

*Proof.* The short exact sequence $0 \longrightarrow N \xrightarrow{p} N \longrightarrow N/pN \longrightarrow 0$ gives rise to a long one:

$$\operatorname{Tor}^1_A(N/pN, M) \longrightarrow N \otimes_A M \xrightarrow{p} N \otimes_A M \longrightarrow N/pN \otimes_A M \longrightarrow 0.$$

Since $M$ is $p$-completely flat, the term on the left vanishes. This implies that $N \otimes M$ is $p$-torsion free, and the result for $N \hat{\otimes} M$ follows from Lemma 7.3.4 below.

If $N$ is $p$-completely flat, then it is $p$-torsion free, and $N \hat{\otimes} M$ is $p$-torsion free by (1). Since $N_1$ and $M_1$ are flat over $A_1$, so is their tensor product. By (2) of Proposition 7.1.2, this is the same as $(N \hat{\otimes} M)_1$, so $N \hat{\otimes} M$ is $p$-completely flat by criterion (3) of Proposition 7.3.2.

To prove (3), first observe that statement (1) implies that $M' := A' \hat{\otimes}_A M$ is $p$-torsion free. Furthermore, $M'_1 \cong A'_1 \otimes_{A_1} M_1$, which is flat (resp. faithfully flat) over $A'_1$ because $M_1$ is flat (resp. faithfully flat) over $A_1$. Then criterion (3) of Proposition 7.3.2 implies that $M'$ is $p$-completely flat (resp. faithfully flat) over $A'$.

Statement (4) follows immediately from criterion (4) of Proposition 7.3.2 and the analogous results for usual flatness. □

**Lemma 7.3.4.** *The $p$-adic completion of a $p$-torsion free abelian group is again $p$-torsion free.*

*Proof.* If $M$ is $p$-torsion free, then for every $n$ we have a commutative diagram with exact rows:

$$\begin{array}{ccccccc}
0 & \longrightarrow & M/pM & \xrightarrow{p^n} & M/p^{n+1}M & \xrightarrow{p} & M/p^{n+1}M \\
 & & \downarrow{\scriptstyle p} & & \downarrow{\scriptstyle \pi} & & \downarrow{\scriptstyle \pi} \\
0 & \longrightarrow & M/pM & \xrightarrow{p^{n-1}} & M/p^nM & \xrightarrow{p} & M/p^nM,
\end{array}$$

where $\pi$ is the canonical projection. Since $\varprojlim$ is left exact, we find that the sequence

$$0 \longrightarrow \varprojlim(M/pM, p) \longrightarrow \varprojlim(M/p^nM, \pi) \xrightarrow{p} \varprojlim(M/p^nM, \pi)$$

is exact. The inverse limit on the left is zero, and the result follows. □

We now investigate the relationship beween $p$-complete flatness and strict exactness. Suppose $A$ is a $p$-adically separated and complete ring. If $N_\bullet$ is a $p$-adic $A$-module with inverse limit $\hat{N}$ and $M$ is any $A$-module, then $\{N_n \otimes_A M \cong (\hat{N} \otimes_A M)_n\}$ forms a $p$-adic $A$-module whose inverse limit is $\hat{N}\hat{\otimes}_A M$. Furthermore $(\hat{N} \otimes M)_n \cong N_n \otimes_A M_n$ for every $n$, by statement (2) of Proposition 7.1.2.

**Proposition 7.3.5.** *Let $K^\bullet$ be a complex of $p$-adically separated and complete $A$-modules and let $M$ be a $p$-completely flat $A$-module.*

*(1) If $K^\bullet$ is strictly acyclic, then so is $K^\bullet\hat{\otimes}_A M$.*

*(2) If $M$ is $p$-completely faithfully flat and $K^\bullet\hat{\otimes}_A M$ is strictly acyclic, then so is $K^\bullet$.*

*(3) If $M$ is $p$-completely faithfully flat, then a homomorphism of complexes $u$ of $p$-adically separated and complete modules is a strict quasi-isomorphism if and only if $u\hat{\otimes}\mathrm{id}_M$ is.*

*In particular, if $\Sigma$ is a strictly short exact sequence of $p$-adically separated and complete $A$-modules, then $\Sigma\hat{\otimes}_A M$ is again strictly short exact.*

*Proof.* If each $K_n^\bullet$ is acyclic, then, since each $M_n$ is flat over $A_n$, so is each $K_n^\bullet \otimes_{A_n} M_n$. It follows from Proposition 7.2.7 that the inverse limit is again acyclic. If now $M$ is $p$-completely faithfully flat, then each $M_n$ is a faithfully flat $A_n$-module, and if each $K_n^\bullet \otimes_A M$ is acyclic, then so is each $K_n^\bullet$. The conclusion again follows from Proposition 7.2.7. Statement (3) follows by taking the mapping cone of $u$, which is again a complex of $p$-adically separated and complete modules. □

**Corollary 7.3.6.** *Suppose that $A$ is a $p$-adically separated and complete ring and $M$ is a $p$-completely faithfully flat $A$-module.*

*(1) Let $N'^{\bullet} \to N^{\bullet}$ be a homomorphism of complexes of $p$-adically separated and complete $A$-modules. Each $N'^{\bullet}_n \otimes_A M \to N^{\bullet}_n \otimes_A M$ is a quasi isomorphism if and only if each each $N'^{\bullet}_n \to N^{\bullet}_n$ is. If this is true for all $n$, then $N'^{\bullet} \to N^{\bullet}$ and $N'^{\bullet}\hat{\otimes}_A M \to N^{\bullet}\hat{\otimes}_A M$ are quasi-isomorphisms.*

*(2) If $N' \to N$ is a homomorphism of $p$-adically separated and complete $A$-modules, then $N' \to N$ is an isomorphism if and only if $N'\hat{\otimes}_A M \to N\hat{\otimes}_A M$ is an isomorphism.*

*Proof.* Statement (1) follows from Proposition 7.3.5 applied to the cone of $N'^{\bullet} \to N^{\bullet}$. Statement (2) follows from (1), or directly from Proposition 7.3.5. because $N' \to N$ is an isomorphism if and only if each $N'_n \to N_n$ is an isomorphism □

**Corollary 7.3.7.** *Suppose $A$ is a $p$-adically separated and complete $p$-torsion free ring. Let $E_{\text{tor}}$ be the $p$-torsion submodule of a $p$-adically separated and complete $A$-module $E$ and let $E^-_{\text{tor}}$ be its closure in $E$. The quotient $E/E^-_{\text{tor}}$ is $p$-torsion free and the inclusion $E^-_{\text{tor}} \to E$ is a strict monomorphism. Furthermore, if $M$ is $p$-completely flat, the natural maps $E^-_{\text{tor}}\hat{\otimes}_A M \to (E\hat{\otimes}_A M)^-_{\text{tor}}$ and $(E/E^-_{\text{tor}})\hat{\otimes}_A M \to (E\hat{\otimes}_A M)/(E\hat{\otimes}_A M_{\text{tor}})^-$ are isomorphisms.*

*Proof.* Since $E/E^-_{\text{tor}}$ is the completion of $E/E_{\text{tor}}$ (see Lemma 7.2.4) and the latter is $p$-torsion free, Lemma 7.3.4 implies that $E/E^-_{\text{tor}}$ is also $p$-torsion free. This implies that $p^n E \cap E^-_{\text{tor}} = p^n E^-_{\text{tor}}$ for all $n$ and hence the strictness of the map $E^-_{\text{tor}} \to E$. Thus the sequence $0 \to E^-_{\text{tor}} \to E \to E/E^-_{\text{tor}} \to 0$ is strictly short exact, and remains so after forming the completed tensor product with $M$, by Proposition 7.3.5. This gives us the exactness of the top row of the following diagram:

$$\begin{array}{ccccccccc}
0 & \longrightarrow & E^-_{\text{tor}}\hat{\otimes}_A M & \longrightarrow & E\hat{\otimes}_A M & \longrightarrow & (E/E^-_{\text{tor}})\hat{\otimes}_A M & \longrightarrow & 0 \\
 & & \downarrow & & \downarrow & & \vdots & & \\
0 & \longrightarrow & (E\hat{\otimes}_A M)^-_{\text{tor}} & \longrightarrow & E\hat{\otimes}_A M & \longrightarrow & (E\hat{\otimes}_A M)/(E\hat{\otimes}_A M)^-_{\text{tor}} & \longrightarrow & 0
\end{array}$$

It is clear that $E_{\text{tor}}\hat{\otimes}_A M$ maps to $(E\hat{\otimes}_A M)_{\text{tor}}$, and hence the same is true for the corresponding closures. This gives the existence of the dashed arrow (the "natural map") in the diagram. The module $(E/E^-_{\text{tor}})\hat{\otimes}_A M$ is $p$-adically separated and complete, and it is $p$-torsion free, by Proposition 7.3.3. Thus the map $E\hat{\otimes}_A M \to (E/E^-_{\text{tor}})\hat{\otimes}_A M$ factors through $(E\hat{\otimes}_A M)/(E\hat{\otimes}_A M)^-_{\text{tor}}$, and so the right vertical arrow is an isomorphism. It follows that the left vertical arrow is also an isomorphism. □

**Proposition 7.3.8.** *If $A \to B$ be a $p$-completely faithfully flat homomorphism of $p$-adically separated and complete $p$-torsion free rings, let $C_B^\bullet$ be the completed augmented Čech-Alexander complex*

$$C_B^\bullet := A \to B \to B\hat{\otimes}_A B \to B\hat{\otimes}_A B\hat{\otimes}_A B \cdots .$$

*(1) If $N$ is a $p$-adically separated and complete $A$-module, then the complex $C_B^\bullet \hat{\otimes} N$, is strictly acyclic, i.e., $C_B^\bullet \otimes N_n$ is acyclic for every $n$.*

*(2) If $N$ is a $p$-adically separated and complete $B$-module endowed with formal descent data $\epsilon\colon (B\hat{\otimes}_A B)\hat{\otimes}_B N \to N\hat{\otimes}_B(B\hat{\otimes}_A B)$, let*

$$\begin{array}{rcl} N^\epsilon & := & \{x \in N : \epsilon(1\hat{\otimes}x) = x\hat{\otimes}1\} \\ (N_n)^\epsilon & := & \{x \in N_n : \epsilon_n(1\hat{\otimes}x) = x\hat{\otimes}1\}. \end{array}$$

*Then:*

*(a) $(N^\epsilon)_n = (N_n)^\epsilon$ for all $n$, and we will denote both by $N_n^\epsilon$. Thus $N^\epsilon \to N$ is a strict monomorphism.*

*(b) The natural maps $B\hat{\otimes}_A N^\epsilon \to N$ and $B \otimes_A N_n^\epsilon \to N_n$ are isomorphisms for all $n$.*

*Proof.* Statement (2) of Proposition 7.1.2 implies that

$$(C_B\hat{\otimes}_A N)_n \cong C_{B_n} \otimes_{A_n} N_n$$

for every $n$. Since each $A_n \to B_n$ is faithfully flat, the acyclicity of each of these complexes is standard. Thus $(C_B\hat{\otimes}_A N)_\bullet$ forms a strictly acyclic complex of $p$-adic $B$-modules, and so Proposition 7.2.7 implies that the limit sequence is also acyclic. This proves (1).

For (2), we start with the fact that each $B\otimes_A (N_n)^\epsilon \to N_n$ is an isomorphism, by usual faithfully flat descent for the homomorphism $A_n \to B_n$. Since the maps $N_n \otimes_{B_n} B_{n-1} \to N_{n-1}$ are isomorphisms, the same is true for the maps $B \otimes_A (N_n)^\epsilon \otimes_{B_n} B_{n-1} \to B \otimes_A (N_{n-1})^\epsilon$. Again invoking faithfully flat descent, we conclude that the maps $(N_n)^\epsilon \otimes_{A_n} A_{n-1} \to (N_{n-1})^\epsilon$ are isomorphisms. Thus $(N_\bullet)^\epsilon$ forms a $p$-adic $A$-module. Then (2) of Proposition 7.1.2 implies that each map $(\varprojlim (N_n)^\epsilon) \otimes_A A_n \to N_n^\epsilon$ is an isomorphism. The left exactness of $\varprojlim$ implies that the natural map $N^\epsilon \to \varprojlim (N_n)^\epsilon$ is an isomorphism, and (2a) follows. We have already seen that each map $B \otimes_A N_n^\epsilon \to N_n$ is an isomorphism, and hence so is the map

$$\varprojlim (B \otimes_A N_n^\epsilon) \to \varprojlim N_n = N.$$

Since $\varprojlim (B \otimes_A N_n^\epsilon) \cong B\hat{\otimes}_A N^\epsilon$, statement (2b) follows. □

A morphism of $p$-torsion free $p$-adic formal schemes $S' \to S$ is said to be $p$-completely flat if for every open affine $U' \subseteq S'$ mapping to an open affine $U \subseteq S$, the corresponding map $\mathcal{O}_S(U) \to \mathcal{O}_{S'}(U')$ is $p$-completely flat. It follows from Proposition 7.3.3 that a $p$-completely flat map $A \to A'$ gives rise to a $p$-completely flat map $\operatorname{Spf} A' \to \operatorname{Spf} A$ and that the family of $p$-completely faithfully flat maps $S'' \to S' \to S$ forms a covering family for a site which we denote by $S_{pcf}$. The objects of (the big version of) $S_{pcf}$ are morphisms of $p$-adic formal schemes $S' \to S$; the small version restricts to such morphisms which are $p$-completely flat. We note that the set of $p$-completely flat maps $U \to S'$ with $U$ affine forms a base for this topology.

## 7.4 $p$-complete quasi-coherence

We now briefly discuss how to work with the notions $p$-adic modules and $p$-complete flatness for sheaves. The equivalence of the two conditions in the following definition is a consequence of Proposition 7.1.2.

**Definition 7.4.1.** *Let $S$ be a $p$-torsion free $p$-adic formal scheme. A sheaf of $\mathcal{O}_S$-modules $E$ on $S_{pcf}$ is $p$-completely quasi-coherent if its value on each affine is $p$-adically separated and complete and the following equivalent conditions are satisfied.*

*(1) For every $p$-completely flat $S$-morphism $U' \to U$ of affine $p$-adic formal schemes which are $p$-completely flat over $S$, the map*

$$\mathcal{O}_S(U')\hat{\otimes}_{\mathcal{O}_S(U)}E(U) \to E(U')$$

*is an isomorphism.*

*(2) For every $p$-completely flat $S$-morphism $U' \to U$ of affine $p$-adic formal schemes which are $p$-completely flat over $S$, the map*

$$\mathcal{O}_{S_n}(U') \otimes_{\mathcal{O}_{S_n}(U)} E_n(U) \to E_n(U')$$

*is an isomorphism for every $n$.*

*In particular, the structure sheaf $\mathcal{O}_S$ is itself $p$-completely quasi-coherent, and if $E$ is $p$-completely quasi-coherent, then each $E_n$ is quasi-coherent as sheaf of $\mathcal{O}_{S_n}$-modules on the scheme $S_n$.*

The following results follow from the discussions above by standard arguments, which we leave to the reader.

**Proposition 7.4.2.** *Let $S$ be a $p$-torsion free $p$-adic formal scheme.*

*(1) If $E$ is a sheaf of $\mathcal{O}_S$-modules and there is a cover $\{U_i \to S\}$ of $S_{pcf}$ such that each $E_{|U_i}$ is $p$-completely quasi-coherent, then $E$ is $p$-completely quasi-coherent.*

*(2) Suppose that $S = \operatorname{Spf} A$ and $M$ is a $p$-adically separated and complete $A$-module. Define a presheaf $\tilde{M}$ on the family of affine $p$-completely flat maps $\operatorname{Spf} A' \to \operatorname{Spf} A$ by $\tilde{M}(\operatorname{Spf} A') := A'\hat{\otimes}_A M$. Then in fact $\tilde{M}$ is a sheaf, and the map $M \to \tilde{M}(S)$ is an isomorphism. Moreover, the functor $M \to \tilde{M}$ defines an equivalence from the category of $p$-adically separated and complete $A$-modules to the category of $p$-complete quasi-coherent sheaves of $\mathcal{O}_S$-modules.*

(3) *If $S$ is affine and $E$ is a $p$-completely quasi-coherent sheaf of $\mathcal{O}_S$-modules on $S_{pcf}$, then $H^i(S_{pcf}, E)$ vanishes for $i > 0$.* $\square$

**Remark 7.4.3.** It will sometimes be convenient for us to view $p$-completely quasi-coherent sheaves geometrically. We do this by adapting Grothendieck's $\mathbf{V}$ construction. Namely, if $Y$ is a $p$-torsion free x$p$-adic formal scheme and $E$ is a $p$-completely quasi-coherent sheaf of $\mathcal{O}_Y$-modules, then for each $n \in \mathbf{N}$, the sheaf $E_n$ is quasi-coherent on $Y_n$ and we can form the scheme $\mathbf{V}E_n := \mathrm{Spec}(S^\bullet E_n)$, affine over $Y_n$. Since formation of $S^\bullet$ is compatible with base change, the system $S^\bullet E_\bullet$ forms a $p$-adic module, and $\varprojlim S^\bullet E_\bullet$ is $p$-adically separated and complete, by Proposition 7.1.2. Thus we can view the collection of $Y$-schemes $\{\mathbf{V}E_n : n \in \mathbf{N}\}$ as the family of reductions of a $p$-adic formal scheme $\mathbf{V}E$, affine over $Y$.

## 7.5 Regular sequences

Various notions of "regular sequences" are used in the literature and in particular in [10]. As is well known, all these notions are equivalent in the noetherian case. We have chosen here to restrict our attention to the standard definition.

**Proposition 7.5.1.** *Let $X \to Y$ be a closed immersion of $p$-torsion free $p$-adic formal schemes. If $\overline{X} \to \overline{Y}$ is a regular immersion, then the same is true of $X \to Y$.*

*Proof.* To check this, we work locally, assuming that $Y = \mathrm{Spf}\, B$, that $X = \mathrm{Spf}\, C$, and that the ideal $I_{\overline{X}/\overline{Y}}$ of $\overline{X}$ in $\overline{Y}$ is generated by a $\overline{B}$-regular sequence. Since $X$ and $Y$ are $p$-torsion free, $I_{\overline{X}/\overline{Y}} = I_{X/Y}/(I_{X/Y} \cap p\mathcal{O}_Y) = I_{X/Y}/pI_{X/Y}$. Choose a sequence $(b_1, \ldots, b_r)$ of elements of $I_{X/Y}$ lifting a $\overline{B}$-regular sequence generating $I_{\overline{X}/\overline{Y}}$. Then $(p, b_1, \ldots, b_r)$ is $B$-regular, and Lemma 7.5.2 below implies that the sequence $(b_1, \ldots, b_r, p)$ is also $B$-regular, and that $C' := B/(b_1, \ldots, b_r)$ is $p$-adically separated and complete and $p$-torsion free. We claim that the surjection $\theta\colon C' \to C$ is an isomorphism. We know that $\overline{\theta}$ is an isomorphism, so if $x$ belongs to the kernel of $\theta$, necessarily $x = px'$ for some $x' \in C'$. Since $C$ is also $p$-torsion free, it follows that $x'$ is also in the kernel of $\theta$, and we see that in fact $x \in p^nC'$ for every $n$. Since $C'$ is $p$-adically separated, this implies that $x = 0$. $\square$

**Lemma 7.5.2.** *Let $B$ be a $p$-torsion free $p$-adically separated and complete ring, let $(b_1, \ldots, b_r)$ be a sequence in $B$ such that $(p, b_1, \ldots, b_r)$ is $B$-regular. Then $(b_1, \ldots, b_r, p)$ is also $B$-regular, and furthermore $B/(b_1, \ldots, b_r)$ is also $p$-torsion free, separated, and complete.*

*Proof.* We first show that $b_1$ is a nonzero divisor in $B$. Namely, if $b_1x$ vanishes in $B$, then it also vanishes in $B/pB$, and since multiplication by $b_1$ on $B/pB$ is injective, it follows that $x = px_1$ for some $x_1 \in B$. Then $pb_1x_1 = b_1px_1 = b_1x = 0$, and since $B$ is $p$-torsion free, this implies that $b_1x_1 = 0$. Repeating this process, we see that $x_1 = px_2$, that $x_2 = px_3$, and so on. Since $B$ is $p$-adically separated, we conclude that $x = 0$. Next we check that $B' := B/b_1B$ is $p$-torsion free. Namely, if $px = b_1y$, then since $b_1$ is a nonzero divisor in $B/pB$, we can write $y = py'$ for some $y'$. Then $px = b_1py'$, hence $x = b_1y'$, so $x$ also maps to zero in $B'$.

It remains to check that $B'$ is also $p$-adically separated and complete. (This is evident in the noetherian case, but requires an argument in general.) Note that for each $n$, the sequence

$$0 \longrightarrow B_n \xrightarrow{b_1} B_n \longrightarrow B'_n \longrightarrow 0$$

is exact, because $B'$ is $p$-torsion free. Then we find a commutative diagram with exact rows:

$$\begin{array}{ccccccccc} 0 & \longrightarrow & B & \xrightarrow{b_1} & B & \longrightarrow & B' & \longrightarrow & 0 \\ & & \downarrow & & \downarrow & & \downarrow & & \\ 0 & \longrightarrow & \varprojlim B_n & \xrightarrow{b_1} & \varprojlim B_n & \longrightarrow & \varprojlim B'_n & \longrightarrow & 0. \end{array}$$

The bottom row is exact because all the transition maps are surjective. Since the first two vertical arrows are isomorphisms, so is the third.

Having proved the proposition when $r = 1$, we proceed by induction on $r$. We have seen that the ring $B/b_1B$ is $p$-adically separated and complete and $p$-torsion free. Since $(b_1, \ldots, b_r)$ was assumed to be $B/pB$-regular, the sequence $(b_2, \ldots, b_r)$ is $B/(p, b_1)$-regular. Thus $(p, b_2, \ldots, b_r)$ is $B/b_1B$-regular. Then the induction hypothesis implies that $(b_2, \ldots, b_r, p)$ is $B/b_1B$-regular, and hence that $(b_1, \ldots, b_r)$ is $B$-regular. The induction hypothesis also implies that $B/(b_1, b_2, \ldots, b_r)$ is $p$-torsion free and separated and complete. □

The following result is due to Micali [27, Theorem 1], with an alternative proofs given by Nowak[29] and C. Huneke [19, 3.1] (with a slightly different hypothesis). Special thanks to Huneke for helpful clarifications about this result. Here is yet another.

**Proposition 7.5.3.** *Let $B$ be a commutative ring and let $(b_1, \ldots, b_r)$ be a $B$-regular sequence generating an ideal $I$. Then for all $n$, the natural map*

$$m_n \colon S^n_B(I) \to I^n$$

*is an isomorphism.*

*Proof.* If $r = 1$, then $I$ is free of rank one, generated by $b_1$, and the result is clear. Let us assume that $r > 1$ and that the result holds for $r-1$. When $n = 1$, there is nothing to prove, and we proceed by induction on $n$. In the ensuing calculations, we use $\odot$ to indicate multiplication in the symmetric algebra $S^\bullet_B(I)$.

For each $n \geq 1$ there is a multilinear function $f_i$ from the $n$th product $I \times I \times \cdots \times I$ to $S_B^{n-1}(I)$ sending $(x_1, \ldots, x_n)$ to $x_i(x_1 \odot x_2 \odot \cdots \odot \hat{x}_i \odot \cdots \odot x_n)$. If $j < i$, we have:

$$\begin{aligned} f_i(x_1, \ldots, x_n) &= x_i(x_1 \odot x_2 \odot \cdots \odot x_j \odot \cdots \odot \hat{x}_i \odot \cdots \odot x_n) \\ &= (x_1 \odot x_2 \odot \cdots \odot x_i x_j \odot \cdots \odot \hat{x}_i \odot \cdots \odot x_n) \\ &= x_j(x_1 \odot x_2 \odot \cdots \odot x_i \odot \cdots \odot \hat{x}_i \odot \cdots \odot x_n) \\ &= x_j(x_1 \odot x_2 \odot \cdots \odot \hat{x}_j \odot \cdots \odot x_i \odot \cdots \odot x_n) \\ &= f_j(x_1, \ldots, x_n). \end{aligned}$$

It follows that $f_i(x_1, \ldots, x_n)$ is unchanged under every transposition, and hence every permutation, of the sequence $(x_1, \ldots, x_n)$. Since $f_i$ is multilinear, it factors through a homomorphism $f\colon S_B^n(I) \to S_B^{n-1}(I)$.

Let $\overline{B} := B/b_1B$, let $\overline{I} := I/(b_1B)$, and let $\overline{b}_i$ be the image of $b_i$ in $\overline{B}$. There is a diagram

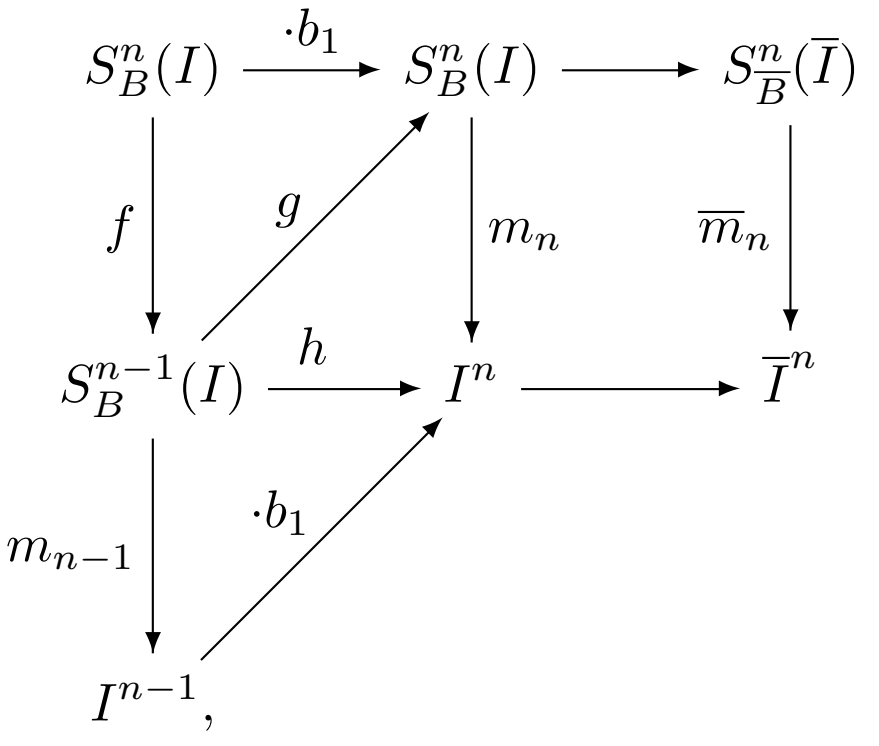

where

$$g(y_1 \odot y_2 \odot \cdots \odot y_{n-1}) := (y_1 \odot y_2 \odot \cdots \odot y_{n-1} \odot b_1)$$

and

$$h(y_1 \odot y_2 \cdots \odot y_{n-1}) := y_1 y_2 \cdots y_{n-1} b,$$

and the remaining maps are given by multiplication, inclusion, and projection. The upper triangle commutes because

$$\begin{aligned} g(f(x_1 \odot x_2 \odot \cdots \odot x_n)) &= g(x_1(x_2 \odot \cdots \odot x_n)) \\ &= x_1 g(x_2 \odot \cdots \odot x_n) \\ &= x_1(x_2 \odot \cdots \odot x_n \odot b_1) \\ &= (x_2 \odot \cdots \odot x_n \odot x_1 b_1) \\ &= b_1(x_2 \odot \cdots \odot x_n \odot x_1) \\ &= b_1(x_1 \odot x_2 \cdots \odot x_n). \end{aligned}$$

The bottom triangle commutes by the definition of $h$, and the central triangle commutes because

$$m_n(g(y_1 \odot \cdots \odot y_{n-1})) = m_n(y_1 \odot \cdots \odot y_{n-1} \odot b_1)$$

$$= y_1 \cdots y_{n-1} b_1$$
$$= h(y_1 \odot \cdots \odot y_{n-1}).$$

The sequence $(\bar{b}_2, \cdots, \bar{b}_r)$ is $\overline{B}$-regular, so the map $\overline{m}_n$ is injective by our induction assumption on $r$, and the map $m_{n-1}$ is injective by our induction assumption on $n$. Since $m_{n-1}$ and $\cdot b_1$ (the bottom slanted arrow) are both injective, $h$ is also injective. The top sequence is right exact because $S^n_{\overline{B}}(\bar{I}) \cong S^n_B(I) \otimes B/(b_1 B)$. To see that $m_n$ is injective, suppose that $x \in S^n_B(I)$ is in its kernel. Since $\overline{m}_n$ is injective and the top sequence is exact in the middle, $x = b_1 y$ for some $y \in S^n_B(I)$. Then $h(f(y)) = m_n(b_1 y) = m_n(x) = 0$, and since $h$ is injective, $f(y) = 0$. But then $x = b_1 y = g(f(y))$ also vanishes. □

# References

[1] Johannes Anschütz and Arthur-César Le Bras. Prismatic Diuedonné theory. *Forum of Mathematics, Pi*, 11:e2:1–92, 2023.

[2] The Stacks Project Authors. *The Stacks Project*. `http://stacks.math.columbia.edu`, 2023.

[3] Pierre Berthelot. Le Théorème de Riemann-Roch. In *SGA 6, Exp. VIII*, volume 225 of *Lecture Notes in Mathematics*. Springer-Verlag, New York, 1971.

[4] Pierre Berthelot. *Cohomologie Cristalline des Schémas de Caractéristique* $p > 0$, volume 407 of *Lecture Notes in Mathematics*. Springer-Verlag, New York, 1974.

[5] Pierre Berthelot. Finitude et pureté cohomologique en cohomologie rigide, avec un appendice par Aise Johan de Jong. *Inv. Math.*, 128:329–377, 1997.

[6] Pierre Berthelot. $\mathcal{D}$-modules arithmétiques II. Descente par Frobenius. *Mém. Soc. Math. France*, 81, 2000.

[7] Pierre Berthelot and Arthur Ogus. *Notes on Crystalline Cohomology*, volume 21 of *Annals of Mathematics Studies*. Princeton University Press, Princeton, 1978.

[8] Bhargav Bhatt and Jacob Lurie. Absolute prismatic cohomology. arXiv:2202.06120v1, Jan 2022.

[9] Bhargav Bhatt, Jacob Lurie, and Akhil Matthew. Revisiting the de Rham Witt complex. *Astérisque*, 424, 2021. arXiv:1804.05501v1.

[10] Bhargav Bhatt and Peter Scholze. Prisms and prismatic cohomology. *Annals of Mathematics*, 196(3):1135–1275, 2022.

[11] James Borger. The basic geometry of Witt vectors, I the affine case. *Algebra and Number Theory*, 2011.

[12] James Borger. The basic geometry of Witt vectors II: Spaces. *Math. Ann.*, 351:877–933, 2011.

[13] Pierre Deligne. Théorie de Hodge II. *Publications Mathématiques de l'I.H.É.S.*, 40:5–57, 1972.

[14] Pierre Deligne and Luc Illusie. Relèvements modulo $p^2$ et décomposition du complexe de de Rham. *Inventiones Mathematicae*, 89:247–270, 1987.

[15] Vladimir Drinfeld. Prismatization. arXiv:2005.04746v5, Jan. 2022.

[16] Vladimir Drinfeld. A stacky approach to crystals. In *Dialogues between physics and mathematics—C. N. Yang at 100*, pages 19–47. Springer, Cham, [2022] ©2022.

[17] Gerd Faltings. Crystalline cohomology and $p$-adic galois representations. In Jun-Ichi Igusa, editor, *Algebraic Analysis, Geometry, and Number Theory*, pages 25–80, Baltimore and London, 1989. The Johns Hopkins University Press.

[18] Robin Hartshorne. On the de Rham cohomology of algebraic varieties. *Publications Mathématiques de l.I.H.É.S.*, 45:1–99, 1976.

[19] Craig Huneke. On the symmetric and Rees algebra of an ideal generated by a $d$-sequence. *Journal of Algebra*, 62:268–275, 1980.

[20] Luc Illusie. Report on crystallline cohomology. In Robin Hartshorne, editor, *Algebraic Geometry Arcata 1974*, volume XXIX of *Proceedings of Symposia in Pure Mathematics*, pages 459–478, Providence, Rhode Island, 1975. American Mathematical Society.

[21] Luc Illusie. Complexe de de Rham Witt et cohomologie cristalline. *Ann. Math. E.N.S.*, 12:501–601, 1979.

[22] André Joyal. $\delta$-anneaux et $\lambda$-anneaux. *C. R. Math. Rep. Acad. Sci. Canada*, 7(4):227–232, 1985.

[23] Nicholas Katz. Algebraic solutions of differential equations ($p$-curvature and the Hodge filtration). *Inventiones Mathematicae*, 18:1–118, 1972.

[24] Nicholas Katz and Tadoa Oda. On the differentiation of de Rham cohomology classes with respect to parameters. *Journal of Mathematics, Kyoto University*, 8:199–213, 1968.

[25] Guitang Lan, Mao Sheng, and Kang Zuo. Semistable Higgs bundles, periodic Higgs bundles and representations of algebraic fundamental groups. *J. Eur. Math. Soc.*, 21:3053–3112, 2019.

[26] Barry Mazur. Frobenius and the Hodge filtration—estimates. *Annals of Mathematics*, 98:58–95, 1973.

[27] Artibano Micali. Sur les algèbres universelles. *Ann. Inst. Fourier Grenoble*, 14,2:33–88, 1964.

[28] Matthew Morrow and Takeshi Tsuji. Generalised representations as $q$-conenections in integral $p$-adic Hodge theory. arXiv:2010.04059v1, 2020.

[29] Krzysztof Jan Nowak. A homological proof of Micali's theorem on the Rees algebra of a regular embedding. *Universitatis Iagellonicae Acta Mathematica*, XXXVII, 1999.

[30] Arthur Ogus. Local cohomological dimension of algebraic varieties. *Annals of Mathematics*, 98:327–365, 1973.

[31] Arthur Ogus. F-isocrystals and de Rham cohomology II—convergent isocrystals. *Duke Mathematical Journal*, 51(4):765–850, 1984.

[32] Arthur Ogus. The convergent topos in characteristic $p$. In *The Grothendieck Festschrift*, volume III, pages 133–162. Birkhauser, 1990.

[33] Arthur Ogus. Higgs cohomology, $p$-curvature, and the Cartier isomorphism. *Compositio Math.*, 140:145–164, 2004.

[34] Arthur Ogus and Vadim Vologodsky. Nonabelian Hodge theory in characteristic $p$. *Publications Mathématiques de l'I.H.E.S.*, 106:1–138, November 2007.

[35] Hidetoshi Oyama. PD Higgs crystals and Higgs cohomolology in characteristic $p$. *Journal of Algebraic Geometry*, 26:735–802, 2917.

[36] Norbert Roby. Lois polynômes et lois formelles en théorie des modules. *Annales Scientifiques de l'École Normale Supérieure 3é série*, 80:213–348, 1963.

[37] Atsushi Shiho. Notes on generalizations of local Ogus-Vologodsky correspondence. *J. Math. Sci. Univ. Tokyo*, 22:798–875, 2015.

[38] Carlos Simpson. Higgs bundles and local systems. *Publ. Math. de l'I.H.É.S.*, 75:5–95, 1992.

[39] Yichao Tian. Finiteness and duality of the cohomology of prismatic crystals. *Journal für die Reine und Angewandte Mathematik (Crelle)*, 2023(800):217–257, 2023.

[40] Daxin Xu. Lifting the Cartier transform of Ogus-Volodgodsky modulo $p^n$. *Memoires de la SMF*, 2019.